\documentclass[english,12pt]{article}   
\usepackage{babel} 
\usepackage{anyfontsize}
\usepackage{colonequals}
\usepackage[compact]{titlesec} 
\usepackage{breakcites}
\usepackage{microtype}
\usepackage{amsthm,amsmath,amssymb,mathabx}
\usepackage{import,enumitem,breqn}
\usepackage{graphicx,float,mathtools,bookmark}
\usepackage{tikz} 
\usepackage{ytableau} 
\usepackage[utf8]{inputenc}
\usepackage[T1]{fontenc}
\usepackage{mathptmx}
\usepackage{autonum}
\usepackage{hyperref}
\usepackage{geometry}
\usepackage{cases}
\usepackage{stmaryrd}

\DeclarePairedDelimiter{\ceil}{\lceil}{\rceil} 
\DeclarePairedDelimiter{\floor}{\lfloor}{\rfloor}
\newtheorem{theorem}{Theorem}[section]
\newtheorem{corollary}[theorem]{Corollary}
\newtheorem{lemma}[theorem]{Lemma}
\newtheorem{proposition}[theorem]{Proposition}
\newtheorem{conjecture}[theorem]{Conjecture}
\newtheorem{definition}[theorem]{Definition}
\theoremstyle{remark}
\theoremstyle{definition} 
\newtheorem{example}[theorem]{Example}
\newtheorem{remark}[theorem]{Remark}
\numberwithin{equation}{subsection}

\renewcommand{\i}{{\mathrm{i}}}
\newcommand{\Ir}{\operatorname{\mathrm{Ir}}}
\newcommand{\AIr}{\operatorname{\mathrm{AIr}}}
\newcommand{\tdet}{\operatorname{\mathrm{\det}}}
\newcommand{\Qu}{\operatorname{\mathrm{Qt}}}

\newcommand{\ct}{{\mathrm{C}}}
\newcommand{\sh}{{\mathrm{S}}}
\newcommand{\del}{{\delta}}
\newcommand{\Del}{{\Delta}}
\newcommand{\rc}{{\circ}}
\newcommand{\fraa}{\frac{1}{2}}

\newcommand{\len}{\mathrm{L}}
\newcommand{\num}{\mathrm{\#}}

\newcommand{\Ima}{\mathrm{Im}}

\newcommand{\ka}{\kappa}
\newcommand{\sig}{\sigma}

\newcommand{\omg}{\omega}
\newcommand{\omgn}{\omega_{\nu}}
\newcommand{\omgd}{\omega_{d}}
\newcommand{\Omg}{\Omega}
\newcommand{\eps}{\epsilon}
\newcommand{\ups}{\upsilon}
\newcommand{\Ups}{\Upsilon}
\newcommand{\Lam}{\Lambda}
\newcommand{\lam}{\lambda}
\newcommand{\gam}{\gamma}

\newcommand{\llo}{l_{1}}
\newcommand{\llt}{l_{2}}
\newcommand{\llr}{l_{3}}

\newcommand{\hbz}{\hat{\Z}}
\newcommand{\hbzt}{\hbz^{2}}
\newcommand{\hbzgez}{\hbz_{\geq 0}}

\newcommand{\dd}{\diamond}
\newcommand{\bop}{\boxplus}
\newcommand{\bom}{\boxminus}

\newcommand{\mor}{\mbox{ or }}
\newcommand{\mand}{\mbox{ and }}
\newcommand{\mif}{\mbox{ if }}

\newcommand{\melse}{\mbox{ else}}
\newcommand{\mfor}{\mbox{ for }}
\newcommand{\mforsome}{\mbox{ for some }}
\newcommand{\mforeach}{\mbox{ for each }}
\newcommand{\mwhen}{\mbox{ when }}

\newcommand{\lcm}{\operatorname{lcm}}

\newcommand{\ord}{\mathrm{ord}}

\newcommand{\abs}[1]{\lvert#1\rvert}
\newcommand{\Li}{\operatorname{Li}}

\newcommand{\Log}{\operatorname{\mathrm{Log}}}
\newcommand{\Tr}[1]{{\mathrm{Tr}}\left(#1\right)}

\newcommand{\fX}{{\mathfrak{X}}}
\newcommand{\fXo}{\fX_{1}}
\newcommand{\fXt}{\fX_{2}}
\newcommand{\fXr}{\fX_{3}}

\newcommand{\fY}{{\mathfrak{Y}}}

\newcommand{\sce}{\succeq}
\newcommand{\sceo}{\sce_{1}}
\newcommand{\scet}{\sce_{2}}
\newcommand{\scer}{\sce_{3}}

\newcommand{\scc}{\succ}
\newcommand{\scco}{\scc_{1}}
\newcommand{\scct}{\scc_{2}}
\newcommand{\sccr}{\scc_{3}}

\newcommand{\geU}{\geq^{U}}
\newcommand{\geUfX}{\geU_{\fX}}

\newcommand{\geUx}{\geU_{x}}

\newcommand{\gUx}{\gU_{x}}

\newcommand{\geUxpm}{\geU_{\xpm}}

\newcommand{\gUxpm}{\gU_{\xpm}}

\newcommand{\gU}{>^{U}}
\newcommand{\gUfX}{\gU_{\fX}}

\newcommand{\xpm}{x^{\pm 1}}
\newcommand{\fXpm}{\fX^{\pm 1}}
\newcommand{\fXopm}{\fXo^{\pm 1}}
\newcommand{\fXtpm}{\fXt^{\pm 1}}

\newcommand{\gefXpm}{\ge_{\fXpm}}
\newcommand{\gefXopm}{\ge_{\fXopm}}

\newcommand{\gfXpm}{>_{\fXpm}}
\newcommand{\gfXopm}{>_{\fXopm}}

\newcommand{\geUfXpm}{\geU_{\fXpm}}

\newcommand{\gUfXpm}{\gU_{\fXpm}}

\newcommand{\qpm}{q^{\pm 1}}
\newcommand{\qq}{q^{\fraa}}

\newcommand{\AfX}{\mathrm{O}_{\fX}}
\newcommand{\AfXo}{\mathrm{O}_{\fXo}}
\newcommand{\AfXt}{\mathrm{O}_{\fXt}}
\newcommand{\AfXr}{\mathrm{O}_{\fXr}}

\newcommand{\geAfX}{\geq_{\AfX}}
\newcommand{\geAfXo}{\geq_{\AfXo}}
\newcommand{\geAfXt}{\geq_{\AfXt}}
\newcommand{\geAfXr}{\geq_{\AfXr}}

\newcommand{\gAfX}{>_{\AfX}}
\newcommand{\gAfXo}{>_{\AfXo}}
\newcommand{\gAfXt}{>_{\AfXt}}
\newcommand{\gAfXr}{>_{\AfXr}}

\newcommand{\StgeAfX}{S(\geAfX,\Q(\fX))}

\newcommand{\StgAfX}{S(\gAfX,\Q(\fX))}

\newcommand{\StgeR}{S(\geq,R)}
\newcommand{\StgR}{S(>,R)}

\newcommand{\StceR}{S(\sce,R)}
\newcommand{\StccR}{S(\scc,R)}
 
\newcommand{\StcefXo}{S(\sceo,\Q(\fXo))}
\newcommand{\StccfXo}{S(\scco,\Q(\fXo))}

\newcommand{\StcefXt}{S(\scet,\Q(\fXt))}
\newcommand{\StccfXt}{S(\scct,\Q(\fXt))}

\newcommand{\StcefXr}{S(\scer,\Q(\fXr))}
\newcommand{\StccfXr}{S(\sccr,\Q(\fXr))}

\newcommand{\StgeUfX}{S(\geUfX,\Q(\fX))}
\newcommand{\StgUfX}{S(\gUfX,\Q(\fX))}

\newcommand{\StgeUx}{S(\geUx,\Q(\fX))}
\newcommand{\StgUx}{S(\gUx,\Q(\fX))}

\newcommand{\StgeUfXpm}{S(\geUfXpm,\Q(\fX))}
\newcommand{\StgUfXpm}{S(\gUfXpm,\Q(\fX))}

\newcommand{\StgeUxpm}{S(\geUxpm,\Q(\fX))}
\newcommand{\StgUxpm}{S(\gUxpm,\Q(\fX))}

\newcommand{\gefX}{\geq_{\fX}}
\newcommand{\gefXo}{\geq_{\fXo}}

\newcommand{\gfX}{>_{\fX}}
\newcommand{\gfXo}{>_{\fXo}}

\newcommand{\gex}{\geq_{x}}

\newcommand{\gx}{>_{x}}

\newcommand{\gzd}{>_{z,d}}

\newcommand{\qr}{q^{\rho}}

\newcommand{\cB}{\mathcal{B}}

\newcommand{\cF}{\mathcal{F}}
\newcommand{\cG}{\mathcal{G}}
\newcommand{\cH}{\mathcal{H}}

\newcommand{\cL}{\mathcal{L}}
\newcommand{\cM}{\mathcal{M}}

\newcommand{\cQ}{\mathcal{Q}}
\newcommand{\cR}{\mathcal{R}}

\newcommand{\cZ}{\mathcal{Z}}

\newcommand{\bC}{\mathbb{C}}

\newcommand{\bP}{\mathbb{P}}
\newcommand{\bQ}{\mathbb{Q}}
\newcommand{\bR}{\mathbb{R}}

\newcommand{\bZ}{\mathbb{Z}}

\newcommand{\R}{\bR}
\newcommand{\Q}{\bQ}
\newcommand{\Z}{\bZ}
\newcommand{\N}{\Zgez}

\newcommand{\Ql}{\Q^{l}}

\newcommand{\Zl}{\Z^{l}}
\newcommand{\Zlo}{\Z^{\llo}}
\newcommand{\Zlt}{\Z^{\llt}}
\newcommand{\Zlr}{\Z^{\llr}}
\newcommand{\Ztl}{\Z^{2l}}

\newcommand{\Zd}{\Z^{d}}
\newcommand{\Qd}{\Q^{d}}

\newcommand{\Zo}{\Z^{1}}
\newcommand{\Zt}{\Z^{2}}
\newcommand{\Zr}{\Z^{3}}

\newcommand{\Zogeo}{\Zo_{\geq 1}}
\newcommand{\Ztgeo}{\Zt_{\geq1}}
\newcommand{\Zrgeo}{\Zr_{\geq1}}

\newcommand{\Zgez}{\Z_{\geq 0}}
\newcommand{\Zgeo}{\Z_{\geq 1}}
\newcommand{\Zget}{\Z_{\geq 2}}

\newcommand{\Zlgez}{\Zl_{\geq 0}}
\newcommand{\Zlgeo}{\Zl_{\geq 1}}

\newcommand{\Zdgez}{\Zd_{\geq 0}}
\newcommand{\Zdgeo}{\Zd_{\geq 1}}

\newcommand{\Zleo}{\Z_{\leq -1}}
\newcommand{\Zlez}{\Z_{\leq 0}}
\newcommand{\Qgo}{\Q_{>0}}
\newcommand{\llq}{>_{q}}
\newcommand{\ggq}{\geq_{q}}
\newcommand{\llqpm}{>_{q^{\pm 1}}  }
\newcommand{\ggqpm}{\geq_{q^{\pm 1}}}
 
\newcommand{\ve}{\vee}
\newcommand{\wcn}{\sqcup}
\newcommand\mdoubleplus{
  \ensuremath{\mathbin{+\mkern-10mu+}}}
\newcommand{\ccn}{\mdoubleplus}
\newcommand{\ehd}{\square}

\newcommand{\ld}{\lhd}
\newcommand{\nld}{\not\ld}
\newcommand{\ldZ}{\ld_{\Z}}
\newcommand{\nldZ}{\nld_{\Z}}
\newcommand{\ldZl}{\ld_{\Zl}}
\newcommand{\nldZl}{\nld_{\Zl}}
\newcommand{\ldZtl}{\ld_{\Ztl}}
\newcommand{\nldZtl}{\nld_{\Ztl}}
\newcommand{\ldZd}{\ld_{\Zd}}

\newcommand{\ldZo}{\ld_{\Zo}}

\newcommand{\ldZr}{\ld_{\Zr}}

\newcommand{\ldZlo}{\ld_{\Zlo}}

\newcommand{\ldZlt}{\ld_{\Zlt}}

\newcommand{\ldZlr}{\ld_{\Zlr}}

\newcommand{\lamo}{\lam_{1}}
\newcommand{\lamt}{\lam_{2}}
\newcommand{\lamr}{\lam_{3}}

\newcommand{\alp}{\alpha}
\newcommand{\alpo}{\alp_{1}}
\newcommand{\alpt}{\alp_{2}}

\newcommand{\bta}{\beta}
\newcommand{\btao}{\bta_{1}}
\newcommand{\btat}{\bta_{2}}

\newcommand{\gamo}{\gam_{1}}
\newcommand{\gamt}{\gam_{2}}
\newcommand{\gamr}{\gam_{3}}

\newcommand{\muo}{\mu_{1}}
\newcommand{\mut}{\mu_{2}}
\newcommand{\mur}{\mu_{3}}

\newcommand{\mn}{m_{1}}
\newcommand{\mt}{m_{2}}
\newcommand{\mr}{m_{3}}
  
\newcommand{\nn}{n_{1}}
\newcommand{\nt}{n_{2}}
\newcommand{\nr}{n_{3}}
 
\newcommand{\kn}{k_{1}}
\newcommand{\kt}{k_{2}}
\newcommand{\kr}{k_{3}}

\newcommand{\thetao}{\theta_{1}}
\newcommand{\thetat}{\theta_{2}}

\newcommand{\cFz}{\cF_{0}}
\newcommand{\cFo}{\cF_{1}}
\newcommand{\cFt}{\cF_{2}}
\newcommand{\cFr}{\cF_{3}}

\newcommand{\fs}{f_{s}}
\newcommand{\gs}{g_{s}}
\newcommand{\hs}{h_{s}}

\newcommand{\rhoo}{\rho_{1}}
\newcommand{\rhot}{\rho_{2}}
\newcommand{\rhor}{\rho_{3}}

\newcommand{\phio}{\phi_{1}}
\newcommand{\phit}{\phi_{2}}
\newcommand{\phir}{\phi_{3}}

\newcommand{\xoo}{x_{1}}

\newcommand{\til}{\tilde}

\newcommand{\tcL}{\til{\cL}}
\newcommand{\tcQ}{\til{\cQ}}

\newcommand{\FD}{\mathrm{F_{D}}}
\newcommand{\RD}{\mathrm{R_{D}}}
\newcommand{\AD}{\mathrm{A_{D}}}

\newcommand{\FZ}{\mathrm{F_{Z}}}
\newcommand{\RZ}{\mathrm{R_{Z}}}
\newcommand{\AZ}{\mathrm{A_{Z}}}

\newcommand{\FC}{\mathrm{F_{C}}}
\newcommand{\RC}{\mathrm{R_{C}}}
\newcommand{\AC}{\mathrm{A_{C}}}
  
\newcommand{\Bo}{\mathrm{Bm}}
\newcommand{\To}{\mathrm{Tp}}

\newcommand{\ci}{C}
\newcommand{\ri}{R}

\newcommand{\Sup}{\Supset}

\newcommand{\oi}{\operatorname{
    \mathrm{I}}}
\newcommand{\tsi}{\operatorname{
    \mathrm{T_{<}}}}
\newcommand{\tei}{\operatorname{
    \mathrm{T_{\leq}}}}
\newcommand{\ofi}{\operatorname{
    \mathrm{F\tsi}}}

\newcommand{\Qti}{Q^{\times}}
\newcommand{\Qr}{\Q^{3}}

\newcommand{\zetao}{\zeta_{1}}
\newcommand{\zetat}{\zeta_{2}}
\newcommand{\zetar}{\zeta_{3}}

\newcommand{\geqx}{\geq^{x}}

\newcommand{\out}{
  \operatorname{\mathrm{O}}}
\newcommand{\cen}{
  \operatorname{\mathrm{C}}}
\newcommand{\tout}{
  \operatorname{\til{\out}}}

\newcommand{\otr}{
  \operatorname{\mathrm{ot}}}
\newcommand{\ctr}{
  \operatorname{\mathrm{ct}}}

\newcommand{\dt}{\tilde{d}}

 \newcommand{\dif}{\operatorname{\mathrm{L}}}
\newcommand{\tdif}{\operatorname{\mathrm{R}}}

\newcommand{\etg}{\gamma}
\newcommand{\ep}{H}

\newcommand{\gqd}{>_{q,d}}
\DeclarePairedDelimiter{\tp}{\llbracket}{\rrbracket}

\newcommand{\tgq}[1]{>_{#1}}
\newcommand{\tgeq}[1]{\geq_{#1}}

\newcommand{\fU}{{\mathfrak{U}}}
\newcommand{\fUq}[1]{\fU_{#1}}

\begin{document}
\author{So Okada} \title{Merged-log-concavity of rational
  functions, almost strictly unimodal sequences, and phase
  transitions of ideal boson-fermion gases}
\maketitle  
\begin{abstract}  
  We obtain some new results on the unimodal sequences of
  the real values of rational functions by polynomials with
  positive integer coefficients.  Thus, we introduce the
  notion of {\it merged-log-concavity} of rational functions.
  Roughly speaking, the notion
  extends Stanley's $q$-log-concavity of polynomials.

  We construct explicit merged-log-concave rational
  functions by $q$-binomial coefficients, Hadamard products,
  and convolutions, extending the Cauchy-Binet formula.
  Then, we obtain the unimodal sequences of rational
  functions by Young diagrams.  Moreover, we consider the
  variation of unimodal sequences by critical points that
  separate strictly increasing, strictly decreasing, and
  hill-shape sequences among almost strictly unimodal
  sequences.  Also, the critical points are zeros of
  polynomials in a suitable setting.
  
  The study above extends the $t$-power series of
  $(\pm t;q)_{\infty}^{\mp 1}$ to some extent by polynomials with
  positive integer coefficients and the variation of
  unimodal sequences.  We then obtain the golden ratio of
  quantum dilogarithms ($q$-exponentials) as a critical
  point.  Additionally, we consider eta products,
  generalized Narayana numbers, and weighted
  $q$-multinomial
  coefficients, which we introduce.
  
  In statistical mechanics, we discuss the grand canonical
  partition functions of some ideal boson-fermion gases with
  or without Casimir energies (Ramanujan summation).  The
  merged-log-concavity gives phase transitions on Helmholtz
  free energies by critical points of the metallic ratios
  including the golden ratio.  In particular, the phase
  transitions implies non-zero particle vacua from zero
  particle vacua as the temperature rises.
 \end{abstract}
 
 \section{Introduction}\label{sec:intro}
 The notion of unimodal sequences includes increasing,
 decreasing, or hill-shape sequences of real numbers.  As
 such, the notion is quite essential whenever one computes.
 Also, there exist some notions of log-concavities for
 polynomials~\cite{But,New,Sta}.  Thus, we introduce the
 notion of {\it merged-log-concavity} for rational
 functions.  Then, to some extent, we study
 merged-log-concave rational functions by polynomials with
 positive integer coefficients and the variation of unimodal
 sequences on continuous parameters.
  
 Let us first recall the following unimodality and
 log-concavity of real numbers and polynomials.
  
 \subsection{Unimodality and log-concavity
   of real numbers}
 Let $\hbz=\Z\cup \{\infty\}$. Also, for $i,j\in \hbz$, let
 $\oi(i,j)=\{k\in \Z \mid i\leq k \leq j\}$ and $\oi(j)=\oi(1,j)$.
 \begin{definition}\label{def:unimodal}
   Suppose a sequence
   $r=\{r_{i}\in \R\}_{i \in \oi(s_{1},s_{2})}$.
      \begin{enumerate}
      \item
        The sequence $r$ is unimodal,
        if
        $r_{s_{1}}\leq \dots \leq r_{\del}\geq r_{\del+1}\geq \dots$ for
        some $\del\in \hbz$ such that $s_{1}\leq \del\leq s_{2}$.  In
        particular, if $\del\in \oi(s_{1}+1,s_{2}-1)$, then
        $r$ is said to be hill-shape.
      \item
        The sequence $r$ is log-concave,
        if $r_{i}^{2}-r_{i+1}r_{i-1}\geq 0$ for each
     $i\in \oi(s_{1}+1,s_{2}-1)$.
   \item If $r_{\mu_{1}}<r_{\mu_{2}}$  for
     some $\mu_{1},\mu_{2}\in \oi(s_{1},s_{2})$,
     then $r$ is said
     to have a slope.
   \item Let  $s_{2}-s_{1}+1\in \hbz$ be the length of $r$.
   \end{enumerate}
\end{definition}

For example, let
  $r=\{r_{i}\in \R\}_{i \in \oi(0,\infty)}$.
 First, suppose that $r$ is unimodal. 
Then, $\del=0$, $\del\in\oi(\infty)$, or
$\del=\infty$ gives a
decreasing sequence $ r_{0}\geq r_{1}\geq r_{2}\geq \dots$,
hill-shape sequence
$ r_{0}\leq \dots \leq r_{\del}\geq r_{\del+1}\geq \dots$, or
increasing sequence
$ r_{0}\leq r_{1}\leq r_{2}\leq \dots$, respectively.
Second, suppose that $r$ is log-concave such that each $r_{i}>0$.
  Then,
$ \frac{r_{1}}{r_{0}}\geq \frac{r_{2}}{r_{1}}\geq\dots$.  Hence,
$r$ is unimodal, since
$\dots \leq r_{\del}\geq r_{\del+1}\geq \dots$ for the smallest
$\del\geq 0$ such that $\frac{r_{\del+1}}{r_{\del}}\leq 1$.
\subsection{\texorpdfstring{$q$}{q}-log-concavity of
  polynomials}\label{sec:intro-q-log}
Unless stated otherwise, let $q$ be an indeterminate.
Then, we take the following binary relations
  on $q$-rational functions $\Q(q)$ 
by $q$-polynomials $\Q[q]$ and Laurent $q$-polynomials
$\Q[q^{\pm 1}]$.

\begin{definition}\label{def:q-pic-poly-q-Laurent-pic-poly}
  Suppose $f,g\in \Q(q)$.
  \begin{enumerate}
  \item Let $f\ggq g$, if $f,g\in \Q[q]$ and $f-g \in \N[q]$.
  \item Let  $f\llq g$, if
    $f\ggq g$ and 
    $f-g\neq 0$.
  \item Let $f\ggqpm g$, if $f,g\in \Q[q^{\pm 1}]$ and
    $f-g\in \N[q^{\pm 1}]$.
  \item Let $f\llqpm g$, if $f\ggqpm g$ and
    $f-g\neq 0$.
    \end{enumerate}
\end{definition}

Sagan introduced the notation $\ggq$ on $\Z[q]$~\cite[Section
1]{Sag}.  Moreover, Stanley gave the following notion
of $q$-log-concavity~\cite[Definition 1.1]{But} \cite[Section 1]{Sag}.  We studied the notion
intensively~\cite{Bre}, but did not extend
the notion to $\Q(q)$.
\begin{definition}\label{def:log-conc}
  Suppose a sequence $f=\{f_{i}(q)\in \N[q]\}_{i\in \Z}$.
  Then, $f$ is called {\it $q$-log-concave}, if
  $f_{i}(q)^{2}-f_{i-1}(q)f_{i+1}(q)\ggq 0$
  for each $i\in \Z$.
\end{definition}
Suppose a $q$-log-concave
$f=\{f_{i}(q)\llq 0\}_{i\in \Z}$ such that
$f_{i}(q)\llq 0$ for each $i\in \Zgez$.  Then,
 for each $r\in \R$,
$\{f_{i}(r)\in \R_{>0}\}_{i\in \Zgez}$ is
log-concave and unimodal.

\subsection{Merged-log-concavity of rational functions}
\label{sec:intro-merged-rational-functions}

We would like a log-concavity of rational
functions that gives polynomials of positive integer
coefficients.   Thus, we introduce the
 merged-log-concavity as a log-concavity of rational
 functions in some generality.
We use the following $q$-analogs.

\begin{definition}\label{def:pochh}
  Assume an indeterminate $a$.
  \begin{enumerate}
  \item For each $n\in \hbzgez$,
    the $q$-Pochhammer symbol
    $(a;q)_{n}$ is $\prod_{i\in\oi(n)}(1-aq^{i-1})$ if $n\geq 1$ and
    $1$ otherwise.  Let $(n)_{q}=(q;q)_{n}$ for our
    convention.
       \item Each $n\in \Zgeo$ defines
          the {\it $q$-number}
         $[n]_{q}=\sum_{ i\in \oi(n)}q^{i-1}$
         and {\it
           $q$-factorial} 
         $[n]!_{q}=\prod_{i\in\oi(n)}[i]_{q}$.
   As special cases,
   $[0]_{q}=0$ and
   $[0]!_{q}=1$.
 \end{enumerate}
 \end{definition}

 Instead of $( \ )$, we use
 $\llbracket \ \rrbracket$ for tuples as follows.
 For instance, this avoids
 confusions between  parentheses of operator
 precedence and tuples, as we consider tuples of tuples.

 \begin{definition}\label{def:cartesian}
   Suppose sets $U_{1},\dots, U_{d}$ for some $d\in \Zgeo$.
   \begin{enumerate}
   \item
     Let
   $\prod_{i\in\oi(d)}U_{i}=U_{1}\times\dots\times U_{d} $ be the
   Cartesian product.
   \item Let
   $a=\tp{a_{1},a_{2},\dots,
     a_{d}}=\tp{a_{i}\in U_{i}}_{i\in\oi(d)}$
    denote a tuple
   $a\in \prod_{i\in\oi(d)}U_{i}$.
 \item
   For a tuple $a\in \prod_{i\in\oi(d)}U_{i}$,
   let $a^{\ve}$ be the
   flip
   $\tp{a_{d-i+1}}_{i\in \oi(d)}\in \prod_{i\in\oi(d)}U_{d-i+1}$.
 \item For each $\alp\in U_{1}$, let $\iota^{d}(\alp)$ be
     the tuple
   $\tp{\alp,\alp,\dots,\alp}\in
   U_{1}^{d}=\prod_{i\in\oi(d)}U_{1}$.
 \end{enumerate}

 In particular, assume $a\in \Zd$. Then, we write $a \geq 0$ if
 each $a_{i}\in \Zgez$, $a=0$ if each $a_{i}=0$, and
 $a\leq 0$ if each $a_{i}\in \Zlez$.  Also, let
 $\abs{a}=\tp{\abs{a_{i}}}_{i\in \oi(d)}\in \Zdgez$.
 Moreover, if $b,c\in \Z^{2d}$, then let
 $\dif(b,c),\tdif(b,c)\in \Zd$ such that
 $\dif(b,c)_{i}=c_{i}-b_{i}$ and
 $\tdif(b,c)_{i}=c_{d+i}-b_{d+i}$ for $i \in\oi(d)$.
   \end{definition}
   For instance, $a\in \Z^{1}$ is a tuple
 $\tp{a_{1}}\in \Z^{1}\neq \Z$.
 We now introduce the
 merged-log-concavity, which is a
 simpler version of that in Definition~\ref{def:merged}.

 \begin{definition} \label{def:merged-simpler} Let
   $u^{-1},\rho,l\in \Zgeo$.  Also, let $w\in \Zgez^{l}$ and
   $a,b\in\Ztl$.
    \begin{enumerate}
    \item In  $\Q(q^{u})$, we define the ring shift
      factor
      \begin{align}
        \Ups(l,w,\rho,a,b,q)
        =
        \begin{dcases}
          \prod_{i\in\oi(l)}
          \frac{(b_{i})^{w_{i}}_{\qr}
          (b_{l+i})^{w_{l-i+1}}_{\qr}}
          {(a_{i})^{w_{i}}_{\qr} 
          (a_{l+i})^{w_{l-i+1}}_{\qr} 
          } \mif a,b\geq 0, \\
          0 \melse.
        \end{dcases}
      \end{align}
     We call $\rho$
     base shift parameter
     of $\Ups(l,w,\rho,a,b,q)$.
    \item For a family
      $\cF=\{\cF_{m}(q^{u})\in \Q(q^{u})\}_{m\in \Zl}$,
       we define the merged determinant
      \begin{align}
          \Delta(\cF)(l,w,\rho,a,b,q)
         =
          \Ups(l,w,\rho,a,b,q)
          \cdot \det\begin{bmatrix}
                  \cF_{\dif(a,b)}&         \cF_{\tdif(a^{\ve},b)}\\
                  \cF_{\dif(a^{\ve},b)} &         \cF_{\tdif(a,b)}
                                          \label{eq:determinant}
                \end{bmatrix}\in \Q(q^{u}).
      \end{align}
      \item Suppose a family
        $f=\{f_{m}\in \Q(q^{u})\}_{m\in \Zl}$ such that
        \begin{align}
          f_{m}&>_{(q^{u})^{\pm 1}}0 \mif m\geq 0,\\
          f_{m}&=0  \melse.
          \end{align}
     Then, we define the parcel
     $\cF=\Lam(l,w,f,q,u)=\{\cF_{m}(q^{u})\in
     \Q(q^{u})\}_{m\in
       \Zl}$ as a family of rational functions such that
     \begin{align}
       \cF_{m}&=
                \begin{dcases}
                  \frac{f_{m}}{\prod_{i\in\oi(l)}(m_{i})^{w_{i}}_{q}}
                  \mif m\geq 0,\\
                  0 \melse.
                \end{dcases}
     \end{align}
     We refer to $l$, $w$, $f$, and $q$ as the
     width, weight, numerator, and base of the parcel
     $\cF$. 
  \item
    \label{fitting-cond-simple}
    For  $a,b\in \Z^{2l}$, we call $\tp{l,a,b}$ fitting,
    if it satisfies the following inequalities:
    \begin{align}
      a&\leq b;\\
      b_{1}\leq \dots \leq b_{l}&<b_{l+1}\leq \dots \leq b_{2l};\\
      0\leq a_{1}\leq \dots \leq a_{l}&<a_{l+1}\leq \dots \leq a_{2l}.
    \end{align}
  \item 
    We call
    $\cF=\Lam(l,w,f,q,u)$ 
    $\rho$-merged-log-concave,
    if each fitting $\tp{l,a,b}$ satisfies
    \begin{align}
      \Delta(\cF)(l,w,\rho,a,b,q)>_{(q^{u})^{\pm 1}}0.
    \end{align}
    \end{enumerate}
\end{definition}

As $\Ups(l,w,\rho,a,b,q)=1$ of $w=0$, merged determinants
generalize $2\times2$ determinants by ring shift factors.  We
put the term ``ring shift factors'', since we study
$\Delta(\cF)(l,w,\rho,a,b,q)\in \Q[q^{\pm u}]$, when
$\cF_{m}\in \Q(q^{u})$ are not necessarily in
$\Q[q^{\pm u}]$.  We use the term ``merged determinants''
for $\Delta(\cF)(l,w,\rho,a,b,q)$, as $\Ups(l,w,\rho,a,b,q)$ and
$\cF_{\dif(a,b)} \cF_{\tdif(a,b)} - \cF_{\dif(a^{\ve},b)}
\cF_{\tdif(a^{\ve},b)}$ merge into $\Delta(\cF)(l,w,\rho,a,b,q)$.

We introduce the notion below to
obtain explicit merged-log-concave parcels.

\begin{definition}\label{def:monomial-index} Let
  $l\in \Zgeo$, $w\in \Zgeo^{l}$, and
  $\gam\in \prod_{i\in\oi(l)}\Qr$. We call the tuple
  $\ups=\tp{l,w,\gam}$ monomial index, if
  \begin{align}
    2\gam_{i,1}&\in \Z  \mforeach i\in\oi(l),
                 \label{cond:monomial-index-integer}\\
    0\leq 2\sum_{j\in\oi(i)}\gam_{j,1}&\leq \sum_{j\in\oi(i)} w_{j}
                                  \mforeach i\in\oi(l).
                                  \label{cond:monomial-index-sum}
  \end{align}
  We call $l$, $w$, and $\gam$ width, weight, and core of
  $\ups$.  We refer to \eqref{cond:monomial-index-integer}
  and~\eqref{cond:monomial-index-sum} as the integer
  monomial condition and the sum monomial condition of
  $\ups$.
 \end{definition}

 We also put the following 
 quadratic polynomials
  and rational numbers.
\begin{definition}\label{def:intro-quad-poly}
  Suppose $l\in \Zgeo$, $\gam\in \prod_{i\in\oi(l)}\Qr$, and
  $\alp\in \Zlgeo$ with
  $\del=\gcd(\alp_{1},\dots,\alp_{l})$.
  \begin{enumerate}
  \item Let
    $t_{\gam_{i}}(z)=
    \gam_{i,1}z^{2}+\gam_{i,2}z+\gam_{i,3}\in \Q[z]$ for each
    $i\in\oi(l)$.
  \item Let $t_{\alp,\gam}:\Zl\to \Q$ such that
    $t_{\alp,\gam}(m)=\sum_{i\in\oi(l)}
    \alp_{i}t_{\gam_{i}}(m_{i})$ for each $m\in \Zl$.
  \item Let $u_{\alp,\gam}=\frac{\delta}{\lam}\in \Q$ for
    the lowest
    $\lam\in \Zgeo$ 
    such that
    $\frac{t_{\alp,\gam}(m)}{u_{\alp,\gam}}\in \Z$ for each
    $m\in \Zlgez$.
  \end{enumerate}
  For simplicity,
  we write $t_{\gam}=t_{\iota^{l}(1),\gam}$ and
  $u_{\gam}=u_{\iota^{l}(1),\gam}$.
\end{definition}
  
Assume a monomial index $\tp{l,w,\gam}$.  For each
$m\in \Zl$, let
\begin{align}
  f_{\gam,m}(q)=
  \begin{dcases}
    q^{t_{\gam}(m)} \mif m\geq 0,\\
    0  \melse.
  \end{dcases}
\end{align}
We then define {\it the monomial parcel}
$\cF_{w,\gam,q}=\Lam(l,w, f_{\gam}(q),q,u_{\gam})$ such that
 each $m\in \Zlgez$ satisfies
\begin{align}
  \cF_{w,\gam,q,m}=\frac{q^{t_{\gam}(m)}}{
  \prod_{i\in \oi(l)}(m_{i})^{w_{i}}_{q}}\in \Q(q^{u_{\gam}}).
\end{align}
Moreover, we have the merged-log-concavity
of monomial parcels as follows.

\begin{theorem}[Theorem~\ref{thm:monomial-poly}]
   \label{thm:monomial-poly-simpler}
   Let $\rho\in \Zgeo$.
   Then, each monomial parcel $\cF_{w,\gam,q}$ is
   $\rho$-merged-log-concave. In particular,  for any
   fitting
   $\tp{l,a,b}$, we have
   \begin{align}
     q^{-t_{\gam}(\dif(a,b))-t_{\gam}(\tdif(a,b))}
         \Delta(\cF_{w,\gam,q})(l,w,\rho,a,b,q)\llq 0.
   \end{align}
\end{theorem}

Unless stated otherwise, we assume the following
simplification until the end of Section~\ref{sec:intro}:
$l=1$; the width of each parcel is $l$; and,
$\cF_{m_{1}}=\cF_{m}$ for each $m=\tp{m_{1}}\in \Zl$ and
parcel $\cF$.  
We also put $\cF_{\gam}=\cF_{w,\gam,q}$ for short,
given a monomial index $\tp{l,w,\gam}$
of $w=\tp{1}$.
For
$\cF_{\gam}$,
the monomial conditions
imply $\gam_{1,1}=\fraa \mor 0$.  Hence, we write
$\tcL_{\tp{\gam_{1,2},\gam_{1,3}}}
=\cF_{\tp{\tp{0,\gam_{1,2},\gam_{1,3}}}}$ and
$\tcQ_{\tp{\gam_{1,2},\gam_{1,3}}}
=\cF_{\tp{\tp{\fraa,\gam_{1,2},\gam_{1,3}}}}$ so that
$\lam\in \Zgez$ satisfies
\begin{align}
\tcL_{\tp{\gam_{1,2},\gam_{1,3}},\lam}& =
\frac{q^{\gam_{1,2} \lam+\gam_{1,3}}}{(\lam)_{q}},\\
\tcQ_{\tp{\gam_{1,2},\gam_{1,3}},\lam}
&=\frac{q^{\frac{\lam^2}{2}
    +\gam_{1,2}\lam+\gam_{1,3}}}{(\lam)_{q}}.
\end{align}
We call
$\tcL_{\tp{\gam_{1,2},\gam_{1,3}}}$ and
$\tcQ_{\tp{\gam_{1,2},\gam_{1,3}}}$ linear and quadratic
monomial parcels, referring to their $q$-powers of numerators.
In particular, we call $\cL=\tcL_{\tp{\fraa,0}}$ {\it
  linear-half monomial parcel} and $\cQ=\tcQ_{\tp{0,0}}$
{\it quadratic-half monomial parcel}.

Clearly, quadratic numerators
 $q^{\frac{\lam^{2}}{2}}$ of
 $\cQ$  give no polynomials with positive integer coefficients by
 \begin{align}
      q^{\frac{\lam^2}{2}}        q^{\frac{\lam^2}{2}}
   -    q^{\frac{(\lam-1)^2}{2}}
   q^{\frac{(\lam+1)^2}{2}}
      =q^{\lam^{2}}-q^{\lam^{2}+1}.
   \label{eq:poly-not-pos-int-coeff}
 \end{align}
 But,
 a key idea
 in Theorem~\ref{thm:monomial-poly-simpler}
 is that
 $\cQ$ gives those by
 merged determinants.
 For example,
    \begin{align}
      \Delta(\cQ)(1,\tp{1},1,\tp{0,1},\tp{1,2},q)
      =&  \frac{(1)_{q} (2)_{q}}{(0)_{q} (1)_{q}}
         (\cQ_{1}^{2}-\cQ_{0} \cQ_{2})\\
      =& 
         (1-q)(1-q^{2})
         \left(\left(\frac{q^{\frac{1^2}{2}}}{1-q}\right)^{2}-
         \frac{q^{\frac{0^{2}}{2}}}{1}\cdot
         \frac{q^{\frac{2^2}{2}}}{(1-q)(1-q^{2})}\right)\\
      =&q.
    \end{align}

    For further discussion, we put the following notation
    for the base shift parameters.
    \begin{definition}\label{def:intro-base-shift}
      For each $\lam\in \Z$ and $\rho\in \Zgeo$, we 
      put the {\it base shift  function}
  \begin{align}
    b_{\lam,\rho}(q)     
    &=
      \begin{dcases}
        \frac{(\lam)_{\qr}}{(\lam)_{q}}
        \mif    \lam\in \Zgez,\\
        0 \melse.
      \end{dcases}
  \end{align}  
  \end{definition}
  Also, we recall the following by tuples.
    \begin{definition}\label{def:gauss-multi}
       Let $d\in \Zgeo$.
       Then,   $i\in \Z$ and $j\in \Zd$ give
       the $q$-multinomial coefficient
       \begin{align}
         {i \brack j}_{q}
         =
         \begin{dcases}
           \frac{(i)_{q}}{\prod_{\lam\in\oi(d)}(j_{\lam})_{q}}
           \mif j\geq 0  \mand \sum_{\lam\in\oi(d)} j_{\lam}=i,\\
           0 \melse.
         \end{dcases}
       \end{align}
       In particular,  $i,j\in \Z$ give the $q$-binomial
       coefficient 
       ${i \brack j}_{q} ={ i \brack \tp{j,i-j}}_{q}={ i
         \brack \tp{i-j,j}}_{q}$.
  \end{definition}
  We then obtain Theorem~\ref{thm:monomial-poly-simpler},
  studying the $q$-multinomial coefficients and base shift
  functions, because
  \begin{align}
    \Delta(\cF)(l,w,\rho,a,b,q)
    &  =\Ups(l,w,\rho,a,b,q)    \cF_{\dif(a,b)} \cF_{\tdif(a,b)}
      - \Ups(l,w,\rho,a,b,q) \cF_{\dif(a^{\ve},b)} \cF_{\tdif(a^{\ve},b)}
    \\ &={ b_{1}\brack a_{1} }^{w_{1}}_{\qr}
         { b_{2} \brack a_{2} }^{w_{1}}_{\qr}
         b_{\dif(a,b)_{1},\rho}(q)^{w_{1}}
         b_{\tdif(a,b)_{1},\rho}(q)^{w_{1}}
         f_{\dif(a,b)}f_{\tdif(a,b)}
    \\&-{ b_{2} \brack a_{1} }^{w_{1}}_{\qr}
    { b_{1}\brack a_{2}}^{w_{1}}_{\qr}
    b_{\dif(a^{\ve},b)_{1},\rho}(q)^{w_{1}}
    b_{\tdif(a^{\ve},b)_{1},\rho}(q)^{w_{1}}
    f_{\dif(a^{\ve},b)}f_{\tdif(a^{\ve},b)}.
  \end{align}
  
  \subsection{Parcel convolutions and
    an extended Cauchy-Binet formula}
  For appropriate parcels $\cF$ and $\cG$, we discuss the
  {\it parcel convolution} $\cF*\cG$, which multiply their
  Toeplitz matrices.  The Cauchy-Binet formula writes the
  minors of a matrix product $AB$ by ones of $A$ and $B$.
  Also, merged determinants extend $2\times 2$ determinants by
  ring shift factors.  Thus, in
  Theorem~\ref{thm:rshift-ext-cb}, we extend general minors
  and the Cauchy-Binet formula to write merged determinants
  of $\cF*\cG$ by those of $\cF$ and $\cG$. This gives the
  merged-log-concavity of parcel convolutions.

  In particular, the multifold convolutions of the linear
  and quadratic monomial parcels $\tcL_{\tp{0,0}}$ and
  $\tcQ_{\tp{-\fraa,0}}$ are merged-log-concave.  Hence,
  their merged determinants give the following polynomials
  with positive integer coefficients.
  
  \begin{theorem} (special cases of
    Theorem~\ref{thm:monom-conv-mult})
    \label{thm:intro-poly-positivity}
      Suppose $d\in \Zgeo$ and $h\in \Zgez$.  The following
      are $q$-polynomials with positive integer coefficients:
        \begin{dmath*}
          [h+1]_{q}
          \left(\sum_{j_{1}\in \Zd}
          {h \brack j_{1}}_{q}\right)^{2}
          -
          [h]_{q}
          \left(
          \sum_{j_{1}\in \Zd}
          {h +1\brack j_{1}}_{q}
          \right)
          \left(
          \sum_{j_{2}\in \Zd}
          {h-1 \brack j_{2}}_{q}
          \right) \llq 0;
        \end{dmath*}
        \begin{dmath*}
          [h+1]_{q}
          \left(\sum_{j_{1}\in \Zd}
            \left(\prod_{i\in\oi(d)} q^{\frac{j_{1,i}(j_{1,i}-1)}{2}}\right)
            {h \brack j_{1}}_{q}\right)^{2}
          -
          [h]_{q}
          \left(
            \sum_{j_{1}\in \Zd}
            \left(
              \prod_{i\in\oi(d)} 
              q^{\frac{j_{1,i}(j_{1,i}-1)}{2}}
            \right)
            {h +1\brack j_{1}}_{q}
          \right)
          \left(
            \sum_{j_{2}\in \Zd}
            \left(
              \prod_{i\in\oi(d)} 
              q^{\frac{j_{2,i}(j_{2,i}-1)}{2}}
            \right)
            {h-1 \brack j_{2}}_{q}
          \right) \llq 0.
        \end{dmath*}
    \end{theorem}

    \subsection{Almost strictly unimodal sequences of
      rational functions and critical
      points}\label{sec:intro-unim-seq-crit-pts}
    
    We recall the following notions of
    {\it almost strictly unimodal
      sequences} \cite[Section 2.2]{Rea} and {\it strictly
      log-concave sequences}.  Then, we consider the
    variation of almost strictly unimodal sequences by the
    merged-log-concavity.

    \begin{definition}\label{def:semi-st-uni-st-log-conc}
      Suppose a sequence
      $r=\{r_{i}\in \R\}_{i\in\oi(s_{1},s_{2})}$.
      \begin{enumerate}
      \item The sequence $r$ is almost strictly
        unimodal, if  we have $\del \in \hbz$
        with $s_{1}\leq \del\leq s_{2}$ such that
        $\{r_{i}\}_{i\leq \delta}$ is strictly increasing,
        $r_{\delta}\geq r_{\delta+1}$, and
        $\{r_{i}\}_{i\geq\delta+1}$ is strictly decreasing.  For
        simplicity, we refer to such $\del$ as the mode of
        $r$.
      \item The sequence $r$ is strictly log-concave, if
        $r_{i}^{2}-r_{i-1}r_{i+1}>0$
        for each $i\in \oi(s_{1}+1,s_{2}-1)$.
      \end{enumerate}
    \end{definition}
    For instance, suppose a strictly log-concave
    $r=\{r_{i}\in \R_{>0}\}_{i\in\oi(s_{1},s_{2})}$.
    Then, $r$ is almost
    strictly unimodal by
    $\frac{r_{s_{1}+1}}{r_{s_{1}}}
    >\frac{r_{s_{1}+2}}{r_{s_{1}+1}}>\dots$.  We also     discuss
    $\{\log(r_{i})\in \R\}_{i\in\oi(s_{1},s_{2})}$, which
    is not necessarily strictly log-concave but almost
    strictly unimodal.

    Let $\cF=\Lam(l,w,f,q,u)$ be $\rho$-merged-log-concave.
    Consider the sequence
    $u(\cF,h)=\{\cF_{m}(h)\in \R_{>0}\}_{m\in \Zgez}$ on the
    continuous parameter $0<h=q^{u}<1$ of $u^{-1}\in \Zgeo$.
    Then, $u(\cF,h)$ is strictly log-concave, because
    \begin{align}
      \cF_{m}(h)^{2}-\cF_{m-1}(h)\cF_{m+1}(h)>0
      \label{ineq:det-positivity}
    \end{align}
    for each $m\in \Zgez$ by 
    $\Delta(\cF)(l,w,\rho,\tp{0,1},\tp{m,m+1},q)
      >_{(q^{u})^{\pm 1}}0$,        
      and
      \begin{align}
      \Ups(l,w,\rho,\tp{0,1},\tp{m,m+1},h^{u^{-1}})
      =
      \left. \frac{(m)^{w_{1}}_{\qr}
      (m+1)^{w_{1}}_{\qr}}{
        (0)^{w_{1}}_{\qr}(1)^{w_{1}}_{\qr}}
        \right|_{q=h^{u^{-1}}}>0.
      \label{ineq:intro-pos-rshift-factor}
    \end{align}
    We thus have the variation of almost strictly unimodal
    sequences $u(\cF,h)$ on $0<h<1$.
      More generally,
    Young diagrams give almost strictly unimodal sequences
    in Section~\ref{sec:semi-st-unim-seq-young}.  For
    instance, the above tuples $\tp{m,m+1}$ correspond to
    the trivial Young diagram
    $\begin{ytableau} 1 \end{ytableau}$, whose number
    coincides with that of
    $\tp{1}=\dif(\tp{m,m+1},\tp{m+1,(m+1)+1})$.
    
    Among almost strictly unimodal sequences,
     there are hill-shape and decreasing sequences.
     Thus, if $u(\cF,h)$ is one of
     these boundary sequences, then we call $h$ {\it
       critical point} of $\cF$.  In particular, since
     $u(\cF,h)$ is almost strictly unimodal, $u(\cF,h)$ is
     hill-shape and decreasing if and only if
    \begin{align}
      \cFz(h)=\cFo(h).\label{eq:hill-shape-and-decreasing}
    \end{align}
    Moreover, since
    $q^{\lam u}(1-q)(\cFz(h)-\cFo(h))\in \Q[q^{u}]$ for some
    $\lam \in \Zgeo$, these $h$ are zeros of
    $q^{u}$-polynomials over $0<q^{u}<1$.  As for the term
    ``critical points'', these $h$ give not only boundary
    sequences, but also the zero discrete derivative
    $0=\cFo(h)-\cFz(h)$.

    \subsection{Metallic ratios as explicit critical points}
    \label{sec:intro-explicit-crit-pts}
    For $n\in \Zgeo$,
    $\frac{-n+\sqrt{n^{2}+4}}{2}:1=
    1:\frac{n+\sqrt{n^{2}+4}}{2}$.
    Hence, we call $\frac{-n+\sqrt{n^{2}+4}}{2}<1$
    {\it metallic ratios} for our convention, instead of
    $\frac{n+\sqrt{n^{2}+4}}{2}>1$
    in~\cite[Section 1]{GilWor}.  For instance,
    $\frac{-1+\sqrt{5}}{2}=0.618\dots$ and
    $\frac{-2+\sqrt{8}}{2}=0.414\dots$ are the golden
    and    silver ratios.
    Then,
        the golden ratio is the critical point of 
        $\cL$, because
        over
        $0<q^{\fraa}=q^{u_{\tp{{0,\fraa,0}}}}<1$,
        it solves
    \begin{align}
      \cL_{0}(q^{\fraa})
      =1=\frac{q^{\fraa}}{1-q}=\cL_{1}(q^{\fraa}),
    \end{align}
    or $q+q^{\fraa}-1=0$.
    
    Over $0<h<1$, if $u(\cF,h)$ changes from a strictly
    decreasing sequence to a two-slope hill-shape or
    strictly increasing sequence, then we say that $\cF$ has
    a phase transition.  Thus, $\cL$ has the phase
    transition by the golden ratio in the figure below.
    \vspace{1em}
    \begin{figure}[H]{}
      \centering  
      \resizebox{130mm}{!}{ 
        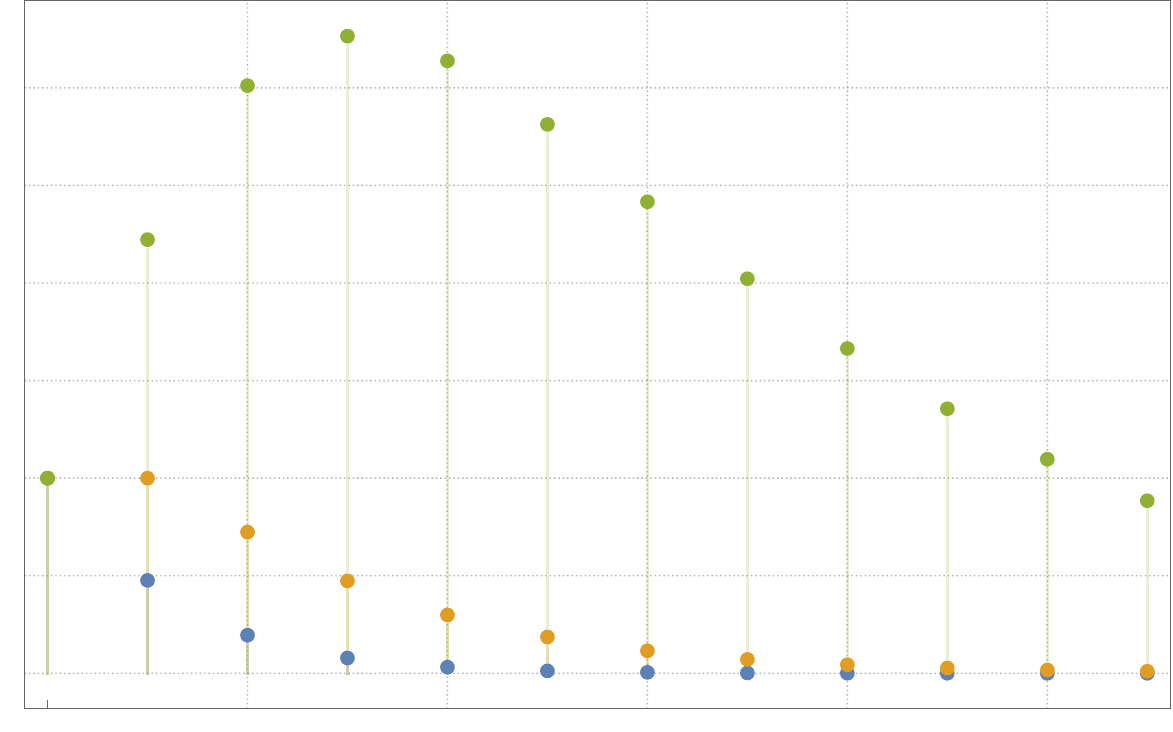}
      \caption{ $\cL_{m}(h)$ of
        $h=0.4$   (bottom),
        $h=\frac{-1+\sqrt{5}}{2}$  (middle),
        and $h=0.8$ (top)}
      \label{fig:pt-linear-half}
\end{figure}

Furthermore, consider the $n$-fold convolutions $\cL^{*n}$ of
$n\in \Zgeo$.  Then, we obtain all the metallic ratios as
the critical points of $\cL^{*n}$, solving
\begin{align}
  \cL^{*n}_{0}(q^{\fraa})
  &=1=\frac{nq^{\fraa}}{1-q}
    = \cL^{*n}_{1}(q^{\fraa})
\end{align}
for $0<\qq<1$.  The same holds for $\cQ^{*n}$, since
$\cQ_{0}=\cL_{0}$ and $\cQ_{1}=\cL_{1}$.
 
\subsection{Characterizations of $\cL$ and $\cQ$ with the
  golden ratio of quantum dilogarithms}
\label{sec:characterization}
As we are interested in $q$-polynomials, we define the
following notion.
\begin{definition}\label{def:intro-ideal}
  Let $\rho\in \Zgeo$.  Suppose a parcel
  $\cF=\Lam(l,w,f,q,u)$.  Then, we call $\cF$ {\it
    $\rho$-ideal}, if $\cF$ is $\rho$-merged-log-concave and
  each $m\in \Zgez$ gives
\begin{align}
  \Delta(\cF)(l,w,\rho,\tp{0,1},\tp{m,m+1},q)\llq 0.
\end{align}  
\end{definition}
 Section~\ref{sec:primal-monom} confirms that
if a monomial parcel $\cF_{\gam}$
 is $\rho$-ideal for some $\rho\in \Zgeo$, then
 $\cF_{\gam}$  is $\rho'$-ideal for any $\rho'\in \Zgeo$.
Hence, we simply call $\cF_{\gam}$ ideal,
if it is $\rho$-ideal for some $\rho\in \Zgeo$.

Consider all the ideal monomial parcels $\cF_{\gam}$.  Then,
$\cL$ is extremal among them by phase transitions,
because if $\cF_{\gam}$ has a phase transition, then 
Section~\ref{sec:primal-monom} proves that $\cL_{m}$ is the maximum
among $\frac{\cF_{\gam,m}}{\cF_{\gam,0}}$
for
each $m\in \Zgez$ and $0<q<1$; i.e., we have the
inequality
\begin{align}
  \cL_{m}=\frac{\cL_{m}}{\cL_{0}}
  = \frac{q^{\frac{m}{2}}}{(m)_{q}}
  \geq \frac{q^{\gam_{1,1}m^{2}+\gam_{1,2}m}}{(m)_{q}}
  = \frac{\cF_{\gam,m}}{\cF_{\gam,0}}
  \label{ineq:monomial-extremal}
\end{align}
in real numbers.
Likewise, we
prove that $\cQ$ is extremal among all the quadratic
$\cF_{\gam}$.

For each $\cF_{\gam}$, we put the generating function
\begin{align}
\cZ_{\cF_{\gam}}(q^{u_{\gam}},t)=\sum_{m\in \Zgez}
\cF_{\gam,m}(q^{u_{\gam}})t^{m}.
\end{align}
For instance,
by $u_{\tp{\tp{\fraa,-\fraa,0}}}=
u_{\tp{\tp{0,0,0}}}=1$,
we
have the generating functions
$\cZ_{\tcL_{\tp{\tp{0,0}}}}(q,t)
=\cZ_{\cF_{\tp{\tp{0,0,0}}}}(q,t)
=\sum_{m\in \Zgez}
\frac{t^{m}}{(m)_{q}}$ and
$\cZ_{\tcQ_{\tp{\tp{-\fraa,0}}}}(q,t)=
\cZ_{\cF_{\tp{\tp{\fraa,-\fraa,0}}}}(q,t)
=\sum_{m\in \Zgez}
\frac{q^{\frac{m(m-1)}{2}}t^{m}}{(m)_{q}}
$ for $\tcL_{\tp{\tp{0,0}}}$ and
$\tcQ_{\tp{\tp{-\fraa,0}}}$.

Suppose $\cF_{\gam,0}=1$. Section~\ref{sec:q-dilog}
confirms that {\it quantum
  dilogarithms}~\cite{FadKas,KonSoi} coincide with the
generating functions $\cZ_{\cF_{\gam}}(q^{u_{\gam}},t)$ of
$\cF_{\gam}$\footnote{The quantum dilogarithms
  $\cZ_{\cF_{\gam}}(q^{u_{\gam}},t)$ of \cite{FadKas} and
  \cite{KonSoi} are of quadratic and linear $\cF_{\gam}$,
  respectively.  By the change of variables
  $t\mapsto (1-q)t$, $\cZ_{\cF_{\gam}}(q^{u_{\gam}},t)$ are
  also said to be $q$-exponentials (c.f. \cite[equations
  (0.7.7) and (0.7.8)]{KoeSwa}).}.  It seems that we have
not obtained the golden ratio for quantum dilogarithms.
Thus, since $\cL$ is the extremal parcel among $\cF_{\gam}$,
it can be seen that we obtain the golden ratio of quantum
dilogarithms as the critical point of $\cL$.  Likewise, the
extremal $\cQ$ among quadratic $\cF_{\gam}$ gives the golden
ratio.

\subsection{On the full version of the
  merged-log-concavity}
We generalize the merged-log-concavity in
Definition~\ref{def:merged-simpler} to the full version in
Definition~\ref{def:merged}.  First, the full version allows
parcels of finitely-many rational functions and
merged determinants of polynomials with
non-negative coefficients for the merged-log-concavity.
Definition~\ref{def:merged-simpler} avoids these
 for simpler
notations.
Second, we generalize the ring shift factor
$\Ups(l,w,\rho,a,b,q)$ into
\begin{align}
  \phi(\qr)^{(\dif(a,b)_{1}+\tdif(a,b)_{1})w_{1}} 
  \frac{[b_{1}]!_{\qr}^{w_{1}}
  [b_{2}]!_{\qr}^{w_{1}}}{
  [a_{1}]!_{\qr}^{w_{1}} [a_{2}]!_{\qr}^{w_{1}}}
\end{align}
for an appropriate $\phi(\qr)\in \Q(\qr)$ such as
$\phi(\qr)=1-\qr$. Thus, we also extend a parcel
$\Lam(l,w,f,q,u)$ into
$\left\{\frac{f_{m}}{
    \phi(q)^{mw_{1}}[m]^{w_{1}}_{q}}\right\}_{m\in \Z}$, whose
denominators are not necessarily $(m)^{w_{1}}_{q}$.  In
particular, the real values of a parcel for $0<q^{u}<1$ depend
on the choice of $\phi$. Third, we formulate the full version
by multivariate rational functions, though we often take
rational functions of rational powers of $q$.  Hence, we
consider the variation of unimodal sequences on multiple
parameters in some general setting.  Fourth, we axiomatize
the partial orders in
Definition~\ref{def:q-pic-poly-q-Laurent-pic-poly} to
discuss positivities of rational functions more generally.

After Section~\ref{sec:intro}, we stick to the full version
to avoid confusion.

\subsection{Comparison with the strong
  \texorpdfstring{$q$}{q}-log-concavity
  and \texorpdfstring{$q$}{q}-log-concavity}
Sagan~\cite[Section 1]{Sag} gave the notion of {\it strong
  $q$-log-concavity} of polynomials.
Section~\ref{sec:weight-zero-merged-q-log-conc} explains
that strongly $q$-log-concave polynomials are essentially
weight-zero merged-log-concave parcels.
Furthermore, in a
suitable setting, strongly $q$-log-concave polynomials over
$(\lam)_{q}$ give merged-log-concave parcels.
However, the converse does not hold, as
$q^{\frac{\lam^{2}}{2}}$ of $\cQ$ are not strongly
$q$-log-concave by
equation~\eqref{eq:poly-not-pos-int-coeff}.

We also compare the $q$-log-concavity and the
merged-log-concavity.  This gives $q$-log-concave
polynomials from weight-zero merged-log-concave parcels.

\subsubsection{On Newton's log-concavities}
More explicitly, we compare the $q$-log-concavity and the
merged-log-concavity by the following {\it Newton's
  log-concavities} on polynomials.  Let $d\in \Zget$.
Then, for $\alp=\{\alp_{\lam}\in \R_{>0}\}_{\lam\in \oi(d)}$,
 consider $p=\{p_{\lam}\in \R_{>0}\}_{\lam\in \oi(0,d)}$ such that
\begin{align}
  \prod_{i\in\oi(d)}(1+\alp_{i} t)
  = \sum_{\lam\in\oi(0,d)}
p_{\lam}t^{\lam}.
\end{align}
Then, it holds that
$\frac{\lam}{\lam+1} \frac{d-\lam}{d+1-\lam}
p_{\lam}^{2}-p_{\lam-1}p_{\lam+1}\geq 0$ for all
$\lam\in\oi(d-1)$~\cite[p241--p243]{New} (see \cite[Lemma
1.1]{Bra} for a proof).  Moreover, since
$\frac{\lam}{\lam+1} \frac{d-\lam}{d+1-\lam}<1$, it holds
that $p_{\lam}^{2}-p_{\lam-1}p_{\lam+1}\geq 0$ for all
$\lam\in\oi(d-1)$, which we call Newton's log-concavity of
$\alp$.  Thus, 
that  of $\alp$ implies the
unimodality of $p$.

First,
let $g(h,d)=\{g(h,d)_{\lam}=h^{\lam}\}_{\lam\in\oi(0,d)}$
a finite geometric
sequence for $0<h<1$.
Then,
the $q$-log-concavity of $(-t;q)_{d}$ extends Newton's
log-concavity of  $g(h,d)$ by $q$-polynomials with positive
integer coefficients as follows.  By the $q$-binomial
theorem, consider $\cB(d,\lam,q)\in \Q[q]$ such that
\begin{align}
  (-t;q)_{d}
  &=\sum_{\lam\in \oi(0,d)}q^{\frac{\lam(\lam-1)}{2}}
    { d \brack \lam}_{q}t^{\lam}
    =\sum_{\lam\in \oi(0,d)}\cB(d,\lam,q)t^{\lam}.
\end{align}
Then, \cite[Theorem 3.2]{Sag} (see \cite{But,Kra}) gives
the  $q$-log-concavity:
\begin{align}
  \cB(d,\lam,q)^{2}-
  \cB(d,\lam-1,q)\cB(d,\lam+1,q)
  &\ggq 0.
    \label{ineq:newton-binom-qlog-quadratic}
\end{align}
By $q=h$, this implies  Newton's log-concavity
of $g(h,d)$, and hence the unimodal sequences
$\{\cB(d,\lam,h)\in \bR_{> 0}\}_{\lam\in \oi(0,d)}$
for $0<h<1$.

Second, we consider the infinite geometric sequences
$g(h,\infty)=\{g(h,d)_{\lam}=h^{\lam}\}_{\lam\in\oi(0,\infty)}$ for
$0<h<1$ by the merged-log-concavity.  This uses the
following equations of Euler~\cite[Chapter 16]{Eul}, which
are also of $q,t\in \bC$ and $|q|,|t|<1$.
\begin{definition}
  \label{def:euler-eq}
  In the ring of formal power series $\Q[[q,t]]$, we call
  $(t;q)^{-1}_{\infty} =\sum_{\lam\in \Zgez} \frac{t^{\lam}}
  {(\lam)_{q}}$ and
  $(-t;q)_{\infty} = \sum_{\lam\in \Zgez}
  \frac{q^{\frac{\lam(\lam-1)}{2}}} {(\lam)_{q}}t^{\lam}$
  Euler's linear and quadratic binomial identities,
  respectively.
\end{definition}
In particular, $\cZ_{\tcQ_{\tp{-\fraa,0}}}(q,t)= (-t;q)_{\infty}$.
Hence, the merged-log-concavity of $(-t;q)_{\infty} $ extends
Newton's log-concavity of finite $g(h,d)$ to infinite $g(h,\infty)$ by
$q$-polynomials with positive integer coefficients.  This
also gives the variation of almost strictly unimodal
sequences
$\{\tcQ_{\tp{-\fraa,0},m}(h)\in \R_{>0}\}_{m\in \Zgez}$ for
$0<h<1$.

Moreover, monomial parcels generalize $(-t;q)_{\infty}$.  In
particular,
$\cZ_{\tcL_{\tp{0,0}}}(q,t) =(t;q)_{\infty}^{-1}$ for the
monomial parcel $\tcL_{\tp{0,0}}$.  Also, monomial parcels
give more merged-log-concave parcels with strongly
$q$-log-concave polynomials by convolutions and certain
Hadamard products, which we discuss in
Section~\ref{sec:hadam}.  Hence, the merged-log-concavity
extends the $t$-power series of $(\pm t;q)_{\infty}^{\mp 1}$ to
some extent by polynomials with positive integer
coefficients and the variation of unimodal sequences.

\subsection{Multimonomial indices and
  eta products}
Euler's binomial identities are often used to compute the
$q$-expansions of {\it eta products} via Jacobi's triple
product identity~\cite{Koh}.  But, we consider multivariate
deformation of the eta products by Euler's binomial
identities and the merged-log-concavity. This approach on
eta products uses the notion of {\it multimonomial indices}
below.

\begin{definition}\label{def:multimonomial-ind}
  Let $d\in \Zgeo$, $w\in \Zogeo$,
  $\alp,\bta\in \Zgeo^{d}$, and
  $\gam\in \prod_{i\in\oi(d)}\Qr$.  We call
  $\tp{d,w,\alp,\bta,\gam}$ multimonomial index, if
  $\tp{1,w,\tp{\gam_{i}}}$ is a monomial index for each
  $i\in\oi(d)$.  We refer to $d$, $w$, $\alp$, $\bta$,
  and $\gam$ as the depth, weight, inner-exponent,
  outer-exponent, and core of $\tp{d,w,\alp,\bta,\gam}$.
\end{definition}

For a parcel $\cF_{w,\gam,q}$, we put the generating
function
\begin{align}
  \cZ_{w,\gam,q}(t)
  =\sum_{m\in \Zgez}\cF_{m}(q^{u_{\gam}})t^{m}.
  \label{eq:intro-gen-functions}
\end{align}
For instance,
$\cZ_{\tp{1},\gam,q}(t)=
\cZ_{\cF_{\gam}}(q^{u_{\gam}},t)$.
Furthermore, we define the following notions of graded monomial products and
monomial convolutions.
\begin{definition}\label{def:monomial-conv}
  Suppose a multimonomial index $\tp{d,w,\alp,\bta,\gam}$.
  Let $z=\tp{z_{i}}_{i\in\oi(d)}$ be a tuple of
  indeterminates.  Then,  we define
  the graded monomial product
    \begin{align}
      M(d,w, \alp,\bta, \gam, q,z)
      =
      \prod_{\lam\in\oi(d)}
      \left(
      \cZ_{w,\tp{\gam_{\lam}},q^{\alp_{\lam}}}(z_{\lam})
      \right)^{\bta_{\lam}}
      \in \Q(q^{u_{\alp,\gam}}) [[z_{1},\dots,z_{d}]].
    \end{align}
    In particular, for an indeterminate $v$, we define the
    monomial convolution
    \begin{align}
      \cM(d,w, \alp,\bta,\gam,q,v)
      =
        M(d,w, \alp,\bta, \gam, q,\iota^{d}(v))
       =\sum_{\lam\in\Zgez}
       \cM(d,w, \alp,\bta,\gam,q,v)_{\lam} v^{\lam}.
    \end{align}
  \end{definition}
  Also, let us recall the eta function and eta products.
 \begin{definition}\label{def:eta-fun-prod}
   For the imaginary unit $\i$, let
   $q=e^{2\pi \i \tau}$ of $\tau\in \bC$ such that $\Ima(\tau)>0$.
   Then,
   $\eta(\tau)= q^{\frac{1}{24}}(q;q)_{\infty}$     
   is the (Dedekind) eta function.
   Moreover,
    when
    $d\in \Zgeo$, $\alp\in \Zgeo^{d}$, and
    $\bta\in \Zd_{\neq 0}$,
    $\ep_{d,\alp,\bta}(\tau)=\prod_{\lam\in\oi(d)}
    \eta(\alp_{\lam}\tau)^{\bta_{\lam}}$
    is an eta product.
  \end{definition}

Then, we put the following notation to describe
some multimonomial indices and multivariate
deformation of eta products.
\begin{definition}\label{def:intro-hat-bta-gam}
  Let $d\in \Zgeo$,
$\alp\in \Zgeo^{d}$,
  $\bta\in \Zd_{\neq 0}$, and
  $\ka\in \Qd$.  Suppose a tuple
  $z=\tp{z_{\lam}}_{\lam\in\oi(d)}$ of indeterminates.
    \begin{enumerate}
    \item
      We put  a tuple
      $\etg(\bta,\ka)\in \prod_{\lam\in\oi(d)}\Qr$
      such that
      \begin{align}
        \etg(\bta,\ka)_{\lam}
        &=
          \begin{dcases}
            \tp{0,\ka_{\lam},-\frac{1}{24}}
            \mif \bta_{\lam}\leq -1,\\
            \tp{\fraa,\frac{\ka_{\lam}}{2},\frac{1}{24}}
            \melse.
          \end{dcases}
      \end{align}
    \item We put a tuple
      $J(z,q,\alp,\bta,\ka)
      =\tp{J(z,q,\alp,\bta,\ka)_{\lam}}_{\lam\in\oi(d)}$
      of indeterminates such that
     \begin{align}
       J(z,q,\alp,\bta,\ka)_{\lam}=
       \begin{dcases}
         q^{-(1-\ka_{\lam})\alp_{\lam}}\cdot z_{\lam} \mif \bta_{\lam}\leq -1,\\
         -q^{-\frac{1-\ka_{\lam}}{2}\alp_{\lam}}\cdot z_{\lam} \melse.
       \end{dcases}
     \end{align}
   \end{enumerate}
\end{definition}

Let $w=\tp{1}$.  Then, $\tp{1,w,\tp{\etg(\bta,\ka)_{\lam}}}$
for each $\lam\in \oi(d)$ is a monomial index. Thus, for each
$q=e^{2\pi \i \tau}$ of $\Ima(\tau)>0$, we obtain the graded
monomial product
$M(d,w, \alp,\abs{\bta}, \etg(\bta,\ka), q,z)$ as a
$J(z,q,\alp,\bta,\ka)$-analog of $\ep_{d,\alp,\bta}(\tau)$,
because Euler's binomial identities imply
\begin{align}
M(d,w, \alp,\abs{\bta}, \etg(\bta,\ka), q,z)= \ep_{d,\alp,\bta}(\tau)    
\end{align}
for the limit $J(z,q,\alp,\bta,\ka)\mapsto \iota^{d}(1)$ as
$z$ varies in $\bC^{d}$.

We prove that all monomial convolutions are generating
functions of merged-log-concave parcels by parcel
convolutions.  We also prove that all graded monomial
convolutions are generating functions of merged-log-concave
parcels of general widths.  Thus, the $J$-analogs
$M(d,w, \alp,\abs{\bta}, \etg(\bta,\ka), q,z)$ of
$\ep_{d,\alp,\bta}(\tau)$ give polynomials with positive integer
coefficients and the variation of almost strictly unimodal
sequences by Young diagrams.  We then conjecture the
unimodality of merged determinants on some eta products.

Furthermore, we define the notion of {\it weighted
  $q$-multinomial coefficients} for an indeterminate $q$.
This comes with weighted $q$-Pascal identities.  Then, we
explicitly describe monomial convolutions and their merged
determinants by the weighted $q$-multinomial
coefficients. This generalizes
Theorem~\ref{thm:intro-poly-positivity} on the change of
variables $q\mapsto q^{\alp_{\lam}}$.  We also propose a
conjecture on the generalized Narayana numbers \cite{Guy} by
merged determinants.

\subsection{Statistical-mechanical phase
  transitions by the
  merged-log-concavity}
\label{sec:intro-stat-pt}

Sections~\ref{sec:intro-stat-pt} and~\ref{sec:casimir}
consider the {\it grand canonical partition functions} of
some ideal boson-fermion gases.  We then use the
mathematical merged-log-concavity.  This gives some
statistical-mechanical phase transitions on the vacua that
have the lowest Helmholtz free energies.  In particular, we
are interested in how particles exist. Then, zero particle
vacua become non-zero particle vacua as the temperature
rises by the phase transitions.  The phase transitions are
some discontinuous changes in the systems of the ideal
boson-fermion gases by polynomials with positive integer
coefficients.

More explicitly, Section~\ref{sec:intro-stat-pt} considers
the phase transitions at the metallic ratios for ideal boson or
fermion gases by the generating functions $\cZ_{\cL}(\qq,t)$
and $\cZ_{\cQ}(\qq,t)$.
Furthermore, by
the monomial convolutions,
Section~\ref{sec:casimir} discusses
ideal (mixed)
boson-fermion gases with or without Casimir energies
(Ramanujan summation of zero-point energies).  In
particular, we obtain the phase transitions by the
merged-log-concavity in either case.

Unless stated otherwise, Section~\ref{sec:intro-stat-pt}
assumes the thermodynamic beta $\bta>0$ and chemical
potential $\mu<0$ with $q=e^{-\bta}$, $\mu'=-\mu \bta>0$, and
$t=e^{-\mu'}$.
Then, we have
\begin{align}
  0<\qq=e^{-\frac{\bta}{2}}, t=e^{-\mu'}<1. \label{ineq:qtone}
\end{align}
The author refers to~\cite[Chapter
1]{KapGal} for fundamental statistical-mechanical notions.

\subsubsection{Ideal boson gases}
\label{sec:intro-ideal-boson-gas}
Let us take the following operators
by the delta function $\del_{\lam,\lam'}$
of $\lam,\lam'\in \Q$ such that
$\del_{\lam,\lam'}= 1$ if $\lam=\lam'$ and 
$\del_{\lam,\lam'}=0$ otherwise.
\begin{definition}\label{def:boson-hamil-numb}
  For $\lam\in \Zgeo$, let
  $a_{b,\lam}, a_{b,\lam}^{\dagger}$ be the bosonic annihilation
  and creation operators that satisfy the commutator
  relations
  $[a_{b,\lam},a_{b,\lam'}^{\dagger}] =\del_{\lam,\lam'}$ and
  $[a^{\dagger}_{b,\lam},a^{\dagger}_{b,\lam'}]
  =[a_{b,\lam},a_{b,\lam'}]=0$.
  Also, for $v\in \Q$ and $\eps_{v,\lam}=\lam-v$,
  let $H_{b,v}$  and $N_{b}$ be the Hamiltonian and number
  operators such that
  \begin{align}
    H_{b,v}&=\sum_{\lam\in\Zgeo}
             \eps_{v,\lam} a_{b,\lam}^{\dagger}a_{b,\lam},
             \label{eq:hamil-bv}
    \\
    N_{b}&=\sum_{\lam\in\Zgeo}a_{b,\lam}^{\dagger}a_{b,\lam}.
  \end{align}
\end{definition}

In particular, we consider the bosonic system $B(1,\fraa)$
of
$H_{b,\fraa}$ and $N_{b}$.
Then,  $B(1,\fraa)$ represents an ideal bosonic gas
with the grand canonical partition function:
\begin{align}
  \cZ_{B(1,\fraa)}(q,t)
  =
  \Tr{e^{-\bta\left(H_{b,\fraa}-\mu N_{b}\right)}}
  =\Tr{e^{-\bta H_{b,\fraa}} \cdot e^{-\mu' N_{b}}}.
\end{align}

Let us write the partition function as the generating
function $\cZ_{\cL}(\qq,t)$.  Now, we have the eigenvalues
$n_{\lam}\in \Zgez$ of
$a_{b,\lam}^{\dagger}a_{b,\lam}$ by
the eivenvectors
$ |n_{\lam}\rangle=\frac{1}{\sqrt{n_{\lam}!}}
(a_{b,\lam}^{\dagger})^{n_{\lam}}|0\rangle$.  Hence, there is
the bosonic Fock space $F_{b}$ that has
$|\nn,\nt,\dots,n_{k},\dots\rangle$  such that
$\sum_{\lam\in \Zgeo}n_{\lam}<\infty$.
In particular, for
$\lam\in \Zgeo$ and $n_{\lam}\in \Zgez$,
\begin{align}
\langle \nn,\dots, n_{k},\dots|
  e^{-\bta H_{b,\fraa}}
  \cdot e^{-\mu' N_{b}}
  | \nn,\dots, n_{k},\dots\rangle
  = e^{-\bta \sum_{\lam} n_{\lam}\eps_{\fraa,\lam}}\cdot
    e^{-\mu' \sum_{\lam} n_{\lam}},  
\end{align}
and
$\left(1-q^{\lam-\fraa}\cdot t\right)^{-1} = \sum_{n_{\lam}}
e^{-\bta n_{\lam} \eps_{\fraa,\lam}} \cdot e^{-\mu' n_{\lam}}$.
Therefore, by Euler's linear binomial identity,
$\sum_{r\in \Zgez}\frac{q^{\frac{r}{2}}}{(r)_{q}}t^{r}
=\cZ_{\cL}(\qq,t)$ gives
  \begin{align}
  \cZ_{B(1,\fraa)}(q,t)
  =\cZ_{\cL}(\qq,t).
  \label{eq:linear-half-ideal-boson}
  \end{align}
  Let us mention that 
  this equality is known in physics~\cite[Chapter 1]{Dim}.

Moreover, for each $n\in \Zgeo$,
suppose that $B(n,\fraa)$ has $n$ sub-systems with
negligible interactions and $B(1,\fraa)$ represents each
sub-system.  Thus, $B(n,\fraa)$ realizes an ideal boson gas
with the grand canonical partition function
$\cZ_{\cL}(\qq,t)^{n}$.
\begin{remark}\label{rmk:ignore}
We ignore the overall factor
$\frac{1}{n!}$ for $\cZ_{\cL}(\qq,t)^{n}$, even
with identical particles on the sub-systems, because we are
interested in the unimodality of the $t$-power series of
$\cZ_{\cL}(\qq,t)^{n}$.  In
Sections~\ref{sec:intro-stat-pt}
and~\ref{sec:casimir} below, when we define a bosonic or
fermionic system that has sub-systems with negligible
interactions, we assume the same ignorance on its grand
canonical partition function.  
\end{remark}

\subsubsection{Ideal fermion gases}
\label{sec:intro-ideal-fermion-gas}
Similarly, we consider ideal fermion gases.
  \begin{definition}\label{def:ferm-hamil-numb}
    Let $a_{f,\lam}, a_{f,\lam}^{\dagger}$ of
    $\lam\in \Zgeo$ be the fermionic annihilation and
    creation operators that satisfy the anti-commutator
    relations
    $\{a_{f,\lam},a_{f,\lam'}^{\dagger}\}
    =\del_{\lam,\lam'}$
    and $\{a^{\dagger}_{f,\lam},a^{\dagger}_{f,\lam'}\}
    =\{a_{f,\lam},a_{f,\lam'}\}=0$.
    Also,
    for   $v\in \Q$,
    let $H_{f,v}$
    and $N_{f}$ be the Hamiltonian and
  number operators such that
  \begin{align}
    H_{f,v}&=\sum_{\lam\in\Zgeo}
             \eps_{v,\lam}a_{f,\lam}^{\dagger}a_{f,\lam},
             \label{eq:hamil-fv}
    \\
    N_{f}&=\sum_{\lam\in\Zgeo}a_{f,\lam}^{\dagger}a_{f,\lam}.
\end{align}
\end{definition}

Let us consider the fermionic system $F(1,\fraa)$ of $H_{f,\fraa}$ and $N_{f}$.  It represents an ideal fermion gas
with
  the grand canonical partition function
 \begin{align}
   \cZ_{F(1,\fraa)}(q,t)=
   \Tr{e^{-\bta H_{f,\fraa}} \cdot e^{-\mu' N_{f}}}
 \end{align}
 such that
 \begin{align}
\cZ_{F(1,\fraa)}(q,t)
   =\cZ_{\cQ}(\qq,t).
   \label{eq:quadratic-half-ideal-fermion}
 \end{align}

 Let us confirm
 equation~\eqref{eq:quadratic-half-ideal-fermion}.  For
 $\lam\in \Zgeo$, $n_{\lam}\in \{0,1\}$ are eigenvalues of
 $a_{f,\lam}^{\dagger}a_{f,\lam}$ on the eivenvectors
 $|0\rangle$ and $\alp_{f,\lam}^{\dagger}|0\rangle$ by the Pauli exclusion
 principle.  Thus, there is the fermionic Fock space $F_{f}$
 that consists of $| \nn,\dots, n_{k},\dots\rangle$ such that
 $\sum_{\lam\in \Zgeo}n_{\lam}<\infty$.  Also,
 $\cZ_{\cQ}(\qq,t) =\sum_{r\in
   \Zgez}\frac{q^{\frac{r^{2}}{2}}}{(r)_{q}} t^{r}$.  Hence,
 equation~\eqref{eq:quadratic-half-ideal-fermion} follows
 from Euler's quadratic binomial identity\footnote{ This
   does not require $|t|=|e^{-\mu'}|<1$~\cite[Lemma
   1.2]{Koh}.  But, we assume $\mu'>0$, i.e., $\mu<0$, for
   simplicity. See also \cite[Fig 1 and 2]{Cow} for $\mu$ in
   high temperatures. }, because
 $n_{\lam}\in \{0,1\}$ and  $\lam\in \Zgeo$ imply
 $1+q^{\lam-\fraa}\cdot t =
 \sum_{n_{\lam}}e^{-\bta n_{\lam}
   \eps_{\fraa,\lam}}\cdot
 e^{-\mu' n_{\lam}}$, and
\begin{align}
\langle \nn,\dots, n_{k},\dots|
  e^{-\bta H_{f,\fraa}}
    \cdot e^{-\mu' N_{f}}
  | \nn,\dots, n_{k},\dots\rangle
  = e^{-\bta \sum_{\lam} n_{\lam}\eps_{\fraa,\lam}}\cdot
  e^{-\mu' \sum_{\lam} n_{\lam}}.  
\end{align}

Moreover, assume that for each $n\in \Zgeo$,
$F(n,\fraa)$
has $n$
sub-systems with negligible interactions and $F(1,\fraa)$
represents each sub-system.  Hence, $F(n,\fraa)$ realizes an
ideal fermion gas with the grand canonical partition
function $\cZ_{\cQ}(\qq,t)^{n}$.

\subsubsection{Phase transitions
  of free energies and almost strictly unimodal sequences}
We discuss phase transitions of $\cZ_{\cL}(\qq,t)^{n}$ and
$\cZ_{\cQ}(\qq,t)^{n}$ by the series
expansions:
\begin{align}
  \cZ_{\cL}(\qq,t)^{n}
  &=\sum_{\lam\in\Zgez}Z_{\cL,n,\lam}(\qq) t^{\lam};
    \label{eq:t-po-1}
    \\\cZ_{\cQ}(\qq,t)^{n}
  &=\sum_{\lam\in\Zgez}Z_{\cQ,n,\lam}(\qq)t^{\lam}.
    \label{eq:t-po-2}
\end{align} 
These $Z_{\cL,n,\lam}(\qq)$ and $Z_{\cQ,n,\lam}(\qq)$
 of  particle numbers $\lam$
are the {\it canonical
  partition functions} of  $B(n,\fraa)$ 
and $F(n,\fraa)$.
In particular,  $n\in \Zgeo$ and $\lam\in \Zgez$ give
positive real values $Z_{\cL,n,\lam}(\qq)=\frac{q^{\frac{\lam}{2}}}{(\lam)_{q}}$
and $Z_{\cQ,n,\lam}(\qq)=\frac{q^{\frac{\lam^{2}}{2}}}{(\lam)_{q}}$
by
inequalities~\eqref{ineq:qtone}. 
Thus, we consider the  sequences:
\begin{align}
  z_{\cL,n}(\qq)
  &=\{z_{\cL,n,\lam}(\qq)
    =\log (Z_{\cL,n,\lam}(\qq))\in \R\}_{\lam\in\Zgez};\\
  z_{\cQ,n}(\qq)
  &=\{z_{\cQ,n,\lam}(\qq)
    =\log (Z_{\cQ,n,\lam}(\qq))\in \R\}_{\lam\in\Zgez}.
\end{align}

We have $\cZ_{\cL}(\qq,t)^{n}=\cZ_{\cL^{*n}}(\qq,t)$ and
$\cZ_{\cQ}(\qq,t)^{n}=\cZ_{\cQ^{*n}}(\qq,t)$.  Hence, by
critical points in the mathematical sense of
Section~\ref{sec:intro-unim-seq-crit-pts},
Theorem~\ref{thm:intro-poly-positivity} gives the following
phase transitions of the sequences
$z_{\cL,n}(\qq)$ and $z_{\cQ,n}(\qq)$.

\begin{theorem}\label{thm:intro-phase-transition}
  (Corollary~\ref{cor:conv-prob-linear-half-quad-half}
  by the terminology of Section~\ref{sec:intro-stat-pt})
  Fix $n\in \Zgeo$.
  \begin{enumerate}
  \item 
    Then, we have the almost strictly unimodal sequences
    $z_{\cL,n}(\qq)$ of $0<\qq<1$ with the critical point $0<c_{\cL,n}<1$
    as follows.
    \begin{enumerate}
    \item For each
       $0<\qq< c_{\cL,n}$, the mode
    $m_{\cL,n}(\qq)=0$ gives the strictly decreasing sequence:
    $z_{\cL,n,m_{\cL,n}(\qq)}(\qq)> z_{\cL,n,1}(\qq)
      > z_{\cL,n,2}(\qq)>\dots$.
  \item If $\qq=c_{\cL,n}$, then
  $m_{\cL,n}(\qq)=0$ gives
  the hill-shape and
  decreasing sequence:
  $z_{\cL,n,m_{\cL,n}(\qq)}(\qq)= z_{\cL,n,1}(\qq)
  > z_{\cL,n,2}(\qq)>\dots$.
\item For each $1>\qq>c_{\cL,n}$,
  $m_{\cL,n}(\qq)\in \Zgeo$ gives the two-slope hill-shape sequence:
    $z_{\cL,n,0}(\qq)<\dots <z_{\cL,n,m_{\cL,n}(\qq)}
    (\qq)
    \geq z_{\cL,n,m_{\cL,n}(\qq)+1}(\qq)>\dots$.      
    \end{enumerate}
  \item Also, we have the same for the almost strictly unimodal
    sequences $z_{\cQ,n}(\qq)$ with the critical point
    $0<c_{\cQ,n}=c_{\cL,n}<1$.
\item In particular, we have the metallic ratio
      \begin{align}
           c_{\cL,n}=c_{\cQ,n}=\frac{-n+\sqrt{n^{2}+4}}{2},
      \end{align}
      which is the golden ratio for $n=1$.
      \end{enumerate}
\end{theorem}

We now recall the
{\it Helmholtz free energies} (or free energies for short)
$A_{\cL,n,\lam}(\qq)$ and $A_{\cQ,n,\lam}(\qq)$ of the
canonical partition functions $ Z_{\cL,n,\lam}(\qq)$ and
$Z_{\cQ,n,\lam}(\qq)$ such that
\begin{align}
  A_{\cL,n,\lam}(\qq)
  &=-\frac{\log (Z_{\cL,n,\lam}(\qq))}{\bta}
    =-\frac{z_{\cL,n,\lam}(\qq)}{\bta},\\
  A_{\cQ,n,\lam}(\qq)
  &=-\frac{\log (Z_{\cQ,n,\lam}(\qq))}{\bta}
  =-\frac{z_{\cQ,n,\lam}(\qq)}{\bta}.
\end{align}
By
inequality~\eqref{ineq:det-positivity}
of 
 $q$-polynomials,
the free energies
satisfy the following inequalities of real numbers:
\begin{align}
  \sum_{i\in \{-1,1\}}
  A_{\cL,n,\lam}(\qq)
  -
  A_{\cL,n,\lam+i}(\qq)
  &<0,\\
\sum_{i\in \{-1,1\}}
  A_{\cQ,n,\lam}(\qq)
  -
  A_{\cQ,n,\lam+i}(\qq)
  &<0.
\end{align}
Therefore,
the merged determinants
of $\cL$ and $\cQ$ give to these inequalities
a  structure
in $\Z[q]$, which
is the smallest $\Z$-ring consisting
of the
base variable $q$.
For example, we obtain the $q$-polynomials
in Theorem~\ref{thm:intro-poly-positivity}
for these inequalities.

Moreover, suppose an ideal monomial parcel $\cF_{\gam}$ with
a phase transition. Also, let $\cF_{\gam,0}=1$.  Then, in
Section~\ref{sec:casimir}, we confirm the generating
function $\cZ_{\cF_{\gam}}(q^{u_{\gam}},t)$ as the grand
canonical partition function of $B(1,1-\gam_{1,2})$ if
$\gam_{1,1}=0$ or $F(1,\fraa-\gam_{1,2})$ if
$\gam_{1,1}=\fraa$.  Hence,
inequality~\eqref{ineq:monomial-extremal} means that
$A_{\cL,1,\lam}(\qq)$ are the lowest among the free energies
$A_{\cF_{\gam},1,\lam}(q^{u_{\gam}})$; i.e., $0<q<1$ and
$\lam\in \Zgez$ satisfy
\begin{align}
  A_{\cL,1,\lam}(\qq)\leq A_{\cF_{\gam},1,\lam}(q^{u_{\gam}}).
\end{align}
Also, $\cQ$ satisfies the same inequalities among quadratic
$\cF_{\gam}$, though
$A_{\cL,1,\lam}(\qq)<A_{\cQ,1,\lam}(\qq)$ for $\lam\in \Zget$.

Let
$A_{\cL,n}(\qq) =\{A_{\cL,n,\lam}(\qq)\}_{\lam\in \Zgez}$ and
$A_{\cQ,n}(\qq) =\{A_{\cQ,n,\lam}(\qq)\}_{\lam\in \Zgez}$.
Then, by $\bta>0$, Theorem~\ref{thm:intro-phase-transition}
implies the following phase transitions on the particle
numbers of the lowest free energies. These phase transitions
give non-zero particle vacua of the free energies as
the temperature rises.

\begin{corollary}\label{cor:transitions}
  Suppose  $n\in \Zgeo$.  Then, 
  the bosonic free energies $A_{\cL,n}(\qq)$  satisfy
  $\min(A_{\cL,n}(\qq)) =A_{\cL,n,m_{\cL,n}(\qq)}(\qq)$ by
  $m_{\cL,n}(\qq)=0$ for $0<\qq\leq c_{\cL,n}$ and
  $m_{\cL,n}(\qq)\in \Zgeo$ for $c_{\cL,n}<\qq<1$.
  Also,
the fermionic free energies
$A_{\cQ,n}(\qq)$ satisfy
  the same.
\end{corollary}

By $1>\qq=e^{-\frac{\bta}{2}}>c_{\cL,n}$, we have
$0<\bta<-2 \log \left(c_{\cL,n}\right)$.
Moreover, $\bta\propto \frac{1}{T}$ for the temperature $T$.
Thus, suppose a high temperature $T$ such that
$1>\qq>c_{\cL,n}$.  Then, we obtain the positive particle
numbers $m_{\cL,n}(\qq), m_{\cQ,n}(\qq)\in \Zgeo$ such that
$ A_{\cL,n,m_{\cL,n}(\qq)}(\qq) =\min(A_{\cL,n}(\qq))$ and
$ A_{\cQ,n,m_{\cL,n}(\qq)}(\qq) =\min(A_{\cQ,n}(\qq))$ in
Corollary~\ref{cor:transitions}.

For instance, by $\bta=\frac{1}{k_{B} T}$ of the Boltzmann
constant $k_{B}$, the temperatures $T_{n}$ of critical
points $c_{\cL,n}=c_{\cQ,n}$ satisfy
\begin{align}
      T_{n}&=\frac{1}{
             -2 \log\left(\frac{-n+\sqrt{n^{2}+4}}{2}\right)
    \cdot ( 1.38\cdots \times 10^{-23})}.
    \end{align}
    Then, by the Kelvin unit $K$, 
    \begin{align}
      T_{1}&= 0.75\cdots \times 10^{23} K,\\
      T_{2}&= 0.41\cdots \times 10^{23} K,\\
      &...
    \end{align}
    By the monomial convolutions in
    Section~\ref{sec:casimir}, we have different critical
    points and the corresponding temperatures on the ideal
    boson and fermion systems $B(n,\ka)$ and $F(n,\ka)$ of
    $n\in \Zgeo$ and $\ka\in \Q_{<1}$, solving
    $nq^{1-\ka}=1-q$ for $0<q^{1-\ka}<1$.

\begin{figure}[H]
  \centering     
  \resizebox{130mm}{!}   
{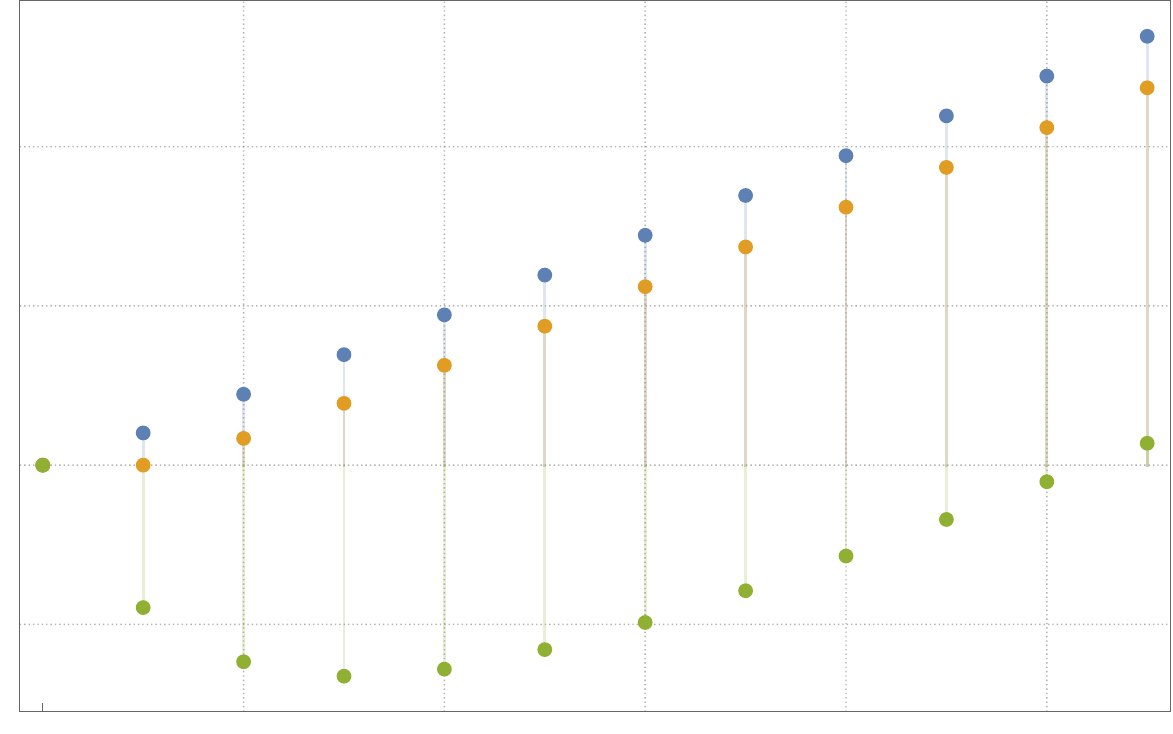} 
  \caption{$A_{\cL,1,\lam}(\qq)$ of
    $\qq=0.4$ (top),
    $\qq=c_{\cL,1}$ (middle), and $\qq=0.8$
    (bottom) }
  \label{fig:pt-linear-half-free-energies}
\end{figure}

    We now reinterpret Figure~\ref{fig:pt-linear-half}.
    Figure~\ref{fig:pt-linear-half-free-energies}
    illustrates the phase transition of the bosonic free
    energies $A_{\cL,1}(\qq)$ with the zero and non-zero
    particle vacua by the golden ration $c_{\cL,1}$ in
    Corollary~\ref{cor:transitions}.  We have a similar
    figure for the fermionic free energies $A_{\cQ,1}(\qq)$.
    \vspace{1em} 
     
Bose-Einstein condensations differ from the phase
transitions in Corollary~\ref{cor:transitions} or
Figure~\ref{fig:pt-linear-half-free-energies}, because
non-zero particle vacua keep appearing for higher
temperatures.  More precisely, for each $\lam \in \Zgeo$, we
obtain the $\lam$-particle vacuum at the temperature that
solves $\qq=1-q^{\lam}$.

It is natural to consider many states at high
temperatures. But, instead of  the infinite product
$(-t\qq;q)_{\infty}$, there is the following finite
analog to the phase transition on
$A_{\cQ,1}(\qq)$. Consider
\begin{align}
  (-t\qq;q)_{2}= 1+(\qq+q^{\frac{3}{2}})t+q^{2}t^{2}
  =\sum_{i\in\oi(0,2)}F(\qq)_{i}t^{i}.        
\end{align}
Then, the sequences $F(\qq)=\{F(\qq)_{i}\in \bR_{>0}\}_{i\in\oi(0,2)}$ of
$0<\qq<1$ are strictly log-concave by the discriminant of
the $t$-polynomial
$(-t\qq;q)_{2}$ (c.f.
inequalities~\eqref{ineq:newton-binom-qlog-quadratic}).  We
have the critical point $c=0.68233\dots$ that solves
$1=\qq+q^{\frac{3}{2}}$. Therefore, $F(\qq)$ become strictly decreasing
for $0<\qq <c$, hill-shape and decreasing for $\qq=c$,
and two-slope hill-shape for $c< \qq<1$.


\section*{Acknowledgments}
The author thanks Professors Tom Braden, Francesco Brenti,
Bruno Carneiro da Cunha, Sergey Galkin, Seikou Kato, Anatol
N. Kirillov, Yoshinori Mishiba, Aníbal L. S. Netto, Takuya
Okuda, Bruce E. Sagan, Fernando A. N. Santos, and Richard
P. Stanley for their helpful communications.  Among them,
Professor Takuya Okuda encouraged the author to study much
further at a very early stage of this work.  Also, Professor
Francesco Brenti responded to the author in an email at a
very early stage of this work, suggesting several
references.  This work was supported by the Research
Institute for Mathematical Sciences, a joint usage/research
center located in Kyoto University.  In particular, the
author thanks Professor Anatol N. Kirillov for reading
preliminary versions of this manuscript and providing
helpful comments in Kyoto during the author's two
several-day visits in 2017 and 2018. The author thanks
Departamento de Matem\'atica of Universidade Federal de
Pernambuco in Recife for their excellent support to his stay
from April 2018 to March 2019.
  
\section{Notations for families
  and rational functions}
\label{sec:prelim-notation}

We now fix some notations to develop our theory of the
merged-log-concavity.

\subsection{Families}

Suppose a family $F$ for an index set $I$.  We indicate each
element of $F$ by $F_{i}$ for each $i\in I$ so that
$F=\{F_{i}\}_{i\in I}$. Also, if $F=G$ for a family
$G=\{G_{i}\}_{i\in J}$, then $I=J$ and each $G_{i}=F_{i}$.
When $I\subset \Z$, we also call $F$ sequence.
In particular,
a tuple $F$ satisfies
 $I=\oi(d)$ of some
$d\in \Zgeo$ in this manuscript.

\subsection{Polynomials, Laurent
  polynomials, and rational functions}
\label{sec:polynomials-Laurent-rational-functions}

Unless stated otherwise, let
$\fX=\{X_{1},\dots, X_{L}\}=\{X_{i}\}_{i\in\oi(L)}$ be a
finite set of indeterminates (symbols) for some
$L\in \Zgeo$.  As we call $\fX$
set, if a set $\fX'$ has the same elements as $\fX$, then
$\fX'=\fX$, even with the indices $i\in\oi(L)$.
 
Let $Q$ be a commutative ring.
Then, $Q[[\fX]]=Q[[X_{1},\dots,X_{L}]]$,
$Q[\fX]=Q[X_{1},\dots, X_{L}]$,
and
$Q[\fXpm] =Q[X_{1}^{\pm 1},\dots, X_{L}^{\pm 1}]$
denote the ring of formal power series, 
multivariate polynomials, 
and Laurent polynomials, respectively.
We write an element $f$ in each of these rings
as
$f= \sum_{j\in \Z^{L}}f_{j_{1},j_{2},\dots,
  j_{L}}X_{1}^{j_{1}} X_{2}^{j_{2}} \dots X_{L}^{j_{L}}$ 
possibly with zero $f_{j_{1},j_{2},\dots, j_{L}}\in Q$.
If $0\neq f \in Q[\fXpm]$, then
$\deg_{X_{i}}f\in \Z$
and $\ord_{X_{i}}f\in \Z$ 
denote the degree and the order of $f$ as the
Laurent $X_{i}$-polynomial
of the ring $Q[(\fX\setminus \{X_{i}\})^{\pm}]$.
Let
$-\deg_{X_{i}}(0)=\ord_{X_{i}}(0)=\infty$ for our convention.

If $Q$ is field, then $Q(\fX)=Q(X_{1},\dots, X_{L})$ is the
field of rational functions.
We often look at real values of rational functions.  Let
$Q\subset \R$, $f\in Q(\fX)$, and
$r=\tp{r_{1},\dots, r_{L}}\in \R^{L}$. Then, we consider
$f(r)=f'(r)\in \R$, provided $f=f'\in Q(\fX)$ and
$f'(r)\in \R$. For a family
$F=\{F_{i}\in \Q(\fX)\}_{i\in I}$ and $r\in \R^{L}$, we
put $F(r)=\{F_{i}(r)\in \R\}_{i\in I}$.

If needed, we write
$\fXo=\{X_{1,i}\}_{i\in\oi(L_{1})},
\fXt=\{X_{2,i}\}_{i\in\oi(L_{2})}, \dots$ for multiple finite
sets of indeterminates.  Also, if
we write
$\fXt\subset Q(\fXo)$, then we assume $\fXt$ algebraically
independent over $Q$, unless stated otherwise.

\section{Squaring orders}
\label{sec:squaring-orders}

We introduce the notion of squaring orders on rational
functions. This is to discuss not only polynomials with positive
integer coefficients, but also  positive real values of
rational functions.  We use the following notion
(see~\cite{GilJer}).

\begin{definition}\label{def:partial-order} Suppose a set $R$.
  \begin{enumerate}
  \item A binary relation $\sce$ on $R$ is called a partial
    order, if $\sce$ satisfies the following conditions.
      \begin{enumerate}
      \item $f\sce f$ for each $f\in R$
        (reflexivity).
      \item $f_{1}\sce f_{2}$
        and $f_{2}\sce f_{3}$ imply $f_{1}\sce f_{3}$
        (transitivity).
    \item $f_{1}\sce f_{2}$ and
      $f_{2}\sce f_{1}$ imply $f_{1}= f_{2}$
      (antisymmetricity).
    \end{enumerate}
    Suppose a binary relation $\scc$ on $R$.
    If $f\scc f$ never holds (irreflexivity) and
    $\succ$ has the transitivity, then $\scc$ is called a
    strict partial order on $R$.     We also refer to a partial order
    and strict partial order as an inequality and 
    strict inequality, if no confusion occurs.
  \item Let $R$ be a ring. Assume a partial order $\sce$ on
    $R$. Then, $R$ is called a partially ordered ring of
    $\sce$ (or $\sce$-poring for short), when $R$ satisfies
    the following conditions.
     \begin{enumerate}
     \item $f_{1}\sce f_{2}$ and $f_{3}\in R$ imply
       $f_{1}+f_{3}\sce f_{2}+f_{3}$ (additivity).
     \item
        $f_{1}\sce 0$ and $f_{2} \sce 0$ imply
        $f_{1}f_{2}\sce 0$ (multiplicativity).
      \end{enumerate}
      Similarly, if a strict partial order $\scc$ on $R$
      satisfies the additivity and multiplicativity, then
      $R$ is called a strictly partially ordered ring of
      $\scc$ (or strict $\scc$-poring).
    \end{enumerate}
\end{definition}

Then, we recall the following properties of porings and strict
porings.  We write down a full proof for the completeness of
this manuscript.

\begin{lemma}\label{lem:succ-succeq}
  The following statements hold, if $R$ is a $\sce$-poring
  and strict $\scc$-poring.
  \begin{enumerate}
  \item\label{c:succ-succeq-sinv1} $f\sce g$ is equivalent to
    $f-g\sce 0$.
  \item\label{c:succ-succeq-stwoadd1} $f_{1}\sce f_{2}$
    and $g_{1}\sce g_{2}$ imply
    $f_{1}+g_{1}\sce f_{2}+g_{2}$.
  \item \label{c:succ-succeq-stwomul1} $f_{1}\sce f_{2}$ and
    $g \sce 0$ imply $f_{1}g\sce f_{2}g$.
  \item
    \label{c:succ-succeq-zero-zero-1}
    $f_{1}\sce f_{2}\sce 0$ and $g_{1}\sce g_{2}\sce 0$
    imply $f_{1} g_{1}\sce f_{2}g_{2}\sce 0$.
  \item\label{c:succ-succeq-sinv2} $f\scc g$ is equivalent
    to $f-g\scc 0$.
  \item \label{c:succ-succeq-stwoadd2} $f_{1}\scc f_{2}$ and
    $g_{1}\scc g_{2}$ imply $f_{1}+g_{1}\scc f_{2}+g_{2}$.
  \item
    \label{c:succ-succeq-stwomul2}
    $f_{1}\scc f_{2}$ and $g \scc 0$ imply
    $f_{1}g\scc f_{2}g$.
  \item
    \label{c:succ-succeq-zero-zero-2}
    $f_{1}\scc f_{2}\scc 0$ and $g_{1}\scc g_{2}\scc 0$
    imply $f_{1} g_{1}\scc f_{2}g_{2}\scc 0$.
  \end{enumerate}
\end{lemma}
\begin{proof}
  Let us prove Claim~\ref{c:succ-succeq-sinv1}. By the
  additivity and $-g\in R$, $f\sce g$ implies
  $f-g\sce g-g=0$. Also, $f-g\sce 0$ gives
  $f=f-g+g\sce 0+g=g$ by $g\in R$.
  
  Let us confirm Claim~\ref{c:succ-succeq-stwoadd1}.  By
  $g_{1}\in R$, $f_{1}+g_{1}\sce f_{2}+g_{1}$.  Also,
  $f_{2}+g_{1}\sce f_{2}+g_{2}$ by $f_{2}\in R$.  Since
  $f_{1}+g_{1}\sce f_{2}+g_{1}\sce f_{2}+g_{2}$,
  Claim~\ref{c:succ-succeq-stwoadd1} follows from the
  transitivity of $\sce$.
  
  Let us prove Claim~\ref{c:succ-succeq-stwomul1}.  By
  Claim~\ref{c:succ-succeq-sinv1}, $f_{1}\sce f_{2}$ implies
  $f_{1}-f_{2}\sce 0$.  Thus,
  $(f_{1}-f_{2})g=f_{1}g-f_{2}g\sce 0$ by $g\sce 0$ and the
  multiplicativity of $\sce$. Then, $f_{1}g\sce f_{2}g$ by
  Claim~\ref{c:succ-succeq-sinv1}.
  
  Let us prove Claim~\ref{c:succ-succeq-zero-zero-1}.  By
  Claim~\ref{c:succ-succeq-stwomul1}, $f_{1}\sce f_{2}$
  gives $f_{1}g_{1}\sce f_{2}g_{1}$.  Since
  $g_{1}\sce g_{2}$ implies $f_{2}g_{1}\sce f_{2}g_{2}$, we
  have $f_{1}g_{1}\sce f_{2}g_{1}\sce f_{2}g_{2}$.  Then,
  $f_{1}g_{1}\sce f_{2}g_{2}$ by the transitivity.  Thus,
  $f_{2}g_{2}\sce 0$ by $f_{2},g_{2}\sce 0$ and the
  multiplicativity.

  We also obtain Claims~\ref{c:succ-succeq-sinv2},
  ~\ref{c:succ-succeq-stwoadd2},
  ~\ref{c:succ-succeq-stwomul2},
  and~\ref{c:succ-succeq-zero-zero-2},
 replacing $\sce$ by $\scc$ in the above.
\end{proof}

Let us consider the following sets of binary relations.

\begin{definition}\label{def:succ-site}
  Suppose a binary relation $\geq$ on a set $R$.  Then, let
  $\StgeR=\{f\in R\mid f\geq g \mor g\geq f \mforsome g\in R\}$.
\end{definition}
  
Then, we introduce squaring orders of a ring by two subrings
and four binary relations.
      
\begin{definition}\label{def:succ-gen}
  Suppose a ring $R$.
  \begin{enumerate}
  \item Assume binary relations $\geq,>$ and $\sce,\scc$ on
    $R$ with the following conditions.
    \begin{enumerate}
    \item $f>g$ implies $f\geq g$ ($>$-$\geq$ implication).
    \item $f \scc g$ implies $f\sce g$ ($\scc$-$\sce$
      implication).
    \item $f\sce 0$ implies $f\geq 0$ (half $\sce$-$\geq$
      implication).
    \item $f\scc 0$ implies $f> 0$ (half $\scc$-$>$
      implication).
    \end{enumerate}
    We refer to the four implications above as the squaring
    implications of $\geq,>$ and $\sce,\scc$.
  \item Moreover, assume the following conditions.
    \begin{enumerate}
    \item $f \sce g \scc h$ or $f \scc g \sce h$ implies
      $f\scc h$ (semi-strict transitivity).
    \item There are the $\geq$-poring
      $\StgeR$ and strict
      $>$-poring $\StgR$ such that $\StgeR=\StgR$
      ($\geq$-$>$-poring equality).
    \item There are the $\sce$-poring $\StceR$ and strict
      $\scc$-poring $\StccR$ such that
      $\StceR=\StccR$ ($\sce$-$\scc$-poring equality).
    \item $\StceR\subset \StgeR$ ($\sce$-$\geq$-poring
      inclusion).
    \end{enumerate}
    Then, we call $\sce,\scc$ squaring orders on
    $\tp{R,\geq,>}$. Also, we call $\scc$ strict squaring order
    of $\sce$.  For simplicity, when we write $\sce,\scc$ as
    squaring orders on $\tp{R,\geq,>}$, we assume that $\scc$
    is a strict squaring order of $\sce$.
  \end{enumerate}
\end{definition}

We use the terminology ``squaring orders'' by the following
implication diagram (``squaring'' or ``square diagram +
ring'' orders):
\begin{center}
  \begin{tikzpicture}[scale=.7]
    \node (a) at (0,2) {$\scc$};
    \node (b) at (-1,1) {$>$};
    \node (c) at (1,1) {$\sce$};
    \node (d) at (0,0) {$\geq$};
    \draw[->] (a) -- (b);
    \draw[->] (a) -- (c);
    \draw[->] (b)--(d);
    \draw[->] (c)--(d);
  \end{tikzpicture}
\end{center}
Also, notice that $\sce$ is
not necessarily ``larger than or equal to'', because
$f\sce g$ and $f\neq g$ is not necessarily the same as
$f\scc g$.
Moreover,
let us confirm the following lemma.
\begin{lemma}\label{lem:half-full-implication}
  Suppose binary relations $\geq,>$ on a ring $R$ with the
  $>$-$\geq$ implication and $>$-$\geq$-poring equality.
  Also, 
  suppose binary relations $\sce, \scc$ on $R$ with the
  $\sce$-$\geq$-poring inclusion and $\sce$-$\scc$-poring
  equality.
  \begin{enumerate}
  \item \label{c:half-full-implications-sce-geq} The half
    $\sce$-$\geq$ implication gives the $\sce$-$\geq$
    implication.
  \item \label{c:half-full-implications-scc->} 
    The half $\scc$-$>$ implication gives the $\scc$-$>$
    implication.
   \end{enumerate}
\end{lemma}
\begin{proof}
  Let us prove Claim~\ref{c:half-full-implications-sce-geq}.
  Let $f\sce g$ so that $f,g\in \StceR$.  Since $\StceR$
  is a $\sce$-poring by the $\sce$-$\scc$-poring equality,
  the additivity of $\sce$ gives
  $f-g\sce 0$.  Hence,
  $f-g\geq 0$ by the half $\sce$-$\geq$ implication.  Also, 
   $f,g \in \StgeR$ by the $\sce$-$\geq$-poring
  inclusion.  Thus, $f\geq g$ holds, since $\StgeR$ is the
  $\geq$-poring.

  Let us prove Claim~\ref{c:half-full-implications-scc->}.
  Let $f\scc g$ so that $f,g\in \StccR$.  Then, the
  additivity of $\scc$ gives $f-g\scc 0$, since $\StccR$ is
  a strict $\scc$-poring by the $\sce$-$\scc$-poring
  equality.  Thus, $f-g> 0$ by the half $\scc$-$>$
  implication.  Also, we have
  $f,g \in \StccR=\StceR\subset \StgeR=\StgR$ by the
  $\sce$-$\scc$-poring equality, $\sce$-$\geq$-poring
  inclusion, and $\geq$-$>$-poring equality. Therefore,
  $f>g$ by the additivity of $>$ on $\StgR$.
\end{proof}

Then, we have the following implications of usual
inequalities by squaring orders.

\begin{proposition}\label{prop:succ}
  Squaring orders $\sce,\scc$ on $\tp{R,\geq,>}$ satisfy the
  following.
  \begin{enumerate}
  \item $f_{1}\sce f_{2}$ and $g_{1}\scc g_{2}$ imply
    $f_{1}+g_{1}\scc f_{2}+g_{2}$.
    \label{c:succ-sce-scc-sum}
  \item $f_{1}\scc f_{2}\sce 0$ and
    $g_{1}\scc g_{2}\sce 0$ imply
    $f_{1} g_{1}\scc f_{2}g_{2}\sce 0$.
    \label{c:succ-scc-sce-mult}
  \item $f_{1}\sce f_{2}\scc 0$ and $g_{1}\scc g_{2}\sce 0$
    imply $f_{1} g_{1}\scc f_{2}g_{2}\sce 0$.
    \label{c:succ-sce-scc-mult}        
  \end{enumerate}
\end{proposition}
\begin{proof}
  Let us confirm Claim~\ref{c:succ-sce-scc-sum}.  Since the
  $\scc$-$\sce$ implication gives $g_{1}\sce g_{2}$, we have
  $g_{1}\in \StceR$.  Thus, since $f_{1}\sce f_{2}$, the
  additivity of $\sce$ implies
  \begin{align}
    f_{1} +g_{1}&\sce f_{2} +g_{1}.\label{ineq:succ-p1}
  \end{align}
  Moreover, $f_{1}\sce f_{2}$ yields $f_{2}\in \StccR$ by
  the $\sce$-$\scc$-poring equality.
  Thus, since $g_{1}\scc g_{2}$, the additivity of $\scc$ gives
  \begin{align}
    f_{2}+g_{1}&\scc f_{2}+g_{2}.\label{ineq:succ-p2}
  \end{align}
  Then, $f_{1}+g_{1}\sce f_{2}+g_{1}\scc f_{2}+g_{2}$ by
  inequalities~\eqref{ineq:succ-p1}
  and~\eqref{ineq:succ-p2}.  Hence,
  Claim~\ref{c:succ-sce-scc-sum} holds by the semi-strict
  transitivity of $\sce,\scc$.
            
  Let us prove Claim~\ref{c:succ-scc-sce-mult}.  By the
  semi-strict transitivity, $f_{1}\scc f_{2}\sce 0$ implies
  $f_{1}\scc 0$. Thus, by Claim~\ref{c:succ-succeq-stwomul2}
  of Lemma~\ref{lem:succ-succeq}, $g_{1}\scc g_{2}$ implies
  \begin{align}
    f_{1}g_{1}&\scc f_{1}g_{2}.\label{ineq:succ-p3}
  \end{align}
  Moreover, 
  $f_{1}\scc f_{2}\sce 0$ gives $f_{1}\sce f_{2}\sce 0$ by the $\scc$-$\sce$ implication.
  Thus, by $g_{2}\sce 0$ and
  Claim~\ref{c:succ-succeq-stwomul1} of
  Lemma~\ref{lem:succ-succeq}, we have
  \begin{align}
    f_{1}g_{2}&\sce f_{2}g_{2}.\label{ineq:succ-p4}
  \end{align}
  Hence, $f_{1}g_{1}\scc f_{1}g_{2} \sce f_{2}g_{2}$ by
  inequalities~\eqref{ineq:succ-p3}
  and~\eqref{ineq:succ-p4}.  Therefore,
  $f_{1}g_{1} \scc f_{2}g_{2}$ by the semi-strict
  transitivity.  Also, $f_{2},g_{2}\sce 0$ yields
  $f_{2}g_{2}\sce 0$ by the multiplicativity.

  Let us prove Claim~\ref{c:succ-sce-scc-mult}.  The
  semi-strict transitivity gives
  $g_{1}\scc 0$ by $g_{1}\scc g_{2}\sce 0$.
  Hence, $g_{1} \sce 0$ holds by the
  $\scc$-$\sce$ implication. Then, since
  $f_{1}\sce f_{2}$,
  Claim~\ref{c:succ-succeq-stwomul1} of
  Lemma~\ref{lem:succ-succeq}  implies
  \begin{align}
    f_{1}g_{1}\sce f_{2}g_{1}.\label{ineq:succ-p5}
  \end{align}
  Also,
  since $g_{1}\scc g_{2}$ and
  $f_{2}\scc 0$,
  Claim~\ref{c:succ-succeq-stwomul2} of
  Lemma~\ref{lem:succ-succeq}  gives
  \begin{align}
    f_{2}g_{1}\scc f_{2}g_{2}.\label{ineq:succ-p6}
  \end{align}
  Therefore, we obtain
  $f_{1}g_{1}\sce f_{2}g_{1}\scc f_{2}g_{2}$ from
  inequalities~\eqref{ineq:succ-p5}
  and~\eqref{ineq:succ-p6}.  This implies
  $f_{1}g_{1}\scc f_{2}g_{2}$ by the semi-strict
  transitivity. Moreover, $f_{2}\scc 0$ yields $f_{2}\sce 0$
  by the $\scc$-$\sce$ implication. Therefore,
  $f_{2}g_{2}\sce 0$ follows from the multiplicativity.
\end{proof}

We state the corollary below for our later discussion.

\begin{corollary}\label{cor:pair-orders}
  Suppose squaring orders $\sce,\scc$ on $\tp{R,\geq,>}$.
  Let $f_{1},f_{2},g_{1},g_{2}\in R$ such that
  \begin{align}
    f_{1}\scc f_{2}\sce 0& \mand
                           g_{1}\sce g_{2},
                           \label{c:pair-orders-scc-sce-sce-sce}\\
    f_{1} \scc 0& \mand   g_{1} \scc 0,
                  \label{c:pair-orders-scc-scc}\\
    g_{2}\scc 0 &\mor g_{2}=0.
                  \label{c:pair-orders-scc-eq}
  \end{align}
  Then, we have
  $f_{1}g_{1} \scc f_{2}g_{2}\sce 0$.
\end{corollary}
\begin{proof}
  Since $f_{2},g_{2}\sce 0$, the multiplicativity of $\sce$
  gives $f_{2}g_{2}\sce 0$.  Hence, let us prove
  \begin{align}
    f_{1}g_{1}\scc f_{2}g_{2}.
    \label{ineq:pair-orders-p1}
  \end{align}
  Conditions~\eqref{c:pair-orders-scc-sce-sce-sce}
  and~\eqref{c:pair-orders-scc-eq} give the following cases:
  \begin{align}
    f_{1}&\scc f_{2}\sce 0 \mand g_{1}\sce g_{2}\scc 0;
           \label{c:pair-orders-p2}\\
    f_{1}&\scc f_{2}\sce 0 \mand g_{1}\sce g_{2}=0.
           \label{c:pair-orders-p3}
  \end{align}

  For Case~\ref{c:pair-orders-p2},
  Claim~\ref{c:succ-sce-scc-mult} in
  Proposition~\ref{prop:succ} yields
  inequality~\eqref{ineq:pair-orders-p1}.  For
  Case~\ref{c:pair-orders-p3}, $g_{2}=0$ implies
  $f_{2}g_{2}=0$. Then,
  inequality~\eqref{ineq:pair-orders-p1} follows, because
  $f_{1}g_{1}\scc 0$ by
  condition~\eqref{c:pair-orders-scc-scc} and the
  multiplicativity.
\end{proof}

The following binary relations give
 explicit squaring orders.  Assume that
a semiring contains the additive unit 0, but not necessarily
the multiplicative unit 1.

\begin{definition}\label{def:ineq} 
  Consider a semiring $U\subset \Q$ and the subset
  $U_{ \geq 0}=\{u\in U\mid u\geq 0\}$.  Assume
  $U\neq \{0\}$ and $U=U_{\geq 0}$.
 Let 
 $\AfX=\{r\in \R^{L}\mid 0<r_{i}<1 \mforeach i\in\oi(L)\}$.  
 Then, we put the following six binary relations on $\Q(\fX)$.
  \begin{enumerate}
  \item
    \begin{enumerate}
      \item 
  $f\geUfX g,      \mif f,g\in \Q[\fX] \mand f-g\in U[\fX]$.
\item  $f\gUfX g, \mif f\geUfX g \mand f-g\neq 0$.
    \end{enumerate}
  \item
    \begin{enumerate}
    \item       $f\geUfXpm g, \mif f,g\in \Q[\fXpm] \mand
         f-g\in U[\fXpm]$.
       \item    $f\gUfXpm g, \mif f\geUfXpm g \mand f-g\neq 0$.
    \end{enumerate}
  \item
    \begin{enumerate}
    \item  $f\geAfX g, \mif f(r),g(r)\in \R
         \mand f(r)-g(r)\geq 0 \mforeach r\in \AfX$.
       \item
         $f\gAfX g, \mif f\geAfX g \mand
         f(r)-g(r)\neq 0 \mforeach r\in \AfX$.
       \end{enumerate}
  \end{enumerate}
  For simplicity, let $\gefX=\gefX^{\Zgez} $, $\gfX=\gfX^{\Zgez}$,
  $\gefXpm =\gefXpm^{\Zgez} $, and $\gfXpm=\gfXpm^{\Zgez} $.
\end{definition}
Also, we use the following notation for our convenience.
\begin{definition}
  \label{def:det-qu}
  Let  $F,G\in \Q(\fX)^{2}$.
  Then, we write $\tdet(F,G)=\det
  \begin{bmatrix}
    F_{1}& F_{2}\\
    G_{1}& G_{2}
  \end{bmatrix}\in \Q(\fX)$.
  Also, if $F_{1}\neq 0$, then we write
  $\Qu(F)=\frac{F_{2}}{F_{1}}\in \Q(\fX)$.
\end{definition}  
By the following lemmas, we  prove that $\StgeAfX$ is a
$\geAfX$- and strict $\gAfX$-poring.

\begin{lemma}\label{lem:fun-n}
  Let $f\in \Q[\fXo]$.  Then, $f(r)=0$ for each
  $r\in A_{\fXo}$ if and only if $f=0\in \Q[\fXo]$.
\end{lemma}
\begin{proof}
  The if part is clear. Thus, let us prove the only if part.
  For $\fXo=\{X_{1,i}\}_{i\in\oi(L_{1})}$, if $L_{1}=1$, then
  it follows from the synthetic division and the infinite
  cardinality of $A_{\fXo}$.  Hence, we use the induction on
  $L_{1}$.  Let
  $\fXt=\{X_{2,i}=X_{1,i+1}\}_{i\in\oi(L_{2})}$ with
  $L_{2}=L_{1}-1\geq 1$. Also, assume $f\neq 0$ in
  $\Q[\fXo]$.  Then, for
  $f=\sum_{j\in \Z^{L_{1}}} f_{j_{1},\dots, j_{L_{1}}}
  X_{1,1}^{j_{1}}\dots X_{1,L_{1}}^{j_{L_{1}}}$, there is
  $\mu\in \Z^{L_{1}}$ such that
  $f_{\muo,\dots, \mu_{L_{1}}}\neq 0$.  Then, in $\Q[\fXt]$,
  $h= \sum_{j\in \Z^{L_{2}}}f_{\muo,j_{1},\dots, j_{L_{2}}}
  X_{2,1}^{j_{1}}\dots X_{2,L_{2}}^{j_{L_{2}}}\neq 0$.  The
  induction gives $v\in A_{\fXt}$ such that $h(v)\neq 0$.
  Also, since
  $g=f(\tp{X_{1},v_{1}\dots, v_{L_{2}}})\in \Q[X_{1}]$
  satisfies $g\neq 0$, the induction gives a real number
  $0<t<1$ such that $g(\tp{t})\neq 0$.  Thus,
  $u=\tp{t,v_{1},\dots, v_{L_{2}}}\in \AfXo$ gives
  $f(u)\neq 0$.
\end{proof}

\begin{lemma}\label{lem:succ-site}
  If $R$ is a strict $\scc$-poring such that
  $\emptyset\neq \StccR\subset R$, then $\StccR= R$.
\end{lemma}
\begin{proof}
  Because some $f\scc g$ implies $f-g\scc 0$ by
  Claim~\ref{c:succ-succeq-sinv2} of
  Lemma~\ref{lem:succ-succeq}, each $h\in R$ satisfies
  $f-g+h\scc h$ by the additivity of $\scc$.
\end{proof}

\begin{proposition}\label{prop:poring-geAfX}
  We have the following porings.
  \begin{enumerate}
  \item \label{c:poring-geAfX-StgeAfX-geAfX} $\StgeAfX$ is the
    $\geAfX$-poring such that
    \begin{align}
      \StgeAfX= \{f\in \Q(\fX) \mid f(r)\in \R \mforeach r\in
      \AfX\}.
       \label{eq:poring-geAfX-StgeAfX-geAfX}
    \end{align}
  \item  $\StgeAfX$ is a strict $\gAfX$-poring.
    \label{c:poring-geAfX-StgeAfX-gAfX}
  \item  $\StgeAfX=\StgAfX$.
    \label{c:poring-geAfX-StgeAfX-StgAfX}
  \end{enumerate}
\end{proposition}
\begin{proof}
  Let us prove Claim~\ref{c:poring-geAfX-StgeAfX-geAfX}.
  First, if $f\in \StgeAfX$, then the reflexivity
  $f \geAfX f$ follows from $f(r)-f(r)=0\geq 0$ for each
  $r\in \AfX$.  In particular,
  equation~\eqref{eq:poring-geAfX-StgeAfX-geAfX} holds.

  Second, if $f_{1}\geAfX f_{2}\geAfX f_{3}$, then
  $f_{1}(r)-f_{2}(r)\geq 0$ and $f_{2}(r)-f_{3}(r)\geq 0$ for each
  $r\in \AfX$.  Thus, since $f_{1}(r)-f_{3}(r)\geq 0$ for each
  $r\in \AfX$, the transitivity $f_{1}\geAfX f_{3}$ holds.

  Third, let us prove the antisymmetricity.  If
  $f\geAfX g\geAfX f$, then $f(r)-g(r)\geq 0$ and
  $g(r)-f(r)\geq 0$ for each $r\in \AfX$.  Thus,
  $f(r)=g(r)\in \R$ for each $r\in \AfX$.  Then, there are
  polynomials $F,G\in \Q[\fX]^{2}$ such that
  $f=\Qu(F)$ and $g=\Qu(G)$,
  $F_{1}(r)\cdot G_{1}(r)\neq 0$, and
  $\tdet(F,G)(r)=0$ for each $r\in \AfX$.
  Then, since $\tdet(F,G)=0\in \Q[\fX]$ by
  Lemma~\ref{lem:fun-n}, the reflexivity $f=g$ holds.
  
  Fourth, let us confirm the additivity and multiplicativity
  of $\geAfX$.  Let $f_{1}\geAfX f_{2}$ and
  $f_{3}\in \StgeAfX$. Then,
  $f_{1}(r)+f_{3}(r)-(f_{2}(r)+f_{3}(r))=f_{1}(r)-f_{2}(r)\geq
  0$ for each $r\in \AfX$.  Thus, the additivity
  $f_{1}+f_{3}\geAfX f_{2}+f_{3}$ holds.  Also, if
  $f_{1},f_{2}\geAfX 0$, then the multiplicativity
  $f_{1}f_{2}\geAfX 0$ follows from $f_{1}(r)f_{2}(r)\geq 0$
  for each $r\in \AfX$.  Hence,
  Claim~\ref{c:poring-geAfX-StgeAfX-geAfX} holds.

  Claim~\ref{c:poring-geAfX-StgeAfX-gAfX} follows as above,
  since $f\gAfX f$ does not hold on $\StgeAfX$.  Also,
   Lemma~\ref{lem:succ-site} yields
 Claim~\ref{c:poring-geAfX-StgeAfX-StgAfX} by $1\gAfX 0$.
\end{proof}

We now introduce the following squaring implications and orders on $\fX$.

\begin{definition}\label{def:succ-on-x}
  We refer to the squaring implications of binary relations
  $\sce,\scc$ and $\geAfX, \gAfX$ on $\Q(\fX)$ as the squaring
  implications of $\sce,\scc$ on $\fX$.  Moreover, we refer
  to squaring orders $\sce,\scc$ on
  $\tp{\Q(\fX),\geAfX, \gAfX}$ as squaring orders on $\fX$.
\end{definition}
 
By the following lemma, we obtain squaring orders on
$\fX$.
\begin{lemma}\label{lem:poring-geUfX-geUfXpm}
  We have the following.
  \begin{enumerate}
  \item \label{c:poring-geUfX-geUfXpm-StgeUfX-geUfX}
    $\StgeUfX$ is the $\geUfX$-poring such that
    $\StgeUfX=\Q[\fX]$.
  \item  $\StgeUfX$ is a strict $\gUfX$-poring.
    \label{c:poring-geUfX-geUfXpm-StgeUfX-gUfX}
  \item  $\StgeUfX=\StgUfX$.
    \label{c:poring-geUfX-geUfXpm-StgeUfX-StgUfX}
  \item \label{c:poring-geUfX-geUfXpm-StgeUfXpm-geUfXpm}
    $\StgeUfXpm$ is the
    $\geUfXpm$-poring such that  $\StgeUfXpm=\Q[\fXpm]$.
  \item  $\StgeUfXpm$ is a strict $\gUfXpm$-poring.
    \label{c:poring-geUfX-geUfXpm-StgeUfXpm-gUfXpm}
  \item  $\StgeUfXpm=\StgUfXpm$. 
    \label{c:poring-geUfX-geUfXpm-StgeUfXpm-StgUfXpm}
  \end{enumerate}
\end{lemma}
\begin{proof}
  Let us prove
  Claim~\ref{c:poring-geUfX-geUfXpm-StgeUfX-geUfX}.  By
  $0\in U[\fX]$, we have the reflexivity $f\geUfX f$ for
  each $f\in \StgeUfX$.  Also, $\StgeUfX=\Q[\fX]$ by the
  reflexivity.  Let $f_{1}\geUfX f_{2} \geUfX f_{3}$.  Since
  $f_{1}-f_{2},f_{2}-f_{3}\in U[\fX]$ and $U$ is a semiring,
  $f_{1}-f_{3}=(f_{1}-f_{2})+(f_{2}-f_{3})\in U[\fX]$.
  Thus, the transitivity $ f_{1}\geUfX f_{3}$ follows.  If
  $f_{1}\geUfX f_{2}\geUfX f_{1}$, then
  $f_{1}-f_{2},f_{2}-f_{1}\in U[\fX]$.  Hence,
  $U=U_{\geq 0}$ yields the antisymmetricity $f_{1}=f_{2}$.
  Moreover, if $f_{1}\geUfX f_{2}$ and $f_{3}\in \StgeUfX$,
  then the additivity $f_{1}+f_{3}\geUfX f_{2}+f_{3}$ holds
  by $(f_{1}+f_{3})-(f_{2}+f_{3})=f_{1}-f_{3}\in U[\fX]$.
  Also, $f_{1},f_{2}\geUfX 0$ implies the multiplicativity
  $f_{1}f_{2}\geUfX0$, since $U$ is a semiring.  Then,
  Claim~\ref{c:poring-geUfX-geUfXpm-StgeUfX-geUfX} holds.

   Claim~\ref{c:poring-geUfX-geUfXpm-StgeUfX-gUfX} follows as above, since
  $f\gUfX g$ demands $f\neq g$. Also,
  Lemma~\ref{lem:succ-site} gives
  Claim~\ref{c:poring-geUfX-geUfXpm-StgeUfX-StgUfX}, because there is
  $f\in U$ such that $f\gUfX 0$ by $U\supsetneq \{0\}$.
  
  Similar arguments hold for
   Claims~\ref{c:poring-geUfX-geUfXpm-StgeUfXpm-geUfXpm},
  ~\ref{c:poring-geUfX-geUfXpm-StgeUfXpm-gUfXpm},
  and~\ref{c:poring-geUfX-geUfXpm-StgeUfXpm-StgUfXpm}.
  \end{proof}

  \begin{proposition}\label{prop:succ-exist}
  The binary relations $\geUfX$, $\gUfX$, $\geUfXpm$, $\gUfXpm$,
  $\geAfX$, $\gAfX $ are squaring orders on $\fX$ such that
  $\gUfX$, $\gUfXpm$, $\gAfX$ are strict squaring orders of
  $\geUfX$, $\geUfXpm$, $\geAfX$, respectively.
    \end{proposition}
    \begin{proof}
      First, let us prove that $\geAfX, \gAfX$ are squaring
      orders on $\fX$.  By Proposition~\ref{prop:poring-geAfX},
      $\geAfX, \gAfX$ on $\fX$ satisfy the
      $\geAfX$-$\gAfX$-poring equality.  Also,
      $\geAfX,\gAfX$ have the squaring implications on $\fX$ by
      Definition~\ref{def:ineq}.  Let us confirm the
      semi-strict transitivity of $\geAfX,\gAfX$.  Assume
      $f_{1}\geAfX f_{2}\gAfX f_{3}$.  Then,
      $f_{1}(r)-f_{2}(r)\geq 0$ and $ f_{2}(r)-f_{3}(r) >0$ for
      each $r\in \AfX$.  Thus,
      $f_{1}(r)-f_{3}(r)=
      (f_{1}(r)-f_{2}(r))+(f_{2}(r)-f_{3}(r))>0$.  Then,
      $f_{1}\gAfX f_{3}$.  Similarly,
      $f_{1}\gAfX f_{2}\geAfX f_{3}$ implies
      $f_{1}\gAfX f_{3}$.  Hence, $\geAfX, \gAfX$ are
      squaring orders on $\fX$ such that $\gAfX$ is a strict
      squaring order of $\geAfX$.

      Second, let us prove that $\geUfX, \gUfX$ are squaring
      orders on $\fX$.  By
      Claims~\ref{c:poring-geUfX-geUfXpm-StgeUfX-geUfX},
      ~\ref{c:poring-geUfX-geUfXpm-StgeUfX-gUfX},
      and~\ref{c:poring-geUfX-geUfXpm-StgeUfX-StgUfX} of
      Lemma~\ref{lem:poring-geUfX-geUfXpm}, $\geUfX, \gUfX$
      give the $\geAfX$-$\gAfX$-poring equality. 
      The squaring implications of $\geUfX$, $\gUfX$ on
      $\fX$ follows from Definition~\ref{def:ineq}.
      Also, we have the
      $\geUfX$-$\geAfX$-poring inclusion by
      Claim~\ref{c:poring-geUfX-geUfXpm-StgeUfX-geUfX}, of
      Lemma~\ref{lem:poring-geUfX-geUfXpm} and
       Claim~\ref{c:poring-geAfX-StgeAfX-geAfX} of
       Proposition~\ref{prop:poring-geAfX}.  Let us check the
      semi-strict transitivity of $\geUfX, \gUfX$.  Let
      $f_{1}\gUfX f_{2}\geUfX f_{3}$.  Then,
      $f_{1}-f_{2}\in U[\fX]$ with $f_{1}- f_{2}\neq 0$ and
      $f_{2}-f_{3}\in U[\fX]$.  Hence, $f_{1}\gUfX f_{3}$,
      since $f_{1}-f_{3}\in U[\fX]$ and
      $f_{1}-f_{3}=(f_{1}-f_{2})+(f_{2}-f_{3}) \neq 0$ by
      $U=U_{\geq 0}$.  Also,
      $f_{1}\geq _{\fX}^{U}f_{2}\gUfX f_{3}$ implies
      $f_{1}\gUfX f_{3}$.  Therefore, $\geUfX,\gUfX$ are
      squaring orders such that $\gUfX$ is a strict squaring
      order of $\geUfX$.

       Similar arguments hold for
       $\geUfXpm,\gUfXpm$.
    \end{proof}

     Squaring orders on $\fX$ are not necessarily
     ones in Proposition~\ref{prop:succ-exist}.
     For example, if
     $\fX=\{X_{1}\}$, then
     convergent Maclaurin series on
     $\abs{X_{1}}<1$ with positive
     coefficients define a strict squaring order on $\fX$.

     We compare squaring orders by the following
     terminology.

\begin{definition}\label{def:succ-compatile}
  Consider fields $\Q(\fXo)\subset \Q(\fXt)$.  Assume
  squaring orders $O_{i}=\{\sce_{i},\scc_{i}\}$ on $\fX_{i}$
  for $i\in\oi(2)$.  Then, $O_{2}$ is said to be
  compatible to $O_{1}$, if the $\sceo$-$\scet$ and
  $\scco$-$\scct$ implications hold.  We write
  $O_{2}\Sup O_{1}$, if $O_{2}$ is compatible to $O_{1}$.
\end{definition}
      
For instance, $\{\geAfX, \gAfX\}\Sup O$ for any squaring
orders $O=\{\sce,\scc\}$ on $\fX$ by
Lemma~\ref{lem:half-full-implication}.
    
\section{Admissible variables}\label{sec:adm}
 
By squaring orders, we introduce the notion of admissible
variables to discuss polynomials and their values in real
numbers.
We call $x\in \Q(\fX)^{l}$ indeterminate,
if each $x_{i}$ is an indeterminate.  First, we extend
Definition~\ref{def:q-pic-poly-q-Laurent-pic-poly} for an
indeterminate $x\in \Q(\fX)^{l}$, whose elements are not
necessarily algebraically independent over $\Q$.
       
\begin{definition}\label{def:ineq-adm} 
  For $l\in \Zgeo$, suppose an indeterminate
  $x\in \Q(\fX)^{l}$.  Let $U\subset \Q$ be a semiring such
  that $U=U_{\geq 0}$ and $U\neq \{0\}$. 
  \begin{enumerate}
  \item Let us define the following subsets in $\Q(\fX)$.
  \begin{enumerate}
  \item The subset $U[x]$ consists of
    $f=\sum_{j\in \Zlgez}
    f_{j_{1},\dots,j_{l}}x_{1}^{j_{1}} \dots x_{l}^{j_{l}}$
    for some finitely many non-zero
    $f_{j_{1},\dots,j_{l}}\in U$.
  \item The subset $U[\xpm]$ consists of
    $f =\sum_{j\in \Zl} f_{j_{1},\dots,j_{l}}x_{1}^{j_{1}}
    \dots x_{l}^{j_{l}} $ for some finitely many non-zero
    $f_{j_{1},\dots,j_{l}}\in U$.
  \end{enumerate}

\item Let us define the following binary relations on
  $\Q(\fX)$. 
  \begin{enumerate}
  \item
    \begin{enumerate}
    \item  $f\geUx g,  \mif f,g\in \Q[x] \mand f-g\in U[x]$.
    \item $f\gUx g, \mif f\geUx g \mand f-g\neq 0$.  
    \end{enumerate}
  \item 
    \begin{enumerate}
      \item $f\geUxpm g, \mif f,g\in \Q[\xpm] \mand f-g\in U[\xpm]$.
      \item $f\gUxpm g, \mif f\geUxpm g \mand f-g\neq 0$.  
    \end{enumerate}
  \end{enumerate}
\end{enumerate}

  For simplicity, if $x=\tp{q}$, then let $\gU_{q}=\gUx$,
  $\geU_{q}= \geUx$, $\gU_{\qpm}=\gUxpm$, and
  $\geU_{\qpm}=\geUxpm$
  Also, if $U=\Zgez$, then we often omit the
  superscript ${}^{U}$ for the binary relations above.
\end{definition}

For an indeterminate $x\in \Q(\fX)$, $\gx$ on $\Q(\fX)$ does
not have to be a squaring order on $\fX$, because $x\gAfX 0$
is not necessarily true.  Thus, we define the following.

\begin{definition}\label{def:admissible}
  Suppose squaring orders $O=\{\sce,\scc\}$ on $\fX$ and an
  indeterminate $x\in \Q(\fX)$.  We call $x$ $O$-admissible
  (or admissible for short), if $x$ and $O$ satisfy the
  following conditions:
  \begin{enumerate}
  \item $f\gx 0$ implies $f\scc 0$ (half $\gx$-$\scc$
    implication);
  \item $1\gAfX x$ (upper condition of $x$ on $\fX$).
  \end{enumerate}
  For $l\in \Zgeo$ and $y\in \Q(\fX)^{l}$, if
  each $y_{i}$ is $O$-admissible, then we call $y$
  $O$-admissible (or admissible for short).
\end{definition}

For example,
for each $X_{i}\in \fX$ and  $d\in \Zgeo$,
$X_{i}^{d}$
is $\{\gefX,\gfX\}$- and
$\{\gefXpm,\gfXpm\}$-admissible. Also, let us state the following
lemmas to obtain $\geUx$- and $\geUxpm$-porings

\begin{lemma}\label{lem:pos-nonneg-implication}
  Suppose squaring orders $O=\{\sce,\scc\}$ on $\fX$ and an
  indeterminate $x\in \Q(\fX)$. Then, the half $\gx$-$\scc$
  implication yields the half $\gex$-$\sce$ implication.
\end{lemma}
\begin{proof}
  Since $f\gex 0$ says $f\gx 0$ or $f=0$,
  the assertion holds by the reflexivity of $\sce$.
\end{proof}

\begin{lemma}\label{lem:adm-succ}
  Assume squaring orders $O=\{\sce,\scc\}$ on $\fX$.  For
  $l\in \Zgeo$, consider an $O$-admissible $x\in \Q(\fX)^{l}$.
  Then, we have the following.
  \begin{enumerate}
  \item \label{c:adm-succ-half-gex-sce} The half $\gex$-$\sce$
    implication holds.
  \item\label{c:adm-succ-half-gx-scc} The half $\gx$-$\scc$
    implication holds.
  \item\label{c:adm-succ-sx1} Let $f\geUx 0$. Then,
    there is
    $\lam\in \Zgeo$ such that $\lam f\sce 0$.  Also,
    $\lam\mu f\sce 0$ for each $\mu\gex 0$.
  \item\label{c:adm-succ-sx2} Let $f\gUx 0$. Then,
    there is
    $\lam\in \Zgeo$ such that $\lam f\scc 0$.  Also,
    $\lam\mu f\scc 0$ for each $\mu\gx 0$.
  \item\label{c:adm-succ-sxpm1} Let $f\geUxpm 0$.
    Then,
    there is a monomial $\lam$ of variables $x_{i}$ such
    that $\lam\gx 0$ and $\lam f\sce 0$.  Also,
    $\lam \mu f \sce 0$ for each $\mu\gex 0$.
  \item\label{c:adm-succ-sxpm2} Let $f\gUxpm 0$. Then, there
    is a 
    monomial $\lam$ of variables $x_{i}$ such that
    $\lam\gx 0$ and $\lam f\scc 0$.  Also, $\lam \mu f \scc 0$
    for each $\mu\gx 0$.
  \end{enumerate}
\end{lemma}
\begin{proof}
  Claim~\ref{c:adm-succ-half-gex-sce} follows from the
  additivity and multiplicativity of $\sce$, since
  $f\geq_{x_{i}}0$ implies $f\sce 0$ by
  Lemma~\ref{lem:pos-nonneg-implication}.
  Claim~\ref{c:adm-succ-half-gx-scc} holds similarly by the
  half $>_{x_{i}}$-$\scc$ implication .

  Let us prove Claim~\ref{c:adm-succ-sx1}.  Because
  $U=U_{\geq 0}\subset \Q$ and $f\in U[x]$, there is
  $\lam\in \Zgeo$ such that $\lam f \gex 0$, Thus,
  $\lam f\sce 0$ by Claim~\ref{c:adm-succ-half-gex-sce}.
  Since $\mu\sce 0$ by Claim~\ref{c:adm-succ-half-gex-sce},
  the latter assertion follows from the multiplicativity of
  $\sce$.

  Let us prove Claim~\ref{c:adm-succ-sx2}.  Because $f\in U[x]$
  and $f\neq 0$, some $\lam\in \Zgeo$ yields $\lam f \gx 0$.
  This gives
   $\lam f\scc 0$ by Claim~\ref{c:adm-succ-half-gx-scc}.  Also,
  $\lam \mu f \scc 0$ by the multiplicativity of $\scc$ and
  Claim~\ref{c:adm-succ-half-gx-scc}.

  Let us prove Claim~\ref{c:adm-succ-sxpm1}.  We have some
  finitely many $f_{j_{1},\dots,j_{l}}\in \Q_{> 0}$ such
  that
  $f=\sum_{j\in \Zl} f_{j_{1},\dots,j_{l}}x_{1}^{j_{1}}
  \dots x_{l}^{j_{l}}$.  Hence, there is $u\in \Zlgeo$ such
  that $u_{i}+j_{i}\geq 0$, whenever $i\in\oi(l)$ and
  $f_{j_{1},\dots,j_{l}}\in \Q_{> 0}$.  This gives
  $v\in \Zgeo$ such that the monomial
  $\lam=v x_{1}^{u_{1}}\dots x_{l}^{u_{l}}\gx 0$ satisfies
  $\lam f \gex 0$.  Thus,
  Claim~\ref{c:adm-succ-half-gex-sce} yields
  $\lam f \sce 0$.  The latter assertion follows from the
  multiplicativity of $\sce$, since $\mu\sce 0$ by
  Claim~\ref{c:adm-succ-half-gex-sce}.
  
  Let us prove Claim~\ref{c:adm-succ-sxpm2}.  As above for
  $\geUxpm$, some $u\in \Zlgeo$ and
  $v\in \Zgeo$ give the monomial
  $\lam=v x_{1}^{u_{1}}\dots x_{l}^{u_{l}}\gx 0$ such that
  $\lam f \gex 0$. Also, $\lam f\neq 0$, since $\Q[\fX]$ is
  an integral domain and $\lam\neq 0$.  Thus, since
  $\lam f \gx 0$, Claim~\ref{c:adm-succ-half-gx-scc} yields
  $\lam f \scc 0$.  The latter assertion holds by the
  multiplicativity of $\scc$ and $\mu\scc 0$ of
  Claim~\ref{c:adm-succ-half-gx-scc}.
\end{proof}

\begin{proposition}\label{prop:adm-poring}
  Suppose squaring orders $O=\{\sce,\scc\}$ on $\fX$.  Let
  $l\in \Zgeo$ and $x\in \Q(\fX)^{l}$ be $O$-admissible.
  Then, we have the following porings.
  \begin{enumerate}
  \item        \label{c:adm-poring-sx1}
    $\StgeUx$ is the $\geUx$-poring
    such that $\StgeUx=\Q[x]$.
  \item  $\StgeUx$ is a strict $\gUx$-poring.
    \label{c:adm-poring-sx2}
  \item  $\StgeUx=\StgUx$.
    \label{c:adm-poring-sx3}
  \item   \label{c:adm-poring-sxpm1}
    $\StgeUxpm$ is the $\geUxpm$-poring
    such that  $\StgeUxpm=\Q[\xpm]$.
  \item   $\StgeUxpm$ is a strict $\gUxpm$-poring.
    \label{c:adm-poring-sxpm2}
  \item  $\StgeUxpm=\StgUxpm$.
    \label{c:adm-poring-sxpm3}
  \end{enumerate}
\end{proposition}
\begin{proof}
  Let us prove Claim~\ref{c:adm-poring-sx1}.  First, each
  $f\in \StgeUx$ satisfies the reflexivity $f\geUx f$, since
  $f-f=0\in U[x]$.  In particular, $\StgeUx=\Q[x]$ by the
  reflexivity.

  Second,  $f\geUx g \geUx h$ implies the
  transitivity $f\geUx h$, since $f-g,g-h\in U[x]$
  gives $(f-g)+(g-h)\in U[x]$ for the semiring $U$.
  
  Third, let us prove the antisymmetricity of $\geUx$.  Let
  $f\geUx g\geUx f$. Then, Claim~\ref{c:adm-succ-sx1} of
  Lemma~\ref{lem:adm-succ} gives $\lamo,\lamt\in \Zgeo$
  such
  that $\lamo(f-g) \sce 0$ and $\lamt(g-f)\sce 0$.  This
  gives $\lamo\lamt(f-g) \sce 0$ and
  $\lamo\lamt(g-f)=-\lamo\lamt (f -g) \sce 0$ by
  $\lamo,\lamt\gex 0$ and Claim~\ref{c:adm-succ-sx1} of
  Lemma~\ref{lem:adm-succ}.  Also, the additivity of $\sce$
  implies $0 \sce \lamo\lamt (f -g)$.  Hence,
  $\lam (f- g )=0$ for $\lam=\lamo\lamt$ by the
  antisymmetricity of $\sce$.  Then, $\lam f=\lam g$, which
  yields $f=g\in \Q[\fX]$.

  Fourth, let us prove the additivity and multiplicativity of
  $\geUx$.  Let $f\geUx g$ and $h\in \StgeUx$.  Then, the
  additivity $f +h\geUx g+h$ follows, since
  $(f+h)-(g+h)=f-g\in U[x]$.  Next, let $f, g\geUx 0$. Then,
  $fg\in U[x]$, since $f,g\in U[x]$ for the semiring $U$.
  Hence, we obtain the multiplicativity $fg \geUx 0$.  Thus,
  Claim~\ref{c:adm-poring-sx1} follows.

  Let us confirm Claim~\ref{c:adm-poring-sx2}.  Thus, let us
  prove the irreflexivity and transitivity of $\gUx$.  For
  each $f\in \StgeUx$, $f \gUx f$ does not hold by $f-f=0$.
  Then, let us prove the transitivity.  Consider
  $f\gUx g \gUx h$.  This implies $f\geUx h$ by the
  transitivity of $\geUx$ and $\gUx$-$\geUx$ implication.
  Moreover, Claim~\ref{c:adm-succ-sx2} of
  Lemma~\ref{lem:adm-succ} gives $\lam\in \Zgeo$ such that
  $\lam(f-g), \lam(g-h)\scc 0$.  Hence,
  $\lam(f-g)+ \lam(g-h)=\lam(f-h) \scc 0$ by the additivity
  of $\scc$.  Therefore, we obtain the transitivity
  $f \gUx h$, since $f-h\neq 0$ by the irreflexivity of
  $\scc$.

  Let us confirm the additivity and multiplicativity of
  $\gUx$ on $\StgeUx$.  The additivity holds as above for
  $\geUx$. Thus, let us prove the
  multiplicativity.  Let $f, g\gUx 0$.  Then, $fg\in U[x]$,
  since $f,g\in U[x]$ and $f,g\neq 0$.  Also, $fg\neq 0$,
  since $\Q[\fX]$ is an integral domain. Hence, the
  multiplicativity $fg \gUx 0$ follows.  Then,
  Claim~\ref{c:adm-poring-sx2} holds.

  Lemma~\ref{lem:succ-site} yields
  Claim~\ref{c:adm-poring-sx3}, since $U\neq \{0\}$ gives some
  $f,g\in U[x]$ such that $f\gUx g$.

  Let us prove Claim~\ref{c:adm-poring-sxpm1}.  The
  reflexivity, transitivity, additivity, and
  multiplicativity of $\geUxpm$ on $\StgeUxpm$ hold as
  above for $\geUx$.
  In particular, $\StgeUxpm=\Q[\xpm]$ by the
  reflexivity.
  Hence, let us confirm the
  antisymmetricity.  Suppose $f\geUxpm g\geUxpm f$.
  By  Claim~\ref{c:adm-succ-sxpm1} of
  Lemma~\ref{lem:adm-succ}, there are
   monomials $\lam,\lam'$ of
  variables $x_{i}$ such that $\lam,\lam'\gx 0$,
  $\lam(f-g)\sce 0$, and $\lam' (g-f)\sce 0$.  Moreover,
  since $\lam,\lam'\gex 0$, we have
  $\lam\lam'(f-g)\sce 0$, and
  $\lam\lam' (g-f)=-\lam\lam'(f-g)\sce 0$ by
  Claim~\ref{c:adm-succ-sxpm1} of
  Lemma~\ref{lem:adm-succ}.
  Hence, $\lam\lam' (f-g) \sce 0$ and
  $0\sce \lam\lam' (f -g)$ by the additivity of $\sce$.
  Therefore, $ \lam\lam' (f-g)=0\in \Q[\fX]$ by the
  antisymmetricity of $\sce$.  In particular, the
  antisymmetricity $f=g$ of $\geUx$ holds, since $\Q[\fX]$
  is an integral domain and $\lam\lam'\neq 0$.  Thus,
  Claim~\ref{c:adm-poring-sxpm1} holds.

  Let us prove Claim~\ref{c:adm-poring-sxpm2}.  The
  irreflexivity holds by $f-f=0$.  As above for $\gUx$, we
  have the additivity and multiplicativity of $\gUxpm$, and
  the transitivity by Claim~\ref{c:adm-succ-sxpm2} of
  Lemma~\ref{lem:adm-succ}.

  Lemma~\ref{lem:succ-site} implies
  Claim~\ref{c:adm-poring-sxpm3}, since $f\gUxpm g$ for some
  $f,g\in U[\xpm]$.
\end{proof}

Moreover, we prove the following to obtain squaring orders
by admissible variables.
\begin{lemma}\label{lem:adm-bounds-poring-inclusions}
  Suppose squaring orders $O=\{\sce,\scc\}$ on $\fX$.
  For $l\in \Zgeo$, let
  $x\in \Q(\fX)^{l}$ be $O$-admissible. Then, we have the
  following.
  \begin{enumerate}
  \item \label{c:adm-bounds} For each $r\in \AfX$ and
    $i\in\oi(l)$, we have
    $x_{i}(r)\in\R$ such that $0<x_{i}(r)<1$.
  \item \label{c:adm-poring-inclusions} We have
    $\StgeUx\subset \StgeAfX$ and $\StgeUxpm\subset \StgeAfX$.
  \end{enumerate}
\end{lemma}
\begin{proof}
  Let us prove Claim~\ref{c:adm-bounds}.  Since $x_{i}>_{x_{i}}0$,
  we have $x_{i}\scc 0$ by the half $>_{x_{i}}$-$\scc$ implication.
  Then, $x_{i}\gAfX 0$ by the half $\scc$-$\gAfX$ implication.
  Hence, $x_{i}(r)\in \R$ and $x_{i}(r)>0$ for each $r\in \AfX$.
  Thus, Claim~\ref{c:adm-bounds} holds, since the upper
  condition of $x_{i}$ implies $1>x_{i}(r)$.

  Claim~\ref{c:adm-poring-inclusions} follows from
  Claim~\ref{c:adm-bounds} and Definitions~\ref{def:ineq}
  and~\ref{def:ineq-adm}.
\end{proof}

By Claim~\ref{c:adm-bounds} of
Lemma~\ref{lem:adm-bounds-poring-inclusions}, for each
$\ka\in \Q$ and $r\in \AfX$, we assume
\begin{align}
  x_{i}(r)^{\ka}=\sqrt[\ka]{x_{i}(r)}\in \R_{>0}.
\end{align}
 
Moreover, we discuss real values of rational functions and
polynomials by the following admissible variables, which are
not necessarily algebraically independent.

\begin{theorem}\label{thm:adm}
  Consider squaring orders $O=\{\sce,\scc\}$ on $\fX$.  For
  $l\in \Zgeo$, let $x\in \Q(\fX)^{l}$ be $O$-admissible.
  Then, $\geUx, \gUx,\geUxpm, \gUxpm$ are squaring orders
  on
  $\fX$ such that $\gUx$ and $ \gUxpm$ are strict squaring
  orders of $\geUx$ and $\geUxpm$, respectively.
\end{theorem}
\begin{proof}
  Let us prove the assertion for $\geUx,\gUx$.  We have the
  $\gAfX$-$\geAfX$ and $\gUx$-$\geUx$ implications by
  Definitions~\ref{def:ineq} and~\ref{def:ineq-adm}.  Also,
  Claims~\ref{c:adm-poring-sx1}, \ref{c:adm-poring-sx2},
  and~\ref{c:adm-poring-sx3} of Proposition~\ref{prop:adm-poring}
  give the $\geUx$-$\gUx$-poring equality. The
  $\geUx$-$\geAfX$-poring inclusion holds by
  Claim~\ref{c:adm-poring-inclusions} of
  Lemma~\ref{lem:adm-bounds-poring-inclusions}.  Hence, let
  us prove the half $\geUx$-$\geAfX$ and $\gUx$-$\gAfX$
  implications and the semi-strict transitivity of
  $\geUx,\gUx$.
  
  First, let us confirm the half $\geUx$-$\geAfX$
  implication.  Let $f\geUx 0$.  Then, $\lam f \sce 0$ for
  some $\lam\in \Zgeo$ by Claim~\ref{c:adm-succ-sx1} of
  Lemma~\ref{lem:adm-succ}.  This implies $\lam f \geAfX 0$
  by the half $\sce$-$\geAfX$ implication.  Therefore, since
  $f \geAfX 0$ by $\lam\in\Zgeo$, the half $\geUx$-$\geAfX$
  implication follows.
  
  Second, let us prove the half $\gUx$-$\gAfX$ implication.
  Let $f\gUx 0$.  This gives $\lam f \succ 0$ for some
  $\lam\in \Zgeo$ by Claim~\ref{c:adm-succ-sx2} of
  Lemma~\ref{lem:adm-succ}. Then, $\lam f \gAfX 0$ by the
  half $\scc$-$\gAfX$ implication.  Thus, since $f \gAfX 0$
  by $\lam\in \Zgeo$, the half $\gUx$-$\gAfX$ implication
  holds.

  Third, let us prove the semi-strict transitivity.  Let
  $f \geUx g \gUx h$.  Then, $f-g,g-h\in U[x]$ yields
  $f-h\in U[x]$, which implies $f\geUx h$. Also, the
  $\geUx$-$\geAfX$ and $\gUx$-$\gAfX$ implications give
  $f \geAfX g \gAfX h$.  Then, since $f \neq h$ by
  $f\gAfX h$, the semi-strict transitivity $f \gUx h$ holds.
  Similarly, $f \gUx g \geUx h$ implies $f\gUx h$.  Thus, we
  obtain the assertion for $\geUx,\gUx$.
  
  The assertion for $\geUxpm,\gUxpm$ follows from a
  parallel argument. In particular, we replace
  Claims~\ref{c:adm-succ-sx1} and~\ref{c:adm-succ-sx2} of
  Lemma~\ref{lem:adm-succ} with
  Claims~\ref{c:adm-succ-sxpm1} and~\ref{c:adm-succ-sxpm2}
  of Lemma~\ref{lem:adm-succ} to prove the half
  $\geUxpm$-$\geAfX$ and $\gUxpm$-$\gAfX$ implications.
\end{proof}

We also obtain the following admissible variables
from given ones.
\begin{proposition}\label{prop:adm-film}
  Suppose squaring orders $O=\{\sce,\scc\}$ on
  $\fX$ and $l\in \Zgeo$. Let
  $x\in \Q(\fX)^{l}$ be $O$-admissible.
  \begin{enumerate} 
  \item\label{c:adm-film-O-admissible-gen} Suppose
    $f=\sum_{j\in \Zlgez} f_{j_{1},\dots, j_{l}}
    x_{1}^{j_{1}}\dots x_{l}^{j_{l}}\in \Q_{\geq 0}[x]$
    with the following conditions:
    \begin{enumerate}
    \item $f\not\in \Q$;
      \label{c:adm-film-O-admissible-gen-non-const}
    \item $f_{j_{1},\dots, j_{l}}\neq 0$ implies
      $f_{j_{1},\dots, j_{l}}\scc 0$;
      \label{c:adm-film-O-admissible-scc}
    \item 
      \label{c:adm-film-O-admissible-bound}
      $0<\sum_{j\in \Zlgez} f_{j_{1},\dots, j_{l}}\leq 1$.
    \end{enumerate}
    Then, $f$ is $O$-admissible.
  \item\label{c:adm-film-O-admissible-monom} If
    $f=\prod_{i\in\oi(l)} x_{i}^{j_{i}}\neq 1$ for some
    $j\in \Zlgez$, then $f$ is $O$-admissible.
  \end{enumerate}
\end{proposition}
\begin{proof}
  Let us prove Claim~\ref{c:adm-film-O-admissible-gen}.
  First, let us confirm
  \begin{align}\label{ineq:adm-film-p1}
    f\scc 0.
  \end{align}
  Since Claim~\ref{c:adm-succ-half-gx-scc} of
  Lemma~\ref{lem:adm-succ} gives
  $x_{1}^{j_{1}}\dots x_{l}^{j_{l}}\scc 0$ for each
  $j\in \Zlgez$, $f_{j_{1},\dots, j_{l}}\neq 0$ implies
  $f_{j_{1},\dots, j_{l}}x_{1}^{j_{1}}\dots
  x_{l}^{j_{l}}\scc 0$ by
  Condition~\ref{c:adm-film-O-admissible-scc} and the
  multiplicativity of $\scc$.  Thus, the additivity of
  $\scc$ yields inequality~\eqref{ineq:adm-film-p1}.
  
  Second, $g>_{f} 0$ implies $g\scc 0$, since we have
  inequality~\eqref{ineq:adm-film-p1}, the additivity and
  multiplicativity of $\scc$, and $1\scc 0$ by
  Claim~\ref{c:adm-succ-half-gx-scc} of
  Lemma~\ref{lem:adm-succ}.  Hence, the half $>_{f}$-$\scc$
  implication follows.

  Third, let us confirm the upper condition of $f$ on $\fX$.
  Condition~\ref{c:adm-film-O-admissible-gen-non-const}
  gives some $j\in \Zlgez$ such that
  $f_{j_{1},\cdots,j_{l}}x_{1}^{j_{1}}\dots
  x_{l}^{j_{l}}\not\in \Q$.  Hence, $1\gAfX f$ by
  Condition~\ref{c:adm-film-O-admissible-bound}, since
  Claim~\ref{c:adm-bounds} of
  Lemma~\ref{lem:adm-bounds-poring-inclusions} says
  $0<x_{i}(r)<1$ for each $i\in\oi(l)$ and $r\in \AfX$.
  Thus, Claim~\ref{c:adm-film-O-admissible-gen} follows.

  Claim~\ref{c:adm-film-O-admissible-gen} and $1\scc 0$
  imply Claim~\ref{c:adm-film-O-admissible-monom}.
\end{proof}

We put the following notion to discuss explicit values of
rational functions over $\AfX$.

\begin{definition}\label{def:fully-adm}
  Suppose squaring orders $O=\{\sce,\scc\}$ on $\fX$.
  Let
  $x\in \Q(\fX)$ be $O$-admissible.  We call $x$ fully
  $O$-admissible by $\fX$ (or fully admissible for short),
  if there is $d\in \Zgeo$ such that
  $x^{\frac{1}{d}}\in \fX$.  Moreover, let $l\in \Zgeo$ and
  $y\in \Q(\fX)^{l}$.  If each $y_{i}$ is fully
  $O$-admissible by $\fX$, then we call $y$ fully
  $O$-admissible by $\fX$ (or fully admissible for short).
\end{definition}

\section{Notations for families and
  some $q$-analogs}
\label{sec:families-q-analogs}
For further discussion, we use the following notations of
families and some $q$-analogs.

\subsection{Families}
\label{sec:families}

 We fix some general
 notations of families.
 First, we take the following
 notation for families of rational functions
 of a filed, which is often $\Q$ or $\R$ in this manuscript.

 \begin{definition}
   \label{def:families-general-operations}
  Let $Q$ be a field.
  Consider families $F=\{F_{i}\in Q(\fX)\}_{i\in I}$
  and $F'=\{F'_{i}\in Q(\fX)\}_{i\in I}$.
  \begin{enumerate}
  \item If $F_{i}=F_{j}$ whenever $i,j\in I$, then $F$ is
    said to be flat.
  \item Let $F^{F'}=\{F_{i}^{F'_{i}}\in Q(\fX)\}_{i \in l}$
    for the exponentiation, if for each $i\in I$, $F_{i}$
    and $F'_{i}$ are not simultaneously zero and
    $F'_{i}\in \Z$.
  \item Let $F\pm F'= \{F_{i}\pm F'_{i}\}_{i\in I}$ for the
    addition and subtraction.
  \item Let $F\rc F'= \{F_{i}F'_{i}\}_{ i \in l}$ for the
    Hadamard product (term-wise product).
  \item Assume $\lam\in Q(\fX)$. Let
    $\lam \pm F=F\pm \lam= \{\lam \pm F_{i}\}_{i\in I}$ for
    the scalar product.  Also, let
    $ \lam F=F \lam= \{\lam F_{i}\}_{i\in I}$ for the scalar
    multiplication.
  \end{enumerate}

  Moreover, if $I$ is a finite set, then we put the
  following.

  \begin{enumerate}[label=(\alph*)]
  \item Let $\len(F)=\num(I)$,
    the number of elements of $I$.
    \item For $h\in Q(\fX)$,
    let $\len_{h}(F)=\num({\{i\in I\mid F_{i}=h\}})$.
  \item Let $\sum F=\sum_{i \in I} F_{i}$ and
    $\prod F=\prod_{i\in I} F_{i}$ for the sum and product of $F$.
  \end{enumerate}
\end{definition}

We then take the following notation for a family of real
numbers.  This generalizes some tuple equality and
inequality in Definition~\ref{def:cartesian}.

\begin{definition}\label{def:families-general-ineq}
  Suppose families $F=\{F_{i}\in \R\}_{i\in I}$ and
  $F'=\{F'_{i}\in \R\}_{i\in I}$. Let $\lam\in \R\cup\{\pm \infty\}$.
  \begin{enumerate}
  \item Let $F\geq F'$ (or $F>F'$), if each
    $F_{i}\geq F'_{i}$ (or $F_{i}>F_{i}'$).
  \item Let $F\geq \lam$ (or $F>\lam$), if each
    $F_{i}\geq \lam$ (or $F_{i}>\lam$).
  \item Let $F\leq \lam$ (or $F<\lam$), if each
    $F_{i}\leq \lam$ (or $F_{i}<\lam$).
  \item Let $F= \lam$, if $F\leq \lam$ and $F\geq \lam$.
  \item We call $F$ positive, negative, or zero, if $F>0$,
    $F<0$, or $F=0$, respectively.  Also, we
    call $F$ non-negative or non-positive, if $F\geq 0$ or
    $F\leq 0$, respectively.
   \item Let $\abs{F}=\{\abs{F_{i}}\in \R_{\geq 0}\}_{i\in I}$.
  \end{enumerate}
\end{definition}

Moreover, we put the following sets of bounded increasing
integer sequences.
\begin{definition}\label{def:bounded-increasing}
  Suppose $d\in \Zgeo$ and $\lam\in \Zt$.
   \begin{enumerate}
   \item Let
     $\tsi(d,\lam) =\{m\in \Zd\mid \lamo\leq
     \mn<\dots <m_{d}\leq
     \lamt\}$ and $\tsi(d,\lamt)=\tsi(d,\tp{1,\lamt})$.
   \item\label{d:bounded-increasing-leq} Let
     $\tei(d,\lam) =\{m\in \Zd\mid \lamo\leq \mn\leq \dots \leq m_{d}\leq
     \lamt\}$ and $\tei(d,\lamt)=\tei(d,\tp{1,\lamt})$.
    \end{enumerate}
  \end{definition}
  
  \subsubsection{Tuples}
  In particular, we put the following notation for tuples.
  \begin{definition}\label{def:tuples-general-operations}
    Let $d,d'\in \Zgeo$ and $\lam\in \tei(2,d)$.  For a field
    $Q$, assume $f\in Q(\fX)$, $m\in Q(\fX)^{d}$, and
    $m'\in Q(\fX)^{d'}$.
  \begin{enumerate}
  \item Let $m\ccn m'\in Q(\fX)^{d+d'}$ be the concatenation
    such that
    \begin{align}
      (m\ccn m')_{i}
      =
      \begin{dcases}
        m_{i} \mfor i\in\oi(d),\\
        m'_{i-d} \mfor i\in\oi(d+1,d+d').
      \end{dcases}
    \end{align}
    For $\mu\in \Zgeo$, let us write the $\mu$-fold
    concatenation
    $m^{\ccn \mu}=m\ccn \cdots \ccn m\in Q(\fX)^{\mu d}$.
  \item Let us write the palindromization
    $m^{\wcn}=m\mdoubleplus m^{\ve}\in Q(\fX)^{2d}$.
  \item Let
    $m(\lamo,\lamt), m(\lamt,\lamo)\in
    Q(\fX)^{\lamt-\lamo+1}$ for the sub-tuples of $m$ such
    that
    \begin{align}
      m(\lamo,\lamt)
      &=\tp{m_{\lamo+i-1}}_{i\in\oi(\lamt-\lamo+1)},\\
      m(\lamt,\lamo)
      &=m(\lamo,\lamt)^{\ve}.
    \end{align}
  \item For $l''=\lamt-\lamo+1$ and $m''\in Q(\fX)^{l''}$,
    let us write the internal addition and subtraction
    $m\pm_{\tp{\lamo,\lamt}}m'' = m(1,\lamo-1) \ccn
    (m(\lamo,\lamt) \pm m''(1,l'')) \ccn m(\lamt+1,l) \in
    Q(\fX)^{l}$.
  \end{enumerate}
\end{definition}

\subsubsection{Some tuple $q$-analogs}
We write the following tuple $q$-analogs by
Definition~\ref{def:pochh}.
\begin{definition}\label{def:g-multi}
  Let $l\in \Zgeo$ and $m,m',w\in \Zlgez$. Suppose an
  indeterminate $x\in \Q(\fX)^{l}$.
  \begin{enumerate}
  \item Let the $x$-Pochhammer symbol
    $(m)_{x}^{w} =\prod_{i\in\oi(l)}(m_{i})_{x_{i}}^{w_{i}}$,
    $x$-number
    $[m]^{w}_{x} =\prod_{i\in\oi(l)} [m_{i}]_{x_{i}}^{w_{i}}$,
    $x$-factorial
    $[m]!_{x}^{w} =\prod_{i\in\oi(l)} [m_{i}]!_{x_{i}}^{w_{i}}$,
    and $x$-binomial coefficient
    ${m \brack m'}^{w}_{x} =\prod_{i\in\oi(l)} {m_{i} \brack
      m'_{i}}_{x_{i}}^{w_{i}}$.
  \item Suppose $x=\iota^{l}(q)$ for an indeterminate
    $q\in\Q(\fX)$.  Let
    $(m)_{q}^{w}=(m)_{x}^{w}$,
    $[m]^{w}_{q}=[m]^{w}_{x}$,
    $[m]!_{q}^{w}=[m]!_{x}^{w}$,
    and ${m \brack m'}^{w}_{q}= {m \brack m'}^{w}_{x}$.
  \end{enumerate}
  In particular, we often
  omit the superscript $w$ above when $w=\iota^{l}(1)$.
\end{definition}

\section{Fitting condition
  and base shift functions}

We introduce {\it fitting tuples} and {\it base shift
  functions}, extending those in
Definitions~\ref{def:merged-simpler}
and~\ref{def:intro-base-shift}.
 
\subsection{Fitting condition}\label{sec:fitting}
We put the following notion of {\it gates} to introduce 
fitting condition.
\begin{definition}\label{def:lhd}
  Let $s\in \hbzt$.
  \begin{enumerate}
  \item We call $s$ gate, if $s_{1}<\infty$ and
    $s_{1}\leq s_{2}$.
  \item Let $s$ be a gate with $l\in \Zgeo$ and $U\subset \hbz$.
    \begin{enumerate}
    \item We call $s$ finite, if $s<\infty$.  If not, we call
    $s$ infinite.
  \item We write $m\ld_{U} s$, if
    $m\in \{ u \in U \mid s_{1}\leq u \leq s_{2}\}$.
  \item We write $m\ld_{U^{l}} s$, if
    $m \in \{ u \in U^{l} \mid s_{1}\leq u \leq s_{2}\}$. 
  \end{enumerate}
\end{enumerate}
\end{definition}    

In this manuscript,
writing $m\ld_{U} s'$ or
$m\ld_{U^{l}} s'$ presumes that $s'$ is a gate.  Also, we
assume that $\{F_{m}\in \Q(\fX)\}_{m\ld_{U}s}$ and
$\{G_{m}\in \Q(\fX) \}_{m\ld_{U^{l}}s}$ denote families, unless
stated otherwise.  Hence, if
$\{F_{m}\in \Q(\fX)\}_{m\ldZ s}=\{F'_{m}\in \Q(\fX) \}_{m\ldZ
  s'}$, then $s=s'$ and each $F_{m}=F'_{m}$.

We use the following notation for computations on fitting
tuples.
\begin{definition}\label{def:sigma-bopbom}
  For $l\in \Zgeo$, let $m\in \Q(\fX)^{l}$ and
  $k,k'\in \Q(\fX)^{2l}$.
  \begin{enumerate}
  \item \label{d:sigma-bopbom-sigma} Let
    $\sig(k)\in \Q(\fX)^{l}$ such that
    $\sig(k)_{i}=\sum k(i+1,2l-i+1)$ for $i\in\oi(l)$.
  \item \label{d:sigma-bopbom-bopbom} In $\Q(\fX)^{l}$, let
    $m\bop k= m+\sig(k)^{\ve}$ and $m\bom k=m-\sig(k)$.
    We
    call $\bop$ and $\bom$ $\sig$-plus and $\sig$-minus.
  \item We call $k$ and $k'$ $\sig$-equivalent (or
    equivalent for short), if $\sig(k)=\sig(k')$.
  \end{enumerate}
\end{definition}
Hence, for $\sig$-equivalent $k$ and $k'$, we have
$(m\bop k)_{i}=m_{i}+\sig(k)_{l-i+1}= m_{i}+\sum
k(l-i+2,l+i)=(m\bop k')_{i}$ and
$(m\bom k)_{i}=m_{i}-\sig(k)_{i}=m_{i}-\sum k(i+1,2l-i+1)
=(m\bom k')_{i}$.

\begin{remark}\label{rmk:operation-order}
  For $l\in \Zgeo$, let $a\in \Q(\fX)^{l}$.  This gives
  expressions such as $a\ccn a +a\ccn a$ or $a \rc a-a$.
  Thus, we take the following precedence: (1) flip,
  palindromization; (2) $\rc$, scalar multiplication; (3)
  $\ccn$, $\bop$, $\bom$; (4) tuple addition/subtraction;
  (5) scalar addition/subtraction.  Other than the
  precedence, we use parentheses and center dots to avoid
  confusion.
\end{remark}

We now introduce the fitting condition and related notion.
\begin{definition}\label{def:fitting}  
  For $l\in \Zgeo$, consider $m,n\in \Zl$ and $k \in \Ztl$
  with
  a gate $s\geq 0$.  Let $\mu=\tp{s,l,m,n,k}$.  We call $l$,
  $m$, $n$, and $k$ width, left ladder, right ladder, and
  support of $\mu$.  Also, we call $k_{1}$ free parameter of
  $\mu$.
  \begin{enumerate}
  \item\label{d:fitting-nu} We define $a=\nu(k)$ and
    $b=\nu(m,n,k)$ in $\Ztl$ such that $a_{i}=\sum k(1,i)$ for
    $i\in\oi(2l)$ and $b=a+m\ccn n$.
  \item
    \label{d:fitting-fitting}
    We call $\mu$ fitting, if its ladders and support satisfy
    the following:
    \begin{align}
      m\ccn n&\ldZtl s;
      \label{ineq:fit-inc}\\
      b_{1}\leq\dots\leq b_{l}<b_{l+1}\leq \dots &\leq b_{2l};
      \label{ineq:fit-top}\\
      0\leq a_{1}\leq \dots\leq a_{l}<a_{l+1}\leq \dots &\leq a_{2l}.
      \label{ineq:fit-bottom}
    \end{align}
    We refer to \eqref{ineq:fit-inc}, \eqref{ineq:fit-top},
    and~\eqref{ineq:fit-bottom} as the inclusion condition,
    upper slope condition, and lower slope condition of
    $\mu$.
  \item\label{d:fitting-wrapped} We call $\mu$ wrapped, if
    $(m\bom k )\ccn (n\bop k) \ldZtl s$.
    \end{enumerate}
\end{definition}

\begin{remark}\label{rmk:free-param}
  Let $\mu=\tp{s,l,m,n,k}$ be fitting.  We call $k_{1}$ free
  parameter, since for each $k_{1}+\lam\in \Zgez$,
  $\mu'=\tp{s,l,m,n,\tp{k_{1}+\lam}\ccn k(2,l)}$ is
  fitting. Also, $\mu$ is wrapped if and only if $\mu'$ is
  wrapped. However, the free parameter is still important,
  since it gives different polynomials with positive integer
  coefficients later by the merged-log-concavity.
\end{remark}

\begin{example}\label{ex:ab}
  Let $l=2$. Suppose a fitting $\tp{s,l,m,n,k}$ with
  $a=\nu(k)$ and $b=\nu(m,n,k)$.  Then, there is the following
  diagram: \def\rgeq{\rotatebox[origin=c]{-90}{$\geq$}}
  \begin{align}
    \begin{array}{lllll}
      &  b_{1}  \leq     &  b_{2}   <    & b_{3}   \leq & b_{4}\\
      & \rgeq \mn     & \rgeq\mt    & \rgeq \nn    & \rgeq \nt\\
      0\stackrel{\kn}{\leq} &  a_{1}   \stackrel{\kt}{\leq}  &  a_{2}
                                                        \stackrel{\kr}{<} & a_{3}
                                           \stackrel{k_{4}}{\leq }  & a_{4},
    \end{array}
  \end{align}
  where numbers along inequalities indicate corresponding
  differences.
\end{example}

Let us the following lemmas for our later discussion.
\begin{lemma}\label{lem:fitting-nonneg}
  Suppose a fitting $\mu=\tp{s,l,m,n,k}$ with $a=\nu(k)$ and
  $b=\nu(m,n,k)$.  Then, we have the following.
  \begin{enumerate}
  \item 
    $k=\tp{a_{1},a_{2}-a_{1},\cdots, a_{2l}-a_{2l-1}}\geq 0$.
    \label{c:fitting-nonneg-k-nonneg}
  \item   $b\geq a\geq k_{1}\geq 0$.
    \label{c:fitting-nonneg-ab-nonneg}
  \end{enumerate}
\end{lemma}
\begin{proof}
  Claim~\ref{c:fitting-nonneg-k-nonneg} follows, since the
  lower slope condition says $\kn=\sum k(1,1)=a_{1}\geq 0$ and
  $k_{i}=\sum k(1,i)-\sum k(1,i-1)=a_{i}-a_{i-1}\geq 0$ for
  $i\in\oi(2,2l)$.  Claim~\ref{c:fitting-nonneg-ab-nonneg}
  holds, since $a\geq k_{1}$ by
  Claim~\ref{c:fitting-nonneg-k-nonneg} and
  $b-a=m\ccn n\geq 0$ by the inclusion condition.
\end{proof}

\begin{lemma}\label{lem:sigma-nu}
  For $l\in \Zgeo$, let $m,n\in \Zl$ and $k\in \Ztl$ with
  $a=\nu(k)$ and $b=\nu(m,n,k)$.  Then, each $i\in\oi(l)$ gives
  $a_{2l-i+1}-a_{i}=\sig(k)_{i}$ and
  $b_{2l-i+1}-b_{i}=n_{l-i+1}+\sig(k)_{i}-m_{i}$.
\end{lemma}
\begin{proof}
  We have
  $a_{2l-i+1}-a_{i}=\sum k(1,2l-i+1)-\sum k(1,i)=\sig(k)_{i}$.
  Also,
  $b_{2l-i+1}-b_{i} =\sum k(1,2l-i+1) +n_{l-i+1}-(\sum
  k(1,i)+m_{i}) =\sig(k)_{i}+n_{l-i+1}-m_{i}$.
\end{proof}

\begin{lemma}\label{lem:slope-conditions-sigma}
  Suppose a fitting $\mu=\tp{s,l,m,n,k}$.  Then, we have
  the
  following inequalities:
  \begin{align}
    \sig(k)_{1}\geq \sig(k)_{2}\geq
    \dots &\geq\sig(k)_{l}=k_{l+1}> 0;
\label{ineq:slope-conditions-sigma-nonneg-odd}\\
    n_{l}+\sig(k)_{1}-\mn\geq
    n_{l-1}+\sig(k)_{2}-\mt\geq \dots&\geq \nn+\sig(k)_{l}-m_{l}>0.
\label{ineq:slope-conditions-sigma-nonneg-even}
  \end{align}
  In particular, if $l=1$, then the upper slope condition of
  $\mu$ is the same as $\nn+k_{l+1}-m_{l}> 0$.
\end{lemma}
\begin{proof}
  We have $\sig(k)_{l}=\sum k(l+1,l+1)=k_{l+1}$.  Then,
inequalities~\eqref{ineq:slope-conditions-sigma-nonneg-odd}
and~\eqref{ineq:slope-conditions-sigma-nonneg-even} follow
  from Lemma~\ref{lem:sigma-nu} and the slope conditions of
  $\mu$.  Let $b=\nu(m,n,k)$.  When $l=1$, the upper slope
  condition is $b_{2}-b_{1}>0$, which is equivalent to
  $\nn+k_{l+1}-m_{l}> 0$ by Lemma~\ref{lem:sigma-nu}.
\end{proof}
 
In particular, we have the following equivalence for the
width-one fitting condition.

\begin{lemma}\label{lem:fitting-one}
  For $l=1$, suppose $m,n\in \Zl$ and $k\in \Ztl$ with a gate
  $s\geq 0$.  Then, $\mu=\tp{s,l,m,n,k}$ is fitting if and only
  if $m\ccn n\ldZtl s$, $k\geq \tp{0,1}$, and $n+\kt> m$.
\end{lemma}
\begin{proof}
  Let $a=\nu(k)$ and $b=\nu(m,n,k)$.  First,
  $m\ccn n\ldZtl s$ is the inclusion condition of $\mu$.
  Second, the lower slope condition of $\mu$ is
  $a_{1}=\kn\geq 0$ and $a_{2}-a_{1}=\kt\geq 1$.  Third, we have
  the equivalence between the upper slope condition and
  $\nn+\kt-\mn>0$ in
  Lemma~\ref{lem:slope-conditions-sigma}.
\end{proof}

In Section~\ref{sec:strictly-fitting}, we discuss a variant
of the fitting condition, demanding strict inequalities on
the slope conditions.

\subsection{Base shift functions}\label{sec:base-change}
We extend 
$b_{\lam,\rho}(q)$  in Definition~\ref{def:intro-base-shift}.
\begin{definition}\label{def:bshift-fun}
  Let $\lamo\in \Zgez$, $\lamt\in \Z$, and $\lamr\in \Zgeo$
  with a gate $s\geq 0$.  Consider an indeterminate
  $x\in \Q(\fX)$ and $0\neq \phi(x)\in \Q(x)$.  In
  $\Q(\fX)$, we define the base shift function
  \begin{align}
    b(s,\lamo,\lamt,\phi,\lamr,x,\fX)
    &=
      \begin{dcases}
        \left(\frac{ \phi(x^{\lamr})^{\lamt}
            [\lamt]!_{x^{\lamr}}} { \phi(x)^{\lamt}
            [\lamt]!_{x}}\right)^{\lamo}
        \mif \lamt\ldZ s,\\
        0 \melse.
      \end{dcases}
  \end{align}
\end{definition}
We state some of their special values.
\begin{lemma}\label{lem:bshift-fun-sp}
  Assume either $\lamo=0$, $\lamt\leq 0$, or $\lamr=1$.  Then,
  we have
  \begin{align}
    b(s,\lamo,\lamt,\phi,\lamr,x,\fX)
    &=
      \begin{dcases}
        1 
        \mif \lamt\ldZ s,\\
        0 \melse.
      \end{dcases}
  \end{align}
\end{lemma}
\begin{proof}
  Assume $\lamt\ldZ s$.  When $\lamr=1$,
  $\phi(x^{\lamr})^{\lamt} [\lamt]!_{x^{\lamr}}=
  \phi(x)^{\lamt} [\lamt]!_{x}$.  When $\lamo=0$ or
  $\lamt=0$,
  $\left(\phi(x^{\lamr})^{\lamt}
    [\lamt]!_{x^{\lamr}}\right)^{\lamo}=
  \left(\phi(x)^{\lamt} [\lamt]!_{x}\right)^{\lamo}
  =1$. Hence, the assertion follows.
\end{proof}

For $l\in \Zgeo$, $x\in \Q(\fX)^{l}$, and
$\phi(x)\in \prod_{i\in\oi(l)}\Q(x_{i})$, we write
\begin{align}
  \phi(x)_{i}=\phi_{i}(x_{i})\in \Q(x_{i}).
\end{align}
We then make the following tuple version of base shift
functions.
\begin{definition}\label{def:bshift-fun-multi}
  Assume $l\in \Zgeo$, $\lamo\in \Zlgez$, $\lamt\in \Z^{l}$,
  and $\lamr\in \Zgeo^{l}$ with a gate $s\geq 0$. For an
  indeterminate $x\in \Q(\fX)^{l}$, let
  $\phi(x)\in \prod_{i\in\oi(l)}\Q(x_{i})$ and
  $\phi_{i}(x_{i})\neq 0$ for each $i\in\oi(l)$.
  \begin{enumerate}
  \item \label{d:bshift-fun-multi-one-tuple} We define the
    base shift function
    \begin{align}
      B(s,l,\lamo,\lamt,\phi,\lamr,x,\fX)
      = \prod_{i\in\oi(l)}
      b(s,\lam_{1,i},\lam_{2,i},\phi_{i},\lam_{3,i},x_{i},\fX) \in \Q(\fX).
    \end{align}
  \item \label{d:bshift-fun-multi-two-tuples} Let
    $m,n\in \Zl$.  We define the base shift function
    \begin{align}
      B(s,l,w,m,n,\phi,\rho,x,\fX)
      =B(s,l,w,m,\phi,\rho,x,\fX)
      B(s,l,w,n,\phi,\rho,x,\fX)\in \Q(\fX).
    \end{align}
  \end{enumerate}
\end{definition}

\subsection{Mediators}

Since a strict squaring order $\scc$ presumes a squaring
order $\sce$, we often denote squaring orders
$\{\sce, \scc\}$ just by $\scc$.  Then, we introduce {\it
  mediators}, which extend $q$-Pochhammer symbols in
Definition~\ref{def:pochh} within our purpose.

\begin{definition}\label{def:mediator}
  Consider a gate $s\geq 0$, $l\in \Zgeo$, $w\in \Zlgez$, and
  $\rho\in \Zgeo^{l}$.  Suppose a $\scc$-admissible
  $x\in \Q(\fX)^{l}$ and
  $\phi(x)\in \prod_{i\in\oi(l)}\Q(x_{i})$.  Let
  $\mu=\tp{s,l,w,\scc,\rho,x,\fX}$.  We call $\phi$
  $\mu$-mediator (or mediator for short), if $\phi$ and
  $\mu$ satisfy the following conditions.
  \begin{enumerate}
  \item For each $i\in\oi(l)$, $\phi(x)_{i}^{w_{i}}\gAfX 0$.
    \label{d:mediator-base-pos}
  \item For each $m\ldZl s$,
    $B(s,l,w,m,\phi,\rho,x,\fX)\scc 0$.
    \label{d:mediator-base-shift-pos}
  \end{enumerate}
  We refer to \ref{d:mediator-base-pos}
  and~\ref{d:mediator-base-shift-pos} as the base positivity
  and base-shift positivity of $\phi$ and $\mu$.
\end{definition}

For a $\scc$-admissible variable
$q\in \Q(\fX)$, we have the squaring order $\llq$ of
$q$-polynomials by Theorem~\ref{thm:adm}.
Then, we introduce the notion of {\it canonical mediators} by
 the following proposition.

\begin{proposition}\label{prop:canon-med}
  Let $l=1$. Suppose a $\scc$-admissible $q\in \Q(\fX)$ and
  $x=\tp{q}$.  Let
  $\mu_{s,w,\rho}=\tp{s,l,w,\llq,\rho,x,\fX}$ for each gate
  $s\geq 0$, $w\in \Zlgez$, and $\rho\in \Zgeo^{l}$.  If
  $\phi(x)\in\Q[x]^{l}$ satisfies $\deg_{q}\phi(x)_{1}=1$, then the
  following statements are equivalent.
  \begin{enumerate}
  \item \label{s:canon-med-for-any-gate-shift-weight} For
    each $\mu_{s,w,\rho}$, $\phi$ is a
    $\mu_{s,w,\rho}$-mediator.
  \item \label{s:canon-med-deg-one} For some
    $u\in \Q_{> 0}$, $\phi(x)_{1}=u(1-q)$.
  \end{enumerate}
\end{proposition}
\begin{proof}
  First, let us prove that
  Statement~\ref{s:canon-med-deg-one} implies
  Statement~\ref{s:canon-med-for-any-gate-shift-weight}.
  The base positivity $\phi(x)_{1}\gAfX 0$ holds, because
  $0<q(r)<1$ for each $r\in \AfX$ by Claim~\ref{c:adm-bounds}
  of Lemma~\ref{lem:adm-bounds-poring-inclusions}.  For each
  $m\ldZl s$, let us prove the base-shift positivity:
  \begin{align}
    B(s,l,w,m,\phi,\rho,x,\fX)\llq 0. \label{ineq:canon-med-p1}
  \end{align}
  Assume $s=\tp{0,\infty}$ without loss of generality.  When
  $m=\tp{0}\ldZl s$, Lemma~\ref{lem:bshift-fun-sp} gives
  $B(s,l,w,m,\phi,\rho,x,\fX)=1$.  Thus, assume otherwise.
  Then, inequality~\eqref{ineq:canon-med-p1} holds by
  \begin{align}
    B(s,l,w,m,\phi,\rho,x,\fX)
    =
    \frac{ (u(1-q^{\rhoo}))^{\mn w_{1}}
    [\mn]!_{q^{\rhoo}}^{w_{1}}}
    { (u(1-q))^{\mn w_{1}}      [\mn]!_{q}^{w_{1}}}
    =
    \frac{(\mn)_{q^{\rhoo}}^{w_{1}}}{(\mn)_{q}^{w_{1}}}
    =
    \prod_{h\in\oi(\mn)}[\rhoo]^{w_{1}}_{q^{h}}.
  \end{align}
  
  Second, let us prove the converse. Thus, let
  $\phi(x)_{1}=\ka_{1}-\ka_{2}q$ for some
  $\ka_{1},\ka_{2}\in \Q$ such that $\ka_{2}\neq 0$ by
  $\deg_{q}\phi(x)_{1}=1$.  Also, let $s=\tp{0,\infty}$.  Then,
  for $w=\tp{1}$, $m=\tp{2}$, and $\rho=\tp{2}$, the base-shift
  positivity of $\phi$ and $\mu$ implies
  \begin{align}
    \frac{(\ka_{1}-\ka_{2}q^2)^{2}(1+q^2)}{
    (\ka_{1}-\ka_{2}q)^{2}(1+q)}\llq0.
  \end{align}
  Hence, $\ka_{1}=\ka_{2}$, since we need
  $\ka_{1}-\ka_{2}(-1)^2=0$ when $q=-1$.  Furthermore,
  $\phi(x)^{w_{1}}_{1} =\ka_{1}(1-q) \gAfX 0$ by the base
  positivity of $\phi$ and $\mu$.  Therefore,
  $\ka_{1}\in \Qgo$, since $0<q(r)<1$ for each $r\in \AfX$ by
  Claim~\ref{c:adm-bounds} of
  Lemma~\ref{lem:adm-bounds-poring-inclusions}.
\end{proof}

Then,
as $\phi(\tp{0})_{1}\neq 0$,
we introduce the following notion by $\phi(\tp{0})_{1}=1$.

\begin{definition}\label{def:canonical-pos}
  Let $l\in \Zgeo$. If
  $\phi(x)=\tp{1-x_{i}}_{i\in\oi(l)} \in \prod_{i\in\oi(l)}\Q(x_{i})$
  for an indeterminate $x\in \Q(\fX)^{l}$, then we call
  $\phi$ canonical $l$-mediator (or canonical mediator
  for
  short).
\end{definition}
In particular, when $l=1$, the canonical $\phi(x)=\tp{1-x}\in \Q(\fX)^{l}$ gives
 $q$-Pochhammer symbols
 $\phi(q)_{1}[m]_{q}=(m)_{q}$ for $m\in \Zgez$.
 Hence,  mediators
 generalize  $q$-Pochhammer symbols
  with  the base
  and base-shift positivities.

\section{Merged-log-concavity}
\label{sec:merged}
We introduce the notions of ring shift factors, merged
determinants, pair-wise positive rational functions, parcels, and the merged-log-concavity.

\begin{definition}\label{def:merged} 
  Suppose a gate $s\geq 0$.  Let $l\in \Zgeo$,
  $w\in \Zlgez$, and $\rho\in \Zgeo^{l}$.  Assume
  squaring orders $O=\{\sce, \scc\}$ on $\fX$.  Consider 
  a $\tp{s,l,w,\scc,\rho,x,\fX}$-mediator $\phi$
  for an
  $O$-admissible $x\in \Q(\fX)^{l}$.
  \begin{enumerate}
  \item \label{d:merged-factor-det} Let $m,n\in \Zl$ and
    $k\in \Ztl$ with $a=\nu(k)$ and $b=\nu(m,n,k)$.
    Suppose $y=x^{\rho}$.
    Then, in
    $\Q(\fX)$, we define the ring shift factor
    \begin{align}
      \Ups(s,l,w,m,n,k,\phi,\rho,x,\fX)
      =
      \begin{dcases}
        \prod(\phi(y)^{\wcn})^{(b-a)\rc w^{\wcn}}
        \cdot
        \frac{[b]!_{y^{\wcn}}^{w^{\wcn}}
        }{[a]!_{y^{\wcn}}^{w^{\wcn}}}
        \mif a,b\geq 0, \\
        0 \melse.
      \end{dcases}
    \end{align}
    Moreover, suppose a family
    $\cF=\{\cF_{m}\in \Q(\fX)\}_{m\in\Zl}$.  Then, in
    $\Q(\fX)$, we define the merged determinant
    \begin{align}
      \Delta(\cF)(s,l,w,m,n,k,\phi,\rho,x,\fX)
      =
      \Ups(s,l,w,m,n,k,\phi,\rho,x,\fX)
      \det \begin{bmatrix}
        \cF_{m}&         \cF_{(n\bop k)^{\ve}}\\
        \cF_{m\bom k} &         \cF_{n^{\ve}}
      \end{bmatrix}.
    \end{align}
\item We call a
   family $\fs=\{f_{s,m}\in \Q(\fX)\}_{m\in \Zl}$ pairwise
   $\tp{s,l,\scc}$-positive (or pairwise positive for short),
   if
   $f_{s,m}f_{s,n}\scc 0$ for each
   $m,n \ldZl s$ and $ f_{s,m}= 0$ for each
   $m \nldZl s$. 
  \item \label{d:merged-parcel} Suppose a pairwise
    $\tp{s,l,\scc}$-positive
    $\fs=\{f_{s,m}\in \Q(\fX)\}_{m\in \Zl}$.  Then, we define
    the parcel
    $\cF=\Lam(s,l,w,\scc,\fs,\phi,\rho,x,\fX) = \{\cF_{m}\in
    \Q(\fX)\}_{m\in \Zl}$ as a family of rational functions
    such that
    \begin{align}
      \cF_{m}
      =
      \begin{dcases}
        \frac{f_{s,m}}{\prod \phi(x)^{m\rc w}\cdot [m]!_{x}^{w}}
        \mif m\ldZl s, \\
        0  \melse.  
      \end{dcases}
    \end{align}
    We refer to $s$, $l$, $w$, $\scc$, $\fs$, $\phi$,
    $\rho$, $x$, and $\fX$ as the gate, width, weight,
    strict squaring order, numerator, mediator, base shift,
    base, and coordinate of $\cF$.  We call them parcel
    parameters of $\cF$.
  \item \label{d:merged-log-concavity} Suppose 
    $\cF=\Lam(s,l,w,\scc,\fs,\phi,\rho,x,\fX)$ with squaring
    orders $O'=\{\sce',\scc'\}\Sup O$.
    \begin{enumerate}
    \item We call $\cF$
      $\tp{s,l,w,\scc',\phi,\rho,x,\fX}$-merged-log-concave (or
      $\scc'$-merged-log-concave for short), if any fitting
      $\tp{s,l,m,n,k}$ gives
      \begin{align}
        \Delta(\cF)(s,l,w,m,n,k,\phi,\rho,x,\fX)\scc' 0.
      \end{align}
    \item Similarly, we call $\cF$
      $\tp{s,l,w,\sce',\phi,\rho,x,\fX}$-merged-log-concave (or
      $\sce'$-merged-log-concave for short), if any fitting
      $\tp{s,l,m,n,k}$ gives
      \begin{align}
        \Delta(\cF)(s,l,w,m,n,k,\phi,\rho,x,\fX)\sce' 0.
      \end{align}
    \end{enumerate}
    We also refer to $\scc'$-merged-log-concavity and
    $\sce'$-merged-log-concavity as {\it strict
      merged-log-concavity} and {\it non-strict
      merged-log-concavity}, if no confusion occurs.
  \end{enumerate}
\end{definition}

Since $\Ups(s,l,w,m,n,k,\phi,\rho,x,\fX)=1$ when $w=\tp{0}$,
merged determinants extend $2\times 2$-determinants by
ring shift factors. Moreover, we discuss 
general minors in Theorem~\ref{thm:rshift-ext-cb}.

Suppose width-$l$ parcels $\cF, \cF'$.
As families of
 rational functions, we consider $\cF=\cF'$ if
 $ \cF_{m}=\cF'_{m} \in \Q(\fX)$ for each $m\in \Zl$, even
 with different parcel parameters.
 We  discuss the change of parcel parameters and
 the merged-log-concavity in
  Claim~\ref{c:merged-param-change-merged} of
 Proposition~\ref{prop:merged-param-change}.

 We simplify the notation of
 ring shift factors, parcels, and merged determinants in
 the following cases.
\begin{definition}\label{def:merged-abbrev}
  Suppose a parcel $\cF=\Lam(s,l,w,\scc,\fs,\phi,\rho,x,\fX)$.
  \begin{enumerate}
  \item \label{d:merged-abbrev-canon} If $\phi$ is the
    canonical mediator, then we write
    \begin{align}
      \Ups(s,l,w,m,n,k,\rho,x,\fX)
      &= \Ups(s,l,w,m,n,k,\phi,\rho,x,\fX),\\
      \Delta(\cF)(s,l,w,m,n,k,\rho,x,\fX)
      &=\Delta(\cF)(s,l,w,m,n,k,\phi,\rho,x,\fX),\\
      \Lam(s,l,w,\scc,\fs,\rho,x,\fX)
      &= \Lam(s,l,w,\scc,\fs,\phi,\rho,x,\fX).
    \end{align}
  \item \label{d:merged-abbrev-triv-shift}
    When $\rho=\iota^{l}(1)$, we write
    \begin{align}
      \Ups(s,l,w,m,n,k,\phi,x,\fX)
      &=\Ups(s,l,w,m,n,k,\phi,\rho,x,\fX),\\
      \Delta(\cF)(s,l,w,m,n,k,\phi,x,\fX)
      &=\Delta(\cF)(s,l,w,m,n,k,\phi,\rho,x,\fX),\\
      \Lam(s,l,w,\scc,\fs,\phi,x,\fX)
      &= \Lam(s,l,w,\scc,\fs,\phi,\rho,x,\fX).
    \end{align}
  \item When $\rho=\iota^{l}(1)$ and $\phi$ is the canonical
    mediator, we write
    \begin{align}
      \Ups(s,l,w,m,n,k,x,\fX)
      &=\Ups(s,l,w,m,n,k,\phi,\rho,x,\fX),\\
      \Delta(\cF)(s,l,w,m,n,k,x,\fX)
      &=\Delta(\cF)(s,l,w,m,n,k,\phi,\rho,x,\fX),\\
      \Lam(s,l,w,\scc,\fs,x,\fX)
      &=\Lam(s,l,w,\scc,\fs,\phi,\rho,x,\fX).
    \end{align}
  \item\label{d:merged-abbrev-weight-zero}
    When $w=\iota^{l}(0)$, we write
    \begin{align}
      \Ups(s,l,m,n,k,\fX)
      &=\Ups(s,l,w,m,n,k,\phi,\rho,x,\fX),\\
      \Delta(\cF)(s,l,m,n,k,\fX)
      &=\Delta(\cF)(s,l,w,m,n,k,\phi,\rho,x,\fX),\\
      \Lam(s,l,\scc,\fs,\fX)
      &=\Lam(s,l,w,\scc,\fs,\phi,\rho,x,\fX).
    \end{align}
  \end{enumerate}
\end{definition}

\section{Preliminary discussions
  on the merged-log-concavity}
\label{sec:prelim}
\subsection{On mediators}
We obtain
the following $q$-Pochhammer symbols
 by canonical mediators.

\begin{proposition}\label{prop:lam-poch}
  Let
  $\cF=\Lam(s,l,w,\scc,\fs,\rho,x,\fX)$.  Then, we have
  \begin{align}
    \cF_{m}
    &=
      \begin{dcases}
        \frac{f_{s,m}}{(m)^{w}_{x}}
        \mfor m\ldZl s,\\
        0 \melse.
      \end{dcases}
      \label{eq:lam-pochh}
  \end{align}
  Moreover, let $m,n\in \Zl$, $k\in \Ztl$,
  $\rho\in \Zgeo^{l}$, and $y=x^{\rho}$ with $a=\nu(k)$ and
  $b=\nu(m,n,k)$.  Then, we have
  \begin{align}
    \Delta(\cF)(s,l,w,m,n,k,\rho,x,\fX)
    &=
      \begin{dcases}
        \frac{(b)_{y^{\wcn}}^{w^{\wcn}}}
        {(a)_{y^{\wcn}}^{w^{\wcn}}}
        ( \cF_{m} \cF_{n^{\ve}}      -  \cF_{m\bom k}
        \cF_{(n\bop k)^{\ve}}) \mif a,b\geq 0,\\
        0 \melse.
      \end{dcases}
      \label{eq:lam-pochh-merged}
  \end{align}
\end{proposition}
\begin{proof}
  Let
  $\phi(x)=\tp{1-x_{i}}_{i\in\oi(l)}\in \prod_{i\in\oi(l)}\Q(x_{i})$.
  First, we prove equation~\eqref{eq:lam-pochh}.  If
  $m\nldZl s$, then $\cF_{m}=0$.
  If $m\ldZl s$, then equation~\eqref{eq:lam-pochh} follows
  from $m\geq 0$ and
  $\prod \phi(x)^{m\rc w}\cdot [m]!_{x}^{w} =\prod_{i\in\oi(l)}
  (1-x_{i})^{m_{i}w_{i}}[m_{i}]!_{x_{i}}^{w_{i}}
  =(m)_{x}^{w}$.
  
  Second, we prove equation~\eqref{eq:lam-pochh-merged}.  If
  $a,b\not\geq 0$, then 
  equation~\eqref{eq:lam-pochh-merged} follows from
  $\Ups(s,l,w,m,n,k,\rho,x,\fX)=0$.  Let $a,b\geq 0$.
  Then, we have
  \begin{align}
    (\phi(y)^{\wcn})^{b\rc w^{\wcn}}\cdot [b]!_{y^{\wcn}}^{w^{\wcn}}
    &=\prod_{i\in\oi(l)}(1-y_{i})^{b_{i}w_{i}}[b_{i}]_{y_{i}}^{w_{i}}
      \cdot
      \prod_{i\in\oi(l)}
      (1-y_{l-i+1})^{b_{i+l}w_{l-i+1}}
      [b_{i+l}]_{y_{l-i+1}}^{w_{l-i+1}}\\
    & =\prod_{i\in\oi(l)}(b_{i})_{y_{i}}^{w_{i}}
      (b_{i+l})_{y_{l-i+1}}^{w_{l-i+1}}\\
    &=(b)_{y^{\wcn}}^{w^{\wcn}}.
  \end{align}
  Similarly,
  $(\phi(y)^{\wcn})^{a\rc w^{\wcn}} \cdot
  [a]!_{y^{\wcn}}^{w^{\wcn}} =(a)_{y^{\wcn}}^{w^{\wcn}}$.
  Hence, we  obtain equation~\eqref{eq:lam-pochh-merged}.
\end{proof}

Also, 
trivial base shifts give
the following invariance 
of the merged-log-concavity
on the choices of mediators.
\begin{proposition}\label{prop:mediator-equiv}
  Consider $\cF=\Lam(s,l,w,\scc,\fs,\phi,x,\fX)$ and
  $\cG=\Lam(s,l,w,\scc,\fs,\psi,x,\fX)$.  For a
  fitting $\mu=\tp{s,l,m,n,k}$, let $a=\nu(k)$ and
  $b=\nu(m,n,k)$.  Then, we have
  \begin{align}
    \frac{\prod (\phi(x)^{\wcn})^{(b-a)\rc w^{\wcn}}}{
    \prod \phi(x)^{m\rc w} \cdot \prod \phi(x)^{n^{\ve}\rc w}}
    &=
      \frac{\prod (\psi(x)^{\wcn})^{(b-a)\rc w^{\wcn}}}{
      \prod \psi(x)^{m\rc w}\cdot \prod \psi(x)^{n^{\ve}\rc w} }
      =1,
      \label{eq:mediator-equiv-phi-psi}\\
    \frac{\prod (\phi(x)^{\wcn})^{(b-a)\rc w^{\wcn}}}{
    \prod \phi(x)^{(m\bom k)\rc w} \cdot
    \prod \phi(x)^{(n\bop k)^{\ve}\rc w}}
    &=\frac{\prod (\psi(x)^{\wcn})^{(b-a)\rc w^{\wcn}}}{
      \prod \psi(x)^{(m\bom k)\rc w}
      \cdot \prod \psi(x)^{(n\bop k)^{\ve}\rc w} }
      =1.
      \label{eq:mediator-equiv-phi-psi-bopbom}
  \end{align}
  In particular,
  $\cF$ is $\scc'$-merged-log-concave if and only if
  $\cG$ is  $\scc'$-merged-log-concave. Similarly,
  $\cF$ is $\sce'$-merged-log-concave if and only if
  $\cG$ is  $\sce'$-merged-log-concave.
\end{proposition}
\begin{proof}
  Let us prove equations~\eqref{eq:mediator-equiv-phi-psi}.
  Because
  $(\phi(x)^{n^{\ve}\rc w})^{\ve}
  =(\phi(x)^{\ve})^{(n^{\ve}\rc w)^{\ve}}
  =(\phi(x)^{\ve})^{n\rc w^{\ve}}$,
  we have
  $\prod \phi(x)^{n^{\ve}\rc w}
  =\prod (\phi(x)^{n^{\ve}\rc w})^{\ve}
  =\prod (\phi(x)^{\ve})^{n\rc w^{\ve}}$.
  This gives
  $\prod \phi(x)^{m\rc w} \cdot \prod \phi(x)^{n^{\ve}\rc w}
    =\prod (\phi(x)\ccn \phi(x)^{\ve})^{(m\ccn n) \rc (w\ccn w^{\ve})}
    =\prod (\phi(x)^{\wcn})^{(m\ccn n) \rc w^{\wcn}}$.
  Thus, we obtain
  equations~\eqref{eq:mediator-equiv-phi-psi} by
  $m\ccn n=b-a$.
  
  Equations~\eqref{eq:mediator-equiv-phi-psi-bopbom} follow
  from equations~\eqref{eq:mediator-equiv-phi-psi}, because
  \begin{align}
    \prod \phi(x)^{(m\bom k)\rc w} \cdot \prod
    \phi(x)^{(n\bop k)^{\ve}\rc w}
    &=\prod_{i\in\oi(l)}
      \phi(x)_{i}^{(m_{i}-\sig(k)_{i}) w_{i}}
      \cdot
      \phi(x)_{i}^{(n_{l-i+1}+\sig(k)_{i}) w_{i}}\\
    &=\prod \phi(x)^{m\rc w} \cdot \prod \phi(x)^{n^{\ve}\rc w}.
  \end{align}
  
  Latter statements hold by
  equations~\eqref{eq:mediator-equiv-phi-psi}
  and~\eqref{eq:mediator-equiv-phi-psi-bopbom}, since
  when $\mu$ is fitting,
  \begin{dmath*}
    \Ups(s,l,w,m,n,k,\phi,x,\fX)\cF_{m}\cF_{n^{\ve}}   
    =
    \frac{\prod (\phi(x)^{\wcn})^{(b-a)\rc w^{\wcn}}}{
      \prod \phi(x)^{m\rc w} \cdot \prod \phi(x)^{n^{\ve}\rc w}}
    \frac{[b]!_{x^{\wcn}}^{w^{\wcn}}
    }{[a]!_{x^{\wcn}}^{w^{\wcn}}}
    \frac{f_{m}f_{n^{\ve}}}{
      [m]!_{x}^{w}[n^{\ve}]!_{x}^{w}},
  \end{dmath*}
  and when $\mu$ is wrapped and fitting,
  \begin{dmath*}
    \Ups(s,l,w,m,n,k,\phi,x,\fX)\cF_{m\bom k}\cF_{(n\bop k)^{\ve}}
    =
    \frac{\prod (\phi(x)^{\wcn})^{(b-a)\rc w^{\wcn}}}{
      \prod \phi(x)^{(m\bom k)\rc w}
      \cdot \prod \phi(x)^{(n\bop k)^{\ve}\rc w} }
    \frac{[b]!_{x^{\wcn}}^{w^{\wcn}}
    }{[a]!_{x^{\wcn}}^{w^{\wcn}}}
    \frac{f_{m\bom k}f_{(n\bop k)^{\ve}}}
    {[m\bom k]!_{x}^{w}[(n\bop k)^{\ve}]!_{x}^{w}}.
  \end{dmath*}
\end{proof}
For parcels with trivial base shifts,
Proposition~\ref{prop:mediator-equiv} gives 
\begin{align}
  \Delta(\cF)(s,l,w,m,n,k,\phi,x,\fX)
  =\Delta(\cG)(s,l,w,m,n,k,\psi,x,\fX)\in \Q(\fX)
\end{align}
for each fitting $\tp{s,l,m,n,k}$.
However, in general, 
\begin{align}
  \cF_{m}(r) = \cG_{m}(r)\in \R
\end{align}
is not true for $r\in \AfX$.  Thus, the choices of
mediators matter for explicit
$\cF_{m}(r),\cG_{m}(r)\in \R$.

\subsection{On 
  coordinates}
\label{sec:choices}
Suppose a parcel $\cF=\Lam(s,l,w,\scc,\fs,\phi,\rho,x,\fXo)$. Then, the coordinate $\fX_{1}$ is important
for us to compute $\cF_{m}(r)\in \R$ of
$r\in \AfXo$.  Hence, we
consider the change of coordinates
$\fXt\subset \Q(\fXo)$
  such that
$\fXt$ has powers of some elements of $\fXo$, and
$\cF$ is still a parcel on $\Q(\fXt)$.

First, we introduce the notion of {\it faithful squaring
  orders}.
\begin{definition}\label{def:faithful}
  Assume $\Q(\fXo)$ for $\fXo=\{X_{1,i}\}_{i\in\oi(L_{1})}$.
  \begin{enumerate}
  \item For $L_{2}\in\oi(L_{1})$, let
    $\lamo \in \tsi(L_{2},L_{1})$ and
    $\lamt\in \Zgeo^{L_{2}}$. Then, we put the set
    $\fX_{1,\lamo,\lamt}=
    \{X_{1,\lam_{1,j}}^{\lam_{2,j}}\}_{j\in\oi(L_{2})}$,
    which we call restricted coordinate of $\fX_{1}$.  Also,
    let $\fX_{1,\lamo}= \fX_{1,\lamo,\iota^{L_{2}}(1)}$.
  \item Assume squaring orders $O_{1}=\{\sceo,\scco\}$ on
    $\fXo$.  We call $O_{1}$ faithful, if for any
    $\fXt= \fX_{1,\lamo,\lamt}$, binary
    relations $\scet=\sce\mid_{\Q(\fXt)}$
    and $\scct=\scc\mid_{\Q(\fXt)}$ give the
    $\scet$-$\scct$-poring equality $\StcefXt=\StccfXt$.
  \end{enumerate}
\end{definition}

We confirm the existence of faithful squaring orders.
\begin{lemma}\label{lem:faithful}
  Assume squaring orders $O_{1}=\{\sceo,\scco\}$ on $\fXo$.
  For $\fXt=\fX_{1,\lamo,\lamt}$, let
  $\scet=\sceo\mid_{\Q(\fXt)}$ and
  $\scct=\scco\mid_{\Q(\fXt)}$.  Then, we have the following.
  \begin{enumerate}
  \item \label{c:faithful-site-restrict-eq}
    We have $\StcefXt=\StcefXo\cap \Q(\fXt)$.
  \item\label{c:faithful-adm-faithful} If there is an
    $O_{1}$-admissible $x\in \Q(\fXo)$, then $O_{1}$ is
    faithful.
  \item \label{c:faithful-existence} Suppose a semiring
    $U\subset \Q$ such that $U=U_{\geq 0}$ and $1\in U$.  If
    $O_{1}$ is either $\{\gefXo^{U},\gfXo^{U}\}$,
    $\{\gefXopm^{U},\gfXopm^{U}\}$, or $\{\geAfXo,\gAfXo\}$,
    then $O_{1}$ is faithful.
  \end{enumerate}
\end{lemma}
\begin{proof}
  Let us prove Claim~\ref{c:faithful-site-restrict-eq}.  If
  $f\in \StcefXo\cap \Q(\fXt)$, then $f\scet f$ by the
  reflexivity of $\sceo$.  Hence, $f\in \StcefXt$. Also, if
  $f\in \StcefXt$, then $f\scet g$ or $g \scet f$ for some
  $g\in \Q(\fXt)$.  Thus, $f\in \StcefXo\cap \Q(\fXt)$.
  
  Let us prove Claim~\ref{c:faithful-adm-faithful}.
  Claim~\ref{c:adm-succ-half-gx-scc} of
  Lemma~\ref{lem:adm-succ} gives $1\scco 0$ by $1\in \Z[x]$.
  This implies $1\scct 0$ by $\{1,0\}\subset
  \Q(\fXt)$. Hence, Claim~\ref{c:faithful-adm-faithful}
  follows from Claim~\ref{c:faithful-site-restrict-eq},
  since $\StcefXt=\StccfXt$ by Lemma~\ref{lem:succ-site}.   
  
  Claim~\ref{c:faithful-existence} follows from
  Claim~\ref{c:faithful-site-restrict-eq}, since each
  $y\in \fXo$ is $O_{1}$-admissible by $1\in U$.
\end{proof}

By Lemma~\ref{lem:faithful}, the squaring orders of parcels are
always faithful, because parcels have admissible
variables. However, in general, we have unfaithful squaring
orders.

\begin{example}\label{rmk:faithful}
  Consider $\fXo=\{X_{1,1},X_{1,2}\}$ and
  $\fXt=\{X_{1,1}\}$.  Let $f \sceo g $ on $\Q(\fX_{1})$,
  if $f,g\in \Z[\fXt]$, $f-g\in \Zgez[\fXt]$, and
  $\ord_{X_{1,1}}(f-g)>0$.  Also, let $f\scco g$, if
  $f \sceo g$ and $f\neq g$.  Then, $O_{1}=\{\sceo,\scco\}$
  is unfaithful  on $\fXo$ as follows.

  Let us confirm that $O_{1}$ has squaring orders on
  $\fXo$.  Thus, we prove that $\StcefXo$ is a
  $\sceo$-poring.  First, $\ord_{X_{1,1}}(0)=\infty>0$ gives
  the reflexivity of $\sceo $.  In particular,
  $\StcefXo=\Z[\fXt]$.  Second, $f\sceo g\sceo h$ implies
  the transitivity $f\sceo h$, since
  $\ord_{X_{1,1}}(f-h)=
  \min(\ord_{X_{1,1}}(f-g),\ord_{X_{1,1}}(g-h))>0$ by
  $f-g,g-h\in\Zgez[\fXt]$.  Third, $f\sceo g\sceo f$ implies
  the antisymmetricity $f=g$ by $f-g,g-f\in \Zgez[\fXt]$.
  Fourth, the additivity and multiplicativity of $\sceo$
  hold, because we have $f+h\sceo g+h$ for $f\sceo g$ and
  $h\in \StcefXo$ by $(f+h)-(g+h)=f-h$, and we have
  $fg\sceo 0$ for $f,g\sceo 0$ by
  $\ord_{X_{1,1}}(fg)=\ord_{X_{1,1}}(f)+\ord_{X_{1,1}}(g)$.
  Hence, $\StcefXo$ is a $\sceo$-poring.
  
  Let us prove that $\StcefXo$ is a strict $\scco$-poring.
  First, $f\scco f$ does not hold by $f=f$.  Second,
  $f \scco g\scco h$ gives the transitivity $f\scco h$,
  because we have $f\neq h$ by $f-g,g-h\in \Zgez[\fXt]$ and
  $f-g,g-h\neq 0$, and we have
  $\ord_{X_{1,1}}(f-h)= \min
  (\ord_{X_{1,1}}(f-g),\ord_{X_{1,1}}(g-h))>0$.  Third, the
  additivity and multiplicativity of $\scco$ hold as above
  for $\sceo$. Thus,
  $\StcefXo$ is a strict $\scco$-poring.

  Let us prove the squaring implications of $O_{1}$ on
  $\fXo$.  First, the $\gAfXo$-$\geAfXo$ and $\scco$-$\sceo$
  implications follow from Definition~\ref{def:ineq} and the
  definition of $O_{1}$. Second, the half $\sceo$-$\geAfXo$
  implication holds, since $f\sceo 0$ implies
  $f \in \Zgez[\fXt]$.  Third, the half $\scco$-$\gAfXo$
  implication holds similarly.
  
  Furthermore, let us prove the semi-strict transitivity of
  $O_{1}$.  If $f\sceo g \scco h$, then
  $f-g,g-h \in \Zgez[\fXt]$ and $g-h\neq 0$.  Hence,
  $f-h\neq 0$. This gives $f\scco h$, because
  $f\sceo g \scco h$ implies $f\sceo g \sceo h$ by the
  $\scco$-$\sceo$ implication, and $f\sceo g \sceo h$ implies
  $f\sceo h$ by the transitivity of $\sceo$.  Similarly, if
  $f\scco g \sceo h$, then $f \scco h$.

  Therefore, $O_{1}$ has squaring orders on $\fXo$,
  because we have the $\sceo$-$\geAfXo$-poring
  inclusion by
  $\StcefXo=\Z[\fXt]$ and the $\sceo$-$\scco$-poring
  equation by Lemma~\ref{lem:succ-site}.
  
  However, for $\fXr=\{X_{1,2}\}$, consider
  $\scer=\sceo|_{\Q(\fXr)}$ and
  $\sccr=\scco|_{\Q(\fXr)}$. Then, $O_{1}$ is not faithful
  by $\emptyset=\StccfXr\neq\StcefXr=\Z$, since
  $\StcefXr=\StcefXo\cap\Q(\fXr)$ by
  Claim~\ref{c:faithful-site-restrict-eq} of
  Lemma~\ref{lem:faithful}, and no $f,g\in \StcefXr$ give
  $f\sccr g$.
\end{example}

The following gives squaring orders as restrictions of
faithful squaring orders.

\begin{lemma}\label{lem:succ-rest}
  Consider faithful squaring orders $O_{1}=\{\sceo,\scco\}$
  on $\fXo$.  For $\fXt=\fX_{1,\lamo,\lamt}$, let
  $\scct=\scco\mid_{\Q(\fXt)}$ and
  $\scet=\sceo\mid_{\Q(\fXt)}$ so that
  $O_{2}=\{\scet,\scct\}$.  Then, we have the following.
  \begin{enumerate}
  \item \label{c:succ-rest-gAfX-geAfX-equiv} Let $f,g\in
    \Q(\fXt)$. Then,
    $f\gAfXo g$ if and only if $f\gAfXt g$.  Also,
    $f\geAfXo g$ if and only if $f\geAfXt g$.
  \item \label{c:succ-rest-existence} We have squaring
    orders $O_{2}$ on $\fXt$ such that $\scct$ is a strict
    squaring order of $\scet$.
  \item \label{c:succ-rest-adm} If $x\in \Q(\fXt)$ is
    $O_{1}$-admissible, then $x$ is $O_{2}$-admissible.
  \end{enumerate}
\end{lemma}
\begin{proof}
  Let $L_{i}=\num(\fX_{i})$ for $i\in\oi(2)$.  Let us prove
  Claim~\ref{c:succ-rest-gAfX-geAfX-equiv}.  First, we
  prove
  the only if part for strict inequalities.  Fix some
  $u\in \R$ such that $0<u<1$. Then, for
  $r=\tp{r_{i}}_{i\in\oi(L_{2})}\in \AfXt$, we put
  ${}_{u,\lamo,\lamt} \sqrt{r}=\tp{\ka_{i}}_{i
    \in\oi(L_{1})}\in \AfXo$ such that
  \begin{align}
    {\ka}_{i}=
    \begin{dcases}
      u &\mif  i\neq \lam_{1,j}
          \mforeach j\in\oi(L_{2}),\\
      \sqrt[\leftroot{0}\uproot{6}\lam_{2,j}]{r_{j}}
        &\mif i=\lam_{1,j}
          \mforsome j\in\oi(L_{2}).
    \end{dcases}
  \end{align}
  Then, $f>_{\AfXo}g$ implies $f(r)>g(r)$ in $\Q(\fXt)$ for
  each $r\in \AfXt$, because $f(r)$ and $g(r)$ are
  $f({}_{u,\lamo,\lamt}\sqrt{r})$ and
  $g({}_{u,\lamo,\lamt}\sqrt{r}) $ in $\Q(\fXo)$.

  Second, let us prove the if part for strict
  inequalities. For $r\in \AfXo$, we put
  $r^{\lamo,\lamt}=
  (r^{\lam_{2,j}}_{\lam_{1,j}})_{j\in\oi(L_{2})}\in\AfXt$.
  Then, $f\gAfXt g$ yields $f(r)>g(r)$ in $\Q(\fXo)$ for
  each $r\in \AfXo$, since $f(r)$ and $g(r)$ are
  $f( r^{\lamo,\lamt})$ and $g( r^{\lamo,\lamt})$ in
  $\Q(\fXt)$.  Therefore, we obtain the equivalence for
  strict inequalities.  The latter equivalence for
  non-strict inequalities holds similarly.
  
  We prove Claim~\ref{c:succ-rest-existence}.  First,
  $\StcefXt$ is a $\scet$-poring and strict $\scct$-poring,
  because $\StcefXo=\StccfXo$ by the $\sceo$-$\scco$-poring
  equality and $\StcefXo\cap \Q(\fXt)=\StcefXt$ by
  Claim~\ref{c:faithful-site-restrict-eq} of
  Lemma~\ref{lem:faithful}.

  Second, let us prove the $\scet$-$\geAfXt$-poring
  inclusion.  Let $f\in \StcefXt$.  This gives $f\scet g$ or
  $g\scet f$ for some $g\in \Q(\fXt)$.  Suppose $f\scet g$
  for simplicity.  Then, since $f\sceo g$, we have
  $f \geAfXo g$ by the half $\sceo$-$\geAfXo$ implication
  and Claim~\ref{c:half-full-implications-sce-geq} of
  Lemma~\ref{lem:half-full-implication}.  Hence, since
  $f\geAfXt g$ by Claim~\ref{c:succ-rest-gAfX-geAfX-equiv},
  the
  $\scet$-$\geAfXt$-poring inclusion holds.
  
  Third, $O_{2}$ satisfies the $\scet$-$\scct$-poring
  equality, because $O_{1}$ is faithful.
  
  Fourth, let us confirm the
  squaring implications of $O_{2}$ on
  $\fXt$.  By Claim~\ref{c:succ-rest-gAfX-geAfX-equiv}, the
  $\gAfXt$-$\geAfXt$ implication follows from the
  $\gAfXo$-$\geAfXo$ implication.  Also, the $\scct$-$\scet$
  implication follows from the $\scco$-$\sceo$ implication.
  Let us confirm the half $\scet$-$\geAfXt$ implication.
  Let $f\scet 0$.  Then, since $f\sceo 0$, we have
  $f\geAfXt 0$ by the half $\sceo$-$\geAfXo$ implication and
  Claim~\ref{c:succ-rest-gAfX-geAfX-equiv}.  Similarly, we
  obtain the half $\scct$-$\gAfXt$ implication.
  
  Fifth, the semi-strict transitivity holds for $O_{2}$,
  as this holds for $O_{1}$.  Therefore, $O_{2}$ consists
  of squaring orders on $\fXt$ such that $\scct$ is a strict
  squaring order of $\scet$.
  
  Let us prove Claim~\ref{c:succ-rest-adm}.  We have the
  half $\gx$-$\scct$ implication, as $O_{2}$ restricts
  $O_{1}$ to $\Q(\fXt)$.  We have the upper condition of $x$
  on $\AfXt$ by $1>_{\AfXo}x$, because for each
  $r\in \AfXt$, $x(r)$ is $x({}_{u,\lamo,\lamt}\sqrt{r})<1$
  in $\Q(\fXo)$.
\end{proof}

We introduce the following notation for
restricted coordinates.
\begin{definition}\label{def:rest-base-numerator}
  Assume $\cF=\Lam(s,l,w,\scco,\fs,\phi,\rho,x,\fXo)$ and
  $\fXt=\fX_{1,\lamo,\lamt}$.
  \begin{enumerate}
  \item We say that $\tp{\cF,\fXt}$ satisfies the base
    condition, if $x\in \Q(\fXt)^{l}$.
  \item We say that $\tp{\cF,\fXt}$ satisfies the numerator
    condition, if $f_{s,m}\in \Q(\fXt)$ for each $m\in \Zl$.
  \end{enumerate}
\end{definition}

Then, the following proposition gives  parcels on restricted
coordinates.

\begin{proposition}\label{prop:rest-parcel}
  Consider $\cF=\Lam(s,l,w,\scco,\fs,\phi,\rho,x,\fXo)$ for
  squaring orders $O_{1}=\{\sceo,\scco\}$ on $\fXo$.
  Let
  $\fXt=\fX_{1,\lamo,\lamt}$. Assume $\tp{\cF,\fXt}$ with
  the base and numerator conditions.  Let
  $O_{2}=\{\scet,\scct\}$ for $\scet=\sceo\mid_{\Q(\fXt)}$ and
  $\scct=\scco\mid_{\Q(\fXt)}$.  Then, we have the following.
  \begin{enumerate}
  \item \label{c:rest-parcel-squaring-order}
    $O_{2}$ has squaring orders on $\fXt$ such that
    $\scct$ is a strict squaring
    order of  $\scet$.
  \item \label{c:rest-parcel-adm}
    $x$ is $O_{2}$-admissible.
  \item \label{c:rest-parcel-pairwise-pos}
    $\fs$ is pairwise $\tp{s,l,\scct}$-positive.
  \item \label{c:rest-parcel-med}
    $\phi$ is a $\tp{s,l,w,\scct,\rho,x,\fXt}$-mediator.
  \item There is a parcel
    $\Lam(s,l,w,\scct,\fs,\phi,\rho,x,\fXt)$.
    \label{c:rest-parcel-existence}
  \end{enumerate}
\end{proposition}
\begin{proof}
  Since $x$ is $O_{1}$-admissible, $O_{1}$ is faithful by
  Claim~\ref{c:faithful-adm-faithful} of
  Lemma~\ref{lem:faithful}.  This gives
  Claim~\ref{c:rest-parcel-squaring-order} by
  Claim~\ref{c:succ-rest-existence} of
  Lemma~\ref{lem:succ-rest}, Claim~\ref{c:rest-parcel-adm}
  follows from the base condition of $\tp{\cF,\fXt}$ and
  Claim~\ref{c:succ-rest-adm} of Lemma~\ref{lem:succ-rest}.
  Claim~\ref{c:rest-parcel-pairwise-pos} holds by the
  numerator condition of $\tp{\cF,\fXt}$.
  
  Let us prove Claim~\ref{c:rest-parcel-med}.  Let
  $\mu_{i}=\tp{s,l,w,\scc_{i},\rho,x,\fX_{i}}$ for
  $i\in\oi(2)$.  First, we confirm the base positivity of
  $\phi$ and $\mut$.  The base condition of
  $\tp{\cF,\fXt}$ gives
  $\phi(x)_{i}\in \Q(\fXt)$.  Also, the base positivity of
  $\phi$ and $\muo$ reads $\phi(x)_{i}^{w_{i}}>_{\AfXo}0$.  Thus,
  Claim~\ref{c:succ-rest-gAfX-geAfX-equiv} of
  Lemma~\ref{lem:succ-rest} yields
  $\phi(x)_{i}^{w_{i}}\gAfXt0$.  Second, the base-shift
  positivity of $\phi$ and $\muo$ implies that of $\phi$ and
  $\mut$, because the base condition of $\tp{\cF,\fXt}$ gives
  $B(s,l,w,m,\phi,\rho,x,\fXo)\in \Q(\fXt)$ by
  Item~\ref{d:bshift-fun-multi-one-tuple} of 
  Definition~\ref{def:bshift-fun-multi}.  In particular, we obtain
  Claim~\ref{c:rest-parcel-med}.
  
  Claim~\ref{c:rest-parcel-existence} follows from
  Claims~\ref{c:rest-parcel-squaring-order},
  \ref{c:rest-parcel-adm}, \ref{c:rest-parcel-pairwise-pos},
  and~\ref{c:rest-parcel-med} and
  Item~\ref{d:merged-parcel}
  of Definition~\ref{def:merged}.
\end{proof}
\begin{remark}\label{rmk:restriction}
  If $\scco=\gx$ and $\sceo=\gex$ in
  Proposition~\ref{prop:rest-parcel}, then by
  $x\in \Q(\fXt)^{l}$, $f \scco g$ is equivalent to
  $f \scct g$, and $f \sceo g$ is equivalent to $f \scet g$.
\end{remark}

By Proposition~\ref{prop:rest-parcel}, we introduce the
following parcels.
\begin{definition}\label{def:rest-parcel}
  Under the assumption of
  Proposition~\ref{prop:rest-parcel}, we define the
  restricted parcel
  $r_{\lamo,\lamt}(\cF) =\Lam(s,l,w,\scct,\fs,\phi,\rho,x,\fXt)$.
\end{definition}

In Definition~\ref{def:rest-parcel}, $\scct$
on $\Q(\fXt)$
depends on
$\scco$ on $\Q(\fXo)$, as $\scct=\scco\mid_{\Q(\fXt)}$.
However, we have the following for some restricted
squaring orders (see Remark~\ref{rmk:restriction} for
$\gex,\gx$).

\begin{proposition}\label{prop:rest-succ-indep}
  Let $\fXo=\{X_{1,i}\}_{i\in\oi(L_{1})}$ and
  $\fXt=\fX_{1,\lamo,\lamt}$.  Then, we have the
  following:
  \begin{align}
    \geAfXo\mid_{\Q(\fXt)}
    =\geAfXt &\mand  \gAfXo \mid_{\Q(\fXt)}
    =\gAfXt;
      \label{eq:rest-succ-indep-geAfX}\\
    \geU_{\fXo}\mid_{\Q(\fXt)}
    =\geU_{\fXt}
    &  \mand     \gU_{\fXo}\mid_{\Q(\fXt)}
      =\gU_{\fXt};
      \label{eq:rest-succ-indep-geUfX}\\
    \geU_{\fXopm}\mid_{\Q(\fXt)} 
    =\geU_{\fXtpm}
    & \mand
      \gU_{\fXopm}\mid_{\Q(\fXt)}
    =\gU_{\fXtpm}.
      \label{eq:rest-succ-indep-geUfXpm}
  \end{align}
\end{proposition}
\begin{proof}
  Lemma~\ref{lem:succ-rest} gives
  equations~\eqref{eq:rest-succ-indep-geAfX}.  Let us
  prove
  equations~\eqref{eq:rest-succ-indep-geUfX}.  First, let
  $f\geU_{\fXo} g$ for $f,g\in \Q(\fXt)$. Then,
  $f,g\in \Q[\fXo]\cap \Q(\fXt)=\Q[\fXt]$.  This also gives
  $f\geU_{\fXt} g$, since $f-g\in U[\fXo]$ implies
  $f-g\in U[\fXt]$.  Second, if $f\geU_{\fXt} g$ for
  $f,g\in \Q[\fXt]$, then $f\geU_{\fXo} g$ by
  $f-g\in U[\fXt]\subset U[\fXo]$.
  Moreover, if $f,g\in \Q(\fXt)$ and $f\neq g$, then
  $f\gU_{\fXo} g$ is the same as $f\gU_{\fXt} g$,
  since
  $f\gU_{\fXo} g$ and $f\gU_{\fXt} g$
  imply
  $f\geU_{\fXo} g$ and $f\geU_{\fXt} g$,
  respectively.
  Thus,
  equations~\eqref{eq:rest-succ-indep-geUfX} follows.  
  Now, equations~\eqref{eq:rest-succ-indep-geUfXpm}
   hold similarly.
\end{proof}

Furthermore, we make the following notion of
{\it optimal
parcel coordinates}. 
\begin{definition}\label{def:optimal} 
  Consider a parcel $\cF=\Lam(s,l,w,\scc,\fs,\phi,\rho,x,\fX)$.
  \begin{enumerate}
  \item Let $\fX=\{X_{i}\}_{i\in\oi(L)}$ and
    $\lamo=\tp{1,2,\dots, L}\in \Zgeo^{L}$.  We call $\fX$
    optimal for $\cF$, if for each $\lamt\in \Zgeo^{L}$ such
    that $\lamt\neq \iota^{L}(1)$, either
    $x\not\in \Q(\fX_{\lamo,\lamt})^{l} $ or
    $f_{s,m}\not\in \Q(\fX_{\lamo,\lamt})$ for some
    $m\in \Zl$.
  \item   We call $\fX$ fully optimal for $\cF$
    (or fully optimal for simplicity), if $\fX$ is optimal
    and $x\in \Q(\fX)^{l}$ is fully admissible.
  \end{enumerate}
\end{definition}

Since $\Q[\fX]$ is a unique factorization domain,
we  use the following notation when our argument
is invariant under the choices of irreducible fractions.

\begin{definition}
  If $f\in \Q(\fX)$, then we put some
  $\Ir(\fX)(f)\in \Q[\fX]^{2}$ such that
  $f=\Qu(\Ir(\fX)(f))$ is an irreducible fraction.
\end{definition}

If $\rho_{1},\rho_{2}\in \Zgeo^{L}$ for
$L\in \Zgeo$, then let
$\lcm(\rho_{1},\rho_{2})=\tp{\lcm(\rho_{1,i},\rho_{2,i})}_{i\in\oi(L)}
\in\Zgeo^{L}$.
We then state the following lemmas to obtain the optimal parcel
coordinate of a parcel. 

\begin{lemma}\label{lem:coord-intersection}
  Assume $\Q(\fX)$ of $\fX=\{X_{i}\}_{i\in\oi(L)}$.
  \begin{enumerate}
  \item \label{c:coord-intersection-subsets} For subsets
    $\fY_{1},\fY_{2}\subset\fX$, let
    $\fY_{3}=\fY_{1}\cap \fY_{2}$. Then,
    $\Q(\fY_{1})\cap\Q(\fY_{2})=\Q(\fY_{3})$.
  \item \label{c:coord-intersection-exponents} For
    $i\in\oi(2)$, suppose $\rho_{i}\in \Zgeo^{L}$ and
    $\fY_{i}=\{X_{j}^{\rho_{i,j}}\}_{j\in\oi(L)}$.  Let
    $\fY_{3}=\{X_{j}^{\lcm(\rho_{1},\rho_{2})_{j}}\}_{j\in\oi(L)}$.
    Then, $\Q(\fY_{1})\cap \Q(\fY_{2})=\Q(\fY_{3})$.
  \item \label{c:coord-intersection-exclusion} Let
    $\rho_{i}\in \Zgeo^{L}$ and
    $\fY_{i}=\{X_{j}^{\rho_{i,j}}\}_{j\in\oi(L)}$ for
    $i\in\oi(2)$. Let $p\in\oi(L)$.  Assume
    $f= \Qu(g)\in \Q(\fY_{1})\cap \Q(\fY_{2})$ such 
    that
    $g\in \Q[\fY_{1}]^{2}$ and
    $\deg_{X_{p}}g_{1},\deg_{X_{p}}g_{2}< \rho_{2,p}$.  Then,
    $ f\in \Q(\fY_{1}\setminus \{X_{p}\})$.
  \end{enumerate}
\end{lemma}
\begin{proof}
  Let us prove Claim~\ref{c:coord-intersection-subsets}.  We
  have $\Q(\fY_{1})\cap \Q(\fY_{2})\supset \Q(\fY_{3})$ by
  $\Q[\fY_{1}]\cap \Q[\fY_{2}]\supset \fY_{3}$.  If
  $f\in \Q(\fY_{1})\cap \Q(\fY_{2})$, then
  $f= \Qu(\Ir(\fY_{1})(f)) =\Qu(\Ir(\fY_{2})(f))$, which
  gives $\tdet(\Ir(\fY_{1})(f),\Ir(\fY_{2})(f)) =0$ in the
  unique factorization domain $\Q[\fX]$.  Hence, each
  irreducible of $\Ir(\fY_{1})(f)_{2}$ has an associate of
  $\Ir(\fY_{2})(f)_{1}$, and vice versa.  Then,
  $f\in \Q(\fY_{3})$, since we have
  $\Ir(\fY_{1})(f)_{2},\Ir(\fY_{2})(f)_{1} \in \Q[\fY_{3}]$.
  Thus, Claim~\ref{c:coord-intersection-subsets} follows.
  Claim~\ref{c:coord-intersection-exponents} holds
  similarly.
  
  Let us prove Claim~\ref{c:coord-intersection-exclusion}.
  If $f=0$, then the assertion is clear. Thus, assume
  $f\neq 0$.  By $f\in \Q(\fY_{2})$,
  $f=\Qu(\Ir(\fY_{2})(f))$.  Also, we have an irreducible
  fraction $f= \Qu(g')$ such that $g'\in \Q[\fY_{1}]^{2}$ and
  $\deg_{X_{p}}g'_{1},\deg_{X_{p}}g'_{2}< \rho_{2,p}$,
  dividing common factors of $g_{1},g_{2}$ if necessarily.
  Then, since $\tdet(\Ir(\fY_{2})(f), g')=0$, each
  irreducible of $\Ir(\fY_{2})(f)_{1}$ has an associate of
  $g'_{1}$. This implies that
  $\deg_{X_{p}}\Ir(\fY_{2})(f)_{1}$ divides
  $\deg_{X_{p}} g'_{1}$. Hence, $\deg_{X_{p}}g'_{1}=0$ by
  the assumption.  Similarly, $\deg_{X_{p}}g'_{2}=0$. Thus,
  Claim~\ref{c:coord-intersection-exclusion} holds.
\end{proof}


\begin{lemma}\label{lem:lcm-base-numerator}
  Let $\cF=\Lam(s,l,w,\scc,\fs,\phi,\rho,x,\fX)$ of
  $\fX=\{X_{i}\}_{i\in\oi(L)}$.  For $i\in\oi(2)$, let
  $\rho_{i}\in \Zgeo^{L}$ and
  $\fY_{i}=\{X_{i}^{\rho_{i,j}}\}_{j\in\oi(L)}$. If both
  $\tp{\cF,\fY_{1}}$ and $\tp{\cF,\fY_{2}}$ satisfy the base and
  numerator conditions, then so does $\tp{\cF,\fY_{3}}$ of
  $\fY_{3}=\{X_{i}^{\lcm(\rho_{1},\rho_{2})_{i}}\}_{i\in\oi(L)}$.
\end{lemma}
\begin{proof}
  The statement holds      by
  Claim~\ref{c:coord-intersection-exponents} of
  Lemma~\ref{lem:coord-intersection}.
\end{proof}


\begin{lemma}\label{lem:rest-coord-smallest-largest}
  Let $\cF=\Lam(s,l,w,\scc,\fs,\phi,\rho,x,\fXo)$ of
  $\fXo=\{X_{1,i}\}_{i\in\oi(L_{1})}$.  Then, we have the
  following.
  \begin{enumerate}
  \item \label{c:rest-coord-smallest-largest-subset}There is the
    smallest non-empty
    subset $\fXt\subset\fXo$ such that $\tp{\cF,\fXt}$ has the base
    and numerator conditions.
  \item \label{c:rest-coord-smallest-largest-exponent}
    Let
    $\fX_{1,\lamo}=\fXt$ in
    Claim~\ref{c:rest-coord-smallest-largest-subset} for
    some $\lamo\in \Zgeo^{L_{2}}$.  Consider the partial
    order $\geq$ on $\Zgeo^{L_{2}}$.  Then, there is the
    largest $\lamt\in \Zgeo^{L_{2}}$ such that
    $\tp{\cF,\fXr}$ of $\fXr=\fX_{1,\lamo,\lamt}$ has the
    base and numerator conditions.
  \end{enumerate}
\end{lemma}
\begin{proof}
  Let us prove
  Claim~\ref{c:rest-coord-smallest-largest-subset}.  Assume
  $\fY_{1},\fY_{2}\subset\fXo$ such that $\tp{\cF,\fY_{1}}$ and
  $\tp{\cF,\fY_{2}}$ satisfy the base and numerator conditions.
  Then, $\tp{\cF, \fY_{1}\cap \fY_{2}}$ has the base and
  numerator conditions by
  Claim~\ref{c:coord-intersection-subsets} of
  Lemma~\ref{lem:coord-intersection}.  Thus,
  Claim~\ref{c:rest-coord-smallest-largest-subset} holds by
  the finiteness of $\fXo$.
  
  Let us prove
  Claim~\ref{c:rest-coord-smallest-largest-exponent}.  If
  there are $u_{1},u_{2}\in \Zgeo^{L_{2}}$ such that
  $\tp{\cF,\fX_{1,\lamo,u_{1}}}$ and
  $\tp{\cF,\fX_{1,\lamo,u_{2}}}$ have the base and numerator
  conditions, then there is
  $u_{3} =\tp{\lcm(u_{1},u_{2})}
  \in \Zgeo^{L_{2}}$ such that
  $\tp{\cF,\fX_{1,\lamo,u_{3}}}$ has the base and numerator
  conditions by Lemma~\ref{lem:lcm-base-numerator} and
  $u_{3}\geq u_{1},u_{2}$.  Hence, let us prove
  Claim~\ref{c:rest-coord-smallest-largest-exponent} by
  contradiction as follows.
  
  Suppose $\mu_{i}\in \Zgez^{L_{2}}$ for $i\in \Zgeo$ and
  $\fY_{i}=\fX_{1,\lamo,\mu_{i}}$ with the following
  three conditions: first, $\muo=\iota^{L_{2}}(1)$; second, $\tp{\cF,\fY_{i}}$
  satisfies the base and numerator conditions for each
  $i\in \Zgeo$; third, some $p\in\oi(L_{2})$ satisfies
  \begin{align}
    \lim_{i\to \infty}\mu_{i,p}=\infty.
    \label{eq:rest-coord-smallest-largest-p1}
  \end{align}
  If there is $g\in \Q(\fX_{1})$ such that
  $g\in \Q(\fY_{i})$ for each $i\in \Zgeo$ and
  $g=\Qu(h)$ for $h\in \Q[\fY_{1}]^{2}$,
  then equation~\eqref{eq:rest-coord-smallest-largest-p1}
  yields $\ka\in \Zgeo$ such that
  $\deg_{X_{1,\lam_{1,p}}} h_{i}< \mu_{\ka,p}$ for
  $i\in\oi(2)$.  This implies
  $g\in \Q(\fY_{1} \setminus \{X_{1,\lam_{1,p}}\})$ by
  Claim~\ref{c:coord-intersection-exclusion} of
  Lemma~\ref{lem:coord-intersection}.  Hence, the base and
  numerator conditions of $\tp{\cF, \fY_{i}}$ give
  $x_{i}, f_{s,m}\in \Q(\fY_{1} \setminus \{X_{1,\lam_{1,p}}\})$ for
  each $i\in\oi(l)$ and $m\ldZl s$.  This contradicts the
  smallest assumption of $\fY_{1}=\fXt$.
\end{proof}

Thus, we obtain the following existence and uniqueness of
optimal parcel coordinates.

\begin{proposition}\label{prop:optimal}
  Let $\cF=\Lam(s,l,w,\scco,\fs,\phi,\rho,x,\fXo)$ of
  $\fXo=\{X_{1,i}\}_{i\in\oi(L_{1})}$.  Then, there is
  $\fXt=\fX_{1,\lamo,\lamt}$ such that $\fXt$ is optimal for
  $r_{\lamo,\lamt}(\cF)=\Lam(s,l,w,\scct,\fs,\phi,\rho,x,\fXt)$.
  In particular, $\fXt$ is uniquely determined by $\cF$ and
  $\fXo$.
\end{proposition}
\begin{proof} 
  Lemma~\ref{lem:rest-coord-smallest-largest} gives the
  existence.  Let us prove the uniqueness.  Let
  $\fXt'=\fX_{1,\lamo',\lamt'}$ for
  $\lamo',\lam'_{2}\in \Zgeo^{L'_{2}}$ such that $\fXt'$ is
  optimal for
  $r_{\lam'_{1},\lam'_{2}}(\cF)=
  \Lam(s,l,w,\scct',\fs,\phi,\rho,x,\fXt')$.
  
  First, let us prove $\lamo=\lamo'$.  Suppose
  $X_{1,\lam_{1,1}}\in \fX_{1,\lamo} \setminus
  \fX_{1,\lamo'}$.  Then, since
  $\Q(\fX_{1,\lamo})\cap \Q(\fX_{1,\lamo'})
  =\Q(\fX_{1,\lamo}\cap \fX_{1,\lamo'})$ by
  Claim~\ref{c:coord-intersection-subsets} of
  Lemma~\ref{lem:coord-intersection},
  $x_{i},f_{s,m}\in \Q(\fX_{1,\lamo}\cap \fX_{1,\lamo'})$
  for each $i\in\oi(l)$ and $m\ldZl s$. Thus, some
  $\ka\in \Z_{>1}$ satisfies
  \begin{align}
    x_{i},f_{s,m}\in \Q\left(
    \{X_{1,\lam_{1,1}}^{\lam_{2,1}\ka}\}\cup\left(
    \fXt\setminus
    \{X_{1,\lam_{1,1}}^{\lam_{2,1}}\} \right)
    \right)
  \end{align}for each $i\in\oi(l)$ and
  $m\ldZl s$. This contradicts
  the optimal property of $\fXt$.
  Second,
  Claim~\ref{c:rest-coord-smallest-largest-exponent}  of
  Lemma~\ref{lem:rest-coord-smallest-largest} gives
  $\lamt=\lamt'$.
\end{proof}

\subsection{
  Merged determinants and
  \texorpdfstring{$q$}{q}-binomial coefficients}
We write merged determinants by 
 $q$-binomial coefficients and base shift functions. This gives
 general non-negativities and positivities of merged
 determinants. 
 We first put the
following notation.

\begin{definition}\label{def:left-right-merged-det}
  Consider a parcel
  $\cF =\Lam(s,l,w,\scc,\fs,\phi,\rho,x,\fX)$ with
  $m,n\in \Z^{l}$, $k\in \Ztl$, $a=\nu(k)$, and
  $b=\nu(m,n,k)$. Let $y=x^{\rho}$. Then, in $\Q(\fX)$, we
  put
  \begin{align}
    &\Delta_{L}(\cF)(s,l,w,m,n,k,\phi,\rho,x,\fX)
    \\&=
    f_{s,m}
    f_{s,n^{\ve}}
    B(s,l,w,m,n^{\ve},\phi,\rho,x,\fX)
    { b \brack a}^{w^{\wcn}}_{y^{\wcn}},\\
    &\Delta_{R}(\cF)(s,l,w,m,n,k,\phi,\rho,x,\fX)
    \\ & =
         f_{s,m\bom k}
         f_{s,(n\bop k)^{\ve}}
         B(s,l,w,m\bom k,(n\bop k)^{\ve},\phi,\rho,x,\fX)
         { b \brack a^{\ve}}^{w^{\wcn}}_{y^{\wcn}}.
  \end{align}
\end{definition}

Also, let us write $\sig$-plus $\bop$ and $\sig$-minus
$\bom$ by $\nu$ and flips.
\begin{lemma}\label{lem:ab-define}
  Let $l\in \Zgeo$, $m,n\in \Zl$, and $k\in \Ztl$ with
  $a=\nu(k)$ and $b=\nu(m,n,k)$.  Then, we have
  \begin{align}
    m\bom k&=(b-a^{\ve})(1,l),
             \label{eq:m_bom_k_ba}\\
    n\bop k&=(b-a^{\ve})(l+1,2l).
             \label{eq:n_bop_k_ba}
  \end{align}
\end{lemma}
\begin{proof}
  First, equation~\eqref{eq:m_bom_k_ba} holds,
  because each
  $i\in\oi(l)$ gives
  $b_{i}-a_{2l-i+1}=\sum k(1,i)+m_{i}- \sum k(1,2l-i+1)
  =m_{i}-  \sum k(i+1,2l-i+1)
  =m_{i}- \sig(k)_{i}$.
  Second, equation~\eqref{eq:n_bop_k_ba} follows, since
  each $i\in\oi(l)$ yields
  $b_{i+l}-a_{l-i+1}
  =      \sum k(1,i+l)+n_{i}   -\sum k(1,l-i+1)
    =n_{i}+  \sum k(l-i+2,l+i)
    =n_{i}+\sig(k)_{l-i+1}$.
\end{proof}

Then,
we obtain the following statement
by $\Delta_{L}$ and $\Delta_{R}$. We use the following
statement to examine the merged-log-concavity 
not only by $q$-binomial
coefficients but also by general non-negativities
and positivities on squaring orders.

\begin{theorem}\label{thm:merged-binom-bshift}
  Suppose a parcel
  $\cF=\Lam(s,l,w,\scc,\fs,\phi,\rho,x,\fX)$. Let
  $\mu=\tp{s,l,m,n,k}$ for $m,n\in \Zl$ and $k\in \Ztl$.
  \begin{enumerate}
  \item \label{c:merged-binom-bshift-left-right-equations}
    We have the following equations:
    \begin{align}
      \Delta_{L}(\cF)(s,l,w,m,n,k,\phi,\rho,x,\fX)
      &=
        \Ups(s,l,w,m,n,k,\phi,\rho,x,\fX)
        \cF_{m} \cF_{n^{\ve}};
        \label{eq:merged-binom-bshift-left}
      \\
      \Delta_{R}(\cF)(s,l,w,m,n,k,\phi,\rho,x,\fX)
      &      =
        \Ups(s,l,w,m,n,k,\phi,\rho,x,\fX)
        \cF_{m\bom k} \cF_{(n\bop k)^{\ve}}.
        \label{eq:merged-binom-bshift-right}
    \end{align}
    In particular, we obtain
    \begin{align} 
      \Delta(\cF)(s,l,w,m,n,k,\phi,\rho,x,\fX)
      &=
        \Delta_{L}(\cF)(s,l,w,m,n,k,\phi,\rho,x,\fX)
      \\&-
      \Delta_{R}(\cF)(s,l,w,m,n,k,\phi,\rho,x,\fX).
      \label{eq:merged-binom-bshift-diff}
    \end{align}
  \item \label{c:merged-binom-bshift-sce} We have the
    following inequalities:
    \begin{align}
      \Delta_{L}(\cF)(s,l,w,m,n,k,\phi,\rho,x,\fX)
      &\sce 0;
        \label{ineq:merged-binom-bshift-sce-left}  \\
      \Delta_{R}(\cF)(s,l,w,m,n,k,\phi,\rho,x,\fX)
      &\sce 0.
        \label{ineq:merged-binom-bshift-sce-right}
    \end{align}
  \item \label{c:merged-binom-bshift-scc-eq} Let $\mu$ be
    fitting.  Then, we obtain
    \begin{align}
      \Delta_{L}(\cF)(s,l,w,m,n,k,\phi,\rho,x,\fX)\scc 0,
      \label{ineq:merged-binom-bshift-scc-left}
    \end{align}
    and  
    \begin{numcases}{\Delta_{R}(\cF)(s,l,w,m,n,k,\phi,\rho,x,\fX)}
      = 0 & \mif \mbox{$\mu$ is unwrapped},
      \label{eq:merged-binom-bshift-eq-right}  \\
      \scc 0 & \mif \mbox{$\mu$ is wrapped.}
      \label{eq:merged-binom-bshift-scc-right}
    \end{numcases}
  \end{enumerate}
\end{theorem}
\begin{proof}
  Let $y=x^{\rho}$, $a=\nu(k)$, and $b=\nu(m,n,k)$.  Let us
  prove
  Claim~\ref{c:merged-binom-bshift-left-right-equations}.
  First, we prove
  equation~\eqref{eq:merged-binom-bshift-left}.  If
  $m\ccn n \nldZtl s$, then $\cF_{m}=f_{s,m}=0$ or
  $\cF_{n^{\ve}}=f_{s,n^{\ve}}=0$.  Also, if
  $a\ccn b\not\geq 0$, then
  ${ b \brack a}^{w^{\wcn}}_{y^{\wcn}}=
  \Ups(s,l,w,m,n,k,\phi,\rho,x,\fX)=0$.  Thus, when
  $m\ccn n\nldZtl s$ or $a\ccn b\not\geq 0$,
  equation~\eqref{eq:merged-binom-bshift-left} follows from
  $0=0$.
     
  Assume $m\ccn n\ldZtl s$ and $a\ccn b\geq 0$.  Then,
  $m=(b-a)(1,l)\geq 0$ gives
  \begin{align}
    &f_{s,m}\cdot B(s,l,w,m,\phi,\rho,x,\fX)
      \cdot
      {b(1,l) \brack a(1,l)}^{w}_{y}
    \\&=
    f_{s,m}\cdot \prod_{i\in\oi(l)}
    \frac{ \phi_{i}(y_{i})^{m_{i}w_{i}}
    [m_{i}]!_{y_{i}}^{w_{i}}}
    { \phi_{i}(x_{i})^{m_{i}w_{i}}      [m_{i}]!_{x_{i}}^{w_{i}}} 
    \cdot \prod_{i\in\oi(l)}
    \frac{[b_{i}]!_{y_{i}}^{w_{i}}
    }{[a_{i}]!_{y_{i}}^{w_{i}} [m_{i}]!_{y_{i}}^{w_{i}}}
    \\&= f_{s,m} \cdot
    \prod_{i\in\oi(l)}
    \frac{\phi_{i}(y_{i})^{m_{i}w_{i}} [b_{i}]!_{y_{i}}^{w_{i}}
    }{\phi_{i}(x_{i})^{m_{i}w_{i}}
    [a_{i}]!_{y_{i}}^{w_{i}} [m_{i}]!_{x_{i}}^{w_{i}}}
    \\& =
    \frac{f_{s,m}}{\prod \phi(x)^{m\rc w}\cdot
    [m]!_{x}^{w}}
    \cdot
    \prod \phi(y)^{(b-a)(1,l)\rc w}\cdot 
    \frac{[b(1,l)]!_{y}^{w}}{[a(1,l)]!_{y}^{w}}
    \\&     =\cF_{s,m}\cdot
    \prod \phi(y)^{(b-a)(1,l)\rc w}\cdot 
    \frac{[b(1,l)]!_{y}^{w}}{[a(1,l)]!_{y}^{w}}.
    \label{eq:merged-binom-bshift-bshift-parcel}
  \end{align}
  Also, $n=(b-a)(l+1,2l) \geq 0$ yields
  \begin{align}
    &f_{s,n^{\ve}}\cdot B(s,l,w,n^{\ve},\phi,\rho,x,\fX)
      \cdot
      {b(l+1,2l) \brack a(l+1,2l)}^{w^{\ve}}_{y^{\ve}}
    \\&    =
    f_{s,n^{\ve}}
    \cdot
    \prod_{i\in\oi(l)}
    \frac{ \phi_{l-i+1}(y_{l-i+1})^{n_{i}w_{l-i+1}}
    [n_{i}]!_{y_{l-i+1}}^{w_{l-i+1}}}
    { \phi_{l-i+1}(x_{l-i+1})^{n_{i}w_{l-i+1}}
    [n_{i}]!_{x_{l-i+1}}^{w_{l-i+1}}}
    \cdot \prod_{i\in\oi(l)}
    \frac{[b_{i+l}]!_{y_{l-i+1}}^{w_{l-i+1}}}
    {[a_{i+l}]!_{y_{l-i+1}}^{w_{l-i+1}} [n_{i}]!_{y_{l-i+1}}^{w_{l-i+1}}}
    \\&
    =
    f_{s,n^{\ve}}
    \cdot
    \prod_{i\in\oi(l)}
    \frac{
    \phi_{l-i+1}(y_{l-i+1})^{n_{i}w_{l-i+1}}
    [b_{i+l}]!_{y_{l-i+1}}^{w_{l-i+1}}}
    { \phi_{l-i+1}(x_{l-i+1})^{n_{i}w_{l-i+1}}
    [a_{i+l}]!_{y_{l-i+1}}^{w_{l-i+1}}
    [n_{i}]!_{x_{l-i+1}}^{w_{l-i+1}}}
    \\&       =
    \frac{f_{s,n^{\ve}}
    }{\prod \phi(x)^{n^{\ve}\rc w}\cdot  [n^{\ve}]!_{x}^{w}}
    \cdot
    \prod (\phi(y)^{\ve})^{(b-a)(l+1,2l)\rc w^{\ve}}
    \cdot
    \frac{[b(l+1,2l)]!_{y^{\ve}}^{w^{\ve}}
    }{[a(l+1,2l)]!_{y^{\ve}}^{w^{\ve}}}
    \\&=
    \cF_{s,n^{\ve}}\cdot
    \prod (\phi(y)^{\ve})^{(b-a)(l+1,2l)\rc w^{\ve}}
    \cdot
    \frac{[b(l+1,2l)]!_{y^{\ve}}^{w^{\ve}}
    }{[a(l+1,2l)]!_{y^{\ve}}^{w^{\ve}}}.
    \label{eq:merged-binom-bshift-bshift-parcel-ve}
  \end{align}
  Furthermore, we have
  \begin{align}
    \frac{[b(1,l)]!_{y}^{w}}{[a(1,l)]!_{y}^{w}}\cdot
    \frac{[b(l+1,2l)]!_{y^{\ve}}^{w^{\ve}}
    }{[a(l+1,2l)]!_{y^{\ve}}^{w^{\ve}}}
    &=   \frac{[b]!_{y^{\wcn}}^{w^{\wcn}}}
      {[a]!_{y^{\wcn}}^{w^{\wcn}}},\\
    \prod \phi(y)^{(b-a)(1,l)\rc w}\cdot
    \prod (\phi(y)^{\ve})^{(b-a)(l+1,2l)\rc w^{\ve}}
    &=\prod (\phi(y)^{\wcn})^{(b-a)\rc w^{\wcn}}.
  \end{align}
  Thus,
  equations~\eqref{eq:merged-binom-bshift-bshift-parcel}
  and~\eqref{eq:merged-binom-bshift-bshift-parcel-ve} give
  equation~\eqref{eq:merged-binom-bshift-left}.

  Second, we prove
  equation~\eqref{eq:merged-binom-bshift-right}.  If
  $(m\bom k)\ccn (n\bop k)\nldZtl s$, then
  $\cF_{s,m\bom k}=f_{s,m\bom k}=0$ or
  $\cF_{s,(n\bop k)^{\ve}}=f_{s,(n\bop k)^{\ve}}=0$.  
  If $a\ccn b\not\geq 0$, then 
  ${ b \brack a^{\ve}}^{w^{\wcn}}_{y^{\wcn}}=
  \Ups(s,l,w,m,n,k,\phi,\rho,x,\fX)=0$.  Thus, when
  $(m\bom k)\ccn (n\bop k)\nldZtl s$ or $a\ccn b\not\geq 0$,
  equation~\eqref{eq:merged-binom-bshift-right} holds by
  $0=0$.

  Hence,  assume a wrapped $\mu$ with $a\ccn b \geq 0$.
  Since Lemma~\ref{lem:ab-define} says
  $(b-a^{\ve})(1,l)=m\bom k$, which is non-negative by
  $m\bom k \ldZl s$, we obtain
  \begin{align}
    &f_{s,m\bom k}
      \cdot
      B(s,l,w,m\bom k,\phi,\rho,x,\fX)
      \cdot
      {b(1,l) \brack a^{\ve}(1,l)}^{w}_{y}
    \\&=
    f_{s,m\bom k}
    \cdot
    \prod_{i\in\oi(l)}
    \frac{ \phi_{i}(y_{i})^{(m\bom k)_{i}w_{i}}
    [(m\bom k)_{i}]!_{y_{i}}^{w_{i}}}
    { \phi_{i}(x_{i})^{(m\bom k)_{i}w_{i}} 
    [(m\bom k)_{i}]!_{x_{i}}^{w_{i}}} 
    \cdot \prod_{i\in\oi(l)}
    \frac{[b_{i}]!_{y_{i}}^{w_{i}}}
    {[a_{2l-i+1}]!_{y_{i}}^{w_{i}}
    [(m\bom k)_{i}]!_{y_{i}}^{w_{i}}}
    \\&    = f_{s,m\bom k} \cdot
    \prod_{i\in\oi(l)}
    \frac{\phi_{i}(y_{i})^{(m\bom k)_{i}w_{i}}
    [b_{i}]!_{y_{i}}^{w_{i}}}{
    \phi_{i}(x_{i})^{(m\bom k)_{i}w_{i}}
    [a_{2l-i+1}]!_{y_{i}}^{w_{i}}
    [(m\bom k)_{i}]!_{x_{i}}^{w_{i}}}
    \\&     =
    \frac{f_{s,m\bom k}}{\prod \phi(x)^{(m\bom k) \rc w}
    \cdot [m\bom k]!_{x}^{w}}
    \cdot
    \prod \phi(y)^{(b-a^{\ve})(1,l)\rc w}\cdot 
    \frac{[b(1,l)]!_{y}^{w}}{[a^{\ve}(1,l)]!_{y}^{w}}
    \\ & =\cF_{s,m\bom k}\cdot
         \prod \phi(y)^{(b-a^{\ve})(1,l)\rc w}
         \cdot
         \frac{[b(1,l)]!_{y}^{w}}{[a^{\ve}(1,l)]!_{y}^{w}}.
         \label{eq:merged-binom-bshift-bshift-parcel-bom}
  \end{align}
  Also, Lemma~\ref{lem:ab-define} says
  $(b-a^{\ve})(l+1,2l)=n\bop k$, which is non-negative by
  $n\bop k\ldZl s$. Hence, we obtain
  \begin{align}
    &f_{s,(n\bop k)^{\ve}}\cdot B(s,l,w,(n\bop k)^{\ve},\phi,\rho,x,\fX)
      \cdot
      {b(l+1,2l) \brack a^{\ve}(l+1,2l)}^{w^{\ve}}_{y^{\ve}}
    \\&
    =
    f_{s,(n\bop k)^{\ve}}
    \cdot
    \prod_{i\in\oi(l)}
    \frac{ \phi_{l-i+1}(y_{l-i+1})^{(n\bop k)_{i}w_{l-i+1}}
    [(n\bop k)_{i}]!_{y_{l-i+1}}^{w_{l-i+1}}}
    { \phi_{l-i+1}(x_{l-i+1})^{(n\bop k)_{i}w_{l-i+1}}
    [(n\bop k)_{i}]!_{x_{l-i+1}}^{w_{l-i+1}}}
    \\ &\cdot       \prod_{i\in\oi(l)}
         \frac{[b_{i+l}]!_{y_{l-i+1}}^{w_{l-i+1}}}
         {[a_{l-i+1}]!_{y_{l-i+1}}^{w_{l-i+1}} [(n\bop
         k)_{i}]!_{y_{l-i+1}}^{w_{l-i+1}}}
    \\&       =
    f_{s,(n\bop k)^{\ve}}
    \cdot
    \prod_{i\in\oi(l)}
    \frac{ \phi_{l-i+1}(y_{l-i+1})^{(n\bop k)_{i}w_{l-i+1}}
    [b_{i+l}]!_{y_{l-i+1}}^{w_{l-i+1}}}
    { \phi_{l-i+1}(x_{l-i+1})^{(n\bop k)_{i}w_{l-i+1}}
    [a_{l-i+1}]!_{y_{l-i+1}}^{w_{l-i+1}}
    [(n\bop k)_{i}]!_{x_{l-i+1}}^{w_{l-i+1}}}
    \\&       =
    \frac{ f_{s,(n\bop k)^{\ve}}}{   \prod \phi(x)^{(n\bop k)^{\ve}\rc w}\cdot
    [(n\bop k)^{\ve}]!_{x}^{w}}
    \cdot
    \prod (\phi(y)^{\ve})^{(b-a^{\ve})(l+1,2l)\rc w^{\ve}}
    \cdot
    \frac{[b(l+1,2l)]!_{y^{\ve}}^{w^{\ve}}}
    {[a^{\ve}(l+1,2l)]!_{y^{\ve}}^{w^{\ve}}}
    \\&=
    \cF_{s,(n\bop k)^{\ve}}\cdot \prod
    (\phi(y)^{\ve})^{(b-a^{\ve})(l+1,2l)\rc w^{\ve}}
    \cdot
    \frac{[b(l+1,2l)]!_{y^{\ve}}^{w^{\ve}}
    }{[a^{\ve}(l+1,2l)]!_{y^{\ve}}^{w^{\ve}}}.
    \label{eq:merged-binom-bshift-bshift-parcel-bop}
  \end{align}
  Moreover, we have
  \begin{align}
    \frac{[b(1,l)]!_{y}^{w}}{[a^{\ve}(1,l)]!_{y}^{w}}
    \cdot
    \frac{[b(l+1,2l)]!_{y^{\ve}}^{w^{\ve}}
    }{[a^{\ve}(l+1,2l)]!_{y^{\ve}}^{w^{\ve}}}
    =   \frac{[b]!_{y^{\wcn}}^{w^{\wcn}}}{
    [a^{\ve}]!_{y^{\wcn}}^{w^{\wcn}}}
    =\frac{[b]!_{y^{\wcn}}^{w^{\wcn}}}{
    [a]!_{y^{\wcn}}^{w^{\wcn}}},
    \label{eq:merged-binom-bshift-qfactorial}
  \end{align} 
  since
  $[a^{\ve}]!_{y^{\wcn}}^{w^{\wcn}}
  =\prod_{i\in\oi(l)}
      [a_{2l-i+1}]_{y_{i}}^{w_{i}}\cdot
      \prod_{i\in\oi(l)}
      [a_{l-i+1}]_{y_{l-i+1}}^{w_{l-i+1}}
    =
    \prod_{i\in\oi(l)}
      [a_{l+i}]_{y_{l-i+1}}^{w_{l-i+1}}\cdot
      \prod_{i\in\oi(l)}
        [a_{i}]_{y_{i}}^{w_{i}}
    =[a]!_{y^{\wcn}}^{w^{\wcn}}$.
  Also, it holds that 
  \begin{align}
    \prod \phi(y)^{(b-a^{\ve})(1,l)\rc w}
    \cdot\prod (\phi(y)^{\ve})^{(b-a^{\ve})(l+1,2l)\rc w^{\ve}}
    = 
    \prod (\phi(y)^{\wcn})^{(b-a)\rc w^{\wcn}}.
    \label{eq:merged-binom-bshift-med-prod}
  \end{align}
  Thus,
  equations~\eqref{eq:merged-binom-bshift-bshift-parcel-bom},
  \eqref{eq:merged-binom-bshift-bshift-parcel-bop},
  \eqref{eq:merged-binom-bshift-qfactorial},
  and~\eqref{eq:merged-binom-bshift-med-prod} yield
  equation~\eqref{eq:merged-binom-bshift-right}.

  Third,  
  equations~\eqref{eq:merged-binom-bshift-left}
  and~\eqref{eq:merged-binom-bshift-right} give
  equation~\eqref{eq:merged-binom-bshift-diff}.
  Hence,
  Claim~\ref{c:merged-binom-bshift-left-right-equations}
  holds.

  Let us prove Claim~\ref{c:merged-binom-bshift-sce}.
  First, we prove
  inequality~\eqref{ineq:merged-binom-bshift-sce-left}.
  Thus, let us confirm the $\sce$-non-negativities of
  $f_{s,m}f_{s,n}$, $B(s,l,w,m,n,\phi,\rho,x,\fX)$, and
  ${b \brack a}^{w^{\wcn}}_{y^{\wcn}}$.

  We have $f_{s,m}=0$ for $m\nldZl s$.  Hence, the
  $\scc$-$\sce$ implication implies
  \begin{align}
    f_{s,m}f_{s,n}\sce 0
    \label{ineq:merged-binom-bshift-pairwise-nonneg}
  \end{align}
  for each $m,n\in \Zl$, because the pairwise positive
  property of $\fs$ gives
  \begin{align}
    f_{s,m}f_{s,n}\scc 0
    \label{ineq:merged-binom-bshift-pairwise-pos}
  \end{align} 
  for each $m,n\ldZl s$.  Similarly, we have
  \begin{align}
    B(s,l,w,m,n,\phi,\rho,x,\fX)\sce 0
    \label{ineq:merged-binom-bshift-nonneg}
  \end{align}
  for each $m,n\in \Zl$, because the base-shift positivity of
  $\phi$ and $\tp{s,l,w,\scc,\rho,x,\fX}$ says
  \begin{align}
    B(s,l,w,m,n,\phi,\rho,x,\fX)\scc 0
    \label{ineq:merged-binom-bshift-pos}
  \end{align}
  for each $m,n\ldZl s$.
 
  Let us confirm
  \begin{align}
    {b \brack a}^{w^{\wcn}}_{y^{\wcn}}\sce 0.
    \label{eq:merged-binom-bshift-gauss-noneg}
  \end{align}
  We have
  ${b \brack a}^{w^{\wcn}}_{y^{\wcn}} =\prod_{i\in\oi(l)}{b_{i}
    \brack a_{i}}^{w_{i}}_{x_{i}^{\rho_{i}}} \cdot{b_{l+i} \brack
    a_{l+i}}^{w_{l-i+1}}_{x_{l-i+1}^{\rho_{l-i+1}}}$.  If
  $b-a\not\geq 0$, $a\not\geq 0$, or $b\not\geq 0$, then,
  ${b \brack a}^{w^{\wcn}}_{y^{\wcn}}=0$.  Hence, suppose
  $b-a\geq 0$, $a\geq 0$, and $b\geq 0$.  Then,
  ${b_{i} \brack a_{i}}^{w_{i}}_{x_{i}^{\rho_{i}}}>_{x_{i}}0$
  and
  ${b_{l+i} \brack
    a_{l+i}}^{w_{l-i+1}}_{x_{l-i+1}^{\rho_{l-i+1}}}>_{x_{l-i+1}}0$
  for each $i\in\oi(l)$.  Moreover, since each $x_{i}$ is
  $\scc$-admissible, the multiplicativity of $\scc$ yields
  \begin{align}  
    {b \brack a}^{w^{\wcn}}_{y^{\wcn}}\scc 0.
    \label{ineq:merged-binom-bshift-gauss-pos}
  \end{align}
  In particular,
  inequality~\eqref{eq:merged-binom-bshift-gauss-noneg}
  follows.
     
  Thus, inequality~\eqref{ineq:merged-binom-bshift-sce-left}
  holds by
  inequalities~\eqref{ineq:merged-binom-bshift-pairwise-nonneg},
  \eqref{ineq:merged-binom-bshift-nonneg},
  and~\eqref{eq:merged-binom-bshift-gauss-noneg}.
  Similarly,
  inequality~\eqref{ineq:merged-binom-bshift-sce-right}
  holds.  Hence, we obtain
  Claim~\ref{c:merged-binom-bshift-sce}.
       
  Let us prove Claim~\ref{c:merged-binom-bshift-scc-eq}.  We
  obtain
  inequality~\eqref{ineq:merged-binom-bshift-scc-left} by
  inequalities~\eqref{ineq:merged-binom-bshift-pairwise-pos},
  \eqref{ineq:merged-binom-bshift-pos},
  and~\eqref{ineq:merged-binom-bshift-gauss-pos} by the
  inclusion condition of $\mu$ and
  Claim~\ref{c:fitting-nonneg-ab-nonneg} of
  Lemma~\ref{lem:fitting-nonneg}.  Similarly, we have
  inequality~\eqref{eq:merged-binom-bshift-scc-right}, since
  $(m\bom k)\ccn (n\bop k)\ldZtl s$ gives $b-a^{\ve}\geq 0$ by
  Lemma~\ref{lem:ab-define}.  Also,
  equation~\eqref{eq:merged-binom-bshift-eq-right} holds,
  since $ (m\bom k)\ccn (n\bop k) \nldZtl s$ gives
  $f_{s,m\bom k}=0$ or $ f_{s,(n\bop k)^{\ve}}=0$.
\end{proof}

 \subsection{Cut and shift operators}
 \label{sec:cut-shift}
 We introduce some
 general operators to give a parcel from another
 one.
  \subsubsection{Cut operators}
  \begin{definition}\label{def:cut}
   Suppose gates $s_{1},s_{2}\geq 0$ such that
   $s_{1,1}\leq s_{2,1}\leq s_{2,2}\leq s_{1,2}$. Let
   $\cF_{1}=\Lam(s_{1},l,w,\scc,f_{1,s_{1}},\phi,\rho,x,\fX)$.
   Then, we define the parcel
   $\cF_{2} =\ct_{s_{1},s_{2}}(\cF_{1})
   =\Lam(s_{2},l,w,\scc,f_{2,s_{2}},\phi,\rho,x,\fX)$ such that
   $f_{2,s_{2},m}=
   f_{1,s_{1},m}$ if $m\ldZl s_{2}$ and
   $f_{2,s_{2},m}=0$ otherwise.
    We call $\ct_{s_{1},s_{2}}$ cut operator.
  \end{definition}

  Then, we have  the following merged-log-concavity
  on cut operators.
 
 \begin{proposition}\label{prop:cut-op}
   Let
   $\cF_{i}=\Lam(s_{i},l,w,\scc,f_{i,s_{i}},\phi,\rho,x,\fX)$ for
   $i\in\oi(2)$ such that
   $\cF_{2}=\ct_{s_{1},s_{2}}(\cF_{1})$.  Suppose a fitting
   $\mu=\tp{s_{2},l,m,n,k}$.
   \begin{enumerate}
   \item \label{c:cut-op-wrapped}
     If $\mu$ is wrapped,
     then
   \begin{align}
     \Delta(\cF_{2})
     (s_{2},l,w,m,n,k,\phi,\rho,x,\fX)
     =\Delta(\cF_{1})
     (s_{1},l,w,m,n,k,\phi,\rho,x,\fX).
   \end{align}
 \item \label{c:cut-op-unwrapped} If $\mu$ is unwrapped,
   then
   \begin{align}
     \Delta(\cF_{2})(s_{2},l,w,m,n,k,\phi,\rho,x,\fX)
     =
     \Delta_{L}(\cF_{1})(s_{1},l,w,m,n,k,\phi,\rho,x,\fX).
   \end{align}

 \item \label{c:cut-op-merged} $\cF_{2}$ is
   $\scc'$-merged-log-concave, if $\cF_{1}$ is
   $\scc'$-merged-log-concave.
   \end{enumerate}
 \end{proposition}
 \begin{proof}
   Let us prove Claim~\ref{c:cut-op-wrapped}. We have
   $s_{1,1}\leq s_{2,1}\leq s_{2,2}\leq s_{1,2}$ and
   $(m\bom k)\ccn (n\bop k)\ldZtl s_{2}$. Thus,
   $f_{2,s_{2},m} f_{2,s_{2},n^{\ve}}
   B(s_{2},l,w,m,n^{\ve},\phi,\rho,x,\fX)
            = f_{1,s_{1},m}
            f_{1,s_{1},n^{\ve}}
     B(s_{1},l,w,m,n^{\ve},\phi,\rho,x,\fX)$ and
     $f_{2,s_{2},m\bom k} f_{2,s_{2},(n\bop k)^{\ve}} 
       B(s_{2},l,w,m\bom k,(n\bop k)^{\ve},\phi,\rho,x,\fX)
       = 
           f_{1,s_{1},m\bom k}
           f_{1,s_{1},(n\bop k)^{\ve}} 
           B(s_{1},l,w,m\bom k,(n\bop k)^{\ve},\phi,\rho,x,\fX)$.
      Claim~\ref{c:cut-op-wrapped} now follows from
     Claim~\ref{c:merged-binom-bshift-left-right-equations}
     of Theorem~\ref{thm:merged-binom-bshift}.
   
     Let us prove Claim~\ref{c:cut-op-unwrapped}.
     Claim~\ref{c:merged-binom-bshift-scc-eq} of
     Theorem~\ref{thm:merged-binom-bshift} yields
     $\Delta(\cF_{2}) (s_{2},l,w,m,n,k,\phi,\rho,x,\fX)
     =\Delta_{L}(\cF_{2}) (s_{2},l,w,m,n,k,\phi,\rho,x,\fX)$.
     This equals to
     $\Delta_{L}(\cF_{2}) (s_{1},l,w,m,n,k,\phi,\rho,x,\fX)$ by
     $m\ccn n\ldZtl s_{1}$.

   Let us prove Claim~\ref{c:cut-op-merged}.
   Thus, we want to confirm
   \begin{align}
     \Delta(\cF_{2})(s_{2},l,w,m,n,k,\phi,\rho,x,\fX)
     &\scc' 0.
       \label{ineq:cut-op-wrapped-s2-scc}
   \end{align}
   Because $\mu$ is fitting, $\tp{s_{1},l,m,n,k}$ is fitting. In
   particular,
   \begin{align}
     \Delta(\cF_{1})
     (s_{1},l,w,m,n,k,\phi,\rho,x,\fX)\scc'0.
     \label{ineq:cut-op-wrapped-s1-scc}
   \end{align}
   Hence, if $\mu$ is wrapped, then
   Claim~\ref{c:cut-op-wrapped} and
   inequality~\eqref{ineq:cut-op-wrapped-s1-scc} imply
   inequality~\eqref{ineq:cut-op-wrapped-s2-scc}.  Moreover,
   $\scc'\Sup \scc$, since $\cFo$ is
   $\scc'$-merged-log-concave.  Thus, if $\mu$ is unwrapped,
   then Claim~\ref{c:cut-op-unwrapped} gives
   inequality~\eqref{ineq:cut-op-wrapped-s2-scc}, because
   $m\ccn n\ldZtl s_{2}$ gives
   $f_{2,s_{2},m} f_{2,s_{2},n^{\ve}} =
     f_{1,s_{1},m} f_{1,s_{1},n^{\ve}} \scc 0$ and
     $B(s_{2},l,w,m,n^{\ve},\phi,\rho,x,\fX)=
     B(s_{1},l,w,m,n^{\ve},\phi,\rho,x,\fX)\scc 0$.
 \end{proof}

 \subsubsection{Shift operators}

 \begin{definition}\label{def:shift}
   Suppose
   $\cF_{1} =\Lam(s_{1},l,w,\scc,f_{1,s_{1}},\phi,\rho,x,\fX)$
   with $h\in \Zgez$ and $s_{2}=s_{1}+h$.  We define the
   parcel
   $\sh_{h}(\cF_{1})
   =\Lam(s_{2},l,w,\scc,f_{2,s_{2}},\phi,\rho,x,\fX)$
   such that
   \begin{align}
     f_{2,s_{2},m}
     &=
       \begin{dcases}
         f_{1,s_{1},m-h}\cdot
         \frac{[m]!_{x}^{w}}{[m-h]!_{x}^{w}}
         \mfor m\ldZl s_{2},\\
         0 \melse.
       \end{dcases}
   \end{align}
   We call $\sh_{h}$ shift operator.  Also, let
   $ f_{2,s_{2}}=\sh_{h}(f_{1,s_{1}})$.
 \end{definition}
 Notice that $S_{h}(\cF_{1})$ is a parcel, because the
 $\scc$-admissibility of $x$ and the pairwise positivity of
 $f_{1,s_{1}}$ imply that each $m,n\ldZl s_{1}+h$
 satisfies
 \begin{align}
   f_{2,s_{2},m}f_{2,s_{2},n}=
   f_{1,s_{1},m-h}f_{1,s_{1},n-h} \cdot
   \frac{[m]!_{x}^{w}}{[m-h]!_{x}^{w}}
   \frac{[n]!_{x}^{w}}{[n-h]!_{x}^{w}}
   \scc 0.
   \label{ineq:shift-well}
 \end{align}

 We do not extend Definition~\ref{def:shift} for $h<0$,
 since inequality~\eqref{ineq:shift-well} does not have to
 hold for $h<0$ and $w\neq \iota^{l}(0)$.

 For example, each $m\ldZl s_{1}+1$ satisfies
 \begin{align}
   \sh_{1}(\cF_{1})_{m}
 =\frac{f_{1,s_{1},m-1}}{\prod \phi(x)^{m\rc w}\cdot [m-1]!_{x}^{w}}.
 \end{align}
 In particular, let $l=1$, $x=\tp{q}$, and $w=\tp{1}$.  Also,
 suppose $f_{1,s_{1},m}=q^{m_{1}}$ for $m\ldZl s_{1}$.
 Then, each $m\ldZl s_{1}+1$ gives
 $\sh_{1}(f_{1,s_{1}})_{m}=[m_{1}]_{q}q^{m_{1}-1}$, which
 coincides with the Jackson derivative of $f_{1,s_{1},m}$.

 We first
 prove $\sh_{h+1}=\sh_{1}\rc \sh_{h}$
 to discuss the merged-log-concavity of
  shift operators.
 \begin{lemma}\label{lem:one-shift}
   Assume
   $\cFo=\Lam(s_{1},l,w,\scc,f_{1,s_{1}},\phi,\rho,x,\fX)$ and
   $h\in \Zgez$.  Then,
   $\sh_{h+1}(\cFo)=\sh_{1}(\sh_{h}(\cFo))$.
 \end{lemma}
 \begin{proof}
   Let  $s_{2}=s_{1}+h$ and
   $s_{3}=s_{4}=s_{2}+1$.
   Consider the following parcels:
   \begin{align}
     \cFt
     &
       =\sh_{h}(\cFo)
       =\Lam(s_{2},l,w,\scc,f_{2,s_{2}},\phi,\rho,x,\fX);\\
     \cFr
     &=\sh_{1}(\sh_{h}(\cFo))
       =\Lam(s_{3},l,w,\scc,f_{3,s_{3}},\phi,\rho,x,\fX);\\
     \cF_{4}
     &=\sh_{h+1}(\cFo)=
       \Lam(s_{4},l,w,\scc,f_{4,s_{4}},\phi,\rho,x,\fX).
   \end{align}
   Let $m\ldZl s_{4}$.  Then, the assertion follows from
   $f_{3,s_{3},m}
     =f_{2,s_{2},m-1}\cdot \frac{[m]!_{x}^{w}
       }{[m-1]!_{x}^{w}}
       =f_{1,s_{1},m-1-h}\cdot \frac{[m-1]!_{x}^{w}
       }{[m-1-h]!_{x}^{w}}
       \cdot \frac{[m]!_{x}^{w}
       }{[m-1]!_{x}^{w}}
       =f_{1,s_{1},m-h-1}\cdot
       \frac{[m]!_{x}^{w}}{[m-h-1]!_{x}^{w}}
       =f_{4,s_{4},m}$.
  \end{proof}

  We then state the following compatibility of the fitting
  condition and shift operators.
  \begin{lemma}\label{lem:fit-shift}
    Let $h\in \Zgez$ and $s_{2}=s_{1}+h$.  For $m,n\in \Zl$
     and $k\in \Ztl$, let $\muo=\tp{s_{1},l,m-h,n-h,k}$ and
    $\mut=\tp{s_{2},l,m,n,k}$.  Then, we have the following.
    \begin{enumerate}
    \item\label{c:fit-shift-fitting} $\mut$ is fitting if
      and only if $\muo$ is fitting.
    \item\label{c:fit-shift-wrapped} $\mut$ is wrapped if
      and only if $\muo$ is wrapped.
    \end{enumerate}
 \end{lemma}
 \begin{proof}
   Let $a_{1}=a_{2}=\nu(k)$, $b_{1}=\nu(m-h,n-h,k)$, and
   $b_{2}=\nu(m,n,k)$.  Let us prove
   Claim~\ref{c:fit-shift-fitting}.  Assume that $\mut$ is
   fitting.  Then, the inclusion condition of $\muo$ holds
   as $(m-h)\ccn (n-h)\ldZtl s_{2}-h$, which is $s_{1}$.
   Also, the slope conditions of $\muo$
   follow from $a_{1}=a_{2}$ and $b_{1}=b_{2}-h$.  Thus,
   $\muo$ is fitting.  The converse holds similarly.

   We obtain Claim~\ref{c:fit-shift-wrapped} by
   $(m\bom k)\ccn (n\bop k)-h=((m-h)\bom k)\ccn
   ((n-h)\bop
   k)$, since $(m\bom k)\ccn (n\bop k)\ldZtl s_{2}$ is
   equivalent to $(m\bom k)\ccn (n\bop k)-h\ldZtl s_{2}-h$.
 \end{proof}

 Moreover, we confirm the following equations by
 $q$-binomial coefficients.
 
 \begin{lemma}\label{lem:med-bshift-binom}
   Let $l\in \Zgeo$, $w\in \Zlgez$, $\rho\in \Zgeo^{l}$, and
   $m,n\in \Zgeo^{l}$.  Consider $a,b\in \Zgez^{2l}$ such
   that
   $b-a=m\ccn n$. Let $x\in \Q(\fX)^{l}$ an indeterminate
   with $y=x^{\rho}$. Suppose
   $\phi(x)\in \prod_{i\in\oi(l)}\Q(x_{i})$ such that
   $\prod \phi(x)\neq 0$.  Then, we have
   \begin{dmath}
     \frac{
       \prod
       \phi(y)^{m\rc w}\cdot
       [m]!_{y}^{w}}
     { \prod \phi(x)^{m\rc w}  \cdot [m]!_{x}^{w}}
     \cdot
     [m]_{x}^{w}
     \cdot{b(1,l) \brack a(1,l)}_{y}^{w}
     =  [b(1,l)]_{y}^{w}\cdot
     \frac{ \prod \phi(y)^{(m-1)\rc w}\cdot
       [m-1]!_{y}^{w}}{ \prod \phi(x)^{(m-1)\rc w}\cdot
       [m-1]!_{x}^{w}}
     \cdot
     \frac{ \prod \phi(y)^{w}}{\prod \phi(x)^{w}}
     \cdot
     {b(1,l)-1 \brack a(1,l)}_{y}^{w}.
     \label{eq:med-bshift-binom-left}
   \end{dmath}
   Similarly, we have
   \begin{dmath}
     \frac{\prod \phi(y)^{n^{\ve}\rc w} \cdot
       [n^{\ve}]!_{y}^{w}}{\prod \phi(x)^{n^{\ve}\rc w}\cdot
       [n^{\ve}]!_{x}^{w}}
     \cdot [n^{\ve}]_{x}^{w}
     \cdot {b(l+1,2l)\brack a(l+1,2l)}_{y^{\ve}}^{w^{\ve}}
     =  [b(l+1,2l)]_{y^{\ve}}^{w^{\ve}}
     \cdot
     \frac{\prod \phi(y)^{(n^{\ve}-1)\rc w}\cdot
       [n^{\ve}-1]!_{y}^{w}}{\prod \phi(x)^{(n^{\ve}-1)\rc w}
       \cdot
       [n^{\ve}-1]!_{x}^{w}}
     \cdot
     \frac{\prod \phi(y)^{w}}{\prod\phi(x)^{w}}
     \cdot
     {b(l+1,2l) -1\brack a(l+1,2l)}_{y^{\ve}}^{w^{\ve}}.
     \label{eq:med-bshift-binom-right}
   \end{dmath}
 \end{lemma}
 \begin{proof}
   Since $(b-a)(1,l)=m\geq 1$, we have
   equation~\eqref{eq:med-bshift-binom-left} by
   \begin{align}
     &  \frac{ \prod_{i\in\oi(l)}\phi_{i}(y_{i})^{m_{i}w_{i}}
       [m_{i}]!_{y_{i}}^{w_{i}}}
       { \prod_{i\in\oi(l)}\phi_{i}(x_{i})^{m_{i}w_{i}}
       [m_{i}]!_{x_{i}}^{w_{i}}}
       \cdot
       \prod_{i\in\oi(l)}[m_{i}]_{x_{i}}^{w_{i}}
       \cdot \frac{[b_{i}]!_{y_{i}}^{w_{i}}
       }{[a_{i}]!_{y_{i}}^{w_{i}}[m_{i}]!_{y_{i}}^{w_{i}}}
     \\&=
     \frac{ \prod  \phi(y)^{w}}{ \prod \phi(x)^{w}}\cdot
     \frac{ \prod_{i\in\oi(l)}\phi_{i}(y_{i})^{(m_{i}-1)w_{i}}
     [m_{i}-1]!_{y_{i}}^{w_{i}}}
     { \prod_{i\in\oi(l)}
     \phi_{i}(x_{i})^{(m_{i}-1)w_{i}}      [m_{i}-1]!_{x_{i}}^{w_{i}}}
     \cdot \frac{[b_{i}]!_{y_{i}}^{w_{i}}
     }{[a_{i}]!_{y_{i}}^{w_{i}}[m_{i}-1]!_{y_{i}}^{w_{i}}}\\&
     =  [b(1,l)]_{y}^{w}\cdot
     \frac{ \prod \phi(y)^{(m-1)\rc w} \cdot
     [m-1]!_{y}^{w}}
     { \prod \phi(x)^{(m-1)\rc w} \cdot
     [m-1]!_{x}^{w}}
     \cdot
     \frac{ \prod \phi(y)^{w}}{\prod \phi(x)^{w}}
     \cdot
     {b(1,l)-1 \brack a(1,l)}_{y}^{w}.
   \end{align}
   Since $(b-a)(l+1,2l)=n\geq 1$,
   equation~\eqref{eq:med-bshift-binom-right} holds
   similarly by
   $\frac{[n^{\ve}]!_{y}^{w}}{[n]!_{y^{\ve}}^{w^{\ve}}}
   =\frac{[n^{\ve}-1]!_{y}^{w}}{[n-1]!_{y^{\ve}}^{w^{\ve}}}$.
 \end{proof}
 Then, we obtain
 the following merged-log-concavity for shift operators.
 
 \begin{proposition}\label{prop:shift}
   Consider
   $\cF_{i}=\Lam(s_{i},l,w,\scc,f_{i,s},\phi,\rho,x,\fX)$ for
   $i\in\oi(2)$ such that $\cFt=\sh_{h}(\cFo)$.  Let
   $y=x^{\rho}$.  Also, let $\mut=\tp{s_{2},l,m,n,k}$ be fitting.
   \begin{enumerate}
   \item \label{c:shift-induction}
     Then,  we have
     \begin{align}
       &\Delta(\cFt)(s_{2},l,w,m,n,k,\phi,\rho,x,\fX)
       \\&=
       \prod_{i\in\oi(h)}[\nu(m-i+1,n-i+1,k)]_{y^{\wcn}}^{w^{\wcn}}
       \\&
       \cdot \left(\frac{\prod \phi(y)^{w}}{\prod \phi(x)^{w}}\right)^{2h}
       \cdot \Delta(\cFo)(s_{1},l,w,m-h,n-h,k,\phi,\rho,x,\fX).
     \end{align}
   \item\label{c:shift-merged}
     Assume a squaring order $\scc'$ such that
     \begin{align}
       \left(\frac{\prod \phi(y)^{w}
       }{\prod \phi(x)^{w}}
       \right)^{2}\scc' 0.
      \label{ineq:shift-s2-1}
     \end{align}
     Then, $\cFt$ is $\scc'$-merged-log-concave, when $\cFo$
     is $\scc'$-merged-log-concave.
  \end{enumerate}
\end{proposition}
\begin{proof}
  Let $a=\nu(k)$ and $b=\nu(m,n,k)$.  Let us prove
  Claim~\ref{c:shift-induction}.  If
  Claim~\ref{c:shift-induction} is true for $h=1$, then
  Lemma~\ref{lem:one-shift} yields
  Claim~\ref{c:shift-induction} for each $h\in\Zgeo$.  Thus,
  suppose $h=1$.  Since
  $f_{2,s_{1},n^{\ve}-1}= f_{2,s_{1},(n-1)^{\ve}}$ and
  $f_{2,s_{1},(n\bop k)^{\ve}-1}= f_{2,s_{1},(n\bop
    k-1)^{\ve}}$,
  Claim~\ref{c:merged-binom-bshift-left-right-equations} of
  Theorem~\ref{thm:merged-binom-bshift} gives
  \begin{align}
    &\Delta(\cFt)(s_{2},l,w,m,n,k,\phi,\rho,x,\fX)
    \\&=
    f_{2,s_{2},m}f_{2,s_{2},n^{\ve}}
    \cdot B(s_{2},l,w,m,n^{\ve},\phi,\rho,x,\fX)          
    { b \brack a}^{w^{\wcn}}_{y^{\wcn}}
    \\   &   -
           f_{2,s_{2},m\bom k}f_{2,s_{2},(n\bop k)^{\ve}}
           \cdot 
           B(s_{2},l,w,m\bom k,(n\bop k)^{\ve},\phi,\rho,x,\fX)          
           { b \brack a^{\ve}}^{w^{\wcn}}_{y^{\wcn}}
    \\&=
    f_{2,s_{1},m-1}[m]_{x}^{w}f_{2,s_{1},(n-1)^{\ve}}
    [n^{\ve}]_{x}^{w}
    \cdot
    B(s_{2},l,w,m,n^{\ve},\phi,\rho,x,\fX)          
    { b \brack a}^{w^{\wcn}}_{y^{\wcn}}
    \\ &
         -
         f_{2,s_{1},m\bom k-1}
         [m\bom k]_{x}^{w}
         f_{2,s_{1},(n\bop k-1)^{\ve}}
         [(n\bop k)^{\ve}]_{x}^{w}
    \\ &\cdot
         B(s_{2},l,w,m\bom k,(n\bop k)^{\ve},\phi,\rho,x,\fX)          
         { b \brack a^{\ve}}^{w^{\wcn}}_{y^{\wcn}}.
         \label{eq:shift-p2}
  \end{align}
  
  Moreover, $\mut$ is fitting by the assumption and
  $\muo=\tp{s_{1},m-1,n-1,k}$ is fitting by
  Claim~\ref{c:fit-shift-fitting} of
  Lemma~\ref{lem:fit-shift}.  In particular,
  $m\ccn n \ldZtl s_{2}$ and $(m-1)\ccn (n-1) \ldZtl s_{1}$
  so that $B(s_{2},l,w,m,n^{\ve},\phi,\rho,x,\fX)\neq0 $ and
  $B(s_{1},l,w,m-1,(n-1)^{\ve},\phi,\rho,x,\fX)\neq0$.  Also, we
  have $b-a=m\ccn n\geq h=1$, and $a,b\geq 0$ by
  Claim~\ref{c:fitting-nonneg-ab-nonneg} in
  Lemma~\ref{lem:fitting-nonneg}.  Thus,
  Lemma~\ref{lem:med-bshift-binom} yields
  \begin{align}
    &[m]_{x}^{w}  [n^{\ve}]_{x}^{w}
      B(s_{2},l,w,m,n^{\ve},\phi,\rho,x,\fX)          
      { b \brack a}^{w^{\wcn}}_{y^{\wcn}}
    \\&=
    [m]_{x}^{w}  [n^{\ve}]_{x}^{w}
    B(s_{2},l,w,m,\phi,\rho,x,\fX)
    B(s_{2},l,w,n^{\ve},\phi,\rho,x,\fX) 
    { b \brack a}^{w^{\wcn}}_{y^{\wcn}}
    \\&=
    \frac{\prod \phi(y)^{m\rc w}
    \cdot [m]!_{y}^{w}}
    {\prod \phi(x)^{m\rc w}
    \cdot [m]!_{x}^{w}}
    \cdot  [m]_{x}^{w}
    \cdot {b(1,l) \brack a(1,l)}_{y}^{w}
    \cdot
    \frac{\prod \phi(y)^{n^{\ve}\rc w}
    \cdot [n^{\ve}]!_{y}^{w}}
    {\prod \phi(x)^{n^{\ve}\rc w}
    \cdot [n^{\ve}]!_{x}^{w}}
    \cdot  [n^{\ve}]_{x}^{w}
    \cdot {b(l+1,2l) \brack a(l+1,2l)}_{y}^{w}
    \\ &=
         [b(1,l)]_{y}^{w}\cdot
         \frac{ \prod \phi(y)^{(m-1)\rc w}\cdot
         [m-1]!_{y}^{w}}{ \prod \phi(x)^{(m-1)\rc w}\cdot
         [m-1]!_{x}^{w}}
         \cdot
         \frac{ \prod \phi(y)^{w}}{\prod \phi(x)^{w}}
         \cdot
         {b(1,l)-1 \brack a(1,l)}_{y}^{w}
    \\& \cdot
    [b(l+1,2l)]_{y^{\ve}}^{w^{\ve}}
    \cdot
    \frac{\prod \phi(y)^{(n^{\ve}-1)\rc w}\cdot
    [n^{\ve}-1]^{w}_{y}}{\prod \phi(x)^{(n^{\ve}-1)\rc w}
    \cdot
    [n^{\ve}-1]^{w}_{x}}
    \cdot
    \frac{\prod \phi(y)^{w}}{\prod\phi(x)^{w}}
    \cdot
    {b(l+1,2l) -1\brack a(l+1,2l)}_{y^{\ve}}^{w^{\ve}}
    \\&=
    [b]_{y^{\wcn}}^{w^{\wcn}}\left( \frac{ \prod\phi(y)^{w}}
    { \prod \phi(x)^{w}}\right)^{2}
    B(s_{1},l,w,m-1,(n-1)^{\ve},\phi,\rho,x,\fX)
    { b-1 \brack a}_{y^{\wcn}}^{w^{\wcn}}.
    \label{eq:shift-p4}
  \end{align}
  
  Similarly, $b-a^{\ve}=(m\bom k)\ccn (n\bop k)$ by
  Lemma~\ref{lem:ab-define}.  If $\mut$ is wrapped, then
  $(m\bom k) \ccn (n\bop k)\geq h=1$. Hence,
  Lemma~\ref{lem:med-bshift-binom} again gives
  \begin{dmath}
    [m\bom k]_{x}^{w}
    [(n\bop k)^{\ve}]_{x}^{w}
    B(s_{2},l,w,m\bom k,(n\bop k)^{\ve},\phi,\rho,x,\fX)
    { b \brack a^{\ve}}^{w^{\wcn}}_{x^{\wcn}}
    =
    [b]_{y^{\wcn}}^{w^{\wcn}}\cdot
    \left(\frac{\prod\phi(y)^{w}}{\prod\phi(x)^{w}}\right)^{2}
    \cdot
    B(s_{1},l,w,m\bom k-1,(n\bop k-1)^{\ve},\phi,\rho,x,\fX)
    \cdot { b-1 \brack a^{\ve}}_{y^{\wcn}}^{w^{\wcn}}.
    \label{eq:shift-p5}
  \end{dmath}
  If $\mut$ is unwrapped, then $\muo$ is unwrapped by
  Claim~\ref{c:fit-shift-wrapped} of
  Lemma~\ref{lem:fit-shift}.  Thus,
  equation~\eqref{eq:shift-p5} still holds by $0=0$.

  Therefore, equations~\eqref{eq:shift-p2},
  \eqref{eq:shift-p4}, and~\eqref{eq:shift-p5} give
  \begin{align}
    &\Delta(\cFt)(s_{2},l,w,m,n,k,\phi,\rho,x,\fX)
    \\&=
    [b]_{y^{\wcn}}^{w^{\wcn}}
    \left(\frac{\prod \phi(y)^{w}}{\prod \phi(x)^{w}}\right)^{2}
    \Delta(\cFo)
    (s_{1},l,w,m-1,n-1,k,\phi,\rho,x,\fX).
    \label{eq:shift-p6}
  \end{align}
  Hence, Claim~\ref{c:shift-induction} holds for
  $h=1$.

  Let us prove Claim~\ref{c:shift-merged}.  Since
  $b=m\ccn n+a\geq h$, $\nu(m-i+1,n-i+1,k)=b-i+1\geq 1$ for
  $i\in\oi(h)$.  In particular,
  $\prod_{i\in\oi(h)}[\nu(m-i+1,n-i+1,k)]_{y^{\wcn}}^{w^{\wcn}}\scc
  0$ by Claim~\ref{c:adm-succ-half-gx-scc} of
  Lemma~\ref{lem:adm-succ}.  Then,
  Claim~\ref{c:shift-merged} holds, because
  Claim~\ref{c:shift-induction} and
  inequality~\eqref{ineq:shift-s2-1} imply
  $\Delta(\cFt)(s_{2},l,w,m,n,k,\phi,\rho,x,\fX) \scc'0$.
\end{proof}

We state the following corollary.
\begin{corollary}\label{cor:shift-ineq}
  Let $\cF_{i}=\Lam(s_{i},l,w,\scc,f_{i,s},\phi,\rho,x,\fX)$ for
  $i\in\oi(2)$ such that
  $\cFt=\sh_{h}(\cFo)$.
  Assume one of the following
  four cases.
  \begin{enumerate}
  \item \label{c:shift-ineq-canonical-med}
    $\phi$ is the canonical mediator.
  \item \label{c:shift-ineq-unit}
    $1\ldZ s_{1}$.
  \item \label{c:shift-ineq-weight-zero}
    $w=\iota^{l}(0)$.
  \item\label{c:shift-ineq-zero-triv-shift}
    $\rho=\iota^{l}(1)$.
  \end{enumerate}
  Then, we have
  \begin{align}
    \frac{ \prod \phi(y)^{w}}{ \prod \phi(x)^{w}}\scc 0.
    \label{ineq:shift-scc}
  \end{align}
  In particular, if $\cFo$ is $\scc'$-merged-log-concave, then
  $\cFt$ is $\scc'$-merged-log-concave.
\end{corollary}
\begin{proof}
  First, suppose Case~\ref{c:shift-ineq-canonical-med}.
  Then, since $\rho\geq 1$, we have
  $\frac{\prod \phi(y)^{w}}{\prod \phi(x)^{w}}
  =\prod_{i\in\oi(l)}
  \frac{(1-x_{i}^{\rho_{i}})^{w_{i}} }{(1-x_{i})^{w_{i}}}
  =\prod_{i\in\oi(l)}[\rho_{i}]^{w_{i}}_{x_{i}}\scc 0$ by
  Claim~\ref{c:adm-succ-half-gx-scc} of
  Lemma~\ref{lem:adm-succ}.  Second, suppose
  Case~\ref{c:shift-ineq-unit} so that
  $m=\iota^{l}(1)\ldZl s_{1}$. Then, since $\phi(x)$ is the
  mediator of $\cFo$, its base-shift positivity gives
  \begin{align}
    0\prec B(s_{1},l,w,m,\phi,\rho,x,\fX)
    =
    \prod_{i\in\oi(l)}  \frac{ \phi_{i}(y_{i})^{m_{i}w_{i}}
    [m_{i}]!_{y_{i}}^{w_{i}}}
    { \phi_{i}(x_{i})^{m_{i}w_{i}}      [m_{i}]!_{x_{i}}^{w_{i}}}
    =     \frac{\prod \phi(y)^{w}   }{\prod \phi(x)^{w}}.
  \end{align}
  Third, inequality~\eqref{ineq:shift-scc} holds for
  other cases, since
  $\frac{\prod \phi(y)^{w}}{\prod \phi(x)^{w}}=1\scc 0$.
  
  The latter statement holds by Claim~\ref{c:shift-merged} of
  Proposition~\ref{prop:shift}
  and the multiplicativity of $\scc$.
\end{proof}

\section{Explicit merged-log-concave
  parcels}
\label{sec:explicit-merged-parcels}

By monomial conditions in
Definition~\ref{def:monomial-index},
we construct explicit merged-log-concave parcels for
arbitrary gates, widths, base shifts, and positive
weights.   Also, we discuss
several conjectures on merged determinants.

\subsection{Base shift functions, pre-parcels,
  and pre-merged determinants}
\label{sec:bshift-pre-parcels-pre-merged-det}
We state the following lemmas on the base shift functions
$b_{\lam,\rho}(q)$ in Definition~\ref{def:intro-base-shift}.
\begin{lemma}\label{lem:bsingle-qnum}
  We have
  \begin{align}
    b_{\lam,\rho}(q)
    =
    \begin{dcases}
      \prod_{h\in\oi(\lam)}[\rho]_{q^{h}}\llq0      \mif \lam\geq 1,\\
      1\llq 0  \mif \lam=0.
    \end{dcases}
  \end{align}
\end{lemma}
\begin{proof}
  The assertion follows from
    $b_{\lam,\rho}(q)
    = \prod_{h\in \oi(\lam)}\frac{1-q^{\rho h}}{1-q^{h}}$ such that
    $\frac{1-q^{\rho h}}{1-q^{h}}
    =[\rho]_{q^{h}}$.
  \end{proof}

\begin{lemma}\label{lem:bsingle-degree}
  If  $\lam\in \Zgez$, then
  $\deg_{q} b_{\lam,\rho}(q)=\frac{(\rho-1)\lam(\lam+1)}{2}$.
\end{lemma}
\begin{proof}
  By Lemma~\ref{lem:bsingle-qnum}, 
  $\deg_{q} b_{\lam,\rho}(q)     
    =
      \deg_{q} q^{\rho-1} q^{2(\rho-1)} \dots       q^{\lam(\rho-1)}
      =\frac{(\rho-1)\lam(\lam+1)}{2}$.
\end{proof}

\begin{lemma}\label{lem:bsingle-quot}
  Consider  $\lam\geq k\geq 0$.  Then, we have
  \begin{align}
    \frac{b_{\lam,\rho}(q)}{b_{\lam-k,\rho}(q)}
    =
    \begin{dcases}
      \sum_{ j_{1},\dots, j_{k}\in\oi(0,\rho-1)}
      q^{\sum_{h\in\oi(k)} j_{h}(\lam-k+h)}
      \llq 0 \mif k\geq 1,\\
      1\llq 0 \melse.
    \end{dcases}
  \end{align}
\end{lemma}
\begin{proof}
  If $k=0$, then Lemma~\ref{lem:bsingle-qnum} yields
  the assertion by
  $b_{\lam,\rho}(q)=b_{\lam-k,\rho}(q)\neq 0$.
  If $k\geq 1$, then Lemma~\ref{lem:bsingle-qnum} gives
  \begin{align}
    \frac{b_{\lam,\rho}(q)}{b_{\lam-k,\rho}(q)} 
    &=\prod_{\ka \in\oi(\lam-k+1,\lam)}
      (1+q^{\ka}+q^{2\ka}+
         \dots+q^{(\rho-1)\ka}).
  \end{align}
  Thus, we have the assertion, 
  choosing $q^{j_{i}(\lam-k+i)}$ for some
  $j_{i}\in\oi(0,\rho-1)$ in each factor.
\end{proof}

We introduce a tuple version of 
$b_{\lam,\rho}(q)$, pre-parcels, and pre-merged determinants.

\begin{definition}\label{def:base-pre-parcel-merged}
  Let $\mu\in \Zgeo$ and $\lam=\floor{\frac{\mu}{2}}$.
  Suppose an indeterminate
  $x\in \Q(\fX)^{\mu}$ and $w,\rho\in \Zgeo^{\mu}$ such that
  $x,w,\rho$ are palindromic.  Let $\phi$ be the canonical
  $\mu$-mediator.
  \begin{enumerate}
  \item
    \label{d:base-pre-parcel-merged-base}
    If $m\in \Z^{\mu}$, then we define the base shift
    function
    \begin{align}
      B(\mu,w,m,\rho,x)
      =\prod_{i\in\oi(\mu)}b^{w_{i}}_{m_{i}, \rho_{i}}(x_{i})
      \in \Q(\fX).
    \end{align}
  \item
    \label{d:base-pre-parcel-merged-parcel}
    We define the pre-parcel
    $U(\mu,w,\rho,x)=\{ U_{a}^{b}(\mu,w,\rho,x)\in
    \Q(\fX)\}_{a,b\in\Z^{\mu}}$ such that
    \begin{align}
      U_{a}^{b}(\mu,w,\rho,x)
      &={ b \brack a}^{w}_{x^{\rho}} 
        B(\mu,w,b-a,\rho,x).
    \end{align}
    We refer to $\mu$, $w$, $\rho$, and $x$ as the width,
    weight, base shift, and base of $U(\mu,w,\rho,x)$.
  \item
    \label{d:base-pre-parcel-merged-merged}
    If $e\in \Zgez^{\lam}$, then we define the pre-merged
    determinant
    \begin{align}
      d(U)_{a}^{b}(\mu,w,\rho,e,x)
      =U_{a}^{b}(\mu,w,\rho,x)- U_{a^{\ve}}^{b}(\mu,w,\rho,x)
      \cdot \prod x(1,\lam)^{e}       \in \Q(\fX).
    \end{align}
    We call $e$ degree shift of $d(U)_{a}^{b}(\mu,w,\rho,e,x)$.
  \item If $x$ is flat and
    $e\in \Zgez$, then let
    \begin{align}
      \dt(U)_{a}^{b}(\mu,w,\rho,e,x)
      =d(U)_{a}^{b}(\mu,w,\rho,\tp{e}\ccn \iota^{\lam-1}(0),x).
    \end{align}      
  \end{enumerate}
\end{definition}

We put the integers below to discuss
pre-parcels and pre-merged determinants by $\deg_{q}$.

\begin{definition}\label{def:folded-diff}
  Let $\mu\in \Zget$, $\lam=\floor{\frac{\mu}{2}}$, and
  $a,b,w,\rho\in \Z^{\mu}$.
  \begin{enumerate}
  \item We put $N_{a}^{b}(\mu,w,\rho)
    =\tp{w_{i}\rho_{i} (b_{\mu-i+1}-b_{i})(a_{\mu-i+1}-a_{i})}_{i\in \oi(\lam)}
    \in \Z^{\lam}$.
    Also, let
    \begin{align}
      n_{a}^{b}(\mu,w,\rho)
      &=\sum N_{a}^{b}(\mu,w,\rho),\\
      o_{a}^{b}(\mu,w,\rho)
      &=\floor*{\frac{n_{a}^{b}(\mu,w,\rho)}{2}}.
    \end{align}
  \item
    \label{d:folded-diff-triv-bshift}
    If $\rho=\iota^{\mu}(1)$, then we put
    \begin{align}
      N_{a}^{b}(\mu,w)&=N_{a}^{b}(\mu,w,\rho),\\
      n_{a}^{b}(\mu,w)&=n_{a}^{b}(\mu,w,\rho),\\
      o_{a}^{b}(\mu,w)&=o_{a}^{b}(\mu,w,\rho).
    \end{align}
  \end{enumerate}
\end{definition}

Then, we have the following degree differences of
width-two pre-parcels.

\begin{lemma}\label{lem:deg-diff-width-two}
  Let $a,b\in \Zt$ such that $a,b-a,b-a^{\ve}\geq 0$.  Consider
  flat $w,\rho \in \Zgeo^{2}$ and
  $x=\iota^{2}(q)\in \Q(\fX)^{2}$. Then, we have
  \begin{align}
    \deg_{q}
    \left({ b \brack a}^{w}_{x^{\rho}}\right)-
    \deg_{q}
    \left({ b \brack a^{\ve}}^{w}_{x^{\rho}}\right)
    &=n_{a}^{b}(2,w,\rho),
      \label{eq:deg-diff-width-two-with-shift}\\
    \deg_{q}( U_{a}^{b}(2,w,\rho,x))-
    \deg_{q}(U_{a^{\ve}}^{b}(2,w,\rho,x))
    &=n_{a}^{b}(2,w).
      \label{eq:deg-diff-width-two-sem-shift}
  \end{align}
\end{lemma}
\begin{proof}
  Let us prove
  equation~\eqref{eq:deg-diff-width-two-with-shift}.
  Since $\deg_{q}{B \brack A}_{q}=A(B-A)$
  for
  integers $B\geq A\geq 0$, 
   flat $w$ and
  $\rho$ give
  \begin{align}
    \deg_{q}
    \left({ b \brack a}^{w}_{x^{\rho}}\right)
    &=
    \deg_{q}\left( {b_{1} \brack a_{1}}^{w_{1}}_{q^{\rhoo}}
    {b_{2} \brack a_{2}}^{w_{1}}_{q^{\rhoo}}\right)
    =\rhoo w_{1}(a_{1}(b_{1}-a_{1})+a_{2}(b_{2}-a_{2})),\\
    \deg_{q}\left({ b \brack a^{\ve}}^{w}_{x^{\rho}}\right)
&     =\rhoo w_{1}(a_{2}(b_{1}-a_{2})+a_{1}(b_{2}-a_{1})).
  \end{align}
  Thus, equation~\eqref{eq:deg-diff-width-two-with-shift}
  holds, as
  $a_{1}b_{1}+a_{2}b_{2} -a_{2}b_{1}-a_{1}b_{2}=
  (b_{2}-b_{1})(a_{2}-a_{1})$.
  
  Let us prove
  equation~\eqref{eq:deg-diff-width-two-sem-shift}.  We have
  \begin{align}
    &\deg_{q}( U_{a}^{b}(2,w,\rho,x))-
      \deg_{q}(U_{a^{\ve}}^{b}(2,w,\rho,x))
    \\&=\deg_{q}
    \left({ b \brack a}^{w}_{x^{\rho}}B(2,w,b-a,\rho,x)\right)
    -
    \deg_{q}
    \left({ b \brack a^{\ve}}^{w}_{x^{\rho}}
    B(2,w,b-a^{\ve},\rho,x) \right).
  \end{align}
  Furthermore, Lemma~\ref{lem:bsingle-degree} yields
  \begin{align}
    &\deg_{q}(B(2,w,b-a,\rho,x))
      -  \deg_{q}(B(2,w,b-a^{\ve},\rho,x))
    \\&=
    \deg_{q}(b_{b_{1}-a_{1},\rhoo}^{w_{1}}(q)
    b_{b_{2}-a_{2},\rhoo}^{w_{1}}(q))
    -
    \deg_{q}(b_{b_{1}-a_{2},\rhoo}^{w_{1}}(q)
    b_{b_{2}-a_{1},\rhoo}^{w_{1}}(q))
    \\&
    =\frac{w_{1}(\rhoo-1)}{2}
    \cdot
    \bigg( (b_{1}-a_{1}) (b_{1}-a_{1}+1)
    +(b_{2}-a_{2}) (b_{2}-a_{2}+1)\\&
    -
    (b_{1}-a_{2}) (b_{1}-a_{2}+1)
    -(b_{2}-a_{1}) (b_{2}-a_{1}+1)\bigg)\\&
    =-\frac{w_{1}(\rhoo-1)}{2}
    \cdot
    ((b_{1}-a_{2})^{2} +(b_{2}-a_{1})^{2}
    -(b_{1}-a_{1})^{2}
    -(b_{2}-a_{2})^{2} )\\&
    =-(\rhoo-1)w_{1}(b_{2}-b_{1})(a_{2}-a_{1}).
  \end{align}
  Thus, equation~\eqref{eq:deg-diff-width-two-with-shift}
  gives equation~\eqref{eq:deg-diff-width-two-sem-shift}.
\end{proof}

We further introduce the following integers and tuples to
discuss 
the pre-parcels and pre-merged determinants of general widths.

\begin{definition}\label{def:center-outer}
  Let $\mu\in \Z_{\geq 3}$  and $a\in \Z^{\mu}$.
  Also, for $\lam=\floor{\frac{\mu}{2}}$, let
  $e\in\Z^{\lam}$.
  \begin{enumerate}
  \item Let $\chi(\mu)=2^{1-\ceil{\frac{\mu}{2}}+\lam}$.
  \item If $\mu$ is even, then we put the following
    tuples:
    \begin{align}
      \out(a)&=a(1,\lam-1)\ccn a(\lam+2,\mu)
               \in \Z^{\mu-\chi(\mu)};\\
      \cen(a)&=a(\lam,\lam+1)
               \in \Z^{\chi(\mu)};\\
      \out(\mu,e)&=e(1, \lam-1)
                  \in \Z^{\lam-\chi(\mu)+1};\\
      \cen(\mu,e)&=e(\lam,\lam)\in \Z^{\chi(\mu)-1}.
    \end{align}
  \item If $\mu$ is odd, then we put
    the following tuples:
    \begin{align}
      \out(a)
      &=a(1,\lam)
        \ccn a(\lam+2,\mu)\in \Z^{\mu-\chi(\mu)};\\
      \cen(a)&=a(\lam+1,\lam+1)\in \Z^{\chi(\mu)};\\
      \out(\mu,e)&=e\in \Z^{\lam-\chi(\mu)+1}.
    \end{align}
  \end{enumerate}
  We call $\out(a)$ and $\cen(a)$ outer and center tuples of
  $a$.  Also, we call $\out(\mu,e)$ and $\cen(\mu,e)$
  outer and center tuples of $e$.
\end{definition}

Then, we have the following lemma for pre-parcels and pre-merged determinants.
\begin{lemma}\label{lem:outer-preparcel-premerged}
  For $\mu\in \Z_{\geq 3}$, suppose $a,b\in \Z^{\mu}$.  Let
  $\lam=\floor{\frac{\mu}{2}}$.
  \begin{enumerate}
  \item
    \label{c:outer-preparcel-premerged-preparcel}
    Then, it holds that
    \begin{align}
      U_{a}^{b}(\mu,w,\rho,x)
      &= U_{\out(a)}^{\out(b)}(\mu-\chi(\mu),
        \out(w),\out(\rho),\out(x))
        \cdot U_{\cen(a)}^{\cen(b)}
        (\chi(\mu),\cen(w),\cen(\rho),\cen(x)),
      \\ 
      U_{a^{\ve}}^{b}(\mu,w,\rho,x)
      &=
        U_{\out(a)^{\ve}}^{\out(b)}
        (\mu-\chi(\mu),\out(w),\out(\rho),\out(x))
        \cdot U_{\cen(a)^{\ve}}^{\cen(b)}
        (\chi(\mu),\cen(w),\cen(\rho),\cen(x)).
    \end{align}
  \item
    \label{c:outer-preparcel-premerged-premerged}
    In particular, if $\mu$ is odd, then
    \begin{dmath*}
      d(U)_{a}^{b}(\mu,\rho, w, e ,x)
      =U_{\cen(a)}^{\cen(b)}
      (\chi(\mu),\cen(w),\cen(\rho),\cen(x))
      \cdot 
      d(U)_{\out(a)}^{\out(b)}
      (\mu-\chi(\mu), \out(v), \out(\rho),\out(\mu,e), \out(x));
    \end{dmath*}
    also, if $\mu$ is even, then
    \begin{dmath*}
      d(U)_{a}^{b}(\mu,\rho, w, e ,x)
      =U_{\cen(a)}^{\cen(b)}(\chi(\mu),\cen(w),\cen(\rho),\cen(x))
      \cdot U_{\out(a)}^{\out(b)}(\mu-\chi(\mu),\out(w),\out(\rho),\out(x))
      -
      \prod x(\lam,\lam)^{\cen(\mu,e)}\cdot
      U_{\cen(a)^{\ve}}^{\cen(b)}(\chi(\mu),\cen(w),\cen(\rho),\cen(x))
      \cdot \prod x(1,\lam-1)^{\out(\mu,e)}
      \cdot
      U_{\out(a)^{\ve}}^{\out(b)}
      (\mu-\chi(\mu), \out(w),\out(\rho),\out(x)).
    \end{dmath*}
  \end{enumerate}
\end{lemma}
\begin{proof}
  Let us prove
  Claim~\ref{c:outer-preparcel-premerged-preparcel}.
  Since $\out(w)$, $\out(\rho)$, $\out(x)$, and $\cen(w)$,
  $\cen(\rho)$, $\cen(x)$ are palindromic, each factor in
  equations of 
  Claim~\ref{c:outer-preparcel-premerged-preparcel} exists.
  Also, $a=\out(a)(1,\lam-\chi(\mu)+1) \ccn \cen(a)\ccn
  \out(a)(\lam-\chi(\mu)+2,\mu-\chi(\mu))$ and
  $b=\out(b)(1,\lam-\chi(\mu)+1) \ccn \cen(b) \ccn
  \out(b)(\lam-\chi(\mu)+2,\mu-\chi(\mu))$.
  Thus, Claim~\ref{c:outer-preparcel-premerged-preparcel}
  follows from
  $\out(a^{\ve})=\out(a)^{\ve}$ and
  $\cen(a^{\ve})=\cen(a)^{\ve}$.
  Claim~\ref{c:outer-preparcel-premerged-premerged} holds by
  Claim~\ref{c:outer-preparcel-premerged-preparcel},
  since $C(a)=C(a)^{\ve}$ for odd $\mu$.
\end{proof}

We put the following notions
to later discuss merged determinants by pre-merged ones.
\begin{definition}\label{def:pre-fitting}
  Let $\mu\in \Zget$ and $a,b\in \Z^{\mu}$.
  Consider a tuple $\omg=\tp{\mu,a,b}$.
  \begin{enumerate}
  \item We call $\omg$ pre-fitting, if we have the
    following  conditions:
    \begin{enumerate}
    \item $a$ is non-negative and increasing;
    \item $b$ is increasing;
    \item $b_{1}<b_{\mu}$ and  $a_{1}<a_{\mu}$;
      \label{c:pre-fitting-end-term} 
    \item $a\leq b$. 
    \end{enumerate} 
    We refer to Condition~\ref{c:pre-fitting-end-term} as
    end slope condition of $\omg$.
  \item We call $\omg$ tempered, if $ a^{\ve}\leq b$.
  \item If $\mu\geq 3$, then we define the outer tuple
    $\tout(\omg)= \tp{\mu-\chi(\mu), \out(a),\out(b)}$.
  \end{enumerate}
\end{definition}

For example, if $\tp{\mu,a,b}$ is pre-fitting and tempered,
then
$0\leq a_{1}\leq a_{2}\leq \dots \leq a_{\mu}\leq
  b_{1}\leq b_{2}\leq \dots \leq b_{\mu}$.

  We use the term ``pre-fitting''  by the following statement.

  \begin{proposition}\label{prop:fitting-pre-fitting}
    Suppose a gate $s\geq 0$, $l\in \Zgeo$, $m,n\in \Zl$,
  and $k\in\Ztl$.  Let $\omg_{1}=\tp{s,l,m,n,k}$, $a=\nu(k)$,
  $b=\nu(m,n,k)$, and $\omg_{2}=\tp{2l,a,b}$.  Then, we have
  the following.
  \begin{enumerate}
  \item\label{c:fitting-pre-fitting-fitting}
    If $\omg_{1}$ is fitting, then $\omg_{2}$ is pre-fitting.
  \item\label{c:fitting-pre-fitting-wrapped}
    If $\omg_{1}$ is wrapped, then $\omg_{2}$ is tempered.
  \end{enumerate}
\end{proposition}
\begin{proof}
  First, let us confirm
  Claim~\ref{c:fitting-pre-fitting-fitting}.
  The slope conditions of $\omg_{1}$
  give the end slope condition of $a$ and $b$.
  Also, $b\geq a\geq 0$ by
  Claim~\ref{c:fitting-nonneg-ab-nonneg} of
  Lemma~\ref{lem:fitting-nonneg}.  Hence, $\omg_{2}$ is
  pre-fitting.
  Second, let us prove Claim~\ref{c:fitting-pre-fitting-wrapped}.
  If $\omg_{1}$ is wrapped, then $b\geq a^{\ve}$, because
  $s\geq 0$ and
  $b-a^{\ve}=(m\bom k)\ccn (n\bop k)\ldZtl s$ by
  Lemma~\ref{lem:ab-define}.
\end{proof}

The converse of Claim~\ref{c:fitting-pre-fitting-fitting} in
Proposition~\ref{prop:fitting-pre-fitting} does not hold,
because $b-a\ldZtl s$ is not necessarily true.  Likewise,
the converse of Claim~\ref{c:fitting-pre-fitting-wrapped} in
Proposition~\ref{prop:fitting-pre-fitting} does not hold.

Moreover, we have the following for outer tuples $\tout(\omg)$.
\begin{lemma}\label{lem:center-outer}
  Let $\mu\in \Z_{\geq 3}$, $a,b \in \Z^{\mu}$, and
  $\omg=\tp{\mu,a,b}$.
  \begin{enumerate}
  \item \label{c:center-outer-pre-fitting} If $\omg$ is
    pre-fitting, then $\tout(\omg)$ is pre-fitting.
  \item \label{c:center-outer-tempered} If $\omg$ is
    tempered, then $\tout(\omg)$ is tempered.
  \end{enumerate}
\end{lemma}
\begin{proof}
  Claim~\ref{c:center-outer-pre-fitting} holds, since $a$
  and $\out(a)$ (or $b$ and $\out(b)$) have the same end
  terms.  Also, Claim~\ref{c:center-outer-tempered} holds, since
  $\out(a^{\ve})=\out(a)^{\ve}$ gives
  $\out(b)\geq \out(a)^{\ve}$.
\end{proof}

Equation~\eqref{eq:deg-diff-width-two-sem-shift} in
Lemma~\ref{lem:deg-diff-width-two} is  independent on base
shifts $\rho$.  This extends to the following pre-parcels
of general widths.

\begin{proposition}\label{prop:pre-parcel-deg}
  Suppose a tempered pre-fitting $\omg=\tp{\mu,a,b}$.  Let
  $x=\iota^{\mu}(q)$.  Then,
  \begin{align}
    \deg_{q}( U_{a}^{b}(\mu,w,\rho,x))-
    \deg_{q}(U_{a^{\ve}}^{b}(\mu,w,\rho,x))
    &=n_{a}^{b}(\mu,w),
      \label{eq:pre-parcel-deg-diff}\\
    n_{a}^{b}(\mu,w)&>0.
                     \label{ineq:pre-parcel-deg-pos}
  \end{align}
\end{proposition}
\begin{proof}
  Let us prove equation~\eqref{eq:pre-parcel-deg-diff} by
  induction on $\mu$.  Suppose $\mu=2$.  Then,
  $w$ and $\rho$ are flat, as they are palindromic.
  Thus, Lemma~\ref{lem:deg-diff-width-two} gives
  equation~\eqref{eq:pre-parcel-deg-diff},
  since we have
  $a,b-a,b-a^{\ve}\geq 0$ for the tempered pre-fitting
  $\omg$.

  Hence, assume $\mu>2$.  Since $\omg$ is tempered,
  $U_{a}^{b}(\mu,w,\rho,x)$ and
  $U_{a^{\ve}}^{b}(\mu,w,\rho,x)$ are non-zero by
  $b-a,b-a^{\ve},a\geq 0$.  Moreover,
  $\tout(\omg)=\tp{\mu-\chi(\mu), \out(a),\out(b)}$ is tempered
  pre-fitting by Lemma~\ref{lem:center-outer}.  Thus, since
  $\cen(a)=\cen(a)^{\ve}$ for odd $\mu$, the induction on
  $\mu$ and Claim~\ref{c:outer-preparcel-premerged-preparcel}
  of Lemma~\ref{lem:outer-preparcel-premerged} give
  \begin{dmath*}
    \deg_{q}( U_{a}^{b}(\mu,w,\rho,x))-
    \deg_{q}(U_{a^{\ve}}^{b}(\mu,w,\rho,x))
    =
    \begin{dcases}
      n_{\out(a)}^{\out(b)}(\mu-\chi(\mu),\out(w)) \mif \mu
      \mbox{ is odd},\\
      n_{\out(a)}^{\out(b)}(\mu-\chi(\mu),\out(w))
      + n_{\cen(a)}^{\cen(b)}(\chi(\mu),\cen(w)) \mif \mu
      \mbox{ is even},
    \end{dcases}
    = n_{a}^{b}(\mu,w).
  \end{dmath*}
  
  Let us prove inequality~\eqref{ineq:pre-parcel-deg-pos}.
  Each $i\in\oi(\floor{\frac{\mu}{2}})$ gives
  $ w_{i}(b_{\mu-i+1}-b_{i})(a_{\mu-i+1}-a_{i})\geq 0$, since
  $a$ and $b$ are increasing for the pre-fitting $\omg$.
  Then, the end slope condition of $\omg$ yields
  $n_{a}^{b}(\mu,w)
  =
  \sum_{i\in\oi(\floor{\frac{\mu}{2}})}
  w_{i}(b_{\mu-i+1}-b_{i})(a_{\mu-i+1}-a_{i})
  \geq w_{1}(b_{\mu}-b_{1})(a_{\mu}-a_{1}) >0$.
\end{proof}

\subsection{Positivity of pre-merged determinants}
\label{sec:positivity-pre-merged-det}
We obtain polynomials with positive integer
coefficients by pre-merged determinants.  
For this, we recall the following notions of
$q$-polynomials.

\begin{definition}\label{def:sym-uni-lc}
  Suppose $f(q)\in \Q[q]$.
  \begin{enumerate}
  \item The polynomial $f$ is a palindromic $q$-polynomial, if
    $f_{\ord_{q}(f)+i}=f_{\deg_{q}(f)-i}$
    for each  $i\in \Z$ such that
    $0\leq i \leq \deg_{q}(f)-\ord_{q}(f)$.
  \item
    \label{d:sym-uni-lc-uni}
    The polynomial $f$ is a unimodal $q$-polynomial, if $f_{i-1}>f_{i}$ for
    some $i\in \Z$ such that $\ord_{q}(f) < i<\deg_{q}(f)$,
    then $f_{j-1}\geq f_{j}$ for each $j\in \Z$ such that
    $i< j \leq \deg_{q}(f)$.
  \item The polynomial $f$ is a log-concave $q$-polynomial, if
    $f_{i}^{2}-f_{i-1}f_{i+1}\geq 0$ for each $i\in \Z$ such
    that $\ord_{q}(f)<i<\deg_{q}(f)$.
  \end{enumerate}
\end{definition}
Hence, $f(q)\in\Q[q]$ is unimodal if and only if
$f=0$ or some integer $i\in \Z$ satisfies
$\ord_{q}(f)\leq i \leq \deg_{q}(f)$ and
$f_{\ord_{q}(f)}\leq f_{\ord_{q}(f)+1}
\leq \dots \leq f_{i} \geq \dots \geq
f_{\deg_{q}(f)-1}\geq f_{\deg_{q}(f)}$.

We then recall the  well-known statement below on palindromic
and unimodal $q$-polynomials.
\begin{proposition}(\cite[Proposition 1]{Sta})
  \label{prop:palindromic-unimodal}
  Let $f(q),g(q)\in \Q_{\geq 0}[q]$ be palindromic and unimodal
  $q$-polynomials.  Then, $f(q)g(q)\in \Q_{\geq 0}[q]$ is a
  palindromic and unimodal $q$-polynomial.
\end{proposition}

The following extends
Proposition~\ref{prop:palindromic-unimodal} for the change
of variables $q\mapsto \qr$ of $\rho\in \Zgeo$.  This is to construct
explicit merged-log-concave parcels for non-trivial base
shifts.  For simplicity, let $a \equiv_{\rho} b$ denote
$a \equiv b \pmod \rho$ for $a,b,\rho\in \Z$.

\begin{proposition}(a relative version of
  Proposition~\ref{prop:palindromic-unimodal}
  on $q$-powers)
  \label{prop:rho-palindromic-unimodal}
  Suppose $\rho\in \Zgeo$.  Let $f(q)\in\Q_{\geq 0}[\qr]$ be a
  palindromic and unimodal $\qr$-polynomial . Also, let
  $h(q)\in \Q_{\geq 0}[q]$ be a palindromic and unimodal
  $q$-polynomial.  Consider
  $g(q)=[\rho]_{q}h(q)$.  Then,
  $f(q)g(q)\in \Q_{\geq 0}[q]$ is a palindromic and unimodal
  $q$-polynomial.
\end{proposition}
\begin{proof}
  First, let us confirm that each product of palindromic
  $q$-polynomials is a palindromic $q$-polynomial.  Thus,
  assume palindromic $q$-polynomials
  $\phi(q),\psi(q)\in \Q[q]$ such that
  $\ord_{q}(\phi)=\ord_{q}(\psi)=0$ for simplicity.  Let
  $\gam(q)=\phi(q)\psi(q) =\sum_{i} \gam_{i}q^{i}$.  It follows that
  $\gam_{i} =\sum_{j_{1}+j_{2}=i} \phi_{j_{1}}\psi_{j_{2}}
  =\sum_{j_{1}+j_{2}=i} \phi_{\deg_{q}(\phi)-j_{1}}\psi_{\deg_{q}(\psi)-j_{2}}
  =\sum_{j_{1}+j_{2}=\deg_{q}(\phi)+\deg_{q}(\psi)-i} \phi_{j_{1}}\psi_{j_{2}}
  =\gam_{\deg_{q}(\phi)+\deg_{q}(\psi)-i}$.
  
  Second, let
  $\lam(q)=[\rho]_{q}=\sum_{0\leq i \leq \rho-1}\lam_{i} q^{i}$ and
  $\mu(q)= f(q)\lam(q) = \sum_{0\leq i\leq \deg_{q}(f)+\rho-1}\mu_{i}q^{i}$.
  Then, let us prove that $\mu(q)$ is a unimodal
  $q$-polynomial, since $\mu(q)$ is a palindromic
  $q$-polynomial by the above.  Assume $\ord_{q}(f)=0$ for
  simplicity.  Hence, let us prove $\mu_{i+1}-\mu_{i}\geq 0$ for
  each $0\leq i \leq \frac{\deg_{q}(f)+\rho-1}{2}-1$.
  
  If $j\neq -1, \rho-1$, then $\lam_{j+1}=\lam_{j}$.  Moreover,
  $ \mu_{i+1}-\mu_{i}
  =\sum_{j}f_{i-j}\lam_{j+1}-\sum_{j}f_{i-j}\lam_{j}
  =\sum_{j}f_{i-j}(\lam_{j+1}-\lam_{j})$.  Hence, we have
  \begin{align}
    \mu_{i+1}-\mu_{i}=f_{i+1}-f_{i+1-\rho}.
    \label{eq:rho-palindromic-unimodal-p1}
  \end{align}
  If $\rho=1$, then $\mu(q)= f(q)$.  Thus, let $\rho>1$.  Suppose
  $\frac{\deg_{q}(f)}{2}-1< i\leq \frac{\deg_{q}(f)+\rho-1}{2}-1$.
  Because $\deg_{q}(f)=\rho k$ for some $k\in \Z$, we have
  $\rho k < 2i+2 \leq \rho (k+1)-1$. Then,
  $i+1\not\equiv_{\rho}0$ by $\rho>1$.  Thus,
  $\mu_{i+1}-\mu_{i}=0$ follows from
  equation~\eqref{eq:rho-palindromic-unimodal-p1}, since
  $f_{i+1}=f_{i+1-\rho}=0$.  If
  $0\leq i\leq \frac{\deg_{q}(f)}{2}-1$, then
  $\mu_{i+1}-\mu_{i}\geq 0$ by
  equation~\eqref{eq:rho-palindromic-unimodal-p1}, since
  $f_{i+1}\geq f_{i+1-\rho}$ for the palindromic and unimodal
  $q$-polynomial $f(q)$.  Thus, $\mu(q)$ is a unimodal
  $q$-polynomial.
  
  Third, the assertion holds by
  Proposition~\ref{prop:palindromic-unimodal}, since
  $\mu(q)$ and $h(q)$ are palindromic
  unimodal $q$-polynomials.
\end{proof}

Let us introduce the following to analyze
 pre-merged determinants.

\begin{definition}\label{def:omega}
  Let $a,b\in \Zt$ and
  $w,\rho\in \Zgeo^{2}$. Then, we define the following.
  \begin{enumerate}
  \item
  $\chi_{a,b}(w,\rho,q)
    ={b \brack a}_{q^{\rhoo}}^{w}[\rhoo]_{q}
    \in \Zgez[q]$.
    \item  $\ka_{a,b}(w,\rho,q)
    =\frac{\prod_{i\in\oi(2)}b_{b_{i}- a_{2},\rhoo}^{w_{i}}(q) }
    {[\rhoo]_{q}}
    \in \Q(q)$.
  \end{enumerate}

\end{definition}

Each $q$-binomial coefficient ${ B \brack A}_{q}$ of
$A,B\in \Z$ is a palindromic and unimodal $q$-polynomial
(see~\cite[Theorem 11]{Sta}).  Thus, we have the
following for
$\chi_{a,b}(w,\rho,q)$
by
Proposition~\ref{prop:rho-palindromic-unimodal}.

\begin{corollary}\label{cor:omega}
  If $b\geq a\geq 0$, then $\chi_{a,b}(w,\rho,q)\llq 0$ is a
  palindromic and unimodal $q$-polynomial.
\end{corollary}
\begin{proof}
  Since $b\geq a\geq 0$,
  ${ b_{1} \brack a_{1}}_{q^{\rhoo}}\llq 0$ and
  ${ b_{2} \brack a_{2}}_{q^{\rhoo}}\llq 0$.  Thus,
  $\chi_{a,b}(w,\rho,q)\llq 0$ by $[\rhoo]_{q}\llq 0$.
  Moreover, ${ b_{1} \brack a_{1}}_{q^{\rhoo}}$ and
  ${ b_{2} \brack a_{2}}_{q^{\rhoo}}$ are palindromic and
  unimodal $q^{\rhoo}$-polynomials.  Hence,
  Proposition~\ref{prop:rho-palindromic-unimodal} implies
  that ${ b \brack a}^{w}_{q^{\rhoo}} $ is a palindromic and
  unimodal $q^{\rhoo}$-polynomial, and that
  $\chi_{a,b}(w, \rho,q)$ is a palindromic and unimodal
  $q$-polynomial.
\end{proof}

Moreover, we have the following  positivity
of  $\ka_{a,b}(w,\rho,q)$.
\begin{lemma}\label{lem:kappa}
  If $a,b\in \Zt$ satisfy $b_{2}>a_{2}$ and
  $b_{1}\geq a_{2}$, then $\ka_{a,b}(w,\rho,q)\llq 0$.
\end{lemma}
\begin{proof}
  By $b_{2}>a_{2}$ and $b_{1}\geq a_{2}$,
  Lemma~\ref{lem:bsingle-qnum} gives
  $\frac{b_{b_{2}-a_{2},\rhoo}^{w_{2}}(q) }{[\rhoo]_{q}}
  =\frac{\prod_{h\in\oi(b_{2}-a_{2})}[\rhoo]^{w_{2}}_{q^{h}}
  }{[\rhoo]_{q}}\llq 0$ and
  $b_{b_{1}-a_{2},\rhoo}^{w_{1}}(q) \llq 0$.  Hence, the
  assertion follows.
\end{proof}

Then, the following positivity holds by width-two pre-merged determinants.

\begin{proposition}\label{prop:width-two-pre-merged}
  Consider a tempered pre-fitting $\omg=\tp{2,a,b}$.  Assume
  a flat $w\in \Zgeo^{2}$, $x=\iota^{2}(q)$, and
  $e\in\oi(0,n_{a}^{b}(2,w))$.  Then, we have
  \begin{align}
    \dt(U)_{a}^{b}(2,w,\rho,e,x)\llq 0.
    \label{ineq:width-two-pre-merged-pos}
  \end{align}
\end{proposition}
\begin{proof}
  Since $\omg$ is tempered pre-fitting, 
  we have
  \begin{align}
    0\leq a_{1}< a_{2}\leq b_{1}<b_{2}.
    \label{ineq:width-two-pre-merged-two-strict}
  \end{align}
  By
  inequalities~\eqref{ineq:width-two-pre-merged-two-strict},
  we prove inequality~\eqref{ineq:width-two-pre-merged-pos}
  in two steps: first, we take trivial base and degree shifts
  $\rhoo=1$ and $e=0$; second, we take general
  $\rhoo\geq 1$ and $e\in\oi(0,n_{a}^{b}(2,w))$.
  
  Let us prove
  inequality~\eqref{ineq:width-two-pre-merged-pos} for
  $\rhoo=1$ and $e=0$.  First, let $w_{1}=1$.  Thus, by the
  induction $b_{2}$, let us prove
  \begin{align}
    \dt(U)_{a}^{b}(2,w,\rho,e,x)=
    {b \brack a}_{q}
    -
    {b \brack a^{\ve}}_{q}
    \llq 0.
    \label{ineq:width-two-pre-merged-triv-baseshift}
  \end{align}
  The smallest possible $b_{2}$ is $2$ by
  inequalities~\eqref{ineq:width-two-pre-merged-two-strict}.
  Then, $a_{1}=0$, $a_{2}=1$, and $b_{1}=1$.  This yields
  inequality~\eqref{ineq:width-two-pre-merged-triv-baseshift} by
  the  direct computation:
  \begin{align}
    \dt(U)_{a}^{b}(2,w,\rho,e,x)=
    {1 \brack 0}_{q}
    {2 \brack 1}_{q}
    -
    {1 \brack 1}_{q}
    {2 \brack 0}_{q}
    =q        \llq 0.
    \label{ineq:width-two-pre-merged-q}
  \end{align}
  
  Let $b_{2}\geq 2$.  By the $q$-Pascal identity, the first
  term of
  inequality~\eqref{ineq:width-two-pre-merged-triv-baseshift}
  gives
  \begin{align}
    {b \brack a}_{q}
    &=
      \left( 
      {b_{1}-1 \brack a_{1}-1}_{q}
      +q^{a_{1}}{b_{1}-1 \brack a_{1}}_{q}
      \right)
      \cdot
      \left( 
      {b_{2}-1 \brack a_{2}-1}_{q}
      +q^{a_{2}}{b_{2}-1 \brack a_{2}}_{q}
      \right)
    \\&       =
    {b_{1}-1 \brack a_{1}-1}_{q}
    {b_{2}-1 \brack a_{2}-1}_{q} 
    +q^{a_{1}}
    {b_{1}-1 \brack a_{1}}_{q}
    {b_{2}-1 \brack a_{2}-1}_{q}
    \\ &          +q^{a_{2}}
         {b_{1}-1 \brack a_{1}-1}_{q}
         {b_{2}-1 \brack a_{2}}_{q}
         +q^{a_{1}+a_{2}}
         {b_{1}-1 \brack a_{1}}_{q}
         {b_{2}-1 \brack a_{2}}_{q}.
  \end{align}
  Similarly,  the second term of
  inequality~\eqref{ineq:width-two-pre-merged-triv-baseshift} gives
  \begin{align}
    {b \brack a^{\ve}}_{q}
    &=
      {b_{1}-1 \brack a_{2}-1}_{q}
      {b_{2}-1 \brack a_{1}-1}_{q} 
      +q^{a_{2}}
      {b_{1}-1 \brack a_{2}
      }_{q}{b_{2}-1 \brack a_{1}-1}_{q}
    \\ &
         +q^{a_{1}}{b_{1}-1 \brack a_{2}-1
         }_{q}{b_{2}-1 \brack a_{1}}_{q}
         +q^{a_{1}+a_{2}}{b_{1}-1 \brack a_{2}
         }_{q}{b_{2}-1 \brack a_{1}}_{q}.
  \end{align}
  Notice that
  inequality~\eqref{ineq:width-two-pre-merged-triv-baseshift}
  is strict.  Thus, comparing the powers of $q$, we prove the
  following:
  \begin{align}
    {b-1 \brack a-1}_{q}
    -
    {b-1 \brack a^{\ve}-1}_{q}
    &\ggq 0;  
      \label{ineq:width-two-pre-merged-four-minuses}
    \\
    {b-1 \brack (a_{1},a_{2}-1)}_{q}
    -
    {b-1 \brack (a_{1},a_{2}-1)^{\ve}}_{q} 
    &\ggq   0;
      \label{ineq:width-two-pre-merged-two-minuses-left}\\
    {b-1 \brack (a_{1}-1,a_{2})}_{q}
    -
    {b-1 \brack (a_{1}-1,a_{2})^{\ve}}_{q}
    &\ggq 0;
      \label{ineq:width-two-pre-merged-two-minuses-right}
    \\
    {b-1 \brack a}_{q}
    -
    {b-1 \brack a^{\ve}}_{q}
    &\llq 0.
      \label{ineq:width-two-pre-merged-two-minuses-center}
  \end{align}
  
  Inequality~\eqref{ineq:width-two-pre-merged-four-minuses}
  holds as follows.  If $a_{1}\geq 1$, then we obtain
  inequality~\eqref{ineq:width-two-pre-merged-four-minuses}
  by the induction.  If $a_{1}=0$, then
  inequality~\eqref{ineq:width-two-pre-merged-four-minuses}
  still holds by $0\ggq 0$. 
  
  Let us prove
  inequality~\eqref{ineq:width-two-pre-merged-two-minuses-left}.
  If $a_{2}-1>a_{1}$, then
  inequality~\eqref{ineq:width-two-pre-merged-two-minuses-left}
  follows from the induction.  When $a_{2}-1=a_{1}$, we
  obtain
  inequality~\eqref{ineq:width-two-pre-merged-two-minuses-left}
  by $0\ggq 0$.

  Assume $b_{1}-1<a_{2}$. Then,
  ${b_{1}-1 \brack a_{2}}_{q}=0$. Hence,
  inequality~\eqref{ineq:width-two-pre-merged-two-minuses-right}
  follows, since ${b-1 \brack \tp{a_{1}-1,a_{2}}}_{q}\ggq 0$.
  Also,
  inequality~\eqref{ineq:width-two-pre-merged-two-minuses-center}
  holds, because ${b-1 \brack a}_{q}\llq 0$ follows from
  $b-1\geq a$ in
  inequalities~\eqref{ineq:width-two-pre-merged-two-strict}.
  
  Hence, suppose $a_{2} \leq b_{1}-1$.  Then,
  inequality~\eqref{ineq:width-two-pre-merged-two-minuses-center}
  follows from the induction.  If $a_{1}\geq 1$, then
  inequality~\eqref{ineq:width-two-pre-merged-two-minuses-right}
  holds by the induction. If $a_{1}=0$, then we still have
  inequality~\eqref{ineq:width-two-pre-merged-two-minuses-right}
  by $0\ggq 0$.  Therefore, we obtain
  inequality~\eqref{ineq:width-two-pre-merged-triv-baseshift}.
  
  Second, if $w_{1}\geq 1$, then since
  inequality~\eqref{ineq:width-two-pre-merged-triv-baseshift}
  implies ${b \brack a}_{q} \llq {b \brack a^{\ve}}_{q}$,
  the $w_{1}$-fold multiplications of both sides give
  inequality~\eqref{ineq:width-two-pre-merged-pos} for
  $\rhoo=1$ and $e=0$.
  
  Let us prove
  inequality~\eqref{ineq:width-two-pre-merged-pos} for
  general $\rhoo\geq 1$ and $e\in\oi(0,n_{a}^{b}(2,w))$.  We have
  \begin{align}
    &\dt(U)_{a}^{b}(2,w,\rho,e,x)
    \\&=
    b_{b_{1}- a_{2},\rhoo}^{w_{1}}(q)
    b_{b_{2}- a_{2},\rhoo}^{w_{1}}(q)\left(
    {b \brack a}^{w}_{q^{\rhoo}}
    \frac{ b_{b_{1}- a_{1},\rhoo}^{w_{1}}(q)
    }{ b_{b_{1}- a_{2},\rhoo}^{w_{1}}(q)}
    -
    q^{e}
    {b \brack a^{\ve}}^{w}_{q^{\rhoo}}
    \frac{ b_{b_{2}- a_{1},\rhoo}^{w_{1}}(q)
    }{ b_{b_{2}- a_{2},\rho}^{w_{1}}(q)}\right)\\
    &=
      \ka_{a,b}(w,\rho,q) \left(
      \chi_{a,b}(w,\rho,q)
      \frac{ b_{b_{1}- a_{1},\rhoo}^{w_{1}}(q)}
      { b_{b_{1}- a_{2},\rhoo}^{w_{1}}(q)}
      -
      q^{e}\chi_{a^{\ve},b}
      (w,\rho,q)
      \frac{ b_{b_{2}- a_{1},\rhoo}^{w_{1}}(q)}
      { b_{b_{2}- a_{2},\rhoo}^{w_{1}}(q)}\right).
      \label{eq:width-two-pre-merged-baseshift-quotients}
  \end{align}
  Also,
  inequalities~\eqref{ineq:width-two-pre-merged-two-strict}
  yield $d=b-a_{1}\in \Zt$ and $k=a_{2}-a_{1}\in \Z$ such that
  \begin{align}
    d &>0, \label{ineq:width-two-pre-merged-d}\\
    k &>0. \label{ineq:width-two-pre-merged-k}
  \end{align}
  Moreover, for $\lam\in \bZ^{k}$, let
  \begin{align}
    E_{1}(\lam,w,d,e,k)
    &= w_{1}\sum_{i\in\oi(k)}\lam_{i}(d_{1}-k+i)\in \Z,\\
    E_{2}(\lam,w,d,e,k)
    &=e + w_{1}\sum_{i\in\oi(k)}\lam_{i}(d_{2}-k+i)\in \Z.
  \end{align}
  Then, since we have $d-k=b-a_{2}\geq 0$ and
  inequality~\eqref{ineq:width-two-pre-merged-k},
  Lemma~\ref{lem:bsingle-quot} gives
  \begin{align}
    \frac{ b_{b_{i}- a_{1},\rhoo}(q)^{w_{1}}}
    { b_{b_{i}- a_{2},\rhoo}(q)^{w_{1}}}
    &=   \frac{b_{d_{i},\rhoo}(q)^{w_{1}}}
      {b_{d_{i}-k,\rhoo}(q)^{w_{1}}}
      =\sum_{\lam\in\bZ^{k}, 0\leq \lam\leq \rhoo-1}
      q^{E_{i}(\lam,w,d,0,k)}
  \end{align}
  for $i\in\oi(2)$.  Hence, assume $\lam\in\Z^{k}$ such that
  \begin{align}
    0\leq \lam\leq \rhoo-1. 
    \label{ineq:width-two-pre-merged-lam}
  \end{align}
  Then, by Lemma~\ref{lem:kappa} and
  equation~\eqref{eq:width-two-pre-merged-baseshift-quotients},
  inequality~\eqref{ineq:width-two-pre-merged-pos} follows
  from
  \begin{align}
    \chi_{a,b}(w,\rho,q)
    q^{ E_{1}(\lam,w,d,e,k)}
    -
    \chi_{a^{\ve},b}(w,\rho,q)
    q^{ E_{2}(\lam,w,d,e,k)}
    \llq 0. 
    \label{ineq:width-two-pre-merged-chi-diff}                  
  \end{align}

  Let us prove
  inequality~\eqref{ineq:width-two-pre-merged-chi-diff}.  We
  have
  $\ord_{q}(\chi_{a,b}(w,\rho,q)) =
  \ord_{q}(\chi_{a^{\ve},b}(w,\rho,q))=0$.  Also,
  Lemma~\ref{lem:deg-diff-width-two} gives
  \begin{align}
    \deg_{q}(\chi_{a,b}(w,\rho,q))
    -\deg_{q}(\chi_{a^{\ve},b}(w,\rho,q))
    &
      =
      \deg_{q}
      \left({b \brack a}_{y^{\rho}}^{w}[\rho]_{q}\right)
      -\deg_{q}
      \left({b \brack
      a^{\ve}}_{y^{\rho}}^{w}[\rho]_{q}\right)
    \\ &
         =
         \deg_{q}
         {b \brack a}_{y^{\rho}}^{w}
         -\deg_{q}
         {b \brack a^{\ve}}_{y^{\rho}}^{w}\\
    &=              n_{a}^{b}(2,w,\rho).
  \end{align}
  Moreover, since
  inequality~\eqref{ineq:width-two-pre-merged-pos} of $e=0$
  gives
  ${b \brack a}^{w}_{q^{\rhoo}} - {b \brack
    a^{\ve}}^{w}_{q^{\rhoo}} >_{q^{\rhoo}}0$, multiplying
  $[\rhoo]_{q}$ on both sides gives
  $ \chi_{a,b}(w,\rho,q) - \chi_{a^{\ve},b}(w,\rho,q) \llq 0$.
  Hence, since $\chi_{a,b}(w,\rho,q)$ and
  $\chi_{a^{\ve},b}(w,\rho,q) $ are palindromic and unimodal 
  $q$-polynomials by Corollary~\ref{cor:omega},
  inequality~\eqref{ineq:width-two-pre-merged-chi-diff}
  follows from  
  \begin{align}
    0\leq  E_{2}(\lam,w,d,e,k) - E_{1}(\lam,w,d,e,k)
    \leq n_{a}^{b}(2,w,\rho).
    \label{ineq:width-two-pre-merged-E-diff}
  \end{align}
  
  Let us prove
  inequalities~\eqref{ineq:width-two-pre-merged-E-diff}.
  Because
  $0\leq e\leq
    n_{a}^{b}(2,w)=w_{1}(b_{2}-b_{1})(a_{2}-a_{1})$,
    inequalities~\eqref{ineq:width-two-pre-merged-d},
  ~\eqref{ineq:width-two-pre-merged-k},
  and~\eqref{ineq:width-two-pre-merged-lam} 
  give
  the right-hand side
  of inequalities~\eqref{ineq:width-two-pre-merged-E-diff}:
  \begin{align}
    E_{2}(\lam,w,d,e,k)
    - E_{1}(\lam, w,d,e,k)
    &=e+w_{1}(d_{2}-d_{1})\sum \lam_{i}
    \\&
    \leq e+w_{1}(d_{2}-d_{1})k (\rhoo -1)
    \\&
    \leq w_{1}\rhoo(b_{2}-b_{1}) (a_{2}-a_{1})
    \\&=n_{a}^{b}(2,w,\rho).
  \end{align}
  Thus, since $d_{2}-d_{1}=b_{2}-b_{1}>0$ and $\lam\geq 0$
  imply the left-hand side,
  inequalities~\eqref{ineq:width-two-pre-merged-E-diff}
  follow.  Then, the assertion holds.
\end{proof}

We put the following remark, specializing
parameters such as the base shift $\rho$.

\begin{remark}
  Suppose $w_{1}=\rhoo=1$ and $e=0$.
  Then,
for integers $0\leq a_{1}<a_{2}\leq b_{1}<b_{2}$,
Proposition~\ref{prop:width-two-pre-merged}
restricts to
  \begin{align}
    \dt(U)_{a}^{b}(2,w,\rho,e,x)=
    {b_{1} \brack a_{1}}_{q}
    {b_{2} \brack a_{2}}_{q}
    -
    {b_{1} \brack a_{2}}_{q}
    {b_{2} \brack a_{1}}_{q}
    \llq 0. \label{ineq:simplest}
  \end{align}
  This does not claim the $q$-log-concavity of
  $q$-binomial coefficients.  Thus, this differs from
  the pioneering
  results~\cite{But,Kra,Sag} on the $q$-log-concavity.
  
  Let us explain more precisely.  As in~\cite[Corollary
  3]{Kra}, for non-negative integers $b_{1}\leq b_{2}$,
  $a_{1}\leq a_{2}$, $k$, and
  $\lam \leq d(a,b,k):=k(2((a_{2}-a_{1})+k)+(b_{2}-b_{1}))$,
  they have obtained
  \begin{align}
    T(a,b,k,\rho,q):=
    {b_{1} \brack a_{2}}_{q}
    {b_{2} \brack a_{1}}_{q}
    -
    q^{\lam}{b_{1} \brack a_{2}+k}_{q}
    {b_{2} \brack a_{1}-k}_{q}
    \ggq0. \label{ineq:BSK}
  \end{align}
  
  Let us confirm that $T(a,b,k,\rho,q)$ is not a monomial,
  if $\deg_{q}T(a,b,k,\rho,q)>0$.  First, if
  ${b_{1} \brack a_{2}}_{q} {b_{2} \brack a_{1}}_{q} {b_{1}
    \brack a_{2}+k}_{q} {b_{2} \brack a_{1}-k}_{q}\neq 0$,
  then,
  $\deg_{q}\left({b_{1} \brack a_{2}}_{q} {b_{2} \brack
      a_{1}}_{q}\right) - \deg_{q}\left({b_{1} \brack
      a_{2}+k}_{q} {b_{2} \brack a_{1}-k}_{q}\right)
  =d(a,b,k)>1$.
  Also, ${b_{1} \brack a_{2}}_{q} {b_{2} \brack a_{1}}_{q}$
  and
  ${b_{1} \brack a_{2}+k}_{q} {b_{2} \brack a_{1}-k}_{q}$
  are palindromic and unimodal $q$-polynomial with
  $ \ord_{q}\left({b_{1} \brack a_{2}}_{q} {b_{2} \brack
      a_{1}}_{q}\right)= \ord_{q}\left({b_{1} \brack
      a_{2}+k}_{q} {b_{2} \brack a_{1}-k}_{q}\right)=0$.
  Hence, $T(a,b,k,\rho,q)$ is not a monomial.  Second, suppose
  ${b_{1} \brack a_{2}}_{q} {b_{2} \brack a_{1}}_{q} {b_{1}
    \brack a_{2}+k}_{q} {b_{2} \brack a_{1}-k}_{q}= 0$.
  Then,
  ${b_{1} \brack a_{2}}_{q} {b_{2} \brack a_{1}}_{q}=0$
  implies
  $ {b_{1} \brack a_{2}+k}_{q} {b_{2} \brack a_{1}-k}_{q}=0$
  by $d(a,b,k)\geq 0$.  Also,
  ${b_{1} \brack a_{2}+k}_{q} {b_{2} \brack a_{1}-k}_{q}=0$
  implies that $T(a,b,k,\rho,q)$ is not a monomial when
  $\deg_{q}T(a,b,k,\rho,q)>0$, as $T(a,b,k,\rho,q)$ is a
  palindromic and unimodal $q$-polynomial.
  
  However, $\dt(U)_{a}^{b}(2,w,\rho,e,x)=q$ in
  inequality~\eqref{ineq:width-two-pre-merged-q}.
  Furthermore,
  unlike inequality~\eqref{ineq:BSK},
  inequality~\eqref{ineq:simplest} is of $\llq 0$. This is
  important for us to obtain polynomials with positive
  integer coefficients.
\end{remark}

We extend Proposition~\ref{prop:width-two-pre-merged} after
the following lemma.

\begin{lemma}\label{lem:pre-merged-term-nonneg-pos}
  Consider a pre-parcel $U_{a}^{b}(\mu,w,\rho,x)$ with
  $\lam=\floor{\frac{\mu}{2}}$ and $e\in \Zgez^{\lam}$.
  Let $x\in \Q(\fX)^{\mu}$ be $\scc$-admissible.  Then, we
  have the following.
  \begin{enumerate}
  \item \label{c:pre-merged-term-pos}
    $U_{a}^{b}(\mu,w,\rho,x)\gx 0$ if and only if
    $b\geq a\geq 0$.
  \item \label{c:pre-merged-term-nonneg}
    $U_{a}^{b}(\mu,w,\rho,x)\gex 0$ and
    $U_{a}^{b}(\mu,w,\rho,x)\cdot \prod x(1,\lam)^{e} \gex 0$.
  \end{enumerate}
\end{lemma}
\begin{proof}
  Let us prove Claim~\ref{c:pre-merged-term-pos}.  If
  $b\geq a\geq 0$ is not true, then
  ${b \brack a}^{w}_{x^{\rho}}=0$.  If $b\geq a\geq 0$, then
  ${b \brack a}^{w}_{x^{\rho}}\gx 0$ and
  $B(\mu,w,b-a,\rho,x)\gx 0$ by
  Lemma~\ref{lem:bsingle-qnum},
  since $\gx$ has the
  multiplicative property by Theorem~\ref{thm:adm}.  
  Thus,
  Claim~\ref{c:pre-merged-term-pos} follows.
  
  Claim~\ref{c:pre-merged-term-nonneg} holds by
  $\prod x(1,\lam)^{e}\gex 0$.
\end{proof}


\begin{theorem}\label{thm:pre-merged}
  Suppose a pre-fitting $\omg=\tp{\mu,a,b}$.  Consider a
  $\scc$-admissible  $x\in \Q(\fX)^{\mu}$ 
  and
  $w,\rho\in \Zgeo^{\mu}$ such that $x$, $w$, and
  $\rho$ are palindromic.
  Let
  $\lam=\floor{\frac{\mu}{2}}$ and
  $e\in \Z^{\lam}$ such that
  \begin{align}
    0\leq e\leq N_{a}^{b}(\mu,w).
    \label{ineq:pre-merged-shift}
  \end{align}
  Then, we obtain
  \begin{align}
    d(U)_{a}^{b}(\mu,w,\rho,e,x)    \gx 0.
    \label{ineq:pre-merged-pos}
  \end{align}
\end{theorem}
\begin{proof}
  Suppose that $\omg$ is tempered.
  Claim~\ref{c:pre-merged-term-pos} of
  Lemma~\ref{lem:pre-merged-term-nonneg-pos} implies that
  $U_{a^{\ve}}^{b}(\mu,w,\rho,x)=0$ as $b\not\geq a^{\ve}$, and
  that
  $d(U)_{a}^{b}(\mu,w,\rho,e,x) =U_{a}^{b}(\mu,w,\rho,x)\gx 0$ as
  $b\geq a\geq 0$ for the pre-fitting $\omg$.
  
  Hence, assume that $\omg$ is tempered.  Let us prove
  inequality~\eqref{ineq:pre-merged-pos} by the induction on
  $\mu$.  When $\mu=2$,
  Proposition~\ref{prop:width-two-pre-merged} gives
  inequality~\eqref{ineq:pre-merged-pos}.  Suppose
  $\mu\geq 3$.  Then, inequalities~\eqref{ineq:pre-merged-shift}
  imply
  $0\leq \out(\mu,e)\leq N_{\out(a)}^{\out(b)}(\mu-\chi(\mu),\out(w))$.
  Also, $\tout(\omg)$ is tempered and pre-fitting by
  Lemma~\ref{lem:center-outer}. Then, the induction implies
  \begin{align}
    d(U)_{\out(a)}^{\out(b)}(\mu-\chi(\mu), \out(w), \out(\rho),
    \out(\mu,e), \out(x)) \gx 0.
    \label{ineq:pre-merged-outer-pos}
  \end{align}
  Moreover, since $b\geq a\geq 0$ implies
  $\cen(b)\geq \cen(a)\geq 0$, Claim~\ref{c:pre-merged-term-pos}
  of Lemma~\ref{lem:pre-merged-term-nonneg-pos} gives
  \begin{align}
    U_{\cen(a)}^{\cen(b)}(\chi(\mu),\cen(w),\cen(\rho),\cen(x)) \gx 0.
    \label{ineq:pre-merged-center-pos}
  \end{align}
  In particular, if $\mu$ is odd, then
  inequalities~\eqref{ineq:pre-merged-outer-pos}
  and~\eqref{ineq:pre-merged-center-pos} give
  inequality~\eqref{ineq:pre-merged-pos} by
  Claim~\ref{c:outer-preparcel-premerged-premerged} of
  Lemma~\ref{lem:outer-preparcel-premerged}.
  
  Hence, let $\mu$ be even.  First, assume that $\cen(a)$ or
  $\cen(b)$ is flat.  Then,
  $0\leq \cen(\mu,e)_{1}=e_{\lam}
  \leq
  \rho_{\lam}w_{\lam}(b_{\lam+1}-b_{\lam})
  (a_{\lam+1}-a_{\lam})=0$.
  Also, since $\cen(x)$, $\cen(w)$, and $\cen(\rho)$ are flat,
  \begin{align}
    &{\cen(b) \brack \cen(a)}^{\cen(w)}_{\cen(x)^{\cen(\rho)}}
      b_{b_{\lam}-a_{\lam},\rho_{\lam}}^{w_{\lam}}(x_{\lam})
      b_{b_{\lam+1}-a_{\lam+1},\rho_{\lam+1}}^{w_{\lam+1}}
      (x_{\lam+1})
    \\&
    =
    {\cen(b) \brack \cen(a)^{\ve}}^{\cen(w)}_{\cen(x)^{\cen(\rho)}}
    b_{b_{\lam}-a_{\lam+1},\rho_{\lam}}^{w_{\lam}}(x_{\lam})
    b_{b_{\lam+1}-a_{\lam},\rho_{\lam+1}}^{w_{\lam+1}}(x_{\lam+1}).
  \end{align}
  Thus,
  inequality~\eqref{ineq:pre-merged-center-pos} implies
  \begin{align}
    \prod x(\lam,\lam)^{\cen(\mu,e)}\cdot
    U_{\cen(a)^{\ve}}^{\cen(b)}(\chi(\mu),\cen(w),\cen(\rho),\cen(x)) =
    U_{\cen(a)}^{\cen(b)}(\chi(\mu),\cen(w),\cen(\rho),\cen(x))\gx 0.
  \end{align}
  Therefore,
  inequality~\eqref{ineq:pre-merged-pos} follows from
  inequality~\eqref{ineq:pre-merged-outer-pos} and
  Claim~\ref{c:outer-preparcel-premerged-premerged} of
  Lemma~\ref{lem:outer-preparcel-premerged}.
  
  Second, suppose
  $0\leq \cen(a)_{1}<\cen(a)_{2}\leq \cen(b)_{1}<\cen(b)_{2}$.
  Then, because
  $0\leq \cen(\mu,e)\leq N_{\cen(a)}^{\cen(b)}(\chi(\mu),\cen(w))$,
  the induction gives
  $
  d(U)_{\cen(a)}^{\cen(b)}
  (\chi(\mu),\cen(w),\cen(\rho),\cen(\mu,e),\cen(x))\gx0$.
  Thus, we have
  \begin{dmath}
    U_{\cen(a)}^{\cen(b)}(\chi(\mu),\cen(w),\cen(\rho),\cen(x))
    \gx
    \prod x(\lam,\lam)^{\cen(\mu,e)}\cdot
    U_{\cen(a)^{\ve}}^{\cen(b)}(\chi(\mu),\cen(w),\cen(\rho),\cen(x)).
    \label{ineq:pre-merged-center-merged}
  \end{dmath}
  Moreover, inequality~\eqref{ineq:pre-merged-center-pos}
  implies
  \begin{dmath}
    U_{\out(a)}^{\out(b)}(\mu-\chi(\mu),\out(w),\out(\rho),\out(x))
    \gx
    \prod x(1,\lam-1)^{\out(\mu,e)}\cdot
    U_{\out(a)^{\ve}}^{\out(b)}
    (\mu-\chi(\mu),\out(w),\out(\rho),\out(x)).
    \label{ineq:pre-merged-outer-merged}
  \end{dmath}
  Therefore, Claim~\ref{c:succ-scc-sce-mult} of
  Proposition~\ref{prop:succ},
  Claim~\ref{c:outer-preparcel-premerged-premerged} of
  Lemma~\ref{lem:outer-preparcel-premerged}, and
  Claim~\ref{c:pre-merged-term-nonneg} of
  Lemma~\ref{lem:pre-merged-term-nonneg-pos} give
  inequality~\eqref{ineq:pre-merged-pos} from
  inequalities~\eqref{ineq:pre-merged-center-merged}
  and~\eqref{ineq:pre-merged-outer-merged}.
\end{proof}
In particular, if $x=\iota^{\mu}(q)$, then
Theorem~\ref{thm:pre-merged} gives
$d(U)_{a}^{b}(\mu,w,\rho,e,x)\llq 0$.

\subsection{Merged-log-concavity by
  identities}\label{sec:monom-identity}
We obtain explicit merged-log-concave parcels by identity
functions.  This uses
the following variants of pre-parcels and pre-merged
determinants
with {\it proper mediators}.

\begin{definition}\label{def:quasi-parcel-merged-proper-med}
  Suppose a gate $s\geq 0$, $l\in \Zgeo$, and
  $w,\rho\in \Zgeo^{l}$.  Let
  $\mu=\tp{s,l,w,\scc,\rho,x,\fX}$ for a $\scc$-admissible
  $x\in \Q(\fX)^{l}$.  Consider a $\mu$-mediator $\phi$.
  \begin{enumerate}
  \item We define the quasi-parcel
    $V(s,l,w,\phi,\rho,x)
    =\{V_{a}^{b}(s,l,w,\phi,\rho,x)\in \Q(\fX)\}_{
      a,b\in \Ztl}$ such that
    \begin{align}
      V_{a}^{b}(s,l,w,\phi,\rho,x)
      &=              { b \brack a}^{w^{\wcn}}_{(x^{\rho})^{\wcn}}
        B(s,2l,w^{\wcn},b-a,\phi^{\wcn},\rho^{\wcn},x^{\wcn},\fX).
    \end{align}
    We refer to $s$, $l$, $w$, $\phi$, $\rho$, and $x$ as
    the gate, width, weight, mediator, base shift, and base
    of the quasi-parcel $V(s,l,w,\phi,\rho,x)$.
  \item
    \label{d:quasi-parcel-merged-proper-med-merged}
    Let $a,b\in \Ztl$ and $e\in \Zlgez$.  We define the
    quasi-merged determinant
    \begin{align}
      d(V)_{a}^{b}(s,l,w,\phi,\rho,e,x)
      =
      V_{a}^{b}(s,l,w,\phi,\rho,x)
      -
      V_{a^{\ve}}^{b}(s,l,w,\phi,\rho,x)
      \cdot \prod x^{e}\in \Q(\fX).
    \end{align}
    We call $e$ degree shift of the quasi-merged determinant
    $d(V)_{a}^{b}(s,l,w,\phi,\rho,e,x)$.
  \item We call $\phi$ proper $\mu$-mediator (or proper
    mediator for short), if
    \begin{align}
      d(V)_{a}^{b}(s,l,w,\phi,\rho,e,x)\scc 0
    \end{align}
    for each fitting $\tp{s,l,m,n,k}$, $a=\nu(k)$,
    $b=\nu(m,n,k)$, and $e\in \Zl$ such that
    $0\leq e\leq N_{a}^{b}(2l,v)$.
  \item If  $x$ is flat and
    $e\in \Zgez$, then let
    \begin{align}
      \dt(V)_{a}^{b}(s,l,w,\phi,\rho,e,x)=
      d(V)_{a}^{b}(s,l,w,\phi,\rho,\tp{e}\ccn \iota^{l-1}(0),x).
    \end{align}
  \end{enumerate}
\end{definition}

Let us compare pre-merged and quasi-merged
determinants.
\begin{lemma}\label{lem:pre-merged-quasi-merged}
  Let $s=\tp{0,\infty}$, $l\in \Zgeo$, and $e\in \Zlgez$.  Consider a
  $\scc$-admissible $x\in \Q(\fX)^{l}$ and the canonical
  $l$-mediator $\phi$.  Then, 
  $    d(V)_{a}^{b}(s,l,w,\phi,\rho,e,x)
    =d(U)_{a}^{b}(2l,w^{\wcn},\rho^{\wcn},e,x^{\wcn})$.
  \end{lemma}
\begin{proof}
  Since $\phi(x)$ is canonical,
  Item~\ref{d:base-pre-parcel-merged-base} of
  Definition~\ref{def:base-pre-parcel-merged} gives
  \begin{align}
    B(2l,w^{\wcn},b-a,\rho^{\wcn},x^{\wcn})
    &=
      B(s,l,w,(b-a)(1,l),((b-a)(l+1,2l))^{\ve},\phi,\rho,x,\fX),\\
    B(2l,w^{\wcn},b-a^{\ve},\rho^{\wcn},x^{\wcn})
    &=   B(s,l,w,(b-a^{\ve})(1,l),
      ((b-a^{\ve})(l+1,2l))^{\ve},\phi,\rho,x,\fX).
  \end{align}
  Thus, the assertion follows, because
  \begin{align}
    U_{a}^{b}(2l,w^{\wcn},\rho^{\wcn},x^{\wcn})
    &=
      { b \brack a}^{w^{\wcn}}_{(x^{\rho})^{\wcn}}
      B(2l,w^{\wcn},b-a,\rho^{\wcn},x^{\wcn})
      =V_{a}^{b}(s,l,w,\phi,\rho,x),\\
    U_{a^{\ve}}^{b}(2l,w^{\wcn},\rho^{\wcn}
    ,x^{\wcn})
    &=
      { b \brack a^{\ve}}^{w^{\wcn}}_{(x^{\rho})^{\wcn}}
      B(2l,w^{\wcn},b-a^{\ve},\rho^{\wcn},
      x^{\wcn})
      =V_{a^{\ve}}^{b}(s,l,w,\phi,\rho,x).
  \end{align}  
\end{proof}
We confirm the existence of proper mediators.
\begin{proposition}\label{prop:proper-med}
  Suppose $s_{1}=\tp{0,\infty}$, $l\in \Zgeo$, and
  $w,\rho\in \Zgeo^{l}$.  Let
  $\muo=\tp{s_{1},l,w,\scc,\rho,x,\fX}$ for a $\scc$-admissible
  $x\in \Q(\fX)^{l}$.  Let
  $\mut=\tp{s_{2},l,w,\scc,\rho,x,\fX}$ for a gate
  $s_{2}\geq 0$.  Consider a $\muo$-mediator $\phi$.
  \begin{enumerate}
  \item If $\phi$ is the canonical $l$-mediator, then
    $\phi$ is a proper $\muo$-mediator.
    \label{c:proper-med-canonical-proper}
  \item If $\rho=\iota^{l}(1)$, then $\phi$ is a proper
    $\muo$-mediator.
    \label{c:proper-med-triv-shift}
  \item If $\phi$ is a proper $\muo$-mediator, then $\phi$
    is a proper $\mut$-mediator.
    \label{c:proper-med-proper-proper}
  \end{enumerate}
\end{proposition}
\begin{proof}
  Let us prove Claim~\ref{c:proper-med-canonical-proper}.
  Consider a fitting $\tp{s_{1},l,m,n,k}$, $a=\nu(k)$,
  $b=\nu(m,n,k)$, and $e\in \Zl$ such that
  $0\leq e\leq N_{a}^{b}(2l,v)$.  Then, $\tp{2l,a,b}$ is pre-fitting
  by Claim~\ref{c:fitting-pre-fitting-fitting} of
  Proposition~\ref{prop:fitting-pre-fitting}.  Thus, $\phi$ is
  proper, since Theorem~\ref{thm:pre-merged} and
  Lemma~\ref{lem:pre-merged-quasi-merged} give
  $d(V)_{a}^{b}(s_{1},l,w,\phi,\rho,e,x)
  =d(U)_{a}^{b}(2l,w^{\wcn},\rho^{\wcn},e,x^{\wcn}) \scc 0$.
  
  Let us prove Claim~\ref{c:proper-med-triv-shift}.  When
  $\rho=\iota^{l}(1)$, Lemma~\ref{lem:bshift-fun-sp} implies
  \begin{align}
    B(s_{1},2l,w^{\wcn},b-a,\phi^{\wcn},\rho^{\wcn},x^{\wcn},\fX)
    =B(s_{1},2l,w^{\wcn},b-a^{\ve},\phi^{\wcn},\rho^{\wcn},
    x^{\wcn},\fX)
    =1.
  \end{align}
  Thus, Claim~\ref{c:proper-med-triv-shift} follows from
  Theorem~\ref{thm:pre-merged}.
  
  For a fitting $\tp{s_{2},l,m,n,k}$ with $a=\nu(k)$ and
  $b=\nu(m,n,k)$, Claim~\ref{c:proper-med-proper-proper}
  follows from
  \begin{align}
    d(V)_{a}^{b}(s_{2},l,w,\phi,\rho,e,x)
    =
    \begin{dcases}
      d(V)_{a}^{b}(s_{1},l,w,\phi,\rho,e,x)\scc 0 \mif b-a^{\ve}\ldZtl s_{2},\\
      { b \brack a}^{w^{\wcn}}_{(x^{\rho})^{\wcn}}
      B(s_{2},2l,w^{\wcn},b-a,\phi^{\wcn},\rho^{\wcn},x^{\wcn},\fX)
      \scc 0 \melse.
    \end{dcases}
  \end{align}
\end{proof}

Let us use the following notation of identity functions.
\begin{definition}\label{def:identity}
  Suppose a gate $s\geq 0$ and $l\in \Zgeo$.  We write the
  family $1_{s,l}=\{1_{s,l,m}\in \Q(\fX)\}_{m\in \Zl}$ such
  that
  \begin{align}
    1_{s,l,m}
    =
    \begin{dcases}
      1 \mfor m\ldZl s,\\
      0 \melse.
    \end{dcases}
  \end{align}
  We call $1_{s,l}$ $\tp{s,l}$-identity (or identity for
  simplicity).
\end{definition}

Let $r\in \Q$.  Then,
$r 1_{s,l}=\{r 1_{s,l,m}\in \Q(\fX)\}_{m\in \Zl}$ is pairwise
$\tp{s,l,\scc}$-positive if and only if $r^{2}\scc0$, because
$(r1_{s,l,m})^{2}\scc 0$  for some $m\ldZl s$
when $r 1_{s,l}$ is pairwise $\tp{s,l,\scc}$-positive.  Hence,
we introduce {\it constant parcels} as follows.
\begin{definition}\label{def:constant-parcel}
  For $r\in \Q$ such that $r^{2}\scc 0$,  we call
  a parcel $\Lam(s,l,w,\scc,r 1_{s,l},\rho,\phi,x,\fX)$ constant
  parcel.
\end{definition}

In particular,
Claim~\ref{c:adm-succ-half-gx-scc} of
Lemma~\ref{lem:adm-succ} gives
$1\scc 0$ by $1\gx 0$.
Hence,  $\Lam(s,l,w,\scc, 1_{s,l},\rho,\phi,x,\fX)$
is a constant parcel.

We obtain 
merged-log-concave constant parcels on arbitrary 
gates, widths, positive weights, and base shifts after
 the following lemma.
\begin{lemma}\label{lem:merged-binom-shift-w-inverse}
  Suppose a parcel
  $\cF=\Lam(s,l,w,\scc, f_{s},\rho,\phi,x,\fX)$.  Assume a
  fitting $\tp{s,l,m,n,k}$ with $a=\nu(k)$ and
  $b=\nu(m,n,k)$. Let $y=x^{\rho}$.  Then, we have
  \begin{align}
    &
      (f_{s,m}f_{s,n^{\ve}})^{-1}\Delta(\cF)(s,l,w,m,n,k,\phi,\rho,x,\fX)
    \\ &=
         B(s,l,w,(b-a)(1,l),((b-a)(l+1,2l))^{\ve},
         \phi,\rho,x,\fX)
         { b \brack a}^{w^{\wcn}}_{y^{\wcn}}
    \\&-
    (f_{s,m}f_{s,n^{\ve}})^{-1}
    f_{s,m\bom k }f_{s,(n\bop k)^{\ve}}
    B(s,l,w,(b-a^{\ve})(1,l),
    ((b-a^{\ve})(l+1,2l))^{\ve},\phi,\rho,x,\fX)
    { b \brack a^{\ve}}^{w^{\wcn}}_{y^{\wcn}}.
  \end{align}
\end{lemma}
\begin{proof}
  Claim~\ref{c:merged-binom-bshift-left-right-equations} of
  Theorem~\ref{thm:merged-binom-bshift} and
  Lemma~\ref{lem:ab-define} yield the assertion,
  since $f_{s,m}$ and $f_{s,n^{\ve}}$ are invertible by  the pairwise positivity.
\end{proof}

\begin{proposition}\label{prop:constant-parcel}
  Consider a constant parcel
  $\cF=\Lam(s,l,w,\scc,r 1_{s,l},\rho,\phi,x,\fX)$
  with   $w\in \Zgeo^{l}$ and a proper mediator $\phi$. Then, $\cF$
  is
  $\scc$-merge-log-concave.
\end{proposition}
\begin{proof}
  Let $r=1$  for simplicity.  For a
  fitting $\tp{s,l,m,n,k}$ with $a=\nu(k)$ and $b=\nu(m,n,k)$,
  Lemma~\ref{lem:merged-binom-shift-w-inverse}
  gives
  $\Delta(\cF)(s,l,w,m,n,k,\phi,\rho,x,\fX)
  =
  d(V)_{a}^{b}(s,l,w,\phi,\rho,\iota^{l}(0),x)$.
  Hence, the statement holds by the properness of $\phi$.
\end{proof}

\begin{example}\label{ex:constant-parcel}
  If
  $\cF =\Lam(s,l,w,\scc,1_{s,l}, \rho,x,\fX)$ with
  $w\in \Zgeo^{l}$, then
    \begin{align}
    \cF_{m}=
    \begin{dcases}
      \frac{1}{(m)_{x}^{w}} \mfor
      m\ldZl s,\\
      0 \melse.
    \end{dcases}
  \end{align}
  Also, Claim~\ref{c:proper-med-canonical-proper} of
  Proposition~\ref{prop:proper-med} and
  Proposition~\ref{prop:constant-parcel} imply that $\cF$ is
  $\scc$-merged-log-concave.
\end{example}

\subsection{Merged-log-concavity
  by functional monomial indices}
\label{sec:merged-functional-monom-proper-med}

We first define {\it $\sig$-difference functions} by $\sig$-plus
and minus. Then, we define {\it functional monomial indices}
by the $\sig$-difference functions to give more
merged-log-concave parcels.

\begin{definition}\label{def:functional-monom-index}
  Assume a gate $s\geq 0$, $l\in \Zgeo$, and
  $w\in \Zgeo^{l}$.  Suppose a function $t:\Zl\to \Q$.
  \begin{enumerate}
  \item For each $m,n\in \Zl$ and $k\in \Ztl$, we define the
    $\sig$-difference function
    \begin{align}
      t_{\Del}(m,n,k)=t(m\bom k)+
      t((n\bop k)^{\ve})-t(m)-t(n^{\ve})\in \Q.
    \end{align}
  \item We call the tuple $\ups=\tp{s,l,w,t}$ functional monomial index, if
    \begin{align}
      t_{\Del}(m,n,k)\in \Z,
      \label{cond:functional-monom-index-integer}\\
      0\leq t_{\Del}(m,n,k)\leq n_{a}^{b}(2l,w)
      \label{cond:functional-monom-index-sum}
    \end{align}
    for each wrapped fitting $\tp{s,l,m,n,k}$ with $a=\nu(k)$ and
    $b=\nu(m,n,k)$.
    
    We call $s,l,w$, and $t$ gate, width,
    weight, and core function of $\ups$.  We refer to
    \eqref{cond:functional-monom-index-integer}
    and~\eqref{cond:functional-monom-index-sum} as the
    integer monomial condition and the sum monomial
    condition of $\ups$.
  \end{enumerate}
\end{definition}

Suppose $x=\iota^{l}(q)$.  Then,
Proposition~\ref{prop:constant-parcel} gives
merged-log-concave constant parcels for functional monomial
indices $\tp{s,l,w,t}$ of the zero function $t$.  Moreover,
the following extends Proposition~\ref{prop:constant-parcel}
for general functional monomial indices.

\begin{proposition}\label{prop:functional-merged}
  Let $\ups=\tp{s,l,w,t}$ be a functional monomial index.
  Assume a $\scc$-admissible $q\in \Q(\fX)$ with
  $x=\iota^{l}(q)$.  For $\mu=\tp{s,l,w,\llq,\rho,x,\fX}$, suppose a
  proper $\mu$-mediator $\phi$.  Consider a parcel
  $\cF=\Lam(s,l,w,\scc,\fs,\phi,\rho,x,\fX)$ such that
  \begin{align}
    f_{s,m}=
    \begin{dcases}
      q^{t(m)}  \mfor m\ldZl s,\\
      0 \melse.
    \end{dcases}
  \end{align}
  Also, suppose a fitting $\tp{s,l,m,n,k}$ with $a=\nu(k)$ and
  $b=\nu(m,n,k)$.
  Then, we have
  \begin{align}
    q^{-t(m)-t(n^{\ve})}
    \Delta(\cF)(s,l,w,m,n,k,\phi,\rho,x,\fX)
    &=\dt(V)_{a}^{b}(s,l,w,\phi,\rho,t_{\Del}(m,n,k),x)
    \label{eq:functional-merged-quasi-merged}
    \\&\llq 0.
    \label{ineq:functional-merged-pos}
  \end{align}
  In particular, $\cF$ is
  $\scc$-merged-log-concave.
\end{proposition}
\begin{proof}
  Lemma~\ref{lem:merged-binom-shift-w-inverse}
  gives
  equation~\eqref{eq:functional-merged-quasi-merged},
  since $(f_{s,m}f_{s,n^{\ve}})^{-1}
  f_{s,m\bom k}f_{s,(n\bop k)^{\ve}}=
  q^{  t_{\Del}(m,n,k)}$.
  
  Let us prove
  inequality~\eqref{ineq:functional-merged-pos}.  The
  monomial conditions of $\ups$ give $e(t,m,n,k)\in \Zl$ such
  that $\sum e(t,m,n,k)=t_{\Del}(m,n,k)$ and
  $0\leq e(t,m,n,k)\leq N_{a}^{b}(2l,w)$.  Then, since
  $x=\iota^{l}(q)$ and $\phi$ is proper,
  inequality~\eqref{ineq:functional-merged-pos} follows from
  \begin{align}
    \dt(V)_{a}^{b}(s,l,w,\phi,\rho,t_{\Del}(m,n,k),x,\fX)
    =d(V)_{a}^{b}(s,l,w,\phi,\rho,e(t,m,n,k),x,\fX)
    \llq 0.
  \end{align}
  
  Moreover, since $f_{s}$ is pairwise
  $\tp{s,l,\scc}$-positive,
  $f_{s,m}f_{s,n^{\ve}}=q^{t(m)+t(n^{\ve})}\scc 0$.  This
  gives the latter statement by
  inequality~\eqref{ineq:functional-merged-pos} and the half
  $\llq$-$\scc$ implication.
\end{proof}

\subsection{Monomial indices and  functional monomial indices}
\label{sec:equiv-monomial-indices}
In general,
we do not know how to explicitly give the core function $t$
of a functional monomial index $\tp{s,l,w,t}$.
However, $3l$ rational numbers specify the core $\gam$ of a
monomial index $\tp{l,w,\gam}$ in
Definition~\ref{def:monomial-index}.  Also, cut operators
give merged-log-concave parcels of any gates from ones of
infinite gates by Proposition~\ref{prop:cut-op}.  Thus, we
look for an equivalence between monomial indices $\tp{l,w,\gam}$
and some functional monomial indices $\tp{s,l,w,t}$ on infinite
gates $s$.  This yields explicit merged-log-concave parcels
by monomial indices.

We use the following wrapped fitting tuples.

\begin{lemma}\label{lem:wrapped-fitting-prescribed}
  Consider an infinite gate $s\geq 0$ and $l\in \Zgeo$.
  For each $k\in \Z^{2l}$ and  $m\in \Zl$,
  let $\mu_{m,k}=\tp{s,l,m,m,k}$.
  \begin{enumerate}
  \item \label{c:wrapped-fitting-prescribed-ab-pos} Suppose
    $k\in \Zgeo^{2l}$ and flat $m\in \Zl$ with $a=\nu(k)$
    and $b=\nu(m,m,k)$.  Then, $b_{j}-b_{i}=a_{j}-a_{i}>0$
    for each $1\leq i<j \leq 2l$.
  \item
    \label{c:wrapped-fitting-prescribed-flat-fitting-wrapped}
    Let $k\in \Zgeo^{2l}$ and flat $m\in\Z^{l}$.
    \begin{enumerate}
    \item If $m \ldZl s$, then $\mu_{m,k}$ is fitting.
      \label{c:wrapped-fitting-prescribed-flat-fitting}
    \item If $m\geq \sig(k)_{1}+s_{1}$, then $\mu_{m,k}$
      is wrapped and fitting.
      \label{c:wrapped-fitting-prescribed-flat-wrapped}
    \end{enumerate}
  \item Suppose $r\in \Zgeo^{l}$ and $\lam\in \Zl$. 
    Then, there are
    $k\in \Zgeo^{2l}$ and flat $m\in \Zl$ with the
    following properties.
    \label{c:wrapped-fitting-prescribed-modulo}
    \begin{enumerate}
    \item For each $i\in\oi(l)$,
      $0<\sig(k)_{i}\equiv_{r_{i}} \lam_{i}$.
      \label{c:wrapped-fitting-prescribed-modulo-k}
    \item $\mu_{m,k}$ is wrapped and fitting.
\label{c:wrapped-fitting-prescribed-modulo-fitting-wrapped}
    \end{enumerate}
  \item Suppose $\lam\in \oi(l)$ and $R\in \Zgeo$. Then,
    there are $k\in \Zgeo^{2l}$ and flat $m\in \Zl$ with
    the following properties.
    \label{c:wrapped-fitting-prescribed-any-num}
    \begin{enumerate}
    \item $\sig(k)_{i}= 2(l-i)+1$ when $i\in\oi(\lam+1,l)$.
      \label{c:wrapped-fitting-prescribed-any-num-right}
    \item $\sig(k)_{i}=R+2(l-i)$ when
      $i\in\oi(\lam)$.
      \label{c:wrapped-fitting-prescribed-any-num-left}
    \item $\mu_{m,k}$ is wrapped and fitting.
      \label{c:wrapped-fitting-prescribed-any-num-fitting-wrapped}
    \end{enumerate}
  \end{enumerate}
\end{lemma}
\begin{proof}
  Claim~\ref{c:wrapped-fitting-prescribed-ab-pos}
  holds, since
  $b_{j}-b_{i}=(a_{j}+\mn)-(a_{i}+\mn)=a_{j}-a_{i}$
  by the flat $m$
  and  $a_{j}-a_{i}=\sum k(i+1,j)>0$ by $k\in \Zgeo^{2l}$.

  Claim~\ref{c:wrapped-fitting-prescribed-ab-pos} gives
  Claim~\ref{c:wrapped-fitting-prescribed-flat-fitting} by
  $a_{1}=k_{1}\geq 0$.
  Let us prove Claim~\ref{c:wrapped-fitting-prescribed-flat-wrapped}.   First, $\mu_{m,k}$ is fitting by
  Claim~\ref{c:wrapped-fitting-prescribed-flat-fitting}, because
  $m\ldZl s$ by $\sig(k)_{1}\geq 0$ and  $s_{2}=\infty$.
  Second,  $\mu_{m,k}$ is wrapped, because
  $m\bop k\ldZl s$ and
  $(m\bom k)_{i}=m_{i}-\sig(k)_{i}
  \geq m_{i}-\sig(k)_{1}\geq s_{1}$ by
  Lemma~\ref{lem:slope-conditions-sigma}.

  Let us prove
  Claim~\ref{c:wrapped-fitting-prescribed-modulo}.  Hence,
  let us obtain $k\in \Zgeo^{2l}$ with
  Property~\ref{c:wrapped-fitting-prescribed-modulo-k}.  Let
  $k(1,l)=\iota^{l}(1)$ and
  $k_{l+1}=p_{1}r_{l}+\lam_{l}\in \Zgeo$ for some
  $p_{1}\in \Zgeo$.  Moreover, $i\in\oi(2,l)$ implies
  $2\leq l-i+2< l+i-1< l+i\leq 2l$.  Thus, we inductively put
  $p_{i}\in \Zgeo$ for $i\in\oi(2,l)$ so that
  $k_{l+i}\in \Zgeo$ and
  \begin{align}
    k_{l+i}
    &=p_{i}r_{l-i+1}+\lam_{l-i+1}-\sum k(l-i+2,l+i-1).
      \label{eq:wrapped-fitting-prescribed-p1}
  \end{align}
  
  First, if $i=l$, then
  $\sig(k)_{l}=k_{l+1}\equiv_{r_{l}} \lam_{l}$.
  Second, let  $i\in\oi(l-1)$.
  Then,
  \begin{align}
    2\leq l-i+1 \leq  l.
    \label{eq:wrapped-fitting-prescribed-p2}
  \end{align}
  This gives $2\leq l-(l-i+1)+2<l+(l-i+1)-1\leq 2l-1$ and
  $(i+1)+1\leq 2l-(i+1)+1$.  Hence, we have
  \begin{align}
    \sum k(l-(l-i+1)+2,l+(l-i+1)-1)
    &=\sum k(i+1,2l-i)
    \\&=\sum k((i+1)+1,2l-(i+1)+1)+k_{i+1}
    \\    &=\sig(k)_{i+1}+k_{i+1}.
  \end{align}
  Therefore, since $l-(l-i+1)+1=i$,
  equation~\eqref{eq:wrapped-fitting-prescribed-p1} and
  inequalities~\eqref{eq:wrapped-fitting-prescribed-p2}
  imply
  \begin{align}
    \sig(k)_{i}
    &=\sum k(i+1,2l-i+1)\\
    &=\sig(k)_{i+1}+k_{i+1}+k_{l+l-i+1}\\
    &=
      \sig(k)_{i+1}+k_{i+1}+p_{l-i+1}
      r_{i}+
      \lam_{i}-\sig(k)_{i+1}-k_{i+1}\\ 
    &\equiv_{r_{i}} \lam_{i}.
  \end{align}
  Then, Property~\ref{c:wrapped-fitting-prescribed-modulo-k}
  holds, since $k\in \Zgeo^{2l}$ implies $\sig(k)>0$.
  Moreover, Claim~\ref{c:wrapped-fitting-prescribed-modulo}
  follows, since
  Claim~\ref{c:wrapped-fitting-prescribed-flat-wrapped}
  gives a flat $m\in \Zl$ with
  Property~\ref{c:wrapped-fitting-prescribed-modulo-fitting-wrapped}.
  
  Let us prove
  Claim~\ref{c:wrapped-fitting-prescribed-any-num}.  Let
  $k\in \Zgeo^{2l}$ such that $k_{i}=1$ for
  $2\leq i\neq \lam+1 \leq 2l$ and $k_{\lam+1}=R$.  Then,
  Property~\ref{c:wrapped-fitting-prescribed-any-num-right}
  holds, since $i\in\oi(\lam+1,l)$ implies
  $\sig(k)_{i}= \sum k(i+1,2l-i+1)= 2(l-i)+1$ by
  $\lam+1<i+1\leq 2l-i+1$.  Also,
  Property~\ref{c:wrapped-fitting-prescribed-any-num-left}
  follows, since $i\in\oi(\lam)$ implies
  $i+1\leq \lam+1\leq 2l-i+1$, and hence
  \begin{align}
    \sig(k)_{i}
    = \sum k(i+1,2l-i+1)
    =(2l-i+1)-i-1+R
    =R+2(l-i).
  \end{align}
  Thus, Claim~\ref{c:wrapped-fitting-prescribed-any-num}
  follows, since
  Claim~\ref{c:wrapped-fitting-prescribed-flat-wrapped}
  gives a flat $m\in \Zl$ with
  Property~\ref{c:wrapped-fitting-prescribed-any-num-fitting-wrapped}.
\end{proof}

Let us rewrite $n_{a}^{b}(2l,w)$.

\begin{lemma}\label{lem:weight-sigma-sum}
  Let $l\in \Zgeo$.  Consider $m,n\in \Zl$ and $k\in \Ztl$ with
  $a=\nu(k)$ and $b=\nu(m,n,k)$.  Then,
  $n_{a}^{b}(2l,w)=\sum_{i\in\oi(l)}
  w_{i}\sig(k)_{i}(n_{l-i+1}+\sig(k)_{i}-m_{i})$.
\end{lemma}
\begin{proof}
  The statement holds by Lemma~\ref{lem:sigma-nu}, since
  $n_{a}^{b}(2l,w)=\sum N_{a}^{b}(2l,w)$ for
  $N_{a}^{b}(2l,w)_{i}
  =w_{i}(b_{2l-i+1}-b_{i})(a_{2l-i+1}-a_{i})$.
\end{proof}

This gives the
$\sig$-differences $t_{\gam,\Del}$ below by
quadratic polynomials $t_{\gam}$ 
in Definition~\ref{def:intro-quad-poly}.
\begin{lemma}\label{lem:sigma-diff-poly}
  Let $l\in \Zgeo$ and $\gam\in \prod_{i\in\oi(l)}\Qr$.  Suppose
  $m,n\in \Zl$ and $k\in \Ztl$. Then, we have the following.
  \begin{enumerate}
  \item \label{c:sigma-diff-poly-quad-sigma}
    $t_{\gam,\Del}(m,n,k) =2\sum_{i\in\oi(l)}
    \gam_{i,1}\sig(k)_{i}(n_{l-i+1}+\sig(k)_{i}-m_{i})$.
  \item\label{c:sigma-diff-poly-integer} If each
    $2\gam_{i,1}\in \Z$, then $t_{\gam,\Del}(m,n,k)\in \Z$.
  \end{enumerate}
\end{lemma}
\begin{proof}
  Claim~\ref{c:sigma-diff-poly-quad-sigma} holds, because we
  have
  \begin{dmath*}
    t_{\gam,\Del}(m,n,k)
    =
    t_{\gam}(m\bom k) +t_{\gam}((n\bop k)^{\ve})
    -t_{\gam}(m)-t_{\gam}(n^{\ve}) 
    =\sum_{i\in\oi(l)}
    t_{\gam,i}(m_{i}-\sig(k)_{i})
    +t_{\gam,i}(n_{l-i+1}+\sig(k)_{i})
    -t_{\gam,i}(m_{i})-t_{\gam,i}(n_{l-i+1}) 
    =\sum_{i\in\oi(l),j\in\oi(3)} 
    \gam_{i,j}(m_{i}-\sig(k)_{i})^{3-j}
    +
    \gam_{i,j}(n_{l-i+1}+\sig(k)_{i})^{3-j}
    -\gam_{i,j} m_{i}^{3-j} - \gam_{i,j}    n_{l-i+1}^{3-j}
    =\sum_{i\in\oi(l)}
    \gam_{i,1}(-2m_{i}\sig(k)_{i}+\sig(k)_{i}^2)
    +\gam_{i,1}(2n_{l-i+1}\sig(k)_{i}+\sig(k)_{i}^2)
  \end{dmath*} 
  Claim~\ref{c:sigma-diff-poly-integer} follows from
  Claim~\ref{c:sigma-diff-poly-quad-sigma}.
\end{proof}
We  discuss the sum monomial condition of monomial indices by
the following.

\begin{lemma}\label{lem:sums-ineq-mon-cond}
  Let $l\in \Zgeo$ and $u,w\in \Ql$.  Then, each
  $h\in\oi(l)$ satisfies
  \begin{align}
    0\leq \sum_{i\in\oi(h)}u_{i}\leq  \sum_{i\in\oi(h)} w_{i}
    \label{ineq:sums-ineq-mon-cond-uw}
  \end{align}
  if and only if each decreasing $\ka\in \Ql_{\geq 0}$ satisfies
  \begin{align}
    0\leq  \sum_{i\in\oi(l)}\ka_{i}u_{i} \leq \sum_{i\in\oi(l)}\ka_{i}w_{i}.
    \label{ineq:sums-ineq-mon-cond-ka-uw}
  \end{align}
\end{lemma}
\begin{proof}
  We obtain inequalities~\eqref{ineq:sums-ineq-mon-cond-uw}
  from inequalities~\eqref{ineq:sums-ineq-mon-cond-ka-uw},
  taking $\ka=\iota^{h}(1)\ccn \iota^{l-h}(0)$.  Hence, let us
  prove inequalities~\eqref{ineq:sums-ineq-mon-cond-ka-uw}
  from inequalities~\eqref{ineq:sums-ineq-mon-cond-uw}.  If
  $l=1$, then this holds by $0\leq u_{1}\leq w_{1}$ and
  $\kn\geq 0$. Thus, let $l\geq 2$.
  
  First, let us prove
  $0\leq \sum_{i\in\oi(l)}\ka_{i}u_{i}$.
  Inequalities~\eqref{ineq:sums-ineq-mon-cond-uw} imply
  $-u_{l}\leq \sum_{i\in\oi(l-1)}u_{i}$.  Also,
  $\ka(1,l-1)-\ka_{l}\in \Ql_{\geq 0}$ is decreasing.  Hence, the
  induction on $l$ gives
  $\sum_{i\in\oi(l)}\ka_{i}u_{i}
    =\sum_{i\in\oi(l-1)} \ka_{i}u_{i}+\ka_{l}u_{l}
    \geq \sum_{i\in\oi(l-1)} \ka_{i}u_{i}-  \sum_{i\in\oi(l-1)} \ka_{l}u_{i}
     =\sum_{i\in\oi(l-1)} (\ka_{i}-\ka_{l})u_{i}
    \geq 0$.
    Second, $0\leq \sum_{i\in\oi(l)} \ka_{i} (w_{i}-u_{i})$ holds as
  above, since
  $0\leq \sum_{i\in\oi(h)} (w_{i}-u_{i})\leq \sum_{i\in\oi(h)}w_{i}$ for
  each $h\in\oi(l)$ by
  inequalities~\eqref{ineq:sums-ineq-mon-cond-uw}.  
\end{proof}

Then, we obtain the following equivalence of monomial
indices and functional monomial indices by polynomials.

\begin{theorem}\label{thm:monomial-ind-equiv}
  For $l\in \Zgeo$, let $t_{i}(z)\in \Q[z]$ for each
  $i\in\oi(l)$. Let
  $t(m)=\sum_{i\in\oi(l)} t_{i}(m_{i})\in \Q$ for each
  $m\in \Zl$.  Suppose an infinite gate $s\geq 0$ and
  $w\in \Zgeo^{l}$. Let $\phi=\tp{s,l,w,t}$. Then, the following
  statements are equivalent.
  \begin{enumerate}
  \item
    \label{s:monomial-ind-equiv-functional}
    $\phi$ is a  functional monomial index.
  \item
    \label{s:monomial-ind-equiv-non-functional}
    $\psi=\tp{l,w,\gam}$ such that $t=t_{\gam}$ is
    a monomial index.
  \end{enumerate}
\end{theorem}
\begin{proof}
  Let us prove
  Statement~\ref{s:monomial-ind-equiv-functional} from
  Statement~\ref{s:monomial-ind-equiv-non-functional}.
  First, the integer monomial condition of $\phi$ follows from
  that of $\psi$ by Claim~\ref{c:sigma-diff-poly-integer} of
  Lemma~\ref{lem:sigma-diff-poly}.  Second, the sum monomial
  condition of $\phi$ holds as follows.  Suppose a fitting
  $\tp{s,l,m,n,k}$. Then,
  Lemma~\ref{lem:slope-conditions-sigma} gives a
  decreasing
  $\tp{\sig(k)_{i}(n_{l-i+1}+\sig(k)_{i}-m_{i})}_{i\in\oi(l)}\in
  \Zlgez$.  Then, by Lemma~\ref{lem:sums-ineq-mon-cond},
  the sum monomial condition of $\psi$ gives
  \begin{align}
    0\leq 2\sum_{i\in\oi(l)}\gam_{i,1}
    \sig(k)_{i}(n_{l-i+1}+\sig(k)_{i}-m_{i}) 
    \leq \sum_{i\in\oi(l)}w_{i}\sig(k)_{i}(n_{l-i+1}+\sig(k)_{i}-m_{i}).
  \end{align}
  Thus, the sum monomial condition of $\phi$ holds by
  Lemma~\ref{lem:weight-sigma-sum} and
  Claim~\ref{c:sigma-diff-poly-quad-sigma} of
  Lemma~\ref{lem:sigma-diff-poly}.  In particular,
  Statement~\ref{s:monomial-ind-equiv-functional} follows.
  
  Let us prove
  Statement~\ref{s:monomial-ind-equiv-non-functional} from
  Statement~\ref{s:monomial-ind-equiv-functional}.  First,
  we prove each $\deg_{z}t_{i}(z)\leq 2$.  Assume
  $d=\max( \deg_{z}t_{1}(z),\dots, \deg_{z}t_{l}(z)) >2$.
  Then, suppose all integers $j_{1}<\dots<j_{u}$ such that
  $t_{j_{i}}(z)
  =\alp_{j_{i},d}z^{d}+\alp_{j_{i},d-1}z^{d-1}+\dots$ and
  $\alp_{j_{i},d}\neq 0$.  Also, for a flat $m\in \Zl$, let
  \begin{align}
    t_{j_{i},\Del}(m,m,k)
    &=
      t_{j_{i}}((m\bom k)_{j_{i}})+
      t_{j_{i}}((m\bop k)^{\ve}_{j_{i}})
      -t_{j_{i}}(m_{j_{i}})-t_{j_{i}}(m^{\ve}_{j_{i}})
    \\&=
    t_{j_{i}}(\mn-\sig(k)_{j_{i}})+
    t_{j_{i}}(\mn+\sig(k)_{j_{i}})
    -t_{j_{i}}(\mn)-t_{j_{i}}(\mn).
  \end{align}
  This gives $\deg_{\mn}t_{j_{i},\Del}(m,m,k)\leq d-2$,
  because
  \begin{align}
    \alp_{j_{i},d-\lam}
    \left((\mn-\sig(k)_{j_{i}})^{d-\lam}+
    (\mn+\sig(k)_{j_{i}})^{d-\lam}
    -\mn^{d-\lam}-\mn^{d-\lam}\right)
  \end{align}
  has neither $\mn^{d-\lam}$ nor
  $\mn^{d-\lam-1}$  for each $\lam\in\oi(0,d-1)$.
  In particular, we have
  \begin{align}
    t_{\Del}(m,m,k)=     2m_{1}^{d-2}
    {d \choose d-2}\sum_{i\in\oi(u)}
    \alp_{j_{i},d}  \sig(k)_{j_{i}}^{2}+\dots.
    \label{eq:monomial-ind-equiv-p3}
  \end{align}
  Moreover, by
  Claim~\ref{c:wrapped-fitting-prescribed-any-num} of
  Lemma~\ref{lem:wrapped-fitting-prescribed}, each
  $R\in \Zgeo$ and $\lam\in \oi(l)$ gives a flat
  $g(R,\lam)\in \Zl$ and wrapped fitting
  $\mu(R,\lam)=\tp{s,l,g(R,\lam),g(R,\lam),h(R,\lam)}$ such that
  $\sig(h(R,\lam))_{i}=R+2(l-i)$ for $i\in\oi(\lam)$ and
  $\sig(h(R,\lam))_{i}=2(l-i)+1$ for $i\in\oi(\lam+1,l)$.  In
  particular, a large enough $R\in \Zgeo$ yields
  \begin{align}
    \sum_{i\in\oi(u)} \alp_{j_{i},d}  \sig(h(R,j_{1}))_{j_{i}}^{2} \neq 0.
    \label{eq:monomial-ind-equiv-p4}
  \end{align}
  Also,
  Claim~\ref{c:wrapped-fitting-prescribed-flat-wrapped} of
  Lemma~\ref{lem:wrapped-fitting-prescribed} gives a wrapped
  fitting $\tp{s,l,g',g',h(R,j_{1})}$ for each flat
  $g'\in \Zl$ such that $g'\geq g(R,j_{1})$.  For
  $a(R,j_{1})=\nu(h(R,j_{1}))$ and
  $b(g',R)=\nu(g',g',h(R,j_{1}))$,
  Lemma~\ref{lem:weight-sigma-sum} implies
  $n_{a(R,j_{1})}^{b(g',R)}(2l,w) =\sum_{i\in\oi(l)}
  w_{i}\sig(h(R,j_{1}))_{i}^{2}$, which is independent of
  $g'$.

  Thus, since $d>2$,
  equation~\eqref{eq:monomial-ind-equiv-p3} and
  inequality~\eqref{eq:monomial-ind-equiv-p4} imply that
  large $R\in \Zgeo$ and $g'\geq g(R,j_{1})$ would violate the
  sum monomial condition
  $0\leq t_{\Del}(g',g',h(R,j_{1}))\leq
  n_{a(R,j_{1})}^{b(g',R)}(2l,w)$ of $\phi$.  This gives
  $\gam\in\prod_{i\in\oi(l)} \Qr$ such that $t=t_{\gam}$, since
  $\deg_{z}t_{i}(z)\leq 2$ for each $i\in\oi(l)$.
  
  Second, let us confirm the integer monomial condition of
  $\psi$. Suppose some $j\in\oi(l)$ such that
  $2\gam_{j,1}\not\in \Z$.  For each $i\in \oi(l)$, suppose
  $y_{i}\geq 1$ and $2\gam_{i,1}=\frac{x_{i}}{y_{i}}$ with
  coprime $x_{i}$ and $y_{i}$.  By
  Claim~\ref{c:wrapped-fitting-prescribed-modulo} of
  Lemma~\ref{lem:wrapped-fitting-prescribed},
  $y=\tp{y_{i}}_{i\in\oi(l)}$ gives a flat $g(y)\in \Zl$ and
  wrapped fitting $\tp{s,l,g(y),g(y),h(y)}$ such that
  $\sig(h(y))_{j}\equiv_{y_{j}}1$ and
  $\sig(h(y))_{i}\equiv_{y_{i}}0$ for each
  $1\leq i \neq j\leq l$.  Then,
  $2\gam_{j,1}\sig(h(y))_{j}^{2}\not\in \Z$ and
  $2\gam_{i,1}\sig(h(y))_{i}^{2}\in \Z$ for each
  $1\leq i \neq j\leq l$.  However, this would violate the integer
  monomial condition of $\phi$, since
  $t_{\Del}(g(y),g(y),h(y))=
  \sum_{i}2\gam_{i,1}\sig(h(y))_{i}^{2}\in \Z$ by
  Claim~\ref{c:sigma-diff-poly-quad-sigma} of
  Lemma~\ref{lem:sigma-diff-poly}.  Thus, $\psi$ has the
  integer monomial condition.
  
  Third, let us prove the sum monomial condition of $\psi$.
  For each $R\in \Zgeo$ and $\lam\in\oi(l)$, let
  $a(R,\lam)=\nu(g(R,\lam))$ and
  $b(R,\lam)=\nu(g(R,\lam),g(R,\lam),h(R,\lam))$ by
  $\mu(R,\lam)$ above.  Then,
  Lemma~\ref{lem:weight-sigma-sum} and
  Claim~\ref{c:sigma-diff-poly-quad-sigma} of
  Lemma~\ref{lem:sigma-diff-poly} yield
  \begin{align}
    &n_{a(R,\lam)}^{b(R,\lam)}(2l,w)
      -t_{\Del}(g(R,\lam),g(R,\lam),w)
    \\&
    =
    \sum_{i\in\oi(\lam)}(w_{i}-2\gam_{i,1})(R+2(l-i))^{2}
    +   \sum_{i\in\oi(\lam+1,l)}(w_{i}-2 \gam_{i,1})(2(l-i)+1)^{2}.
  \end{align}
  Thus, because $w>0$ gives
  $\sum_{i\in\oi(\lam)} w_{i}\sig(h(R,\lam))_{i}^{2} =
  \sum_{i\in\oi(\lam)} w_{i}(R+(l-i))^{2} >0$, we have
  \begin{align}
    \lim_{R\to\infty}
    \frac{n_{a(R,\lam)}^{b(R,\lam)}(2l,w)
    -t_{\Del}(g(R,\lam),g(R,\lam),w)}{
    \sum_{i\in\oi(\lam)} w_{i}\sig(h(R,\lam))_{i}^{2}}
    =
    1-\frac{\sum_{i\in\oi(\lam)}2\gam_{i,1}}
    {\sum_{i\in\oi(\lam)}w_{i}}.
  \end{align}
  This limit has to be non-negative by
  $n_{a(R,\lam)}^{b(R,\lam)}(2l,w)
  -t_{\Del}(g(R,\lam),g(R,\lam),w)\geq 0 $ in the sum monomial
  condition of $\phi$.  Hence, $\sum_{i\in\oi(\lam)}w_{i}>0$
  implies that each $\lam\in\oi(l)$ gives
  \begin{align}
    \sum_{i\in\oi(\lam)}2\gam_{i,1}\leq \sum_{i\in\oi(\lam)}w_{i}.
  \end{align}
  Similarly,
  $\lim_{R\to\infty} \frac{t_{\Del}(g(R,\lam),g(R,\lam),w)}{
    \sum_{i\in\oi(\lam)} w_{i}\sig(h(R,\lam))_{i}^{2}}
  =\frac{\sum_{i\in\oi(\lam)}2\gam_{i,1}}{\sum_{i\in\oi(\lam)}w_{i}}$.
  This limit has to be non-negative by
  $t_{\Del}(g(R,\lam),g(R,\lam),w)\geq 0$ in the sum monomial
  condition of $\phi$.  Thus, each $\lam\in\oi(l)$ gives
  \begin{align}
    0\leq \sum_{i\in\oi(\lam)}2\gam_{i,1}.
  \end{align}
  Therefore, $\psi$ has the sum monomial condition.
\end{proof}

\begin{remark}\label{rmk:monomial-poly}
  Suppose the notation in
  Theorem~\ref{thm:monomial-ind-equiv} with $l=1$.  If
  $\phi$ is a functional monomial index, then
  Theorem~\ref{thm:monomial-ind-equiv} and the monomial
  conditions of $\psi$ imply $t=t_{\gam}$ and
  $\gam_{1,1}\geq 0$.
  Proposition~\ref{prop:vanishing-const-neg-quad} also gives
  this inequality in some general setting.
\end{remark}

\begin{remark}\label{rmk:without-tuple-flips}
  The following variant of merged determinants would give
  fewer monomial indices for an analog of
  Theorem~\ref{thm:monomial-ind-equiv}.  For simplicity, let
  $\cF=\Lam(s,l,w,\scc,\fs,x,\fX)$. Also, assume a wrapped
  fitting $\tp{s,l,m,n,k}$ with $a=\nu(k)$ and $b=\nu(m,n,k)$.

  Then, $\Delta(\cF)(s,l,w,m,n,k,x,\fX)$ takes flips
  $x^{\ve}$, $w^{\ve}$, $n^{\ve}$, and $(n\bop k)^{\ve}$.
  Hence, without flips, consider the following variant
  \begin{align}
    \Delta'(\cF)(s,l,w,m,n,k,x,\fX)
    =
    \frac{(b)_{x\ccn x}^{w\ccn w}}{(a)_{x\ccn x}^{w\ccn w}}\left(
    \frac{f_{s,m}}{(m)_{x}^{w}}\cdot
    \frac{f_{s,n}}{(n)_{x}^{w}}
    -
    \frac{f_{s,m\bom k}}{(m\bom k)_{x}^{w}}\cdot
    \frac{f_{s,n\bop k}}{(n\bop k)_{x}^{w}}
    \right)\in \Q(\fX).
  \end{align}
  Moreover, suppose a monomial index $\tp{l,w,\gam}$.  Then,
  $t_{\gam,\Del}(m,n,k)$ takes $n^{\ve}$ and
  $(n\bop k)^{\ve}$. Instead, let
  \begin{align}
    t'_{\gam,\Del}(m,n,k)=t_{\gam}(m\bom k)+
    t_{\gam}(n\bop k)-t_{\gam}(m)-t_{\gam}(n).
  \end{align}
  
  In general, $t'_{\gam,\Del}(m,n,k)$ has non-zero
  $\gam_{i,2}$ terms, unlike $t_{\gam,\Del}(m,n,k)$ in
  Claim~\ref{c:sigma-diff-poly-quad-sigma} of
  Lemma~\ref{lem:sigma-diff-poly}.  Moreover, let $s$ be
  infinite.  Then, as in the proof of
  Theorem~\ref{thm:monomial-ind-equiv}, we would demand
  $\tp{\gam_{i,2}}_{i\in\oi(l)}$ to be palindromic for the
  integer monomial condition of $\tp{s,l,w,t_{\gam}}$.  This
  leaves fewer monomial indices for us to compute
  polynomials with positive integer coefficients.
  
\end{remark}

\subsection{Monomial  parcels}
\label{sec:monomial}
By monomial indices, we introduce monomial parcels, which
are explicit merged-log-concave parcels.  In particular,
they extend Proposition~\ref{prop:constant-parcel} by flat
bases.

\begin{definition}\label{def:monomial-parcel}
  Suppose an infinite gate $s\geq 0$ and monomial index
  $\tp{l,w,\gam}$.  Let $q\in \Q(\fX)$.
  \begin{enumerate}
  \item We define the $t$-monomials
    $\Psi_{s,\gam,q}=\{\Psi_{s,\gam,q,m}\in \Q(\fX)\}_{m\in \Zl}$
    such that
    \begin{align}
      \Psi_{s,\gam,q,m}
      &=
        \begin{dcases}
          q^{t_{\gam}(m)}
          \mwhen  m\ldZl s,\\
          0 \melse.
        \end{dcases}
    \end{align}
  \item Let $q$ be $\scc$-admissible.
    For each $m\ldZl s$,
    assume
    \begin{align}
      \Psi_{s,\gam,q,m}=q^{t_{\gam}(m)} \scc 0.
    \end{align}
    Suppose a proper $\tp{s,l,w,\llq,\rho,x,\fX}$-mediator
    $\phi$ for $x=\iota^{l}(q)$.  Then, we define the monomial
    parcel
    \begin{align}
      \cF=\Lam(s,l,w,\scc,\Psi_{s,\gam,q},\phi,\rho,x,\fX).
    \end{align}
    In particular, we call $\cF$ linear or quadratic,
    when each $\gam_{i,1}=0$ or otherwise.
  \end{enumerate}
\end{definition}
Explicitly, each $m\ldZl s$ gives
\begin{align}
  \cF_{m}
  = \frac{ \Psi_{s,\gam,q,m}}{\prod \phi(x)^{m\rc w}
  \cdot [m]!_{q}^{w}}  
  = \frac{q^{t_{\gam}(m)}}{\prod \phi(x)^{m\rc w}
  \cdot [m]!_{q}^{w}},
\end{align}
which is $\frac{q^{t_{\gam}(m)}}
{(m)^{w}_{q}}$ when
$\phi$ is the canonical mediator. Moreover,
we  obtain the following $q$-polynomials with positive integer
coefficients by quasi-merged determinants.

\begin{theorem}\label{thm:monomial-poly}
  Consider a monomial parcel
  $\cF=\Lam(s,l,w,\scc,\Psi_{s,\gam,q},\phi,\rho,x,\fX)$.  For each
  fitting $\tp{s,l,m,n,k}$ with $a=\nu(k)$ and
  $b=\nu(m,n,k)$, we have
  \begin{align}
    q^{-t_{\gam}(m)-t_{\gam}(n^{\ve})}
    \Delta(\cF)(s,l,w,m,n,k,\phi,\rho,x,\fX)
    =\dt(V)_{a}^{b}(s,l,w,\phi,\rho,t_{\gam,\Del}(m,n,k),x)         
    \llq 0.
  \end{align}
  In particular, $\cF$ is $\scc$-merged-log-concave.
\end{theorem}
\begin{proof}
  Statements follow from
  Proposition~\ref{prop:functional-merged} and
  Theorem~\ref{thm:monomial-ind-equiv}.
\end{proof}

If $l=1$ and $w=\tp{1}$, then the monomial conditions of
$\tp{l,w,\gam}$ imply $\gam_{1,1}=0$ or
$\gam_{1,1}=\fraa$.
Suppose $s=\tp{0,\infty}$. Also,
 let
$\gam_{1}=\tp{\tp{0,0,0}}$ and
$\gam_{2}=\tp{\tp{\fraa,-\fraa,0}}$
so that we have the linear and quadratic
$\cF_{i}=\Lam(s,l,w,\scc,\Psi_{s,\gam_{i},q},\rho,x,\fX)$ of
$i\in\oi(2)$.  Then,
for an indeterminate
$t$,
$(\pm t;q)^{\mp 1}_{\infty}$ are
$\sum_{m\in \Zgeo^{l}}\cF_{i,m}t^{m_{1}}$ of $i\in\oi(2)$
by Euler's binomial identities.

\begin{example}
  For $s=\tp{0,\infty}$, $l=1$, $w=\tp{1}$, and
  $\gam=\tp{\tp{0,0,0}}$,  consider
  $\cF=\Lam(s,l,w,\scc,\Psi_{s,\gam,q},x,\fX)$. 
  Let  $m=n=\tp{2}$ and $k=\tp{1,1}$. Then, we have
  $a=\nu(k)=\tp{1,2}$,
  $b=\nu(m,n,k)=\tp{3,4}$,
  $m\bom k=\tp{1}$,
  $n\bop k=\tp{3}$.
  Hence, we
  obtain
  \begin{align}
    \Delta(\cF)(s,l,w,m,n,k,x,\fX)
    &=\frac{ (3)_{q}(4)_{q}}{  (1)_{q}(2)_{q}}
      \left(\frac{1}{(2)_{q}}\cdot \frac{1}{(2)_{q}}
      -
      \frac{1}{(1)_{q}}\cdot \frac{1}{(3)_{q}}\right)
    \\&
    ={ 3 \brack 1 }_{q}{ 4 \brack 2 }_{q}-
    { 3 \brack 2 }_{q}{ 4 \brack 1 }_{q}
    \\&
    =q^6 + q^5 + 2q^4 + q^3 + q^2\llq 0.
  \end{align}
\end{example}

\begin{example}
  Let $s=\tp{0,\infty}$, $l=2$, $w=\iota^{l}(1)$, and
  $\gam=\iota^{l}(\tp{0,0,0})$. Suppose
  $\cF=\Lam(s,l,w,\scc,\Psi_{s,\gam,q},x,\fX)$. 
  Moreover, let   $m=n=\iota^{l}(3)$ and $k=\iota^{2l}(1)$
  so that
  $a=\nu(k)=\tp{1,2,3,4}$,
  $b=\nu(m,n,k)=\tp{4,5,6,7}$,
  $m\bom k=\tp{0,2}$,
  $n\bop k=\tp{4,6}$.
  Then, we have
  \begin{align}
    \Delta(\cF)(s,l,w,m,n,k,x,\fX)
    &=\prod_{i\in \oi(4)}\frac{  (i+3)_{q}}
      {  (i)_{q}} \cdot
      \Bigg(\frac{1}{(3)_{q}(3)_{q}}\cdot
      \frac{1}{(3)_{q}(3)_{q}}
      -
      \frac{1}{(0)_{q}(2)_{q}}
      \cdot \frac{1}{(4)_{q}(6)_{q}}\Bigg)
    \\&
    ={ 4 \brack 1 }_{q}{ 5 \brack 2 }_{q}
    { 6 \brack 3 }_{q}{ 7 \brack 4 }_{q}
    -
    { 4 \brack 4 }_{q}{ 5 \brack 3 }_{q}
    { 6 \brack 2 }_{q}{ 7 \brack 1 }_{q},
  \end{align}
  which is
  \begin{dmath*}
    q^{30} + 4q^{29} + 13q^{28} + 34q^{27} + 76q^{26} + 151q^{25} + 273q^{24} + 452q^{23} + 695q^{22} + 999q^{21} + 1346q^{20} + 1710q^{19} + 2052q^{18} + 2330q^{17} + 2506q^{16} + 2557q^{15} + 2470q^{14} + 2262q^{13} + 1958q^{12} + 1600q^{11} + 1229q^{10} + 886q^9 + 593q^8 + 368q^7 + 208q^6 + 106q^5 + 47q^4 + 18q^3 + 5q^2 + q \llq 0.
  \end{dmath*}
\end{example}
\subsection{On the density of merged determinants}\label{sec:density}

We use the following notation to discuss merged determinants
of some monomial parcels.

\begin{definition}\label{def:density}
  Consider $f\in \Q[z]$.
  \begin{enumerate}
  \item Let us call $f$ $z$-dense, if $f_{i}\neq 0$ for each
    $i\in \Z$ such that $\ord_{z}(f)\leq i\leq \deg_{z}(f)$.
  \item Suppose $g\in \Q[z]$.  Let us write $f\gzd g$, if
    $f-g$ is $z$-dense and $f-g>_{z}0$.
  \end{enumerate}
\end{definition}

In particular,
 $\rho\in \Zgeo$ give
different
$z^{\rho}$-densities.
Hence, we consider
a relative version of the assertion that
the sequence $\{f_{i}\}_{i\in\Z}$
of a non-zero $f\in \Q[z]$ 
has no internal zeros.
We first state
the following density transitivity by base shift.

\begin{lemma}\label{lem:density-bshift}
  Let $\rho,\lam\in \Zgeo$.  If $f >_{z^{\rho},d}0$, then
  $f b_{\lam,\rho}(z)\gzd0$.
\end{lemma}
\begin{proof}
  If $\rho=1$, then the statement follows from $f\gzd0$
  and
  $b_{\lam,\rho}(z)=1$ in
  Lemma~\ref{lem:bshift-fun-sp}.
  Hence, let $\rho\geq 2$.  Also, we assume
  $\ord_{z}(f)=0$,
  replacing $f$ by $z^{-\ord_{z}(f)} f>_{z^{\rho},d}0$.

  First, suppose $\lam=1$.  Then,
  Lemma~\ref{lem:bsingle-qnum} gives
  $b_{\lam,\rho}(z)=[\rho]_{z}$.  Thus, the statement
  follows,
  because we have $(f[\rho]_{z})_{i}>f_{\rho r}>0$ for
  each
  $i\in\oi(0,\deg_{z}(f)+\rho-1)$ and $i=\rho r+\mu$ such
  that
  $r\geq 0$ and $\mu\in\oi(0,\rho-1)$.
  
  Second, assume $\lam\geq 2$.  Then, $\rho\geq 2$ and
  Lemma~\ref{lem:bsingle-degree} give
  $\deg_{z}(b_{\lam-1,\rho}(z))
    =\frac{(\rho-1)\lam (\lam-1)}{2}\geq \lam-1$.
    Therefore, $\ord_{z}(f b_{\lam-1,\rho}(z))=0$ and
  $fb_{\lam-1,\rho}(z)\gzd 0$ yield
  \begin{align}
    fb_{\lam,\rho}(z)
    =f b_{\lam-1,\rho}(z) [\rho]_{z^{\lam}}
    =\sum_{i\in\oi(0,\rho-1)} f b_{\lam-1,\rho}(z)
    z^{\lam i}
    \gzd 0.
  \end{align}
\end{proof}

Moreover, the binary relation
$\gzd$ has the following poring.
\begin{lemma}\label{lem:density-poring}
  $\Q[z]$ is a strict $\gzd$-poring.
\end{lemma}
\begin{proof}
  First, since $0\gzd 0$ is false, the irreflexivity holds.
  Second, suppose $f_{1}\gzd f_{2} \gzd f_{3}$.  Then,
  $\ord_{z}(f_{j})\leq \ord_{z}(f_{j+1})$,
  $\deg_{z}(f_{j}) \geq \deg_{z}(f_{j+1})$, and
  $f_{j,i}> f_{j+1,i}$ when
  $i\in \oi(\ord_{z}(f_{j}),\deg_{z}(f_{j}))$ and
  $j\in \oi(2)$.  Hence,
  $\ord_{z}(f_{1})\leq \ord_{z}(f_{3})$,
  $\deg_{z}(f_{1})\geq \deg_{z}(f_{3})$, and
  $f_{1,i}>f_{3,i}$ when
  $i\in \oi(\ord_{z}(f_{1}),\deg_{z}(f_{1}))$.  The
  transitivity $f_{1}\gzd f_{3}$ follows.  Third, if
  $f_{1}\gzd f_{2}$ and $f_{3}\in \Q[z]$, then the additivity
  follows from
  $(f_{1}+f_{3})-(f_{2}+f_{3}) =f_{1}-f_{2}\gzd 0$.  Fourth,
  suppose $f_{1},f_{2}\gzd 0$.  When
  $i\in \oi(\ord_{z}(f_{1})+\ord_{z}(f_{2}),
  \deg_{z}(f_{1})+\deg_{z}(f_{2}))$, we have
  $j_{k}\in \oi(\ord_{z}(f_{k}), \deg_{z}(f_{k})) $ for
  $k\in \oi(2)$ such that $j_{1}+j_{2}=i$.  Thus, the
  multiplicativity $f_{1}f_{2} \gzd 0$ follows from
  $(f_{1}f_{2})_{i}\geq f_{1,j_{1}}f_{2,j_{2}}>0$.
\end{proof}

Hence, we have the following $q$-density of pre-parcels.

\begin{proposition}\label{prop:density-pre-parcels}
  Let $\mu\in \Zgeo$ and $x=\iota^{\mu}(q)$. Consider
  $a,b\in \Z^{\mu}$ such that $b\geq a\geq 0$.  Then, we have
  $U_{a}^{b}(\mu,w,\rho,x)>_{q,d} 0$.
\end{proposition}
\begin{proof}
  Since
  ${b_{i} \brack a_{i}}_{q^{\rho_{i}}}>_{q^{\rho_{i}},d}0$, we
  have
  $U_{\tp{a_{i}}}^{\tp{b_{i}}}(1,\tp{1},\tp{\rho_{i}},\tp{q})
  ={b_{i} \brack a_{i}}_{q^{\rho_{i}}}
  b_{b_{i}-a_{i},\rho_{i}}(q)>_{q,d}0$ by
  Lemma~\ref{lem:density-bshift}.  Therefore,
  Lemma~\ref{lem:density-poring} gives the assertion, since
  $U_{a}^{b}(\mu,w,\rho,x)
  =\prod_{i\in\oi(\mu)} U_{\tp{a_{i}}}^{\tp{b_{i}}}(1,\tp{1},
  \tp{\rho_{i}},\tp{q})^{w_{i}}$.
\end{proof}

To discuss pre-merged determinants, 
we put the following notion of {\it strict} or {\it almost
  strict} pre-fitting tuples, which are tempered.
\begin{definition}
  Suppose $\ups=\tp{\mu,a,b}$ for $\mu\in \Zget$ and
  $a,b\in \Z^{\mu}$.
  \begin{enumerate}
  \item We call $\ups$ strictly pre-fitting, if
    $0<a_{1}<\dots <a_{\mu}\leq b_{1}<\dots <b_{\mu}$.
  \item We call $\ups$ almost strictly pre-fitting, if
    $0\leq a_{1}<\dots <a_{\mu}\leq b_{1}<\dots <b_{\mu}$.
  \end{enumerate}
\end{definition}

Then, we conjecture the following $q$-density of pre-merged
determinants.  This implies the $q$-density of merged
determinants of some width-two monomial parcels by
Lemma~\ref{lem:pre-merged-quasi-merged} and
Theorem~\ref{thm:monomial-poly}.

\begin{conjecture}\label{conj:density}
  Let $\mu=2$, $w=\iota^{\mu}(1)$, and $x=\iota^{\mu}(q)$.  Suppose a
  flat $\rho\in \Zgeo^{\mu}$ and strictly pre-fitting
  $\ups=\tp{\mu,a,b}$.
  Let $e\in\oi(0,n_{a}^{b}(\mu,w))$.  Then,
  \begin{align}
    \dt(U)_{a}^{b}(\mu,w,\rho,e,x)>_{q,d} 0.
  \end{align}
\end{conjecture}

\begin{example}
  Conjecture~\ref{conj:density} does not extend to almost
  strictly pre-fitting tuples. For example, if
  $w=\rho=\iota^{2}(1)$, $a=\tp{0,2}$, $b=\tp{2,3}$, and
  $e=1$, then
  \begin{align}
    \dt(U)_{a}^{b}(2,w,\rho,e,x)=
    { 2 \brack 0 }_{q}{ 3 \brack 2 }_{q}-
    q { 2 \brack 2 }_{q}{ 3 \brack 0 }_{q}
    =
    q^2+1\not>_{q,d} 0.
  \end{align}
\end{example}

Furthermore, we prove
the following $q$-density of higher-width pre-merged
determinants, assuming Conjecture~\ref{conj:density}.

\begin{theorem}\label{thm:pre-merged-density}
  Assume Conjecture~\ref{conj:density}.  Let $\mu\in \Zget$
  and $x=\iota^{\mu}(q)$. Consider a pre-fitting
  $\tp{\mu,a,b}$ such that $a>0$.  Let
  $e\in\oi(0,n_{a}^{b}(\mu,w))$.  Then, we have
  \begin{align}
    \dt(U)_{a}^{b}(\mu,w,\rho,e,x)\gqd 0.
  \end{align}
\end{theorem}
\begin{proof}
  If $a^{\ve}\not\leq b$, then
  $ U_{a^{\ve}}^{b}(\mu,w,\rho,x)=0$.  This gives the assertion
  by Proposition~\ref{prop:density-pre-parcels}.  Hence, let
  \begin{align}
    a^{\ve}\leq b.
    \label{ineq:pre-merged-density-a-ve-leq-b}
  \end{align}

  If $\mu=2$, then the assertion becomes
  Conjecture~\ref{conj:density}.  Assume an odd $\mu\geq 3$.
  Then, Claim~\ref{c:outer-preparcel-premerged-premerged}
  of
  Lemma~\ref{lem:outer-preparcel-premerged} gives
  $\dt(U)_{a}^{b}(\mu,w,\rho,e,x) ={\cen(b) \brack
    \cen(a)}^{\cen(w)}_{\cen(x)^{\cen(\rho)}}
  b_{\cen(b)-\cen(a),\cen(\rho)}(q) \cdot
  U_{\out(a)}^{\out(b)}(\mu-1,\out(\rho),\out(w),e,\out(x))$.
  Also,
  $0\leq e\leq n_{a}^{b}(\mu,w) =
  n_{\out(a)}^{\out(b)}(\mu-1,\out(w))$.  Therefore, since
  ${\cen(b) \brack \cen(a)}^{\cen(w)}_{\cen(x)^{\cen(\rho)}}
  b_{\cen(b)-\cen(a),\cen(\rho)}(q)>_{q,d}0$ by
  Proposition~\ref{prop:density-pre-parcels}, the assertion
  holds by
  Lemma~\ref{lem:density-poring} and the induction on
  $\mu$
  
  Suppose an even $\mu\geq 2$.  Let $\mu=2l$. Consider
  $E\in \Zl$ such that $0\leq E \leq N_{a}^{b}(\mu,w)$ and
  $\sum E=e$.  Then, 
  Claim~\ref{c:outer-preparcel-premerged-premerged} of
  Lemma~\ref{lem:outer-preparcel-premerged} implies
  \begin{dmath}
    \dt(U)_{a}^{b}(\mu,w,\rho,e,x)
    =
    U_{\cen(a)}^{\cen(b)}(2,\cen(w),\cen(\rho),\cen(x))
    U_{\out(a)}^{\out(b)}(\mu-2,\out(w),\out(\rho),\out(x))
    -
    \prod x(l,l)^{\cen(\mu,E)}\cdot
    U_{\cen(a)^{\ve}}^{\cen(b)}(2,\cen(w),\cen(\rho),\cen(x))
    \cdot \prod x(1,l-1)^{\out(\mu,E)}
    \cdot
    U_{\out(a)^{\ve}}^{\out(b)}(\mu-2,\out(w),\out(\rho),\out(x)).
    \label{eq:pre-merged-density}
  \end{dmath}
  Also, since $\out(a)^{\ve}\leq \out(b)$ by
  inequality~\eqref{ineq:pre-merged-density-a-ve-leq-b}, the
  induction gives
  \begin{dmath}
    U_{\out(a)}^{\out(b)}(\mu-2,\out(w),\out(\rho),\out(x))
    \gqd
    \prod x(1,l-1)^{\out(\mu,E)}
    \cdot
    U_{\out(a)^{\ve}}^{\out(b)}(\mu-2,\out(w),\out(\rho),\out(x))
    \gqd 0.
    \label{ineq:pre-merged-density-outer}
  \end{dmath}
  
  First, assume that $\cen(a)$ or $\cen(b)$ is flat.  Then,
  $N_{a}^{b}(\mu,w)_{l}= w_{l}(b_{l+1}-b_{l}) (a_{l+1}-a_{l})
  =0$ gives $\cen(\mu,E)=\tp{0}$, which implies
  $\prod x(l,l)^{\cen(\mu,E)}\cdot
  U_{\cen(a)^{\ve}}^{\cen(b)}(2,\cen(w),\cen(\rho),\cen(x)) =
  U_{\cen(a)}^{\cen(b)}(2,\cen(w),\cen(\rho),\cen(x))\gqd 0$
  by Proposition~\ref{prop:density-pre-parcels}.  Hence,
  equation~\eqref{eq:pre-merged-density} and
  inequalities~\eqref{ineq:pre-merged-density-outer} give
  the assertion by
Claim~\ref{c:succ-succeq-stwomul2} of
  Lemma~\ref{lem:succ-succeq} and
  Lemma~\ref{lem:density-poring}.

  Second, let $\cen(a)_{1}<\cen(a)_{2}$ and
  $\cen(b)_{1}<\cen(b)_{2}$.
  Then, since $\prod x(l,l)^{\cen(\mu,E)}\gqd 0$ and
  $\cen(b)\geq \cen(a)^{\ve}$ by
  inequality~\eqref{ineq:pre-merged-density-a-ve-leq-b},
  the induction gives
  $U_{\cen(a)}^{\cen(b)}(2,\cen(w),\cen(\rho),\cen(x)) \gqd \prod
  x(l,l)^{\cen(\mu,E)}\cdot
  U_{\cen(a)^{\ve}}^{\cen(b)}(2,\cen(w),\cen(\rho),\cen(x))
  >_{q,d} 0$.
  Thus, equation~\eqref{eq:pre-merged-density} and
  inequalities~\eqref{ineq:pre-merged-density-outer} yield
  the assertion by
  Claim~\ref{c:succ-succeq-zero-zero-2} of
  Lemma~\ref{lem:succ-succeq} and
  Lemma~\ref{lem:density-poring}.
\end{proof}

\subsection{On the almost log-concavity,
  unimodality,
  and palindromicity of pre-merged
  determinants}\label{sec:monomial-conj}

Let us discuss the log-concavity, unimodality, and
palindromicity by pre-merged determinants.
We use the
following terminology to avoid conjecturing upon
Conjecture~\ref{conj:density}.
\begin{definition}\label{def:almost}
  Let $f\in \Q[z]$.  Assume $\ka\in\Zgeo$ and integers
  $\lamo<\dots <\lam_{\ka}$ such that
  $f=\sum_{i\in\oi(\kappa)}f_{\lam_{i}}z^{\lam_{i}}$ and
  $f_{\lam_{i}}\neq 0$ whenever
  $\ord_{z}(f)\leq \lam_{i}\leq \deg_{z}(f)$.
  \begin{enumerate}
  \item Let us call $f$ almost palindromic $z$-polynomial,
    if each $f_{\lam_{i}}=f_{\lam_{\ka-i+1}}$.
  \item Let us call $f$ almost unimodal $z$-polynomial, if
    $f_{\lam_{i-1}}>f_{\lam_{i}}$ for some
    $i\in\oi(2,\ka-1)$ implies
        $f_{\lam_{j-1}}\geq f_{\lam_{j}}$ for each
    $j\in\oi(i+1,\ka)$.
  \item Let us call $f$ almost log-concave $z$-polynomial,
    if $f_{\lam_{i}}^{2}-f_{\lam_{i-1}}f_{\lam_{i+1}}\geq 0$
    for each $i\in\oi(2,\ka-1)$.
  \end{enumerate}
\end{definition}

\subsubsection{On the almost log-concavity}

\begin{conjecture}
  \label{conj:log-concavity-of-pre-merged-width-two}
  Let $\mu= 2$, $x=\iota^{\mu}(q)$, and
  $\del,\lam\in \Zgeo$.  Then, there is
  $H_{\del,\lam}\in \Zgeo$ such that
  the pre-merged determinant $\dt(U)_{a}^{b}(\mu,
  \iota^{\mu}(H_{\del,\lam}),\iota^{\mu}(t), e,x)$ is
   an almost log-concave $q$-polynomial for
  each $t\in\oi(\lam)$,
  strictly pre-fitting $\tp{\mu,a,b}$ with $b\leq \del$, and
  $e\in\oi(0,n_{a}^{b}(\mu,\iota^{\mu}(H_{\del,\lam})))$.
\end{conjecture}

\begin{example}
  The $q$-polynomial
  $\dt(U)_{\tp{1,2}}^{\tp{3,4}}(2,\iota^{2}(1),\iota^{2}(1),0,\iota^{2}(q))
    =
    {3 \brack 1 }_{q}{4 \brack 2 }_{q}
    -{3 \brack 2 }_{q}{4 \brack 1 }_{q}
    =q^{6} + q^{5} + 2q^{4} + q^{3} + q^{2}$
    is not log-concave.
    However,
    $\dt(U)_{\tp{1,2}}^{\tp{3,4}}(2,\iota^{2}(3),\iota^{2}(1),0,\iota^{2}(q))=
    {3 \brack 1 }^{3}_{q}{4 \brack 2 }^{3}_{q}
    -{3 \brack 2 }^{3}_{q}{4 \brack 1 }^{3}_{q}
    =q^{18} + 6q^{17} + 24q^{16} + 67q^{15} + 150q^{14} + 273q^{13} + 422q^{12} + 555q^{11} + 633q^{10} + 622q^{9} + 531q^{8} + 387q^{7} + 241q^{6} + 123q^{5} + 51q^{4} + 15q^{3}
    + 3q^{2}$ is log-concave.
  One can check that  setting $H_{20,10}=3$
  supports Conjecture~\ref{conj:log-concavity-of-pre-merged-width-two}.
\end{example}

We conjecture the following analog for pre-parcels, though
$q$-binomial coefficients are not necessarily log-concave.

\begin{conjecture}\label{conj:log-concavity-of-pre-parcel}
  Let $\lam,\del\in \Zgeo$.
  Then, there is $h_{\del,\lam}\in \Zgeo$ such that
  $U_{\tp{a_{1}}}^{\tp{b_{1}}}
    (1,\tp{h_{\del,\lam}},\tp{t},\tp{q})
    ={ b_{1}\brack a_{1}}_{q^{t}}^{h_{\del,\lam}}b_{b_{1}-a_{1},t}
    (q)^{h_{\del,\lam}}$
    is a log-concave $q$-polynomial
    for any $t\in\oi(\lam)$ and $1\leq a_{1}< b_{1}\leq \del$.
  \end{conjecture}

\begin{example}
  The pre-parcel
  $U_{\tp{2}}^{\tp{4}}(1,\tp{1},\tp{1},\tp{q})
  ={4 \brack 2 }_{q}
    =q^{4} + q^{3} + 2 q^{2} + q + 1.
  $ is not a log-concave $q$-polynomial.
  However,
  $U_{\tp{2}}^{\tp{4}}(1,\tp{3},\tp{1},\tp{q})
    {4 \brack 2}^{3}_{q}
    =q^{12} + 3 q^{11} + 9 q^{10} + 16 q^{9} + 27 q^{8} + 33
    q^{7} +
    38 q^{6} + 33q^{5} + 27 q^{4} + 16 q^{3} + 9 q^{2} +
    3q + 1
  $ is a log-concave $q$-polynomial.
  One can check that setting $h_{30,20}=3$ supports
  Conjecture~\ref{conj:log-concavity-of-pre-parcel}.
\end{example}

We also conjecture the following almost log-concavity
for higher-width pre-merged determinants.
\begin{conjecture}
  \label{conj:log-concavity-of-pre-merged-width-four}
  Let $\mu\in \Z_{\geq 4}$ and $\lam\in \Zgez$. If
  $\tp{\mu,a,b}$ is  almost strictly pre-fitting, then
  $\dt(U)_{a}^{b}(\mu,\iota^{\mu}(1),\iota^{\mu}(\lam),
    o_{a}^{b}(\mu,\iota^{\mu}(1)),
    \iota^{\mu}(q))$
  is an almost log-concave $q$-polynomial.
\end{conjecture} 

\begin{example}\label{ex:non-palindromic}
  We have $o_{\tp{0,2,3}}^{\tp{3,5,7}}(3,\iota^{3}(1))=
  o_{\tp{0,1,2,3}}^{\tp{3,4,5,7}}(4,\iota^{4}(1))=6$.
  Then,
  the following is not a log-concave $q$-polynomial:
  \begin{dmath*}
    \dt(U)_{\tp{0,2,3}}^{\tp{3,5,7}}(3,\iota^{3}(1),\iota^{3}(1),6,
    \iota^{3}(q))=
    q^{18}+ 2q^{17}+ 5q^{16}+ 9q^{15}+ 15q^{14}+ 21q^{13}+ 28q^{12}+ 33q^{11}+ 37q^{10}+ 38q^{9}+ 37q^{8}+ 33q^{7}+ 28q^{6}+ 21q^{5}+ 15q^{4}+ 9q^{3}+ 5q^{2} + 2q + 1.
  \end{dmath*}
  However,
  the following is a log-concave $q$-polynomial:
  \begin{dmath*}
    \dt(U)_{\tp{0,1,2,3}}^{\tp{3,4,5,7}}(4,\iota^{4}(1),\iota^{4}(1),6,
    \iota^{4}(q))=
    q^{21}+ 3q^{20}+ 8q^{19}+ 17q^{18}+ 31q^{17}+ 50q^{16}+ 74q^{15}+ 98q^{14}+ 121q^{13}+ 138q^{12}+ 147q^{11}+ 146q^{10}+ 137q^{9}+ 119q^{8}+ 97q^{7}+ 73q^{6}+ 50q^{5}+ 31q^{4}+ 17q^{3}+ 8q^{2} + 3q + 1.
  \end{dmath*}
\end{example}

\subsubsection{On the almost unimodality}

\begin{example}
  A pre-merged determinant does not have to be almost
  unimodal, because
  $\dt(U)_{\tp{0,2}}^{\tp{2,5}}(2,\iota^{2}(1), \iota^{2}(1),
  3,\iota^{2}(q))=
    {2 \brack 0}_{q}{5 \brack 2}_{q}
    -q^{3}{2 \brack 2}_{q}{5 \brack 0}_{q}
    =1+q+2q^{2}+q^{3}+2q^{4}+q^{5}+q^{6}$.
    Even with the trivial degree shift, 
   the following is not unimodal either:
  \begin{dmath}
\dt(U)_{\tp{0,6}}^{\tp{11,14}}(2,\iota^{2}(1), \iota^{2}(1),
  0,\iota^{2}(q)) =
    {11 \brack 0}_{q}{14 \brack 6}_{q}
    -{11 \brack 6}_{q}{14 \brack 0}_{q}
    =
    q^{48}+q^{47}+2 q^{46}+3 q^{45}+5 q^{44}+7 q^{43}+11 q^{42}+14 q^{41}+20 q^{40}+25 q^{39}+33 q^{38}+40 q^{37}+51 q^{36}+59 q^{35}+71 q^{34}+81 q^{33}+94 q^{32}+103 q^{31}+115 q^{30}+122 q^{29}+132 q^{28}+136 q^{27}+141 q^{26}+140 q^{25}+141 q^{24}+135 q^{23}+130 q^{22}+120 q^{21}+111 q^{20}+98 q^{19}+87 q^{18}+73 q^{17}+62 q^{16}+49 q^{15}+39 q^{14}+29 q^{13}+22 q^{12}+15 q^{11}+10 q^{10}+6 q^9+4 q^8+2 q^7+q^6.
    \label{eq:dense-not-unimodal}
  \end{dmath}
\end{example}
However, we conjecture the following  almost
unimodality for higher-width pre-merged determinants.

\begin{conjecture}\label{conj:unimodal-width-three}
  Let $\mu\in \Z_{\geq 3}$ and $\lam\in \Zgeo$ Consider an
  almost strictly pre-fitting $\tp{\mu,a,b}$.
  Then,
  $\dt(U)_{a}^{b}(\mu,\iota^{\mu}(1),\iota^{\mu}(\lam),0,\iota^{\mu}(q))$
   is an almost unimodal $q$-polynomial.
\end{conjecture}

\begin{example}
  By equation~\eqref{eq:dense-not-unimodal},
  $\dt(U)_{\tp{0,6}}^{\tp{11,14}}(2,\iota^{2}(1),
  \iota^{2}(1),0,\iota^{2}(q))$ is not  unimodal.  But,
  $\dt(U)_{\tp{0,2,6}}^{\tp{11,12,14}}(3, \iota^{3}(1),
  \iota^{3}(1),0,\iota^{3}(q))$ is unimodal, as it is
  \begin{dmath*}
    q^{68} + 2q^{67} + 5q^{66} + 9q^{65} + 17q^{64} +
    28q^{63} + 47q^{62} + 72q^{61} + 111q^{60} + 161q^{59}
    + 233q^{58} + 322q^{57} + 443q^{56} + 588q^{55} + 775q^{54} +
    993q^{53} + 1262q^{52} + 1565q^{51} + 1924q^{50} + 2315q^{49} + 2761q^{48} + 3230q^{47} + 3744q^{46} + 4263q^{45} + 4809q^{44} + 5335q^{43} + 5862q^{42} + 6338q^{41} + 6786q^{40} + 7153q^{39} + 7465q^{38} + 7671q^{37} + 7802q^{36} + 7813q^{35} + 7742q^{34} + 7552q^{33} + 7286q^{32} + 6917q^{31} + 6492q^{30} + 5992q^{29} + 5465q^{28} + 4897q^{27} + 4333q^{26} + 3762q^{25} + 3223q^{24} + 2705q^{23} + 2238q^{22} + 1810q^{21} + 1441q^{20} + 1118q^{19} + 853q^{18} + 631q^{17} + 458q^{16} + 320q^{15} + 219q^{14} + 143q^{13} + 91q^{12} + 54q^{11} + 31q^{10} + 16q^{9} + 8q^{8} + 3q^{7} + q^{6}.
  \end{dmath*}
\end{example}

Moreover, we conjecture the following almost unimodality for width-two pre-merged determinants.

\begin{conjecture}\label{conj:unimodal-width-two}
  Let $\mu=2$.  Then,
  $\dt(U)_{a}^{b}(\mu,\iota^{\mu}(1),\iota^{\mu}(1),0, \iota^{\mu}(q))$
  is an almost unimodal $q$-polynomial, provided
  $\tp{\mu,a,b}$ is strictly pre-fitting.
\end{conjecture}

Conjecture~\ref{conj:unimodal-width-two} does not contradict
equation~\eqref{eq:dense-not-unimodal}, as $1\leq a$ in
Conjecture~\ref{conj:unimodal-width-two}.

\subsubsection{On the almost palindromicity}
\begin{conjecture}\label{conj:palindromic-width-three}
  Let $\mu\in\Z_{\geq 3}$. 
  Then,
  $\dt(U)_{a}^{b}(\mu,\iota^{\mu}(1),\iota^{\mu}(1), 0,
    \iota^{\mu}(q))$
    is not an almost palindromic $q$-polynomial,
    provided $\tp{\mu,a,b}$ is almost strictly
  pre-fitting.
\end{conjecture}

In particular,
$\dt(U)_{\tp{a_{1},a_{2}}}^{\tp{b_{1},b_{2}}}(2,
\iota^{2}(1),\iota^{2}(1), 0, \iota^{2}(q))$ would not be palindromic
when $a_{2}-a_{1}\geq 2$ and $b_{2}-b_{1}\geq 2$, because these
$a_{1},a_{2},b_{1},b_{2}$ give rise to non-palindromic cases
of $\mu=3$ in Conjecture~\ref{conj:palindromic-width-three}.
Hence, when $\mu=2$, we conjecture the following, which
also provides infinitely many
almost palindromic unimodal $q$-polynomials.

\begin{conjecture}\label{conj:non-palindromic-seq}  
  For $\lam\in \Zgez$, let
  $M(\lam)
  =\{\mu\in \Z^{4}\mid 0\leq \mu(1)<\mu(2)\leq \mu(3)<\mu(4)=\lam\}$ and
  $N(\lam)
  =\{\mu\in M(\lam)\mid 
      \dt(U)_{\mu(1,2)}^{\mu(3,4)}(2,\iota^{2}(1),
      \iota^{2}(1), 0,\iota^{2}(q)) 
      \mbox{ is not almost palindromic.} \}$.
  Then,
  $O(\lam)=\frac{\num(N(\lam))}{\num(M(\lam))}$
  satisfies 
  $0=O(3) <O(4)  <O(5) <\dots $.
\end{conjecture}

On base shifts,
we have the following width-two pre-merged
determinants.
\begin{example}
  The pre-merged determinant
  $\dt(U)_{\tp{0,1}}^{\tp{1,4}}(2,\iota^{2}(1),\iota^{2}(1),
    1,\iota^{2}(q))
    =q^{3} + q^{2} + 1$
    is palindromic.  However,
  $\dt(U)_{\tp{0,1}}^{\tp{1,4}}(2,\iota^{2}(1),\iota^{2}(2),
    1,\iota^{2}(q))=
    q^{13} + 2q^{12} + 2q^{11} + 4q^{10} + 5q^{9} + 5q^{8}
    +       6q^{7} + 6q^{6} + 5q^{5} + 4q^{4} + 4q^{3} +
    2q^{2} + q + 1$
    is not palindromic.
\end{example}

Still, we conjecture the following palindromic transitivity
on base shifts 
by width-three pre-merged
determinants.

 \begin{conjecture}\label{conj:transitivity-palindromicity}
   For $\mu\in\Z_{\geq 3}$, let
    $\tp{\mu,a,b}$ be almost strictly pre-fitting.
  Suppose $\lam\in \Zgeo$ and
  $e\in\oi(0,n_{a}^{b}(\mu,\iota^{\mu}(1)))$.  Then,
  $\dt(U)_{a}^{b}(\mu,\iota^{\mu}(1),\iota^{\mu}(\lam),e,
    \iota^{\mu}(q))$
  is an almost palindromic $q$-polynomial if and only
  if
  $
    \dt(U)_{a}^{b}(\mu,\iota^{\mu}(1),  \iota^{\mu}(\lam+1),e,\iota^{\mu}(q))$
  is an almost palindromic $q$-polynomial.
\end{conjecture}


\section{Hadamard products of parcels}
\label{sec:hadam}

We introduce external and internal Hadamard products of
parcels.  They construct merged-log-concave parcels of
higher widths and weights from given ones.

\subsection{External Hadamard products}
\label{sec:external}

External Hadamard products increase
the widths of parcels.  
To introduce the products,
we use the following truncations of
fitting tuples.

\begin{definition}\label{def:truncations}
  Assume a fitting $\mu=\tp{s,\llo,\mn,\nn,\kn}$. Let
  $\lam\in \tei(2,\llo)$ and $\llt=\lamt-\lamo+1$.
  \begin{enumerate}
  \item We define the truncation
    $t(\lam,\mu)= \tp{s,\llt,\mt,\nt,\kt}$ such that
    \begin{align}
      \mt&=\mn(\lamo,\lamt),
           \label{eq:truncations-m}\\
      \nt&=\nn(\llo-\lamt+1,\llo-\lamo+1),
           \label{eq:truncations-n}
      \\
      k_{2,1}&=\sum \kn(1,\lamo),
               \label{eq:truncations-k-left-sum}\\
      k_{2,i}&=k_{1,\lamo+i-1} \mfor      i\in\oi(2,\llt),
               \label{eq:truncations-k-left}\\
      k_{2,\llt+1}&=\sum \kn(\lamt+1,2\llo-\lamt+1),
                    \label{eq:truncations-k-right-sum}\\
      k_{2,\llt+i}&=k_{1,2\llo-\lamt+i} \mfor i\in\oi(2,\llt).
                    \label{eq:truncations-k-right}
    \end{align}
  \item
    \label{d:truncations-off-centered}
    If $\lamo=1$, then we define the off-centered truncation
    $\otr(\llt,\mu)=   t(\lam,\mu)$.
  \item
    If $\lamt=\llo$, then we define the centered truncation
    $\ctr(\llt,\mu)=  t(\lam,\mu)$.
  \end{enumerate}
\end{definition}

Let us state
the following properties of
these truncations.

\begin{proposition}\label{prop:fitting-trunc}
  For a fitting $\muo=\tp{s,\llo,\mn,\nn,\kn}$, consider
  $\mut=\tp{s,\llt,\mt,\nt,\kt}= t(\lam,\mu)$.  Let
  $a_{i}=\nu(k_{i})$ and $b_{i}=\nu(m_{i},n_{i},k_{i})$ for
  $i\in\oi(2)$.  Then, we have the following.
  \begin{enumerate}
  \item\label{c:fitting-trunc-ab-equations}
    There exist the following equations:
    \begin{align}
      a_{2}(1,\llt)
      &=a_{1}(\lamo,\lamt);
        \label{eq:fitting-trunc-a-left}\\
      a_{2}(\llt+1,2\llt)
      &=a_{1}(2\llo-\lamt+1,2\llo-\lamo+1);
        \label{eq:fitting-trunc-a-right}\\
      b_{2}(1,\llt)
      &=b_{1}(\lamo,\lamt);
        \label{eq:fitting-trunc-b-left}\\
      b_{2}(\llt+1,2\llt)
      &=b_{1}(2\llo-\lamt+1,2\llo-\lamo+1).
        \label{eq:fitting-trunc-b-right}
    \end{align}
  \item \label{c:fitting-trunc-fitting} $\mut$ is fitting.
  \item \label{c:fitting-trunc-sigma-equations}
    There exist the following  equations:
    \begin{align}
      \sig(\kt)
      &=\sig(\kn)(\lamo,\lamt);
        \label{eq:fitting-trunc-sigma}\\
      \mt\bom \kt
      &=(\mn\bom \kn)(\lamo,\lamt);
        \label{eq:fitting-trunc-bom}\\
      (\nt\bop \kt)^{\ve}
      &=(\nn\bop \kn)^{\ve}(\lamo,\lamt).
        \label{eq:fitting-trunc-bop}
    \end{align}
  \item \label{c:fitting-trunc-wrapped} If $\muo$ is
    wrapped, then $\mut$ is wrapped.
  \end{enumerate}
\end{proposition}
\begin{proof}
  Let us prove Claim~\ref{c:fitting-trunc-ab-equations}.
  Since equation~\eqref{eq:truncations-k-left-sum} gives
  $a_{2,1}=\sum \kn(1,\lamo)=a_{1,\lamo}$,
  equation~\eqref{eq:fitting-trunc-a-left} follows from
  equation~\eqref{eq:truncations-k-left} by
  $\lamo+\llt-1= \lamt$. Also,
  equation~\eqref{eq:truncations-k-right-sum} gives
  $a_{2,\llt+1}=a_{1,2\llo-\lamt+1}$.  Hence,
  equation~\eqref{eq:truncations-k-right} yields
  equation~\eqref{eq:fitting-trunc-a-right}.  Since
  $b_{2}=a_{2}+\mt\ccn \nt$,
  equations~\eqref{eq:fitting-trunc-b-left}
  and~\eqref{eq:fitting-trunc-b-right} follow from
  equations~\eqref{eq:truncations-m},
  \eqref{eq:truncations-n}, \eqref{eq:fitting-trunc-a-left},
  and~\eqref{eq:fitting-trunc-a-right}.
  
  Let us prove Claim~\ref{c:fitting-trunc-fitting}.  Since
  $\mn,\nn\ldZlo s$, equations~\eqref{eq:truncations-m}
  and~\eqref{eq:truncations-n} imply the inclusion condition
  $\mt,\nt\ldZlt s$ of $\mut$.  Moreover, $\mut$ satisfies
  the slope conditions by
  equations~\eqref{eq:fitting-trunc-a-left},
  \eqref{eq:fitting-trunc-a-right},
  \eqref{eq:fitting-trunc-b-left},
  and~\eqref{eq:fitting-trunc-b-right}, because
  $a_{2,\llt}=a_{1,\lamt}< a_{1,2\llo-\lamt+1}=a_{2,\llt+1}$
  and
  $b_{2,\llt}=b_{1,\lamt}<
  b_{1,2\llo-\lamt+1}=b_{2,\llt+1}$.
                
  Let us prove Claim~\ref{c:fitting-trunc-sigma-equations}.
  First, we obtain equation~\eqref{eq:fitting-trunc-sigma},
  because by equations~\eqref{eq:truncations-k-left},
  \eqref{eq:truncations-k-right-sum},
  and~\eqref{eq:truncations-k-right}, each
  $i\in\oi(\llt)$ satisfies
  \begin{align}
    \sig(\kt)_{i}
    &= \sum \kt(i+1,2\llt-i+1)\\
    &= \sum \kt(i+1,\llt+(\llt-i+1))\\
    &= \sum \kn(i+\lamo,2\llo-\lamt+(\llt-i+1))\\
    &=\sum \kn((i+\lamo-1)+1,2\llo-(i+\lamo-1)+1)\\
    &=\sig(\kn)_{i+\lamo-1}.
  \end{align}
  Second, equation~\eqref{eq:fitting-trunc-bom} holds,
  because for each $i\in\oi(\llt)$,
  equation~\eqref{eq:fitting-trunc-sigma} gives
  \begin{align}
    (\mt\bom \kt)_{i}
    &=m_{2,i}-\sig(\kt)_{i}\\
    &=m_{1,i+\lamo-1}-\sig(\kn)_{i+\lamo-1}\\
    &=(\mn\bom \kn)_{i+\lamo-1}.
  \end{align}
  Third, equation~\eqref{eq:fitting-trunc-bop} holds as follows.
  For each $i\in\oi(\llt)$,
  equation~\eqref{eq:fitting-trunc-sigma} gives
  \begin{align}
    (\nt\bop \kt)_{i}
    &=n_{2,i}+\sig(\kt)_{\llt-i+1}\\
    &=n_{1,\llo-\lamt+i}+\sig(\kn)_{\llt-i+1+\lamo-1}\\
    &=n_{1,\llo-\lamt+i}+\sig(\kn)_{\llo-(\llo-\lamt+i)+1}\\
    &=(\nn\bop \kn)_{\llo-\lamt+i}.
  \end{align}
  Thus, for each $i\in\oi(0,\llt-1)$, we obtain
  \begin{align}
    (\nn\bop \kn)^{\ve}_{\lamo+i}
    &=(\nn\bop \kn)_{\llo-(\lamo+i)+1}\\
    &=(\nn\bop \kn)_{\llo-\lamt+\llt-i}\\
    &=(\nt\bop \kt)_{\llt-i}\\
    &=(\nt\bop \kt)^{\ve}_{i+1}.
  \end{align}
  Therefore, Claim~\ref{c:fitting-trunc-sigma-equations}
  follows.
  
  Claim~\ref{c:fitting-trunc-wrapped} follows from
  Claim~\ref{c:fitting-trunc-sigma-equations}.
\end{proof}

Now, we have the following corollary for
off-centered/centered truncations.

\begin{corollary}\label{cor:trunc}
  Let $\ka\in \Zgeo^{3}$ such that
  $\ka_{3}=\ka_{1}+\ka_{2}$.  For a fitting
  $\mur=\tp{s,\ka_{3},\mr,\nr,\kr}$, consider
  \begin{align}
    \muo
    &=\otr(\ka_{1},\mur)
      =\tp{s,\ka_{1},\mn,\nn,\kn},\\
    \mut
    &=\ctr(\ka_{2},\mur)
      =\tp{s,\ka_{2},\mt,\nt,\kt}.
  \end{align}
  Let $a_{i}=\nu(k_{i})$ and $b_{i}=\nu(m_{i},n_{i},k_{i})$
  for $i\in\oi(3)$.  Also, suppose
  $u_{i}\in \Zgez^{\ka_{i}}$ and an indeterminate
  $x_{i}\in \Q(\fX)^{\ka_{i}}$ for $i\in\oi(3)$.  Then, we
  have the following.
  \begin{enumerate}
  \item  $\muo$ is fitting.
    \label{c:trunc-off-centered-fitting}
  \item\label{c:trunc-off-centered-eq}
    There exist the following equations:  
    \begin{align}
      \mn
      &=\mr(1,\ka_{1});
        \label{eq:trunc-off-centered-m}\\
      \nn
      &= \nr
        (\ka_{3}-\ka_{1}+1,\ka_{3});
        \label{eq:trunc-off-centered-n}\\
      a_{1}
      &=a_{3}(1,\ka_{1}) \ccn
        a_{3}(2\ka_{3}-\ka_{1}+1,2\ka_{3});
             \label{eq:trunc-off-centered-a}\\
      b_{1}
      &=b_{3}(1,\ka_{1})  \ccn
        b_{3}(2\ka_{3}-\ka_{1}+1,2\ka_{3});
             \label{eq:trunc-off-centered-b}\\
      \sig(\kn)
      &=
        \sig(\kr)(1,\ka_{1});
        \label{eq:trunc-off-centered-sigma}\\
      \mn\bom \kn
         &=(\mr\bom \kr)(1,\ka_{1});
           \label{eq:trunc-off-centered-bom}\\
      (\nn\bop \kn)^{\ve}
         &= (\nr\bop \kr)^{\ve}(1,\ka_{1}).
           \label{eq:trunc-off-centered-bop}
    \end{align}
  \item
    \label{c:trunc-off-centered-binom}
    Moreover, we have
    \begin{align}
      {b_{1} \brack  a_{1}}^{u_{1}^{\wcn}}_{x_{1}^{\wcn}}
      &=   {b_{3}(1,\ka_{1}) \brack
        a_{3}(1,\ka_{1})}^{u_{1}}_{x_{1}}
        {b_{3}(2\ka_{3}-\ka_{1}+1,2\ka_{3})\brack
        a_{3}(2\ka_{3}-\ka_{1}+1,
        2\ka_{3})}^{u_{1}^{\ve}}_{x_{1}^{\ve}},
        \label{eq:trunc-off-centered-binom-prod}\\
      {b_{1} \brack
      a^{\ve}_{1}}^{u_{1}^{\wcn}}_{x_{1}^{\wcn}}
      &=  {b_{3}(1,\ka_{1}) \brack
        a_{3}( 2\ka_{3}-\ka_{1}+1, 2\ka_{3})^{\ve}}^{u_{2}}_{x_{2}}
        {b_{3}(2\ka_{3}-\ka_{1}+1,2\ka_{3}) \brack
        a_{3}(1,\ka_{1})^{\ve}}^{u_{1}^{\ve}}_{x_{1}^{\ve}}.
        \label{eq:trunc-off-centered-binom-prod-wed}
    \end{align}
  \item  $\mut$ is fitting.  \label{c:trunc-centered-fitting}
  \item\label{c:trunc-centered-eq}
    There exist the following equations:  
    \begin{align}
      \mt&=\mr(\ka_{1}+1,\ka_{3});
           \label{eq:trunc-centered-m}\\
      \nt&=\nr(1,\ka_{3}-\ka_{1});
           \label{eq:trunc-centered-n}\\
      a_{2}&=a_{3}(\ka_{1}+1,2\ka_{3}-\ka_{1});
             \label{eq:trunc-centered-a}\\
      b_{2}&=b_{3}(\ka_{1}+1,2\ka_{3}-\ka_{1});
             \label{eq:trunc-centered-b}\\
      \sig(\kt)
         &=\sig(\kr)(\ka_{1}+1,\ka_{3});
           \label{eq:trunc-centered-sigma}\\
      \mt\bom \kt
         &=      (\mr\bom \kr)(\ka_{1}+1,\ka_{3});
           \label{eq:trunc-centered-bom}\\
      (\nt\bop \kt)^{\ve}
         &= (\nr\bop \kr)^{\ve}(\ka_{1}+1,\ka_{3}).
           \label{eq:trunc-centered-bop}
    \end{align}
  \item \label{c:trunc-centered-binom}
    We have
    \begin{align}
      {b_{2} \brack a_{2}}^{u_{2}^{\wcn}}_{x_{2}^{\wcn}}
      &=
        {b_{3}(\ka_{3}+1,\ka_{3}) \brack
        a_{3}(\ka_{1}+1,\ka_{3})}^{u_{2}}_{x_{2}}
        {b_{3}(\ka_{3}+1,2\ka_{3}-\ka_{1}) \brack
        a_{3}(\ka_{3}+1,
        2\ka_{3}-\ka_{1})}^{u_{2}^{\ve}}_{x_{2}^{\ve}},
        \label{eq:trunc-centered-binom-prod}\\
      {b_{2} \brack a^{\ve}_{2}}^{u_{2}^{\wcn}}_{x_{2}^{\wcn}}
      &=
        {b_{3}(\ka_{1}+1,\ka_{3}) \brack
        a_{3}(\ka_{3}+1,
        2\ka_{3}-\ka_{1})^{\ve}}^{u_{2}}_{x_{2}}
        {b_{3}(\ka_{3}+1,2\ka_{3}-\ka_{1}) \brack
        a_{3}(\ka_{1}+1,
        \ka_{3})^{\ve}}^{u_{2}^{\ve}}_{x_{2}^{\ve}}.
        \label{eq:trunc-centered-binom-prod-wed}
    \end{align}
  \item \label{c:trunc-binom}
    Moreover, we have
    \begin{align}
      {b_{3} \brack
      a_{3}}^{(u_{1}\ccn u_{2})^{\wcn}}_{
      (x_{1}\ccn x_{2})^{\wcn}}
      &=
        {b_{1} \brack a_{1}}^{u_{1}^{\wcn}}_{x_{1}^{\wcn}}
        {b_{2} \brack a_{2}}^{u_{2}^{\wcn}}_{x_{2}^{\wcn}},
        \label{eq:trunc-binom-prod}
      \\
      {b_{3} \brack a^{\ve}_{3}}^{(u_{1}\ccn u_{2})^{\wcn}}_{
      (x_{1}\ccn x_{2})^{\wcn}}
      &=
        {b_{1} \brack a^{\ve}_{1}}^{u_{1}^{\wcn}}_{x_{1}^{\wcn}}
        {b_{2} \brack a^{\ve}_{2}}^{u_{2}^{\wcn}}_{x_{2}^{\wcn}}.
        \label{eq:trunc-binom-prod-wed}
    \end{align}
  \item
    \label{c:trunc-wrapped}
    If $\mur$ is wrapped, then $\muo$ and
    $\mut$ are wrapped.
  \end{enumerate}
\end{corollary}
\begin{proof}
  Let $\lamo=1$ and $\lamt=\ka_{1}$ in
  Proposition~\ref{prop:fitting-trunc}.  Then,
  Claim~\ref{c:trunc-off-centered-fitting} follows from
  Claim~\ref{c:fitting-trunc-fitting} of
  Proposition~\ref{prop:fitting-trunc}, since
  $\lamt-\lamo+1=\ka_{1}$.  Also, we have
  Claim~\ref{c:trunc-off-centered-eq} as follows.
  Item~\ref{d:truncations-off-centered} of
  Definition~\ref{def:truncations} gives
  equations~\eqref{eq:trunc-off-centered-m}
  and~\eqref{eq:trunc-off-centered-n}.
  Equations~\eqref{eq:trunc-off-centered-a}
  and~\eqref{eq:trunc-off-centered-b} hold by
  Claim~\ref{c:fitting-trunc-ab-equations} of
  Proposition~\ref{prop:fitting-trunc}.  Moreover, we obtain
  equations~\eqref{eq:trunc-off-centered-sigma},
  \eqref{eq:trunc-off-centered-bom},
  and~\eqref{eq:trunc-off-centered-bop} by
  Claim~\ref{c:fitting-trunc-sigma-equations} of
  Proposition~\ref{prop:fitting-trunc}.
  Claim~\ref{c:trunc-off-centered-binom} follows from
  equations~\eqref{eq:trunc-off-centered-a}
  and~\eqref{eq:trunc-off-centered-b}.
  
  Similarly, let $\lamo=\ka_{3}-\ka_{2}+1$ and
  $\lamt=\ka_{3}$ in Proposition~\ref{prop:fitting-trunc}.
  Then, we obtain Claims~\ref{c:trunc-centered-fitting},
  \ref{c:trunc-centered-eq},
  and~\ref{c:trunc-centered-binom} as above, since
  $\lamt-\lamo+1=\ka_{2}$.

  Claim~\ref{c:trunc-binom} follows from
  Claims~\ref{c:trunc-off-centered-binom}
  and~\ref{c:trunc-centered-binom}.  Furthermore,
  Claim~\ref{c:trunc-wrapped} holds by
  Claim~\ref{c:fitting-trunc-wrapped} of
  Proposition~\ref{prop:fitting-trunc}.
\end{proof}

We state the following compatibility on squaring
orders.

\begin{lemma}\label{lem:succ-on-coord-unions}
  Suppose $\fXr=\fXo\cup \fXt$.  Then, squaring orders
  $\{\geAfXr,\gAfXr\}$ on $\fXr$ are
  compatible to squaring orders $\{\geAfXo,\gAfXo\}$ on
  $\fXo$.
\end{lemma}
\begin{proof}
  If $f>_{\AfXo}0$, then $f>_{A_{\fXr}}0$ by
  $\fXr=\fXo\cup \fXt$.  A similar argument holds for
  $\geAfXo$ and $\geq_{A_{\fXr}}$.
\end{proof}

Then,  we define the external Hadamard
products of parcels as parcels by the following proposition.

\begin{proposition}\label{prop:ext-hadam}
  Suppose parcels
  $\cF_{i}
  =\Lam(s,l_{i},\scc_{i},w_{i},f_{i,s},\phi_{i},\rho_{i},x_{i},\fX_{i})$
  for $i\in\oi(2)$.  For $\fXr=\fXo\cup \fXt$, consider
  $f_{3,s}=\{f_{3,s,m}=
  f_{1,s,m(1,\llo)}f_{2,s,m(\llo+1,\llr)}\in \Q(\fXr)\}_{m\in
    \Zlr}$.  Let $\llr=\llo+\llt$, $w_{3}=w_{1}\ccn w_{2}$,
  $x_{3}=x_{1}\ccn x_{2}$, $\rhor=\rhoo\ccn \rhot$, and
  $\phir(x_{3})=\phio(x_{1})\ccn \phit(x_{2})$.  For
  $O_{i}=\{\sce_{i},\scc_{i}\}$ of $i\in \oi(2)$, assume
  squaring orders $O_{3}=\{\scer,\sccr\}$ on $\fXr$ such
  that $O_{3}\Sup O_{1},O_{2}$.
  \begin{enumerate}
  \item
    \label{c:ext-hadam-parcel}
    Then, there exists a parcel
    $\cFr =\Lam(s,\llr,\sccr,w_{3},f_{3,s},
    \phir,\rhor,x_{3},\fXr)$.
  \item
    \label{c:ext-hadam-merged}
    For a fitting $\mur= \tp{s,\llr,\mr,\nr,\kr}$, let
    $\muo=\tp{s,\llo,\mn,\nn,\kn}= \otr(\llo,\mur)$ and
    $\mut=\tp{s,\llt,\mt,\nt,\kt}= \ctr(\llt,\mur)$.
    Then, we have
    \begin{dmath*} 
    \Delta(\cFr)(s,\llr,w_{3},\mr,\nr,
    \kr,\phir,\rhor,x_{3},\fXr) = \prod_{i\in\oi(2)}
    \Delta_{L}(\cF_{i})(s,l_{i},w_{i},m_{i},n_{i},
    k_{i},\phi_{i},\rho_{i},x_{i},\fX_{i}) - \prod_{i\in\oi(2)}
    \Delta_{R}(\cF_{i})(s,l_{i},w_{i},m_{i},n_{i},
    k_{i},\phi_{i},\rho_{i},x_{i},\fX_{i}).
    \end{dmath*}
  \end{enumerate}
\end{proposition}
\begin{proof} 
  Let us prove Claim~\ref{c:ext-hadam-parcel}.  First, we
  confirm that $x_{3}$ is $\sccr$-admissible.  If
  $g>_{x_{3,i}}0$ for some $i\in\oi(\llr)$, then $g\scco 0$
  or $g\scct 0$ by $x_{3}=x_{1}\ccn x_{2}$ and the half
  $>_{x_{j,i}}$-$\scc_{j}$ implications of $j\in\oi(2)$.
  This yields $f\sccr0 $ by $O_{3}\Sup O_{1}, O_{2}$.  Thus,
  each $i\in\oi(\llr)$ gives the half $>_{x_{3,i}}$-$\sccr$
  implication.  Also, each $i\in\oi(\llr)$ satisfies the
  upper condition of $x_{3,i}$ on $\fXr$ by
  $\fXr=\fXo\cup \fXt$ and
  Lemma~\ref{lem:succ-on-coord-unions}.  Therefore, $x_{3}$
  is $\sccr$-admissible.

  Second, we prove that $f_{3,s}$ is pairwise
  $\tp{s,\llr,\sccr}$-positive.  If $m,n\ldZlr s$, then
  $f_{3,s,m}f_{3,s,n}= f_{1,s,m(1,\llo)}
  f_{2,s,m(\llo+1,\llr)} f_{1,s,n(1,\llo)}
  f_{2,s,n(\llo+1,\llr)} \sccr 0$, since
  $m(1,\llo),n(1,\llo) \ldZlo s$ and
  $m(\llo+1,\llr),n(\llo+1,\llr)\ldZlt s$.  Also, if
  $m\nld_{\Zlr} s$, then
  $f_{3,s,m}=f_{1,s,m(1,\llo)} f_{2,s,m(\llo+1,\llr)} =0$,
  since $m(1,\llo) \nld_{\Zlo}s$ or
  $m(\llo+1,\llr)\nld_{\Zlt}s$.  Thus, $f_{3,s}$ is pairwise
  $\tp{s,\llr,\sccr}$-positive.
  
  Third, we prove that $\phir$ is a $\lamr$-mediator for
  $\lamr=\tp{s,\llr,w_{3},\sccr,\rhor,x_{3},\fXr}$.  Since
  $\cFo$ and $\cFt$ are parcels, $\phir$ and $\lamr$ have
  the base positivity by
  $\phir(x_{3})= \phio(x_{1}) \ccn \phit(x_{2})$ and
  Lemma~\ref{lem:succ-on-coord-unions}.  Moreover, let
  $m\ldZlr s_{3}$.  Then, $m(1,\llo)\ldZlo s$ and
  $m(\llo+1,\llr) \ldZlt s$ give
  \begin{align}
    B(s,\llr,w_{3},m,\phir,\rhor,x_{3},\fXr)
    &=
      B(s,\llo,w_{1},m(1,\llo),\phio,\rhoo,x_{1},\fXo)
    \\ &\cdot
         B(s,\llt,w_{2},m(\llo+1,\llr),\phit,\rhot,x_{2},\fXt)
    \\& \sccr 0. 
  \end{align}
  Thus, $\phir$ and $\lamr$ have the base-shift positivity.
  In particular, since $\phir$ is a $\lamr$-mediator,
  Claim~\ref{c:ext-hadam-parcel} follows.
  
  Let us prove Claim~\ref{c:ext-hadam-merged}.  Let
  $a_{3}=\nu(k_{3})$ and $b_{3} =\nu(m_{3},n_{3},k_{3})$.
    Then,
  Claim~\ref{c:merged-binom-bshift-left-right-equations} of
  Theorem~\ref{thm:merged-binom-bshift} gives
  \begin{dmath*}
    \Delta(\cFr)(s,\llr,w_{3},\mr,\nr,\kr,
    \phir,\rhor,x_{3},\fXr)
    =f_{3,s,\mr} f_{3,s,\nr^{\ve}} 
    B(s,\llr,w_{3},\mr,\nr^{\ve},
    \phir,\rhor,x_{3},\fXr)
    {b_{3} \brack a_{3}}^{w_{3}^{\wcn}}_{x_{3}^{\wcn}}
    -
    f_{3,s,\mr \bom \kr}
    f_{3,s,(\nr\bop \kr)^{\ve}}
    B(s,\llr,w_{3},\mr\bom \kr,
    (\nr\bop \kr)^{\ve},
    \phir,\rhor,x_{3},\fXr)
    {b_{3} \brack
      a_{3}^{\ve}}^{w_{3}^{\wcn}}_{x_{3}^{\wcn}}.
  \end{dmath*}
  Hence, Claims~\ref{c:trunc-off-centered-eq}
  and~\ref{c:trunc-centered-eq}
  of
  Corollary~\ref{cor:trunc} give
  Claim~\ref{c:ext-hadam-merged} by
   Item~\ref{d:bshift-fun-multi-two-tuples} of  
   Definition~\ref{def:bshift-fun-multi} and
   Claim~\ref{c:trunc-binom} of Corollary~\ref{cor:trunc}.
 \end{proof}


\begin{definition}\label{def:ext-hadam}
  Under the assumption of Proposition~\ref{prop:ext-hadam},
  we define the external Hadamard product $\cFo\ehd \cFt$ as
  the parcel
  $\cFr=\Lam(s,\llr,\sccr,w_{3},f_{3,s},
  \phir,\rhor,x_{3},\fXr)$.
\end{definition}

Hence, $\cF_{3,m}=\cF_{1,m(1,\llo)}\cF_{2,m(\llo+1\llr)}$
for $m\in \Zlr$, because $m\ld_{\Zlr}s$ implies
\begin{align}
  \cF_{3,m}
  &=\frac{f_{3,s,m}}{
    \prod \phi(x_{3})^{m\rc w_{3}}
    \cdot [m]!_{x_{3}}^{w_{3}}}\\
  &=\frac{f_{1,s,m(1,\llo)}}{
    \prod \phi(x_{1})^{m(1,\llo)\rc w_{1}}
    \cdot [m(1,\llo)]!_{x_{1}}^{w_{1}}}
    \cdot
    \frac{f_{2,s,m(\llo+1,\llr)}}{
    \prod \phi(x_{2})^{m(\llo+1,\llr)\rc w_{2}}
    \cdot [m(\llo+1,\llr)]!_{x_{2}}^{w_{2}}}.
\end{align}

Moreover,
 external
Hadamard products have
 the following merged-log-concavity.

\begin{theorem}\label{thm:merged-ext-hadam}
  Let
  $\cF_{i}
  =\Lam(s,l_{i},\scc_{i},w_{i},f_{i,s},\phi_{i},\rho_{i},x_{i},\fX_{i})$
  for $i\in\oi(3)$ such that $\cFr=\cFo\ehd \cFt$.  Assume
  squaring orders $O_{i}'=\{\sce_{i}',\scc'_{i}\}$ on
  $\fX_{i}$ for $i\in\oi(3)$ such that
  $O'_{3}\Sup O'_{1}, O'_{2}$.  Let $\cFt$ be
  $\scet'$-merged-log-concave.  Then, we have the following.
  \begin{enumerate}
  \item
    \label{c:merged-ext-hadam-scc}
    $\cFr$ is $\sccr'$-merged-log-concave, if
    $\cFo$ is $\scco'$-merged-log-concave.
  \item
    \label{c:merged-ext-hadam-sce}
    $\cFr$ is
    $\scer'$-merged-log-concave, if $\cFo$ is
    $\sceo'$-merged-log-concave.
  \end{enumerate}
\end{theorem}
\begin{proof}
  We prove Claim~\ref{c:merged-ext-hadam-scc}.  Suppose a
  fitting $\mur=\tp{s,\llr,\mr,\nr,\kr}$. Then,
  Claims~\ref{c:trunc-off-centered-fitting}
  and~\ref{c:trunc-centered-fitting} of
  Corollary~\ref{cor:trunc} give fitting
  $\tp{s,\llo,\mn,\nn,\kn}=\otr(\llo,\mur)$ and
  $\tp{s,\llt,\mt,\nt,\kt}=\ctr(\llt, \mur)$.  Also, for
  $i\in\oi(2)$, let
  \begin{align}
    L_{i}
    &=
      \Delta_{L}(\cF_{i})(s,l_{i},w_{i},m_{i},n_{i},
      k_{i},\phi_{i},\rho_{i},x_{i},\fX_{i}),\\
    R_{i}
    &=
      \Delta_{R}(\cF_{i})(s,l_{i},w_{i},m_{i},n_{i},
      k_{i},\phi_{i},\rho_{i},x_{i},\fX_{i}).
  \end{align}
  
  Then, by the merged-log-concavity of $\cFo$ and $\cFt$, 
  we obtain
  \begin{align}
    L_{1}\scco' R_{1} \mand L_{2} \scet' R_{2}.
  \end{align}
  Moreover, we have $L_{1} \scco 0$ and $R_{1}\sceo 0$ by
  Claim~\ref{c:merged-binom-bshift-scc-eq} of
  Theorem~\ref{thm:merged-binom-bshift} for $\cFo$ on
  $O_{1}$.  Also, we have $L_{2} \scct 0$, and
  $R_{2}\scct 0$ or $R_{2}=0$ by that for $\cFt$ on $O_{2}$.
  Hence, because $O_{3}'\Sup O_{1}', O_{2}'$ with
  $O'_{1}\Sup O_{1}$ and $O'_{2}\Sup O_{2}$, we obtain
  \begin{align}
    L_{1}&\sccr' R_{1} \scer' 0 \mand
           L_{2}\scer' R_{2},\\
    L_{1}&\sccr' 0 \mand L_{2}\sccr' 0,\\
    R_{2}&\sccr' 0 \mor R_{2}=0.  
  \end{align}
  In particular, Corollary~\ref{cor:pair-orders} gives
  Claim~\ref{c:merged-ext-hadam-scc}.

  Claim~\ref{c:merged-ext-hadam-sce} holds by
  Claim~\ref{c:succ-succeq-zero-zero-1} of
  Lemma~\ref{lem:succ-succeq}, since
  $L_{i}\scer' R_{i}\scer ' 0$ for $i\in\oi(2)$ in the above.
\end{proof}

\begin{remark}\label{rmk:not-ext-hadam}
  Under the assumption of
  Theorem~\ref{thm:merged-ext-hadam}, suppose $\fXo=\fXt$,
  $O_{1}=O_{2}=O_{3}$, and canonical mediators $\phio,\phit$
  for simplicity.  Then, Theorem~\ref{thm:merged-ext-hadam}
  gives higher-width merged-log-concave parcels from
  lower-width ones. But, a width-two monomial parcel does not
  have to be the external Hadamard product of  width-one
  monomial parcels.

  For instance, let $\llo=\llt=1$ and $\llr=2$.  Also,
  assume monomial indices $\tp{l_{i},w_{i},\gam_{i}}$ for
  $i\in \oi(3)$ such that $\gam_{3,1,1}>0>\gam_{3,2,1}$.
  Then, $x_{1}=x_{2}=\iota^{\llo}(q)$ and
  $x_{3}=\iota^{\llr}(q)$ imply
  \begin{align}
    \Lam(s,\llr,w_{3},\scc,\Psi_{s,\gamr,q},x_{3},\fXr)
    & \neq
      \Lam(s,\llo,w_{1},\scc,\Psi_{s,\gamo,q},
      x_{i},\fXo)
      \ehd
      \Lam(s,\llt,w_{2},\scc,\Psi_{s,\gamt,q},
       x_{2},\fXt),
    \end{align}
    because $\gam_{1,1,1},\gam_{2,1,1}\geq 0$ in the sum
    monomial conditions of $\tp{\llo,w_{1},\gamo}$ and
    $\tp{\llt,w_{2},\gamt}$.
\end{remark}

We also define the following multifold products 
for our later discussion.
\begin{definition}\label{def:multi-fold-ext-hadam}
  Let $\cF =\Lam(s,l,\scc,w,\fs,\phi,\rho,x,\fX)$ and
  $d\in \Zgeo$. Then, we write the $d$-fold external Hadamard
  product
  $\cF^{\ehd d}= \Lam(s,d l,\scc,w^{\ccn d},\gs, \phi^{\ccn
    d},\rho^{\ccn d},x^{\ccn d},\fX)$ such that
  $g_{s,m}=\prod_{\lam\in\oi(d)} f_{s,m((\lam-1)l+1,\lam l)}$ for
  each $m\in \Z^{d l}$.
\end{definition}

\subsection{Reduction to width-one
  merged-log-concave
  parcels}\label{sec:reduct-to-width-one}

Unlike external Hadamard products, the following notation
gives width-one parcels by higher-width ones.

\begin{definition}\label{def:diag}
  For  $l\in \Zgeo$, suppose
  $\cF=\{\cF_{m} \in \Q(\fX)\}_{m\in \Zl}$.  Then, let us
  write the diagonal  $D(\cF)$ as the family
  $\{D(\cF)_{\tp{i}}=\cF_{\iota^{l}(i)} \in \Q(\fX)\}_{\tp{i}\in
    \Zo}$.
\end{definition}

\begin{proposition}\label{prop:diag-merg}
  Suppose a parcel
  $\cF=\Lam(s,\llo,\scc,w_{1},\fs,\phio,\rhoo, x_{1},\fX)$
  with flat $x_{1}=\iota^{\llo}(q)$, $\phio(x_{1})$, and
  $\rhoo$.  For $\llt=1$, let $w_{2}=\tp{\sum w_{1}}$ and
  $\rhot=\tp{\rho_{1,1}}$ in $\Zlt$. Also, let $x_{2}=\tp{q}$ and
  $\phit(x_{2})=\tp{\phio(x_{1})_{1}}$ in $\Q(\fX)^{\llt}$.
  Then, we have the following.
  \begin{enumerate}
  \item\label{c:diag-merg-parcel}
    There exists the width-one parcel 
    $D(\cF)=\Lam(s,\llt,\scc,w_{2},
    \gs,\phit,\rhot,x_{2},\fXt)$
    of $\gs=D(\fs)$.
  \item \label{c:diag-merg-scc}$D(\cF)$ is
    $\scc$-merged-log-concave, if $\cF$ is
    $\scc$-merged-log-concave.
  \item \label{c:diag-merg-sce}$D(\cF)$ is
    $\sce$-merged-log-concave, if $\cF$ is
    $\sce$-merged-log-concave.
  \end{enumerate}
\end{proposition}
\begin{proof}
  Suppose $m_{i}=\iota^{\llo}(i)$ for $i\in \Z$.  Let us prove
  Claim~\ref{c:diag-merg-parcel}.  Thus, we confirm the
  existence of a parcel
  $\cG= \Lam(s,\llt,\scc,w_{2}, \gs,\phit,\rhot,x_{2},\fX)$.
  First, $\gs$ is pairwise $\tp{s,\llt,\scc}$-positive,
  since
  $g_{s,\tp{i}}g_{s,\tp{j}}=f_{s,m_{i}}f_{s,m_{j}}\scc 0$ if
  $\tp{i},\tp{j}\ldZlt s$ and
  $g_{s,\tp{i}}g_{s,\tp{j}}=f_{s,m_{i}}f_{s,m_{j}}=0$
  otherwise.  Second, $q$ is $\scc$-admissible, since $\cF$
  is a parcel.  Third, $\phit$ is a
  $\tp{s,\llt,w_{2},\scc,\rhot,x_{2},\fX}$-mediator by the
  following. Since $\phio(x_{1})$ is flat, we have the base
  positivity:
  \begin{align}
    \phit(x_{2})^{w_{2,1}}_{1}
    =\phio(x_{1})^{\sum w_{1}}_{1}
    =\prod_{\lam\in\oi(\llo)}
    \phio(x_{1})^{ w_{1,\lam}}_{\lam}
    \gAfX 0.
  \end{align}
  Also, since $\rhoo$ and $x_{1}$ are flat, we have the
  base-shift positivity: each $\tp{i}\ldZlt s$ satisfies
  \begin{align}
    B(s,\llt,w_{2},\tp{i},\phit,\rhot,x_{2},\fX)
    &=b(s,\sum w_{1},i,\phi_{1,1},\rho_{1,1},q,\fX)\\
    &= \frac{ \phi_{1,1}(q^{\rho_{1,1}})^{i\sum w_{1}}
      \cdot [i]!_{q^{\rho_{1,1}}}^{\sum w_{1}}}
      { \phi_{1,1}(q)^{i\sum w_{1}} \cdot     [i]!_{q}^{\sum w_{1}}}\\
    &=\prod_{\lam\in\oi(\llo)}
      \frac{ \phi_{1,\lam}(q^{\rho_{1,\lam}})^{i w_{1,\lam}}
      \cdot [i]!_{q^{\rho_{1,\lam}}}^{ w_{1,\lam}}}
      { \phi_{1,\lam}(q)^{i w_{1,\lam}} \cdot     [i]!_{q}^{ w_{1,\lam}}}\\
    &=\prod_{\lam\in\oi(\llo)}b(s, w_{1,\lam},i,\phi_{1,\lam},
      \rho_{1,\lam},x_{1,\lam},\fX)\\
    &=B(s,\llo,w_{1},m_{i},\phio,\rhoo,x_{1},\fX)\\
    &\scc 0.
  \end{align}
  Hence, we have a parcel
  $\cG= \Lam(s,\llt,\scc,w_{2}, \gs,\phit,\rhot,x_{2},\fX)$.
  
  Claim~\ref{c:diag-merg-parcel} follows, because
  $\tp{i}\ldZlt s$ implies
  \begin{align}
    D(\cF)_{\tp{i}}
    &=\cF_{m_{i}}\\
    &=\frac{f_{s,m_{i}}}{\prod \phio(x_{1})^{m_{i}\rc w_{1}}\cdot
      [m_{i}]!_{x_{1}}^{w_{1}}}\\
    &=\frac{g_{s,\tp{i}}}{ \prod_{\lam\in\oi(\llo)}
      \phi_{1,1}(x_{1,1})^{i w_{1,\lam}}\cdot
      [i]!_{x_{1,1}}^{w_{1,\lam}}}\\
    &=\frac{g_{s,\tp{i}}}{
      \phi_{1,1}(x_{1,1})^{i \sum w_{1}}\cdot
      [i]!_{x_{1,1}}^{\sum w_{1}}}\\
    &=\frac{g_{s,\tp{i}}}{ \prod
      \phit(x_{2})^{ \tp{i}\rc  w_{2}}\cdot
      [\tp{i}]!_{x_{2}}^{w_{2}}}\\
    &=\cG_{\tp{i}}.
  \end{align}
   
  Let us prove Claim~\ref{c:diag-merg-scc}.  Suppose a
  fitting $\mut=\tp{s,\llt,\tp{i},\tp{j},\kt}$.  Let
  $a_{2}=\nu(\kt)$ and $b_{2}=\nu(\tp{i},\tp{j},\kt)$.  Also,
  suppose $\kn\in \Zlo$ such that $k_{1,1}=k_{2,1}$,
  $k_{1,\llo+1}=k_{2,2}$, and $k_{1,\lam}=0$ otherwise.  Let
  $a_{1}=\nu(\kn)$ and $b_{1}=\nu(m_{i},m_{j},\kn)$.
  
  First, we prove that $\muo=\tp{s,\llo,m_{i},n_{j},\kn}$ is
  fitting.  Since $\mut$ is fitting,
  Lemma~\ref{lem:fitting-one} implies
  \begin{align}
    \tp{i},\tp{j}&\ldZlt s,\label{ineq:diag-1}\\
    k_{2,1}\geq 0, k_{2,2}&\geq 1,\label{ineq:diag-2}\\
    j+k_{2,2}&>i. \label{ineq:diag-3}
  \end{align}
  Condition~\eqref{ineq:diag-1} gives $m_{i},m_{j}\ldZlo
  s$. Next, condition~\eqref{ineq:diag-2} yields
  $k_{2,1}= a_{1,1}= \dots = a_{1,\llo}<a_{1,\llo+1}=
  k_{2,1}+k_{2,2}= \dots =a_{1,2\llo}$.  Also,
  condition~\eqref{ineq:diag-3} implies
  $i+k_{2,1}=b_{1,1}=\dots =b_{1,\llo} <j+k_{2,1}+k_{2,2}=
  b_{1,\llo+1}=\dots =b_{1,2\llo}$.  Hence, $\muo$ is
  fitting.
  
  Second, let us prove the following equations:
  \begin{align}
    D(\cF)_{\tp{i}\bom \kt}&=\cF_{m_{i}\bom \kn};
                          \label{eq:diag-01}\\
    D(\cF)_{\tp{j}\bop \kt}&= \cF_{m_{j}\bop \kn};
                          \label{eq:diag-02}\\
    [a_{2}]!_{(x_{2}^{\rhot})^{\wcn}}^{w_{2}^{\wcn}}
                           &=
                             [a_{1}]!_{
                             (x_{1}^{\rhoo})^{\wcn}}^{w_{1}^{\wcn}};
                             \label{eq:diag-03}\\
    [b_{2}]!_{(x_{2}^{\rhot})^{\wcn}}^{w_{2}^{\wcn}}
                           &=
                             [b_{1}]!_{
                             (x_{1}^{\rhoo})^{\wcn}}^{w_{1}^{\wcn}};
                             \label{eq:diag-04}\\
    \prod (\phit(x_{2}^{\rhot})^{\wcn})^{(b_{2}-a_{2})
    \rc w_{2}^{\wcn}} 
                           &=
                             \prod (\phio(x_{1}^{\rhoo})^{\wcn})^{(b_{1}-a_{1})
                             \rc w_{1}^{\wcn}}.
                             \label{eq:diag-05}\\
  \end{align}
  Equations~\eqref{eq:diag-01} and~\eqref{eq:diag-02} hold,
  because $\tp{i}\bom \kt=\tp{i-k_{2,2}}$ and
  $m_{i}\bom \kn=\iota^{\llo}(i-k_{2,2})$, and
  $\tp{j}\bop \kt=\tp{j+k_{2,2}}$ and
  $m_{j}\bop \kn=\iota^{\llo}(i+k_{2,2})$.  We obtain
  equation~\eqref{eq:diag-03}, because
  $a_{2}=\tp{k_{2,1},k_{2,1}+k_{2,2}}$ and
  $w_{2}=\tp{\sum w_{1}}$ yield
  \begin{align}
    [a_{1}]!_{(x_{1}^{\rhoo})^{\wcn}}^{w_{1}^{\wcn}}
    =
    \prod_{\lam\in\oi(\llo)}
    [k_{2,1}]!_{q^{\rho_{1,1}}}^{w_{1,\lam}}
    \cdot
    \prod_{\lam\in\oi(\llo)}
    [k_{2,1}+k_{2,2}]!_{q^{\rho_{1,1}}}^{w_{1,\llo-\lam+1}}
    =[a_{2}]!_{(x_{2}^{\rhot})^{\wcn}}^{w_{2}^{\wcn}}.
  \end{align}
  Similarly, we obtain equation~\eqref{eq:diag-04} by
  $b_{2}=\tp{i+k_{2,1},j+k_{2,1}+k_{2,2}}$.  
  Equation~\eqref{eq:diag-05} follows, since
  \begin{align}
    \prod (\phio(x_{1}^{\rhoo})^{\wcn})^{(b_{1}-a_{1})
    \rc w_{1}^{\wcn}}
    &= \prod_{\lam\in\oi(\llo)}
      \phi_{1,1}(q^{\rho_{1,1}})^{i w_{1,\lam}}
      \cdot
      \prod_{\lam\in\oi(\llo)}
      \phi_{1,1}(q^{\rho_{1,1}})^{j w_{1,\llo-\lam+1}}\\
    &=\phi_{2,1}(q^{\rho_{2,1}})^{i \sum w_{1}}
      \cdot \phi_{2,1}(q^{\rho_{2,1}})^{j \sum w_{1}}\\
    &=
      \prod (\phit(x_{2}^{\rhot})^{\wcn})^{(b_{2}-a_{2})
      \rc w_{2}^{\wcn}}.
  \end{align}
  
  Third, equations~\eqref{eq:diag-01} and~\eqref{eq:diag-02}
  give
  $D(\cF)_{\tp{i}} D(\cF)_{\tp{j}^{\ve}}- D(\cF)_{\tp{i}\bom \kt}
  D(\cF)_{(\tp{j}\bop \kt)^{\ve}} = \cF_{m_{i}}
  \cF_{m_{j}^{\ve}} - \cF_{m_{i}\bom \kn} \cF_{(m_{j}\bop
    \kn)^{\ve}}$.  Moreover, equations~\eqref{eq:diag-03},
  ~\eqref{eq:diag-04}, and~\eqref{eq:diag-05} yield
  $\Ups(s,\llt,w_{2},\tp{i},\tp{j},\kt,\phit,\rhot,x_{2},\fX) =
  \Ups(s,\llo,w_{1},m_{i},m_{j},\kn,\phio,\rhoo,x_{1},\fX)$.
  Therefore, Claim~\ref{c:diag-merg-scc} holds, because the
  $\scc$-merged-log-concavity of $\cF$ implies
  \begin{align}
    \Delta(D(\cF))(s,\llt,w_{2},\tp{i},\tp{j},\kt,\phit,
    \rhot,x_{2},\fX) =
    \Delta(\cF)(s,\llo,w_{1},m_{i},m_{j},\kn,\phio,
    \rhoo,x_{1},\fX) \scc 0.
  \end{align}
  
  Claim~\ref{c:diag-merg-sce} holds similarly.
\end{proof}

\subsection{Internal Hadamard
  products}\label{sec:internal}
Internal Hadamard products of parcels
 give higher-weight merged-log-concave parcels
 from lower-weight ones.
 In particular, we obtain
higher-weight strictly merged-log-concave parcels
from weight-zero non-strictly ones.
To define the products, 
we state the following.
\begin{proposition}\label{prop:int-hadam}
  Consider  $\llo,\llt\in \Zgeo$
  with $\lam\in \tei(2,\llo)$
  and $\llt=\lamt-\lamo+1$.  Let
  $\cF_{i}
  =\Lam(s,l_{i},w_{i},\scc_{i},f_{i,s},\phi_{i},
  \rho_{i},x_{i},\fX_{i})$ for $i\in\oi(2)$ such
  that
  \begin{align}
    x_{2}
    &=x_{1}(\lamo,\lamt),\\
    \phit(x_{2})
    &=\phio(x_{1})(\lamo,\lamt),\\
    \rhot
    &=\rhoo(\lamo,\lamt).
  \end{align}
  Let $O_{i} =\{\sce_{i},\scc_{i}\}$ on $\fX_{i}$ for
  $i\in\oi(2)$.  Suppose squaring orders
  $O_{3}=\{\scer,\sccr\}$ on $\fXr=\fXo\cup \fXt$ such that
  $O_{3} \Sup O_{1}, O_{2}$.  Also, let
  \begin{align}
    w_{3}
    &=w_{1}+_{(\lamo,\lamt)}w_{2}\in \Zgez^{\llo},\\
    f_{3,s}
    &=\{f_{3,s,m}=f_{1,s,m}
      f_{2,s,m(\lamo,\lamt)}\in \Q(\fXr)\}_{m\in \Zlo}.
  \end{align}
  Then, there is a parcel
  \begin{align}
    \cFr
    &=\Lam(s,\llo,w_{3},\sccr,f_{3,s},\phio,\rhoo,x_{1},\fXr).
  \end{align}
\end{proposition}
\begin{proof}
  Let us prove that $x_{1}$ is $\sccr$-admissible.  Since
  $\cFo$ is a parcel, $x_{1}$ is $\scco$-admissible.  This
  gives the half $>_{x_{1,i}}$-$\sccr$ implication for each
  $x_{1,i}$ by $O_{3}\Sup O_{1}$.  Therefore, $x_{1}$ is
  $\sccr$-admissible by
  Lemma~\ref{lem:succ-on-coord-unions}, because we have the
  upper condition of each $x_{1,i}$ on $\fXo$.
  
  Let us confirm that $f_{3,s}$ is pairwise
  $\tp{s,\llo,\sccr}$-positive.  For each $m,n\ldZlo s$, 
    $f_{3,s,m}f_{3,s,n}= f_{1,s,m} f_{2,s,m(\lamo,\lamt)}
  f_{1,s,n} f_{2,s,n(\lamo,\lamt)} \sccr 0$, since
  $f_{1,s,m} f_{1,s,n} \scco 0$ and
  $f_{2,s,m(\lamo,\lamt)}f_{2,s,n(\lamo,\lamt)}\scct 0$ with
  $O_{3} \Sup O_{1}, O_{2}$.  If $m\nld_{\Zlo}s$, then
  $f_{1,s,m}=0$ gives
  $f_{3,s,m}=f_{1,s,m} f_{2,s,m(\lamo,\lamt)}= 0$.
  
  Let us prove that $\phio$ is a $\mu$-mediator for
  $\mu=(s,\llo,w_{3},\sccr,\rhoo,x_{1},\fXr)$.  Now,
  $\phio(x_{1})$ and $\phit(x_{2})$ are 
  mediators of $\cFo$
  and $\cFt$. Hence, we have 
  $\phio(x_{1})_{i}^{w_{1,i}}>_{\AfXo}0$
  for each $i\in\oi(\llo)$, and
  $\phio(x_{1})_{i}^{w_{2,i-\lamo+1}}\gAfXt0$ for each
  $i\in\oi(\lamo,\lamt)$.  Hence, $\phio$ and $\mu$ have the
  base positivity, since
  Lemma~\ref{lem:succ-on-coord-unions} implies
  \begin{align}
    \phio(x_{1})_{i}^{w_{3,i}}
    =
    \begin{dcases}
      \phio(x_{1})_{i}^{w_{1,i}}
      \phio(x_{1})_{i}^{w_{2,i-\lamo+1}}
      >_{A_{\fXr}}0
      \mforeach i\in\oi(\lamo,\lamt),\\
      \phio(x_{1})_{i}^{w_{1,i}}
      >_{A_{\fXr}}0 \melse.
    \end{dcases}
  \end{align}
  Also, $\phio$ and $\mu$ have the base-shift
  positivity, because each $m\ldZlo s$ satisfies
  \begin{dmath*}
    B(s,\llo,w_{3},m,\phio,\rhoo,x_{1},\fXr)
    =
    \prod_{i\in\oi(\llo)}\frac{ \phi(x_{i}^{\rho_{1,i}})^{w_{3,i}m_{i}}
      [m_{i}]!_{x_{i}^{\rho_{1,i}}}^{w_{3,i}}}{ 
      \phi(x_{i})^{w_{3,i}m_{i}}
      [m_{i}]!_{x_{i}}^{w_{3,i}}}
    =
    B(s,\llo,w_{1},m,\phio,\rhoo,x_{1},\fXo)
    \cdot B(s,\llt,w_{2},m(\lamo,\lamt),
    \phit,\rhot,x_{2},\fXt)
    \sccr 0.
  \end{dmath*}
  Hence,  $\phio$ is a $\mu$-mediator. This confirms
  the assertion.
\end{proof}

 Then, we have the following parcels.

\begin{definition}\label{def:int-hadam}
  Under the assumption in Proposition~\ref{prop:int-hadam},
  we define the internal Hadamard product
  $\cFo\dd_{\lam} \cFt$ as the parcel
  $\cFr=
  \Lam(s,\llo,w_{3},\sccr,f_{3,s},\phio,\rhoo,x_{1},\fXr)$.
\end{definition}

In particular, the diagonal reduction gives
\begin{align}
  D(\cF \ehd \cG)= D(\cF \dd_{\lam} \cG),  
\end{align}
whenever both sides are defined.

The existence of $\cF\dd_{\lam}\cG$ does not
imply that of $\cG\dd_{\lam}\cF$. But we have the following
commutativity.

\begin{proposition}\label{prop:int-hadam-comm}
  Let
  $\cF_{i}
  =\Lam(s,l_{i},w_{1},\scc_{i},f_{i,s},\phi_{i},\rho_{i},x_{i},\fX_{i})$
  for $i\in\oi(2)$.  Assume the parcel
  $\cFr=\cFo\dd_{\lam}\cFt$.
  Then, the parcel
  $\cFr'=\cFt\dd_{\lam}\cFo$
  exists if and only if $\llo=\llt$.  Furthermore, if
  $\llo=\llt$, then we have $\cFr=\cFr'$.
\end{proposition}
\begin{proof}
  Since there is the parcel $\cFr$,
  Proposition~\ref{prop:int-hadam} gives $\llo\geq \llt$.
  Also, we need $\llt\geq \llo$ to define the parcel
  $\cFr'$.  On the other hand, if $\llo=\llt$, then
  $\lam=(1,\llo)$.  This gives the parcel $\cFr'=\cFr$ by
  Proposition~\ref{prop:int-hadam}.
\end{proof}

If there is $\scc$-admissible $x\in \Q(\fX)^{l}$, then
$1\scc 0$ by Claim~\ref{c:adm-succ-half-gx-scc} of
Lemma~\ref{lem:adm-succ}.  Hence, we define the following
multifold products.
\begin{definition}\label{def:self-int-hadam}
  Suppose $\cF=\Lam(s,l,w,\scc,\fs,\phi,\rho,x,\fX)$ and
  $d\in \Zgez$.  We define the $d$-fold internal Hadamard
  product $\cF^{\dd d}$ such that
  \begin{align}
    \cF^{\dd d}
    =
    \begin{dcases}
      \Lam(s,l,d w,\scc,\fs^{d},\phi,\rho,x,\fX) \mif d\geq 1,\\
      \Lam(s,l,\scc,1_{s,l},\fX) \mif d=0.
    \end{dcases}
  \end{align}
\end{definition}

As we are interested in explicit real values of parcels,
we put the following remark on their coordinates.

\begin{remark}\label{rmk:int-hadam-non-optimal}
  Suppose
  $\cF_{i}=\Lam(s,l_{i},w_{i},\scc_{i},f_{i,s},\phi_{i},
  x_{i},\fX_{i})$ for $i\in\oi(3)$ such that
  $\cFr =\cFo\dd_{\lam}\cFt$.  Then, $\fXr$ does not have to
  be optimal for $\cFr$, even when $\fXo$ and $\fXt$ are
  optimal for $\cFo$ and $\cFt$, respectively.
  For example, let $s=\tp{0,\infty}$ and $\fXo=\{\qq\}$.
  Consider the monomial parcel
  $\cFo=\cFt=
    \Lam(s,1,\tp{1},\scc,\Psi_{s,\tp{\tp{0,\fraa,0}},q},x,\fX)$.
    Then, $\fX$ is optimal for 
     $\cFo$ and $\cFt$,
    but not for
    $\cFr=\cFo^{\dd 2}
          =\Lam(s,1,\tp{2},\scc,\Psi_{s,\tp{\tp{0,1,0}},q},x,\fX)$.
  \end{remark}

  We have the following merged-log-concavity of internal
Hadamard products.
\begin{theorem}\label{thm:merged-int-hadam}
  Let
  $\cF_{i} =\Lam(s,l_{i},w_{i},\scco,f_{i,s},
  \phi_{i},\rho_{i},x_{i},\fX_{i})$ for $i\in\oi(3)$ such that
  $\cFr=\cFo \dd_{\lam} \cFt$.  Consider squaring
  orders $O_{i}'=\{\sce_{i}',\scc'_{i}\}$ on $\fX_{i}$ for
  $i\in\oi(3)$ such that $O'_{3} \Sup O'_{1},O'_{2}$.  Let
  $\cFt$ be $\scet'$-merged-log-concave.  Then, we have the
  following.
  \begin{enumerate}
  \item \label{c:merged-int-hadam-scc} $\cFr$ is
    $\sccr'$-merged-log-concave, if $\cFo$ is
    $\scco'$-merged-log-concave.
  \item \label{c:merged-int-hadam-non-sce} $\cFr$ is
    $\scer'$-merged-log-concave, if $\cFo$ is
    $\sceo'$-merged-log-concave.
  \end{enumerate}
\end{theorem}
\begin{proof}
  Let us prove Claim~\ref{c:merged-int-hadam-scc}.  Assume a
  fitting $\muo=\tp{s,\llo,\mn,\nn,\kn}$.  Then,
  Claim~\ref{c:fitting-trunc-fitting} of
  Proposition~\ref{prop:fitting-trunc} gives the fitting
  $\mut=\tp{s,\llt,\mt,\nt,\kt}=t(\lam, \muo)$.
  Also, we put
   $\mur=\tp{s,\llr,\mr,\nr,\kr}=\muo$.  Let
  $a_{i}=\nu(k_{i})$, $b_{i}=\nu(m_{i},n_{i},k_{i})$, and
  $y_{i}=x_{i}^{\rho_{i}}$ for $i\in\oi(3)$.  Moreover,
  consider $L_{i},R_{i}\in \Q(\fX_{i})$ for $i\in\oi(3)$ such
  that
  \begin{align}
    L_{i}&=\Delta_{L}(\cF_{i})(s,l_{i},w_{i},m_{i},n_{i},
           k_{i},\phi_{i},\rho_{i},x_{i},\fX_{i})\\
         &=f_{i,s,m_{i}} f_{i,s,n_{i}^{\ve}}
           B(s,l_{i},w_{i},m_{i},n_{i}^{\ve},\phi_{i},\rho_{i},x_{i},\fX_{i})
           { b_{i} \brack a_{i}}^{w_{i}^{\wcn}}_{y_{i}^{\wcn}},\\
    R_{i}
         &=\Delta_{R}(\cF_{i})(s,l_{i},w_{i},m_{i},n_{i},
           k_{i},\phi_{i},\rho_{i},x_{i},\fX_{i})\\
         &=    f_{i,s,m_{i}\bom k_{i}} f_{i,s,(n_{i}\bop k_{i})^{\ve}}
           B(s,l_{i},w_{i},m_{i}\bom k_{i},
           (n_{i}\bop k_{i})^{\ve},\phi_{i},\rho_{i},x_{i},\fX_{i})
           {b_{i} \brack
           a_{i}^{\ve}}^{w_{i}^{\wcn}}_{y_{i}^{\wcn}}.
  \end{align}
  Then, let us prove
  \begin{align}
    L_{1}L_{2}&=L_{3},\label{eq:merged-int-hadam-L}\\
    R_{1}R_{2}&=R_{3}.\label{eq:merged-int-hadam-R}
  \end{align}
  
  First, let us confirm
  \begin{dmath}
    B(s,\llr,w_{3},\mr,\nr^{\ve},
    \phir,\rhor,x_{3},\fXr)
    =\prod_{i\in \oi(2)}
    B(s,l_{i},w_{i},m_{i},n_{i}^{\ve},
    \phi_{i},\rho_{i},x_{i},\fX_{i}).   
    \label{eq:merged-int-hadam-base-shift-mn}
  \end{dmath}
  Each $i\in\oi(\lamo,\lamt)$ gives
  \begin{dmath*}
    b(s,w_{3,i},m_{3,i},\phi_{3,i},\rho_{3,i},x_{3,i},\fXr)
    =b(s,w_{3,i},m_{1,i},\phi_{1,i},\rho_{1,i},x_{1,i},\fXr)
    =\frac{\phio(y_{1,i})^{m_{1,i}w_{3,i}}
      [m_{1,i}]_{y_{1,i}}^{w_{3,i}}}{
      \phio(x_{1,i})^{m_{1,i}w_{3,i}}[m_{1,i}]_{x_{1,i}}^{w_{3,i}}}
    =\frac{\phio(y_{1,i})^{m_{1,i}w_{1,i}}
      [m_{1,i}]_{y_{1,i}}^{w_{1,i}}}{
      \phio(x_{1,i})^{m_{1,i}w_{1,i}}
      [m_{1,i}]_{x_{1,i}}^{w_{1,i}}}
    \cdot
    \frac{\phit(y_{2,i-\lamo+1})^{m_{2,i-\lamo+1}
        w_{2,i-\lamo+1}}
      [m_{2,i-\lamo+1}
      ]_{y_{2,i-\lamo+1}}^{w_{2,i-\lamo+1}}}{
      \phit(x_{2,i-\lamo+1}
      )^{m_{2,i-\lamo+1}w_{2,i-\lamo+1}}
      [m_{2,i-\lamo+1}
      ]_{x_{2,i-\lamo+1}}^{w_{2,i-\lamo+1}}}
    =
    b(s,w_{1,i},m_{1,i},\phi_{1,i},\rho_{1,i},x_{1,i},\fXo)
    \cdot
    b(s,w_{2,i-\lamo+1},m_{2,i-\lamo+1},
    \phi_{2,i-\lamo+1},\rho_{2,i-\lamo+1},x_{2,i-\lamo+1},\fXt).
  \end{dmath*}
  Similarly,
  since $\nt^{\ve}=\nn^{\ve}(\lamo,\lamt)$,
  each $i\in\oi(\lamo,\lamt)$ yields
  \begin{dmath*}
    b(s,w_{3,i},(n^{\ve}_{3})_{i},
    \phi_{3,i},\rho_{3,i},x_{3,i},\fXr)
    =b(s,w_{3,i},(n^{\ve}_{1})_{i},
    \phi_{1,i},\rho_{1,i},x_{1,i},\fXr)
    =
    b(s,w_{1,i},(n^{\ve}_{1})_{i},
    \phi_{1,i},\rho_{1,i},x_{1,i},\fXo)
    \cdot
    b(s,w_{2,i-\lamo+1},(n^{\ve}_{2})_{i-\lamo+1},
    \phi_{2,i-\lamo+1},\rho_{2,i-\lamo+1},x_{2,i-\lamo+1},\fXt).
  \end{dmath*}
  Hence, equation~\eqref{eq:merged-int-hadam-base-shift-mn}
  follows.
  
  Second, let us prove
  \begin{dmath}
    B(s,\llo,w_{3},\mr\bom \kr, (\nr\bop
    \kr)^{\ve},\phir,\rhor,x_{3},\fXr) = \prod_{i\in
      \oi(2)}B(s,l_{i},w_{i},m_{i}\bom k_{i}, (n_{i}\bop
    k_{i})^{\ve},\phi_{i},\rho_{i},x_{i},\fX_{i}).
    \label{eq:merged-int-hadam-base-shift-bombop}
  \end{dmath}
  Claim~\ref{c:fitting-trunc-sigma-equations} of
  Proposition~\ref{prop:fitting-trunc} gives
  $\mt\bom \kt= (\mn\bom \kn)(\lamo,\lamt)$ and
  $(\nt\bop \kt)^{\ve}= (\nn\bop \kn)^{\ve}(\lamo,\lamt)$.
  Thus,
  equation~\eqref{eq:merged-int-hadam-base-shift-bombop}
  follows, since each $i\in\oi(\lamo,\lamt)$ satisfies the
  following equations:
  \begin{dgroup*}
    \begin{dmath*}
      b(s,w_{3,i},(\mr\bom
      \kr)_{i},\phi_{3,i},\rho_{3,i},x_{3,i},\fXr)
      =b(s,w_{3,i},(\mn\bom
      \kn)_{i},\phi_{1,i},\rho_{1,i},x_{1,i},\fXr)
      =
      b(s,w_{1,i},(\mn\bom \kn)_{i},
      \phi_{1,i},\rho_{1,i},x_{1,i},\fXo)
      \cdot
      b(s,w_{2,i-\lamo+1},(\mt\bom \kt
      )_{i-\lamo+1},
      \phi_{2,i-\lamo+1},\rho_{2,i-\lamo+1},
      x_{2,i-\lamo+1},\fXt);
    \end{dmath*}
    \begin{dmath*}
      b(s,w_{3,i},(\nr\bop \kr)^{\ve}_{i},
      \phi_{3,i}, \rho_{3,i},x_{3,i},\fXr)
      =b(s,w_{3,i},(\nn\bop \kn)^{\ve}_{i},
      \phi_{1,i}, \rho_{1,i},x_{1,i},\fXr)
      =
      b(s,w_{1,i},(\nn\bop
      \kn)^{\ve}_{i},\phi_{1,i},\rho_{1,i},x_{1,i},\fXo)
      \cdot
      b(s,
      w_{2,i-\lamo+1},(\nt\bop \kt)^{\ve}_{i-\lamo+1},
      \phi_{2,i-\lamo+1},
      \rho_{2,i-\lamo+1},x_{2,i-\lamo+1},\fXt).
    \end{dmath*}
  \end{dgroup*}
  
  Third, since $w_{3}=w_{1}+_{\tp{\lamo,\lamt}}w_{2}$,
  Claim~\ref{c:fitting-trunc-ab-equations} of
  Proposition~\ref{prop:fitting-trunc} implies
  \begin{align}
    { b_{1} \brack a_{1}}^{w_{1}^{\wcn}}_{y_{1}^{\wcn}}
    { b_{2} \brack a_{2}}^{w_{2}^{\wcn}}_{y_{2}^{\wcn}}
    &=
      { b_{3} \brack
      a_{3}}^{w_{3}^{\wcn}}_{y_{3}^{\wcn}},
      \label{eq:merged-int-hadam-binom}
    \\
    { b_{1} \brack a^{\ve}_{1}}^{w_{1}^{\wcn}}_{y_{1}^{\wcn}}
    { b_{2} \brack a^{\ve}_{2}}^{w_{2}^{\wcn}}_{y_{2}^{\wcn}}
    &=
      { b_{3} \brack a^{\ve}_{3}
      }^{w_{3}^{\wcn}}_{y_{3}^{\wcn}}.
      \label{eq:merged-int-hadam-binom-wed}
  \end{align}
  
  Therefore,
  equations~\eqref{eq:merged-int-hadam-base-shift-mn}
  and~\eqref{eq:merged-int-hadam-binom} give
  equation~\eqref{eq:merged-int-hadam-L}.  Also,
  equations~\eqref{eq:merged-int-hadam-base-shift-bombop}
  and~\eqref{eq:merged-int-hadam-binom-wed} yield
  equation~\eqref{eq:merged-int-hadam-R}.

  Let us prove that $\cFr$ is $\sccr'$-merged-log-concave.
  We have $L_{1} \scco' R_{1}$ and $L_{2} \scet' R_{2}$ by
  the merged-log-concavity of $\cFo$ and $\cFt$.  Hence,
  because $O'_{3}\Sup O'_{1}, O'_{2}$,
  Claim~\ref{c:merged-binom-bshift-scc-eq} of
  Theorem~\ref{thm:merged-binom-bshift} gives
  \begin{align}
    L_{1} &\sccr' R_{1}\scer' 0 \mand L_{2} \scer' R_{2},\\
    L_{1}&\sccr'0 \mand   L_{2}\sccr'0,\\
    R_{2}&\sccr'0 \mor   R_{2}=0.
  \end{align}
  Then, equations~\eqref{eq:merged-int-hadam-L}
  and~\eqref{eq:merged-int-hadam-R} imply
  $L_{3}\sccr' R_{3}$ by Corollary~\ref{cor:pair-orders}.
  In particular, Claim~\ref{c:merged-int-hadam-scc} holds.
  
  Claim~\ref{c:merged-int-hadam-non-sce} follows from
  Claim~\ref{c:succ-succeq-zero-zero-1} of
  Lemma~\ref{lem:succ-succeq}, because $L_{1} \sceo' R_{1}$
  and $L_{2} \scet' R_{2}$.
\end{proof}

In particular, we obtain higher-weight strictly
merged-log-concave parcels from weight-zero non-strict ones
by
Theorems~\ref{thm:merged-ext-hadam}
and~\ref{thm:merged-int-hadam}.
 
\begin{corollary}\label{cor:ext-int-hadam}
  Let $\lam\in \Zgeo$.  For each $i\in\oi(\lam)$, let
  $\cF_{i}=\Lam(s,l_{i},\scc_{i},f_{i,s},\fX_{i})$ be
  $\sce_{i}$-merged-log-concave.  Let
  $\fX=\cup_{i\in\oi(\lam)} \fX_{i}$. Consider squaring
  orders $O=\{\sce,\scc\}$ on $\fX$ such that $O$ is
  compatible to each $\{\sce_{i},\scc_{i}\}$.  Let
  \begin{align}
    \cG= \Lam(s,l,\scc,\gs,\fX) = \cFo\ehd
    \cFt\ehd \dots \ehd \cF_{\lam}.
  \end{align}
    Then, we have the following.
  \begin{enumerate}
  \item
    \label{c:ext-int-hadam-more-deform}
    For each $w\in \Zlgez$, if
    $\Lam(s,l,w,\scc,k_{s},\phi,\rho,x,\fX)$ is a
    $\sce$-merged-log-concave parcel, then there is the
    $\scc$-merged-log-concave parcel
    $\Lam(s,l,w,\scc,k_{s} \gs,\phi,\rho,x,\fX)$.
  \item
    \label{c:ext-int-hadam-identity}
    For each $w\in \Zlgeo$, if $\phi$ is a proper
    $\tp{s,l,w,\scc,\rho,x,\fX}$-mediator, then there is a
    $\scc$-merged-log-concave parcel
    $\Lam(s,l,w,\scc,\gs,\phi,\rho,x,\fX)$.
  \end{enumerate}  
\end{corollary}
\begin{proof}
  Claim~\ref{c:ext-int-hadam-more-deform} follows from
  Claim~\ref{c:merged-int-hadam-scc} of
  Theorem~\ref{thm:merged-int-hadam}, because
  $\cG=\Lam(s,l,\tp{0},\scc,\gs,\phi,\rho,x,\fX)$ is
  $\sce$-merged-log-concave by
  Claim~\ref{c:merged-ext-hadam-sce} of
  Theorem~\ref{thm:merged-ext-hadam}.

  Let us prove Claim~\ref{c:ext-int-hadam-identity}.  Since
  $1\scc 0$ by the $\scc$-admissible $x$ and
  Claim~\ref{c:adm-succ-half-gx-scc} of
  Lemma~\ref{lem:adm-succ}, we have the constant parcel
  $\Lam(s,l,w,\scc,1_{s,l},\phi,\rho,x,\fX)$.  This is
  $\scc$-merged-log-concave by
  Proposition~\ref{prop:constant-parcel}.  Hence,
  Claim~\ref{c:ext-int-hadam-more-deform} gives
  Claim~\ref{c:ext-int-hadam-identity}.
\end{proof}

\begin{remark}\label{rmk:weight-difference}
  Assume an infinite gate $s\geq 0$ with $l=1$ and
  $w\in \Zgeo^{l}$.  In Claim~\ref{c:ext-int-hadam-identity}
  of Corollary~\ref{cor:ext-int-hadam},
  $\Lam(s,l,w,\scc,\fs,x,\fX)$ is $\scc$-merged-log-concave,
  if $\Lam(s,l,\scc,\fs,\fX)$ is $\sce$-merged-log-concave.
  However, the converse is not necessarily true.  For
  example,
  $\Lam(s,l,w,\llq, \Psi_{s,\tp{\tp{\fraa,0,0}},q},x, \fX)$ is
  $\llq$-merged-log-concave, but
  $\Lam(s,l,\llq,\Psi_{s,\tp{\tp{\fraa,0,0}},q},\fX)$ is not
  $\ggq$-merged-log-concave by
  equation~\eqref{eq:poly-not-pos-int-coeff}.  Thus, we
  study nontrivial-weight merged-log-concave parcels
  to obtain polynomials with positive integer coefficients
  in this manuscript.
\end{remark}

\section{Weight-zero
  merged-log-concavity, strong
  \texorpdfstring{$q$}{q}-log-concavity,
  and
  \texorpdfstring{$q$}{q}-log-concavity}
\label{sec:weight-zero-merged-q-log-conc}

We discuss the merged-log-concavity by the strong
$q$-log-concavity and $q$-log-concavity.  In particular, the
strong $q$-log-concavity of polynomials corresponds to the
weight-zero merged-log-concavity in a suitable
setting. Thus, Corollary~\ref{cor:ext-int-hadam} implies
that strongly $q$-log-concave polynomials give higher-weight
merged-log-concave rational functions.  Additionally, we
prove some analogs of conjectures in
Section~\ref{sec:explicit-merged-parcels} by $q$-numbers and
give conjectures on weight-zero parcels.

\subsection{Strong
  \texorpdfstring{$q$}{q}-log-concavity
  and merged-log-concavity}
\label{sec:st-q-log-c-merged}
Let us recall the strong $q$-log-concavity.

\begin{definition}(\cite[Section
  1]{Sag})\label{def:st-q-log-conc} A sequence
  $f=\{f_{m}\in \N[q]\}_{m\in \Z}$ is called strongly
  $q$-log-concave, if $n\geq m$ implies
  $f_{m}f_{n}-f_{m-1}f_{n+1}\ggq 0$.
\end{definition}

We then formulate the following notion
 by tuples and
squaring orders to compare
 the strong $q$-log-concavity with the merged-log-concavity.

\begin{definition}\label{def:succ-strong-log-conc}
  Let $l\in \Zgeo$ and $\cF=\Lam(s,l,\scc,\fs,\fX)$.  Suppose
  squaring orders $O'=\{\sce',\scc'\} \Sup \{\sce,\scc\}$.
  \begin{enumerate}
  \item Let us call $\cF$ strongly $\scc'$-log-concave, if
    $\cF_{m}\cF_{n^{\ve}}-\cF_{m-1}\cF_{(n+1)^{\ve}}\scc' 0$
    whenever $n,m\ldZl s$, $n^{\ve}\geq m$, and
    $m\ccn (n+1)\in \Ztl$ is increasing.
  \item Let us call $\cF$ strongly $\sce'$-log-concave, if
    $\cF_{m}\cF_{n^{\ve}}-\cF_{m-1}\cF_{(n+1)^{\ve}}\sce' 0$
    whenever $n,m\ldZl s$, $n^{\ve}\geq m$,  and
    $m\ccn (n+1)\in \Ztl$ is increasing.
  \end{enumerate}
\end{definition}

We use the following telescoping lemma.

\begin{lemma}\label{lem:telescoping}
  Let $l\in \Zgeo$ and $\cF=\Lam(s,l,\scc,\fs,\fX)$.
  Then, $\cF$ is strongly $\scc'$-log-concave if and only if
  \begin{align}
    \cF_{m}\cF_{n^{\ve}}-\cF_{m-k}\cF_{(n+k)^{\ve}}\scc' 0
    \label{ineq:telescoped}
  \end{align}
  whenever $m,n\ldZl s$, $n^{\ve}\geq m$, $k\in \Zgeo$, and
  $m\ccn (n+k)$ is increasing.
\end{lemma}
\begin{proof}
  The if part holds by $k=1$.  Thus, we prove the only if
  part, assuming $k=1$ and
  inequality~\eqref{ineq:telescoped}. Consider $m,n\ldZl s$
  such that $n^{\ve}\geq m$ and $m\ccn(n+1)$ is increasing.
  
  First, suppose $i\in \Zgez$ such that $m-i\geq s_{1}$ and
  $n+i\leq s_{2}$. Then, $m-i,n+i\ldZl s$, because
  $(n+i)^{\ve}=n^{\ve}+i\geq m-i$ by
  $n^{\ve}\geq m$.  Thus, since
  $(m-i)\ccn (n+i+1)$ is increasing,
  $\cF_{m-i}\cF_{(n+i)^{\ve}}-\cF_{m-i-1}\cF_{(n+i+1)^{\ve}}
  \scc' 0$.
  
  Second, suppose $i\in \Zgez$ such that either
  $m-i< s_{1} $ or $n+i> s_{2}$.  If $m-i<s_{1}$, then
  $m-i-1<s_{1}$.  Also, if $n+i>s_{2}$, then
  $n+i+1>s_{2}$. Hence,
  $\cF_{m-i}\cF_{(n+i)^{\ve}}-\cF_{m-i-1}\cF_{(n+i+1)^{\ve}}=0$.
  In particular, each $k\in \Zgeo$ yields
  \begin{align}
    \cF_{m}\cF_{n^{\ve}}-\cF_{m-k}\cF_{(n+k)^{\ve}}
    &=
      \sum_{i\in\oi(0,k-1)}
      (\cF_{m-i}\cF_{(n+i)^{\ve}}-\cF_{m-i-1}\cF_{(n+i+1)^{\ve}})\\
    &\sce'
      \cF_{m}\cF_{n^{\ve}}-\cF_{m-1}\cF_{(n+1)^{\ve}}\\
    &\scc' 0.
  \end{align}
\end{proof}
We also use the following special fitting tuples.
\begin{lemma}\label{lem:plus-one-fitting}
  For some $l\in \Zgeo$,
  suppose
  $n,m\ldZl s$ such that $n^{\ve}\geq m$ and
  $m\ccn (n+1)\in \Ztl$ is increasing.  
  Then, $\mu=\tp{s,l,m,n,k}$ is fitting
  for $k=\iota^{l}(0)\ccn \tp{1}\ccn \iota^{l-1}(0)$.
\end{lemma}
\begin{proof}
  Let $a=\nu(k)$
  and $b=\nu(m,n,k)$.  Then,
  $a=\iota^{l}(0)\ccn \iota^{l}(1)$ and $b=m\ccn (n+1)$.
  Hence, $\mu$ satisfies the slope conditions,
  since $n^{\ve}\geq m$ implies
  $b_{l}=m_{l}<\nn+1=b_{l+1}$. 
\end{proof}

Then, we state the  comparison below.
\begin{proposition}\label{prop:lc-lc}
  Let $l\in \Zgeo$ and $\cF=\Lam(s,l,\scc,\fs,\fX)$.
  Then, we have the following.
  \begin{enumerate}
  \item \label{c:lc-lc-scc-non-equiv} $\cF$ is strongly
    $\scc'$-log-concave, if $\cF$ is
    $\scc'$-merged-log-concave.
  \item \label{c:lc-lc-sce-non-equiv} $\cF$ is strongly
    $\sce'$-log-concave, if $\cF$ is
    $\sce'$-merged-log-concave.
  \end{enumerate}
  Moreover, assume $l=1$. Then, we have the following.
  \begin{enumerate}[label=(\alph*)]
  \item \label{c:lc-lc-scc-equiv}$\cF$ is
    $\scc'$-merged-log-concave if and only if $\cF$ is
    strongly $\scc'$-log-concave.
  \item \label{c:lc-lc-sce-equiv} $\cF$ is
    $\sce'$-merged-log-concave if and only if $\cF$ is
    strongly $\sce'$-log-concave.
  \end{enumerate}
\end{proposition}
\begin{proof}
  Let us prove Claim~\ref{c:lc-lc-scc-non-equiv}.  Consider
  $n,m\ldZl s$ such that $n^{\ve}\geq m$ and
  $m\ccn (n+1)\in \Ztl$ is increasing.  Also, let
  $\mu=\tp{s,l,m,n,k}$ for
  $k=\iota^{l}(0)\ccn \tp{1}\ccn \iota^{l-1}(0)$.  Then,
  we have
  $\Delta(\cF)(s,l,m,n,k,\fX)\scc' 0$ by
  Lemma~\ref{lem:plus-one-fitting}.  
  Claim~\ref{c:lc-lc-scc-non-equiv} now follows, since
  $m\bom k=m-1$ and $n\bop k=n+1$ imply
  $\Delta(\cF)(s,l,m,n,k,\fX)= f_{s,m}f_{s,n^{\ve}}- f_{s,m-1}
  f_{s,(n+1)^{\ve}}$.  Claim~\ref{c:lc-lc-sce-non-equiv}
  holds similarly.
  
  Let us prove Claim~\ref{c:lc-lc-scc-equiv}.  First, we
  prove the only if part.  By Lemma~\ref{lem:fitting-one},
  $\cF$ is $\scc'$-merged-log-concave if and only if
  \begin{align}
    \Delta(\cF)(s,l,m,n,k,\fX)=
    f_{s,m}f_{s,n}-
    f_{s,m- \kt} f_{s,n+\kt}
    \scc' 0
    \label{eq:lc-lc-merged-simple}
  \end{align}
  whenever
  \begin{align}
    m,n\ldZl s, \ k=\tp{\kn,\kt}\geq \tp{0,1}, \mand
    n+\kt> m.
    \label{cond:lc-lc-fitting-simple}
  \end{align}
  Thus, 
  inequality~\eqref{eq:lc-lc-merged-simple} gives the strong
  $\scc'$-log-concavity of $\cF$,
  since
conditions~\eqref{cond:lc-lc-fitting-simple} hold by
$\kt=1$ and $n,m\ldZ s$ such that $n\geq m$.

  Second, we prove the if part.  Consider $n,m,k$
  that satisfy conditions~\eqref{cond:lc-lc-fitting-simple}.
  If $n\geq m$, then
  inequality~\eqref{eq:lc-lc-merged-simple} holds by
  Lemma~\ref{lem:telescoping}.  If $n<m$, then let
  $m'=n$, $n'=m$, and $\kt'=\nn+\kt-\mn\in \Zgeo$.
  In particular, we obtain
  $n'>m'$, 
  $m'-\kt'=n-(n+\kt-m)=m-\kt$, and
  $n'+\kt'=m+(n+\kt-m)=n+\kt$.
  Then, Lemma~\ref{lem:telescoping} implies
  inequality~\eqref{eq:lc-lc-merged-simple}, since
  $f_{s,m}f_{s,n}- f_{s,m- \kt} f_{s,n+\kt}
  =f_{s,m'}f_{s,n'}- f_{s,m'- k'_{2}} f_{s,n'+k'_{2}}$.
  Hence, Claim~\ref{c:lc-lc-scc-equiv} holds.
  Claim~\ref{c:lc-lc-sce-equiv} follows similarly.
\end{proof}

In particular, we have the following explicit correspondence
between the strong $q$-log-concavity
 and the merged-log-concavity.

\begin{corollary}\label{cor:merged-qstr}
  Let $l=1$.  Assume a gate $s\geq 0$ and $\llq$-admissible
  $q\in \Q(\fX)$.
  \begin{enumerate}
  \item \label{c:merged-qstr-from-qstr} Consider a strongly
    $q$-log-concave
    $f_{1,s}=\{f_{1,s,m}\in \N[q]\}_{m \in \Z}$ such that
    $f_{1,s,m}\llq 0$ if $m\ldZ s$ and
    $f_{1,s,m}=0$ otherwise.
    \begin{enumerate}
    \item \label{c:merged-qstr-from-qstr-to-non-strict-merged}
      Then, there is the $\ggq$-merged-log-concave
      $\cF=\Lam(s,l,\llq,f_{2,s},\fX)$ such that
      $f_{2,s,m}=f_{1,s,\mn}$ for each
      $m=\tp{\mn}\in \Zl$.
    \item \label{c:merged-qstr-from-qstr-to-strict-merged}
      Moreover, if
      \begin{align}
        f_{1,s,m}f_{1,s,n}- f_{1,s,m-1}f_{1,s,n+1}\llq 0
        \label{ineq:merged-qstr-from-qstr-to-strict-merged-llq}
      \end{align}
      whenever $n,m\ldZ s$ and $n\geq m$, then $\cF$ is
      $\llq$-merged-log-concave.
    \end{enumerate}
  \item \label{c:merged-qstr-from-merged}Conversely, consider
    a $\ggq$-merged-log-concave 
    $\cF=\Lam(s,l,\llq,f_{2,s},\fX)$ 
    such that $f_{2,s}=\{f_{2,s,m}\in \N[q]\}_{m\in \Zl}$.
    \begin{enumerate}
    \item Then, there is the strongly $q$-log-concave
      $f_{1,s}=\{f_{1,s,m}=f_{2,s,\tp{m}}\}_{m \in
        \Z}$.
    \item Furthermore, if $\cF$ is
      $\llq$-merged-log-concave, then we have
      $f_{1,s,m}f_{1,s,n}- f_{1,s,m-1}f_{1,s,n+1}\llq 0$
      whenever  $n,m\ldZ s$ and $n\geq m$.
    \end{enumerate}
  \end{enumerate}
\end{corollary}
\begin{proof}
  Let us prove
  Claim~\ref{c:merged-qstr-from-qstr-to-non-strict-merged}.
  Let $n,m\in \Zl$ such that $n\geq m$.  Then,
  $f_{2,s,m}f_{2,s,n}-f_{2,s,m-1}f_{2,s,n+1}\ggq0$, since
  $f_{1,s}$ is strongly $q$-log-concave.  Hence, $\cF$ is
  strongly $\ggq$-log-concave, since $n\geq m$ if and only if
  $n^{\ve}\geq m$ and $m\ccn (n+1)$ is increasing.  Thus,
  $\cF$ is $\ggq$-merged-log-concave by
  Claim~\ref{c:lc-lc-sce-equiv} of
  Proposition~\ref{prop:lc-lc}.

  Since
  inequality~\eqref{ineq:merged-qstr-from-qstr-to-strict-merged-llq}
  implies that $\cF$ is strongly $\llq$-log-concave,
  Claim~\ref{c:merged-qstr-from-qstr-to-strict-merged}
  follows from
  Claim~\ref{c:lc-lc-scc-equiv} of
  Proposition~\ref{prop:lc-lc}.

  Claim~\ref{c:merged-qstr-from-merged} follows similarly
  from
  Claims~\ref{c:lc-lc-scc-equiv} and~\ref{c:lc-lc-sce-equiv}
  of
  Proposition~\ref{prop:lc-lc}.
\end{proof}

\subsection{\texorpdfstring{$q$}{q}-log-concavity
  and merged-log-concavity}
As 
in Section~\ref{sec:st-q-log-c-merged},
we formulate the following
notion by tuples and squaring orders
to compare the $q$-log-concavity in
Definition~\ref{def:log-conc} 
 with the merged-log-concavity.

\begin{definition}\label{def:succ-log-conc}
  Let $l\in \Zgeo$ and $\cF=\Lam(s,l,\scc,\fs,\fX)$.  Suppose
  squaring orders $\{\sce',\scc'\} \Sup \{\sce,\scc\}$.
  \begin{enumerate}
  \item Let us call $\cF$ $\scc'$-log-concave, if
    $\cF_{m}\cF_{m}-\cF_{m-1}\cF_{m+1}\scc' 0$ whenever
    $m\ldZl s$ and $m\ccn (m+1)\in \Ztl$ is increasing.
  \item Let us call $\cF$ $\sce'$-log-concave, if
    $\cF_{m}\cF_{m}-\cF_{m-1}\cF_{m+1}\sce' 0$
    whenever $m\ldZl s$ and $m\ccn (m+1)\in \Ztl$ is
    increasing.
  \end{enumerate}
\end{definition}
Then, we have the following comparison.

\begin{proposition}\label{prop:lc-lc-non-strong}
  Let $l\in \Zgeo$ and $\cF=\Lam(s,l,\scc,\fs,\fX)$.  Then,
  we have the following.
  \begin{enumerate}
  \item $\cF$ is $\scc'$-log-concave, if $\cF$ is
    $\scc'$-merged-log-concave.
  \item $\cF$ is $\sce'$-log-concave, if $\cF$ is
    $\sce'$-merged-log-concave.
  \end{enumerate}
\end{proposition}
\begin{proof}
  The assertions hold by
  Claims~\ref{c:lc-lc-scc-non-equiv} and~\ref{c:lc-lc-sce-non-equiv}
  of Proposition~\ref{prop:lc-lc},
  since the strong $\scc'$-log-concavity and the strong
  $\sce'$-log-concavity imply the $\scc'$-log-concavity
  and the $\sce'$-log-concavity, respectively.
\end{proof}
This gives the $q$-log-concavity
from the merged-log-concavity as follows.

\begin{corollary}\label{cor:merged-qstr-non-strong}
  Let $l=1$. Suppose that $\cF=\Lam(s,l,\llq,f_{1,s},\fX)$
  is $\ggq$-merged-log-concave and
  $f_{1,s}=\{f_{1,s,m}\in \N[q]\}_{m\in \Zl}$.  Then, we have
  the following.
  \begin{enumerate}
  \item The sequence
    $f_{2,s}=\{f_{2,s,m} =f_{1,s,\tp{m}}\in \N[q]\}_{m \in \Z}$
    is $q$-log-concave.
  \item If $\cF$ is $\llq$-merged-log-concave,
    then 
    $f_{2,s,m}^{2}- f_{2,s,m-1}f_{2,s,m+1}\llq 0$
    for each $m\ldZ s$.
  \end{enumerate}
\end{corollary} 
\begin{proof}
  The assertions follow from 
  Proposition~\ref{prop:lc-lc-non-strong}.
\end{proof}

\subsection{Some analogs of conjectures in
  Section~\ref{sec:explicit-merged-parcels}}
\label{sec:weight-zero-conjec}

We put the following weight-zero parcels.

\begin{definition} \label{def:chi} Assume a gate $s\geq 1$,
  $l\in \Zgeo$, and $q\in \Q(\fX)$.  Let
  \begin{align}
    \chi_{s,l,q,m}=
    \begin{dcases}
      [m]_{q}\in \Q(\fX) \mfor m\ldZl s,\\
      0 \melse.&
    \end{dcases}
  \end{align}
  Moreover, if $q$ is $\scc$-admissible,  then we call
  $\Lam(s,l,\scc,\chi_{s,l,q},\fX)$ $q$-number parcel.
\end{definition}
We prove the density and unimodality on
the merged determinants of width-one
$q$-number parcels. Then, this confirms some analogs of
Conjectures~\ref{conj:density}
and~\ref{conj:unimodal-width-two} for the density and
unimodality on the merged determinants of positive-weight
parcels.
Moreover, we put a conjecture on
higher-width $q$-number parcels in Section~\ref{sec:subsec-conj}.

Let us put the following integers.

\begin{definition}\label{def:min-qnum}
  Suppose $\mn,\nn,\lam\in \Z$. Then, let
  $I(\lam,\mn,\nn)
  =\min(\lam+1,\mn,\nn,\mn+\nn-1-\lam)\in \Z$.    
  \end{definition}
  We compute these
  integers $I(\lam,\mn,\nn)$
  for the products of
$q$-numbers as follows.
  \begin{lemma}\label{lem:chi-lem}
  Let $\mn,\nn\in \Zgeo$ and $\lam\in \bZ$.
  \begin{enumerate}
  \item \label{c:chi-lem-min} If $\mn\leq \nn$, then 
    \begin{numcases}{I(\lam,\mn,\nn)=}
      \lam+1  \mif \lam\in\oi(0,\mn-1),
      \label{eq:chi-lem-a}\\
      \mn \mif \lam\in\oi(\mn-1,\nn-1),
      \label{eq:chi-lem-b}\\
      \mn+\nn-1-\lam \mif \lam\in\oi(\nn-1,\mn+\nn-2).
      \label{eq:chi-lem-c}
    \end{numcases}
  \item \label{c:chi-lem-prod} We have
    $[\mn]_{q}[\nn]_{q}= \sum_{\lam\in\oi(0,\mn+\nn-2)}
    I(\lam,\mn,\nn) q^{\lam}$.
  \item
    \label{c:chi-lem-ineq}
    Assume $\nn+\kt>\mn$ for some
    $\kt\in \Zgeo$.
    \begin{enumerate}
    \item \label{c:chi-lem-ineq-nonst}
      Then,
      $I(\lam,\mn,\nn)
        \geq
        I(\lam,\mn-\kt,\nn+\kt)$.
    \item
      \label{c:chi-lem-ineq-st}
      Also, if $\mn-\kt\geq 1$,  then there is
      $\lam\in \Z$ such that
      \begin{align}
        \mn+\nn-2
        &\geq \nn-1 \geq \lam\geq \mn-\kt\geq 0,
          \label{ineq:chi-lem-ineq-st-param}\\
        I(\lam,\mn,\nn)
        &>
          I(\lam,\mn-\kt,\nn+\kt).
          \label{ineq:chi-lem-lem-ineq-st}
      \end{align}
    \end{enumerate}
  \end{enumerate}
\end{lemma}
\begin{proof}
  Let us prove Claim~\ref{c:chi-lem-min}.  First, if
  $\lam\in\oi(0,\mn-1)$, then equation~\eqref{eq:chi-lem-a}
  follows from
  $\lam+1\leq \mn\leq \nn\leq \nn+\mn-(\lam+1)= \mn+\nn-1-\lam$.
  Second, if $\lam\in\oi(\mn-1,\nn-1)$, then
  equation~\eqref{eq:chi-lem-b} holds by
  $\mn\leq \lam+1\leq \nn$ and $\mn+\nn-(1+\lam)\geq \mn$.  Third, if
  $\lam\in\oi(\nn-1,\mn+\nn-2)$, then
  equation~\eqref{eq:chi-lem-c} follows from
  $\lam+1\geq \nn\geq \mn \geq \mn+\nn-1-\lam$.
  
  Let us prove Claim~\ref{c:chi-lem-prod}.  We write the
  product as
  \begin{align}
    [\mn]_{q}[\nn]_{q}
    =\sum_{t_{1}\in\oi(0,\mn-1), t_{2}\in\oi(0,\nn-1)} q^{t_{1}+t_{2}}.
    \label{eq:chi-lem-prod-two-indices}
  \end{align}
  Assume $\nn\geq \mn$ without loss of generality.  Let
  $\lam\in \Z$.  Thus, we consider $\tp{t_{1},t_{2}}$ in
  equation~\eqref{eq:chi-lem-prod-two-indices} such that
  $\lam=t_{1}+t_{2}$.  First, let $\lam\in\oi(0,\mn-1)$.
  Then, to obtain $q^{\lam}$, we have $\lam+1$ choices of
  $\tp{t_{1},t_{2}}$ such that
  $\tp{\lam,0}, \tp{\lam-1,1},\dots, \tp{0,\lam}$, since
  $\lam\leq \nn-1$.  Second, if $\lam\in\oi(\mn-1,\nn-1)$, then
  there are $\mn$ choices
  $\tp{\mn-1,\lam-(\mn-1)}, \dots, \tp{0,\lam}$ by
  $\mn\geq 1$.  Third, if $\lam\in\oi(\nn-1,\mn+\nn-2)$, then
  there are $\mn+\nn-1-\lam$ choices
  $\tp{\mn-1,\lam-(\mn-1)},\dots, \tp{\lam-(\nn-1),\nn-1}$,
  since $0\leq \nn-\mn\leq \lam-\mn+1\leq \nn-1$ and
  $\lam-\nn+1\leq \mn-1$.  Thus, in the both sides of
  Claim~\ref{c:chi-lem-prod}, coefficients of $q^{\lam}$
  agree by Claim~\ref{c:chi-lem-min}. 
  
  Let us prove Claim~\ref{c:chi-lem-ineq-nonst}.  If
  $\mn\leq \nn$, then $\mn-\kt<\mn\leq \nn$.
  If $\mn>\nn$, then $\mn>\nn>\mn-\kt$.
  Hence, in both cases, we have
  \begin{align}
    I(\lam,\mn,\nn)
    \geq         \min(\lam+1,\mn-\kt,\nn+\kt,\mn+\nn-1-\lam)
    =I(\lam,\mn-\kt,\nn+\kt).
  \end{align}
  
  Let us prove Claim~\ref{c:chi-lem-ineq-st}.
  First,
   there is
   $\lam\in \Z$ in
   inequalities~\eqref{ineq:chi-lem-ineq-st-param},
   because
  $\mn+\nn-2\geq \nn-1$
  by 
  $\mn\geq \kt+1\geq 1$, and
  $\nn-1\geq \mn-\kt$
  by $\nn>\mn-\kt$.
  
  Second, let us prove inequality~\eqref{ineq:chi-lem-lem-ineq-st}.
  By inequalities~\eqref{ineq:chi-lem-ineq-st-param}, we have
  $\nn\geq \lam+1$, which implies
  $\mn+\nn-1-\lam \geq \mn$.  Thus,
  \begin{align}
    I(\lam,\mn,\nn)    \geq \min(\lam+1,\mn).
  \end{align}
  Moreover, inequalities~\eqref{ineq:chi-lem-ineq-st-param}
  give $\lam+1\geq \mn-\kt+1> \mn-\kt$.  Also, $\nn>\mn-\kt$
  yields $\nn+\kt>\mn-\kt$.  Thus, since
  $\mn+\nn-1-\lam\geq \mn>\mn-\kt$, we have
  \begin{align}
    \min(\lam+1,\mn)>\mn-\kt=I(\lam,\mn-\kt,\nn+\kt).
  \end{align}
\end{proof}

Suppose
$\cF=\Lam(s,l,\scc,\chi_{s,l,q},\fX)$ for $l=1$.  Then, 
$\cF$
is $\ggq$-merged-log-concave
by
Corollary~\ref{cor:merged-qstr} and the strong
$q$-log-concavity of $q$-numbers~\cite[Lemma 2.1]{Sag}.
Moreover, the
following proposition states the $\llq$-merged-log-concavity
of $\cF$.
Hence, the
following proposition
confirms the analogs of
Conjectures~\ref{conj:density}
and~\ref{conj:unimodal-width-two} for
$\Delta(\cF)(s,l,m,n,k,\fX)$, since a log-concave
$q$-polynomial of positive coefficients is $q$-dense and
unimodal.

\begin{proposition}\label{prop:chi-lc}
  Let $l=1$.  Suppose $\cF=\Lam(s,l,\scc,\chi_{s,l,q},\fX)$
  and a fitting $\mu=\tp{s,l,m,n,k}$.  Then, we have the
  following.
  \begin{enumerate}
  \item \label{c:chi-lc-merged}
    $\Delta(\cF)(s,l,m,n,k,\fX)\llq 0$.
  \item \label{c:chi-lc-palin} $\Delta(\cF)(s,l,m,n,k,\fX)$
    is a palindromic $q$-polynomial.
  \item \label{c:chi-lc-qlc} $\Delta(\cF)(s,l,m,n,k,\fX)$ is
    a log-concave $q$-polynomial.
  \end{enumerate}
\end{proposition}
\begin{proof}
  Since $\cF$ is of weight-zero, we have
  \begin{align}
    \Delta(\cF)(s,l,m,n,k,\fX) =\chi_{s,l,q,m}\chi_{s,l,q,n}-
    \chi_{s,l,q,m-\kt}\chi_{s,l,q,n+\kt}.
    \label{eq:chi-lc-merged-fitting}
  \end{align}
  Suppose that $\mu$ is unwrapped.  Then,
  $\Delta(\cF)(s,l,m,n,k,\fX)= \chi_{s,l,q,m}\chi_{s,l,q,n}$.
  This gives Claim~\ref{c:chi-lc-merged} by
  $m_{1},n_{1}\geq s_{1}\geq 1$.  Moreover,
  we obtain
  Claim~\ref{c:chi-lc-palin} by
  Proposition~\ref{prop:rho-palindromic-unimodal} and
  Claim~\ref{c:chi-lc-qlc} by Claims~\ref{c:chi-lem-min}
  and~\ref{c:chi-lem-prod} of Lemma~\ref{lem:chi-lem}.
  Hence,  assume that $\mu$ is wrapped.
  
  Let us prove Claim~\ref{c:chi-lc-merged}.  Since $\mu$ is
  wrapped,
   Lemma~\ref{lem:fitting-one} implies
  \begin{align}
    \nn+\kt>\mn>\mn-\kt\geq s_{1}\geq 1.
    \label{ineq:chi-lc-fitting-wrapped}
  \end{align}
  Thus, Claim~\ref{c:chi-lc-merged} holds by
  equation~\eqref{eq:chi-lc-merged-fitting}, since
  Claims~\ref{c:chi-lem-prod} and~\ref{c:chi-lem-ineq} of
  Lemma~\ref{lem:chi-lem} give
  \begin{dmath*}
    [\mn]_{q}[\nn]_{q}
    =\sum_{\lam\in\oi(0,\mn+\nn-2)}
    I(\lam,\mn,\nn)
    q^{\lam} \llq
    \sum_{\lam\in\oi(0,\mn+\nn-2)}
    I(\lam,\mn-\kt,\nn+\kt)
    q^{\lam} =[\mn-\kt]_{q}[\nn+\kt]_{q}.
  \end{dmath*}
  
  Let us prove Claim~\ref{c:chi-lc-palin}.  By
  Proposition~\ref{prop:rho-palindromic-unimodal},
  $\chi_{s,l,q,m}\chi_{s,l,q,n}$ and 
  $\chi_{s,l,q,m-\kt}\chi_{s,l,q,n+\kt}$ are palindromic
  $q$-polynomials.  Also,
  $\ord_{q}(\chi_{s,l,q,m}\chi_{s,l,q,n})
  =0=\ord_{q}(\chi_{s,l,q,m-\kt}\chi_{s,l,q,n+\kt})$ and
  $\deg_{q} (\chi_{s,l,q,m}\chi_{s,l,q,n})=\mn+\nn-2 =
  \deg_{q}(\chi_{s,l,q,m-\kt}\chi_{s,l,q,n+\kt})$.  Thus,
  Claim~\ref{c:chi-lc-palin} holds by
  equation~\eqref{eq:chi-lc-merged-fitting}, since
  $\Delta(\cF)(s,l,m,n,k,\fX)$ is a difference of
  palindromic $q$-polynomials with the same orders and
  degrees.
  
  Let us prove Claim~\ref{c:chi-lc-qlc}.  For each
  $\lam\in \Z$, let us define the following integers:
  \begin{align}
    P(\lam,\mn,\nn,\kt)
    &=I(\lam,\mn,\nn)-I(\lam,\mn-\kt,\nn+\kt);\\
    \del(P)(\lam,\mn,\nn,\kt)
    &=
      P(\lam,\mn,\nn,\kt)^{2}-
      P(\lam-1,\mn,\nn,\kt)       P(\lam+1,\mn,\nn,\kt).
  \end{align}
  Thus, for $\lam\in\oi(0,\mn+\nn-2)$, we prove
  \begin{align}
    \del(P)(\lam,\mn,\nn,\kt)\geq 0.
    \label{ineq:chi-lc-qlc-pointwise}
  \end{align}
  
  First, assume $\mn\leq \nn$.  Then,
  Claim~\ref{c:chi-lem-min} of Lemma~\ref{lem:chi-lem} and
  inequalities~\eqref{ineq:chi-lc-fitting-wrapped} give
  $I(\lam,\mn,\nn)$, $I(\lam,\mn-\kt,\nn+\kt)$, and
  $P(\lam,\mn,\nn,\kt)$ in Figure~\ref{fig:lc}.
  \begin{figure}[H]
    \centering     
    \includegraphics[width=14cm]      
    {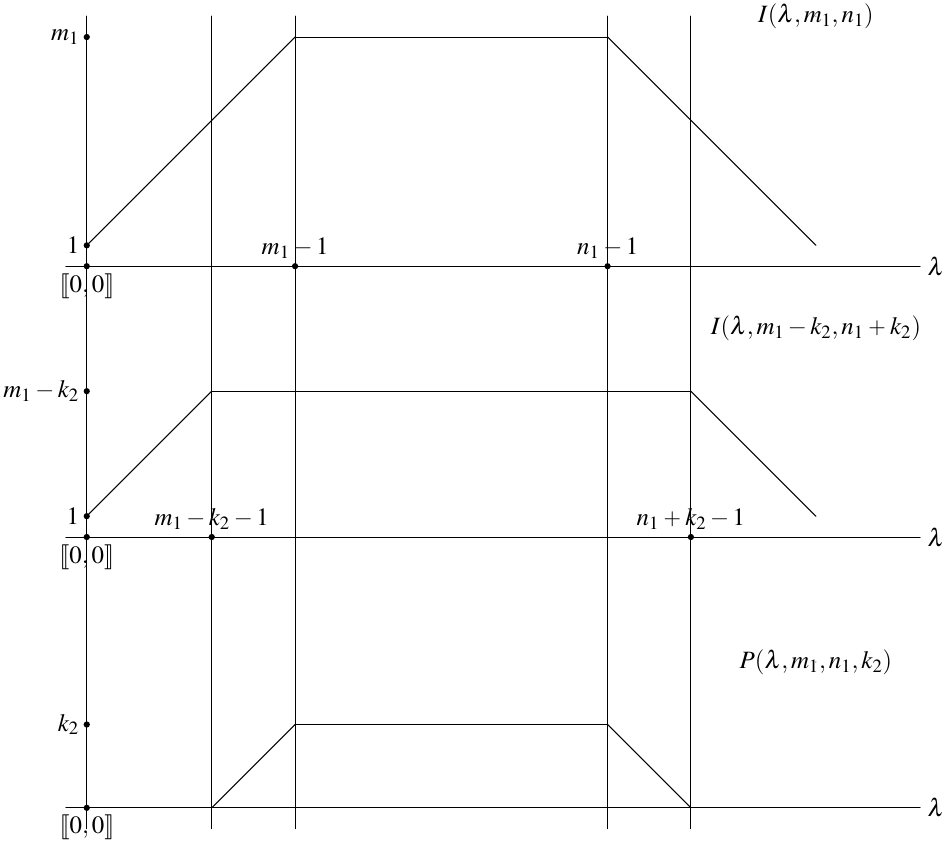}
    \caption{$P(\lam,\mn,\nn,\kt)$ when $\mn\leq \nn$}
    \label{fig:lc}
  \end{figure}
  In particular, we have the following explicit numbers:
  \begin{align}
    P(\lam,\mn,\nn,\kt)=
    \begin{dcases}
      0  &\mif \lam\in\oi(0,\mn-\kt-1),\\
      \lam+1-\mn+\kt &\mif   \lam\in\oi(\mn-\kt-1,\mn-1),\\
      \kt  &\mif   \lam\in\oi(\mn-1,\nn-1), \\
      \nn+\kt-1-\lam  &\mif   \lam\in\oi(\nn-1,\nn+\kt-1),\\
      0  &\mif \lam\in\oi(\nn+\kt-1,\mn+\nn-2).
    \end{dcases}
  \end{align}
  Hence, $\mn< \nn$ gives the following:
  \begin{align}
    \del(P)(\lam,\mn,\nn,\kt) =
    \begin{dcases}
      0  &\mif \lam\in\oi(0,\mn-\kt-1),\\
      1 &\mif \lam\in\oi(\mn-\kt,\mn-2),\\
      \kt &\mif \lam=\mn-1,\\
      0  &\mif \lam\in\oi(\mn,\nn-2),\\
      \kt  &\mif \lam=\nn-1,\\
      1  &\mif \lam\in\oi(\nn,\nn+\kt-2),\\
      0 &\mif \lam\in\oi(\nn+\kt-1,\mn+\nn-2).
    \end{dcases}
  \end{align}
  If $\mn=\nn$, then we have the following:
  \begin{align}
    \del(P)(\lam,\mn,\nn,\kt)=
    \begin{dcases}
      0 &\mif \lam\in\oi(0,\mn-\kt-1),\\
      1 &\mif \lam\in\oi(\mn-\kt,\mn-2),\\
      2\kt-1 &\mif \lam=\mn-1=\nn-1,\\
      1 &\mif \lam\in\oi(\nn,\nn+\kt-2),\\ 
      0 &\mif \lam\in\oi(\nn+\kt-1,\mn+\nn-2).
    \end{dcases}
  \end{align}
  Therefore, inequality~\eqref{ineq:chi-lc-qlc-pointwise}
  follows from $\kt\geq 1$ in
  inequalities~\eqref{ineq:chi-lc-fitting-wrapped}.
  
  Second, suppose $\mn>\nn$.  Then,
  Claim~\ref{c:chi-lem-min} of Lemma~\ref{lem:chi-lem} and
  inequalities~\eqref{ineq:chi-lc-fitting-wrapped} give
  $I(\lam,\mn,\nn)$, $I(\lam,\mn-\kt,\nn+\kt)$, and
  $P(\lam,\mn,\nn,\kt)$ in Figure~\ref{fig:lc2}.
  \begin{figure}[H]
    \centering 
    \includegraphics[width=14cm]  
    {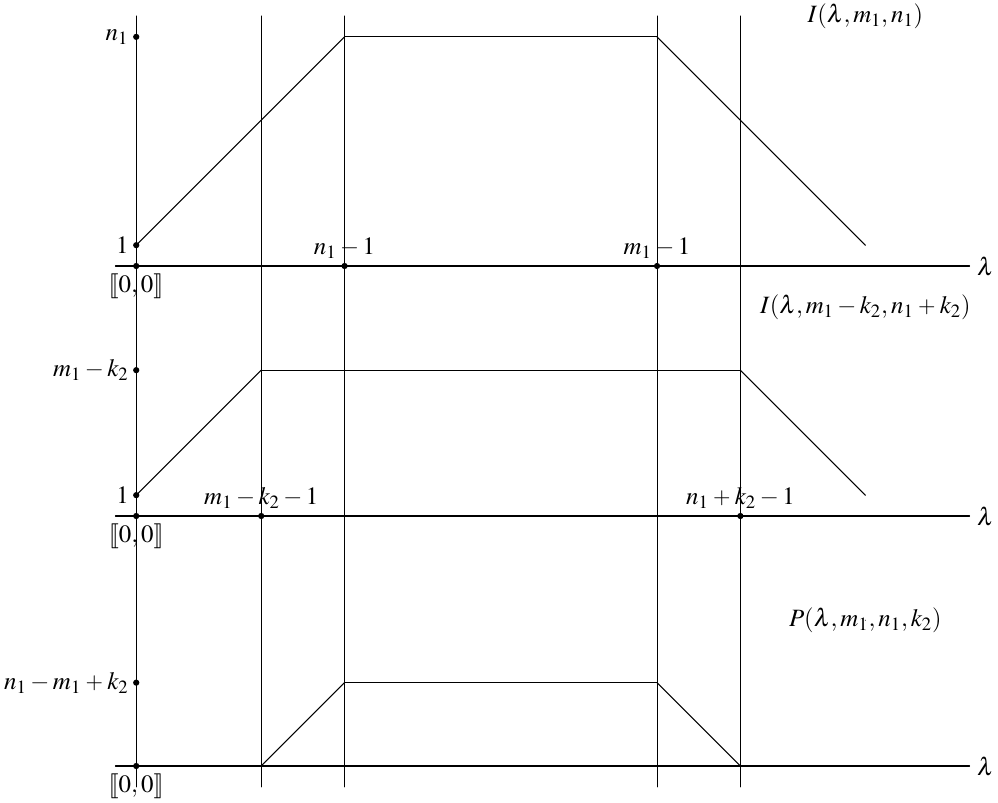}
    \caption{$P(\lam,\mn,\nn,\kt)$ when
      $\mn> \nn$}
    \label{fig:lc2}
  \end{figure}
  Thus, we have the following:
  \begin{align}
    P(\lam,\mn,\nn,\kt)=    
    \begin{dcases}
      0 &\mif  \lam\in\oi(0,\mn-\kt-1),\\
      \lam+1-\mn+\kt &\mif \lam\in\oi(\mn-\kt-1,\nn-1),\\
      \nn-\mn+\kt &\mif \lam\in\oi(\nn-1,\mn-1),\\
      \nn+\kt-1-\lam &\mif \lam\in\oi(\mn-1,\nn+\kt-1),\\
      0  &\mif \lam\in\oi(\nn+\kt-1,\mn+\nn-2).
    \end{dcases}
  \end{align}
  Then, we obtain the following:
  \begin{align}
    \del(P)(\lam,\mn,\nn,\kt)
    =
    \begin{dcases}
      0 & \mif \lam\in\oi(0,\mn-\kt-1),\\
      1 &\mif \lam \in \oi(\mn-\kt,\nn-2),\\
      \nn+\kt-\mn &\mif \lam=\nn-1,\\
      0  &\mif \lam\in\oi(\nn,\mn-2),\\
      \nn+\kt-\mn  &\mif \lam=\mn-1,\\
      1  &\mif \lam\in\oi(\mn, \nn+\kt-2),\\
      0 &\mif \lam\in\oi(\nn+\kt-1,\mn+\nn-2).
    \end{dcases}
  \end{align}
  Since $\nn+\kt>\mn$ in
  inequalities~\eqref{ineq:chi-lc-fitting-wrapped},
  inequality~\eqref{ineq:chi-lc-qlc-pointwise} holds.
\end{proof}

\subsection{Conjectures on weight-zero
  parcels}\label{sec:subsec-conj}
We state the following conjectures, motivated by
Proposition~\ref{prop:chi-lc}.

\begin{conjecture}\label{conj:qnum-uni}
  If $l\in \Zget$ and
  $\cF=\Lam(s,l,\scc,\chi_{s,l,q},\fX)$, then
  $\Delta(\cF)(s,l,m,n,k,\fX)$ is a unimodal palindromic
  $q$-polynomial for each fitting $\tp{s,l,m,n,k}$.
\end{conjecture}


\begin{example}\label{ex:non-lc-qnum}
  Let $l=2$, $m=\tp{5,2}$, $n=\tp{2,5}$, and
  $k=\tp{0,0,1,0}$. Then, since $m\bom k=\tp{4,1}$ and
  $n\bop k=\tp{3,6}$, $\cF=\Lam(s,l,\scc,\chi_{s,l,q},\fX)$
  gives
  \begin{align}
    \Delta(\cF)(s,l,m,n,k,\fX)
    &= [5]_{q}[2]_{q}[5]_{q}[2]_{q}
      -
      [4]_{q}[1]_{q}[6]_{q}[3]_{q}
    \\&
    =q^{9} + 2q^{8} + 3q^{7} + 5q^{6} +
    6q^{5} + 5q^{4} + 3q^{3} + 2q^{2} + q,
  \end{align}
  which  not a log-concave, but unimodal
  $q$-polynomial.
\end{example}

Let us recall  {\it the $q$-Starling
  polynomials of the first and second kind} $c(\ka,\lam,q)$ and
$S(\ka,\lam,q)$ in $\Z[q]$.
For each $\ka\in \Zgez$ and $\lam\in \Z$,
the delta function $\del_{i,j}$ defines
\begin{align}
  c(\ka,\lam,q)
  &  =
  \begin{dcases}
    c(\ka-1,\lam-1,q)
    +[\ka-1]_{q}c(\ka-1,\lam,q) \mif \ka\geq 1,\\
    \del_{\ka,\lam} \mif \ka=0,
  \end{dcases}\\
  S(\ka,\lam,q)
  &=
  \begin{dcases}
    S(\ka-1,\lam-1,q)+[\lam]_{q}S(\ka-1,\lam,q) \mif \ka \geq 1,\\
    \del_{\ka,\lam} \mif \ka=0.
  \end{dcases}
\end{align}
Moreover, we define the parcels below by
$c(\ka,\lam,q), S(\ka,\lam,q)\llq 0$ for $\ka\geq \lam \geq 1$.
\begin{definition}\label{def:stirling-on-tuples}
  Let $l\in \Zgeo$ and $\ka\in \Zgeo$. For the gate
  $s=\tp{1,\ka}$, we put
   families
   $c_{s,l,q}=\{c_{s,l,q,m}
   =\prod_{i\in\oi(l)} c(\ka,m_{i},q) \in\Z[q]\}_{m\in \Zl}$
   and $S_{s,l,q}=\{S_{s,l,q,m}=\prod_{i\in\oi(l)} S(\ka,m_{i},q) \in\Z[q]\}_{m\in \Zl}$.
   Furthermore, if  $q$ is $\llq$-admissible on $\fX$, then
    we call
    $\Lam(s,l,\llq,c_{s,l,q},\fX)$ and
    $\Lam(s,l,\llq,S_{s,l,q},\fX)$
    $q$-Stirling parcels of 
    the first and second kinds.
\end{definition}

Sequences
$\{c(\ka,\lam,q)\}_{\lam\in \Z}$ and
$\{S(\ka,\lam,q)\}_{\lam\in \Z}$ are strongly
$q$-log-concave~\cite[Theorems 2.4 and 2.5]{Sag}.  Hence,
$q$-Stirling parcels $\Lam(s,l,\llq,c_{s,l,q},\fX)$ and
$\Lam(s,l,\llq,S_{s,l,q},\fX)$ are $\ggq$-merged-log-concave
by Theorem~\ref{thm:merged-ext-hadam} and
Corollary~\ref{cor:merged-qstr}.  Moreover, we
conjecture the following.

\begin{conjecture}\label{conj:stir1}
  Suppose $q$-Stirling parcels
  $\cF=\Lam(s,l,\llq,c_{s,l,q},\fX)$ and
  $\cG=\Lam(s,l,\llq,S_{s,l, q},\fX)$.  Then,
   for each fitting $\tp{s,l,m,n,k}$,
  $\Delta(\cF)(s,l,m,n,k,\fX)$ and
  $\Delta(\cG)(s,l,m,n,k,\fX)$ are log-concave
  $q$-polynomials.
\end{conjecture}
\begin{example}\label{ex:q-st}
  For $l=1$ and $s=\tp{1,3}$, let
  $\cF=\Lam(s,l,\llq,c_{s,l,q},\fX)$ and
  $\cG=\Lam(s,l,\llq,S_{s,l,q},\fX)$.  Suppose $m=\tp{2}$,
  $n=\tp{1}$, and $k=\tp{0,2}$, which give $m\bom k=\tp{0}$
  and $n\bop k=\tp{3}$.  Then, we have the following
  non-palindromic, but log-concave
  $q$-polynomials:
  \begin{align}
    \Delta(\cF)(s,l,m,n,k,\fX)
    &= c(3,2,q)c(3,1,q)
      =q^{2}+3q+2;\\
    \Delta(\cG)(s,l,m,n,k,\fX)
    &= S(3,2,q)S(3,1,q)
      =q+2.
  \end{align}

\end{example}

Let us recall the notion of {\it strong $q$-log-convexity}.

\begin{definition}\cite[Section 1]{CWY}
  A sequence $f=\{f_{m}\in \R[q]\}_{m\in \Zgez}$ is called
  strongly $q$-log-convex, if $n\geq m\geq 1$ implies
  $-(f_{m}f_{n}-f_{m-1}f_{n+1})\ggq0$.
\end{definition}

\begin{remark}
  Let $\cF=\Lam(s,l,w,\scc,\phi,\rho,\fs,x,\fX)$ be a parcel.
  Then, there is the notion of $\scc'$-merged-log-convexity
  (or $\sce'$-merged-log-convexity) of $\cF$ such that each
  wrapped fitting $\tp{s,l,m,n,k}$ gives
  $-\Delta(\cF)(s,l,w,m,n,k,\phi,\rho,x,\fX)\scc' 0$ (or $\sce' 0$).
\end{remark}

Moreover, we put parcels by {\it the Ramanujan polynomials}
$R_{\lam+1}(q)$ and {\it Bessel polynomials} $B_{\lam}(q)$
in $\Z[q]$ such that each $\lam\in\Zgez$ satisfies
\begin{align}
  R_{\lam+1}(q)
  &=
    \begin{dcases}
      \lam(1+q)R_{\lam}(q)+
      q^{2}\frac{ d R_{\lam}(q)}{d q} \mif \lam\in \Zgeo,\\
      1 \mif \lam=0,
  \end{dcases}\\
  B_{\lam}(q)
  &=\sum_{\ka\in\oi(0,\lam)}\frac{(\lam+\ka)!}{
    (\lam-\ka)!\ka!}\frac{q^{\ka}}{2^{\ka}}.
\end{align}

\begin{definition}\label{def:rama-bessel}
  Suppose $l\in \Zgeo$ and a gate $s\geq 0$.  Then, let
  $R_{s,l,q}=\{R_{s,l,q,m}\in \Z[q]\}_{m\in \Zl}$ and
  $B_{s,l,q}=\{B_{s,l,q,m}\in \Z[q]\}_{m\in \Zl}$ be families
  such that
    \begin{align}
      R_{s,l,q,m}
      &=
        \begin{dcases}
          \prod_{i\in\oi(l)}   R_{m_{i}}(q) \mfor m\ldZl s,\\
          0 \melse,
        \end{dcases}\\
      B_{s,l,q,m}
      &=
        \begin{dcases}
          \prod_{i\in\oi(l)}   B_{m_{i}}(q) \mfor m\ldZl s,\\
          0 \melse.
        \end{dcases}
    \end{align}
\end{definition}
Sequences $\{R_{\lam}(q)\}_{\lam\in \Zgeo}$ and
$\{B_{\lam}(q)\}_{\lam\in \Zgez}$ are strongly
$q$-log-convex~\cite[Corollaries 3.2 and 3.3]{CWY}.  But, we
conjecture the following log-concavity.
\begin{conjecture}\label{conj:rama-bessel}
  Suppose $\cF=\Lam(s,l,\llq,R_{s,l,q},\fX)$ and
  $\cG=\Lam(s,l,\llq,B_{s,l,q},\fX)$ for $s=\tp{1,\infty}$.
  Then,
   $\Delta(\cF)(s,l,m,n,k,\fX)$ and
  $\Delta(\cG)(s,l,m,n,k,\fX)$ are log-concave
  $q$-polynomials for each fitting $\tp{s,l,m,n,k}$.
\end{conjecture}

\begin{example}\label{ex:rama}
  For $l=1$ and $s=\tp{1,\infty}$, let
  $\cF=\Lam(s,l,\llq,R_{s,l,q},\fX)$ and
  $\cG=\Lam(s,l,\llq,B_{s,l,q},\fX)$.
  If $m=\tp{3}$, $n=\tp{2}$, and $k=\tp{0,2}$, then
  $m\bom k=\tp{1}$ and $n\bop k=\tp{4}$ give the
  following
  non-palindromic, but log-concave $q$-polynomials:
  \begin{align}
    -\Delta(\cF)(s,l,m,n,k,\fX)
    &=-(R(3,q)R(2,q)-R(1,q)R(4,q))\\
    &=12 q^3 + 18 q^2 + 12q + 4;\\
    -\Delta(\cG)(s,l,m,n,k,\fX)
    &=-( B(3,q)B(2,q)-B(1,q)B(4,q))\\
    &=60q^{5} + 120q^4 + 72q^3 + 19q^2 + 2q.
  \end{align}
\end{example}

\section{Almost strictly unimodal sequences and Young diagrams}
\label{sec:semi-st-unim-seq-young}

By the
merged-log-concavity,
we discuss almost strictly unimodal sequences.
In particular, Young diagrams give
infinite-length such sequences.

\subsection{Strict log-concavity
  and merged-log-concavity}
\label{sec:st-lc-merged-for-almost strictly-unim}

Suppose a strictly log-concave  sequence
$r=\{r_{i}\in \R_{>0}\}_{i\ldZ s}$.  Then,
$j-1,j+2\ldZ s$ give
\begin{align}
  r_{j}^{2}-r_{j-1}r_{j+1}&>0,\\
  r_{j+1}^{2}-r_{j}r_{j+2}&>0.
                            \label{ineq:st-lc-one-minus-one}
\end{align}
 Hence, $r$ is almost strictly
unimodal by
$\frac{r_{j}}{r_{j-1}}> \frac{r_{j+1}}{r_{j}}>
\frac{r_{j+2}}{r_{j+1}}$.
However, inequality~\eqref{ineq:st-lc-one-minus-one}
do not always extend for merged-log-concave parcels.

For instance, consider a $\scc'$-merged-log-concave
$\cF=\Lam(s,l,w,\scc,\fs,\phi,\rho,x,\fX)$. Then,
fitting $\tp{s,l,m,m,k}$ and $\tp{s,l,m\bop k, m\bop k,k}$
yield
\begin{align}
  \Ups(s,l,w,m,m,k,\phi,\rho,x,\fX)
  (\cF_{m}\cF_{m^{\ve}}-
  \cF_{m\bom k} \cF_{(m \bop k)^{\ve}})
  &\scc' 0,\\
  \Ups(s,l,w,m\bop k,m\bop k,k,\phi,\rho,x,\fX)
  (\cF_{m\bop k}\cF_{(m\bop k)^{\ve}}-
  \cF_{(m\bop k)\bom k}\cF_{((m\bop k)\bop k)^{\ve}})
  &\scc' 0.
\end{align}
However, in general, we do not have
$\cF_{(m\bop  k)\bom k}=\cF_{m}$,
which corresponds to
\begin{align}
  r_{(j+1)-1}=r_{j}
  \label{eq:one-minus-one}
\end{align}
in inequality~\eqref{ineq:st-lc-one-minus-one}.

\subsection{Fitting paths}
\label{sec:fitting-paths}

We define the following {\it fitting paths}
to discuss almost strictly unimodal sequences
by the merged-log-concavity.

\begin{definition}\label{def:fitting-paths-len-one}
  Suppose a gate $\theta$ such that
  $\thetat-\thetao=1$.  Consider a sequence
  $P=\{P_{i}=\tp{s,l,m_{i},n_{i},k_{i}}\}_{i\ldZ \theta}$
  such that $P_{\thetao}$ and $P_{\thetat}$ are fitting tuples,
  $k_{\thetao}$ and $k_{\thetat}$ are $\sig$-equivalent, and
  $m_{\thetat}=(n_{\thetao}\bop k_{\thetao})^{\ve}$.
  \begin{enumerate}
  \item We call $P$ fitting path of type zero, if
    $m_{i}=m_{i}^{\ve}=n_{i}$ for each $i\ldZ\theta$.
  \item Assume $m_{\thetao}=n_{\thetao}$.
    \begin{enumerate}
    \item We call $P$ fitting path of type 1-1, if
      $m_{\thetat}=n_{\thetat}$.
    \item We call $P$ fitting path of type 1-2, if
      $m_{\thetat}=n_{\thetat}^{\ve}$.
    \end{enumerate}
  \item Assume 
    $m_{\thetao}^{\ve}=n_{\thetao}$.
    \begin{enumerate}
    \item We call $P$ fitting path of type 2-1, if
      $m_{\thetat}=n_{\thetat}$.
    \item We call $P$ fitting path of type 2-2, if
      $m_{\thetat}= n_{\thetat}^{\ve}$.
    \end{enumerate}
  \end{enumerate}
\end{definition}
We define fitting paths of general lengths.

\begin{definition}\label{def:fitting-paths-gen-len}
  Let $\theta$ be a gate.  Consider a sequence
  $P=\{P_{i}=\tp{s,l,m_{i},n_{i},k_{i}}\}_{i\ldZ \theta}$ such
  that $Q_{j}=\{P_{i}\}_{i\ldZ \tp{j,j+1}}$ is a fitting
  path for each $j\in\oi(\thetao,\thetat-1)$.
  \begin{enumerate}
  \item We call $P$ fitting path. 
  \item We call $P$ fitting path of the type $A$, if
    $Q_{j}$
    has the same type $A$ for all
    $j\in\oi(\thetao,\thetat-1)$.
  \end{enumerate}
  We refer to $s$, $l$, $\theta$, and
  $\{k_{i,1}\}_{i\ldZ \theta}$ of $P$ as the gate, width,
  length, free-path parameters of $P$.  Moreover, if
  $\thetat-\thetao=n\in \Zgez$, then $P$ is said to be
  length-$n$ or finite-length; otherwise, $P$ is said to be
  $P$ infinite-length.
  
  We also call a sequence
  $P=\{P_{i}=\tp{s,l,m_{i},n_{i},k_{i}}\}_{i\ldZ \theta}$
  length-zero fitting path, if $\thetat-\thetao=0$ and
  $P_{\thetao}$ is fitting.
\end{definition}

On the terminology of free-path parameters of $P$, suppose
$k_{i,1}+\lam_{i}\in \Zgez$ for each $i\ldZ \theta$.  Then,
$p_{i,\lam_{i}}=\tp{s,l,m_{i},n_{i},\tp{k_{i,1}+\lam_{i}}\ccn
  k(2,l)}$ is always fitting as in
Remark~\ref{rmk:free-param}.  Also,
$\{p_{i,\lam_{i}}\}_{i\ldZ \theta}$ is a fitting path by
$\sig(k_{i})= \sig(\tp{k_{i,1}+\lam_{i}}\ccn k(2,l))$.

\begin{example}\label{ex:paths-same-tuples}
  Let $s=\tp{0,\infty}$, $l=1$, $\mn=\tp{1}$, and
  $\kn=\tp{0,1}$.  Then, $\mu_{1}= \tp{s,l,\mn,\mn,\kn}$ is
  fitting, because $\mn\ldZl s$,
  $\nu(\kn)=\tp{0,1}$,
  and $\nu(\mn,\mn,\kn)=\nu(\kn)+\mn\ccn \mn=\tp{1,2}$.
  Also, let $\mt=(\mn\bop \kn)^{\ve}=\tp{2}$.  Then,
  $\mu_{2}=\tp{s,l,\mt,\mt,\kt}=\tp{s,l,(\mn\bop k)^{\ve},
    (\mn\bop k)^{\ve},\kn}$ is fitting, because 
  $\mt\ldZl s$ and
  $\nu(\mt,\mt,\kt)=\nu(\kt)+\mt\ccn \mt=\tp{2,3}$.
  Hence, $\{\mu_{i}\}_{i\ldZ\tp{1,2}}$ is a
  fitting path of type 1-1, and of type zero by
  $l=1$.
\end{example}

\begin{example}\label{ex:paths-equiv-tuples}
  The following different but equivalent  supports
  give a fitting path.
  Let $l=2$, $s=\tp{0,\infty}$, $\mn=\tp{0,3}$, and
  $\kn=\tp{0,4,4,4}$.  Then,
  $\mu_{1}=\tp{s,l,\mn,\mn^{\ve},\kn}$ is fitting, since
  \begin{align}
    \nu(\kn)
    &=\tp{0,4,8,12},\\
    \nu(\mn,\nn,\kn)
    &=\nu(\kn)+\mn\ccn \mn^{\ve}=\tp{0,7,11,12}.
  \end{align} 
  
  Should we continue with $\kn$, then either
  $\tp{s,l,(\mn^{\ve}\bop \kn)^{\ve}, \mn^{\ve}\bop
    \kn,\kn}$ or
  $\tp{s,l,(\mn^{\ve}\bop \kn)^{\ve}, (\mn^{\ve}\bop
    \kn)^{\ve},\kn}$ has to be fitting.  However, neither of
  them is fitting, because
  $\mn^{\ve}\bop \kn =\tp{3,0}+\tp{4,12}=\tp{7,12}$ implies
  \begin{align}
    \nu((\mn^{\ve}\bop \kn)^{\ve}, \mn^{\ve}\bop \kn,\kn)
    &=\nu(\kn)+\tp{12,7}\ccn \tp{7,12}=\tp{12,11,15,24},\\
    \nu((\mn^{\ve}\bop \kn)^{\ve},
    (\mn^{\ve}\bop \kn)^{\ve},\kn)
    &=\nu(\kn)+\tp{12,7}\ccn \tp{12,7}=\tp{12,11,20,19}.
  \end{align}
  
  Hence, consider $\kt=\tp{0,7,4,1}\neq \kn$. Then, not only
  $\kn$ and $\kt$ are equivalent, but also
  $\mu_{2}=\tp{s,l,\mt,\nt,\kt}=\tp{s,l,(\mn^{\ve}\bop
    \kn)^{\ve}, \mn\bop \kn,\kt}$ is fitting, because
\begin{align}
    \nu(\kt)
    &=\tp{0,7,11,12},\\
    \nu(\mt,\nt,\kt)
    &=\nu(\kt)+\tp{12,7}\ccn \tp{7,12}
      =\tp{12,14,18,24}.
\end{align}
Thus, $\{ \mu_{i}\}_{i\in\oi(2)}$ is a type 2-2 fitting path
with different but equivalent  supports $\kn, \kt$.
\end{example}

Moreover, we introduce the following notion of parcels.

\begin{definition}\label{def:static-parcel} 
  Let $\cF=\Lam(s,l,w,\scc,\fs,\phi,\rho,x,\fX)$.  We call
  $\cF$ static at $r\in \AfX$, if each $m\in \Zl$ gives
  $\cF_{m}(r)=\cF_{m^{\ve}}(r)$.  We simply call $\cF$
  static, if $\cF$ is static at every $r\in \AfX$.
\end{definition}

For instance, we have
the following static monomial parcels (see
Remark~\ref{rmk:without-tuple-flips}).
\begin{lemma}\label{lem:static-parcel}
  Suppose $\cF=\Lam(s,l,w,\scc,\fs,\phi,\rho,x,\fX)$ and
  $r\in \AfX$ such that $x$, $\phi(x)$, and $w$ are
  palindromic,
  and $f_{s,m}(r)=f_{s,m^{\ve}}(r)$ for each $m\in \Zl$.
  Then, $\cF$ is static at $r$.
\end{lemma}
\begin{proof}
  Since
  $\cF_{m}=\frac{f_{m}}{\prod \phi(x)^{m\rc w}\cdot
    [m]!_{x}^{w}}$ for
  each $m\ldZl s$, the statement follows.
\end{proof}

\begin{proposition}\label{prop:static-mono}
  Suppose a monomial index $\tp{l,w,\gam}$ with
  palindromic
  $w$, $\tp{\gam_{i,1}}_{i\in\oi(l)}$, and
  $\tp{\gam_{i,2}}_{i\in\oi(l)}$.  Then, a monomial parcel
  $\Lam(s,l,w,\scc,\Psi_{s,\gam,q},\rho,x,\fX)$ is
  static.
\end{proposition}
\begin{proof}
  The statement follows from
  Definition~\ref{def:monomial-parcel} and
  Lemma~\ref{lem:static-parcel}.
\end{proof}

To discuss further, we prove the following 
on $\sig$-sum and $\sig$-minus.

\begin{lemma}\label{lem:bombop-wedge}
  Let $l\in \Zgeo$, $m\in \Zl$, and $k,k'\in \Ztl$.  Then, we
  have the following equations.
  \begin{align}
    (m\bop k)^{\ve}&=m^{\ve}\bom (-k),
                     \label{eq:bop-bom}\\
    (m\bom k)^{\ve}&=m^{\ve} \bop (-k),
                     \label{eq:bom-bop}\\
    m\bop k\bop k'&=m\bop (k+k'),
                    \label{eq:bopbop-bop}\\
    m\bom k\bom k'&=m\bom (k+k'),
                    \label{eq:bombom-bom}\\
    (m\bop k)^{\ve}\bom k&=m^{\ve},
                           \label{eq:bopwedbom-wed}\\
    (m\bom k)^{\ve}\bop k&=m^{\ve},
                           \label{eq:bomwedbop-wed}\\
    m\bom k \bop k'&= m\bop k' \bom k.
                     \label{eq:bombop-bopbom}
  \end{align}
\end{lemma}
\begin{proof}
  We have
  $(m\bop k)^{\ve} =m^{\ve}+(\sig(k)^{\ve})^{\ve}=
  m^{\ve}-\sig(-k)=m^{\ve}\bom (-k)$.  This gives
  equation~\eqref{eq:bop-bom}. Then,
  equation~\eqref{eq:bom-bop} follows, since
  equation~\eqref{eq:bop-bom} implies
  $(m\bom k)^{\ve}=((m^{\ve}\bop (-k))^{\ve})^{\ve}
  =m^{\ve}\bop (-k)$.
  
  We obtain equation~\eqref{eq:bopbop-bop}, since
  $m\bop k \bop k'= m\bop k+\sig(k')^{\ve}
  =m+\sig(k)^{\ve}+\sig(k')^{\ve} =m+\sig(k+k')^{\ve}=m\bop
  (k+k')$.  Likewise, 
  equation~\eqref{eq:bombom-bom} follows.

  We obtain equation~\eqref{eq:bopwedbom-wed}, because
  equations~\eqref{eq:bop-bom} and~\eqref{eq:bombom-bom}
  give
  $(m\bop k)^{\ve}\bom k =m^{\ve}\bom (-k)\bom k =m^{\ve}$.
  Similarly, equation~\eqref{eq:bomwedbop-wed} follows from
  equations~\eqref{eq:bom-bop} and~\eqref{eq:bopbop-bop}.
  Since
  $m\bom k \bop k'=m-\sig(k)+\sig(k')^{\ve}
  =m+\sig(k')^{\ve}-\sig(k) =m \bop k' \bom k$, 
  equation~\eqref{eq:bombop-bopbom} holds.
\end{proof}
We now realize the following analog of
equation~\eqref{eq:one-minus-one}.
\begin{proposition}
  \label{prop:static-parcel-one-minus-one}
  Consider a fitting path
  $P=\{P_{i}=\tp{s,l,m_{i},n_{i},k_{i}}\}_{i\ldZ \theta}$ and
  static parcel $\cF=\Lam(s,l,w,\scc,\fs,\phi,\rho,x,\fX)$ at
  $r\in \AfX$.  Then, each $j\in\oi(\thetao+1,\thetat)$ 
  satisfies
\begin{align}
  \cF_{m_{j}\bom k_{j}}(r)=\cF_{n_{j-1}}(r)=\cF_{m_{j-1}}(r).
\end{align}
\end{proposition}
\begin{proof}
  Since
  $m_{j}=(n_{j-1}\bop k_{j-1})^{\ve}$ or
$m_{j}=(n_{j-1}^{\ve}\bop k_{j-1})^{\ve}$,
$m_{j}\bom k_{j}=n_{j-1}$ or $n_{j-1}^{\ve}$ by
 Claims~\ref{eq:bop-bom}
 and~\ref{eq:bombom-bom} of Lemma~\ref{lem:bombop-wedge}.
Also, $m_{j-1}=n_{j-1}$ or $n_{j-1}^{\ve}$ by
Definition~\ref{def:fitting-paths-len-one}.
\end{proof}

\subsection{Positivity of ring shift factors}
We generalize
inequality~\eqref{ineq:intro-pos-rshift-factor} by
ring shift factors in Definition~\ref{def:merged}.
Let us first
obtain some $\gAfX$-positivities of mediators and
$q$-numbers by
the following lemma.

\begin{lemma}\label{lem:gAfX-closedness}
  The set $C=\{f\in \Q(\fX)\mid f\gAfX 0\}$ is
  a group under multiplications.
\end{lemma}
\begin{proof}
  Clearly, $1\in C$. If $g\in C^{2}$, then
  $\Qu(g)(r)=\Qu(g(r))>0$ for each $r\in \AfX$.
\end{proof} 

\begin{proposition}\label{prop:AfX-pos}
  Suppose a $\mu=\tp{s,l,w,\scc,\rho,x,\fX}$-mediator $\phi$.
  \begin{enumerate}
  \item \label{c:AfX-pos-med-qfact} Let $m\in \Zlgez$ and
    $\lam\in \Zgeo$.  Then, 
    \begin{align}
      \prod \phi(x)^{m\rc w}& \gAfX 0,
                       \label{ineq:AfX-pos-med-qfact-med}\\
      [m]!_{x^{\lam}}^{w}&\gAfX 0.
                           \label{ineq:AfX-pos-med-qfact-qfact}
    \end{align}
  \item \label{c:AfX-pos-med-bshift} Let $m\ldZl s$. Then,
        $\prod \phi(x^{\rho})^{m\rc w}\gAfX 0$.
  \item \label{c:AfX-pos-med-bshift-fitting} Suppose a
    fitting $\tp{s,l,m,n,k}$ with $a=\nu(k)$ and
    $b=\nu(m,n,k)$. Then, we have
    $\prod (\phi(x^{\rho})^{\wcn})^{(b-a)\rc w^{\wcn}}
    \gAfX 0$.
  \end{enumerate}
\end{proposition}
\begin{proof}
  Let us prove Claim~\ref{c:AfX-pos-med-qfact}.  First,
  inequality~\eqref{ineq:AfX-pos-med-qfact-med} follows from
  Lemma~\ref{lem:gAfX-closedness} and the base positivity of
  $\phi$ and $\mu$.
  Second, we prove
  inequality~\eqref{ineq:AfX-pos-med-qfact-qfact}.  By
  $\lam\in \Zgeo$,
  each $i\in\oi(l)$ implies
  $[m_{i}]!^{w_{i}}_{x^{\lam}_{i}} >_{x_{i}} 0$,
  and so
  $[m]!_{x^{\lam}}^{w}\scc 0$  by
   the half $\gx$-$\scc$
  implication in Claim~\ref{c:adm-succ-half-gx-scc} of
  Lemma~\ref{lem:adm-succ}.  Hence, 
  inequality~\eqref{ineq:AfX-pos-med-qfact-qfact} follows
  from the half $\scc$-$\gAfX$ implication.
  
  Let us prove Claim~\ref{c:AfX-pos-med-bshift}.  The
  base-shift positivity of $\phi$ and $\mu$
  and the compatibility $\gAfX \Sup \scc$
  imply
  $B(s,l,w,m,\phi,\rho,x,\fX)\gAfX 0$.
  Then, by
  Definition~\ref{def:bshift-fun-multi}, we have
  \begin{align}
    B(s,l,w,m,\phi,\rho,x,\fX)
    &=\prod_{i\in\oi(l)}
      \frac{ \phi(x_{i}^{\rho_{i}})^{w_{i}m_{i}}
      [m_{i}]!_{x_{i}^{\rho_{i}}}^{w_{i}}}
      { \phi(x)^{w_{i}m_{i}}      [m_{i}]!_{x_{i}}^{w_{i}}}
   =
      \frac{ \prod \phi(x^{\rho})^{m\rc w}\cdot
      [m]!_{x^{\rho}}^{w}}{
      \prod \phi(x)^{m\rc w}\cdot
      [m]!_{x}^{w}}
    \gAfX 0.
  \end{align}
  Therefore, Claim~\ref{c:AfX-pos-med-bshift} follows from
  Lemma~\ref{lem:gAfX-closedness} and
  Claim~\ref{c:AfX-pos-med-qfact}.
  
  Claim~\ref{c:AfX-pos-med-bshift} gives
  Claim~\ref{c:AfX-pos-med-bshift-fitting},
since
  $\prod (\phi(x^{\rho})^{\wcn})^{(b-a)\rc w^{\wcn}}
    = \prod (\phi(x^{\rho})^{\wcn})^{(m\ccn n)\rc w^{\wcn}}
    = \prod \phi(x^{\rho})^{m\rc w}
    \cdot \prod \phi(x^{\rho})^{n^{\ve}\rc w}$
     for $m,n^{\ve}\ldZl s$.
\end{proof}
We now generalize
inequality~\eqref{ineq:intro-pos-rshift-factor}
as follows.
\begin{corollary}\label{cor:rshift-factor-pos}
  Suppose a fitting $\tp{s,l,m,n,k}$.  Then,
  \begin{align}
    \Ups(s,l,w,m,n,k,\phi,\rho,x,\fX) \gAfX 0.    
  \end{align}
\end{corollary}
\begin{proof}
  Let $a=\nu(k)$ and $b=\nu(m,n,k)$. Then, 
  $a,b\geq 0$ by the slope conditions.  Thus,
  Claim~\ref{c:AfX-pos-med-qfact} of
  Proposition~\ref{prop:AfX-pos}  implies
  $[b]!_{(x^{\rho})^{\wcn}}^{w^{\wcn}}=[b(1,l)]!_{x^{\rho}}^{w}
    \cdot [b(2l,l+1)]!_{x^{\rho}}^{w}\gAfX 0$ and
    $[a]!_{(x^{\rho})^{\wcn}}^{w^{\wcn}}
    =[a(1,l)]!_{x^{\rho}}^{w}
    \cdot [a(2l,l+1)]!_{x^{\rho}}^{w}
    \gAfX 0$.
    This gives the assertion by
    Item~\ref{d:merged-factor-det} of
    Definition~\ref{def:merged},
    Lemma~\ref{lem:gAfX-closedness}, and
    Claim~\ref{c:AfX-pos-med-bshift-fitting} of
    Proposition~\ref{prop:AfX-pos}.
\end{proof}
\subsection{Almost strictly unimodal sequences
  by the merged-log-concavity}
\label{sec:semi-st-unim-seq-by-merged}

We discuss almost strictly unimodal sequences by fitting
paths and the merged-log-concavity.  First, squaring orders
give the following polynomials of rational functions.
\begin{lemma}\label{lem:geAfX-gAfX-scc}
  For $f\in \Q(\fX)$, consider $\Ir(\fX)(f)\in \Q[\fX]^{2}$.
  \begin{enumerate}
  \item \label{c:geAfX-gAfX-scc-geAfX} $f\geAfX 0$
    if and only if  there exists
    $f(r)=\Qu(\Ir(\fX)(f)(r))\in \R_{\geq 0}$ for each $r\in \AfX$.
  \item \label{c:geAfX-gAfX-scc-gAfX} $f\gAfX 0$ if and
    only
    if there exists
    $f(r)=\Qu(\Ir(\fX)(f)(r))\in
    \R_{>0}$ for each $r\in \AfX$.
  \item \label{c:geAfX-gAfX-scc-scc} Suppose a squaring
    order $\scc$ on $\fX$.  If $f^{2}\scc 0$, then
    $(\prod \Ir(\fX)(f))(r)\neq 0$ for each $r\in \AfX$.
  \end{enumerate}
\end{lemma}
\begin{proof}
  Let us prove Claim~\ref{c:geAfX-gAfX-scc-geAfX}.  Assume
  $f\geAfX 0$.  Then, $\Ir(\fX)(f)_{1}(r)\neq 0$ for each
  $r\in \AfX$, since $f(r)\in \R$ must exist.  Thus, the only
  if part holds. The if part is clear.
  Claim~\ref{c:geAfX-gAfX-scc-gAfX} holds similarly.
  Claim~\ref{c:geAfX-gAfX-scc-scc} follows from
  Claim~\ref{c:geAfX-gAfX-scc-gAfX}, since
  $f^{2}=\frac{\Ir(\fX)(f)_{2}^{2}}{\Ir(\fX)(f)_{1}^{2}}\gAfX
  0$ by $\gAfX \Sup \scc$.
\end{proof}

We use the following notion
for the total positivity or the total negativity of
parcels.
\begin{definition}
  Let $\cF=\Lam(s,l,w,\scc,\fs,\phi,\rho,x,\fX)$.  Then, we call
  $\cF$ $\tp{s,l,\fX}$-positive (or
  $\tp{s,l,\fX}$-negative), provided $\cF_{m}(r)>0$ (or
  $<0$) for each $m\ldZl s$ and $r\in \AfX$.
\end{definition}
\begin{lemma}\label{lem:total-pos-neg}
  Let $\cF=\Lam(s,l,w,\scc,\fs,\phi,\rho,x,\fX)$.  Then, we have
  the following.
  \begin{enumerate}
  \item
    \label{c:total-pos-neg-pos}
    $\cF$ is $\tp{s,l,\fX}$-positive if and only if
    $\cF_{m}(r)>0$ for some $m\ldZl s$ and $r\in \AfX$.
  \item
    \label{c:total-pos-neg-neg}
    $\cF$ is $\tp{s,l,\fX}$-negative if and only if
    $\cF_{m}(r)<0$ for some $m\ldZl s$ and $r\in \AfX$.
  \item
    \label{c:total-pos-neg-pos-or-neg}
    $\cF$ is either $\tp{s,l,\fX}$-positive or negative.
  \end{enumerate}
\end{lemma}
\begin{proof}
  Let us prove Claim~\ref{c:total-pos-neg-pos}.  The only if
  part is clear.  So, let us prove the if part.  First, for
  each $\lam\in \AfX$, we prove
  \begin{align}
    \cF_{m}(\lam)>0.
    \label{ineq:total-pos-neg-pos-AfX}
  \end{align}
  Claim~\ref{c:AfX-pos-med-qfact} of
  Proposition~\ref{prop:AfX-pos} gives
  \begin{align}
    \prod \phi(x)^{m\rc w} \cdot [m]!_{x}^{w}& \gAfX 0.
                                    \label{ineq:total-pos-neg-med-pos}
  \end{align}
  Thus,
  $\cF_{m}(r)>0$ implies
  \begin{align}
    f_{s,m}(r)>0.
    \label{ineq:total-pos-neg-point-wise-pos-numerator}
  \end{align}
  Also, $f_{s,m}^{2}\gAfX 0$, since $f_{s,m}^{2}\scc 0$
  for
  the pairwise positive $\fs$. Thus,
  $\Ir(\fX)(f_{s,m})\in \Q[\fX]^{2}$ satisfies
  $(\prod \Ir(\fX)(f_{s,m}))(\lam)\neq 0$ 
  for each $\lam\in \AfX$
  by
  Claim~\ref{c:geAfX-gAfX-scc-scc} in
  Lemma~\ref{lem:geAfX-gAfX-scc}.   Then,
  $f_{s,m}(\lam)>0$ for each $\lam\in \AfX$ by
  inequality~\eqref{ineq:total-pos-neg-point-wise-pos-numerator},
  since $\Ir(\fX)(f_{s,m})_{1}$ and $\Ir(\fX)(f_{s,m})_{2}$
  are continuous on $\AfX$.  Therefore,
  inequality~\ref{ineq:total-pos-neg-pos-AfX}
  follows from
  inequality~\eqref{ineq:total-pos-neg-med-pos}.  
  
  Second, each $n\ldZl s$ gives $f_{s,m}f_{s,n}\gAfX 0$
  for
  the pairwise positive $\fs$.  Moreover,
  $f_{s,m}(r)f_{s,n}(r)>0$ implies that the signs of
  $f_{s,m}(r)$ and $f_{s,n}(r)$ coincide.  Therefore,
  $\cF_{n}(r)>0$ by
  inequality~\eqref{ineq:total-pos-neg-med-pos}. This gives
  Claim~\ref{c:total-pos-neg-pos}.
  Claim~\ref{c:total-pos-neg-neg} holds similarly.
  Claim~\ref{c:total-pos-neg-pos-or-neg} follows
  from Claims~\ref{c:total-pos-neg-pos}
  and~\ref{c:total-pos-neg-neg}, since $\cF_{m}(r)\neq 0$
  for
  each $m\ldZ s$ and $r\in \AfX$.
\end{proof}

We now introduce the following notation for real numbers
along fitting paths.

\begin{definition}\label{def:path-parcel-sequences}
  Let $\theta$ be a gate.
  Consider the pair $\zeta=\tp{P,\cF}$ of a fitting path
  $P=\{\tp{s,l,m_{i},n_{i},k_{i}}\}_{i\ldZ \theta}$ and parcel
  $\cF=\Lam(s,l,w,\scc,\fs,\phi,\rho,x,\fX)$.  Let
  $r\in \AfX$.
  \begin{enumerate}
  \item We define the extended gate
    $e(\theta)=\tp{\thetao-1,\thetat+1}\in\hbzt$.
  \item We define the path-parcel sequence
    $u(\zeta,r)=\{u(\zeta,r)_{i}\in \R\}_{i\ldZ e(\theta)}$ such
    that
    \begin{align}
      u(\zeta,r)_{i}
      &=
        \begin{dcases}
          \cF_{m_{\thetao}\bom k_{\thetao}}(r)
          &\mif i=\thetao-1,\\
          \cF_{n_{i}}(r)
          &\mif                    i\in\oi(\thetao,\thetat),\\
          \cF_{n_{\thetat\bop k_{\thetat}}}(r)
          &\mif i=\thetat+1<\infty.
        \end{dcases}
    \end{align}
  \end{enumerate}
\end{definition}
By the ends
 of parcel-path sequences,
we put the notion of {\it wrapped fitting paths}.
\begin{definition}\label{def:wrapped-paths}
  Suppose a fitting path
  $P=\{\tp{s,l,m_{i},n_{i},k_{i}}\}_{i\ldZ \theta}$.  We call
  $P$ wrapped, if $P$ satisfies the following conditions:
  \begin{align}
    m_{\thetao}&\bom k_{\thetao}\ldZl s;
                 \label{inc:wrapped-path-lower}\\
    n_{\thetat}&\bop k_{\thetat} \ldZl s
                 \mwhen \thetat<\infty.
                 \label{inc:wrapped-path-upper}
  \end{align}
  We refer to \eqref{inc:wrapped-path-lower}
  and~\eqref{inc:wrapped-path-upper} as the lower inclusion
  condition and the upper inclusion condition of $P$.
\end{definition}

We also state the following lemma  to
obtain almost strictly unimodal sequences by
the merged-log-concavity.

\begin{lemma}\label{lem:st-log-conc-sst-unim}
  If $r=\{r_{i}\in \R_{>0}\}_{i \ldZ s}$ is strictly
  log-concave, then $r$ is almost strictly unimodal.  If
  $r=\{r_{i}\in \R_{<0}\}_{i \ldZ s}$ is strictly
  log-concave, then $-r=\{-r_{i}\in \R_{>0}\}_{i\ldZ s}$ is
  almost strictly unimodal.
\end{lemma}
\begin{proof}
  Let us prove the first statement.  Assume $j-1,j+1\ldZ s$.
  Then, the strict log-concavity of $r$ gives
  $r_{j}^{2}-r_{j-1}r_{j+1}>0$.  Then, $r_{j-1},r_{j}>0$
  implies $\frac{r_{j}}{r_{j-1}}>\frac{r_{j+1}}{r_{j}}$.
  Hence, if $\frac{r_{j}}{r_{j-1}}\leq 1$, then
  $r_{j-1}\geq r_{j}>r_{j+1}>\dots$.  Therefore, $r$ is almost
  strictly unimodal.  The second statement holds similarly.
\end{proof}

By the half $\succ$-$\gAfX$ implication, the
$\scc$-merged-log-concavity gives the
$\gAfX$-merged-log-concavity. Hence, we obtain the
following  short almost strictly unimodal
sequences by the $\gAfX$-merged-log-concavity.

\begin{lemma}\label{lem:short-unim-seq}
  Consider a length-zero fitting path
  $P=\{\tp{s,l,m_{i},n_{i},k_{i}}\}_{i\ldZ \theta}$ and
  $r\in \AfX$.  Let
  $\cF=\Lam(s,l,w,\scc,\fs,\phi,\rho,x,\fX)$ be
  $\gAfX$-merged-log-concave and static at $r$.  Let
  $\zeta=\tp{P,\cF}$. Then, we have the following.
  \begin{enumerate}
  \item \label{c:short-unim-seq-st-lc}
    $u(\zeta,r)$ is strictly log-concave.
  \item  \label{c:short-unim-seq-almost-st-unim}
    $u(\zeta,r)$ or
    $-u(\zeta,r)$ is positive and
    almost strictly unimodal, if $P$ is wrapped.
  \end{enumerate}
\end{lemma}
\begin{proof}
  Let us prove Claim~\ref{c:short-unim-seq-st-lc}.  Suppose
  $\theta=\tp{1,1}$ without loss of generality.  Then, 
  $e(\theta)=\tp{0,2}$. Thus, we want to prove the strict
  log-concavity of the following sequence:
  \begin{align}
    u(\zeta,r)_{0}&=\cF_{\mn\bom \kn}(r);
                 \label{eq:short-unim-seq-first-pp-term}
    \\
    u(\zeta,r)_{1}&=\cF_{\nn}(r);
                 \label{eq:short-unim-seq-second-pp-term}\\
    u(\zeta,r)_{2}&=\cF_{\nn\bop \kn}(r).
                 \label{eq:short-unim-seq-third-pp-term}
  \end{align}

  By the $\gAfX$-merged-log-concavity of $\cF$, we have
  $\Delta(\cF)(s,l,w,\mn,\nn,\kn,\phi,\rho,x,\fX) =
  \Ups(s,l,w,\mn,\nn,\kn,\phi,\rho,x,\fX) (\cF_{\mn}
  \cF_{\nn^{\ve}} - \cF_{\mn\bom \kn}
  \cF_{(\nn\bop    \kn)^{\ve}}) \gAfX 0$.  Thus,
  Lemma~\ref{lem:gAfX-closedness} and
  Corollary~\ref{cor:rshift-factor-pos} yield
  $\cF_{\mn} \cF_{\nn^{\ve}} -\cF_{\mn\bom \kn}
  \cF_{(\nn\bop \kn)^{\ve}}\gAfX 0$.  Also, $\mn=\nn^{\ve}$
  or $\mn=\nn$.  Hence,
  $\cF_{\nn}(r) \cF_{\nn}(r) -\cF_{\mn\bom \kn}(r)
  \cF_{\nn\bop \kn}(r)\gAfX 0$, since $\cF$ is static at
  $r$.  Then, Claim~\ref{c:short-unim-seq-st-lc} follows,
  since equations~\eqref{eq:short-unim-seq-first-pp-term},
  \eqref{eq:short-unim-seq-second-pp-term},
  and~\eqref{eq:short-unim-seq-third-pp-term} give
  $u(\zeta,r)_{1}^{2}-u(\zeta,r)_{0}u(\zeta,r)_{2}>0$.

  Claim~\ref{c:short-unim-seq-almost-st-unim} holds by
  Claim~\ref{c:total-pos-neg-pos-or-neg} of
  Lemma~\ref{lem:total-pos-neg}
  and Lemma~\ref{lem:st-log-conc-sst-unim},
  because
   $\mn\bom \kn\ldZl s$ and $\nn\bop \kn\ldZl s$ imply
   $u(\zeta,r)_{0}u(\zeta,r)_{2}\neq 0$ by
  equations~\eqref{eq:short-unim-seq-first-pp-term}
  and~\eqref{eq:short-unim-seq-third-pp-term}.
\end{proof}

Moreover, we obtain the following almost strictly unimodal sequences
by general fitting paths.

\begin{theorem}\label{thm:path-parcel-seq}
  Consider a fitting path
  $P=\{P_{i}=\tp{s,l,m_{i},n_{i},k_{i}}\}_{i\ldZ \theta}$.
  Assume that $\cF=\Lam(s,l,w,\scc,\fs,\phi,\rho,x,\fX)$ is
  $\gAfX$-merged-log-concave and static at some $r\in \AfX$.
  Let $\zeta=\tp{P,\cF}$.  Then, we have the following.
  \begin{enumerate}
  \item \label{c:path-parcel-seq-lc} $u(\zeta,r)$ is
    strictly log-concave.
  \item \label{c:path-parcel-seq-unim} $u(\zeta,r)$ or
    $-u(\zeta,r)$ is positive and almost strictly unimodal, if $P$ is
    wrapped.
  \end{enumerate}
\end{theorem}
\begin{proof}
  Let us prove Claim~\ref{c:path-parcel-seq-lc}.  Let
  $\lam\ldZ \theta$.  Consider the length-zero fitting path
  $Q_{\lam}=\{P_{i}\}_{i\ldZ \tp{\lam,\lam}}$. For
  $\ka_{\lam}=\tp{Q_{\lam},\cF}$, let us prove
  the following equations:
  \begin{align}
    u(\ka_{\lam},r)_{\lam-1}&=u(\zeta,r)_{\lam-1};
      \label{eq:path-parcel-seq-left}\\
    u(\ka_{\lam},r)_{\lam}
    &= u(\zeta,r)_{\lam};
      \label{eq:path-parcel-seq-middle}\\
    u(\ka_{\lam},r)_{\lam+1}
    &= u(\zeta,r)_{\lam+1}.
      \label{eq:path-parcel-seq-right}
  \end{align}
  In particular, these equations imply
  Claim~\ref{c:path-parcel-seq-lc} by
  Claim~\ref{c:short-unim-seq-st-lc} of
  Lemma~\ref{lem:short-unim-seq} on $u(\ka_{\lam},r)$.
  
  First, let us prove
  equation~\eqref{eq:path-parcel-seq-left}.  If
  $\lam-1\nldZ\theta$, then
  $u(\ka_{\lam},r)_{\lam-1}=\cF_{m_{\lam}\bom
    k_{\lam}}(r)=u(\zeta,r)_{\lam-1}$.  If
  $\lam-1\ldZ \theta$, then
  Proposition~\ref{prop:static-parcel-one-minus-one} gives
  $u(\ka_{\lam},r)_{\lam-1} =\cF_{m_{\lam}\bom k_{\lam}}(r)
  =\cF_{n_{\lam-1}}(r) =u(\zeta,r)_{\lam-1}$.

  Second, equation~\eqref{eq:path-parcel-seq-middle} holds,
  since $\lam\ld \theta$ implies
  $u(\ka_{\lam},r)_{\lam} =\cF_{n_{\lam}}(r)
  =u(\zeta,r)_{\lam}$.
  
  Third, let us prove
  equation~\eqref{eq:path-parcel-seq-right}.  If
  $\lam+1\nldZ\theta$, then
  equation~\eqref{eq:path-parcel-seq-right} follows from
  $u(\ka_{\lam},r)_{\lam+1}=\cF_{n_{\lam}\bop
    k_{\lam}}(r)=u(\zeta,r)_{\lam+1}$.  If
  $\lam+1\ldZ\theta$, then
  $m_{\lam+1}=(n_{\lam}\bop k_{\lam})^{\ve}$, which is
  $n_{\lam+1}$ or $n_{\lam+1}^{\ve}$. Hence,
  $u(\ka_{\lam},r)_{\lam+1} = \cF_{n_{\lam}\bop k_{\lam}}(r)
  =\cF_{n_{\lam+1}}(r) =u(\zeta,r)_{\lam+1}$ by the staticity
  of $\cF$.
  
  Let us prove Claim~\ref{c:path-parcel-seq-unim}.  Since
  $P$ is wrapped, we have
  $u(\zeta,r)_{\thetao-1}= \cF_{m_{\thetao}\bom
    k_{\thetao}}(r)\neq 0$ and
  $u(\zeta,r)_{\thetat+1}= \cF_{n_{\thetat}\bop
    k_{\thetat}}(r)\neq 0 \mif \thetat<\infty$.  Therefore,
  $u(\zeta,r)$ or $-u(\zeta,r)$ is almost strictly unimodal by
  Claim~\ref{c:total-pos-neg-pos-or-neg} of
  Lemma~\ref{lem:total-pos-neg}
  and Lemma~\ref{lem:st-log-conc-sst-unim}.
\end{proof}

\subsection{Infinite-length fitting paths}
\label{sec:infinite-paths}

By the following, we obtain infinite-length
fitting paths, which give infinite-length almost strictly
unimodal sequences by Theorem~\ref{thm:path-parcel-seq}.

\begin{definition}\label{def:inf-paths-sem-yd}
  Let $l\in \Zgeo$.  Let $\lam\in \Zgez^{3}$ such that
  $l\geq \lamo\geq 1$ and $\lamt\geq 1$.  Suppose infinite gates
  $s\geq 0$ and $\theta=\tp{1,\infty}$.  Let
  $t=\{t_{i}\in \Zgez\}_{i\ldZ \theta}$.
  
  \begin{enumerate}
  \item\label{d:inf-paths-sem-yd-flat-mn}
    If $\lamo=l$, then
    we define the sequence
    $P_{s,l,\lam,t}=\{P_{s,l,\lam,t,i}
    =\tp{s,l,m_{i},n_{i},k_{i}}\}_{i\ldZ
      \theta}$ such that
    \begin{align}
      m_{i}=n_{i}
      &=  \iota^{l}(\lamt+\lamr+s_{1})+(i-1)
        \iota^{l}(\lamt),\\
      k_{i}
      &=  \iota^{l}(0)\ccn \tp{\lamt}\ccn\iota^{l-1}(0)
        +\tp{t_{i}}\ccn \iota^{2l-1}(0).
    \end{align}
  \item\label{d:inf-paths-sem-yd-non-flat-mn} If $\lamo< l$,
    then we define the sequence
    $P_{s,l,\lam,t}=\{P_{s,l,\lam,t,i}
    =\tp{s,l,m_{i},n_{i},k_{i}}\}_{i\ldZ \theta}$ such that
    \begin{align}
      m_{i}=n_{i}
      &=
        \iota^{\lamo}(\lamt)\ccn
        \iota^{l-\lamo}(0)+
        \iota^{l}(\lamr+s_{1})+
        (i-1) \iota^{l}(2\lamt),\\
      k_{i}
      &=
        \iota^{\lamo}(0)\ccn\tp{\lamt} \ccn
        \iota^{l-\lamo-1}(0)
        \ccn\tp{\lamt}\ccn
        \iota^{\lamo-1}(0)\ccn\tp{\lamt} \ccn \iota^{l-\lamo-1}(0)
      \\& +\tp{t_{i}}\ccn \iota^{2l-1}(0).
    \end{align}
  \end{enumerate}
  When each $t_{i}=0$, let $P_{s,l,\lam}=P_{s,l,\lam,t}$
  for simplicity.
\end{definition}
For Item~\ref{d:inf-paths-sem-yd-flat-mn} of
Definition~\ref{def:inf-paths-sem-yd}, we have
the following.
\begin{proposition}\label{prop:inf-paths-sem-yd-flat-mn}
  Assume $l\in \Zgeo$ and $\lam\in \Zr$ such that
  \begin{align}
    \lamo=l, \  \lamt\geq 1, \ \lamr\geq 0.
    \label{cond:inf-paths-sem-yd-flat-mn-triple}
  \end{align}
  Then,
  $P_{s,l,\lam,t}=\{P_{s,l,\lam,t,i}=\tp{s,l,m_{i},n_{i},k_{i}}\}_{i\ldZ
    \theta}$ is an infinite-length fitting path of type zero
  such that
  \begin{align}
    \nu(k_{i})&=\iota^{l}(0)\ccn \iota^{l}(\lamt)+t_{i}.
               \label{eq:inf-paths-sem-yd-flat-mn-a}
  \end{align}
\end{proposition}
\begin{proof}
  For simplicity, let $t_{i}=0$.  First, we prove that each
  $P_{s,l,\lam,i}$ is fitting.  We have
  $m_{1},n_{1}\geq s_{1}$ in
  Item~\ref{d:inf-paths-sem-yd-flat-mn} of
  Definition~\ref{def:inf-paths-sem-yd} by
  condition~\eqref{cond:inf-paths-sem-yd-flat-mn-triple}.
  Thus, $m_{i},n_{i}\ldZl s$, since $s$ is an infinite gate.
  Moreover, since
  $k_{i}= \iota^{l}(0)\ccn \tp{\lamt}\ccn\iota^{l-1}(0)$, we
  obtain equation~\eqref{eq:inf-paths-sem-yd-flat-mn-a}.
  This gives
  $\nu(m_{i},n_{i},k_{i}) =\iota^{l}(\lamt+\lamr+s_{1})\ccn
  \iota^{l}(2\lamt+\lamr+s_{1})+(i-1) \iota^{2l}(\lamt)$.  Hence,
  each $P_{s,l,\lam,i}$ is fitting by
  condition~\eqref{cond:inf-paths-sem-yd-flat-mn-triple}.
  
  Second, we prove that
  $Q_{i}=\{Q_{i,\ka}=P_{s,l,\lam,i+\ka}\}_{\ka\ldZ
    \tp{0,1}}$ is a fitting path of type zero for each
  $i\ldZ \theta$.  Since $\sig(k_{i})=\iota^{l}(\lamt)$, we have
  $ (n_{i}\bop k_{i})^{\ve} = \iota^{l}(\lamt+\lamr+s_{1})+i
  \iota^{l}(\lamt)=m_{i+1}$.  Hence, each $Q_{i}$ is a fitting
  path of type zero by $m_{i+\ka}=n_{i+\ka}=m_{i+\ka}^{\ve}$
  for $\ka\ldZ \tp{0,1}$.
\end{proof}

\begin{example}\label{ex:inf-paths-sem-yd-flat-mn}
  Let $s=\tp{0,\infty}$, $l\in \Zgeo$, and
  $\lam=\tp{l,1,0}$.  Then, $\theta=\tp{1,\infty}$ provides the
  infinite-length fitting path
  $P_{s,l,\lam}=\{P_{s,l,\lam,i}=\tp{s,l,m_{i},n_{i},k_{i}}\}_{i\ldZ
    \theta}$ of type-zero such that
  $m_{i}=n_{i}=\iota^{l}(i)$ and
  $k_{i}=\iota^{l}(0)\ccn\tp{1} \ccn \iota^{l-1}(0)$ for
  $i \ldZ \theta$.
\end{example}
For Item~\ref{d:inf-paths-sem-yd-non-flat-mn} of
Definition~\ref{def:inf-paths-sem-yd}, we have the
following.
\begin{proposition}\label{prop:inf-paths-sem-yd-non-flat-mn}
  Suppose $l\in \Zgeo$ and $\lam\in \Zr$ such that
  \begin{align}
    1\leq \lamo< l, \  \lamt\geq 1, \ \lamr\geq 0.
    \label{cond:inf-paths-sem-yd-non-flat-mn-triple}
  \end{align}
  Consider
  $P_{s,l,\lam,t}=\{P_{s,l,\lam,t,i}
  =\tp{s,l,m_{i},n_{i},k_{i}}\}_{i\ldZ \theta}$.  Let
  $L(l,\lamo)=\min(\lamo-1,l-\lamo-1)$ and
  $H(l,\lamo)=\max(\lamo-1,l-\lamo-1)$.  Then, the following
  statements hold.
  \begin{enumerate}
  \item \label{c:inf-paths-sem-yd-non-flat-mn-sig} We have
    \begin{align}
      \sig(k_{i})_{l-j}=
      \begin{dcases}
        \lamt &\mif  j\in\oi(0,L(l,\lamo)),\\
        2\lamt &\mif j\in\oi(L(l,\lamo)+1,H(l,\lamo)),\\
        3\lamt&\mif  j\in\oi(H(l,\lamo)+1,l-1).
      \end{dcases}
    \end{align}
  \item
    \label{c:inf-paths-sem-yd-non-flat-mn-type-1-1-a}
    $P_{s,l,\lam}$ is an infinite-length fitting path of
    type 1-1 such that each $i\ldZ \theta$ satisfies
    \begin{align}
      \nu(k_{i})
      &=\iota^{\lamo}(0)\ccn \iota^{l-\lamo}(\lamt)
        \ccn\iota^{\lamo} (2\lamt)\ccn \iota^{l-\lamo}(3\lamt)+t_{i}.
        \label{eq:inf-paths-sem-yd-non-flat-mn-a}
    \end{align}
  \end{enumerate}
\end{proposition}
\begin{proof}
  For simplicity, assume each $t_{i}=0$.  Let us prove
  Claim~\ref{c:inf-paths-sem-yd-non-flat-mn-sig}.  If
  $\lamo-1>0$ and $l-\lamo-1>0$, then the statement is
  clear. Hence, we suppose otherwise.  First, let
  $\lamo-1=0$. If $l-\lamo-1=0$, then
  $k_{i}=\iota^{\lamo}(0)\ccn\tp{\lamt}
  \ccn\tp{\lamt}\ccn\tp{\lamt}$.  Hence,
  Claim~\ref{c:inf-paths-sem-yd-non-flat-mn-sig} holds,
  since $L(l,\lamo)=H(l,\lamo)=0$.  If $l-\lamo-1>0$, then
  $k_{i}=\iota^{\lamo}(0)\ccn\tp{\lamt} \ccn
  \iota^{l-\lamo-1}(0)\ccn\tp{\lamt} \ccn\tp{\lamt}\ccn
  \iota^{l-\lamo-1}(0)$ gives
  Claim~\ref{c:inf-paths-sem-yd-non-flat-mn-sig} by
  $L(l,\lamo)=0<H(l,\lamo)=l-\lamo-1$.  Second,
  let $\lamo-1>0$ and
  $l-\lamo-1=0$.  Then,
  $k_{i}=\iota^{\lamo}(0)\ccn\tp{\lamt}\ccn\tp{\lamt}\ccn
  \iota^{\lamo-1}(0)\ccn\tp{\lamt}$ implies
  Claim~\ref{c:inf-paths-sem-yd-non-flat-mn-sig} by
  $L(l,\lamo)=0<H(l,\lamo)=\lamo-1$.
  
  Let us prove
  Claim~\ref{c:inf-paths-sem-yd-non-flat-mn-type-1-1-a}.
  Thus, we prove that each
  $P_{s,l,\lam,i}=\tp{s,l,m_{i},n_{i},k_{i}}$ is fitting.
  Since $s$ is an infinite gate, $m_{i},n_{i}\ldZl s$.
  Also,
we obtain
  equation~\eqref{eq:inf-paths-sem-yd-non-flat-mn-a},
  ignoring $\ccn \iota^{\lamo-1}(0)$ or
  $\ccn \iota^{l-\lamo-1}(0)$ in $k_{i}$ when $\lamo=1$ or
  $\lamo=l-1$.  Then,
  since
  $\nu(m_{i}, n_{i},k_{i})
    =\iota^{l}(\lamt+\lamr+s_{1})\ccn
      \iota^{l}(3\lamt+\lamr+s_{1})
      +(i-1) \iota^{2l}(2\lamt)$,
       $P_{s,l,\lam,i}$ is fitting for each
  $i\ldZ \theta$ by
  condition~\eqref{cond:inf-paths-sem-yd-non-flat-mn-triple}.
  
  Let us  prove that
  $Q_{i}=\{Q_{i,\ka}=P_{s,l,\lam,i+\ka}\}_{\ka\ldZ \tp{0,1}}$ is a
  fitting path of type 1-1 for each $i\ldZ \theta$.
  First, assume $L(l,\lamo)=l-\lamo-1$.  Then,
  since $\lamo-1-(l-\lamo-1)=2\lamo-l\geq 0$,
  Claim~\ref{c:inf-paths-sem-yd-non-flat-mn-sig}
  gives
  $\sig(k_{i})
    =
      \iota^{l-\lamo}(3\lamt)
      \ccn
      \iota^{2\lamo-l}(2\lamt)
      \ccn
      \iota^{l-\lamo} (\lamt)$.
      Moreover, $l-\lamo+2\lamo-l=\lamo$ implies
  \begin{align}
    n_{i}\bop k_{i}
    =\iota^{l-\lamo}(2\lamt)\ccn \iota^{\lamo}(3\lamt)
    +\iota^{l}(\lamr+s_{1})
    +(i-1) \iota^{l}(2\lamt).
  \end{align}
  In particular, $(n_{i}\bop k_{i})^{\ve} = m_{i+1}$.
  Thus,
  $Q_{i}$ is a fitting path of type 1-1 for each $i\ldZ \theta$.
  
  Second, assume $L(l,\lamo)=\lamo-1$.
  Then,  since $l-\lam-1-(\lamo-1)=l-2\lamo\geq 0$,
  Claim~\ref{c:inf-paths-sem-yd-non-flat-mn-sig} yields
  $\sig(k_{i})
    =
      \iota^{\lamo}(3\lamt)
      \ccn
      \iota^{l-2\lamo}(2\lamt)
      \ccn
      \iota^{\lamo} (\lamt)$.
  Furthermore,  $\lamo+l-2\lamo=l-\lamo$ 
  gives
  \begin{dmath*}
    n_{i}\bop k_{i}
    =\iota^{l-\lamo}(2\lamt)\ccn \iota^{\lamo}(3\lamt)
    +\iota^{l}(\lamr+s_{1})
    +(i-1) \iota^{l}(2\lamt).
  \end{dmath*}
  Therefore, since
  $(n_{i}\bop k_{i})^{\ve}
  = m_{i+1}$, $Q_{i}$ is a fitting path of type 1-1 for each
  $i\ldZ \theta$.
\end{proof}

\begin{example}
  Let $s=\tp{0,\infty}$, $l=3$, and $\lam=\tp{1,1,0}$.  Then,
$\theta=\tp{1,\infty}$  gives
  the infinite-length
  fitting path
  $P_{s,l,\lam}=\{P_{s,l,\lam,i}
  =\tp{s,l,m_{i},n_{i},k_{i}}\}_{i\ldZ
    \theta}$ of   type 1-1 such that
  $m_{i}=n_{i}=\tp{1,0,0}+2i-2$ and
  $k_{i}=\tp{0,1,0,1,1,0}$ for 
  each $i\ldZ \theta$.
\end{example}

\subsection{Triplet scalings and sums of fitting paths}
\label{sec:scalings-sums-paths}

To construct more fitting
paths,
we first introduce the following scalings and sums.

\begin{definition}\label{def:fitting-scaling-sum}
  Suppose a gate $s$ and $l\in \Zgeo$.
  \begin{enumerate}
  \item Assume a tuple
    $T=\tp{s,l,\alp,\bta,\gam}$ such that
    $\alp,\bta\in \Zl$ and $\gam\in \Ztl$.  Let
    $\lam\in \Z$.  Then, we put the triplet scaling
    $\lam \odot T =\tp{s,l,\lam \alp,\lam \bta,\lam \gam}$.
  \item For each $i\in\oi(2)$, assume a tuple
    $P_{i}=\tp{s,l,m_{i},n_{i},k_{i}}$ such that
    $m_{i},n_{i}\in \Zl$ and $k_{i}\in \Ztl$.  Then, we put
    the triplet sum
    $P_{1}\dotplus P_{2} =\tp{s,l, \mn+\mt, \nn+\nt, \kn+\kt}$.
  \end{enumerate}
\end{definition}
In particular, we obtain
the following fitting tuples.

\begin{lemma}
  \label{lem:fitting-tuples-triplet-scaling-sums}
  Consider fitting tuples
  $P_{j}=\tp{s,l,m_{j},n_{j},k_{j}}$
  for $j\in\oi(2)$.  Suppose $\lam\in \Zt$ such that
  \begin{align}
    \lamo,\lamt\geq 0 \mand \lamo+\lamt >0.
    \label{ineq:fitting-tuples-triplet-scalings-sums}
  \end{align}
  Let
  $P_{3}= \lamo\odot P_{1}\dotplus \lamt \odot P_{2}
  =\tp{s,l,\mr,\nr,\kr}$.  Then, $P_{3}$ is fitting,
  provided $\mr,\nr\ldZl s$.
\end{lemma}
\begin{proof}
  Let $a_{i}=\nu(k_{i})$ and $b_{i}=\nu(m_{i},n_{i},k_{i})$
  for  $i\in\oi(3)$.  Then,   we have
  \begin{align}
    a_{3}&=\lamo \nu(\kn)+\lamt \nu(\kt)\\
         &=\lamo a_{1}+\lamt a_{2},\\
    b_{3}&=a_{3}+\mr\ccn \nr\\
         &=a_{3}+(\lamo \mn+\lamt \mt)\ccn(\lamo \nn+\lamt \nt)\\
         &=\lamo a_{1}+\lamt a_{2}+\lamo \mn \ccn
           \lamo \nn+\lamt \mt\ccn \lamt \nt\\
         &=\lamo  b_{1}+\lamt b_{2}. 
  \end{align}
  Thus, the slope conditions of $P_{1}$ and
  $P_{2}$ imply those of $P_{3}$ by
  inequalities~\eqref{ineq:fitting-tuples-triplet-scalings-sums}.
  This gives the statement by $\mr,\nr\ldZl s$.
\end{proof}

Moreover, we introduce the following notion
on tuple sequences by the scalings and sums above.
\begin{definition}\label{def:fit-sum}
  Suppose gates $s$ and $\theta$.  Let $l\in \Zgeo$.
  \begin{enumerate}
  \item Consider a sequence
    $T=\{T_{i} =\tp{s,l,\alp_{i},\bta_{i},\gam_{i}}\}_{i\ldZ
      \theta}$ with $\alp_{i},\bta_{i}\in \Zl$ and
    $\gam_{i}\in \Ztl$.  Let $\lam\in \Z$.  Then, we define
    the triplet scaling
    $\lam\odot T =\{ \lam\odot T_{i} \}_{i\ldZ \theta}$.
  \item For each $j\in\oi(2)$, consider a sequence
    $P_{j}=\{P_{j,i}
    =\tp{s,l,m_{j,i},n_{j,i},k_{j,i}}\}_{i\ldZ \theta}$ with
    $m_{j,i},n_{j,i}\in \Zl$ and $k_{j,i}\in \Ztl$.  Then, we
    define the triplet sum
    $ P_{1}\dotplus P_{2} =\{ P_{1,i}\dotplus P_{2,i}\}_{i\ldZ \theta}$.
  \end{enumerate}
\end{definition}

For example, suppose an infinite gate $s\geq 0$ and
$\lam=\tp{\lamo,\lamt,0}\in \Zr$ such that $s_{1}=0$ and
$\lamo,\lamt\in \Zgeo$. This gives
$P_{s,l,\lam}=\lamt\odot P_{s,l,\tp{\lamo,1,0}}$.

We now prove the following strict inequality on the ladders
of fitting paths.
\begin{lemma}\label{lem:fitting-path-mn-positivity}
  Assume a fitting path
  $P=\{\tp{s,l,m_{i},n_{i},k_{i}}\}_{ i\ldZ \tp{1,
      2}}$. Then, we have $\nn=\mn< \mt^{\ve} $ or
  $\nn=\mn^{\ve} <\mt^{\ve} $.  In particular, we have
  $\sum \mn=\sum \nn<\sum\mt=\sum \nt$.
\end{lemma}
\begin{proof}
  Since $P$ is a fitting path, $\nn=\mn$ or $\nn=\mn^{\ve}$.
  If $\nn=\mn$, then $\sig(k)>0$ gives
  $\mn<\mn\bop \kn = \nn\bop \kn =\mt^{\ve}$.  If
  $\nn=\mn^{\ve}$, then
  $\mn^{\ve}< \mn^{\ve}+\sig(\kn)^{\ve} =\nn+\sig(\kn)^{\ve}
  =\nn\bop \kn =m_{2}^{\ve}$.  The latter statement follows
  from $\mt=\nt$ or $\mt^{\ve}=\nt$.
\end{proof}

Then, we obtain fitting paths by triplet
scalings and sums.

\begin{proposition}\label{prop:fitting-paths-scalings-sums}
  Assume $\lam\in \Ztgeo$ such that
  \begin{align}
    \lamo,\lamt\geq 0 \mand \lamo+\lamt >0.
    \label{ineq:fitting-paths-scalings-sums}
  \end{align}
  Let $\theta=\tp{1,2}$. Consider fitting paths
  $P_{j}=\{P_{j,i}
  =\tp{s,l,m_{j,i},n_{j,i},k_{j,i}}\}_{ i\ldZ
    \theta}$ for $j\in\oi(2)$ such that $P_{1}$ and
  $P_{2}$ have the same type $A$.  Also, suppose
  $P_{3}= \lamo\odot P_{1}\dotplus \lamt\odot P_{2}
           =\{P_{3,i}
           =\tp{s,l,m_{3,i},n_{3,i},k_{3,i}}\}_{ i\ldZ \theta}$
           with $m_{3,2},n_{3,2}\ldZl s$.
Then, $P_{3}$ is a fitting
  path of the type $A$.
\end{proposition}
\begin{proof}
  Lemma~\ref{lem:fitting-tuples-triplet-scaling-sums}
  implies that
   $P_{3,2}$ is fitting.
  Hence, let us prove that $P_{3,1}$ is fitting.  Thus, we
  confirm $m_{3,1}, n_{3,1} \ldZl s$.  By
  inequalities~\eqref{ineq:fitting-paths-scalings-sums}, we  
  have
  \begin{align}
    m_{3,1}&\geq m_{1,1} \mor m_{2,1},
             \label{ineq:fitting-paths-scalings-sums-m-lower-bound}\\
    n_{3,1}&\geq n_{1,1} \mor n_{2,1}.
             \label{ineq:fitting-paths-scalings-sums-n-lower-bound}
  \end{align}
  Also, since $P_{1}$ and $P_{2}$ are fitting paths,
  Lemma~\ref{lem:fitting-path-mn-positivity} yields
  \begin{align}
    m_{3,1}&=\lamo m_{1,1}+\lamt m_{2,1}
             <
             \begin{dcases}
               \lamo m_{1,2}+\lamt m_{2,2}=m_{3,2} \mif
               \mbox{A is type 2-1 or 2-2},\\
               \lamo
               m^{\ve}_{1,2}+\lamt m^{\ve}_{2,2}=m_{3,2}^{\ve}
               \mif \mbox{A is type 1-1 or 1-2}.
             \end{dcases}
  \end{align}
  Thus, $m_{3,1}\ldZl s$ by
  $m_{3,2}\ldZl s$ and
  inequality~\eqref{ineq:fitting-paths-scalings-sums-m-lower-bound}.
  Also, $n_{3,1}\ldZl s$ by
  $m_{3,2}\ldZl s$ and
  inequality~\eqref{ineq:fitting-paths-scalings-sums-n-lower-bound},
  because
  Lemma~\ref{lem:fitting-path-mn-positivity} gives
  $n_{3,1}=\lamo n_{1,1}+\lamt n_{2,1}
    <\lamo m^{\ve}_{1,2}+\lamt m^{\ve}_{2,2}
    =m_{3,2}^{\ve}$.
  Hence, $P_{3,1}$ is fitting by
  Lemma~\ref{lem:fitting-tuples-triplet-scaling-sums}.
  
  Let us prove that $P_{3}$ is a fitting path of the
  type $A$.  We have
  \begin{align}
    m_{3,2}=(n_{3,1}\bop k_{3,1})^{\ve},
    \label{eq:fitting-paths-scalings-sums-mn}
  \end{align}
  since
  $n_{3,1}\bop k_{3,1} =(\lamo n_{1,1} +\lamt n_{2,1})\bop
  (\lamo k_{1,1}+\lamt k_{2,1}) =\lamo ( n_{1,1}\bop
  k_{1,1}) +\lamt (n_{2,1}\bop k_{2,1}) =\lamo
  m_{1,2}^{\ve}+\lamt m_{2,2}^{\ve} =m_{3,2}^{\ve}$.
  
  For instance, assume that $P_{1}$ and $P_{2}$ have the
  type 1-1.  Then, we have $m_{3,1}=n_{3,1}$ by
  $m_{1,1}=n_{1,1}$ and $m_{2,1}=n_{2,1}$, and
  $m_{3,2}=n_{3,2}$ by $m_{1,2}=n_{1,2}$ and
  $m_{2,2}=n_{2,2}$.  Thus, $P_{3}$ is of type 1-1 by
  equation~\eqref{eq:fitting-paths-scalings-sums-mn}. Similar
  discussions hold for the other types.
\end{proof}

\subsection{Infinite-length fitting paths
  of Young diagrams}
\label{sec:infinite-paths-young}
We prove that Young diagrams
 with box counting
give explicit
infinite-length fitting paths. This uses the following
notation.

\begin{definition}\label{def:young}
  Suppose $p,\lam\in \Zgeo$.
  \begin{enumerate}
  \item $\xi\in \Zgeo^{p}$ is called a Young diagram (or a
    partition of $\sum \xi$), if $\xi$ is decreasing.
  \item We write $Y(p,\lam)$ for the set of Young diagrams
    $\xi\in \Zgeo^{p}$ such that $\xi_{1}\leq \lam$.
  \end{enumerate}
\end{definition}

By Definition~\ref{def:families-general-operations},
we use the following notation for the box counting of
$Y(p,\lam)$.

\begin{definition}\label{def:young-box-counting}
  Let $\xi\in Y(p,l)$, $\lam\in \Z$, and $\lam'\in \Zgeo$.
  \begin{enumerate}
  \item Let
    $\len_{\leq \lam}(\xi)=\sum_{i\leq \lam}\len_{i}(\xi)\in \Zgez$.
  \item   Let 
    $\len(\lam,\xi)
    = 2\len_{\leq \lam-1}(\xi)+\len_{\lam}(\xi)\in \Zgez$.
  \item Let 
    $z(\lam',\lam,\xi)
    =\lam+\tp{\len_{\leq 0}(\xi),
      \len_{\leq 1}(\xi),\dots,\len_{\leq \lam'-1}
    (\xi)}\in \Z^{\lam'}$.
  \end{enumerate}
\end{definition}

Moreover, we define the following sequence of 
fitting tuples.

\begin{definition}\label{def:young-fitting-paths}
  Let $h\in\Zgez$ and $\xi\in Y(p,l)$.  Consider
  $\lam_{i}\in \Zgez^{3}$ for $i\in\oi(p)$ such that
  $\lamo=\tp{\xi_{1},1,h}$ and $\lam_{i}=\tp{\xi_{i},1,0}$ for
  $i\in\oi(2,p)$.  Also, suppose an infinite gate $s\geq 0$.
  Let $l\in \Zgeo$, $\theta=\tp{1,\infty}$, and
  $t=\{t_{i}\in \Zgez\}_{i\ldZ \theta}$.  Then, we define the
  sequence 
  \begin{align}
    \bP_{s,l,\xi,t,h}  = P_{s,l,\lam_{1},t} \dotplus
    P_{s,l,\lam_{2}} \dotplus \dots  \dotplus P_{s,l,\lam_{p}}.    
  \end{align}
    We call $h$ base height of $\bP_{s,l,\xi,t,h}$.  If each
  $t_{i}=0$, then we also write $\bP_{s,l,\xi,h}$ for
  $\bP_{s,l,\xi,t,h}$.
\end{definition}

Then,
we realize the following infinite-length fitting paths
$\bP_{s,l,\xi,t,h}$.

\begin{theorem}\label{thm:fitting-path-young}
  Consider a Young diagram $\xi\in Y(p, l)$.  Then,
  $\bP_{s,l,\xi,t,h} =\{\tp{s,l,m_{i},n_{i},k_{i}}\}_{i\ldZ
    \theta}$ is an infinite-length fitting path of type 1-1.
  Moreover, each $i\ldZ \theta$ satisfies
  \begin{align}
    m_{i}=n_{i}
    &= p+s_{1}+h    -z(l,0,\xi)+
      (i-1)\iota^{l}(\len(l,\xi)),      
      \label{eq:fitting-path-young-mn}\\
    \nu(k_{i})
    &=z(l,0,\xi)\ccn
      z(l,\len(l,\xi),\xi)+t_{i}.
      \label{eq:fitting-path-young-k}
  \end{align}
\end{theorem}
\begin{proof}
  The former statement holds by
  Propositions~\ref{prop:inf-paths-sem-yd-flat-mn},
  ~\ref{prop:inf-paths-sem-yd-non-flat-mn},
  and~\ref{prop:fitting-paths-scalings-sums}.

  Let us prove
  equations~\eqref{eq:fitting-path-young-mn}
  and~\eqref{eq:fitting-path-young-k}.  Suppose
  $s_{1}=h=0$
  and each $t_{i}=0$ for simplicity.  Thus, let
  $\lam_{j}=\tp{\xi_{j},1,0}$ for $j\in\oi(p)$.  Also,
  consider the fitting path
  $P_{s,l,\lam_{j}}
  =\{\tp{s,l,m_{\lam_{j},i},n_{\lam_{j},i},k_{\lam_{j},i}}\}_{
    i \ldZ \theta}$ for $j\in\oi(p)$. Then,
  Propositions~\ref{prop:inf-paths-sem-yd-flat-mn}
  and~\ref{prop:inf-paths-sem-yd-non-flat-mn} give the
  following:
  \begin{align}
    m_{\lam_{j},i}=n_{\lam_{j},i}
    &=  
      \begin{dcases}
        \iota^{l}(1)+(i-1) \iota^{l}(1) \mif \xi_{j}=l,
        \\
        \iota^{\xi_{j}}(1)\ccn\iota^{l-\xi_{j}}(0)+
        (i-1) \iota^{l}(2) \mif \xi_{j}<l;
      \end{dcases}
      \label{eq:fitting-path-young-cases-for-mn}
    \\
    \nu(k_{\lam_{j},i})
    &=
      \begin{dcases}
        \iota^{l}(0)\ccn \iota^{l}(1) \mif \xi_{j}=l,\\
        \iota^{\xi_{j}}(0)\ccn \iota^{l-\xi_{j}}(1)
        \ccn\iota^{\xi_{j}}(2)\ccn \iota^{l-\xi_{j}}(3)  \mif \xi_{j}<l.
      \end{dcases}
    \label{eq:fitting-path-young-cases-for-a}
  \end{align}
  Thus, equation~\eqref{eq:fitting-path-young-mn}
  holds, since
  equation~\eqref{eq:fitting-path-young-cases-for-mn}
  implies
  $m_{1,\ka}=n_{1,\ka}= \sum_{j\in \oi(p)}m_{\lam_{j},1,\ka}=
  \sum_{u\in \oi(\ka,l)}\len_{u}(\xi)$
  for each $\ka\in\oi(l)$. 
  Also,
  equation~\eqref{eq:fitting-path-young-k} holds, as
  equation~\eqref{eq:fitting-path-young-cases-for-a}
  gives
  $\nu(k_{i}) =\iota^{l}(0)\ccn \iota^{l}(\len_{l}(\xi))
  +\sum_{\ka\in\oi(l-1)} \iota^{\ka}(0)\ccn
  \iota^{l-\ka}(\len_{\ka}(\xi)) \ccn\iota^{\ka}(2
  \len_{\ka}(\xi))\ccn \iota^{l-\ka}(3 \len_{\ka}(\xi))$.
\end{proof}

Theorem~\ref{thm:fitting-path-young} gives not only
infinitely many polynomials with positive integer
coefficients by Young diagrams and the merged-log-concavity,
but also infinite-length almost strictly unimodal sequences
by Theorem~\ref{thm:path-parcel-seq}.  Also, $z(l,0,\xi)$ and
$p-z(l,0,\xi)$ in Theorem~\ref{thm:fitting-path-young} count
boxes of Young diagrams in $p\times l$ rectangles as
in the following examples.

\begin{example}\label{ex:f-y-1}
  Let $l=4$, $p=2$, $\xi=\tp{3,1}$, $s=\tp{0,\infty}$, and
  $h=0$.  Then,
  there exists the infinite-length fitting
  path
  $\bP_{s,l,\xi,h} =\{\tp{s,l,m_{i},n_{i},k_{i}}\}_{i\ldZ
    \theta}$ such that each $i\ldZ \theta$ satisfies
  \begin{align}
    m_{i}=n_{i}&=\tp{2,1,1,0}+(i-1)\iota^{l}(4),
                 \label{ex:f-y-1-mn}
    \\
    \nu(k_{i})&=\tp{0,1,1,2,4,5,5,6}.
               \label{ex:f-y-1-k}
  \end{align}

   To explain these equations
   by $p  -z(l,0,\xi)$ and  $z(l,0,\xi)$,
   consider the following Young diagram of $\xi$:
  \begin{align}
    \begin{ytableau}
      ~ & \none & \none \\
      ~ & ~       & ~ &\none
    \end{ytableau}.
  \end{align}
  Then, 
  $p-z(l,0,\xi)$ corresponds to the following
  box counting:
  \begin{align}
    \begin{ytableau}
      ~ & \none & \none \\
      2 & 1       & 1 &\none[0]
    \end{ytableau},
  \end{align}
  where $2,1,1,0$ indicate the numbers of boxes in the
  vertical direction.  Hence, we obtain
  equation~\eqref{ex:f-y-1-mn} by
  Theorem~\ref{thm:fitting-path-young}, since
  $\len_{\leq l-1}(\xi)=2$ and $\len_{l}(\xi)=0$ gives
  $\len(l,\xi) =4$.
   
  Moreover, in the $p\times l$ rectangle,
  $\xi$ gives the following complementary Young diagram:
  \begin{align}
    \begin{ytableau}
      \none & ~ & ~ &~\\
      \none & \none  &\none &~
    \end{ytableau}.
  \end{align}
  Then, $z(l,0,\xi)$ corresponds to the
  following box counting:
  \begin{align}
    \begin{ytableau}
      \none[0] & 1 & 1 &2\\
      \none & \none  &\none &~
    \end{ytableau},
  \end{align}
  where $0,1,1,2$ are the numbers of boxes in the vertical direction.
  Thus, equation~\eqref{ex:f-y-1-k} follows.  
\end{example}
When $\xi_{1}=l$, we have the following example.
\begin{example}\label{ex:f-y-2}
  Let $l=3$, $p=2$, $\xi=\tp{3,1}$, $s=\tp{0,\infty}$, and
  $h=0$.  Then, there is the infinite-length fitting path
  $\bP_{s,l,\xi,h} =\{\tp{s,l,m_{i},n_{i},k_{i}}\}_{i\ldZ
    \theta}$ such that each $i\ldZ \theta$ satisfies
  \begin{align}
    m_{i}=n_{i}
    &=\tp{2,1,1}+(i-1)\iota^{l}(3),
      \label{ex:f-y-2-mn}\\
    \nu(k_{i})
    &=\tp{0,1,1,3,4,4}.
      \label{ex:f-y-2-k}
  \end{align}
  
  To explain these equations, we have the following Young
  diagram of $\xi$:
  \begin{align}
    \begin{ytableau}
      ~ & \none &\none\\
      ~ & ~       & ~ 
    \end{ytableau}.
  \end{align}
  Thus, $p-z(l,0,\xi)$ corresponds to the
  following box counting:
  \begin{align}
    \begin{ytableau}
      ~ & \none & \none \\
      2 & 1       & 1
    \end{ytableau},
  \end{align}
  where $2,1,1$ are the numbers of boxes in the
  vertical direction.  This gives
  equation~\eqref{ex:f-y-2-mn} by
  Theorem~\ref{thm:fitting-path-young}, since
  $\len_{\leq l-1}(\xi)=1$ and $\len_{l}(\xi)=1$.

  Furthermore, in the $p\times l$ rectangle, $\xi$ has the
  following complementary Young diagram:
  \begin{align}
    \begin{ytableau}
      \none & ~ & ~ \\
      \none & \none  &\none
    \end{ytableau}.
  \end{align}
  Then, $z(l,0,\xi)$ corresponds to the following box
  counting:
  \begin{align}
    \begin{ytableau}
      \none[0] & 1 & 1\\
      \none & \none  &\none 
    \end{ytableau},
  \end{align}
  where $0,1,1$ are the numbers of boxes in the
  vertical direction.  Hence, we obtain
  equation~\eqref{ex:f-y-2-k}.
\end{example}

\begin{remark}\label{rmk:Young-succ-log-conc}
  Let $l\in \Zgeo$, $p=1$, $\xi=\tp{l}$, $s=\tp{0,\infty}$, and
  $h=0$.  Then, there is the infinite-length fitting path
  $\bP_{s,l,\xi,h}
  =\{\tp{s,l,m_{i},n_{i},k_{i}}\}_{i\ldZ \theta}$
  such that $m_{i}=n_{i} =i \iota^{l}(1)$ and
  $\nu(k_{i}) =\iota^{l}(0)\ccn \iota^{l}(1)$ for each $i\ldZ \theta$.

  Moreover, suppose a parcel $\cF=\Lam(s,l,\scc,\fs,\fX)$
  and $i\ldZ \theta$.  Then, since $\sig(k)=\iota^{l}(1)$, the
  $\scc'$-log-concavity
  $\cF_{m_{i}}\cF_{m_{i}}-\cF_{m_{i}-1}\cF_{m_{i}+1}\scc' 0$
  coincides with the $\scc'$-merged-log-concavity
  $\cF_{m_{i}}\cF_{n_{i}^{\ve}}-\cF_{m_{i}\bom k_{i}}
    \cF_{(n_{i}\bop k_{i})^{\ve}}\scc' 0$.
 \end{remark}

 We now obtain infinite-length almost strictly
 unimodal sequences by the following lemma.

\begin{lemma}\label{lem:low-inclusion}
  Suppose
  $P=\bP_{s,l,\xi,t,h} =\{\tp{s,l,m_{i},n_{i},k_{i}}\}_{i\ldZ
    \theta}$.
  Then, $P$ is wrapped, when
  \begin{align}
    h\geq 3 \len_{\leq l-1}(\xi).
    \label{ineq:low-inclusion-3}
  \end{align}
\end{lemma}
\begin{proof}
  Since $s$ is infinite, $P$ has the upper inclusion
  condition.  Let us confirm that
  inequality~\eqref{ineq:low-inclusion-3} gives the lower
  inclusion condition $\mn\bom \kn\geq s_{1}$ of $P$.  Let
  $a_{1}=\nu(\kn)$ and $p=\len_{\leq l-1}(\xi)+\len_{l}(\xi)$. Then, 
  Theorem~\ref{thm:fitting-path-young} gives
  \begin{align}
    \sig(\kn)_{1}
    &=a_{1,2l}-t_{1}
      =\len(l,\xi)+\len_{\leq l-1}(\xi)
      =3 \len_{\leq l-1}(\xi)+\len_{l}(\xi)
      =2 \len_{\leq l-1}(\xi)+p.
  \end{align}
  Thus, since 
  $\mn\geq s_{1}+h+p-\len_{\leq l-1}(\xi)$
  by Theorem~\ref{thm:fitting-path-young},
   we obtain
  \begin{align}
    \mn\bom \kn
    &\geq s_{1}+h+p-\len_{\leq l-1}(\xi)- \sig(\kn)_{1}
      =s_{1}+h-3\len_{\leq l-1}(\xi).
  \end{align}
  In particular, inequality~\eqref{ineq:low-inclusion-3}
  gives the
  lower inclusion condition of $P$.
\end{proof}

\begin{proposition}\label{prop:unim-wrapped-three}
  Let
  $P=\bP_{s,l,\xi,h} =\{\tp{s,l,m_{i},n_{i},k_{i}}\}_{i\ldZ
    \theta}$.  Suppose that
  $\cF=\Lam(s,l,w,\scc,\fs,\phi,\rho,x,\fX)$ is
  $\scc$-merged-log-concave and static at some $r\in
  \AfX$. Let $\zeta=\tp{P,\cF}$. Then, we have the following.
  \begin{enumerate}
  \item \label{c:unim-wrapped-three-lc} $u(\zeta,r)$ is
    infinite-length and strictly log-concave.
  \item \label{c:prop-wrapped-three-lc-three} $u(\zeta,r)$ or
    $-u(\zeta,r)$ is infinite-length and almost strictly
    unimodal, if $h\geq 3 \len_{\leq l-1}(\xi)$.
  \end{enumerate}
\end{proposition}
\begin{proof}
  Claim~\ref{c:unim-wrapped-three-lc} is of
  Theorems~\ref{thm:path-parcel-seq}
  and~\ref{thm:fitting-path-young}.  Then,
  Lemma~\ref{lem:low-inclusion} gives
  Claim~\ref{c:prop-wrapped-three-lc-three}.
\end{proof}

\subsection{On the fitting condition
  and infinite-length fitting paths}
\label{sec:strictly-fitting}

The fitting condition in Definition~\ref{def:fitting} is a
key notion for the merged-log-concavity.  To better
understand the notion, we discuss the following variant that
requires strict inequalities on the slope conditions of
fitting tuples. We then prove that the variant has no
infinite-length fitting paths with widths greater than one.

\begin{definition}\label{def:strictly-fitting}
  Let $l\in \Zgeo$,  $m,n\in \Zl$, and
  $k \in \Ztl$. Suppose a gate $s\geq 0$.  Let
  $a=\nu(k)$ and $b=\nu(m,n,k)\in \Ztl$.  We
  call $\mu=\tp{s,l,m,n,k}$ strictly fitting, if $\mu$ satisfies
  the inclusion condition with the following inequalities:
  \begin{align}
    b_{1}<\dots< b_{l}<b_{l+1}<\dots< b_{2l};
    \label{ineq:strictly-fitting-s2}\\
    0< a_{1}<\dots< a_{l}<a_{l+1}<\dots < a_{2l}.
    \label{ineq:strictly-fitting-s3}
  \end{align}
  We refer to
  \eqref{ineq:strictly-fitting-s2}
  and~\eqref{ineq:strictly-fitting-s3} as
  the strict upper slope condition and 
  strict lower slope condition of $\mu$.
\end{definition}

For example, Example~\ref{ex:paths-equiv-tuples} gives a
fitting path of strictly fitting tuples.
Also, some results
 in
 Sections~\ref{sec:prelim},\ref{sec:explicit-merged-parcels},
 and~\ref{sec:hadam} hold by strictly fitting tuples.
 But,  we state the following lemmas to confirm
 that there are no infinite-length fitting paths of strictly fitting tuples.

\begin{lemma}\label{lem:pref-upp-bounds}
  Let $l\in \Zget$.
  Suppose a strictly fitting $\mu=\tp{s,l,m,n,k}$.
  \begin{enumerate}
  \item\label{c:pref-upp-bounds-sig-diff-minus-one}
    Then, $\nn-n_{l}< \sig(k)_{1}-\sig(k)_{l}-1$.
  \item \label{c:pref-upp-bounds-sig-diff-half}
    If $m=n$, then
    $\nn-n_{l}<\frac{\sig(k)_{1}-\sig(k)_{l}}{2}$.
  \end{enumerate}
\end{lemma}
\begin{proof}
  Let $a=\nu(k)$ and $b=\nu(m,n,k)$.  Let us prove
  Claim~\ref{c:pref-upp-bounds-sig-diff-minus-one}.  By
  $l\in \Zget$, the strict upper slope condition of $\mu$
  gives $b_{2l}=n_{l}+a_{2l}>b_{l+1}=\nn+a_{l+1}$. In
  particular, $\nn-n_{l}<a_{2l}-a_{l+1}$.  Moreover, since
  $a_{l}-a_{1}>0$ by the strict lower slope condition of
  $\mu$, $a_{2l}-a_{l+1} < a_{l}-a_{1}+ a_{2l}-a_{l+1}$.
  Therefore,
  Claim~\ref{c:pref-upp-bounds-sig-diff-minus-one} follows
  from
  \begin{align}
    a_{l}-a_{1}+         a_{2l}-a_{l+1}
    =\kt+\dots+k_{l}+k_{l+2}+\dots +k_{2l}
    =\sig(k)_{1}-\sig(k)_{l}.
    \label{eq:pref-upp-bounds-a-sig-diff}
  \end{align}
  
  Let us prove Claim~\ref{c:pref-upp-bounds-sig-diff-half}.
  The strict upper slope condition of $\mu$ implies
  $b_{1}=\nn+a_{1}<b_{l}=n_{l}+a_{l}$ and
  $b_{l+1}=\nn+a_{l+1}<b_{2l}=n_{l}+a_{2l}$.  This gives
  $\nn-n_{l}<a_{l}-a_{1}$ and $\nn-n_{l}<a_{2l}-a_{l+1}$.
  Then, we have $2(\nn-n_{l})<a_{l}-a_{1}+a_{2l}-a_{l+1}$.
  Thus, Claim~\ref{c:pref-upp-bounds-sig-diff-half}
  follows from equations~\eqref{eq:pref-upp-bounds-a-sig-diff}.
\end{proof}

\begin{lemma}\label{lem:pref-low-bounds}
  Let $l\in \Zget$.  Assume a fitting path
  $P=\{ P_{i}=\tp{s,l,m_{i}, n_{i},k_{i}}\}_{i\in\oi(2)}$
  such that each $P_{i}$ is strictly fitting.  Let
  $a_{i}=\nu(k_{i})$ and $b_{i}=\nu(m_{i},n_{i},k_{i})$ for
  $i\in\oi(2)$. Then, each $j\in\oi(l-1)$ satisfies
  \begin{align}
    n_{1,l-j}-n_{1,l-j+1}
    +\sig(\kn)_{j+1}-\sig(\kn)_{j}
    +a_{2,j+1}-a_{2,j}
    &>0,
      \label{ineq:pref-low-bounds-n-sig-a-diff}
    \\        
    n_{1,j}- n_{1,j+1}
    &>1.
      \label{ineq:pref-low-bounds-n-diff}
  \end{align}
\end{lemma}
\begin{proof}
  Let us prove
  inequality~\eqref{ineq:pref-low-bounds-n-sig-a-diff}.
  Since $P$ is a fitting path,
  $\mt=(\nn\bop \kn)^{\ve}=\nn^{\ve}+\sig(\kn)$.  Then,
  $m_{2,j}=n_{1,l-j+1}+\sig(\kn)_{j}$ and
  $m_{2,j+1}=n_{1,l-j}+\sig(\kn)_{j+1}$. Hence,
  \begin{align}
    b_{2,j+1}-b_{2,j}
    &= (m_{2,j+1}+a_{2,j+1})
      -(m_{2,j}+a_{2,j})\\
    &=
      (n_{1,l-j}+\sig(\kn)_{j+1}+a_{2,j+1})
      -(n_{1,l-j+1}+\sig(\kn)_{j}+a_{2,j})\\
    &=n_{1,l-j}-n_{1,l-j+1}
      +\sig(\kn)_{j+1}-\sig(\kn)_{j}
      +a_{2,j+1}-a_{2,j}.
  \end{align}
  This yields
  inequality~\eqref{ineq:pref-low-bounds-n-sig-a-diff},
  since $b_{2,j+1}-b_{2,j}>0$ by the strict
  upper slope  condition.
  
  Let us prove
  inequality~\eqref{ineq:pref-low-bounds-n-diff}.  Since
  $a_{2,2l-j+1}-a_{2,2l-j}>0$ by the strict
  lower slope
  condition, the equivalence of $\kn$ and $\kt$ gives
  \begin{align}
    0&=\sum \kt(j+2,2l-j)- \sum \kt(j+1,2l-j+1)
    +k_{2,j+1}+k_{2,2l-j+1}
    \\& = \sig(\kt)_{j+1}-\sig(\kt)_{j}
    +k_{2,j+1}+k_{2,2l-j+1}
    \\&=\sig(\kn)_{j+1}-\sig(\kn)_{j}
    +a_{2,j+1}-a_{2,j}+a_{2,2l-j+1}-a_{2,2l-j}
     \\&> \sig(\kn)_{j+1}-\sig(\kn)_{j}
    +a_{2,j+1}-a_{2,j}.
\end{align}
  Thus,
  inequality~\eqref{ineq:pref-low-bounds-n-diff}
  follows from
  inequality~\eqref{ineq:pref-low-bounds-n-sig-a-diff}, as
  $n_{1,l-j}-n_{1,l-j+1}
    >
      -(\sig(\kn)_{j+1}-\sig(\kn)_{j}
      +a_{2,j+1}-a_{2,j})>0$.
\end{proof}
Hence, we obtain the following finiteness of fitting paths of
strictly fitting tuples.
\begin{proposition}\label{prop:pre-fitting}
  Let $l\in \Zget$. Assume a fitting path
  $P=\{P_{i}=\tp{s,l,m_{i},n_{i},k_{i}}\}_{i\ldZ \theta}$ such
  that each $P_{i}$ is strictly fitting.  Then,
  $\thetat-\thetao\leq 2$.
\end{proposition}
\begin{proof}
  Let $a_{i}=\nu(k_{i})$ and $b_{i}=\nu(m_{i},n_{i},k_{i})$
  for $i\ldZ \theta$.  Suppose $\thetao=1$ for
  simplicity.  We prove the statement by
  contradiction. Thus, let $\thetat\geq 4$.

  Claim~\ref{c:pref-upp-bounds-sig-diff-minus-one} of
  Lemma~\ref{lem:pref-upp-bounds} gives
  $n_{i,1}-n_{i,l}\leq \sig(k_{i})_{1}-\sig(k_{i})_{l}-2$
  for  $i\ldZ \theta$.  Therefore, we have
  \begin{align}
    m_{2,1}-m_{2,l}
    &\geq -(\sig(\kn)_{1}-\sig(\kn)_{l}-2)
      +(\sig(\kn)_{1}-\sig(\kn)_{l})
      =2,
      \label{ineq:pre-fitting-m-diff}
  \end{align}
  because $\mt=(\nn\bop \kn)^{\ve}$ implies
  $m_{2,1}-m_{2,l}
    =(n_{1,l}+\sig(\kn)_{1})
      -(n_{1,1}+\sig(\kn)_{l})
      =-(n_{1,1}-n_{1,l})+(\sig(\kn)_{1}-\sig(\kn)_{l})$.

  Furthermore, by Lemma~\ref{lem:pref-low-bounds}, each
  $i\in\oi(\thetat-1)$ satisfies
  \begin{align}
    n_{i,1}-n_{i,l}\geq 2.
    \label{ineq:pre-fitting-n-diff}
  \end{align}
  Hence, inequality~\eqref{ineq:pre-fitting-m-diff} yields
  $\mt=\nt$ by $2\leq \thetat-1$, since either $\mt=\nt$ or
  $\mt^{\ve}=\nt$ for the fitting path $P$.
  Then, 
  $\mr=(\nt\bop \kt)^{\ve} =(\mt\bop \kt)^{\ve}$ gives
  \begin{align}
    m_{3,1}
    &=m_{2,l}+\sig(\kt)_{1}
      =(n_{1,1}+\sig(\kn)_{l})+\sig(\kn)_{1},\\
    m_{3,l}
    &=m_{2,1}+\sig(\kt)_{l}
      =(n_{1,l}+\sig(\kn)_{1})+\sig(\kn)_{l}.
  \end{align}
  Then, inequality~\eqref{ineq:pre-fitting-n-diff} implies
  $m_{3,1}-m_{3,l}\geq 2$, and hence
  \begin{align}
    \mr=\nr.
    \label{eq:pre-fitting-m-n}
  \end{align}
  
  However, Claim~\ref{c:pref-upp-bounds-sig-diff-half} of
  Lemma~\ref{lem:pref-upp-bounds} yields
  \begin{align}
    m_{3,1}-m_{3,l}
    &=(n_{2,l}+\sig(\kt)_{1})
      -(n_{2,1}+\sig(\kt)_{l})\\
    &=-(n_{2,1}-n_{2,l})+(\sig(\kt)_{1}-\sig(\kt)_{l})\\
    &>-\frac{\sig(\kt)_{1}-\sig(\kt)_{l}}{2}
      +(\sig(\kt)_{1}-\sig(\kt)_{l})\\
    &=\frac{\sig(\kt)_{1}-\sig(\kt)_{l}}{2}.
  \end{align}
  This contradicts
  Claim~\ref{c:pref-upp-bounds-sig-diff-half} of
  Lemma~\ref{lem:pref-upp-bounds} and
  equation~\eqref{eq:pre-fitting-m-n}.
\end{proof}

\section{Variation of almost strictly unimodal
  sequences and real algebraic sets
  of critical points}
\label{sec:semi-st-unim-seq-crit-pts}

Among almost strictly unimodal sequences,
we have the following boundary sequences: first,
increasing and hill-shape sequences; second, decreasing
and hill-shape sequences; third, increasing and
{\it asymptotically hill-shape sequences}, which we explain
for the merged-log-concavity. By these
boundary sequences, we consider critical points on the variation
of almost strictly unimodal sequences.  We also obtain
real algebraic sets of the critical points in a suitable
setting.

\subsection{Increasing, hill-shape, and decreasing sequences}
\label{sec:inc-hill-dec}

As hill-shape sequences 
have  at least three terms in Definition~\ref{def:unimodal},
we state
the following classification on almost strictly unimodal sequences.

\begin{lemma}\label{lem:uni-type}
  Let $r=\{r_{i}\in \R \}_{i \ldZ s}$ be almost strictly
  unimodal with $s_{2}-s_{1}\geq 2$.  Let $\del$ be the
  mode of $r$.
  \begin{enumerate}
  \item \label{c:uni-type-st-increasing}If $\del=s_{2}$, then $r$
    is  strictly increasing.
  \item  \label{c:uni-type-claim-st-decreasing}
    If $\del= s_{1}$ and $r_{\del}>r_{\del+1}$,
    then $r$ is strictly decreasing.
  \item  \label{c:uni-type-hill-shape}
    In other cases, $r$ is  hill-shape.
  \end{enumerate}
\end{lemma}
\begin{proof}
  Let us prove Claim~\ref{c:uni-type-st-increasing}.  If
  $s_{2}<\infty$, then Claim~\ref{c:uni-type-st-increasing}
  follows, since $r=\{r_{i}\}_{i\leq \delta}$ is strictly
  increasing.  If $s_{2}=\infty$, then
  Claim~\ref{c:uni-type-st-increasing} follows, as
  $r_{\del}$ and $r_{\del+1}$ are non-existent.

  Claim~\ref{c:uni-type-claim-st-decreasing} holds, since
  $\{r_{i}\}_{i>\delta}$ is
  strictly decreasing.

  Let us prove
  Claim~\ref{c:uni-type-hill-shape}.  If $s_{1}<\del<s_{2}$,
  then $\del<\infty$ gives hill-shape sequences.  If
  $\del=s_{1}$ and $r_{\del}=r_{\del+1}$, then
  $s_{2}-s_{1}\geq 2$ gives hill-shape sequences.
\end{proof}

Then, we have the criteria below for some boundary
sequences.

\begin{proposition}\label{prop:inc-dec-hill}
  Consider an almost strictly unimodal
  $r=\{r_{i}\in \R \}_{i \ldZ s}$ such that
  $s_{2}-s_{1}\geq 2$.  Then, we have the following.
  \begin{enumerate}
  \item \label{c:inc-dec-hill-dec-hill} $r$ is decreasing and
    hill-shape if and only if
    $r_{s_{1}}=r_{s_{1}+1}$.
  \item \label{c:inc-dec-hill-inc-hill} $r$ is increasing and
    hill-shape if and only if $r$ satisfies
    the following conditions:
    \begin{align}
      s_{2}&<\infty;
             \label{ineq:inc-dec-hill-inc-hill-infty} \\
      r_{s_{2}-1}&=r_{s_{2}}.
                   \label{eq:inc-dec-hill-inc-hill-plateau}
    \end{align}
  \end{enumerate}
\end{proposition}
\begin{proof}
  Assume the mode $\del$ of $r$.  Let us prove
  Claim~\ref{c:inc-dec-hill-dec-hill}.  First, the if part
  follows, because $s_{2}-s_{1}\geq 2$ gives
  $r_{s_{1}}=r_{s_{1}+1}>r_{s_{1}+2}>\dots$, which is
  hill-shape and decreasing.
  Second, the only if part holds as follows.
   If $s_{1}<\del\leq s_{2}$, then
  $r$ is not decreasing by
  $r_{s_{1}}<s_{s_{1}+1}$.  If
   $s_{1}=\del$,  then 
    $r_{\del}=r_{\del+1}$
  by Claim~\ref{c:uni-type-claim-st-decreasing}
  of Lemma~\ref{lem:uni-type}.
  
  Let us prove Claim~\ref{c:inc-dec-hill-inc-hill}. First,
  the if part holds, since $s_{2}-s_{1}\geq 2$ gives
  $\dots <r_{s_{2}-2}<r_{s_{2}-1}=r_{s_{2}}$, which is
  increasing and hill-shape.  Second, let us prove the only
  if part.  If $s_{2}=\infty$, then $r$ is not increasing when
  $\del<\infty$, or strictly increasing when $\del=\infty$ by
  Claim~\ref{c:uni-type-st-increasing}. Thus,
  inequality~\ref{ineq:inc-dec-hill-inc-hill-infty} holds.
  Moreover, since the strictly decreasing part
  $\{r_{i}\}_{i>\delta}$ has to be trivial, $\del=s_{2}-1$ by
  Claim~\ref{c:uni-type-st-increasing} of
  Lemma~\ref{lem:uni-type}.  Hence,
  equation~\eqref{eq:inc-dec-hill-inc-hill-plateau} follows
  for the increasing $r$.
\end{proof}

\subsection{Critical points on almost strictly unimodal
  sequences} \label{sec:crit-value} We introduce the notion
of critical points on almost strictly unimodal sequences of
merged-log-concave parcels.  This uses the  {\it
  merged pairs} below.
\begin{definition}\label{def:strong-pair}
  Suppose a wrapped fitting path
  $P=\{\tp{s,l,m_{i},n_{i},k_{i}}\}_{i\ldZ \theta}$.  Assume a
  parcel family
  $\cF=\{
  \cF_{i}=\Lam(s,l_{i},w_{i},\scc_{i},f_{i,s},\phi_{i},\rho_{i},
  x_{i},\fX) \}_{i\ldZ \chi}$ for a finite gate $\chi$ such that
  each $\cF_{i}$ satisfies the following  conditions:
  \begin{enumerate}
  \item $\cF_{i}$ is static;
  \item $\cF_{i}$ is
    $\gAfX$-merged-log-concave;
  \item $\cF_{i}$ is $\tp{s,l,\fX}$-positive.
  \end{enumerate}
  Then, we call the tuple $\tp{P,\cF}$
  $\tp{\theta,\chi,\fX}$-merged pair.  For each
  $i\ldZ \chi$, we also call the tuple $\tp{P,\cF_{\chi_{i}}}$
  $\tp{\theta,\fX}$-merged pair for simplicity.
  \end{definition}
 
  By Theorem~\ref{thm:path-parcel-seq}, a
  $\tp{\theta,\chi,\fX}$-merged pair $\tp{P,\cF}$ gives the almost
  strictly unimodal sequences
  $u(\tp{P,\cF_{i}},r)$ of
  $i\ldZ \chi$ and $r\in \AfX$.  Hence, we define the following
  critical points on these
  $u(\tp{P,\cF_{i}},r)$.

\begin{definition}\label{def:crit-pts}
  Suppose a $\tp{\theta,\chi,\fX}$-merged pair $\zeta=\tp{P,\cF}$.
  Let $\ka_{i}=\tp{P,\cF_{i}}$ and $r\in \AfX$.
  \begin{enumerate}
  \item We call $r$ front critical point of
    $\zeta$, if $u(\ka_{i},r)$ is decreasing and
    hill-shape for each $i\ldZ \chi$.
  \item We call $r$ rear critical point of
    $\zeta$, if $u(\ka_{i},r)$ is increasing and
    hill-shape for each $i\ldZ \chi$.
  \end{enumerate}
\end{definition}

We use the term ``front critical point'', since
$u(\ka_{i},r)$ is decreasing and hill-shape when its front
terms satisfy
$u(\ka_{i},r)_{\thetao-1}= u(\ka_{i},r)_{\thetao}$ by
Proposition~\ref{prop:inc-dec-hill}.  Similarly, we use the
term ``rear critical point'', since $u(\ka_{i},r)$ of
$\thetat<\infty$ is increasing and hill-shape when its rear
terms satisfy
$u(\ka_{i},r)_{\thetat}= u(\ka_{i},r)_{\thetat+1}$.

Rear critical points exist only for
 finite-length fitting paths.
But
for an infinite gate $\theta$, there is
a $\tp{\theta,\fX}$-merged pair $\lam=\tp{P,\cG}$ with $b_{1},b_{2},b_{3}\in \AfX$
such that $u(\lam,r)$ is hill-shape if
$b_{1}\leq r<b_{2}$, and  strictly increasing
if $b_{2}\leq r\leq b_{3}$  (see
Section~\ref{sec:non-canonical-mediator-with-pt} for such an example).  Hence, we use the
following notion.

\begin{definition}\label{def:asympt-hill-shape}
  Consider a strictly increasing sequence
  $r=\{r_{i}\in \R\}_{i\ldZ \theta}$ for an
  infinite gate $\theta$.  Then, $r$ is said to be
  asymptotically hill-shape, if
  $\lim_{i\to\infty}\frac{r_{i+1}}{r_{i}}=1$.
\end{definition}

Then, we introduce the notion of {\it asymptotic
  critical points} for merged pairs.

\begin{definition}\label{def:asympt-crit-pts}
  Suppose a $\tp{\theta,\chi,\fX}$-merged pair
  $\zeta=\tp{P,\cF}$ of an
  infinite gate $\theta$.  We call $r\in \AfX$ asymptotic
  critical point of $\zeta$, if each $i\ldZ \chi$
  satisfies
  $\lim_{j\to\infty}\frac{u(\tp{P,\cF_{i}},r)_{j+1}}{u(\tp{P,\cF_{i}},r)_{j}}
    =1$.
\end{definition}

Hence, 
asymptotic critical points give
asymptotically hill-shape
sequences.

\begin{proposition}\label{prop:path-parcel-asymp-hill-shape}
  For a $\tp{\theta,\chi,\fX}$-merged pair $\tp{P,\cF}$,
  assume an asymptotic critical point $r\in \AfX$.  Let
  $\zeta_{i}=\tp{P,\cF_{i}}$.  Then, $u(\zeta_{i},r)$ is
  asymptotically hill-shape for each $i\ldZ \chi$.
\end{proposition}
\begin{proof}
  We prove that each $u(\zeta_{i},r)$ is strictly
  increasing.  Assume some
  $j\ldZ \chi$ such that $u(\zeta_{j},r)$ is hill-shape or
  decreasing.  Then, since $P$ is infinite-length, there
  exists $\lam \ldZ \theta$ such that
  $u(\zeta_{j},r)_{\lam}>   u(\zeta_{j},r)_{\lam+1}>\dots$.
  Furthermore, since $u(\zeta_{j},r) >0$ is strictly log-concave
  by Theorem~\ref{thm:path-parcel-seq}, we have
  $1\geq \frac{u(\zeta_{j},r)_{\lam+1}}{u(\zeta_{j},r)_{\lam}}>
  \frac{u(\zeta_{j},r)_{\lam+2}}{u(\zeta_{j},r)_{\lam+1}}>\dots$
  against  $r$ being an asymptotic critical point.
\end{proof}

\subsection{Parcel ratios
  and vanishing constraints
  of parcel numerators}
\label{sec:vanishing}

Suppose a merged pair $\tp{P,\cF}$ with an infinite-length
$P=\{\tp{s,l,m_{i},n_{i},k_{i}}\}_{i\ldZ \theta}$. Then, we discuss
parcel ratios $\frac{\cF_{n_{i}\bop k_{i}}}{\cF_{n_{i}}}$ by
almost strictly unimodal sequences and polynomials.  This
gives certain vanishing constraints on parcel numerators.

We use the  lemma below on infinite-length fitting
paths.

\begin{lemma}\label{lem:fit-pos-inf}
  Assume an infinite-length fitting path
  $P=\{\tp{s,l,m_{i},n_{i},k_{i}}\}_{ i\ldZ \theta}$. Then,
  $\lim_{i\to\infty}n_{i,j}=\infty$
  for each
  $j\in\oi(l)$.
  In particular, $s$ is an infinite gate.
\end{lemma}
\begin{proof}
  Let $i\ldZ \theta$ such that $i\geq \thetao+1$.  Then,
  $n_{i}$ is $n_{i-1}\bop k_{i-1}$ or
  $(n_{i-1}\bop k_{i-1})^{\ve}$.  Thus,
  $n_{i}\geq i-\thetao\in \Zgeo$ by $\sig(k)\geq 1$ and
  $n_{\thetao}\geq 0$. So, we have the former assertion.
  This
  also gives the latter assertion by the inclusion condition
  of each fitting tuple of $P$.
\end{proof}

We also have the following lemma on parcel ratios.

\begin{lemma}\label{lem:path-parcel-ratio}
  Consider a merged pair $\tp{P,\cF}$ such that
  $P=\{\tp{s,l,m_{i},n_{i},k_{i}}\}_{i\ldZ \theta}$ and
  $\cF = \Lam(s,l,w,\scc,\fs,\phi,\rho,x,\fX)$.  Let $\psi$ be the
  $l$-canonical mediator and $\ka=k_{\thetao}$.  For each
  $j\ldZ \theta$, let
  \begin{align}
    G(j,x)=
    \prod_{i\in\oi(l),\lam\in\oi(\sig(\ka)^{\ve}_{i})}
    (1-x_{i}^{n_{j,i}+\lam})^{w_{i}}.
  \end{align}
  Then, we have the following.
  \begin{enumerate}
  \item \label{c:path-parcel-ratio-by-canonical-med}
    If $j\ldZ \theta$, then
    \begin{align}
      \frac{\cF_{n_{j}\bop k_{j}}}{\cF_{n_{j}}}
      =\frac{f_{s,n_{j}\bop \ka}}{f_{s,n_{j}}}
      \cdot
      \frac{\prod \psi(x)^{\sig(\ka)^{\ve}\rc w}}{\prod
      \phi(x)^{\sig(\ka)^{\ve}\rc w}}
      \cdot
      \frac{1}{G(j,x)}.
    \end{align}
  \item \label{c:path-parcel-ratio-one-limit} If $P$ is
    infinite-length, then
    $\lim_{j\to \infty}G(j,x)|_{x=\tp{x_{1}(r),\dots, x_{l}(r)}}
    =1$ for each $r\in \AfX$.
  \end{enumerate}
\end{lemma}
\begin{proof}
  Let us prove
  Claim~\ref{c:path-parcel-ratio-by-canonical-med}.  For
  equivalent supports $k_{j}$ and $\ka$,
  $\sig(k_{j})=\sig(\ka)$ and
  $n_{j}\bop k_{j}=n_{j}\bop \ka$.  Then, since
  $n_{j}\bop k_{j}-n_{j}=\sig(\ka)^{\ve}$, we obtain
  \begin{align}
    \frac{\cF_{n_{j}\bop k_{j}}}{\cF_{n_{j}}}
    &=\frac{\prod \phi(x)^{n_{j}\rc w} \cdot
      [n_{j}]!_{x}^{w}}{f_{s,n_{j}}} \cdot
      \frac{f_{s,n_{j}\bop k_{j}}}{\prod \phi(x)^{(n_{j}\bop k_{j})\rc w} \cdot
      [n_{j}\bop k_{j}]!_{x}^{w}}\\
    &=\frac{f_{s,n_{j}\bop \ka}}{f_{s,n_{j}}}
      \cdot \frac{1}{\prod \phi(x)^{\sig(\ka)^{\ve}\rc w}}
      \cdot \frac{      [n_{j}]!_{x}^{w}}{  [n_{j}\bop \ka]!_{x}^{w}}\\
    &=\frac{f_{s,n_{j}\bop \ka}}{f_{s,n_{j}}}
      \cdot \frac{\prod \psi(x)^{\sig(\ka)^{\ve}\rc w}}{\prod
      \phi(x)^{\sig(\ka)^{\ve}\rc w}}
      \cdot \frac{      [n_{j}]!_{x}^{w}}{
      \prod_{i\in\oi(l)}(1-x_{i})^{\sig(\ka)^{\ve}_{i}\cdot w_{i}} \cdot
      [n_{j}\bop \ka]!_{x}^{w}}.
  \end{align}
  Thus, Claim~\ref{c:path-parcel-ratio-by-canonical-med}
  holds, since
  \begin{align}
    \prod_{i\in\oi(l)}(1-x_{i})^{\sig(\ka)^{\ve}_{i}\cdot w_{i}} \cdot
    \frac{  [n_{j}\bop \ka]!_{x}^{w}}{      [n_{j}]!_{x}^{w}}
    &=
      \prod_{i\in\oi(l)} (1-x_{i})^{\sig(\ka)^{\ve}_{i}\cdot w_{i}} \cdot
      [n_{j,i}+1]^{w_{i}}_{x_{i}}\dots
      [n_{j,i}+\sig(\ka)^{\ve}_{i}]^{w_{i}}_{x_{i}}\\
    &=
      \prod_{i\in\oi(l)} (1-x_{i}^{n_{j,i}+1})^{w_{i}}\dots
      (1-x_{i}^{n_{j,i}+\sig(\ka)^{\ve}_{i}})^{w_{i}}.
  \end{align}

  Claim~\ref{c:path-parcel-ratio-one-limit} holds by
  Lemma~\ref{lem:fit-pos-inf}, since $0<x_{i}(r)<1$ by
  Claim~\ref{c:adm-bounds} of
  Lemma~\ref{lem:adm-bounds-poring-inclusions}.
\end{proof}

We put the notion below to
discuss asymptotic critical points by
polynomials.
\begin{definition}\label{def:tame} 
  Assume $\cF = \Lam(s,l,w,\scc,\fs,\phi,\rho,x,\fX)$ and an
  infinite-length fitting path
  $P=\{\tp{s,l,m_{i},n_{i},k_{i}}\}_{i\ldZ \theta}$.
  \begin{enumerate}
  \item 
  We call
  $\cF$ tame along $P$ by $t\in \Q(\fX)$, if $t\geAfX 0$ and
  \begin{align}
    t(r)=\lim_{i\to\infty}\frac{f_{s,n_{i}\bop
    k_{i}}}{f_{s,n_{i}}}(r)    
  \end{align}
  for each $r\in \AfX$.
  We call this $t$ tame factor of $\cF$ along $P$.
\item Let $\cF$ be tame along $P$ by $t\in \Q(\fX)$.  Then,
  we put
  $\AIr(\fX)(\zeta),\AIr(\fX)(t,\zeta)\in \Q[\fX]^{2}$ such that
  $\AIr(\fX)(\zeta)=\Ir(\fX)(\prod \phi(x)^{\sig(\ka)^{\ve}\rc w})$
  and
  $\AIr(\fX)(t,\zeta)=\Ir(\fX)(t \cdot \prod \psi(x)^{\sig(\ka)^{\ve}\rc
    w})$.
\end{enumerate}
\end{definition}

On merged pairs and almost strictly unimodal sequences, we
then prove the following limit properties, which are
independent to mediators.

\begin{proposition}\label{prop:path-parcel-ratio-unim-poly}
  Consider a merged pair $\zeta=\tp{P,\cF}$ with
  $\cF = \Lam(s,l,w,\scc,\fs,\phi,\rho,x,\fX)$ and an infinite-length
  $P=\{\tp{s,l,m_{i},n_{i},k_{i}}\}_{i\ldZ \theta}$.
   Let $\psi$
  be the $l$-canonical mediator and
  $\ka=k_{\thetao}$.
  \begin{enumerate}
  \item \label{c:path-parcel-ratio-zero-unim}
    Suppose 
    $\lim_{j\to\infty}\frac{f_{s,n_{j}\bop k_{j}}}{f_{s,n_{j}}}(r)=0$
    for some $r\in \AfX$.
    Then, we have
    \begin{align}
      \lim_{j\to\infty}\frac{\cF_{n_{j}\bop k_{j}}}{\cF_{n_{j}}}(r)
      =0.
      \label{eq:path-parcel-ratio-zero-lim}
    \end{align}
    In particular, $u(\zeta,r)$ is hill-shape or decreasing
    such that
    \begin{align}
      \lim_{i\to \infty}u(\zeta,r)_{i}=0.
      \label{eq:path-parcel-ratio-unim-lim}
    \end{align}
    
  \item \label{c:path-parcel-ratio-no-inf-lim}
    For each $r\in \AfX$,
    $\lim_{j\to\infty}\frac{f_{s,n_{j}\bop k_{j}}}{f_{s,n_{j}}}(r)=\infty$
    does not hold.

  \item \label{c:path-parcel-ratio-lim-poly}
    Let $\cF$ be tame along $P$  by
    $t\in \Q(\fX)$ and $r\in \AfX$.  Then, we have
    \begin{enumerate}
    \item
      \label{c:path-parcel-ratio-lim-poly-nonneg}
      $\lim_{j\to\infty}\frac{\cF_{n_{j}\bop
          k_{j}}}{\cF_{n_{j}}}(r) =t(r)\cdot \frac{\prod
        \psi(x(r))^{\sig(\ka)^{\ve}\rc w}}{\prod
        \phi(x(r))^{\sig(\ka)^{\ve}\rc w}}\in \bR_{\geq 0}$;
    \item
      \label{c:path-parcel-ratio-poly-exp}
      $\tdet(\AIr(\fX)(\zeta),\AIr(\fX)(t,\zeta))(r)=0$ if and
      only if
      $\lim_{j\to\infty}\frac{\cF_{n_{j}\bop
          k_{j}}}{\cF_{n_{j}}}(r) =1$.
    \end{enumerate}
  \end{enumerate}
\end{proposition}
\begin{proof}
  Let us prove Claim~\ref{c:path-parcel-ratio-zero-unim}.
  First, we obtain
  equation~\eqref{eq:path-parcel-ratio-zero-lim} from
  Claim~\ref{c:AfX-pos-med-qfact} of
  Proposition~\ref{prop:AfX-pos} and
  Lemma~\ref{lem:path-parcel-ratio}.  Second, let us obtain
  equation~\eqref{eq:path-parcel-ratio-unim-lim}.  The
  path-parcel sequence $u(\zeta,r)$ is almost strictly unimodal
  by Theorem~\ref{thm:path-parcel-seq}, and of positive real
  numbers by Definition~\ref{def:strong-pair}.  Therefore,
  equation~\eqref{eq:path-parcel-ratio-zero-lim} gives
  $\lam\ldZ \theta$ such that
  $\{u(\zeta,r)_{j}=\cF_{n_{j}}(r)\}_{j\geq \lam}$ is strictly
  decreasing. This gives $v\in \R_{\geq 0}$ such that
  $\lim_{j\to\infty}u(\zeta,r)_{j}=v$.  If $v>0$, then
  $\lim_{j\to\infty}\frac{u(\zeta,r)_{j+1}}{u(\zeta,r)_{j}}=1$ against
  equation~\eqref{eq:path-parcel-ratio-zero-lim}.

  Let us prove
  Claim~\ref{c:path-parcel-ratio-no-inf-lim}. Let us assume
  otherwise. Then, Lemma~\ref{lem:path-parcel-ratio} and
  Claim~\ref{c:AfX-pos-med-qfact} of
  Proposition~\ref{prop:AfX-pos} imply
  $\lim_{j\to\infty}\frac{u(\zeta,r)_{j+1}}{u(\zeta,r)_{j}} =\infty$.
  However, because $u(\zeta,r)$ is strictly log-concave by
  Theorem~\ref{thm:path-parcel-seq}, each $j\ldZ \theta$
  satisfies
  $\infty >\frac{u(\zeta,r)_{j+1}}{u(\zeta,r)_{j}}
  >\frac{u(\zeta,r)_{j+2}}{u(\zeta,r)_{j+1}}>\dots$.

  Since $t\geAfX 0$,
  Claim~\ref{c:path-parcel-ratio-lim-poly-nonneg} holds by
  Lemma~\ref{lem:path-parcel-ratio} and
  Claim~\ref{c:AfX-pos-med-qfact} of
  Proposition~\ref{prop:AfX-pos}.  Let us prove
  Claim~\ref{c:path-parcel-ratio-poly-exp}.
  Claim~\ref{c:AfX-pos-med-qfact} of
  Proposition~\ref{prop:AfX-pos} gives
  $\prod \phi(x)^{\sig(\ka)^{\ve}\rc w}\gAfX 0$.  Thus,
  Claim~\ref{c:path-parcel-ratio-lim-poly-nonneg} implies
  $\lim_{j\to\infty}\frac{\cF_{n_{j}\bop k_{j}}}{\cF_{n_{j}}}(r)
  =1$ if and only if
  $t(r)\cdot\prod \psi(x(r))^{\sig(\ka)^{\ve}\rc w}=\prod
  \phi(x(r))^{\sig(\ka)^{\ve}\rc w}$.  Since
  $t\cdot\prod \psi(x)^{\sig(\ka)^{\ve}\rc w}\geAfX 0$ by
  Claim~\ref{c:AfX-pos-med-qfact} of
  Proposition~\ref{prop:AfX-pos}, we obtain the assertion by
  Claims~\ref{c:geAfX-gAfX-scc-geAfX}
  and~\ref{c:geAfX-gAfX-scc-gAfX} of
  Lemma~\ref{lem:geAfX-gAfX-scc}.
\end{proof}

We use the following infinite-length fitting paths.

\begin{definition}\label{def:paths-w-one}
  Assume $\lam\in \Zgeo$, $h\in \Z$, $l=1$, and
  $s=\tp{0,\infty}$.  Then, let
  $\pi(\lam,h)=\bP_{s,l,\iota^{\lam}(1),h}$.
  Also, let $\pi(\lam)=\pi(\lam,0)$ for simplicity. 
\end{definition}

We then have the following explicit description
of $\pi(\lam,h)$.
\begin{lemma}\label{lem:width-one-fititng-path}
  Let $\lam\in \Zgeo$, $h\in \Z$, $l=1$, $s=\tp{0,\infty}$,
  and $\theta=\tp{1,\infty}$.
  \begin{enumerate}
  \item
    \label{c:width-one-fititng-path-fitting-tuples}
    We have
    $\pi(\lam,h) =\{\tp{s,l,m_{i},n_{i},k_{i}}\}_{i\ldZ
      \theta}$ with $a_{i}=\nu(k_{i})$ and
    $b_{i}=\nu(m_{i},n_{i},k_{i})$ such that each
    $i\ldZ \theta$ satisfies
    \begin{align}
      a_{i}=k_{i}&=\tp{0,\lam}\in \Ztl,\\
      b_{i}&=\tp{\lam i+h,\lam(i+1)+h}\in \Ztl,\\
      m_{i}=n_{i}&=\tp{\lam i+h}\in\Zl.
    \end{align}
  \item
    \label{c:width-one-fititng-path-zero}
    We have $m_{1}\bom k_{1}=\tp{h}\in \Zl$.
  \item\label{c:wrapped-iff} $\pi(\lam,h)$ is wrapped if and
    only if $h\in \Zgez$.  In particular, $\pi(\lam)$ is
    wrapped.
  \end{enumerate}
\end{lemma}
\begin{proof}
  We have $\iota^{\lam}(1)\in Y(\lam,l)$ and
  $a_{i}=\nu(k_{i})=k_{i}$ by $a_{i,1}=0$.  Thus,
  Claim~\ref{c:width-one-fititng-path-fitting-tuples} holds
  by Theorem~\ref{thm:fitting-path-young}, since
  $z(l,0,\iota^{\lam}(1)) =\tp{\len_{0}(\iota^{\lam}(1))}=\tp{0}$
  and $\len(l,\iota^{\lam}(1))=\tp{\lam}$.
  Claim~\ref{c:width-one-fititng-path-zero} follows from
  $\sig(k_{i})=\tp{\lam}$.  We then obtain
  Claim~\ref{c:wrapped-iff}, because $\pi(\lam,h)$ is wrapped
  if and only if $\tp{h}\ldZo \tp{0,\infty}$ by
  Claim~\ref{c:width-one-fititng-path-zero}, and
  $\pi(\lam)=\pi(\lam,0)$.
\end{proof}

Moreover,  we use the  notation below for $\pi(\lam)$.
\begin{definition}
  Suppose a parcel
  $\cF=\Lam(s,l,w,\scc,f_{s},\phi,\rho,x,\fX)$ and
  $n\in \Zgeo$.  Then, we write $\Omg_{n}(\cF)$ for the
  tuple $\tp{\pi(n),\cF}$. In particular, let
  $\Omg(\cF)=\Omg_{1}(\cF)$.
\end{definition}

Then, we obtain the following vanishing constraint on
parcel numerators.
\begin{proposition}
  \label{prop:vanishing-inf-many-zero}
  Let $p\in \Z$.
  For each
  $h\in \Zgez$, suppose $u_{h}\in \Zget$,
  $v_{h}\in \Zleo$,
  $\lam_{h,0,0}\in \Z$,
  and finitely
  many non-zero $\lam_{h,i,j}\in \Z$ of $i\in \Zgeo$ and
  $j\in \Z_{\neq 0}$ with the following conditions:
  \begin{enumerate}
  \item \label{a:vanishing-inf-many-zero-uv-existence}
    $i\leq u_{h}$ and $j\geq v_{h}$, if $i\in \Zgeo$ and
    $j\in \Zleo$ satisfy $\lam_{h,i,j}\neq 0$;
  \item \label{a:vanishing-inf-many-zero-uv-ineq}
    $u_{h+1}\geq u_{h}$ and $v_{h+1}\leq v_{h}$;
  \item \label{a:vanishing-inf-many-zero-lam-lim}
    $\lim_{h\to\infty}\lam_{h,u_{h},v_{h}}=p$;
  \item \label{a:vanishing-inf-many-zero-lim-est}
    $\lim_{h\to\infty} \sum_{i\in \Zgeo,j\in \Z_{\neq 0},\tp{i,j}\neq
      \tp{u_{h}, v_{h}}}
    \abs{\lam_{h,i,j}} y^{h} =\lim_{h\to
      \infty}\lam_{h,0,0}y^{h}=0$, if $0<y<1$.
  \end{enumerate}
  Let $s=\tp{0,\infty}$, $l=1$, and $q$ be
  $\scc$-admissible.  Consider a $\tp{s,l,\fX}$-positive
  and
  $\scc'$-merged-log-concave
  $\cF= \Lam(s,l,w,\scc,\fs,\phi,\rho,x,\fX)$ such that
   each $\tp{h}\ldZl s$ satisfies
  \begin{align}
    f_{s,\tp{h}}=
    \lam_{h,0,0}
    +\sum_{i\in \Zgeo,j\in \Z_{\neq 0}}
    \lam_{h,i,j}q^{j h^i}\in \Q(\fX).
  \end{align}
   Then, we obtain
  \begin{align}
    \lam_{h,u_{h},v_{h}}= 0    
  \end{align}
  for infinitely many $h\in \Zgez$.
\end{proposition}
\begin{proof}
  Suppose $r\in \AfX$ and $h\in \Z_{\geq 2}$.
  First, let us prove
  \begin{align}
    \lim_{h\to\infty}\lam_{h,0,0}
    q(r)^{-v_{h} h^{u_{h}}}=0.
    \label{eq:vanishing-inf-many-zero-lim-lam-q-power}
  \end{align}
  Claim~\ref{c:adm-bounds} of
  Lemma~\ref{lem:adm-bounds-poring-inclusions} gives
  \begin{align} 
    0<q(r)<1. \label{ineq:vanishing-inf-many-zero-q-bounds}
  \end{align}
  Also, since  $-v_{h}\in \Zgeo$,
  each $i\in \Zgeo$ satisfies
  \begin{align}
    -v_{h} h^{i}\geq h.
    \label{ineq:vanishing-inf-many-zero-vh-bounds}
  \end{align}
  Thus,
  inequalities~\eqref{ineq:vanishing-inf-many-zero-q-bounds}
  and~\eqref{ineq:vanishing-inf-many-zero-vh-bounds} and
  Assumption~\ref{a:vanishing-inf-many-zero-lim-est} give
  equation~\eqref{eq:vanishing-inf-many-zero-lim-lam-q-power}.
  
  Second, let $p(i,j,h)=j h^i-v_{h}h^{u_{h}}$ for
  $i\in \Zgeo$ and $j\in \Z_{\neq 0}$.  When
  $\lam_{h,i,j}\neq 0$ and $\tp{i,j}\neq \tp{u_{h},v_{h}}$, we
  prove
  \begin{align}
    p(i,j,h) \geq h.
    \label{ineq:vanishing-inf-many-zero-jhvh-bounds}
  \end{align}
  Assumption~\ref{a:vanishing-inf-many-zero-uv-existence}
  gives the following three cases:
  \begin{itemize}
    \item if $j\in\Zgeo$,
  then
  inequality~\eqref{ineq:vanishing-inf-many-zero-vh-bounds}
  implies $p(i,j,h)\geq -v_{h}h^{u_{h}}\geq h$; 
\item if $j\in\oi(v_{h},-1)$ and $i\in\oi(u_{h}-1)$, then
  $h^{u_{h}-i}-1\geq 1$ and
  inequality~\eqref{ineq:vanishing-inf-many-zero-vh-bounds}
  imply
  $p(i,j,h) \geq v_{h}h^{i}-v_{h}h^{u_{h}}
  =-v_{h}h^{i}(h^{u_{h}-i}-1) \geq -v_{h}h^{i}\geq h$;
    \item 
      if $j\in\oi(v_{h}+1,-1)$ and $i=u_{h}$,
      then $j-v_{h}\geq 1$ yields
   $p(i,j,h)
   =h^{u_{h}}(j-v_{h})\geq h$.
   \end{itemize}

   Third, let
   $\pi(1) =\{\tp{s,l,m_{i},n_{i},k_{i}}\}_{i\ldZ \theta}$ with
  $U(h,q)=q^{-v_{h} h^{u_{h}}}f_{s,n_{h}}$ and
  $V(h,q)=
  U(h,q)-\lam_{h,0,0}
  q^{-v_{h} h^{u_{h}}}
  -\lam_{h,u_{h},v_{h}}$.    
  Then,
 Assumption~\ref{a:vanishing-inf-many-zero-lim-est}
 yields
  \begin{align}
    \lim_{h\to\infty}        V(h,q(r))
    &=0,
      \label{ineq:vanishing-inf-many-zero-lim-V-zero}
  \end{align}
  because
  Claim~\ref{c:width-one-fititng-path-fitting-tuples} of
  Lemma~\ref{lem:width-one-fititng-path} and
  inequalities~\eqref{ineq:vanishing-inf-many-zero-q-bounds}
  and~\eqref{ineq:vanishing-inf-many-zero-jhvh-bounds} give
  \begin{align}
    \abs{V(h,q(r))}
    &\leq
      \sum_{i\in \Zgeo,j\in \Z_{\neq 0},\tp{i,j}\neq \tp{u_{h},v_{h}}}
      \abs{\lam_{h,i,j}             q(r)^{p(i,j,h)}}
      \leq
      \sum_{i\in \Zgeo,j\in \Z_{\neq 0},\tp{i,j}\neq \tp{u_{h},v_{h}}}
      \abs{\lam_{h,i,j}} q(r)^{h}.
  \end{align}
  
  Fourth, let us prove the assertion by contradiction.
  Then, Assumption~\ref{a:vanishing-inf-many-zero-lam-lim}
  implies
  $\lim_{h\to\infty}\lam_{h,u_{h},v_{h}}=p\neq 0$.
  Thus,
  $\lim_{h\to\infty}U(h,q(r))=p\neq 0$ by equations~\eqref{eq:vanishing-inf-many-zero-lim-lam-q-power} and~\eqref{ineq:vanishing-inf-many-zero-lim-V-zero}.
In particular,
  \begin{align}
    \lim_{h\to\infty}\frac{U(h+1,q(r))}{U(h,q(r))}=1.
    \label{eq:vanishing-inf-many-zero-U-quot-one}
  \end{align}
  Moreover,
  inequality~\eqref{ineq:vanishing-inf-many-zero-vh-bounds}
  and Assumption~\ref{a:vanishing-inf-many-zero-uv-ineq}
  imply
  \begin{align}
    v_{h}h^{u_{h}}
    -v_{h+1}(h+1)^{u_{h+1}}
    \geq
    v_{h}h^{u_{h}}
    -v_{h}(h+1)^{u_{h}}
    \geq -v_{h}u_{h}
    h^{u_{h}-1}\geq h,
  \end{align}
  since
  $-v_{h}(h+1)^{u_{h}}
  =-v_{h }h^{u_{h}}-v_{h}u_{h}h^{u_{h}-1} -\dots$.
  Then,
  inequalities~\eqref{ineq:vanishing-inf-many-zero-q-bounds}
  yield
  \begin{align}
    \lim_{h\to\infty}
    \frac{q(r)^{v_{h+1} (h+1)^{u_{h+1}}}}{q(r)^{v_{h} h^{u_{h}}}}
    =    \lim_{h\to\infty}
    \frac{1}{q(r)^{v_{h} h^{u_{h}}-v_{h+1} (h+1)^{u_{h+1}}}}
    =\infty.
  \end{align}
  Therefore,
  because
  $\frac{f_{s,n_{h}\bop k_{h}}}{f_{s,n_{h}}} =
  \frac{f_{s,n_{h+1}}}{f_{s,n_{h}}}$,  
  equation~\eqref{eq:vanishing-inf-many-zero-U-quot-one}
  gives
  \begin{align}
    \lim_{h\to\infty}\frac{f_{s,n_{h}\bop k_{h}}}{f_{s,n_{h}}}(r)
    =
    \lim_{h\to\infty}
    \frac{q(r)^{v_{h+1} (h+1)^{u_{h+1}}}}
    {q(r)^{v_{h} h^{u_{h}}}}
    \cdot
    \lim_{h\to\infty}\frac{U(h+1,q(r))}{U(h,q(r))}
    =\infty.
  \end{align}
  This
  contradicts Claim~\ref{c:path-parcel-ratio-no-inf-lim} in
  Proposition~\ref{prop:path-parcel-ratio-unim-poly},
  since $\Omg(\cF)$ is a merged
  pair by Claim~\ref{c:wrapped-iff} of
  Lemma~\ref{lem:width-one-fititng-path}.  
\end{proof} 

Also, we obtain the following vanishing constraint on parcel
numerators without assuming $t$-monomials (see
Remark~\ref{rmk:monomial-poly}).  Compared to
Proposition~\ref{prop:vanishing-inf-many-zero}, the
following allows $j\in \Q_{\neq 0}$ and takes
$\lam_{i,j}\in \Q$, which is independent to $m\in \Zl$.

\begin{proposition}
  \label{prop:vanishing-const-neg-quad}
  Assume finitely many non-zero
  $\lam_{0,0}, \lam_{i,j}\in \Q$ for $i\in \Zgeo$ and
  $j\in \Q_{\neq 0}$.  Let $s=\tp{0,\infty}$, $l=1$, and $q$ be
  $\scc$-admissible.  Consider a $\tp{s,l,\fX}$-positive
  and
  $\scc'$-merged-log-concave
  $\cF= \Lam(s,l,w,\scc,\fs,\phi,\rho,x,\fX)$ such that each
  $\tp{h}\ldZl s$ satisfies
  \begin{align}
    f_{s,\tp{h}}
    &=
      \lam_{0,0}
      +\sum_{i\in \Zgeo,j\in \Q_{\neq 0}}
      \lam_{i,j}q^{j h^i}\in \Q(\fX).
  \end{align}
    Then, we have
    \begin{align}
    \lam_{i,j}=0
  \end{align}
  for each $i\in \Zget$ and $j\in \Q_{<0}$.
\end{proposition}
\begin{proof}
  We prove the assertion by contradiction.  Hence, suppose
  the largest $u\in \Zget$ such that there exists
  $j\in \Q_{<0}$ with $\lam_{u,j}\neq 0$, and the smallest
  $v\in \Q_{<0}$ such that $\lam_{u,v}\neq 0$.  Let
  $\pi(1) =\{\tp{s,l,m_{i},n_{i},k_{i}}\}_{i\ldZ \theta}$,
  $h\in \Z_{\geq 2}$, and $r\in \AfX$.

     First, suppose $p(i,j,h)=j h^i-vh^{u}$ for $i\in \Zgeo$
     and $j\in \Q_{\neq 0}$.  Then, by
     Claim~\ref{c:width-one-fititng-path-fitting-tuples} of
     Lemma~\ref{lem:width-one-fititng-path}, we consider
  \begin{align}
    U(h,q)&=q^{-v h^{u}}f_{s,n_{h}}
            =\lam_{0,0}
            q^{-v h^{u}}
            +\lam_{u,v}+  
            \sum_{i\in \Zgeo, j\in \Q_{\neq 0}, \tp{i,j}\neq \tp{u, v}}
            \lam_{i,j}q^{p(i,j,h)}.
  \end{align}
  Since $u\in \Zget$ and $v\in \Q_{<0}$, 
  $\lim_{h\to\infty}-v h^{u}=\infty$.
  Hence, if $\lam_{i,j}\neq 0$ and $\tp{i,j}\neq \tp{u,v}$,  then
  $\lim_{h\to\infty}p(i,j,h)= \infty$,
  as we  have the only three cases:
  $j>0$;  $0>j> v$ and $i=u$; and, $0>j$ and
   $1\leq i<u$.
   In particular, since $0<q(r)<1$ by
  Claim~\ref{c:adm-bounds} of
  Lemma~\ref{lem:adm-bounds-poring-inclusions}, we have
  \begin{align}
    \lim_{h\to\infty}U(h,q(r))=\lam_{u,v}.
    \label{eq:vanishing-const-neg-quad-lim-U}
  \end{align}
  
  Second,  $u\in \Zget$ and
  $v\in \Q_{<0}$
  imply
  $vh^{u} -v(h+1)^{u}\geq -vu h^{u-1}\geq -v h>0$,
  since $-v(h+1)^{u}= -v h^{u}-vuh^{u-1}-\dots$.  
  Hence, we have
  \begin{align}
    \lim_{h\to\infty}
    \frac{ q(r)^{v (h+1)^{u}}}{q(r)^{v h^{u}}}
    =
    \lim_{h\to\infty}
    \frac{1}{q(r)^{v h^{u}-v (h+1)^{u}}}
    =\infty.
  \end{align}
  Then, since $\lam_{u,v}\neq 0$,
  equation~\eqref{eq:vanishing-const-neg-quad-lim-U} yields
  \begin{align}
    \lim_{h\to\infty}\frac{f_{s,n_{h}\bop k_{h}}}{f_{s,n_{h}}}(r)
    =
    \lim_{h\to\infty}
    \frac{ q(r)^{v (h+1)^{u}}}{q(r)^{v h^{u}}}
    \cdot \lim_{h\to\infty}\frac{U(h+1,q(r))}{U(h,q(r))}
    =\infty.
  \end{align}
  This contradicts
  Claim~\ref{c:path-parcel-ratio-no-inf-lim} in
  Proposition~\ref{prop:path-parcel-ratio-unim-poly}, since
  $\Omg(\cF)$ is a merged pair by Claim~\ref{c:wrapped-iff}
  of Lemma~\ref{lem:width-one-fititng-path}.
\end{proof}

\subsection{Merged pairs
  and critical points}
\label{sec:variation}

By front, rear, and asymptotic critical points of merged
pairs, we discuss the variation of almost strictly unimodal
sequences. Also, we obtain real algebraic sets not only of
front and rear critical points, but also of asymptotic
critical points by tame factors.  Moreover, we introduce the
notion of semi-phase transitions and phase transitions by
these critical points.  In particular, we obtain front phase
transitions by monomial parcels of general lengths.
 
First,
we introduce the notion of {\it path-parcel differences}.
\begin{definition}\label{def:p-variety-inf}
  Suppose a $\tp{\theta,\chi,\fX}$-merged pair
  $\zeta=\tp{P,\cF}$ such that
  $\cF_{i} = \Lam(s,l,w_{i},\scc_{i},f_{i,s},\phi_{i},
  \rho_{i},x_{i},\fX)$ for $i\ldZ \chi$ and
  $P=\{\tp{s,l,m_{i},n_{i},k_{i}}\}_{i\ldZ \theta}$.  Assume the
   canonical $l$-mediator $\psi$. Then, we define the following
  families of rational functions:
  \begin{enumerate}
  \item
    $\FD(\zeta)=\{  \cF_{i,m_{\thetao}\bom
      k_{\thetao}}-\cF_{i,n_{\thetao}}\in \Q(\fX)\}_{i\ldZ \chi}$;
  \item
    $\RD(\zeta)=\{ \cF_{i,n_{\thetat}}-\cF_{i,n_{\thetat}\bop
      k_{\thetat}} \in \Q(\fX)\}_{i\ldZ\chi}$,
    if $P$ is finite-length;
  \item
    $\AD(\zeta)=\{ \prod \phi_{i}(x_{i})^{\sig(\ka)^{\ve}\rc
      w} - t_{i} \cdot \prod \psi(x_{i})^{\sig(\ka)^{\ve}\rc
      w}\in \Q(\fX)\}_{i\ldZ\chi}$,
    if each $\cF_{i}$ is tame along $P$ by
    $t_{i}\in \Q(\fX)$.
  \end{enumerate}
  We call $\FD(\zeta),\RD(\zeta)$, and $\AD(\zeta)$ front, rear, and
  asymptotic path-parcel differences.

  Also, for a
  $\tp{\theta,\fX}$-merged pair $\ka=\tp{P,\cF}$,
  suppose the
  $\tp{\theta,\tp{1,1},\fX}$-merged pair
  $\ka'=\tp{P,\tp{\cF}}$.  Then, let
  $\FD(\ka)=\FD(\ka')_{1}\in \Q(\fX)$,
  $\RD(\ka)=\RD(\ka')_{1}\in \Q(\fX)$, and
  $\AD(\ka)=\AD(\ka')_{1}\in \Q(\fX)$ for our convenience.
\end{definition}
Then, we obtain the following variation of
almost strictly unimodal sequences.
\begin{theorem}
  \label{thm:path-parcel-variation}
  Suppose a $\tp{\theta,\fX}$-merged pair $\zeta=\tp{P,\cF}$
  with $r\in \AfX$.
  \begin{enumerate}
  \item 
    \label{c:path-parcel-variation-front}
    Then, we have the following.
    \begin{enumerate}
    \item \label{c:path-parcel-variation-front-pos}
      $u(\zeta,r)$ is strictly decreasing, if
      $\FD(\zeta)(r)>0$.
    \item \label{c:path-parcel-variation-front-zero}
      $u(\zeta,r)$ is hill-shape and decreasing, if
      $\FD(\zeta)(r)=0$.
    \item \label{c:path-parcel-variation-front-negative}
      $u(\zeta,r)$ is two-slope hill-shape or strictly
      increasing, if $\FD(\zeta)(r)<0$.
    \end{enumerate}
  \item Suppose that $P$ is finite-length.
    Then, we have the following.
    \begin{enumerate}
    \item \label{c:path-parcel-variation-rear-pos}
      $u(\zeta,r)$ is two-slope hill-shape or strictly
      decreasing, if $ \RD(\zeta)(r)>0$.
    \item \label{c:path-parcel-variation-rear-zero}
      $u(\zeta,r)$ is hill-shape and increasing, if
      $\RD(\zeta)(r)=0$.
    \item \label{c:path-parcel-variation-rear-negative}
      $u(\zeta,r)$ is strictly increasing, if $ \RD(\zeta)(r)<0$.
    \end{enumerate}
  \item Suppose that $\cF$ is tame along $P$. Then, we have the following.
    \begin{enumerate}
    \item \label{c:path-parcel-variation-asym-pos}
      $u(\zeta,r)$ is two-slope hill-shape or strictly
      decreasing, if $\AD(\zeta)(r)>0$.
    \item \label{c:path-parcel-variation-asym-zero}
      $u(\zeta,r)$ is asymptotically hill-shape, if
      $\AD(\zeta)(r)=0$.
    \item \label{c:path-parcel-variation-asym-negative}
      $u(\zeta,r)$ is strictly increasing, if $\AD(\zeta)(r)<0$.
    \end{enumerate}
  \end{enumerate}
\end{theorem}
\begin{proof}
  Let $P=\{\tp{s,l,m_{i},n_{i},k_{i}}\}_{i\ldZ \theta}$.  Then,
  the first and second terms of $u(\zeta,r)$ are
  $\cF_{m_{\thetao}\bom k_{\thetao}}(r)$ and
  $\cF_{n_{\thetao}}(r)$. Thus,
  Claim~\ref{c:path-parcel-variation-front-pos} holds by
  Claim~\ref{c:inc-dec-hill-dec-hill} of
  Proposition~\ref{prop:inc-dec-hill}, because $u(\zeta,r)$ is
  almost strictly unimodal by
  Theorem~\ref{thm:path-parcel-seq}.
  Claims~\ref{c:path-parcel-variation-front-zero}
  and~\ref{c:path-parcel-variation-front-negative} hold
  similarly.
  
  Let us prove Claim~\ref{c:path-parcel-variation-rear-pos}.
  The penultimate and last terms of $u(\zeta,r)$ are
  $\cF_{n_{\thetat}}(r)$ and
  $\cF_{n_{\thetat}\bop k_{\thetat}}(r)$. Thus,
  Claim~\ref{c:path-parcel-variation-rear-pos} holds by
  Claim~\ref{c:inc-dec-hill-inc-hill} of
  Proposition~\ref{prop:inc-dec-hill} and
  Theorem~\ref{thm:path-parcel-seq}.  We obtain
  Claims~\ref{c:path-parcel-variation-rear-zero}
  and~\ref{c:path-parcel-variation-rear-negative} similarly.
  
  Claims~\ref{c:path-parcel-variation-asym-pos}
  and~\ref{c:path-parcel-variation-asym-negative} follow
  from 
Claim~\ref{c:AfX-pos-med-qfact}
  of Proposition~\ref{prop:AfX-pos},
  Claim~\ref{c:path-parcel-ratio-lim-poly-nonneg} of
  Proposition~\ref{prop:path-parcel-ratio-unim-poly}, and
  Theorem~\ref{thm:path-parcel-seq}.  Also,
  Claim~\ref{c:path-parcel-variation-asym-zero} holds by
  Proposition~\ref{prop:path-parcel-asymp-hill-shape}.
\end{proof}

Second, 
we state the following lemma
for the real algebraic sets of
critical points.

\begin{lemma}
  \label{lem:front-rear-path-parcel-diff-poly}
  Consider $\cF= \Lam(s,l,w,\scc,\fs,\phi,\rho,x,\fX)$.  Let
  $r\in \AfX$ and $m,n\ldZl s$.
  Then, we have the following.
  \begin{enumerate}
  \item
    \label{c:front-rear-path-parcel-diff-poly-existence}
    $\Ir(\fX)(\cF_{m})_{1}(r)\neq 0$.
  \item
    \label{c:front-rear-path-parcel-diff-poly-diff}
    $(\cF_{m}-\cF_{n})(r)=0$ if and only if
    $\tdet(\Ir(\fX)(\cF_{m}),\Ir(\fX)(\cF_{n}))(r)=0$.
  \end{enumerate}
\end{lemma}
\begin{proof}
  Let us prove
  Claim~\ref{c:front-rear-path-parcel-diff-poly-existence}.
  Since $f_{s,m}$ is pairwise positive, we have
  $f_{s,m}^{2}\scc 0$.
  This implies $f_{s,m}^{2}\gAfX 0$.
  Also,
  $\prod \phi(x)^{m\rc w} \cdot [m]!_{x}^{w}\gAfX 0$
  by Proposition~\ref{prop:AfX-pos}.
  Then,
    $\cF_{m}^{2}\gAfX 0$ by
    $\cF_{m}=\frac{f_{s,m}}{ \prod \phi(x)^{m\rc w}\cdot [m]!_{x}^{w}}$
  and Lemma~\ref{lem:gAfX-closedness}.
  Hence, we obtain
  Claim~\ref{c:front-rear-path-parcel-diff-poly-existence}
  by Claim~\ref{c:geAfX-gAfX-scc-scc} of
  Lemma~\ref{lem:geAfX-gAfX-scc}.
  Claim~\ref{c:front-rear-path-parcel-diff-poly-diff}
  follows from
  Claim~\ref{c:front-rear-path-parcel-diff-poly-existence}.
\end{proof}

We also put the following  zero sets
 by polynomial determinants.

\begin{definition}\label{def:real-alg-sets}
  Suppose a $\tp{\theta,\chi,\fX}$-merged pair
  $\zeta=\tp{P,\cF}$ such that
  $\cF=\{ \cF_{i}\}_{i\ldZ \chi}$ and
  $P=\{\tp{s,l,m_{i},n_{i},k_{i}}\}_{i\ldZ \theta}$.  Let
  $\lam_{i}=\tp{P,\cF_{i}}$ for $i\ldZ\chi$.  Then, we define
  the following:
  \begin{enumerate} 
  \item
    $\FZ(\zeta)=\{r\in \AfX \mid \tdet(\Ir(\fX)
    (\cF_{i,m_{\thetao}\bom k_{\thetao}}),
    \Ir(\fX)(\cF_{i,n_{\thetao}}))(r)=0 \mforeach i\ldZ
    \chi\}$;
  \item
    $\RZ(\zeta)=\{r\in \AfX\mid
    \tdet(\Ir(\fX)(\cF_{i,n_{\thetat}}),\Ir(\fX)(
    \cF_{i,n_{\thetat}\bop k_{\thetat}}))(r)=0 \mforeach
    i\ldZ\chi\}$, if $P$ is finite-length;
  \item
    $\AZ(\zeta)=\{r\in \AfX \mid
    \tdet(\AIr(\fX)(\lam_{i}),\AIr(\fX)(t_{i},\lam_{i}))(r)=0
    \mforeach i\ldZ\chi\}$, if each $\cF_{i}$ is tame along
    $P$ by $t_{i}\in \Q(\fX)$.
  \end{enumerate}
\end{definition}

In Definition~\ref{def:real-alg-sets}, $\FZ(\zeta)$ ignore
choices of $\Q[\fX]$-polynomials
$\Ir(\fX)( \cF_{i,m_{\thetao}\bom k_{\thetao}} )$ and
$\Ir(\fX)(\cF_{i,n_{\thetao}})$, since
$ \Ir(\fX)( \cF_{i,m_{\thetao}\bom k_{\thetao}} )_{1}(r)\cdot
\Ir(\fX)(\cF_{i,n_{\thetao}})_{1}(r)\neq 0$ for each
$r\in \AfX$ by Claim~\ref{c:geAfX-gAfX-scc-geAfX} of
Lemma~\ref{lem:geAfX-gAfX-scc}.  The same holds for
$\RZ(\zeta)$ and $\AZ(\zeta)$.

We thus confirm the real algebraic sets  of critical
points.
\begin{theorem}\label{thm:variety}
  Consider a $\tp{\theta,\chi,\fX}$-merged pair
  $\zeta=\tp{P,\cF}$.  Then, $\FZ(\zeta)$, $\RZ(\zeta)$, and
  $\AZ(\zeta)$ are real algebraic sets of front, rear, and
  asymptotic critical points of $\zeta$, respectively.
  More explicitly,
  $r\in \AfX$ belongs to $\FZ(\zeta)$, $\RZ(\zeta)$, or
  $\AZ(\zeta)$ if and only if $\FD(\zeta)(r)=0$,
  $\RD(\zeta)(r)=0$, or $\AD(\zeta)(r)=0$, respectively.
\end{theorem}
\begin{proof}
  Statements hold by
  Claim~\ref{c:front-rear-path-parcel-diff-poly-diff} of
  Lemma~\ref{lem:front-rear-path-parcel-diff-poly},
  Claim~\ref{c:path-parcel-ratio-poly-exp} of
  Proposition~\ref{prop:path-parcel-ratio-unim-poly},
  and
  Theorem~\ref{thm:path-parcel-variation}.
\end{proof}

Third, we introduce the notion of {\it front, rear, and
  asymptotic semi-phase and phase transitions} of merged
pairs.  In
Section~\ref{sec:non-canonical-mediator-without-pt}, there
is a merged pair with a semi-phase transition but not with a
phase transition.  Also, in Section~\ref{sec:triv-no-crit},
there is a merged pair with asymptotic critical points but
not with a semi-phase transition.

\begin{definition}\label{def:phase-transitions}
  Suppose a $\tp{\theta,\chi,\fX}$-merged pair $\zeta=\tp{P,\cF}$.
  \begin{enumerate}
  \item We say that $\zeta$ has a front, rear, or asymptotic
    semi-phase transition, if
    $\emptyset \neq \FZ(\zeta)\neq \AfX$,
    $\emptyset \neq \RZ(\zeta)\neq \AfX$, or
    $\emptyset \neq \AZ(\zeta)\neq \AfX$, respectively.
  \item We say that $\zeta$ has a front, rear, or asymptotic
    phase transition between $r_{1}$ and $r_{2}$, if
    each $j\in \chi$
    satisfies
    $\FD(\zeta)_{j}(r_{1})\cdot \FD(\zeta)_{j}(r_{2})<0$,
    $\RD(\zeta)_{j}(r_{1})\cdot \RD(\zeta)_{j}(r_{2})<0$, or
    $\AD(\zeta)_{j}(r_{1})\cdot \AD(\zeta)_{j}(r_{2})<0$,
    respectively.
  \end{enumerate}
\end{definition}

For instance, suppose a $\tp{\theta,\chi,\fX}$-merged pair
$\zeta=\tp{P,\cF}$.  If there is a front phase transition of
$\zeta$ between $r_{1}$ and $r_{2}$, then
Theorem~\ref{thm:path-parcel-variation} gives a path
$L=\{L(u)\}_{1\leq u \leq 2}\subset \AfX$ between
$r_{1}=L(1)$ and $r_{2}=L(2)$ such that
$\FD(\zeta)_{j}(L(u))\in \R$ exists and changes the sign along
$L$ for each $j\in \chi$.  Also, we obtain the same for a rear
or asymptotic phase transition of $\zeta$ between $r_{1}$ and $r_{2}$ by $\RD(\zeta)$ or $\AD(\zeta)$, respectively.

\begin{remark}\label{rmk:phase-trans}
  Strictly decreasing sequences turn into two-slope
  hill-shape or strictly increasing sequences by front phase
  transitions, which give
  statistical-mechanical phase transitions in
  Sections~\ref{sec:intro-stat-pt} and~\ref{sec:casimir}.
\end{remark}

We obtain front phase transitions on monomial parcels
by  the lemma below.

\begin{lemma}\label{lem:BoTo-crit-points}
  Let $T\in \fX$.  Assume
  $\Bo(T),\To(T)\in \Q(\fX)^{2}$ with the
  following conditions:
  \begin{enumerate}
  \item     \label{cond:BoTo-gAfX}
    $\Bo_{1}(T),\Bo_{2}(T),\To_{1}(T),\To_{2}(T)
    \gAfX 0$;
  \item    \label{cond:BoTo-at-zero}  
    $\Qu(\Bo)(0)>\Qu(\To)(0)$ in $\R$;
  \item    \label{cond:BoTo-at-one}
    $\Qu(\Bo)(1)<\Qu(\To)(1)$ in $\R$.

  \end{enumerate}
  Also, let $F_{i}(T) =\frac{\To_{i}(T)}{\Bo_{i}(T)}$ for
  $i\in\oi(2)$.
  Then, there are real numbers
  $0<r_{0}<r_{1}<r_{2}<1$ such
  that   $F_{1}(r_{0})<F_{2}(r_{0})$,
  $F_{1}(r_{1})=F_{2}(r_{1})$, and
  $F_{1}(r_{2})>F_{2}(r_{2})$.
\end{lemma}
\begin{proof}
  Conditions~\ref{cond:BoTo-at-zero}
  and~\ref{cond:BoTo-at-one} give $0<r_{0}<r_{1}<r_{2}<1$
  such that $\Qu(\Bo)(r_{0})>\Qu(\To)(r_{0})$,
  $\Qu(\Bo)(r_{1})=\Qu(\To)(r_{1})$, and
  $\Qu(\Bo)(r_{2})<\Qu(\To)(r_{2})$.  Now,
  Condition~\ref{cond:BoTo-gAfX} implies
  the assertion.
\end{proof}
 
\begin{proposition}\label{prop:gen-width-monom-pt}
  Let $\tp{l,w,\gam}$ be a monomial index with palindromic
  $w$, $\tp{\gam_{i,1}}_{i\in \oi(l)}$, and
  $\tp{\gam_{i,2}}_{i\in \oi(l)}$.  For $s=\tp{0,\infty}$, assume
  a wrapped fitting path
  $P=\{\tp{s,l,m_{i},n_{i},k_{i}}\}_{i\ldZ \theta}$ and
  monomial
  parcel $\cF=\Lam(s,l,w,\scc,\Psi_{s,\gam,q},\rho,x,\fX)$.  Let
  $\zeta=\tp{P,\cF}$.  Then, we have the following.
  \begin{enumerate}
  \item \label{c:gen-width-monom-pt-merged-pair}
    $\zeta$ is a merged pair.
  \item \label{c:gen-width-monom-pt-equiv} If $q$ is
    fully
    admissible by $\fX$, then the following statements are
    equivalent.
    \begin{enumerate}
    \item \label{s:gen-width-monom-pt-equiv-semi-pt}
      $\zeta$ has a front semi-phase transition;
    \item \label{s:gen-width-monom-pt-equiv-ineq}
      $t_{\gam}(m_{\thetao})-t_{\gam}(m_{\thetao}\bom
      k_{\thetao})\in \Qgo;$
    \item \label{s:gen-width-monom-pt-equiv-full-pt}
      $\zeta$ has a front phase transition.
    \end{enumerate}
    Moreover, if $\fX=\{X_{1}\}$ and
    $t_{\gam}(m_{\thetao})-t_{\gam}(m_{\thetao}\bom
    k_{\thetao})\in \Qgo$, then $\zeta$ has a front phase
    transition at the unique front critical point.
  \end{enumerate}
\end{proposition}
\begin{proof}
  Let $\thetao=1$ and each $\gam_{i,3}=0$ for simplicity.
  Let $r\in \AfX$.  Let us prove
  Claim~\ref{c:gen-width-monom-pt-merged-pair}.  First,
  $\cF$ is static by Proposition~\ref{prop:static-mono}.
  Second, $\cF$ is $\gAfX$-merged-log-concave by
  Theorem~\ref{thm:monomial-poly}.  Third, $\cF$ is
  $\tp{s,l,\fX}$-positive, since
  $\cF_{\mn}(r) = \frac{q(r)^{t_{\gam}(\mn)}}
  {(\mn)^{w}_{q}|_{q=q(r)}} \in \R_{>0}$ by $0<q(r)<1$ and
  $\mn\ldZl s$.

  Let us prove Claim~\ref{c:gen-width-monom-pt-equiv}.  Let
  $T\in \fX$ and $h\in \Zgeo$ such that $T=q^{\frac{1}{h}}$
  for the fully admissible $q$.  First, we prove
  Statement~\ref{s:gen-width-monom-pt-equiv-ineq} from
  Statement~\ref{s:gen-width-monom-pt-equiv-semi-pt}.  Since
  $P$ is wrapped, $\mn\bom \kn\ldZl s$.  Also, since $\cF$
  is static, $\cF_{\nn}(r)=\cF_{\mn}(r)$.  Thus, there
  exists $c\in \AfX$ such that
\begin{align}
  \cF_{\mn\bom \kn}(c)
  =\frac{q(c)^{t_{\gam}(\mn\bom \kn)}}
  {(\mn\bom \kn)^{w}_{q}|_{q=q(c)}}
  =\frac{q(c)^{t_{\gam}(\mn)}}
  {(\mn)^{w}_{q}|_{q=q(c)}}
  =\cF_{\nn}(c).
  \label{eq:gen-width-monom-pt-eq}
\end{align}
Moreover, we obtain
$0<\frac{(\mn)^{w}_{q}|_{q=q(c)}}
{(\mn\bom\kn)^{w}_{q}|_{q=q(c)}} <1$ by $\sig(\kn)>0$ and $w>0$.
Hence, equation~\eqref{eq:gen-width-monom-pt-eq} implies
Statement~\ref{s:gen-width-monom-pt-equiv-ineq}, since
otherwise
$q(c)^{t_{\gam}(m_{\thetao})-t_{\gam}(m_{\thetao}\bom
  k_{\thetao})}\geq 1$.
  
Second, let us prove
Statement~\ref{s:gen-width-monom-pt-equiv-full-pt} from
Statement~\ref{s:gen-width-monom-pt-equiv-ineq}.  Let
$\Bo_{1}(T)=(\mn\bom \kn)^{w}_{q}$,
$\Bo_{2}(T)=(\mn)^{w}_{q}$,
$\To_{1}=q^{t_{\gam}(\mn\bom\kn)}$, and
$\To_{2}(T)=q^{t_{\gam}(\mn)}\in \Q(\fX)$.  Then,
Lemma~\ref{lem:BoTo-crit-points} gives a front phase
transition of $\zeta$.
  
Third, Statement~\ref{s:gen-width-monom-pt-equiv-full-pt}
implies Statement~\ref{s:gen-width-monom-pt-equiv-semi-pt}
by Definition~\ref{def:phase-transitions}.  The latter
statement holds by $\fX=\{T\}$.  Thus,
Claim~\ref{c:gen-width-monom-pt-equiv} holds.
\end{proof}

In Sections ~\ref{sec:critt-pts-pt}
and~\ref{sec:primal-monom}, we further discuss phase
transitions by width-one monomial parcels.

\subsection{Ideal merged pairs}
\label{sec:ideal}
Merged pairs concern real values of merged
determinants. Also, there is a $\tp{\theta,\fX}$-merged
pair $\tp{P,\cG}$ such that $\cG$ is not $\llq$-merged-log-concave, but its merged
determinants along $P$ give $q$-polynomials (see
Section~\ref{sec:ideal-prop-merged-pair}).  Thus, we
introduce the notion of {\it ideal merged pairs} to obtain
polynomials with positive integer coefficients.
This extends
Definition~\ref{def:intro-ideal}.

\begin{definition}\label{def:ideal}
  Suppose a $\tp{\theta,\chi,\fX}$-merged pair
  $\zeta=\tp{P,\cF}$ for
  $P=\{\tp{s,l,m_{i},n_{i},k_{i}}\}_{i\ldZ \theta}$ and
  $\cF=\{ \cF_{i}=
  \Lam(s,l,w_{i},\scc_{i},f_{i,s},\phi_{i},\rho_{i},x_{i},\fX)
  \}_{i\ldZ \chi}$.  We call $\zeta$ ideal, if
  \begin{align}
    \Delta(\cF_{i})
    (s,l,w_{i},m_{j},n_{j},k_{j},\phi_{i},\rho_{i},x_{i},\fX)
    >_{x_{i}}0
  \end{align}
  for each $i\ldZ \chi$ and $j\ldZ \theta$.  If
  $\chi_{1}=\chi_{2}$ and $\zeta$ is ideal, then we simply call
  $\tp{P,\cF_{\chi_{1}}}$ ideal.
\end{definition}
For example, suppose an ideal $\tp{\theta,\fX}$-merged
pair $\tp{P,\cF}$ such that
$\cF=
  \Lam(s,l,w,\scc,f_{s},\phi,\rho,x,\fX)$ and
   $x=\tp{q}$.  Then,
$\Delta(\cF)(s,l,w,m_{j},n_{j},k_{j},\phi,\rho,x,\fX)\llq
0$ for each $j\ldZ \theta$.  Moreover,
these $q$-polynomials
give almost strictly unimodal sequences $u(\tp{P,\cF},r)$
for each $r\in \AfX$ in Theorem~\ref{thm:path-parcel-seq}.

\subsection{Comparison of fitting paths}
\label{sec:comparison-fitting-paths}

We compare fitting paths by the following {\it finer fitting paths}.

\begin{definition}\label{def:fit-finer}
  Assume a fitting path
  $P_{i}=\{\tp{s_{i},l,m_{i,j},n_{i,j},k_{i,j}}\}_{j\ldZ 
    \theta_{i}}$ for each $i\in\oi(2)$.
  \begin{enumerate}
  \item We call $P_{1}$ equivalent to $P_{2}$, if
    $\theta_{1,2}-\theta_{1,1}=\theta_{2,2}-\theta_{2,1}$,
    and each $j\in\oi(0,\theta_{1,2}-\theta_{1,1})$ satisfies
        $m_{1,\theta_{1,1}+j}=m_{2,\theta_{2,1}+j}$ and
    $n_{1,\theta_{1,1}+j}=n_{2,\theta_{2,1}+j}$.
      We write
    $P_{1}\equiv P_{2}$, when $P_{1}$ is equivalent to $P_{2}$.
    
  \item We call $P_{1}$ finer than or equivalent to $P_{2}$,
    if any $i_{1},i_{2}\ldZ \thetat$ such that
    $i_{1}<i_{2}$ give $j_{1},j_{2}\ldZ \thetao$
    such that $j_{1}<j_{2}$ with
    $m_{1,j_{\lam}}=m_{2,i_{\lam}}$ and
    $n_{1,j_{\lam}}=n_{2,i_{\lam}}$
    for each $\lam\in\oi(2)$.  We write $P_{1}\vdash P_{2}$,
    when
    $P_{1}$ is finer than or equivalent to $P_{2}$.
  \end{enumerate}
\end{definition}
In particular, we discuss the existence of finest fitting paths.
We first prove the antisymmetricity of $\vdash$ on
the equivalence relation $\equiv$. This uses
the following binary relation.
\begin{definition}\label{def:lessdot}
  Suppose $m,m'\in \Zl$ for $l\in \Zgeo$.
  Let  $m\lessdot m'$, if 
  $m<m'$ or $m^{\ve}<m'$.
\end{definition}
\begin{lemma}\label{lem:lessdot-st-po}
  The binary relation $\lessdot$ is a strict partial order
  on $\Zl$.
\end{lemma}
\begin{proof}
  First, let us prove the irreflexivity.  If $m\in \Zl$, then
  $m<m$ does not hold. Also, $m<m^{\ve}$ does not hold
  either, because $m<m^{\ve}$ implies
  $m_{1}<m^{\ve}_{1}=m_{l}$
  and $m_{l}<m^{\ve}_{l}=m_{1}$.

  Second, let us prove the transitivity.  Let
  $m_{1},m_{2},m_{3}\in \Zl$ such that
  $m_{1}\lessdot m_{2}\lessdot m_{3}$.  This means
  $m_{1}<m_{2}$ or $m_{1}<m_{2}^{\ve}$, and
  $m_{2}<m_{3}$ or
  $m_{2}<m_{3}^{\ve}$.  Thus, $m_{1}<m_{3}$ or
  $m_{1}<m_{3}^{\ve}$, since taking ${}^{\ve}$
  preserves $<$.
\end{proof}

Then, we have the antisymmetricity.
\begin{proposition}\label{prop:vdash-poset}
  Consider fitting paths
  $P_{i}= \{\tp{s_{i},l,m_{i,j},n_{i,j},k_{i,j}}\}_{j\ldZ
    \theta_{i}}$ for $i\in\oi(2)$.  If
  $P_{1}\vdash P_{2}$ and $P_{2}\vdash P_{1}$, then
  $P_{1}\equiv P_{2}$. In particular,
  $\vdash$ is a partial
  order on fitting paths on the equivalence relation $\equiv$.
\end{proposition}
\begin{proof}
  By $P_{1}\vdash P_{2}$, there is an order-preserving map
  $f:\{i: i\ldZ \thetat\}\mapsto \{i: i\ldZ\thetao\}$ such
  that $m_{1,f(i)}=m_{2,i}$ for each $i\ldZ\thetat$.  Also,
  by $P_{2}\vdash P_{1}$, there is an order-preserving map
  $g:\{i: i\ldZ\thetao\}\mapsto \{i:i\ldZ \thetat\}$ such
  that $m_{2,g(i)}=m_{1,i}$ for each $i\ldZ\thetao$.
  
  First, suppose $f(i)<u<f(i+1)$
  for some $u\ldZ\thetao$.  Then,
  $m_{1,u}\lessdot m_{1,f(i+1)}=m_{2,i+1}$ by
  Lemma~\ref{lem:fitting-path-mn-positivity}.  Hence,
  $g(u)<i+1$; otherwise,
    $m_{2,i+1}\lessdot m_{2,g(u)}=m_{1,u}$
    by Lemma~\ref{lem:fitting-path-mn-positivity}
    or
    $m_{2,i+1}= m_{2,g(u)}=m_{1,u}$,
    both of which
    are against
    Lemma~\ref{lem:lessdot-st-po}.  Likewise,
  $m_{2,i}=m_{1,f(i)}\lessdot m_{1,u}$ gives $i<g(u)$, which
  can not occur by $g(u)\in \Z$. Hence, $f(i+1)-f(i)=1$.
  
  Second, assume $f(\theta_{2,1})>\theta_{1,1}$.  Then,
  Lemmas~\ref{lem:fitting-path-mn-positivity}
  and~\ref{lem:lessdot-st-po} imply
  $g(\theta_{1,1})<\theta_{2,1}$, which can not happen.
  Thus, $f(\theta_{2,1})=\theta_{1,1}$. Hence,
  $P_{1}\equiv P_{2}$.
  
  We obtain the latter statement by the reflexivity and
  transitivity of $\vdash$.
\end{proof}

We also put the  notion of {\it restricted fitting paths}.

\begin{definition}\label{def:rest-path}
  Suppose a fitting path  $P=\{P_{i}\}_{i\ldZ \thetao}$.
  For a gate $\thetat$ such that
  $\theta_{1,1}\leq \theta_{2,1}\leq \theta_{2,2}\leq \theta_{1,2}$, we define the
  restricted fitting path
  $r_{\thetao,\thetat}(P)= \{P_{i}\}_{i\ldZ \thetat}$.
\end{definition}

Then, we write width-one fitting paths
by $\pi(\lam,h)$.
\begin{lemma}\label{lem:finest}
  For $l=1$, consider a fitting path
  $P=\{\tp{s,l,m_{i},n_{i},k_{i}}\}_{i\ldZ \thetao}$. Then, we have
  $P\equiv r_{\thetao,\thetat}(\pi(\lam,h))$
  for
  $\tp{\lam}=\sig(k_{\theta_{1,1}})\in \Zl$,
  $\tp{h}=m_{\theta_{1,1}}-\sig(k_{\theta_{1,1}})\in \Zl$, and
  $\thetat=\tp{1,\theta_{1,2}-\theta_{1,1}+1}$. 
\end{lemma}
\begin{proof}
  By $l=1$, each $m_{i}=n_{i}$.  Thus, the assertion
  follows, since equivalent supports $k_{i}$ give 
  $m_{i}=\tp{\lam (i-\theta_{1,1}+1)+h}\in \Zl$ for
  $i \ldZ \thetao$.
\end{proof}

We now characterize $\pi(1)$
as the finest fitting
path for width-one wrapped fitting paths.
\begin{proposition}\label{prop:existence-finest-width-one}
  We have the following.
  \begin{enumerate}
  \item \label{c:finest-vdash} If $P$ is a width-one
    wrapped
    fitting path, then $\pi(1)\vdash P$.
  \item\label{c:finest-unique}
    If
    there exists a wrapped fitting path $P'$
    such that
    $P' \vdash P$ for each
    width-one wrapped fitting path $P$,
    then $\pi(1)\equiv P'$.
  \end{enumerate}
\end{proposition}
\begin{proof}
  Let us prove Claim~\ref{c:finest-vdash}.  
  For some gates $\thetao$ and $\thetat$,
  $P\equiv r_{\thetao,\thetat}(\pi(\lam,h))$ 
  by Lemma~\ref{lem:finest}. Hence,
  Claim~\ref{c:finest-vdash} holds, since
   $h \in \Zgez$ for the wrapped $P$.

  Claim~\ref{c:finest-unique} holds by
  Proposition~\ref{prop:vdash-poset} and
  Claim~\ref{c:finest-vdash}.
\end{proof}

We also prove the non-existence of finest fitting
paths for higher-width wrapped ones by  the following lemma.
\begin{lemma}\label{prop:path-eq-vdash}
  Assume fitting paths
  $P= \{\tp{s,l,m_{i},n_{i},k_{i}}\}_{i\ldZ \theta}$.  Let
  $u_{1},u_{2}\ldZ \theta$ such that $u_{1}<u_{2}$.  Then, there
  is $\mu\in \Zgez^{2}$ such that $\muo\geq 1$,
  $\sum \mu= u_{2}-u_{1}$, and
    \begin{numcases}{
        \mu_{1} \sig(k_{u_{1}})+
        \mu_{2} \sig(k_{u_{1}})^{\ve}
        + n_{u_{1}}^{\ve}
        =}
      m_{u_{2}}, \mor
      \label{eq:path-eq-vdash-gen-diff-no-ve}
      \\
      m_{u_{2}}^{\ve}.
      \label{eq:path-eq-vdash-gen-diff-with-ve}
    \end{numcases}
\end{lemma}
\begin{proof} 
  Since $(n_{u_{2}-1}\bop k_{u_{2}-1})^{\ve}
  =m_{u_{2}}$,
  we have
  \begin{align}
    \sig(k_{u_{2}-1}) + n_{u_{2}-1}^{\ve} &=
                                        m_{u_{2}}.
    \label{eq:path-eq-vdash-gen-first-with-ind}
  \end{align}
  Hence, we prove the assertion by the induction.
  Let $u_{2}-u_{1}\geq 2$. Also, assume
  $\lam\in \Zgez^{2}$ such that $\lam_{1}\geq 1$,
  $\sum \lam=u_{2}-u_{1}-1$, and
  \begin{numcases}{
      \lam_{1} \sig(k_{u_{1}})+
      \lam_{2} \sig(k_{u_{1}})^{\ve}
      + n_{u_{1}}^{\ve}
      =}
    m_{u_{2}-1}, \mor
    \label{eq:path-eq-vdash-gen-diff-no-ve-ind}
    \\
    m_{u_{2}-1}^{\ve}.
    \label{eq:path-eq-vdash-gen-diff-with-ve-ind}
  \end{numcases}

  First, suppose
  equation~\eqref{eq:path-eq-vdash-gen-diff-no-ve-ind}.  If
  $n_{u_{2}-1}^{\ve}= m_{u_{2}-1}$, then we obtain
  equation~\eqref{eq:path-eq-vdash-gen-diff-no-ve}, adding
  equations~\eqref{eq:path-eq-vdash-gen-first-with-ind}
  and~\eqref{eq:path-eq-vdash-gen-diff-no-ve-ind}.
Moreover, by equation~\eqref{eq:path-eq-vdash-gen-first-with-ind},
we have
\begin{align}
  \sig(k_{u_{2}-1})^{\ve} + n_{u_{2}-1}=
  m_{u_{2}}^{\ve}.
  \label{eq:path-eq-vdash-gen-first-with-ind-ve}
\end{align}
Hence,
$n_{u_{2}-1}= m_{u_{2}-1}$ gives
  equation~\eqref{eq:path-eq-vdash-gen-diff-with-ve} by  equations~\eqref{eq:path-eq-vdash-gen-diff-no-ve-ind}
  and~\eqref{eq:path-eq-vdash-gen-first-with-ind-ve}.
  
  Second, suppose
  equation~\eqref{eq:path-eq-vdash-gen-diff-with-ve-ind}.
  If $n_{u_{2}-1}= m_{u_{2}-1}^{\ve}$, then
  equation~\eqref{eq:path-eq-vdash-gen-diff-with-ve} follows
  from
  equations~\eqref{eq:path-eq-vdash-gen-diff-with-ve-ind}
  and~\eqref{eq:path-eq-vdash-gen-first-with-ind-ve}.
  Also,
  $n_{u_{2}-1}^{\ve}= m_{u_{2}-1}^{\ve}$ gives
  equation~\eqref{eq:path-eq-vdash-gen-diff-no-ve} 
  by equations~\eqref{eq:path-eq-vdash-gen-first-with-ind}
  and~\eqref{eq:path-eq-vdash-gen-diff-with-ve-ind}.
\end{proof}

\begin{proposition}
  \label{prop:non-existence-finest-higher-widths}
  If $l\geq 2$,  then there is no width-$l$
  fitting path
  $P$ such that
  $P\vdash P'$ for each width-$l$
  wrapped fitting path $P'$.
\end{proposition}
\begin{proof}
  Let us prove by contradiction. Hence, suppose a fitting
  path $P= \{\tp{s,l,m_{i},n_{i},k_{i}}\}_{i\ldZ \theta}$ such
  that $P\vdash P'$ for each
  width-$l$ wrapped fitting path $P'$.
  
  Let $\xi=\tp{l}\in Y(1,l)$, $s'=\tp{0,\infty}$,
  and $\theta'=\tp{1,\infty}$.  Then, there is a fitting
  path
  $P'= \bP_{s,l,\xi,0}
  =\{\tp{s',l,m'_{i},n'_{i},k'_{i}}\}_{i\ldZ \theta'}$ such
  that $m'_{i}=n'_{i}=i \iota^{l}(1)$ by
  Theorem~\ref{thm:fitting-path-young} and
  Lemma~\ref{lem:low-inclusion}.
  Furthermore, $P\vdash P'$ implies $m'_{1}=m_{i_{1}}$ and
  $m'_{2}=m_{i_{2}}$ for some $i_{1},i_{2}\ldZ \theta$ such
  that $i_{1}<i_{2}$.  Then, 
  Proposition~\ref{prop:path-eq-vdash} gives
  $\mu\in \Zgez^{2}$ such that $\sum \mu=i_{2}-i_{1}$,
  $\mu_{1}\geq 1$, and
  $\mu_{1}\sig(k_{1})+\mu_{2}\sig(k_{1})^{\ve}
  =m_{i_{2}}-m_{i_{1}}
  =m'_{2}-m'_{1}=\iota^{l}(1)$.  Hence,
$\sig(k_{1})>0$ implies
  $i_{2}-i_{1}=1$ and
  $\sig(k_{1})=\iota^{l}(1)$.  Also, we
  obtain flat $m_{j}=n_{j}=(j-i_{1}+1)\iota^{l}(1)$ for
  $j\ldZ \theta$.
  
  However,
 by
  Theorem~\ref{thm:fitting-path-young} and
  Lemma~\ref{lem:low-inclusion},
  $\xi=\tp{l,1}\in Y(2,l)$ gives a
  wrapped fitting path
  $P'= \bP_{s,l,\xi,3}
  =\{\tp{s',l,m'_{i},n'_{i},k'_{i}}\}_{i\ldZ \theta'}$ such
  that
  $m'_{i}=n'_{i}=\tp{5}\ccn \iota^{l-1}(4)+(i-1)\iota^{l}(3)$
  for $i\ldZ \theta'$.  Then, $P\vdash P'$ can not
  hold, since each $m_{i}$ of $P$ is flat.
\end{proof}

By Propositions~\ref{prop:existence-finest-width-one}
and~\ref{prop:non-existence-finest-higher-widths}, we call
$\pi(1)$ finest fitting path for our convenience.

\subsection{On phase transitions
  of  width-one monomial parcels}
\label{sec:critt-pts-pt}


We obtain some phase transitions of width-one monomial
parcels by canonical mediators and the finest fitting path
$\pi(1)$.  This uses the following lemma on tame factors.

\begin{lemma}\label{lem:mono-tame}
  Let $l=1$ and $\mu=\tp{l,w,\gam}$ be a monomial index.
  Consider the $t$-monomials $\Psi_{s,\gam,q}$.  Assume
  $r\in \AfX$ such that $0<q(r)<1$. Then, for $i\in \Zgez$, we
  have
  \begin{align}
    \lim_{i\to\infty}\frac{\Psi_{s,\gam,q(r),\tp{i+1}}}{
    \Psi_{s,\gam,q(r),\tp{i}}}
    =\begin{dcases}
       q(r)^{\gam_{1,2}} &\mif \gam_{1,1}=0,\\
       0 &\melse.
     \end{dcases}
  \end{align}
\end{lemma}
\begin{proof}
  The monomial conditions of $\mu$ yield
  $\gam_{1,1}\geq 0$.   Hence, the
  assertion holds by
  \begin{align}
    \lim_{i\to\infty}\frac{\Psi_{s,\gam,q(r),\tp{i+1}}}{
    \Psi_{s,\gam,q(r),\tp{i}}}
    =
    \lim_{i\to\infty}\frac{q(r)^{\gam_{1,1}(i+1)^{2}+\gam_{1,2}(i+1)}}{
    q(r)^{\gam_{1,1}i^2+\gam_{1,2}i}}
    =  q(r)^{\gam_{1,2}}   \lim_{i\to\infty}q(r)^{(2i+1)\gam_{1,1}}.
  \end{align}
\end{proof}
We then state the phase transitions.

\begin{proposition}\label{prop:monom-front-pt}
  Suppose a monomial parcel $\cF=\Lam(s,l,w,\scc,\Psi_{s,\gam,q},\phi,x,\fX)$ such that
  $s=\tp{0,\infty}$, $l=1$, and
  $\phi(x)=\tp{1-q}\in \Q(\fX)^{l}$.  Then, for the merged pair
  $\zeta=\Omg(\cF)$, we have the following.
  \begin{enumerate}
  \item \label{c:monom-front-pt-no-rear}
    $\zeta$ has no rear critical points.
  \item \label{c:monom-front-pt-no-asymp}
    $\zeta$ has no asymptotic critical points.
  \item \label{c:monom-front-pt-equiv} If $q$ is fully
    admissible by $\fX$, then the following statements
    are equivalent:
    \begin{enumerate}
    \item $\zeta$ has a front semi-phase transition;
    \item $\gam_{1,1}+\gam_{1,2}>0$;
    \item $\zeta$ has a front phase transition.
    \end{enumerate}
    Moreover, if $\fX=\{X_{1}\}$ and
    $\gam_{1,1}+\gam_{1,2}>0$, then $\zeta$ has a single phase
    transition at the front critical point.
  \end{enumerate}
\end{proposition}
\begin{proof}
  Claim~\ref{c:monom-front-pt-no-rear} holds, since $\pi(1)$
  is infinite-length.
  
  Let us prove Claim~\ref{c:monom-front-pt-no-asymp}.  Let
  $r\in \AfX$. Then, Claim~\ref{c:adm-bounds} of
  Lemma~\ref{lem:adm-bounds-poring-inclusions} implies
  $0<q(r)<1$.  Hence, by Lemma~\ref{lem:mono-tame}, $\cF$ is
  tame along $P$ by $t\in \Q(\fX)$ such that
  $t(r)= q(r)^{\gam_{1,2}}$ if $ \gam_{1,1}=0$, and $0$
  otherwise.  Let
  $\pi(1) =\{\tp{s,l,m_{i},n_{i},k_{i}}\}_{i\ldZ \theta}$.  Then,
  $\sig(k_{\thetao})=\tp{1}\in \Zl$ yields
  $\AD(\zeta) = \prod \phi(x)^{\sig(k_{\thetao})^{\ve}\rc w} - t \cdot \prod
  \phi(x)^{\sig(k_{\thetao})^{\ve}\rc w} =(1-q)^{w_{1}}(1-t)$.
  Thus, $\AD(\zeta)=0$ has no solutions over $\AfX$.

  The equivalence of
  Claim~\ref{c:monom-front-pt-equiv} follows from
  Claim~\ref{c:gen-width-monom-pt-equiv} of
  Proposition~\ref{prop:gen-width-monom-pt}, since
  $t_{\gam}(m_{\thetao})-t_{\gam}(m_{\thetao}\bom
  k_{\thetao})=\gam_{1,1}+\gam_{1,2}$.
  Then,
Claims~\ref{c:monom-front-pt-no-rear}
  and~\ref{c:monom-front-pt-no-asymp} give
  the latter statement of
  Claim~\ref{c:monom-front-pt-equiv}.
\end{proof}

\section{Explicit critical
  points,   phase transitions, and
  merged determinants}
\label{sec:examples} 
 We
 use the following notation for simplicity.
\begin{definition}
  Consider a $\tp{\theta,\chi}$-merged pair $\zeta=\tp{P,\cF}$ 
  such that
  $\fX=\{X_{1}\}$.  Let
  $\tp{h}\in \AfX$ for 
  a real number $0<h<1$.
  Then, we write $u(\zeta,h)$ for the path-parcel sequence
    $u(\zeta,\tp{h})$.
  Also, we call $h$ front, rear, or asymptotic
  critical point, if $\tp{h}$ is a front, rear, or asymptotic
    critical point of $\zeta$, respectively.
  \end{definition}
 We often write
 $\FC(\zeta)$, $\RC(\zeta)$, and $\AC(\zeta)$ to indicate a front,
rear, or asymptotic critical point $h$, respectively.

We put the following  parcel $\cL$, which appears in
Section~\ref{sec:intro-merged-rational-functions} by a
different notation.
In particular,  $\cL$ is $\scc$-merged-log-concave
and  fully optimal.

\begin{definition}\label{def:linear-half}
  Let $s=\tp{0,\infty}$, $l=1$, $w=\tp{1}$, and $\fX=\{\qq\}$.
  Then, we define the linear-half monomial parcel
  $\cL=\Lam(s,l,w,\scc,\Psi_{s,\tp{\tp{0,\fraa,0}},q},x,\fX)$.
\end{definition}

\subsection{Golden angle as a critical point}
\label{sec:golden-angle}

For $l=1$, $w=\tp{2}$, $s=\tp{0,\infty}$, and $\fX=\{q\}$,
let
$\cF=\Lam(s,l,w,\llq,\Psi_{s,\tp{\tp{0,1,0}},q},x,\fX)$ such that
$\cF_{\lam} =\frac{q^{\lamo}}{(\lamo)^{2}_{q}}$ for
$\lam\ldZl s$.
In particular,  $\fX$ is 
the
fully optimal coordinate for $\cF$.
Then, by Theorem~\ref{thm:monomial-poly} and
Proposition~\ref{prop:gen-width-monom-pt}, we consider
the ideal merged pair $\zeta=\Omg(\cF)$ for a critical point
and phase transition, since
 the path-parcel sequences
$u(\zeta,r)$
of $r\in \AfX$ are
almost
strictly unimodal 
by
Theorem~\ref{thm:path-parcel-seq}.

\subsubsection{On critical points and phase transitions}
\label{sec:golden-angle-crit}
 Let
$\pi(1)=\{\tp{s,l,m_{i},n_{i},k_{i}}\}_{i\ldZ \theta}$ with
$m_{0}=n_{0}=\tp{0}$ and $k_{0}=a_{0}=b_{0}=\tp{0,1}$.
Then, a real number $0<q<1$ is a front critical point of
$\zeta$ if and only if
\begin{align}
  \cF_{m_{0}}(q)=1=\frac{q}{(1-q)^2}=\cF_{\mn}(q).  
\end{align}
Therefore,  the unique front critical
point
$\FC(\zeta)$ of $\zeta$ is the golden angle: i.e.,
\begin{align}
  \FC(\zeta) =\frac{3-\sqrt{5}}{2}=0.381966\dots.
\end{align}

\begin{figure}[H]
  \centering
  \resizebox{130mm}{!}   
  {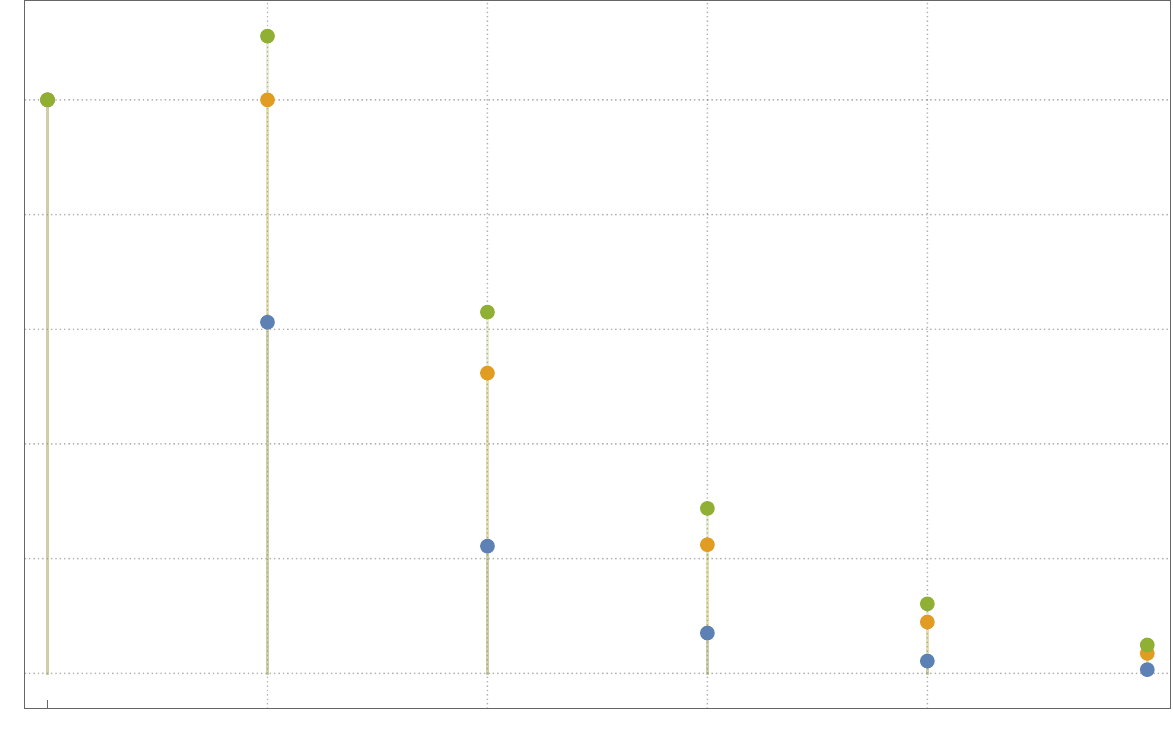}  
  \caption{$\cF_{m_{i}}(q)$ of
     $q=0.3$ (bottom), $\FC(\zeta)$ (middle),
    and $0.4$ (top)}
  \label{fig:golden-angle}
\end{figure}

By Proposition~\ref{prop:monom-front-pt}, $\zeta$ has neither
rear nor asymptotic critical points. However, $\zeta$ has the front
phase transition at $\FC(\zeta)$ such that $u(\zeta,0.3)$ is
strictly decreasing, $u(\zeta,\FC(\zeta))$ is decreasing and
hill-shape, and $u(\zeta, 0.4)$ is two-slope hill-shape.  For
each $i\in \oi(0,5)$, Figure~\ref{fig:golden-angle} plots the
bottom point for $\cF_{m_{i}}(0.3)$, the middle point for
$\cF_{m_{i}}(\FC(\zeta))$, and the top point for
$\cF_{m_{i}}(0.4)$.  \vspace{1em}

\subsubsection{Polynomials with positive integer
  coefficients of an ideal merged pair}
\label{sec:golden-angle-poly}
For $i\in \Zgeo$,
$\zeta$ gives the following $q$-polynomials with positive
coefficients:
\begin{align}
  \Delta(\cF)(s,l,w,m_{i},n_{i},k_{i},x,\fX)
  &=
    \frac{(i)_{q}^{2}\cdot (i+1)_{q}^2}{(0)_{q}^{2}\cdot (1)_{q}^2}
    \cdot 
    ( \cF_{1,m_{i}} \cF_{1,n_{i}}    -  \cF_{1,m_{i-1}}
    \cF_{1,n_{i+1}})\\
  &=
    \frac{(i)_{q}^{2}\cdot (i+1)_{q}^2}{(0)_{q}^{2}\cdot (1)_{q}^2}
    \cdot 
    \left(
    \left(\frac{q^{i}}{(i)^{2}_{q}}\right)^{2}
    -\frac{q^{i-1}}{(i-1)^{2}_{q}}
    \frac{q^{i+1}}{(i+1)^{2}_{q}}
    \right).
\end{align}
Explicitly, some of them are 
\begin{align}
  \Delta(\cF)(s,l,w,\mn,\nn,\kn,x,\fX)
  &=2 q^{3}+q^{4},\\
  \Delta(\cF)(s,l,w,\mt,\nt,\kt,x,\fX)
  &=2 q^{6}+2 q^{7}+q^{8},\\
  \Delta(\cF)(s,l,w,\mr,\nr,\kr,x,\fX)
  &=2 q^{9}+2 q^{10}+2 q^{11}+q^{12},\\
  \Delta(\cF)(s,l,w,m_{4},n_{4},k_{4},x,\fX)
  &=2 q^{12}+2 q^{13}+2 q^{14}+2 q^{15}+q^{16},\\
  \Delta(\cF)(s,l,w,m_{5},n_{5},k_{5},x,\fX)
  &=2 q^{15}+2 q^{16}+2 q^{17}+2 q^{18}+2 q^{19}+q^{20}.
\end{align}

\subsubsection{Golden angle from golden ratio
  as critical points}
\label{sec:golden-angle-ratio}
The linear-half monomial parcel $\cL$ has the golden ratio
as the unique front critical point of $\Omg(\cL)$ in
Section~\ref{sec:intro-explicit-crit-pts}.
Moreover, $\cL$
gives $\cF$, because
\begin{align}
  \cL^{\dd 2}
  &=\Lam(s,l,w,\scc,\Psi_{s,\tp{\tp{0,1,0}},q},x, \{\qq\}),
\end{align}
and 
$\cF= r_{\tp{1},\tp{2}} ( \cL^{\dd 2})$  by the parcel restriction
$r_{\tp{1},\tp{2}}$ (or $\cF=\cL^{\dd 2}$ in $\Q(\qq)$).

\subsection{A non-canonical mediator with
  phase transitions}
\label{sec:non-canonical-mediator-with-pt} 
Consider $s=\tp{0,\infty}$, $l=1$, $w=\tp{2}$, $\rho=\tp{1}$, and
$\fX=\{q\}$.
Also,
let
\begin{align}
  \phi(x)=\tp{\frac{29}{20}-5q+5q^2}\in \Q(x)^{l}.
\end{align}
This $\phi$ is a non-canonical
$\tp{s,l,w,\llq,\rho,x,\fX}$-mediator,
because $\phio(x_{1})^{w_{1}}\gAfX 0$ by
\begin{align}
  \frac{29}{20}-5q+5q^2=
  5\left(q-\fraa\right)^{2}+\frac{1}{5},
  \label{eq:non-canon-quad-med}
\end{align}
and $B(s,l,w,m,\phi,\rho,x,\fX)\llq 0$ for each $m\ldZl s$ by
$\rho=\tp{1}$ and Lemma~\ref{lem:bshift-fun-sp}.  Thus, the analogs
$\phio(x_{1})^{\lamo} [\lamo]!_{q}$ of $q$-Pochhammer
symbols $(\lamo)_{q}$ give
the monomial parcel
$\cF=\Lam(s,l,w,\llq,\Psi_{s,\tp{\tp{0,1,0}},q},\phi,\rho,x,\fX)$
such that each $\lam\ldZl s$ satisfies
\begin{align}
  \cF_{\lam}
  =\frac{q^{\lamo}}{
  \left(\frac{29}{20}-5q+5q^2\right)^{2\lamo}\cdot
   [\lamo]!_{q}^{2}}.
\end{align}
This $\cF$ is $\llq$-merged-log-concave by
Theorem~\ref{thm:monomial-poly}.  Hence, we have the ideal
merged pair $\zeta=\Omg(\cF)$, which gives the same merged
determinants in Section~\ref{sec:golden-angle} by
Proposition~\ref{prop:mediator-equiv}.

\subsubsection{On critical points and
  phase transitions}
\label{sec:non-canonical-crit}
Let $P=\pi(1)=\{\tp{s,l,m_{i},n_{i},k_{i}}\}_{i\ldZ \theta}$ with
$m_{0}=n_{0}=\tp{0}$ and $k_{0}=a_{0}=b_{0}=\tp{0,1}$.
First, unlike the canonical mediator in
Section~\ref{sec:golden-angle}, the non-canonical $\phi$ of
$\zeta$ gives exactly the two front critical points by
equation~\eqref{eq:non-canon-quad-med}.
More explicitly, we have the
front critical points $\FC(\zeta)_{1}=0.253594\dots$ and
$\FC(\zeta)_{2}=0.884885\dots$ that solve
\begin{align}
  \cF_{m_{0}}(q)
  =1=
  \frac{q}
  {\left(\frac{29}{20}-5q+5q^2\right)^2}=\cF_{\mn}(q).
\end{align}
Second, there are no rear critical points for the
infinite-length $P$.  Third, $\cF$ is tame along $P$ by $q$
in Lemma~\ref{lem:mono-tame}. Thus, an asymptotic critical
point $0<q<1$ solves
\begin{align}
  \AD(\zeta)
  =\left(\frac{29}{20}-5q+5q^2\right)^2-q(1-q)^2=0.
\end{align}
For $v(q)=q(1-q)^2$, $v'(q)=(1-q)(1-3q)$ and
$v\left(\frac{1}{3}\right)
>\phio\left(\frac{1}{3}\right)^{2}$.
Then,  equation~\eqref{eq:non-canon-quad-med}
gives  the two asymptotic critical points 
$\AC(\zeta)_{1}=0.30813\dots$ and   
$\AC(\zeta)_{2}=0.63486\dots$.
Hence, we obtain the
following table of phase transitions.
\begin{align}
  \begin{array}{|c|c|c|c|c|c|c|c|c|c|c|c|c|c|c|}
    \hline
    & 0&  &\FC(\zeta)_{1}  & & \AC(\zeta)_{1} &
    & \AC(\zeta)_{2} &   & \FC(\zeta)_{2} &  & 1 \\ \hline
    \FD(\zeta)  &+& + & 0 & - & - & - & - & -& 0 & +& +\\ \hline
    \AD(\zeta) &+& + &+ & +& 0 & - & 0 &+ & + & +& + \\ \hline
  \end{array}
  \label{table:pt-non-canon-quad-med}
\end{align}
In particular, the non-canonical mediator gives not only the front
phase transitions at $\FC(\zeta)_{1}$ and $\FC(\zeta)_{2}$, but
also the asymptotic phase transitions at $\AC(\zeta)_{1}$ and
$\AC(\zeta)_{2}$ between the front phase transitions.

\begin{figure}[H]
  \centering 
  \resizebox{130mm}{!}   
  {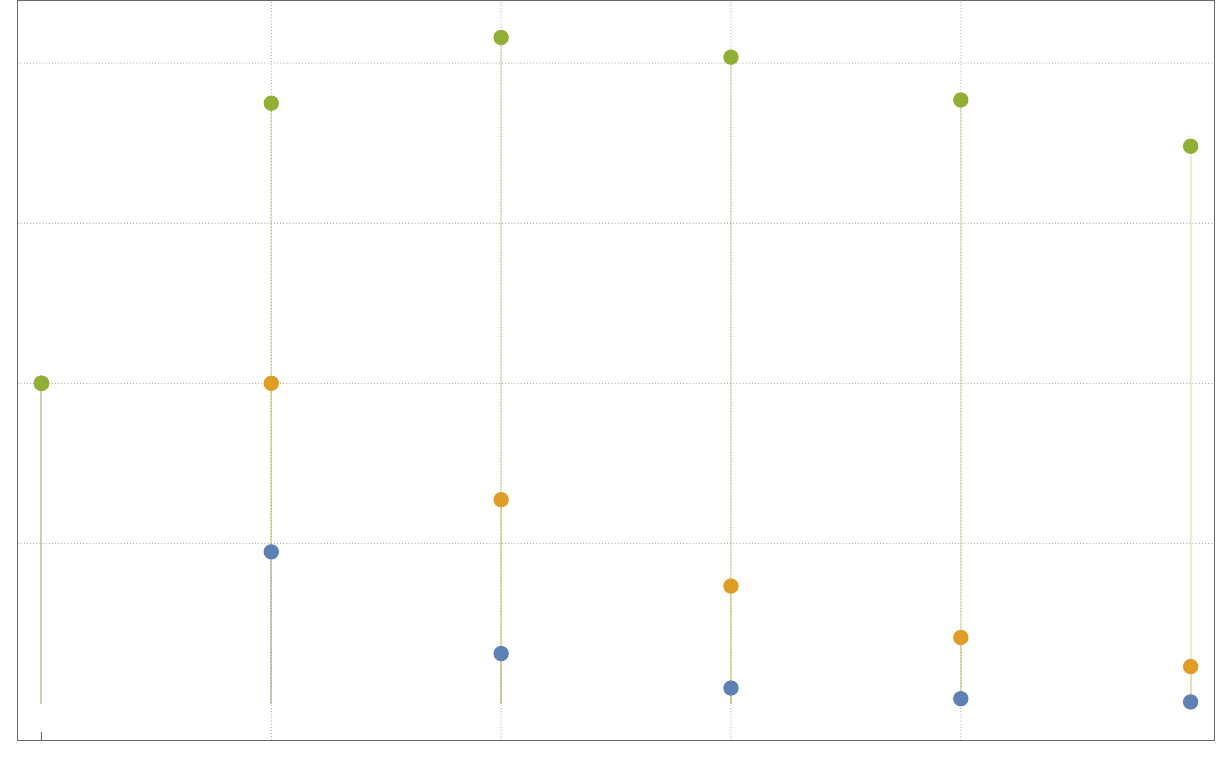}
  \caption{$\cF_{m_{i}}(q)$ of
     $q=0.2$ (bottom),
    $\FC(\zeta)_{1}$ (middle), and $0.3$ (top)}
  \label{fig:non-canonical-med-left-front-pt}
\end{figure}

First, the front phase transition at $\FC(\zeta)_{1}$ gives the
strictly decreasing $u(\zeta,0.2)$, decreasing and hill-shape
$u(\zeta,\FC(\zeta)_{1})$, and two-slope hill-shape
$u(\zeta,0.3)$. For each $i\in \oi(0,5)$,
Figure~\ref{fig:non-canonical-med-left-front-pt} plots
$\cF_{m_{i}}(q)$ of $q=0.2$ for the bottom point,
$q=\FC(\zeta)_{1}$ for the middle point, and $q=0.3$ for the
top point.
\vspace{1em}

\begin{figure}[H]
  \centering
  \resizebox{130mm}{!}   
  {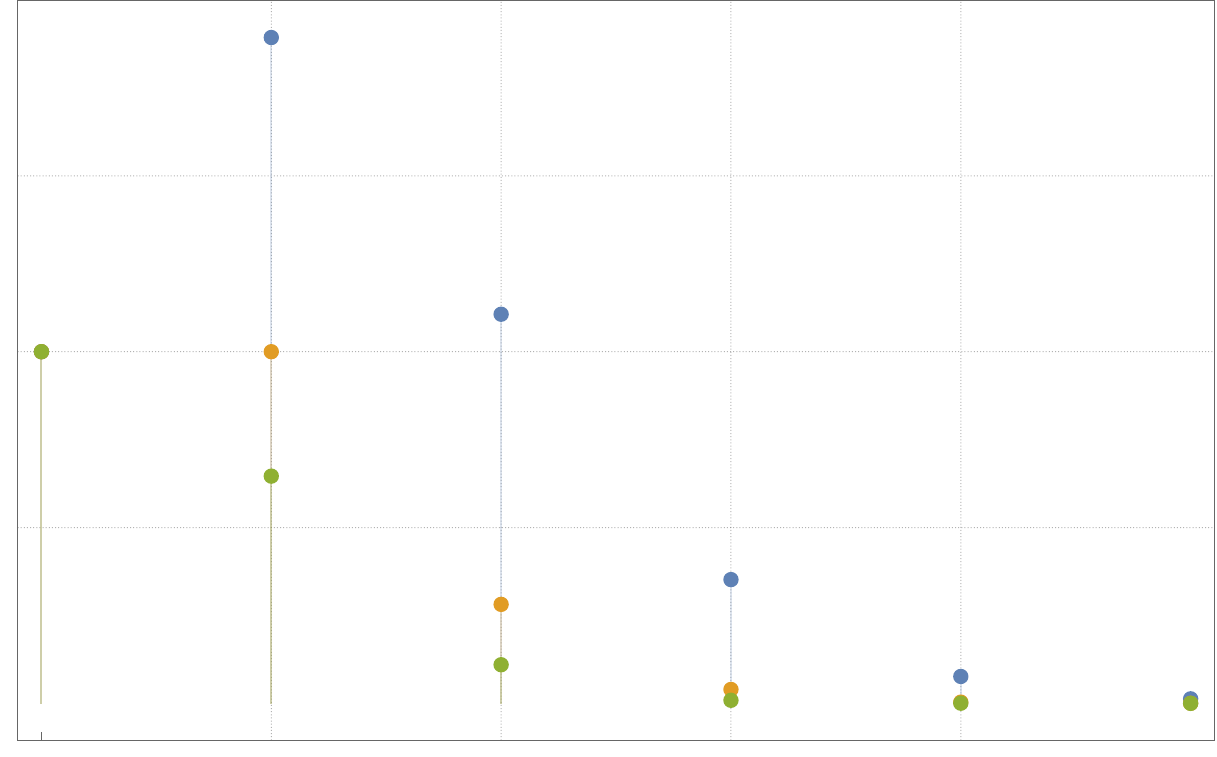}
  \caption{$\cF_{m_{i}}(q)$ of
     $q=0.95$ (bottom),
    $\FC(\zeta)_{2}$ (middle), and $0.8$ (top)}
  \label{fig:non-canonical-med-right-front-pt}
\end{figure}

Second, the front phase transition at $\FC(\zeta)_{1}$ gives
the two-slope hill-shape $u(\zeta,0.8)$, decreasing and
hill-shape $u(\zeta,\FC(\zeta)_{2})$, and strictly decreasing
$u(\zeta,0.95)$.
Figure~\ref{fig:non-canonical-med-right-front-pt} puts
$\cF_{m_{i}}(q)$ of $q=0.95$ for the bottom point,
$q=\FC(\zeta)_{2}$ for the middle point, and $q=0.8$ for the
top point.
\vspace{1em}

\begin{figure}[H]
  \centering
  \resizebox{130mm}{!}{
    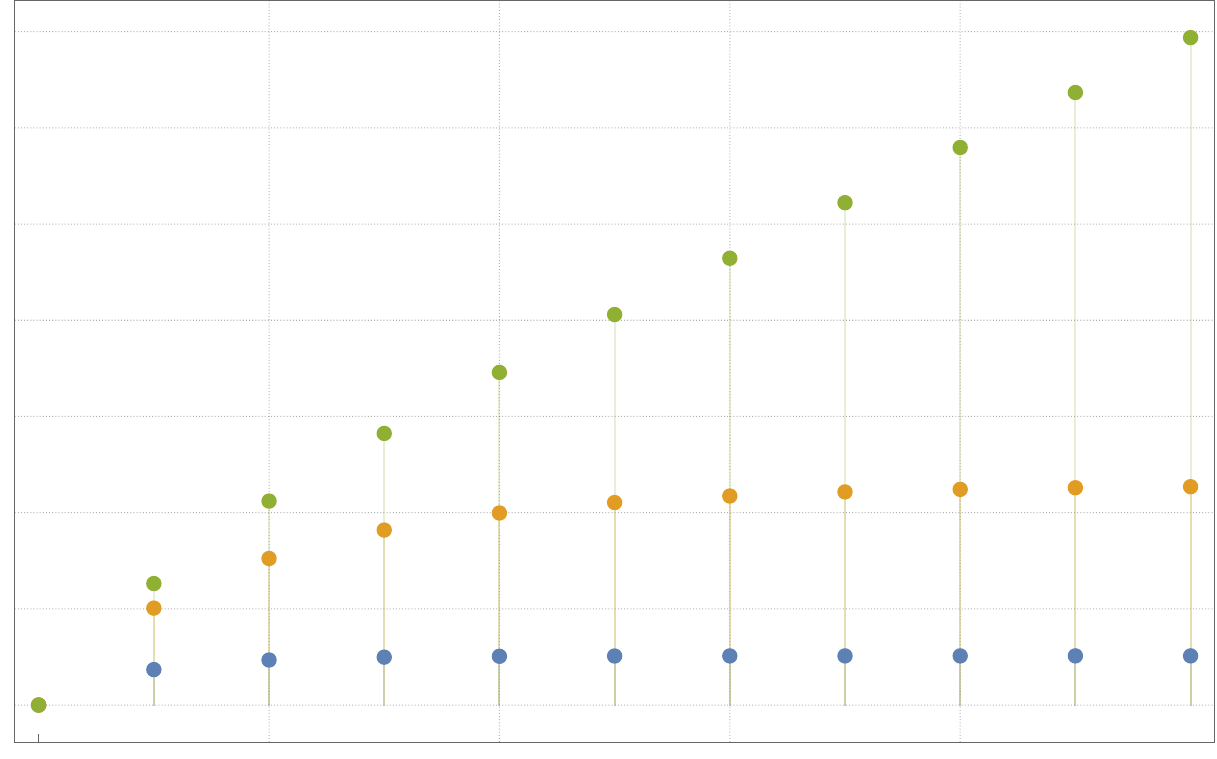}
  \caption{$\log(\cF_{m_{i}}(q))$ of 
      $q=\AC(\zeta)_{1}$ (bottom),
    $\AC(\zeta)_{2}$  (middle),
    and $0.5$ (top)}
  \label{fig:non-canonical-med-asymp-hill-seq}
\end{figure}

Third, Proposition~\ref{prop:path-parcel-asymp-hill-shape}
gives the asymptotically hill-shape sequences
$u(\zeta,\AC(\zeta)_{1})$ and $u(\zeta,\AC(\zeta)_{2})$.  Also, we
obtain the strictly increasing sequence $u(\zeta,q)$ for
$\AC(\zeta)_{1} <q=0.5<\AC(\zeta)_{2}$.
Figure~\ref{fig:non-canonical-med-asymp-hill-seq} plots
$\log(\cF_{m_{i}}(q))$ of $q=\AC(\zeta)_{1} $ for the bottom
point, $q=\AC(\zeta)_{2}$ for the middle point, and $q=0.5$ for
the top point.

\vspace{1em}

\subsection{A non-canonical mediator
  only with a
  semi-phase transition}
\label{sec:non-canonical-mediator-without-pt}

We have a parcel $\cF$ with a front semi-phase transition,
but without a front phase transition.
Consider
$s=\tp{0,\infty}$, $l=1$, $w=\tp{1}$, $\rho=\tp{1}$, and
$\fX=\{q\}$.  Also, let
\begin{align}
  \phi(x)
  =\tp{\frac{3q^{2}}{2} -\frac{q}{2}+\frac{3}{8}}
  \in \Q(x)^{l}.
\end{align}
Then, $\phi$ is a $\tp{s,l,w,\llq,\rho,x,\fX}$-mediator as in Section~\ref{sec:non-canonical-mediator-with-pt},  
since
$\frac{3q^{2}}{2}-\frac{q}{2}+\frac{1}{8}=
\frac{3}{2} \left(q-\frac{1}{6}\right)^{2}+\frac{1}{3}$
and $\rho=\tp{1}$.  Hence, the analogs
$\phio(x_{1})^{\lamo} [\lamo]!_{q}$ of
 $q$-Pochhammer
symbols 
$(\lamo)_{q}$ give
the $\llq$-merge-log-concave monomial parcel
\begin{align}
\cF=\Lam(s,l,w,\llq,\Psi_{s,\tp{\tp{0,1,0}},q},\phi,\rho,x,\fX)  
\end{align}
such that
$\cF_{\lam} =\frac{q^{\lamo}}{ \phio(x_{1})^{\lamo}
  [\lamo]!_{q}}$ for $\lam\ldZl s$.  This also gives the
ideal merged pair $\zeta=\Omg(\cF)$.

\subsubsection{On critical points and
  phase transitions}
Let $P=\pi(1)=\{\tp{s,l,m_{i},n_{i},k_{i}}\}_{i\ldZ \theta}$
with $m_{0}=n_{0}=\tp{0}$ and $k_{0}=a_{0}=b_{0}=\tp{0,1}$.  Then,
$\zeta$ has the single front critical point
$\FC(\zeta)=\fraa$ that solves
\begin{align}
  \cF_{m_{0}}(q)
  =1=
  \frac{q}{\frac{3q^{2}}{2}
  -\frac{q}{2}+\frac{3}{8}}=\cF_{\mn}(q).
\end{align} 
As $\FZ(\zeta)\neq \AfX$, $\zeta$ has the front semi-phase
transition at $\FC(\zeta)$.  But, $\zeta$ has no front phase
transitions, because $\cF_{m_{0}}(q)\geAfX \cF_{m_{1}}(q)$
by
$\frac{3q^{2}}{2} -\frac{q}{2}+\frac{3}{8}-q
=\frac{3}{2}(q-\fraa)^{2}\geq 0$.

There are no rear critical points for the infinite-length
$P$.  By Lemma~\ref{lem:mono-tame}, $\cF$ is tame along $P$
by $q$. But, there are no asymptotic critical points, since
$\phi_{1}(q) \gAfX \frac{1}{4}$ implies
\begin{align}
  \AD(\zeta)
  =\left(\frac{3q^{2}}{2}
  -\frac{q}{2}+\frac{3}{8}\right)-q(1-q)\gAfX 0.
\end{align}

\subsubsection{Polynomials with positive integer
  coefficients of an ideal merged pair}
\label{sec:non-canonical-medaitor-without-pt-poly}
The merged determinants of
$\zeta$ give the following
$q$-polynomials with positive coefficients:
\begin{align}
  \Delta(\cF)(s,l,w,\mn,\nn,\kn,x,\fX)
  &=q^{3};\\
  \Delta(\cF)(s,l,w,\mt,\nt,\kt,x,\fX)
  &=q^{6};\\
  \Delta(\cF)(s,l,w,\mr,\nr,\kr,x,\fX)
  &=q^{9};\\
  \dots.
\end{align}

\subsection{A weight-zero parcel with critical points and without semi-phase transitions}
\label{sec:triv-no-crit}
For $s=\tp{1,\infty}$, $l=1$, and $\fX=\{q\}$, consider the
weight-zero $q$-number parcel
$\cF=\Lam(s,l,\llq,\chi_{s,l,q},\fX)$. Also, there is the fitting path
$P=\pi(\lam,h) =\{\tp{s,l,m_{i},n_{i},k_{i}}\}_{i\ldZ \theta}$
with $m_{0}=n_{0}= \mn\bom \kn= \tp{h}\geq 1$.  Hence, there is the merged pair $\zeta=\tp{P,\cF}$. 

First, $\zeta$ has no front critical points by
$\cF_{m_{1}}-\cF_{m_{0}}\llq 0$.
Second, $\zeta$ has no rear critical points either, since $P$
is infinite-length.  Third, $u(\zeta,r)$ for each
$r\in \AfX$ is asymptotically hill-shape, because $0<q<1$
implies
\begin{align}
  \lim_{i\to \infty}\frac{\chi_{s,l,q,m_{i+1}}}{\chi_{s,l,q,m_{i}}}
  =\lim_{i\to\infty}\frac{1+q+\dots+q^{m_{i+1,1}-1}}{
  1+q+\dots+q^{m_{i,1}-1}}
  =
  \lim_{i\to\infty}
  \frac{1-q^{m_{i+1,1}}}{1-q^{m_{i,1}}}
  =1.
\end{align}
Therefore, $\zeta$ has asymptotic critical points without semi-phase
transitions.
\subsection{A finite-length merged pair with
  a rear phase transition}
\label{sec:crit-second}

For $s=\tp{1,3}$, $l=1$, and $\fX=\{q\}$, we have the $\ggq$-merged-log-concave $q$-Stirling parcel
$\cF=\Lam(s,l,\llq,S_{s,l,q},\fX)$.  Then,
Claim~\ref{c:ext-int-hadam-identity} of
Corollary~\ref{cor:ext-int-hadam} gives
 the
 finite-length $\llq$-merged-log-concave parcel
 \begin{align}
   \cG=\Lam(s,l,w,\llq,S_{s,l,q},x,\fX)   
 \end{align}
for
 $w=\tp{1}$ and $x=\tp{q}$.
 Moreover, 
$\theta=\tp{1,1}$ gives
  the fitting path $P=r_{\tp{1,\infty},\theta}(\pi(1,1))
=\{\tp{s,l,m_{i},n_{i},k_{i}}\}_{i\ldZ
  \theta}
$
such that $\mn=\nn=\tp{2}$ and $\kn=\tp{0,1}$. Therefore,
we obtain the finite-length ideal merged pair $\zeta=\tp{P,\cG}$.  Let
$m_{0}=n_{0}=\tp{1}$, $\mt=\nt=\tp{3}$, and $k_{0}=\kt=\kn$ for
our convenience.

\subsubsection{On critical points and phase transitions}
\label{sec:crit-second-crit}
We have no front critical points of $\zeta$, because
\begin{align}
\cG_{\mt}(q)=
  \frac{q+2}{(1-q)(1-q^2)}
    \gAfX 
  \frac{1}{1-q}
  =\cG_{\mn}(q).
\end{align}
Since $P$ is finite-length, $\zeta$ has no asymptotic critical
points either.  However, $\zeta$ has the rear critical point
$\RC(\zeta)=0.86676\dots$ that solves
\begin{align}
  \cG_{\mt}(q) =\frac{q+2}{(1-q)(1-q^2)}
  =\frac{1}{(1-q)(1-q^2)(1-q^3)}
  =\cG_{\mr}(q).
\end{align}
 
\begin{figure}[H]
  \centering
  \resizebox{130mm}{!}{
    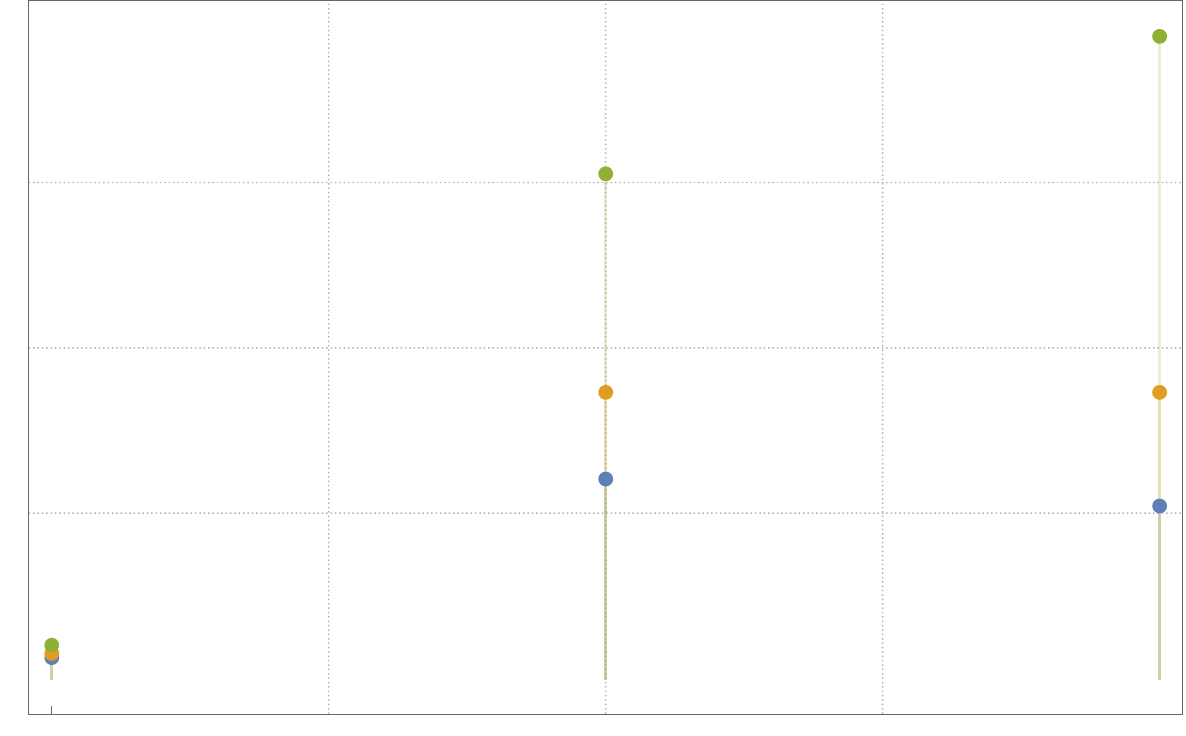}
  \caption{ $\cG_{m_{i}}(q)$ of
    $q=0.84$ (bottom), $\RC(\zeta)$ (middle), and $0.9$ (top)}
  \label{fig:stirling-second-kind-rear-pt}
\end{figure}

In particular, Theorem~\ref{thm:path-parcel-variation} gives
the rear phase transition of $\zeta$ such that
$u(\zeta,0.84)$ is two-slope hill-shape,
$u(\zeta,\RC(\zeta))$ is increasing and hill-shape, and
$u(\zeta,0.9)$ is strictly increasing.  In
Figure~\ref{fig:stirling-second-kind-rear-pt}, each
$i\in\oi(3)$ puts the bottom point for $\cG_{m_{i}}(0.84)$,
the middle point for $\cG_{m_{i}}(\RC(\zeta))$, and the
top point for
$\cG_{m_{i}}(0.9)$.

\subsubsection{Polynomials with positive integer
  coefficients of an ideal merged pair}
\label{sec:crit-second-pol}
The ideal merged pair $\zeta$ gives the following
$q$-polynomials with positive coefficients:
\begin{align}
  \Delta(\cG)(s,l,w,m_{0},n_{0},k_{0},x,\fX)
  &=\frac{(1)_{q}(2)_{q}}{(0)_{q}(1)_{q}}
  \cdot 
  \left(\frac{1}{(1-q)}\right)^{2}  
  \\    &=1+q;
  \\   \Delta(\cG)(s,l,w,\mn,\nn,\kn,x,\fX)
  &=\frac{(2)_{q}(3)_{q}}{(0)_{q}(1)_{q}}
  \cdot 
  \Bigg(\left(   \frac{q+2}{(1-q)(1-q^{2})}\right)^{2}
  \\&
  -
  \frac{1}{(1-q)} \cdot
  \frac{1}{(1-q)(1-q^{2})(1-q^{3})}\Bigg)
  \\ 
    &=3+7 q+9 q^2+5 q^3+q^4;\\
  \Delta(\cG)(s,l,w,\mt,\nt,\kt,x,\fX)
  & =\frac{(3)_{q}(4)_{q}}{(0)_{q}(1)_{q}}
  \cdot 
  \left(\frac{1}{(1-q)(1-q^{2})(1-q^{3})}\right)^{2}\\
  &=1+q+q^2+q^3.
\end{align}

\subsection{A higher-width parcel with
  a phase transition and conjectures}
\label{sec:high-width-w-pt}
 We discuss a higher-width parcel,
 unlike the examples above  in 
 Section~\ref{sec:examples}.
Let $l=3$, $s=\tp{0,\infty}$, $w=\tp{1}$, $\fX=\{\qq\}$, and
$\scc=>_{\qq}$. Also, let
$\gam=\tp{\tp{\fraa,0,0},\tp{-\fraa,0,0},
  \tp{\fraa,0,0}}$.
Then,  we obtain
 the width-three monomial parcel
$\cF=\Lam(s,l,w,\scc,\Psi_{s,\gam,q},x,\fX)$ such that
each $\lam\ldZr s$ satisfies
\begin{align}
  \cF_{\lam}
  =\frac{q^{\frac{\lamo^2-\lamt^{2}+\lamr^2}{2}}}
  { (\lamo)_{q}  (\lamt)_{q}   (\lamr)_{q}}.
\end{align}
Since $\gam$ is palindromic,  $\cF$ is static.

Moreover, let $\xi=\tp{2,1}$, $h=4$, and $\theta=\tp{1,\infty}$.
Then, Theorem~\ref{thm:fitting-path-young} gives the
infinite-length fitting path
$P=\bP_{s,l,\xi,h} =\{\tp{s,l,m_{i},n_{i},k_{i}}\}_{i\ldZ
  \theta}$ such that each $i\ldZ \theta$ gives
\begin{align}
  m_{i}&=n_{i}=\tp{6,5,4}+(i-1)\cdot\iota^{l}(4),\\
  a_{i}&=\nu(k_{i})=\tp{0,1,2,4,5,6},\\
  b_{i}&=\nu(m_{i},n_{i},k_{i})=\tp{6,6,6,10,10,10}
         +(i-1)\cdot \iota^{2l}(4),\\
  \sig(k_{i})&=\tp{6,4,2}.
\end{align}
Then, we have the merged pair $\zeta=\tp{P,\cF}$.  Let
$m_{0}=n_{0}=\tp{2,1,0}$, $k_{0}=\kn$, $a_{0} =\nu(k_{0})$,
and $b_{0}=\nu(m_{0},n_{0},k_{0})$ for our convenience.

\subsubsection{On critical points and phase transitions}\label{sec:high-width-w-pt-crit}
As $s$ is infinite, $\zeta$ has no rear critical points.
Moreover, $\cF$ is tame along $P$ by $0$,
because $0<\qq<1$ gives
\begin{align}
  \lim_{i\to\infty}
  \frac{\Psi_{s,\gam,q,m_{i+1}}}{\Psi_{s,\gam,q,m_{i}}}
  =
  \lim_{i\to\infty}\frac{
  q^{\frac{(4(i+1)+2)^2-(4(i+1)+1)^{2}+(4(i+1))^2}{2}}}
  {q^{\frac{(4i+2)^2-(4i+1)^{2}+(4i)^2}{2}}}
  =\lim_{i\to\infty}q^{16i+12}
  =0.
\end{align}
Thus, $\zeta$ has no asymptotic critical points, since
$u(\zeta,r)$ is hill-shape or decreasing for each
$r\in \AfX$ by Claim~\ref{c:path-parcel-ratio-zero-unim} of
Proposition~\ref{prop:path-parcel-ratio-unim-poly}.

\begin{figure}[H]
  \centering
  \resizebox{130mm}{!}{
    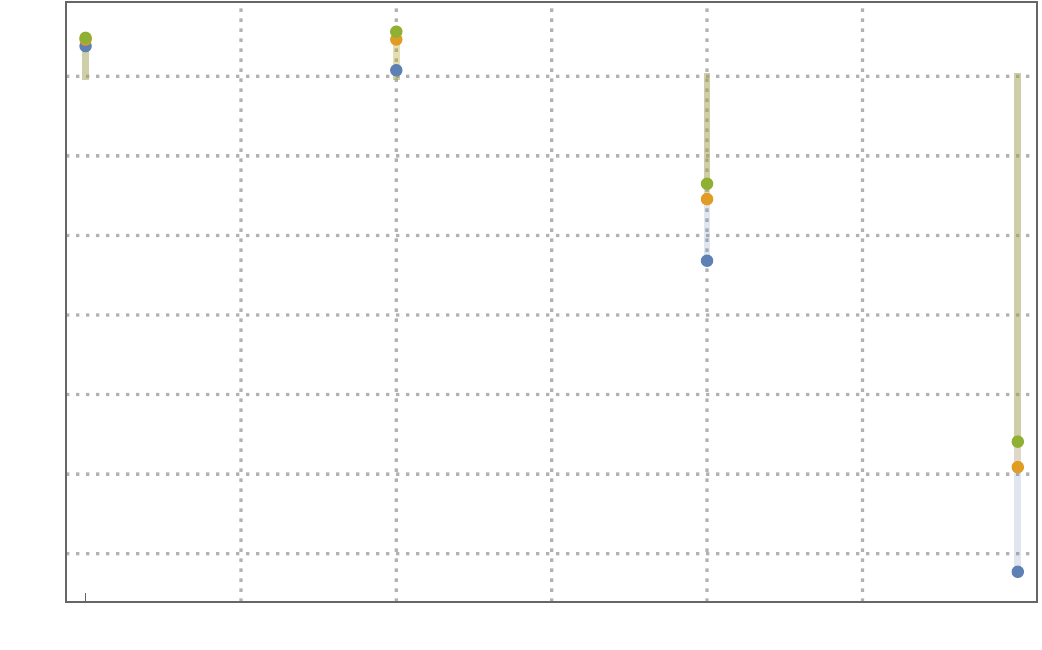}
  \caption{$\log(\cF_{m_{i}}(q))$ of 
      $q=0.8$ (bottom),
      $\FC(\zeta)$ (middle), and $0.83$ (top)}
  \label{fig:phase-transition-width-three} 
\end{figure}

We have $t_{\gam}(\mn)-t_{\gam}(m_{0})=12>0$.
Also, there is the front critical point
$\FC(\zeta)= 0.82439....$ that solves
\begin{align}
  \cF_{m_{0}}(\qq) 
  =\frac{q^{\frac{3}{2}}}{(2)_{q}(1)_{q}(0)_{q}}
  =\frac{q^{\frac{27}{2}}}{(6)_{q}(5)_{q}(4)_{q}}
  =\cF_{\mn}(\qq).
\end{align}
Hence, $\zeta$ has the unique front phase transition at
$\FC(\zeta)$ by Claim~\ref{c:gen-width-monom-pt-equiv} of
Proposition~\ref{prop:gen-width-monom-pt}.  Moreover,
$u(\zeta,0.8)$ is strictly decreasing, $u(\zeta,\FC(\zeta))$ is
decreasing and hill-shape, and $u(\zeta,0.83)$ is two-slope 
hill-shape.  For each $i\in\oi(0,3)$,
Figure~\ref{fig:phase-transition-width-three} gives the
bottom point for $\log(\cF_{m_{i}}(0.8))$, the middle point
for $\log(\cF_{m_{i}}(\FC(\zeta)))$, and the top point for
$\log(\cF_{m_{i}}(0.83))$.
Figure~\ref{fig:phase-transition-width-three} takes the log
scale to avoid point collisions.

\vspace{1em}

\subsubsection{Ideal property of a merged pair}
\label{sec:ideal-prop-merged-pair}
The parcel $\cF$ is not $\llq$-merged-log-concave, since
the following is not a $q$-polynomial:
\begin{dmath*}
  \Delta(\cF)(s,l,w,\tp{3,3,3},\tp{3,4,3},\tp{0,0,0,1,0,1},x,\fX)
  =
  q^{\frac{37}{2}} + 3q^{\frac{35}{2}}+ 7q^{\frac{33}{2}}+ 13q^{\frac{31}{2}}+ 19q^{\frac{29}{2}}+ 23q^{\frac{27}{2}}+ 24q^{\frac{25}{2}}+ 23q^{\frac{23}{2}}+ 20q^{\frac{21}{2}}+ 17q^{\frac{19}{2}}+ 12q^{\frac{17}{2}}+ 7q^{\frac{15}{2}}+ 3q^{\frac{13}{2}} +
  q^{\frac{11}{2}}.
\end{dmath*}
However,  $\zeta$ is ideal, because
each $i\ldZ \theta$ gives
\begin{align}
  \begin{split}
    &\Delta(\cF)(s,l,w,m_{i},n_{i},k_{i},x,\fX)\\
    & =
    \frac{
      (4i+2)_{q}(4i+2)_{q}(4i+2)_{q}
      (4i+6)_{q}(4i+6)_{q}(4i+6)_{q}
    }
    {(0)_{q}(1)_{q}(2)_{q}(4)_{q}(5)_{q}(6)_{q}}
    \\&
    \cdot 
    \Bigg(
    \bigg(\frac{q^{\frac{(4i+2)^2-(4i+1)^{2}+(4i)^2}{2}}}{
    (4i+2)_{q}  (4i+1)_{q}   (4i)_{q}}\bigg)^2
    \\&
    -
    \frac{q^{\frac{(4i-2)^2-(4i-3)^{2}+(4i-4)^2}{2}}}{
    (4i-2)_{q}  (4i-3)_{q}   (4i-4)_{q}}\cdot
    \frac{q^{\frac{(4i+6)^2-(4i+5)^{2}+(4i+4)^2}{2}}}{
    (4i+6)_{q}  (4i+5)_{q}   (4i+4)_{q}}
    \Bigg) 
  \end{split}
\end{align}
such that
\begin{align}
\left(q^{\frac{(4i+2)^2-(4i+1)^{2}+(4i)^2}{2}}\right)^{2}
&=q^{16 i^2 + 8i + 3}\llq 0,\\
q^{\frac{(4i-2)^2-(4i-3)^{2}+(4i-4)^2}{2}}
  \cdot
  q^{\frac{(4i+6)^2-(4i+5)^{2}+(4i+4)^2}{2}}
  &=q^{16 i^2 + 8i + 19}\llq 0.
\end{align}

\subsubsection{Polynomials with positive integer
  coefficients of an ideal merged pair}
\label{sec:high-width-w-pt-poly}
For instance,
$\Delta(\cF)(s,l,w,m_{0},n_{0},k_{0},x,\fX)$ is
\begin{dmath*}
  q^3+3 q^4+6 q^5+10 q^6+15 q^7+20 q^8+23 q^9+24 q^{10}+23 q^{11}+20 q^{12}+15 q^{13}+10 q^{14}+6 q^{15}+3 q^{16}+q^{17}.
\end{dmath*}
Then, $\Delta(\cF)(s,l,w,\mn,\nn,\kn,x,\fX)$ is
\begin{dmath*}
  q^{27}+5 q^{28}+19 q^{29}+58 q^{30}+158 q^{31}+388 q^{32}+885 q^{33}+1890 q^{34}+3828 q^{35}+7390 q^{36}+13688 q^{37}+24412 q^{38}+42089
  q^{39}+70327 q^{40}+114182 q^{41}+180469 q^{42}+278185 q^{43}+418794 q^{44}+616576 q^{45}+888721 q^{46}+1255398 q^{47}+1739429 q^{48}+2365848 q^{49}+3160960
  q^{50}+4151230 q^{51}+5361659 q^{52}+6814044 q^{53}+8524865 q^{54}+10503235 q^{55}+12748773 q^{56}+15249837 q^{57}+17982091 q^{58}+20907732 q^{59}+23975445
  q^{60}+27121214 q^{61}+30270096 q^{62}+33338786 q^{63}+36239100 q^{64}+38881893 q^{65}+41181562 q^{66}+43060365 q^{67}+44452792 q^{68}+45309075 q^{69}+45598218
  q^{70}+45309614 q^{71}+44453850 q^{72}+43061901 q^{73}+41183518 q^{74}+38884194 q^{75}+36241661 q^{76}+33341512 q^{77}+30272892 q^{78}+27123986 q^{79}+23978109
  q^{80}+20910214 q^{81}+17984335 q^{82}+15251804 q^{83}+12750444 q^{84}+10504608 q^{85}+8525955 q^{86}+6814878 q^{87}+5362272 q^{88}+4151661 q^{89}+3161248
  q^{90}+2366030 q^{91}+1739536 q^{92}+1255456 q^{93}+888749 q^{94}+616588 q^{95}+418798 q^{96}+278186 q^{97}+180469 q^{98}+114182 q^{99}+70327 q^{100}+42089
  q^{101}+24412 q^{102}+13688 q^{103}+7390 q^{104}+3828 q^{105}+1890 q^{106}+885 q^{107}+388 q^{108}+158 q^{109}+58 q^{110}+19 q^{111}+
  5 q^{112}+q^{113}.
\end{dmath*}
Also, $\Delta(\cF)(s,l,w,\mt,\nt,\kt,x,\fX)$ is
\begin{dmath*}
  q^{83}+5 q^{84}+19 q^{85}+58 q^{86}+158 q^{87}+390
  q^{88}+899 q^{89}+1951 q^{90}+\dots.  
\end{dmath*}

\subsubsection{Conjectures}
\label{sec:high-width-w-pt-conj}
Some leading coefficients of
$\Delta(\cF)(s,l,w,m_{i},n_{i},k_{i},x,\fX)$ for
$i\in\oi(2)$ coincide in the above.  Also,
these
$\Delta(\cF)(s,l,w,m_{i},n_{i},k_{i},x,\fX)$  are not palindromic, but log-concave
$q$-polynomials.  Thus, we make the following conjectures.
\begin{conjecture}
  For each $i\in \Zgez$,  consider the $q$-polynomial
  \begin{align}
    f_{i}=\sum_{j}f_{i,j}q^{j}=
    \Delta(\cF)(s,l,w,m_{i},n_{i},k_{i},x,\fX).
  \end{align}
  For
  $\ord_{q}(f_{i})<j_{2}<\dots<\deg_{q}(f_{i})$,
  let
  $f^{+}_{i} =\tp{f_{i,\ord_{q}(f_{i})},f_{i,j_{2}},\dots,
    f_{i,\deg_{q}(f_{i})}}$ denote the non-zero positive
  coefficients of $f_{i}$.  Then, each
  $i\in \Zgeo$ satisfies
  \begin{align}
    f^{+}_{i+1}(2+4(i-1),5+4(i-1))
    =      f^{+}_{i}(2+4(i-1),5+4(i-1)).
  \end{align}
\end{conjecture}

For example,
$f^{+}_{1}=\tp{1,5,19,58,158,388,885,\dots}$ and
$f^{+}_{2}=\tp{1,5,19,58,158,390,899,\dots}$
in Section~\ref{sec:high-width-w-pt-poly}.
  Thus, we have
\begin{align}
  f^{+}_{1}(2,5)
  =\tp{5,19,58,158}
  =    f^{+}_{2}(2,5).
\end{align}

We put the
following notation to state another conjecture
on log-concave $q$-polynomials.

\begin{definition}\label{def:poly-w-const}
  For a Laurent polynomial $f\in \Q[z^{\pm 1}]$, let
  \begin{align}
    C_{z}(f)=
    \begin{dcases}
      z^{-\ord_{z}(f)} f \mif f\neq 0,\\
      0 \melse.
    \end{dcases}
  \end{align}
\end{definition}

\begin{conjecture}\label{conj:high-width-w-pt}
  Let $i\in \Zgez$.
  \begin{enumerate}
  \item \label{c:high-width-w-pt-merged-det}
    $\Delta(\cF)(s,l,w,m_{i},n_{i},k_{i},x,\fX)$
    is a log-concave $q$-polynomial.
  \item \label{c:high-width-w-pt-merged-det-shifted-diff}
    The following is a log-concave
    $q$-polynomial:
    \begin{align}
      C_{q}(\Delta(\cF)(s,l,w,m_{i+1},n_{i+1},k_{i+1},x,\fX))
      -
      C_{q}(\Delta(\cF)(s,l,w,m_{i},n_{i},k_{i},x,\fX))\llq 0.
    \end{align}
  \end{enumerate}
\end{conjecture}



Claim~\ref{c:high-width-w-pt-merged-det}
of Conjecture~\ref{conj:high-width-w-pt} is
analogous to
Conjectures~\ref{conj:log-concavity-of-pre-merged-width-two}
and~\ref{conj:log-concavity-of-pre-merged-width-four}.

\section{Parcel convolutions}\label{sec:width-one}
We study some width-one parcels by convolutions.
\subsection{Convolution indices}
Consider a parcel
$\cF=\Lam(s,l,w,\scc,f_{s},\phi,\rho,x,\fX)$ with $l=1$ and
$s=\tp{0,\infty}$.  Then, multiplying generating functions
of $\cF$ gives a sequence of rational functions, which,
however, is not necessarily a parcel
$\Lam(s,l,w,\scc,g_{s},\phi,\rho,x,\fX)$ for some
$g_{s}=\{g_{s,m}\in \Q(\fX)\}_{m\in \Zl}$.  Hence, we
introduce the notion of {\it convolution indices}.

\begin{definition}\label{def:con-ind}
  Assume the following:
  \begin{enumerate}[label=(\alph*)]
  \item \label{a:con-ind-squaring} squaring orders
    $O_{i}=\{\sce_{i},\scc_{i}\}$ on $\fX$ for 
    $i\in\oi(3)$ such that $O_{3}\Sup O_{1}, O_{2}$;
  \item\label{a:con-ind-adm}
    $l=1$, $\lam\in \Qr_{>0}$, and
    $x=\tp{\tp{q^{\lam_{i}}}}_{i\in \oi(3)}
    \in \prod_{i\in \oi(3)}\Q(\fX)^{l}$;
  \item \label{a:con-ind-parcels}    
    $\cF_{i}
    =\Lam(s_{i},l,w,\scc_{i},f_{i,s_{i}},\phi,\rho_{i},x_{i},\fX)$
    for $i\in\oi(2)$;
  \item \label{a:con-ind-gates}
    $\rho=\tp{\rho_{i}}_{i\in \oi(3)}\in \prod_{i\in\oi(3)}\Zl$,
    $s=\tp{s_{1},s_{2},s_{1}+s_{2}}$,
    $O=\tp{O_{i}}_{i\in \oi(3)}$, and
    $o\in \Zgeo^{2}$.
  \end{enumerate}
 
  Then, we call $\ci=\tp{s,l,w,O,\phi,\rho,x,\fX,q,\lam,o}$
  convolution index of $\tp{\cFo,\cFt}$, when
  $\ci$ satisfies the
  following conditions.
  \begin{enumerate}
  \item 
    We have
    \begin{align}
      o_{1}\lamo&=o_{2}\lamt=\lamr,
                  \label{eq:con-ind-exp-equator}\\
      o_{1}^{-1}\rhoo&=o_{2}^{-1}\rhot  =\rhor.
                       \label{eq:con-ind-base-shift-equator}  
    \end{align} 
    We call equations~\eqref{eq:con-ind-exp-equator}
    and~\eqref{eq:con-ind-base-shift-equator}
    $\tp{\lam,o}$-exponent equator and
    $\tp{\rho,o}$-base-shift equator.
  \item
    $\phi$ is a
    $\tp{s_{i},l,w,\scc_{i},\tp{o_{i}},x_{i},\fX}$-mediator
    for each $i\in \oi(2)$. We call this
    $\tp{\ci,o}$-mediator condition.
  \item $\phi$ is a
    $\tp{s_{3},l,w,\sccr,\rhor,x_{3},\fX}$-mediator.  We
    call this $\tp{\ci,\rhor}$-mediator condition.
  \end{enumerate}
\end{definition}

In particular, we have the following on convolution indices.
\begin{lemma}\label{lem:con-ind-sp-cases}
  Suppose a
  convolution index 
  $\ci=\tp{s,l,w,O,\phi,\rho,x,\fX,q,\lam,o}$ 
  of $\tp{\cFo,\cFt}$.
  \begin{enumerate} 
  \item
    \label{c:con-ind-sp-cases-adm}
    Then, $x_{i}$ is $O_{i}$-admissible for each
    $i\in\oi(2)$.
  \item
    \label{c:con-ind-sp-cases-exp-adm}
    Provided the $\tp{\lam,o}$-exponent equator, 
    $x_{3}$ is $O_{3}$-admissible.
  \item \label{c:con-ind-sp-cases-equators-gates}
    Suppose
    that $\lam$ and $\rho$ are flat.
    \begin{enumerate}
    \item
      \label{c:con-ind-sp-cases-equators-gates-equators}
      If $o=\iota^{2}(1)$, then
      we have
      the $\tp{\lam,o}$-exponent and
      $\tp{\rho,o}$-base-shift equators
      and the $\tp{\ci,o}$-mediator condition.
    \item
      \label{c:con-ind-sp-cases-equators-gates-equators-gates}
      If one of $s_{1}$ and $s_{2}$ is infinite or
      $\tp{0,0}$, then we have
      the $\tp{\ci,\rhor}$-mediator condition.
    \item \label{c:con-ind-sp-cases-weight-zero}
      Similarly, if $w=\tp{0}$,
      then we have
      the $\tp{\ci,\rhor}$-mediator condition.
    \end{enumerate}
  \item \label{c:con-ind-sp-cases-mediators}
    Suppose
    that $\phi$ is the canonical $l$-mediator.
    \begin{enumerate}
    \item \label{c:con-ind-sp-cases-o-mediators}
      Then, the
      $\tp{\ci,o}$-mediator condition holds.
    \item \label{c:con-ind-sp-cases-rho-mediator}
      Provided the $\tp{\lam,o}$-exponent equator,
      the $\tp{\ci,\rhor}$-mediator condition holds.
    \end{enumerate} 
  \end{enumerate}
\end{lemma}
\begin{proof}
  Claim~\ref{c:con-ind-sp-cases-adm} holds, since the
  parcel $\cF_{i}$ has the squaring order $O_{i}$ and the
  base $x_{i}$.
  
  Let us prove
  Claim~\ref{c:con-ind-sp-cases-exp-adm}.  We
  have $o\in \Zgeo^{2}$ for the convolution index $\ci$.
  Then, since there is the $\tp{\lam,o}$-exponent equator,
  $x_{3}=\tp{q^{\lamr}}$ is $O_{1}$- and $O_{2}$-admissible
  by Claim~\ref{c:adm-film-O-admissible-monom} of
  Proposition~\ref{prop:adm-film} and
  Claim~\ref{c:con-ind-sp-cases-adm}.  Thus,
  Claim~\ref{c:con-ind-sp-cases-exp-adm} holds by the
  compatibility $O_{3}\Sup O_{1}, O_{2}$ in $\ci$.

  Let us prove
  Claim~\ref{c:con-ind-sp-cases-equators-gates-equators}.
  By $o=\iota^{2}(1)$, the flatness of $\lam$ and $\rho$ yields
  the $\tp{\lam,o}$-exponent equator and
  $\tp{\rho,o}$-base-shift equators.
  Let
  $\mu_{i}=\tp{s_{i},l,w,\scc_{i},\tp{o_{i}},x_{i},\fX}$
  for $i\in \oi(2)$.
  By Claim~\ref{c:con-ind-sp-cases-adm},
  Claim~\ref{c:adm-succ-half-gx-scc} of
  Lemma~\ref{lem:adm-succ} implies $1\scc_{i} 0$
  for $i\in \oi(2)$.
  This gives
  the base-shift positivity  of $\phi$ and $\mu_{i}$
  by
  Lemma~\ref{lem:bshift-fun-sp} and $o=\iota^{2}(1)$.
   Thus,
   the $\tp{\ci,o}$-mediator condition holds,
  since we have
  the base positivity  of $\phi$ and $\mu_{i}$ 
  in the parcel $\cF_{i}$.

  Let us prove
Claim~\ref{c:con-ind-sp-cases-equators-gates-equators-gates}.
  Suppose that $s_{1}$ is infinite or $\tp{0,0}$.
  Then,
  $m\ldZl s_{1}+s_{2}$ implies $m\ldZl s_{1}$ if $s_{1}$ is
  infinite, or $m\ldZl s_{2}$ otherwise.
  Thus, by the flatness of $\rho$ and $\lam$,  $m\ldZl s_{1}+s_{2}$
  gives
  $B(s_{1}+s_{2},l,w,m,\phi,\rhor,x_{3},\fX)
    = B(s_{i},l,w,m,\phi,\rho_{i},x_{i},\fX)\scc_{i} 0$
if $m\ldZl s_{i}$
  Therefore, the $\tp{\ci,\rhor}$-mediator condition holds
  by the compatibility $O_{3}\Sup O_{1}, O_{2}$ in $\ci$.

  Let us prove
  Claim~\ref{c:con-ind-sp-cases-weight-zero}.
  Let $\mu_{3}=\tp{s_{3},l,w,\sccr,\rhor,x_{3},\fX}$.
  Notice that $x_{3}$ is $\sccr$-admissible
  by the flat $\lam$ and the compatibility
  $O_{3}\Sup O_{1}$.
  First,
  Claim~\ref{c:con-ind-sp-cases-adm}
  and
  Claim~\ref{c:adm-succ-half-gx-scc} of
  Lemma~\ref{lem:adm-succ} imply $1\scco 0$,
  which gives $1\sccr 0$ by
  the compatibility.
  Hence, the base-shift positivity of $\phi$ and $\mu_{3}$
  follows from Lemma~\ref{lem:bshift-fun-sp}.
  Moreover, Claim~\ref{c:con-ind-sp-cases-weight-zero}
  follows,
  since we have
  the base positivity  of $\phi$ and $\mu_{3}$
  by $w=\tp{0}$.

  Claim~\ref{c:con-ind-sp-cases-o-mediators} follows
  from
  Claim~\ref{c:con-ind-sp-cases-adm},
  Claim~\ref{c:adm-succ-half-gx-scc} of
  Lemma~\ref{lem:adm-succ}, and
  Proposition~\ref{prop:canon-med}, and so does
  Claim~\ref{c:con-ind-sp-cases-rho-mediator} from
  Claim~\ref{c:con-ind-sp-cases-exp-adm}.
\end{proof}

We now obtain parcels by convolution indices.

\begin{proposition}\label{prop:con-ind-parcel}
  Consider
  $\cF_{i}=\Lam(s_{i},l,w,\scc_{i},f_{i,s_{i}},\phi,\rho_{i},x_{i},\fX)$
  for $i\in\oi(2)$ with a convolution index
  $\ci=\tp{s,l,w,O,\phi,\rho,x,\fX,q,\lam,o}$.  Assume
  $\cH=\{ \cH_{\mn}=\sum_{\mt\in \Zl} \cF_{1,\mt}\cdot
  \cF_{2,\mn-\mt} \in \Q(\fX)\}_{\mn\in \Zl}$, and
  $f_{3,s_{3}}=\{ f_{3,s_{3},m}=\prod \phi(x_{3})^{m\rc w} \cdot
  [m]!_{x_{3}}^{w}\cdot \cH_{m} \in \Q(\fX)\}_{m\in \Zl}$.  Let
  $\tau_{i}=\tp{o_{i}}\in \Zgeo^{l}$ for $i\in\oi(2)$.  Then,
  we
  have the following.
  \begin{enumerate}
  \item \label{c:cond-inv-parcel-taking-zero} When
    $m\nld_{\Zl}s_{3}$, $\cH_{m}=0$.
  \item \label{c:cond-inv-parcel-numerator}
    When $\mn\in \Zl$, 
    $f_{3,s_{3},\mn}=
    \sum_{\mt\in \Zl}
    {\mn\brack\mt}_{x_{3}}^{w}
    \cdot B(s_{1},l,w,\mt,\phi,\tau_{1},x_{1},\fX)
    \cdot B(s_{2},l,w,\mn-\mt,\phi,\tau_{2},x_{2},\fX)
    \cdot f_{1,s_{1},\mt}f_{2,s_{2},\mn-\mt}$.
  \item $f_{3,s_{3}}$ is pairwise
    $\tp{s_{3},l,\sccr}$-positive.
    \label{c:cond-inv-parcel-numerator-pairwise-positive}
  \item \label{c:cond-inv-parcel-parcel}
    $\cH$ is a parcel such that
    $\cH
    =\Lam(s_{3},l,w,\sccr,f_{3,s_{3}},\phi,\rhor,x_{3},\fX)$.
  \end{enumerate}
\end{proposition}
\begin{proof}
  Claim~\ref{c:cond-inv-parcel-taking-zero} holds,
  because $\cF_{1,m}=0$ for $m\nldZl s_{1}$ and
  $\cF_{2,m}=0$ for $m\nldZl s_{2}$.
  
  Let us prove Claim~\ref{c:cond-inv-parcel-numerator}.
  Since $x_{3}=\tp{q^{\lamr}}$,
  we
  have
  $x_{3}=\tp{q^{o_{1}\lamo}} =x_{1}^{\tau_{1}}$
 and $x_{3}=\tp{q^{o_{2}\lamt}}=x_{2}^{\tau_{2}}$
  by
  the $\tp{\lam,o}$-exponent equator.  Thus,
  for each
  $i\in \oi(2)$ and
  $m\ldZl s_{i}$, we have
  \begin{align}
    \frac{\prod \phi(x_{3})^{m\rc w}\cdot
    [m]!_{x_{3}}^{w}}{
    \prod \phi(x_{i})^{m\rc w}\cdot[m]!_{x_{i}}^{w}}
    =
    \frac{\prod \phi(x_{i}^{\tau_{i}})^{m\rc w}\cdot
    [m]!_{x_{i}^{\tau_{i}}}^{w}}{
    \prod \phi(x_{i})^{m\rc w}\cdot[m]!_{x_{i}}^{w}}
    =    B(s_{i},l,w,m,\phi,\tau_{i},x_{i},\fX).
  \end{align}
  Hence, if $\mn\ldZl s_{3}$, then
  \begin{align}
    f_{3,s_{3},\mn}
    &=
      \prod \phi(x_{3})^{\mn\rc w} \cdot  [\mn]!_{x_{3}}^{w}
      \cdot \sum_{\mt\in \Zl} \cF_{1,\mt}\cF_{2,\mn-\mt} 
    \\&=
    \prod \phi(x_{3})^{\mn\rc w} \cdot  [\mn]!_{x_{3}}^{w}
    \\  &       \cdot \sum_{\mt\ldZl s_{1}, \mn-\mt\ldZl s_{2}}
          \frac{f_{1,s_{1},\mt}}{\prod \phi(x_{1})^{\mt\rc w}
          \cdot [m_{2}]!_{x_{1}}^{w}}
          \cdot \frac{f_{2,s_{2},\mn-\mt}}{
          \prod \phi(x_{2})^{(\mn-\mt)\rc w}\cdot
          [\mn-\mt]!_{x_{2}}^{w}}
    \\&=
    \sum_{\mt\ldZl s_{1},\mn-\mt\ldZl s_{2}}
    \frac{   [\mn]!_{x_{3}}^{w}}{ [m_{2}]!_{x_{3}}^{w}
    [\mn-\mt]!_{x_{3}}^{w}}
    \\&
    \cdot \frac{\prod \phi(x_{3})^{\mt\rc w}\cdot
    [m_{2}]!_{x_{3}}^{w}}{
    \prod \phi(x_{1})^{\mt\rc w}\cdot[m_{2}]!_{x_{1}}^{w}}
    \cdot \frac{\prod \phi(x_{3})^{(\mn-\mt)\rc w}
    \cdot [\mn-\mt]!_{x_{3}}^{w}}{
    \prod \phi(x_{2})^{(\mn-\mt)  \rc w}\cdot
    [\mn-\mt]!_{x_{2}}^{w}}
    \cdot f_{1,s_{1},\mt}f_{2,s_{2},\mn-\mt}
    \\&=
    \sum_{\mt\ldZl s_{1},\mn-\mt\ldZl s_{2}}
    {\mn\brack\mt}_{x_{3}}^{w}
    \\&
    \cdot B(s_{1},l,w,\mt,\phi,\tau_{1},x_{1},\fX)
    \cdot B(s_{2},l,w,\mn-\mt,\phi,\tau_{2},x_{2},\fX)
    \cdot f_{1,s_{1},\mt}f_{2,s_{2},\mn-\mt}.
  \end{align}
  If $\mn \nldZl s_{3}$, then both
  sides of
  Claim~\ref{c:cond-inv-parcel-numerator}
  are
  zero by Claim~\ref{c:cond-inv-parcel-taking-zero} and
  the
  pairwise positivities of $f_{1,s_{1}}$ and $f_{2,s_{2}}$.
  
  Let us prove
Claim~\ref{c:cond-inv-parcel-numerator-pairwise-positive}.
For  $m_{1,1},m_{2,1}\ldZl s_{3}$, we obtain
  \begin{align}
    \prod_{i\in \oi(2)}f_{3,s_{3},m_{i,1}}
    &=
      \sum_{m_{1,2},m_{2,2}\ldZl s_{1},
      m_{1,1}-m_{1,2},m_{2,1}-m_{2,2}\ldZl s_{2}}
    \\ & \prod_{i\in \oi(2)}
      {m_{i,1}\brack m_{i,2}}_{x_{3}}^{w}
    \cdot
         B(s_{1},l,w,m_{i,2},\phi,\tau_{1},x_{1},\fX)
         \cdot  B(s_{2},l,w,m_{i,1}-m_{i,2},\phi,\tau_{2},x_{2},\fX)
    \\ &  \cdot f_{1,s_{1},m_{i,2}}
         \cdot f_{2,s_{2},m_{i,1}-m_{i,2}}.
\label{eq:cond-inv-parcel-numerator-prod-expanded}
  \end{align}
  In equation~\eqref{eq:cond-inv-parcel-numerator-prod-expanded},
  $m_{i,1},m_{i,1}-m_{i,2}\geq 0$. This gives
  $\prod_{i\in \oi(2)}
      {m_{i,1}\brack m_{i,2}}_{x_{3}}^{w}
      \sccr0$
      by the $\sccr$-admissibility of $x_{3}$
      in
  Claim~\ref{c:con-ind-sp-cases-exp-adm} of
  Lemma~\ref{lem:con-ind-sp-cases}.
  Also,  in
equation~\eqref{eq:cond-inv-parcel-numerator-prod-expanded},
  the compatibility $\sccr \Sup \scco, \scct$
  implies
  $\prod_{i\in \oi(2)} B(s_{1},l,w,m_{i,2},\phi,\tau_{1},x_{1},\fX)
  \cdot  B(s_{2},l,w,m_{i,1}-m_{i,2},\phi,\tau_{2},x_{2},\fX)\sccr 0$
  and
  $\prod_{i\in \oi(2)} f_{1,s_{1},m_{i,2}}\cdot f_{2,s_{2},m_{i,1}-m_{i,2}}\sccr 0$
  by the pairwise positivities
  of $f_{1,s_{1}}$ and $f_{2,s_{2}}$.
  Hence,
    Claim~\ref{c:cond-inv-parcel-numerator-pairwise-positive}
  holds by
equation~\eqref{eq:cond-inv-parcel-numerator-prod-expanded},
  since $f_{3,s_{3},n}=0$ for $n\nldZl s_{3}$ by
  Claim~\ref{c:cond-inv-parcel-taking-zero}.

  Claim~\ref{c:cond-inv-parcel-parcel} follows
  from
  Claim~\ref{c:cond-inv-parcel-numerator-pairwise-positive}
  and
  the $\tp{\ci,\rhor}$-mediator condition.
\end{proof}

By Proposition~\ref{prop:con-ind-parcel}, we introduce
the parcel convolutions below.

\begin{definition}\label{def:conv-parcel}
  Suppose
  $\cF_{i}
  =\Lam(s_{i},l,w,\scc_{i},f_{i,s_{i}},\phi,\rho_{i},x_{i},\fX)$
  for $i\in\oi(2)$ with a convolution index
  $\ci=\tp{s,l,w,O,\phi,\rho,x,\fX,q,\lam,o}$.
  Then, we define the  parcel convolution
  \begin{align}
    \cFo*\cFt
    =\Lam(s_{3},l,w,\sccr,f_{3,s_{3}},\phi,\rhor,x_{3},\fX)
  \end{align}
  such that
  $f_{3,s_{3},m}=
    \prod \phi(x_{3})^{m\rc w}\cdot
    [m]!_{x_{3}}^{w}\cdot
    \sum_{\nn+\nt=m}\cF_{1,\nn}\cF_{2,\nt}\in \Q(\fX)$
  for  $m\in \Zl$.
\end{definition}

In particular, we have the following parcels by
identity functions.

\begin{lemma}\label{lem:conv-self-zero}
  Suppose $s_{1}=\tp{0,0}$, $l=1$,
  $\cF= \Lam(s_{1},l,w,\scc,1_{s_{1}},\phi,\rho,x,\fX)$,
  and $\cG=\Lam(s_{2},l,w,\scc,f_{s_{2}},\phi,\rho,x,\fX)$.
  Then,
  $\cF*\cG=\cG$.
\end{lemma}
\begin{proof}
  Let  $x=\tp{q}$ and
  $O=\tp{O_{i}=\{\sce,\scc\}}_{i\in \oi(3)}$.
  Then, 
  $\tp{\cF,\cG}$ has
  the convolution index
  $\tp{\tp{s_{1},s_{2},s_{2}},l,w,O,\phi, \iota^{3}(\rho),\iota^{3}(x),\fX,q,
    \iota^{3}(1),\iota^{2}(1)}$
  by Claims~\ref{c:con-ind-sp-cases-equators-gates-equators}
  and~\ref{c:con-ind-sp-cases-equators-gates-equators-gates}
  of Lemma~\ref{lem:con-ind-sp-cases}.
  Thus,
  the
  assertion  holds by
  $\cF_{\tp{0}}=\frac{1_{s',\tp{0}}}{\prod \phi(x)^{\tp{0}\rc w}\cdot
    [\tp{0}]!_{ x}^{w}} =1$.
\end{proof}
 
Hence, we have the convolutions below by
Lemma~\ref{lem:conv-self-zero} and
Claim~\ref{c:con-ind-sp-cases-equators-gates}
of Lemma~\ref{lem:con-ind-sp-cases}.
\begin{definition}\label{def:conv-self}
  Let $l=1$, $\lam\in \Zgez$, and
  $\cF=\Lam(s,l,w,\scc,f_{s},\phi,\rho,x,\fX)$.  Assume that
  $s=\tp{0,0}$ or infinite, or  $w=\tp{0}$.
  Then, we define the
  $\lam$-fold parcel convolution
    \begin{align}
      \cF^{* \lam}
      &=
        \Lam(\lam s,l,w,\scc,g_{\lam s},\phi,\rho,x,\fX)
    \end{align}
    such that $g_{\lam s}=1_{\tp{0,0}} $
    when $\lam=0$, and 
    \begin{align}
      g_{\lam s, m}=    
      \prod \phi(x)^{m\rc w}\cdot
      [m]!_{x}^{w}\cdot
      \sum_{ \nn+\dots+n_{\lam}=m}
      \cF_{\nn}\cdot \cdots \cdot \cF_{n_{\lam}}
      \mfor m\in \Zl \mwhen \lam> 0.
    \end{align}
\end{definition}
 
\subsection{Extension of the
  Cauchy-Binet formula}
\label{sec:rshift}
The Cauchy-Binet formula describes minors of a matrix
product $AB$ by those of $A$ and $B$.  Hence, we extend the
Cauchy-Binet formula to obtain the merged-log-concavity of
parcel convolutions, as merged determinants extend
$2\times 2$ determinants. This uses the following notation with
Definition~\ref{def:bounded-increasing}.
\begin{definition}\label{def:minors}
  Let $d\in \Zgeo$, $\lam\in \Zgeo^{3}$, and $Q$ be a
  commutative ring.
  \begin{enumerate}
  \item Consider a $\lamo\times \lamt$-matrix
    $A=\tp{A_{i,j}}_{i\in\oi(\lamo), j\in\oi(\lamt)}\in
    M_{\lamo,\lamt}(Q)$.  Suppose $\alp \in \tsi(d,\lamo)$
    and $\bta \in \tsi(d,\lamt)$.  Then, we write the minor
    \begin{align}
      A(\alp,\bta)
      &=\det
        \begin{bmatrix}
          A_{\alpo,\btao}&  A_{\alpo,\btat} &\dots  &   A_{\alpo,\bta_{d}}\\
          A_{\alpt,\btao} &  A_{\alpt,\btat} &\dots  &          A_{\alpt,\bta_{d}}\\
          \dots &\dots &\dots &\dots\\
          A_{\alp_{d},\btao} & A_{\alp_{d},\btat} &\dots &A_{\alp_{d},\bta_{d}}
        \end{bmatrix}.
    \end{align}
  \item Suppose $\alp \in \tsi(d,\lamo)$,
    $\bta \in \tsi(d,\lamr)$.  $A\in M_{\lamo,\lamt}(Q)$,
    and $B\in M_{\lamt,\lamr}(Q)$. Then, let
    $\tsi(d,\lamt,\alp,\bta,A,B) =\{\gam\in
    \tsi(d,\lamt)\mid A(\alp,\gam)B(\gam,\bta)\neq 0\}$.
  \end{enumerate}  
\end{definition}

We recall the Cauchy-Binet formula, which holds by summation
reordering in the Leibniz formula of $A(\alp,\gam)$ and
$B(\gam,\bta)$.
\begin{theorem}[the Cauchy-Binet formula]\label{thm:cb}
  Let $d\in \Zgeo$, $\lam\in \Zrgeo$, and $Q$ be a
  commutative ring.  Consider $A\in M_{\lamo,\lamt}(Q)$,
  $B\in M_{\lamt,\lamr}(Q)$, $\alp \in \tsi(d,\lamo)$, and
  $\bta \in \tsi(d,\lamr)$.  Then, 
  $AB(\alp,\bta)
    =
      \sum_{ \gam\in \tsi(d,\lamt)}
      A(\alp,\gam)B(\gam,\bta)$.
\end{theorem}

  Let $\Qti$
  denote the units of a ring $Q$.
We then generalize the ring shift factors in
Definition~\ref{def:merged}.

\begin{definition}\label{def:rshift-ind-1st}
  Let $Q$ be a commutative ring and $d\in \Zgeo$.  Let
  $ \Theta(\ka,u)\in Q$ for $\ka\in\oi(3)$ and $u\in\Z$.
  \begin{enumerate}
  \item For $i\in\oi(2)$, suppose $\phi_{i}\in \Zgeo$,
    $\gam_{i}\in \tsi(d,\phi_{i})$, and $\mu_{i}\in\oi(3)$.
  \begin{enumerate}[label=(\alph*)]
  \item For each $u\in \Z$, we define
    the ring shift factor
    \begin{align}
      \Theta(\mu_{1},\mu_{2},\gam_{1},\gam_{2},u)=
      \begin{dcases}
        \frac{ \prod_{i\in\oi(d)}
        \Theta(\mu_{2},\gam_{2,i}-\gam_{1,1}+u)}
        { \prod_{i\in\oi(d)}  \Theta(\mu_{1},\gam_{1,i}-\gam_{1,1}+u)}
        \mif
        \prod_{i\in\oi(d)}
        \Theta(\mu_{1},\gam_{1,i}-\gam_{1,1}+u)\in \Qti,\\
        0 \melse.
      \end{dcases}
    \end{align}
  \item For each $u\in \Z$ and
    $C\in M_{\phi_{1},\phi_{2}}(Q)$, let
    \begin{align}
      C(\mu_{1},\mu_{2},\gam_{1},\gam_{2},\Theta,u)
      &=
        \Theta(\mu_{1},\mu_{2},\gam_{1},\gam_{2},u)
        \cdot C(\gam_{1},\gam_{2}).
    \end{align}
  \end{enumerate}
  If $\Theta(1,u)=\Theta(2,u)=\Theta(3,u)$ for each
  $u\in \Z$, then we also write $\Theta(u)=\Theta(\ka,u)$ for
  $\ka\in \oi(3)$,
  $\Theta(\gam_{1},\gam_{2},u)=
  \Theta(\mu_{1},\mu_{2},\gam_{1},\gam_{2},u)$, and
  $C(\gam_{1},\gam_{2},\Theta,u)
  =C(\mu_{1},\mu_{2},\gam_{1},\gam_{2},\Theta,u)$.
\item Suppose $\lam\in \Zrgeo$, $A\in M_{\lamo,\lamt}(Q)$,
  $B\in M_{\lamt,\lamr}(Q)$, $\alp \in \tsi(d,\lamo)$, and
  $\bta \in \tsi(d,\lamr)$.  Then, we call
  $\ri=\tp{\Theta,d,\lam,\alp,\bta,A,B,Q}$ ring shift index,
  when $\ri$ satisfies the following conditions:
  \begin{enumerate}
  \item
    \label{cond:rshift-ind-1st-pos-indices} 
    $\Theta(\ka,u)\in \Qti$ for each $\ka\in\oi(3)$ and $u\geq 0$;
  \item
    \label{cond:rshift-ind-1st-pos-indices-non-zero-det}
    $\prod_{i\in\oi(d)}\Theta(\ka,\gam_{i}-\alpo +u) \in \Qti$ for
    each $\ka\in\oi(3)$,
    $\gam\in \tsi(d,\lamt,\alp,\bta,A,B)$, and $u\geq 0$;
  \item
    \label{cond:rshift-ind-1st-indices-non-zero-det}
    $\prod_{i\in\oi(d)}\Theta(\ka,\bta_{i}-\gamo+u) \in \Qti$ for each
    $\ka\in\oi(2,3)$, $\gam\in \tsi(d,\lamt,\alp,\bta,A,B)$,
    and $u\geq \min(0,\gamo-\alpo )$.
  \end{enumerate}        
  We call
  Conditions~\ref{cond:rshift-ind-1st-pos-indices},
  ~\ref{cond:rshift-ind-1st-pos-indices-non-zero-det},
  and~\ref{cond:rshift-ind-1st-indices-non-zero-det}
  $\Theta$-nonsingularity,
  $\oi(3)$-product nonsingularity, and
  $\oi(2,3)$-product nonsingularity of $R$,
  respectively.
\end{enumerate}
\end{definition}

There is a ring shift index
$\tp{\Theta,d,\lam,\alp,\bta,A,B,Q}$ such that each
$\Theta(u)=1$. Then, $AB(\alp,\bta,\Theta,u)=AB(\alp,\bta)$,
$A(\alp,\gam,\Theta,u)=A(\alp,\gam)$, and
$B(\gam,\bta,\Theta,u) =B(\gam,\bta)$.  Hence, ring shift
indices extend minors by ring shift factors.

We state the following lemma to extend the Cauchy-Binet formula.
\begin{lemma}\label{lem:rshift}
  Consider a ring shift index
  $\tp{\Theta,d,\lam,\alp,\bta,A,B,Q}$.  Let
  $\gam\in \tsi(d,\lamt,\alp,\bta,A,B)$. Then, we have the
  following.
  \begin{enumerate}
  \item \label{c:rshift-double-tuples}
    $\Theta(\mu_{1},\mu_{2},\alp,\gam,u)\in \Qti$ whenever
    $\mu_{1},\mu_{2}\in\oi(3)$ and $u\geq 0$.
  \item \label{c:rshift-triple-tuples}
    $\Theta(\mu_{1},\mu_{2},\gam,\bta,u) \in \Qti$ whenever
    $\mu_{1},\mu_{2}\in\oi(2,3)$ and $u\geq \min(0,\gamo-\alpo )$.
  \end{enumerate}
\end{lemma}
\begin{proof}
  Claim~\ref{c:rshift-double-tuples} holds, since we have
  $\prod_{i\in\oi(d)} \Theta(\mu_{1},\alp_{i}-\alpo +u)\in \Qti$ and
  $ \prod_{i\in\oi(d)}\Theta(\mu_{2},\gam_{i}-\alpo +u)\in \Qti$ by the
  $\Theta$- and $\oi(3)$-product nonsingularities.
  
  Let us prove Claim~\ref{c:rshift-triple-tuples}.  If
  $0=\min(0,\gamo-\alpo )$, then
  $\prod_{i\in\oi(d)}\Theta(\mu_{1},\gam_{i}-\gamo+u) \in \Qti$ for
  each $u\geq \min(0,\gamo-\alpo )$ by the
  $\Theta$--nonsingularity; otherwise, the same holds by the
  $\oi(3)$--product nonsingularity.  Thus,
  Claim~\ref{c:rshift-triple-tuples} follows from the
  $\oi(2,3)$-product nonsingularity.
\end{proof}

\begin{theorem}[a ring shift extension of the Cauchy-Binet formula]
  \label{thm:rshift-ext-cb}
  Consider a ring shift index
  $\tp{\Theta,d,\lam,\alp,\bta,A,B,Q}$.
  Then, for each $u\in \Zgez$, we have
  \begin{align}
    AB(1,3,\alp,\bta,\Theta,u)
    &=\sum_{\gam\in \tsi(d,\lamt)}
      A(1,2,\alp,\gam,\Theta,u)
      B(2,3,\gam,\bta,\Theta,\gamo-\alpo +u).
      \label{eq:rshift-ext-cb}
  \end{align}
\end{theorem}
\begin{proof}
  Let
  $\gam\in \tsi(d,\lamt)$ and $u\geq 0$.  First, suppose
  $A(\alp,\gam)B(\gam,\bta)= 0$.  Then,
  \begin{align}
    &A(1,2,\alp,\gam,\Theta,u)
      B(2,3,\gam,\bta,\Theta,\gamo-\alpo +u)
    \\&=
    \Theta(1,2,\alp,\gam,u)
    A(\alp,\gam) 
    \Theta(2,3,\gam,\bta,\gamo-\alpo +u) B(\gam,\bta)
    \\&=0.
    \label{eq:rshift-ext-cb-zero-minor}
  \end{align}
  
  Second, suppose $A(\alp,\gam)B(\gam,\bta)\neq 0$.
  Then,
  Lemma \ref{lem:rshift} yields
  \begin{align}
    A(\alp,\gam)\cdot B(\gam,\bta)
    &=\Theta(1,2,\alp,\gam,u)^{-1} \cdot
      \Theta(2,3,\gam,\bta,\gamo-\alpo +u)^{-1}\\
    & \cdot A(1,2,\alp,\gam,\Theta,u) \cdot
      B(2,3,\gam,\bta,\Theta,\gamo-\alpo +u).
      \label{eq:rshift-ext-cb-minor-prod}
  \end{align}
  Moreover, we have
  \begin{dmath}
    \Theta(1,3,\alp,\bta,u)
    \cdot
    \Theta(1,2,\alp,\gam,u)^{-1}
    \cdot
    \Theta(2,3,\gam,\bta,\gamo-\alpo +u)^{-1}
    =
    \frac{ \prod_{i\in\oi(d)}
      \Theta(3,\bta_{i}-\alpo +u)}{
      \prod_{i\in\oi(d)}\Theta(1,\alp_{i}-\alpo +u)}
    \cdot
    \frac{
      \prod_{i\in\oi(d)}\Theta(1,\alp_{i}-\alpo +u)}
    { \prod_{i\in\oi(d)}\Theta(2,\gam_{i}-\alpo +u)}
    \cdot
    \frac{
      \prod_{i\in\oi(d)}
      \Theta(2,\gam_{i}-\gamo+\gamo-\alpo +u)}
    { \prod_{i\in\oi(d)}
      \Theta(3,\bta_{i}-\gamo+\gamo-\alpo +u)}
    =
    \frac{ \prod_{i\in\oi(d)}\Theta(3,\bta_{i}-\alpo +u)}{
      \prod_{i\in\oi(d)}\Theta(1,\alp_{i}-\alpo +u)}
    \cdot        \frac{
      \prod_{i\in\oi(d)}\Theta(1,\alp_{i}-\alpo +u)}
    { \prod_{i\in\oi(d)}\Theta(2,\gam_{i}-\alpo +u)}
    \cdot
    \frac{
      \prod_{i\in\oi(d)}\Theta(2,\gam_{i}-\alpo +u)}
    { \prod_{i\in\oi(d)}\Theta(3,\bta_{i}-\alpo +u)}
    =1. \label{eq:rshift-ext-cb-rshift-factor-prod}
  \end{dmath}
  
  Therefore,
  equations
  ~\eqref{eq:rshift-ext-cb-zero-minor},
   ~\eqref{eq:rshift-ext-cb-minor-prod},
  and~\eqref{eq:rshift-ext-cb-rshift-factor-prod} and
  Theorem~\ref{thm:cb} yield
  \begin{align}
    &AB(1,3,\alp,\bta,\Theta,u)
    \\&=\Theta(1,3,\alp,\bta,u)
    \cdot AB(\alp,\bta)
    \\&=
    \sum_{\gam \in \tsi(d,\lamt)}
    \Theta(1,3,\alp,\bta,u)
    \cdot
    A(\alp,\gam)B(\gam,\bta)
    \\&=
    \sum_{
    \gam \in
    \tsi(d,\lamt,\alp,\bta,A,B)} 
    \Theta(1,3,\alp,\bta,u)
    \cdot
    \Theta(1,2,\alp,\gam,u)^{-1}\cdot
    \Theta(2,3,\gam,\bta,\gamo-\alpo +u)^{-1}
    \\& \cdot A(1,2,\alp,\gam,\Theta,u)\cdot
    B(2,3,\gam,\bta,
    \Theta,\gamo-\alpo +u)
    \\&    =
    \sum_{
    \gam \in \tsi(d,\lamt,\alp,\bta,A,B)} 
    A(1,2,\alp,\gam,\Theta,u)
    B(2,3,\gam,\bta,\Theta,\gamo-\alpo +u)
    \\&=
    \sum_{
    \gam \in \tsi(d,\lamt)}
    A(1,2,\alp,\gam,\Theta,u)
    B(2,3,\gam,\bta,\Theta,\gamo-\alpo +u).
  \end{align}
\end{proof}
In particular, the ring shift extension reduces to the
Cauchy-Binet formula with trivial ring shift factors
$1=\Theta(1,3,\alp,\bta,u) = \Theta(1,2,\alp,\gam,u) =
\Theta(2,3,\gam,\bta,\gamo-\alpo +u)$ for each
$\gam\in \tsi(d,\lamt,\alp,\bta,A,B)$ and $u\geq 0$.

\begin{remark}
  If $\Theta(u)=[u]!_{q}$ for each $u\geq 0$, then
  $\Theta(u)$ gives a ring shift factor in
  Definition~\ref{def:merged}.  Moreover, consider
  $\Theta(1,u)=[u]!_{q}$, $\Theta(2,u)=[u]!_{q}^{2}$, and
  $\Theta(3,u)=[u]!_{q}^{2}$. Then,
  Proposition~\ref{prop:width-two-pre-merged} and
  Theorem~\ref{thm:rshift-ext-cb} give $q$-polynomials with
  positive integer coefficients by convolutions of
  $\{\frac{1}{[m]_{q}}\}_{m\in \Zgez}$ and
  $\{\frac{1}{[m]^{2}_{q}}\}_{m\in \Zgez}$.  However, we
  discuss a convolution of parcels with the same weights in
  this article.
\end{remark}
Consider the $t$-power series
$\sum_{\lam\in \Zgez}\frac{1} {(\lam)_{q}}t^{\lam}$ and
$\sum_{\lam\in \Zgez}\frac{q^{\frac{\lam^{2}}{2}}}
{(\lam)_{q}}t^{\lam}$.  Then, we conjecture the following polynomial
analog of Schoenberg's P\'{o}lya frequency
on real numbers~\cite[Theorem
3]{Scho} (see \cite{Bra,Edr}).  Proposition~\ref{prop:width-two-pre-merged}
gives the $d=2$ case of the conjecture.

\begin{conjecture}\label{conj:schoenberg-poly}
  Let $h,\lam,\rho,w\in \Zgeo$ and $d\in \Z_{\geq 3}$.  Consider
  $M_{w,h},N_{w,h}\in M_{h,h}(\Q(\qq))$ such that
  $M_{w,h,i,j}= \frac{1}{(j-i)^{w}_{q}}$ and
  $N_{w,h,i,j}=
  \frac{q^{\frac{(j-i)^{2}}{2}}}{(j-i)^{w}_{q}}$ if
  $0\leq j-i$, and $M_{w,h,i,j}=N_{w,h,i,j}=0$ otherwise.
  Moreover, for $\alp,\bta\in \tsi(d,h)$ such that
  $\bta\geq \alp$, let
  \begin{align}
    F(w,\rho,\lam,\alp,\bta)= \frac{
    \prod_{i\in\oi(d)}(\bta_{i}-\alpo +\lam)^{w}_{\qr}
    }{
    \prod_{i\in\oi(d)}
    (\alp_{i}-\alpo +\lam)^{w}_{\qr}}.       
  \end{align}
  Then, we have
  \begin{align}
    F(w,\rho,\lam,\alp,\bta)M_{w,h}(\alp,\bta)
    &\llq 0,\\
    F(w,\rho,\lam,\alp,\bta)N_{w,h}(\alp,\bta)
    &>_{q^{\fraa}} 0.
  \end{align}
\end{conjecture}

\subsection{Fitting tuples and
  strictly increasing sequences}
\label{sec:fit-inc}

We introduce the following notion to discuss fitting tuples
and strictly increasing sequences by the ring shift
extension of the Cauchy-Binet formula.

\begin{definition}\label{def:conv-key}
  Suppose a gate $s\geq 0$ and  $l\in \Zgeo$.  Let
  $u \in \Z$,
  $m,n\in \Zl$, and $k\in \Ztl$.  Then, we define
  $\omgn(\tp{s,l,m,n,k},u )
  =\tp{s,l,\alp,\bta}$ such that
  $\alp=\nu(k)-k_{1}+u \in \Ztl$ and
  $\bta=m \ccn n+\alp\in\Ztl$.
  \end{definition}

We then obtain strictly increasing
sequences from fitting tuples.

\begin{lemma}\label{lem:inc-from-fitting}
  For a fitting $\mu=\tp{s,l,m,n,k}$ and $u \in \Zgez$,
  let $\tp{s,l,\alp,\bta}=\omgn(\mu,u)$.
  \begin{enumerate}
  \item
    \label{c:inc-from-fitting-pre-fitting}
    Then, $\tp{2l,\alp,\bta}$ is pre-fitting and
    $\bta\geq \alp\geq u$.
  \item
    \label{c:inc-from-fitting-inc}
    If $l=1$, then $\alp,\bta\in \tsi(2l,\tp{u,h})$ for some
    $h\in \Zgeo$.
  \end{enumerate}
\end{lemma}
\begin{proof}
  Claim~\ref{c:inc-from-fitting-pre-fitting} follows from
  Claim~\ref{c:fitting-pre-fitting-fitting} of
  Proposition~\ref{prop:fitting-pre-fitting}
  and Claim~\ref{c:fitting-nonneg-ab-nonneg} of
  Lemma~\ref{lem:fitting-nonneg}.
  Claim~\ref{c:inc-from-fitting-pre-fitting} gives
  Claim~\ref{c:inc-from-fitting-inc}, since $\alp,\bta$ are
  strictly increasing when $l=1$.
\end{proof}

The following is
to obtain fitting tuples from increasing sequences.

\begin{definition}\label{def:fit-from-inc}
  For a gate $s\geq 0$ and $l\in \Zgeo$,
  suppose $\mu=\tp{s,l,\alp,\bta}$
  such that $\alp,\bta\in \Ztl$.
  Let  $u\in \Z$.
  Then, we define
  $\omgd(\mu,u)
  =\tp{s,l,m_{\alp,\bta},n_{\alp,\bta},k_{\alp,\bta,u}}$
  such that
  $m_{\alp,\bta}
  =(\bta-\alp)(1,l)\in \Zl$,
  $n_{\alp,\bta}
      =(\bta-\alp)(l+1,2l)\in \Zl$, and
      $k_{\alp,\bta,u}
      =\tp{u}\ccn (\alp(2,2l)-\alp(1,2l-1))\in \Ztl$.
      When $l=1$,  we also write
    $\omgd(s,\alp,\bta,u)=
    \omgd(\mu,u)$
    for our convenience.
\end{definition}

Then, we prove the following reciprocity between $\omgn$
and $\omgd$.

\begin{proposition}\label{prop:recip}
  Assume a gate $s\geq 0$ and $l\in \Zgeo$.  Let
  $m,n\in \Zl$, $k\in \Ztl$, and $\mu=\tp{s,l,m,n,k}$.  Also,
  let $\alp,\bta\in \Ztl$ and
  $\mu'=\tp{s,l,\alp,\bta}$. Then, we have the following.
  \begin{enumerate}
  \item \label{c:eps-omg}
    $\omgd(\tp{s,l,\omgn(\mu,\alpo )_{3},
      \omgn(\mu,\alpo )_{4}},\kn)=\mu$.
  \item  \label{c:omg-eps}
    $\omgn(\omgd(\mu',\kn),\alpo )=\mu'$.
  \item \label{c:omg-eps-fit-st-inc}
    Suppose  $l=1$, $\alpo\in \Zgez$, $k_{1}\in \Zgez$, and
    $\mu=\omgd(\mu',\kn)$.  Then, $\mu'$ satisfies
    $\alp-\bta\ldZtl s$ and
    $\alp,\bta\in \tsi(2l,\tp{\alpo,h})$ for some
    $h\in \Zgeo$ if and only if $\mu$ is fitting.
  \end{enumerate}
\end{proposition}
\begin{proof}
  Let us prove Claim~\ref{c:eps-omg}.  Let
  $\omgn(\mu,\alpo)=\tp{s,l,\alp',\bta'}$.  Then, we
  have 
  \begin{align}
    &\tp{\kn}\ccn (\alp'(2,2l)-\alp'(1,2l-1))\\
    &=\tp{\kn}\ccn ((\nu(k)-\kn+\alpo )(2,2l)
      -(\nu(k)-\kn+\alpo )(1,2l-1))\\
    &=\tp{\kn}\ccn (\nu(k)(2,2l)-\nu(k)(1,2l-1))\\
    &=k.
  \end{align}
  Hence, 
  Claim~\ref{c:eps-omg} holds
  by $\bta'-\alp'=m\ccn n$.

  Let us prove Claim~\ref{c:omg-eps}.  Let
  $\omgd(\mu',\kn) =\tp{s,l,m',n',k'}$. Then, 
  $k'=\tp{\kn}\ccn (\alp(2,2l)-\alp(1,2l-1))$ gives
  $\nu(k')-\kn+\alpo 
    =\tp{\kn,\kn+\alpt -\alpo,
      \dots, \kn+\alp_{2l}-\alpo }
      -\kn+\alpo 
      =\alp$.
  Thus, Claim~\ref{c:omg-eps} follows from
  $m'\ccn n'=\bta-\alp$.

  Let us prove Claim~\ref{c:omg-eps-fit-st-inc}.  The if
  part follows from Claim~\ref{c:omg-eps} and
  Claim~\ref{c:inc-from-fitting-inc} of
  Lemma~\ref{lem:inc-from-fitting}.  Let us confirm the only
  if part.  First, $\alp-\bta\ldZtl s$ gives the inclusion
  condition of $\mu$.  Second, $k_{1}\in \Zgez$ and
  $\alp\in \tsi(2l,\tp{\alpo,h})$ imply the lower slope
  condition of $\mu$.  Third,
  $\bta\in \tsi(2l,\tp{\alpo,h})$ gives the upper slope
  condition of $\mu$.
\end{proof}

\subsubsection{Merged determinants by Toeplitz
matrices and ring shift indices}
We consider the following Toeplitz
matrices.
\begin{definition}\label{def:toeplitz-l-one}
  For $l=1$, suppose
  $\cF=\{\cF_{m}\in \Q(\fX)\}_{m\in \Zl}$.
  For $h\in \Zgeo$ and $i, j\in\oi(h)$, we write
  a matrix $M_{\cF,h}\in M_{h,h}(\Q(\fX))$ such that
  its $\tp{i,j}$-element $M_{\cF,h,i,j}$ is $\cF_{\tp{j-i}}$.
\end{definition}

In particular, we have the following by $\omgd$.
\begin{lemma}\label{lem:non-zero-det-fitting} 
  Suppose $l=1$ and
  $\cF=\Lam(s,l,w,\scc,\fs,\phi,\rho,x,\fX)$.  Let $d=2$.
  Consider $h\in \Zget$ and $\alp,\bta\in \tsi(d,h)$.  Then,
  we have
  the following.
  \begin{enumerate}
  \item \label{c:non-zero-det-fitting-inclusion}
    $M_{\cF,h}(\alp,\bta)\neq 0$
    implies $\bta-\alp\ldZd s$.
  \item \label{c:non-zero-det-fitting-iff} Let
    $\cF$
    be $\scc'$-merged-log-concave.
    Then, $\omgd(s,\alp,\bta,u)$ is
    fitting for each $u\in\Zgez$ if and only if
    $M_{\cF,h}(\alp,\bta) \neq 0$.
  \end{enumerate}
\end{lemma}
\begin{proof}
  Let us prove
  Claim~\ref{c:non-zero-det-fitting-inclusion}.
  Assume $\bta_{1}-\alp_{1}<s_{1}$.  Then,
  $\cF_{\tp{\bta_{1}-\alp_{1}}}=0$.  Also,
  $\cF_{\tp{\bta_{1}-\alp_{2}}}=0$, since
  $\bta_{1}-\alp_{2}<s_{1}$ by $\alp_{2}>\alp_{1}$.
  Hence, we obtain
   the contradiction
  $0\neq M_{\cF,h}(\alp,\bta) =\cF_{\tp{\bta_{1}-\alp_{1}}}
  \cF_{\tp{\bta_{2}-\alp_{2}}} -
  \cF_{\tp{\bta_{2}-\alp_{1}}} \cF_{\tp{\bta_{1}-\alp_{2}}}
  = 0$.  The same contradiction occurs, when
  $\bta_{1}-\alp_{1}>s_{2}$, $\bta_{2}-\alp_{2}<s_{1}$,
  or
  $\bta_{2}-\alp_{2}>s_{2}$.

  Let us prove Claim~\ref{c:non-zero-det-fitting-iff}
  for
  $\omgd(s,\alp,\bta,u)
  =\tp{s,l,m_{\alp,\bta},n_{\alp,\bta},k_{\alp,\bta,u}}$.
  First, we prove the if part.  Since
  $k_{\alp,\bta,u}=\tp{u,\alp_{2}-\alp_{1}}$,
  $k_{\alp,\bta,u}\geq \tp{0,1}$ by $u\geq 0$ and
  $\alp_{2}>\alp_{1}$.  Moreover, since
  $m_{\alp,\bta}\ccn n_{\alp,\bta}=\bta-\alp$,
  $\bta_{2}>\bta_{1}$ gives
  $n_{\alp,\bta}+ k_{\alp,\bta,u,2}=
  \tp{\bta_{2}-\alp_{1}}>\tp{\bta_{1}-\alp_{1}}
  =m_{\alp,\bta}$.  Thus, Lemma~\ref{lem:fitting-one} and
  Claim~\ref{c:non-zero-det-fitting-inclusion} give the if
  part.  Second, we prove the only if part.  We have
  $m_{\alp,\bta}\bom
  k_{\alp,\bta,u}=\tp{\bta_{1}-\alp_{1}-(\alp_{2}-\alp_{1})}
  =\tp{\bta_{1}-\alp_{2}}$ and
  $n_{\alp,\bta}\bop
  k_{\alp,\bta,u}=\tp{\bta_{2}-\alp_{2}+(\alp_{2}-\alp_{1})}=\tp{\bta_{2}-\alp_{1}}$.
  Then,
  $M_{\cF,h}(\alp,\bta) = \cF_{m_{\alp,\bta}}
  \cF_{n_{\alp,\bta}^{\ve}}- \cF_{m_{\alp,\bta}\bom
    k_{\alp,\bta,u}} \cF_{(n_{\alp,\bta}\bop
    k_{\alp,\bta,u})^{\ve}}$.  Thus, the only if part holds
  by the $\scc'$-merged-log-concavity of $\cF$.
\end{proof}

Also, we use the following notation.
\begin{definition}
  Let $l\in \Zgeo$ and $w\in \Zlgez$.  Consider an
  indeterminate $y\in \Q(\fX)^{l}$ and
  $\phi(y)\in \prod_{i\in\oi(l)}\Q(y_{i})$.  Then, for each
  $u\in \Zl$, we write
  \begin{align}
    \til{\Theta}(l,w,y,\phi,u)=
    \begin{dcases}
      \prod \phi(y)^{u\rc w} \cdot
      [u]!_{y}^{w} \mif u\geq 0,\\
      0 \melse.
    \end{dcases}
  \end{align}
\end{definition}

We then realize merged determinants by Toeplitz
matrices and ring shift indices.
\begin{proposition}\label{prop:conv-key}
  Let $l=1$. Suppose
  $\cF_{i}
  =\Lam(s_{i},l,w,\scc_{i},f_{i,s_{i}},\phi,\rho_{i},x_{i},\fX)$
  for $i\in\oi(3)$ such that $\cF_{i}$ is
  $\scc'_{i}$-merged-log-concave for each $i\in\oi(2)$.
  Let
  $\ci=\tp{s,l,w,O,\phi,\rho,x,\fX,q,\lam,o}$ be the convolution
  index of $\tp{\cFo,\cFt}$ and $\cFr=\cFo*\cFt$.
  Let $y=x_{3}^{\rhor}$ and 
  $\Theta(u)=\til{\Theta}(l,w,y,\phi,\tp{u})$ for
  each $u\in \Z$.
   Also, for $d=2$, $h\in \Zget$, and
  $\alp,\bta\in \tsi(d,h)$,
  let $\ri=
  \tp{\Theta,d,\iota^{3}(h),\alp,\bta, M_{\cFo,h},
    M_{\cFt,h},\Q(\fX)}$.
  Then, we have the following.
  \begin{enumerate} 
  \item \label{c:conv-key-rshift-index}
    $\ri$ is a ring shift index.
  \item \label{c:conv-key-merged-det} For each
    $\gam\in \tsi(d,h)$ and $u\geq 0$, we have
    \begin{align}
      M_{\cFo,h}(\alp,\gam, \Theta,u)
      &=\Delta(\cFo)
        (s_{1},l,w,m_{\alp,\gam},n_{\alp,\gam},
        k_{\alp,\gam,u},\phi,\rhoo,x_{1},\fX),
        \label{eq:conv-key-merged-det-cFo}\\
      M_{\cFt,h}(\gam,\bta, \Theta,u)
      &=\Delta(\cFt)
        (s_{2},l,w,m_{\gam,\bta},n_{\gam,\bta},
        k_{\gam,\bta,u},\phi,\rhot,x_{2},\fX),
        \label{eq:conv-key-merged-det-cFt}\\
      M_{\cFr,h}(\alp,\bta,\Theta,u)
      & =\Delta(\cFr)(s_{3},l,w,m_{\alp,\bta},
        n_{\alp,\bta},k_{\alp,\bta,u},\phi,
        \rhor,x_{3},\fX).
        \label{eq:conv-key-merged-det-cFr}
    \end{align}
  \end{enumerate}
\end{proposition}
\begin{proof}
  Let us prove Claim~\ref{c:conv-key-rshift-index}.  First,
  we confirm the $\Theta$-nonsingularity.  The base positivity
  of $\phi$ in $\cFr$ gives $\prod \phi(x_{3})^{w}\gAfX 0$.  Thus,
  $\prod\phi(y)^{ w}\neq 0$ by $\rho_{3}\in \Zgeo^{l}$.  Moreover,
  $q^{\lamo}\gAfX 0$ by Claim~\ref{c:adm-bounds} of
  Lemma~\ref{lem:adm-bounds-poring-inclusions}.  This
  implies
  $y_{1}=q^{\lamr\rho_{3,1}}=q^{\lamo o_{1}\rho_{3,1}}
  \gAfX 0$,
  and so $[\tp{u}]!_{y}^{w}\neq 0$ for each $u\geq 0$.
  Hence, we have the $\Theta$-nonsingularity.

  Second, let us prove the $\oi(3)$- and
  $\oi(2,3)$-product
  nonsingularities.  Assume
  $\gam\in \tsi(d,h,\alp,\bta, M_{\cFo,h},M_{\cFt,h})$. Then,
    Claim~\ref{c:non-zero-det-fitting-inclusion} of
  Lemma~\ref{lem:non-zero-det-fitting} implies
  $\gamo-\alpo \geq 0$ by $M_{\cFo,h}(\alp,\gam) M_{\cFt,h}(\gam,\bta) \neq 0$.  Thus, we obtain the
  $\oi(3)$-product nonsingularity by the
  $\Theta$-nonsingularity and $\gamt>\gamo$.  Also, we have
  $0=\min(0,\gamo-\alpo)$, and $\btao-\gamo\geq 0$ by
  $M_{\cFo,h}(\alp,\gam) M_{\cFt,h}(\gam,\bta) \neq 0$ and
  Claim~\ref{c:non-zero-det-fitting-inclusion} of
  Lemma~\ref{lem:non-zero-det-fitting}.  Thus, the
  $\oi(2,3)$-product nonsingularity follows from the
  $\Theta$-nonsingularity and $\btat>\btao$. In particular,
  Claim~\ref{c:conv-key-rshift-index} holds.

  Let us prove Claim~\ref{c:conv-key-merged-det}.
  Hence,
  we prove
  equation~\eqref{eq:conv-key-merged-det-cFo}.
  Let $a_{\alp,\gam,u}=\nu(k_{\alp,\gam,u})$ and
  $b_{\alp,\gam,u}=\nu(m_{\alp,\gam},n_{\alp,\gam},
  k_{\alp,\gam,u})$. Then, 
  $k_{\alp,\gam,u}=\tp{u,\alpt-\alpo}$ and
  $m_{\alp,\gam}\ccn n_{\alp,\gam}=\gam-\alp$ give
  \begin{align}
    a_{\alp,\gam,u}
    &=\tp{u,\alpt -\alpo +u},
      \label{eq:conv-key-a}\\
    b_{\alp,\gam,u}
    &=a_{\alp,\gam,u}+ m_{\alp,\gam}\ccn n_{\alp,\gam}
      =\tp{\gamo-\alpo +u, \gamt-\alpo +u},
      \label{eq:conv-key-b}\\
    m_{\alp,\gam}\bom                 k_{\alp,\gam,u}
    &=\tp{\gamo-\alpo}-\tp{\alpt -\alpo}
      =\tp{\gamo-\alpt},
      \label{eq:conv-key-m-bom-k}\\
    n_{\alp,\gam}\bop                 k_{\alp,\gam,u}
    &=\tp{\gamt-\alpt}+\tp{\alpt -\alpo}
      =\tp{\gamt-\alpo}.
      \label{eq:conv-key-n-bop-k}
  \end{align}  
  
  First, suppose $\gamo-\alpo+u<0$.  Then,
  $\Theta(\alp,\gam,u)= \frac{\prod_{i\in\oi(d)}\Theta(\gam_{i}-\alpo
    +u)} { \prod_{i\in\oi(d)}\Theta(\alp_{i}-\alpo +u)}=0$.  Also, we
  have
  $\Ups(s_{1},l,w,m_{\alp,\gam},n_{\alp,\gam},
  k_{\alp,\gam,u},\phi,
  \rhoo,\xoo,\fX)=0$ by equation~\eqref{eq:conv-key-b}.
  Hence, equation~\eqref{eq:conv-key-merged-det-cFo}
  holds
  by $0=0$.  Second, suppose $\gamo-\alpo+u\geq 0$.
  Then,
  since $\alpt-\alpo+u\geq 0$ by $u\geq 0$,
  equations~\eqref{eq:conv-key-a},
  \eqref{eq:conv-key-b},
  \eqref{eq:conv-key-m-bom-k},
  and~\eqref{eq:conv-key-n-bop-k} give
  \begin{dmath*}
    M_{\cFo,h}(\alp,\gam,\Theta,u)
    =
    \Theta(\alp,\gam,u)
    M_{\cFo,h}(\alp,\gam)
    =
    \frac{
      \Theta(\gamt-\alpo +u)
      \Theta(\gamo-\alpo +u)}
    {\Theta(\alpt -\alpo +u)\Theta(u)}
    \cdot (M_{\cFo,h,\alpo,\gamo}
    M_{\cFo,h,\alpt,\gamt}
    -
    M_{\cFo,h,\alpo,\gamt}
    M_{\cFo,h,\alpt,\gamo})
    =
    \frac{
      \prod (\phi(y)^{\wcn})^{\tp{\gamt-\alpo +u,
        \gamo-\alpo +u}\rc w^{\wcn}}
      [\tp{\gamt-\alpo +u}]!_{y}^{w}
      [\tp{\gamo-\alpo +u}]!_{y}^{w}}
    { \prod (\phi(y)^{\wcn})^{
        \tp{\alpt -\alpo +u, u}\rc w^{\wcn}}
      [\tp{\alpt -\alpo +u}]!_{y}^{w}
      [\tp{u}]!_{y}^{w}}
    \cdot (\cF_{1,\tp{\gamo-\alpo}}
    \cF_{1,\tp{\gamt-\alpt}}
    -
    \cF_{1,\tp{\gamt-\alpo}}
    \cF_{1,\tp{\gamo-\alpt}})
    =
    \frac{
      \prod (\phi(y)^{\wcn})^{b_{\alp,\gam,u}\rc w^{\wcn}}
      [b_{\alp,\gam,u}]!_{y^{\wcn}}^{w^{\wcn}}}
    {  \prod (\phi(y)^{\wcn})^{a_{\alp,\gam,u}\rc w^{\wcn}}
      [a_{\alp,\gam,u}]!_{y^{\wcn}}^{w^{\wcn}}}
    \cdot (\cF_{1,m_{\alp,\gam}}
    \cF_{1,n_{\alp,\gam}^{\ve}}
    -
    \cF_{1,m_{\alp,\gam}\bom k_{\alp,\gam,u}}
    \cF_{1,(n_{\alp,\gam}\bop k_{\alp,\gam,u})^{\ve}}).
  \end{dmath*}
  Furthermore, we have
  \begin{align}
    \frac{
    \prod (\phi(y)^{\wcn})^{b_{\alp,\gam,u}\rc w^{\wcn}}
    [b_{\alp,\gam,u}]!_{y^{\wcn}}^{w^{\wcn}}}
    {  \prod (\phi(y)^{\wcn})^{a_{\alp,\gam,u}\rc w^{\wcn}}
    [a_{\alp,\gam,u}]!_{y^{\wcn}}^{w^{\wcn}}}
    =
    \Ups(s_{1},l,w,m_{\alp,\gam},
    n_{\alp,\gam},k_{\alp,\gam,u},\phi,
    \rhoo,\xoo,\fX),
  \end{align} 
  since
  $y=x_{3}^{\rhor} =\tp{q^{\lamr}}^{o_{1}^{-1}\rhoo}
  =\tp{q^{\lamo o_{1} o_{1}^{-1}\rho_{1,1}}}
  =\tp{q^{\lamo\rho_{1,1}}} =x_{1}^{\rhoo}$ in the convolution
  index $\ci$.  Thus, we obtain
  equation~\eqref{eq:conv-key-merged-det-cFo}.  Similarly,
  equations~\eqref{eq:conv-key-merged-det-cFt}
  and~\eqref{eq:conv-key-merged-det-cFr} hold.
  \end{proof}

\subsection{Merged-log-concavity}
\label{sec:merged-conv}
If there is a squaring order $\sccr'$ compatible to
$\scco', \scct'$, and $\sccr$ in
Proposition~\ref{prop:conv-key}, then
Theorem~\ref{thm:rshift-ext-cb} and
Proposition~\ref{prop:conv-key} give the
$\scer'$-merged-log-concavity of $\cFr$.  But, we
need the
$\sccr'$-merged-log-concavity of $\cFr$ for
 polynomials with positive integer coefficients.
Hence, we introduce the following notation
to obtain fitting tuples for
$M_{\cFo,h}(\alp,\gam, \Theta,u)$
and
$M_{\cFt,h}(\gam,\bta, \Theta,u)$ in
Proposition~\ref{prop:conv-key}.

\begin{definition}
  Suppose gates $s_{1},s_{2}\geq 0$.  For each $u\in \Z$
  and
  $\alp,\bta\in \Zt$, we define
  $\ofi(s_{1},s_{2},\alp,\bta,u)=\{\gam\in
  \tsi(2,\tp{\alpo,\btat})\mid \omgd(s_{1},\alp,\gam,u) \mand
  \omgd(s_{2},\gam,\bta,u) \mbox{ are fitting.}\}$.
\end{definition}

Then, we prove the following existence of fitting tuples.

\begin{lemma}\label{lem:conv-fitting}
  Let $l=1$ and $\alpo\in \Zgeo$.  Assume the following:
  \begin{enumerate}
  \item \label{a:conv-fitting-gates}
    $s_{3}=s_{1}+s_{2}$
    for gates $s_{1},s_{2},s_{3}\geq 0$;
  \item \label{a:conv-fitting-fitting-tuple} a fitting
    $\mu=\tp{s_{3},l,m,n,k}$;
  \item \label{a:conv-fitting-omg-fitting-tuple}
    $\tp{s_{3},l,\alp,\bta}=\omgn(\mu,\alpo)$;
  \item \label{a:conv-fitting-gamma}
    $\gam=\tp{\max(s_{1,1}+\alpo,\btao-s_{2,2}),
      \min(s_{1,2}+\alpt,\btat-s_{2,1})}\in \Z^{2}$.
  \end{enumerate}
  Then, we have the following.
  \begin{enumerate}[label=(\alph*)]
  \item
    \label{c:conv-fitting-inequalities}
    There exist the following inequalities:
    \begin{align}
      \gamo
      &\leq \min( s_{1,2}+\alpo,\btao-s_{2,1});
        \label{ineq:conv-fitting-max-alpo-btao-min}\\
      \max(s_{1,1}+\alpt,\btat-s_{2,2})
      &\leq \gamt;
        \label{ineq:conv-fitting-max-alpt-btat-min}\\
      \gamo  &< \gamt.
               \label{ineq:conv-fitting-gamo-gamt}
    \end{align}
  \item
    \label{c:conv-fitting-in-ofi}
    For each $u\in\Zgez$, there exists
    $\gam\in\ofi(s_{1},s_{2},\alp,\bta,u)$.
  \end{enumerate}
\end{lemma}
\begin{proof}
  We use the following inequalities to prove
  Claims~\ref{c:conv-fitting-inequalities}
  and~\ref{c:conv-fitting-in-ofi}.  First, the inclusion
  condition of
  Assumption~\ref{a:conv-fitting-fitting-tuple} and
  $\bta-\alp=m \ccn n$ of
  Assumption~\ref{a:conv-fitting-omg-fitting-tuple} imply
  \begin{align}
    s_{3,1}&\leq \btao-\alpo \leq s_{3,2},
             \label{ineq:conv-fitting-s3-btao-alpo-s3}\\
    s_{3,1}&\leq \btat-\alpt \leq s_{3,2}.
             \label{ineq:conv-fitting-s3-btat-alpt-s3}
  \end{align}
  Second, Assumption~\ref{a:conv-fitting-fitting-tuple}
  and
  Claim~\ref{c:inc-from-fitting-inc} of
  Lemma~\ref{lem:inc-from-fitting} give some $h\in \Zgeo$
  such that
  \begin{align}
    1&\leq \alpo <\alpt \leq h,
       \label{ineq:conv-fitting-1-alpo-alpt-h}\\
    1&\leq \btao<\btat\leq h.
       \label{ineq:conv-fitting-1-btao-btat-h}
  \end{align}
  
  Let us prove Claim~\ref{c:conv-fitting-inequalities}.
  First, we prove
  inequality~\eqref{ineq:conv-fitting-max-alpo-btao-min}.
  Then, by Assumption~\ref{a:conv-fitting-gates} and
  inequalities~\eqref{ineq:conv-fitting-s3-btao-alpo-s3}, we
  have
  \begin{align}
    (s_{1,2}+\alpo )-(\btao-s_{2,2})
    &=s_{1,2}+s_{2,2}+\alpo -\btao
      =s_{3,2}-(\btao-\alpo )
      \geq 0,\\
    (\btao-s_{2,1})-
    (s_{1,1}+\alpo) 
    &=\btao-\alpo -(s_{1,1}+s_{2,1})
      =\btao-\alpo -s_{3,1}
      \geq 0.
  \end{align}
  Therefore, we obtain
  inequality~\eqref{ineq:conv-fitting-max-alpo-btao-min},
  because Assumption~\ref{a:conv-fitting-gates} implies
  \begin{align}
    (s_{1,2}+\alpo )-(s_{1,1}+\alpo )=s_{1,2}-s_{1,1}
    &\geq 0,\\
    (\btao-s_{2,1})-(\btao-s_{2,2})=s_{2,2}-s_{2,1}&\geq 0.
  \end{align}
  Second,
  inequality~\eqref{ineq:conv-fitting-max-alpt-btat-min}
  holds similarly by
  Assumption~\ref{a:conv-fitting-gates}
  and inequalities~\eqref{ineq:conv-fitting-s3-btat-alpt-s3}.
Third, we prove
inequality~\eqref{ineq:conv-fitting-gamo-gamt}.  By
Assumption~\ref{a:conv-fitting-gates} and
  inequalities~\eqref{ineq:conv-fitting-1-alpo-alpt-h}
  and~\eqref{ineq:conv-fitting-1-btao-btat-h},
   we have
  \begin{align}
    (s_{1,2}+\alpt )-(s_{1,1}+\alpo )
    &= (s_{1,2}-s_{1,1})+(\alpt -\alpo)> 0,\\
    (\btat-s_{2,1})-(\btao-s_{2,2})
    &= (s_{2,2}-s_{2,1})+ (\btat-\btao)>0.
  \end{align}
  Furthermore, Assumption~\ref{a:conv-fitting-gates} and
  inequalities~\eqref{ineq:conv-fitting-s3-btao-alpo-s3}
  and~\eqref{ineq:conv-fitting-1-btao-btat-h} give
  \begin{align}
    (\btat-s_{2,1})-(s_{1,1}+\alpo )
    = \btat-\alpo -s_{3,1}
    >\btao-\alpo -s_{3,1}
    \geq 0.
  \end{align}
  If $s_{3,2}=\infty$, then
  $\gamo=s_{1,1}+\alpo$
  if $s_{2,2}=\infty$ or
  $\gamt=\btat-s_{2,1}$ if $s_{1,2}=\infty$.
  Thus, suppose $s_{3,2}<\infty$.
  Then, Assumption~\ref{a:conv-fitting-gates} and
  inequalities~\eqref{ineq:conv-fitting-s3-btao-alpo-s3}
  and~\eqref{ineq:conv-fitting-1-alpo-alpt-h} imply
  \begin{align}
    (s_{1,2}+\alpt )-(\btao-s_{2,2})=
    s_{3,2}-\btao+\alpt 
    >s_{3,2}-\btao+\alpo 
    =s_{3,2}-(\btao-\alpo )
    \geq 0.
  \end{align}
  Thus, inequality~\eqref{ineq:conv-fitting-gamo-gamt}
  follows.
  
  Let us prove Claim~\ref{c:conv-fitting-in-ofi}.  First,
  we prove
  \begin{align}
    \gam\in \tsi(2,\tp{\alpo,\btat}).
    \label{inc:conv-fitting-gam-oi}
  \end{align}
  Since $s_{1,1},s_{2,1}\geq 0$ by
  Assumption~\ref{a:conv-fitting-gates},
  Assumption~\ref{a:conv-fitting-gamma} implies
  $\alpo \leq \gamo$ and $\gamt\leq \btat$.  Therefore,
  inclusion~\eqref{inc:conv-fitting-gam-oi} holds by
  inequality~\eqref{ineq:conv-fitting-gamo-gamt}.

  Second, we confirm that
  $\omgd(s_{1},\alp,\gam,u)=
  \tp{s_{1},l,m_{\alp,\gam},n_{\alp,\gam}, k_{\alp,\gam,u}}$
  is fitting.  By Assumption~\ref{a:conv-fitting-gamma} and
  inequalities~\eqref{ineq:conv-fitting-max-alpo-btao-min}
  and~\eqref{ineq:conv-fitting-max-alpt-btat-min},
  each $i\in \oi(2)$ gives  inequalities
  $s_{1,1}+\alp_{i} \leq \gam_{i}\leq s_{1,2}+\alp_{i}$ and
  $\bta_{i}-s_{2,2} \leq \gam_{i}\leq \bta_{i}-s_{2,1}$, which
  are equivalent to
  \begin{align}
    s_{1,1}&\leq \gam_{i}-\alp_{i} \leq s_{1,2},
             \label{ineq:conv-fitting-s1-gam-alpo-s1}\\
    s_{2,1}&\leq \bta_{i}-\gam_{i}\leq s_{2,2}.
             \label{ineq:conv-fitting-s2-bta-gam-s2}
  \end{align}
  Then,
  inequalities~\eqref{ineq:conv-fitting-s1-gam-alpo-s1}
  imply
  $m_{\alp,\gam}\ccn n_{\alp,\gam}=
  \tp{\gam_{i}-\alp_{i}}_{i\in\oi(2)} \ldZtl s_{1}$.
  Moreover,
  inequality~\eqref{ineq:conv-fitting-gamo-gamt} gives
  $n_{\alp,\gam,1}+k_{\alp,\gam,u,2} =\gamt-\alpo
  >\gamo-\alpo =m_{\alp,\gam,1}$.  Hence,
  $\omgd(s_{1},\alp,\gam,u)$ is fitting by
  inequality~\eqref{ineq:conv-fitting-1-alpo-alpt-h} and
  Lemma~\ref{lem:fitting-one}.
  
  Third, we prove that
  $\omgd(s_{2},\gam,\bta,u)
  =\tp{s_{2},l,m_{\gam,\bta},n_{\gam,\bta},
    k_{\gam,\bta,u}}$ is fitting.  By
  inequalities~\eqref{ineq:conv-fitting-s2-bta-gam-s2},
  we have
  $m_{\gam,\bta}\ccn n_{\gam,\bta}
  =\tp{\bta_{i}-\gam_{i}}_{i\in \oi(2)} \ldZtl s_{2}$.  Also,
  inequality~\eqref{ineq:conv-fitting-1-btao-btat-h} gives
  $n_{\gam,\bta,1}+k_{\gam,\bta,u,2} =\btat-\gamo
  >\btao-\gamo =m_{\gam,\bta,1}$.  Therefore,
  $\omgd(s_{2},\gam,\bta,u)$ is fitting by
  inequality~\eqref{ineq:conv-fitting-gamo-gamt}
  and
  Lemma~\ref{lem:fitting-one}.
  Claim~\ref{c:conv-fitting-in-ofi} now follows.
\end{proof}

Moreover, we rewrite
$\ofi(s_{1},s_{2},\alp,\bta,u)$ as follows.
\begin{lemma}\label{lem:eq-ofi-fitting}
  Let $\alp,\bta\in \tsi(2,h)$ for some
  $h\in \Zgeo$.  Then,
  for each  $u\in \Zgez$,
  $\gam\in \ofi(s_{1},s_{2},\alp,\bta,u)$ if and only if
  $\gam \in \tsi(2,h)$ and both $\omgd(s_{1},\alp,\gam,u)$
  and $\omgd(s_{2},\gam,\bta,u)$ are fitting.
\end{lemma}
\begin{proof}
  The only if part is clear. Hence, let us prove the if
  part.  The fitting
  $\omgd(s_{1},\alp,\gam,u)$ and
  $\omgd(s_{2},\gam,\bta,u)$
  give
  $\alpo \leq \gamo$ and
  $\gamt\leq \btat$,
  because
  $m_{\alp,\gam}=\tp{\gamo-\alpo} \ldZl s_{1}$
  and
  $n_{\gam,\bta}=\tp{\btat-\gamt}\ldZl s_{2}$
  for $s_{1}, s_{2}\geq 0$.
  Thus,
  $\gam \in \tsi(2,\tp{\alpo,\btat})$, since
  $\gamo<\gamt$ by 
  $\gam \in \tsi(2,h)$.
\end{proof}
Then, we have the following merged-log-concavity of parcel
convolutions. Its weight-zero case
applies to the convolutions of strongly
$\llq$-log-concave polynomials, which have
convolution indices by Claim~\ref{c:con-ind-sp-cases-weight-zero} of
Lemma~\ref{lem:con-ind-sp-cases}.

\begin{theorem}\label{thm:conv}
  Let $l=1$. Consider
  $\cF_{i}=\Lam(s_{i},l,w,\scc_{i},
  f_{i,s_{i}},\phi,\rho_{i},x_{i},\fX)$ for $i\in\oi(3)$ such that
  $\cFr=\cFo*\cFt$. 
    For a
    fitting $\mu=\tp{s_{3},l,m,n,k}$ and $\alpo \in \Zgeo$,
    let
    $\tp{s_{3},l,\alp,\bta}=\omgn(\mu, \alpo)$.
    Then, we obtain
  \begin{align}
    & \Delta(\cFr)(s_{3},l,w,m,n,k,\phi,\rhor,x_{3},\fX)
    \\& =
    \sum_{ \gam\in \ofi(s_{1},s_{2},\alp,\bta,\kn)}
    \Delta(\cFo) (s_{1},l,w,
    m_{\alp,\gam},n_{\alp,\gam},k_{
    \alp,\gam,\kn},\phi,\rhoo,x_{1},\fX) \\
    & \cdot
      \Delta(\cFt) (s_{2},l,w, m_{\gam,\bta},n_{\gam,\bta},
      k_{\gam,\bta,\gamo-\alpo +\kn},\phi,\rhot,x_{2},\fX).
      \label{eq:conv-ext-cb}
  \end{align}
  Moreover,
 let $\cF_{i}$ be
 $\scc'_{i}$-merged-log-concave for $i\in\oi(2)$
 with
 $O_{i}=\{\sce_{i}, \scc_{i}\}$ for $i\in\oi(3)$ and
    $O'_{i}=\{\sce'_{i}, \scc'_{i}\}$ for $i\in\oi(2)$.
    If there are squaring orders
  $O'_{3}=\{\scer',\sccr'\} \Sup O_{1}', O_{2}', O_{3}$,
  then $\cFr$ is $\sccr'$-merged-log-concave by
  \begin{align}
    \Delta(\cFr)(s_{3},l,w,m,n,k,\phi,\rhor,x_{3},\fX)
    \sccr'0.
    \label{ineq:conv-positivity-positivity}
  \end{align}
\end{theorem}
\begin{proof}
  Let us prove equation~\eqref{eq:conv-ext-cb}.  Hence,
  consider the convolution index
  $\ci=\tp{s,l,w,O,\phi,\rho,x,\fX,q,\lam,o}$ of
  $\tp{\cFo,\cFt}$.  Moreover, let $y=x_{3}^{\rhor}$ and
  $\Theta(u)=\til{\Theta}(l,w,y,\phi,\tp{u})$ for each
  $u\in \Z$.  Then, we have a ring shift index
  $\ri=\tp{\Theta,d,\iota^{3}(h),\alp,\bta,
    M_{\cFo,h},M_{\cFt,h},\Q(\fX)}$ for some $h\in \Zgeo$ by
  Claim~\ref{c:conv-key-rshift-index} of
  Proposition~\ref{prop:conv-key}.
  
  We have
  $M_{\cFr,h}(\alp,\bta,\Theta,\kn)
  =\Delta(\cFr)(s_{3},l,w,m,n,k,\phi,\rhor,x_{3},\fX)$ by the
  convolution index $\ci$ and
  Claim~\ref{c:conv-key-merged-det} of
  Proposition~\ref{prop:conv-key}, because
  $\omgd(s_{3},\alp,\bta,\kn) =\tp{s_{3},l,m,n,k}$ by
  Claim~\ref{c:eps-omg} of Proposition~\ref{prop:recip}.
  Furthermore, if
  $M_{\cFo,h}(\alp,\gam)M_{\cFt,h}(\gam,\bta)\neq 0$, then
  Claim~\ref{c:non-zero-det-fitting-inclusion} of
  Lemma~\ref{lem:non-zero-det-fitting} implies
  $\gamo-\alpo \geq 0$.  In particular,
  Claim~\ref{c:conv-key-merged-det} of
  Lemma~\ref{lem:non-zero-det-fitting} implies
  \begin{align}
    M_{\cFo,h}(\alp,\gam, \Theta,\kn)
    &= \Delta(\cFo)
      (s_{1},l,w,m_{\alp,\gam},n_{\alp,\gam},
      k_{\alp,\gam,\kn},\phi,\rhoo,x_{1},\fX),\\
    M_{\cFt,h}(\gam,\bta,\Theta,\gamo-\alpo +\kn)
    &=\Delta(\cFt)
      (s_{2},l,w,m_{\gam,\bta},n_{\gam,\bta},
      k_{\gam,\bta,\gamo-\alpo +\kn},\phi,\rhot,x_{2},\fX).
  \end{align}
  We have $M_{\cFr,h}= M_{\cFo,h}M_{\cFt,h}$, since
  $\cF_{3,\mn}=\sum_{\mt\in \Zl} \cF_{1,\mt}\cdot \cF_{2,\mn-\mt}$
  for each $\mn\in \Zl$.  Hence,
Theorem~\ref{thm:rshift-ext-cb} with
the ring shift index $\ri$  give
  \begin{dmath*}
    \Delta(\cFr)(s_{3},l,w,m,n,k,\phi,\rhor,x_{3},\fX)
    =M_{\cFr,h}(\alp,\bta,\Theta,\kn) = \sum_{\gam\in \tsi(2,h,
      \alp,\bta, M_{\cFo,h},M_{\cFt,h})}
    M_{\cFo,h}(\alp,\gam,\Theta,\kn) \cdot
    M_{\cFt,h}(\gam,\bta,\Theta,\gamo-\alpo +\kn)
    = \sum_{  \gam\in \tsi(2,h, \alp,\bta, M_{\cFo,h},M_{\cFt,h})}
    \Delta(\cFo) (s_{1},l,w,m_{\alp,\gam},n_{\alp,\gam},
    k_{\alp,\gam,\kn},\phi,\rhoo,x_{1},\fX) \cdot
    \Delta(\cFt) (s_{2},l,w,m_{\gam,\bta},n_{\gam,\bta},
    k_{\gam,\bta,\gamo-\alpo +\kn},\phi,\rhot,x_{2},\fX).
  \end{dmath*}
  This gives equation~\eqref{eq:conv-ext-cb} by
  Claim~\ref{c:non-zero-det-fitting-iff} of
  Lemma~\ref{lem:non-zero-det-fitting} and
  Lemma~\ref{lem:eq-ofi-fitting}.
  
  Let us prove
  inequality~\eqref{ineq:conv-positivity-positivity}.  Let
  $\lam=\tp{ \max(s_{1,1}+\alpo,\btao-s_{2,2}),
    \min(s_{1,2}+\alpt,\btat-s_{2,1})}\in \Zt$.  Then,
  $\lam\in \ofi(s_{1},s_{2},\alp,\bta,\kn)$ by
  Claim~\ref{c:conv-fitting-in-ofi} of
  Lemma~\ref{lem:conv-fitting}.  Also, the
  merged-log-concavity of $\cF_{1}$ and $\cF_{2}$ says
  \begin{align}
    \Delta(\cFo) (s_{1},l,w,
    m_{\alp,\lam},n_{\alp,\lam},k_{
    \alp,\lam,\kn},\phi,\rhoo,x_{1},\fX)
    &\scco'0,\\
    \Delta(\cFt)  (s_{2},l,w,
    m_{\lam,\bta},n_{\lam,\bta},
    k_{\lam,\bta,\lam_{1}-\alpo +\kn},
    \phi,\rhot,x_{2},\fX)
    &\scct'0.   
  \end{align}
  Thus, equation~\eqref{eq:conv-ext-cb} and the
  compatibility $O_{3}' \Sup O_{1}',O_{2}'$ imply
  inequality~\eqref{ineq:conv-positivity-positivity}. 
\end{proof}

\begin{example}\label{ex:conv-thm}
  We explicitly compute both sides of
  equation~\eqref{eq:conv-ext-cb} in Theorem~\ref{thm:conv}.
  Let $l=1$, $s_{1}=s_{2}=s_{3}=\tp{0,\infty}$, $w=\tp{1}$,
  $x_{1}=x_{2}=x_{3}=\tp{q}$, $\fX=\{q\}$,
  $\scco=\scct=\sccr=\llq$, and $\gam=\tp{\tp{0,0,0}}$.
  Then, suppose parcels
  $\cF_{i}=\Lam(s_{i},l,w,\scc_{i},\Psi_{s,\gam,q},x_{i},\fX)$
  for $i\in \oi(2)$.

  Moreover, let $\rhoo=\rhot=\rhor=\tp{1}$,
  $\lam=\iota^{3}(1)$, and $o=\iota^{2}(1)$. Then,
  Claims~\ref{c:con-ind-sp-cases-equators-gates-equators}
  and~\ref{c:con-ind-sp-cases-equators-gates-equators-gates}
  of Lemma~\ref{lem:con-ind-sp-cases} give the convolution
  index $\ci=\tp{s,l,w,O,\phi,\rho,x,\fX,q,\lam,o}$.  Thus,
  consider
  \begin{align}
    \cFr=\cFo*\cFt= \Lam(s_{3},l,w,\sccr,f,x_{3},\fX)
  \end{align}
  such that
  $f_{m}=\prod \phi(x_{3})^{m\rc w}\cdot
  [m]!_{x_{3}}^{w}\cdot \sum_{i+j=m}\cF_{1,i}\cF_{2, j}$ for
  each $m\ldZl s$.
  
  Furthermore, Lemma~\ref{lem:fitting-one} gives a
  fitting
  $\mu=\tp{s_{3},l,m,n,k}$ such that $m=n=\tp{1}$ and
  $k=\tp{0,1}$.  Then, $\alp=\tp{1,2}$ and $\beta=\tp{2,3}$
  satisfy $\tp{s_{3},l,\alp,\bta} =\omgn(\mu, \alpo )$.
  Thus, $\ofi(s_{1},s_{2},\alp,\bta,\kn)$ consists of
  $\gamo=\tp{1,2}$, $\gamt=\tp{1,3}$, and $\gamr=\tp{2,3}$.
  In particular,
  $\gam_{2}=\tp{\max(s_{1,1}+\alpo,\btao-s_{2,2}),
    \min(s_{1,2}+\alpt,\btat-s_{2,1})}$ as in
  Claim~\ref{a:conv-fitting-gamma} of
  Lemma~\ref{lem:conv-fitting}.
  
  The right-hand side of
  equation~\eqref{eq:conv-ext-cb} in
  Theorem~\ref{thm:conv} sums the following
  $q$-polynomials:
  \begin{align}
    &\Delta(\cFo) (s_{1},l,w,
      m_{\alp,\gamo},n_{\alp,\gamo},k_{
      \alp,\gamo,\kn},x_{1},\fX)  \\
    &\cdot \Delta(\cFt)
      (s_{2},l,w, m_{\gamo,\bta},n_{\gamo,\bta},
      k_{\gamo,\bta,\gam_{1,1}-\alpo +\kn},x_{2},\fX)   \\
    &=\Delta(\cFo) (s_{1},l,w,
      \tp{0},\tp{0},\tp{0,1},x_{1},\fX) \cdot
      \Delta(\cFt)  (s_{2},l,w,\tp{1},\tp{1},\tp{0,1},x_{2},\fX)
    \\&=q;\\
    &\Delta(\cFo) (s_{1},l,w,
      m_{\alp,\gamt},n_{\alp,\gamt},
      k_{\alp,\gamt,\kn},x_{1},\fX)  \\
    &\cdot      \Delta(\cFt)  (s_{2},l,w, m_{\gamt,\bta},
      n_{\gamt,\bta},
      k_{\gamt,\bta,\gam_{2,1}-\alpo +\kn},x_{2},\fX) \\
    &=\Delta(\cFo) (s_{1},l,w,
      \tp{0},\tp{1},\tp{0,1},x_{1},\fX)  
      \cdot      \Delta(\cFt)
      (s_{2},l,w,\tp{1},\tp{0},\tp{0,2},x_{2},\fX)\\
    &=q+1;\\
    &\Delta(\cFo) (s_{1},l,w,
      m_{\alp,\gamr},n_{\alp,\gamr},k_{
      \alp,\gamr,\kn},x_{1},\fX) \\
    &\cdot  \Delta(\cFt)  (s_{2},l,w,
      m_{\gamr,\bta},n_{\gamr,\bta},
      k_{\gamr,\bta,\gam_{3,1}-\alpo +\kn},x_{2},\fX)    \\
    &=\Delta(\cFo) (s_{1},l,w,
      \tp{1},\tp{1},\tp{0,1},x_{1},\fX)  
      \cdot   \Delta(\cFt)
      (s_{2},l,w,\tp{0},\tp{0},\tp{1,1},x_{2},\fX)\\
    &=q.
  \end{align}
  Then, we obtain $3q+1$. Therefore, this coincides
  with the left-hand side:
  \begin{align}
    \Delta(\cFr)(s_{3},l,w,m,n,k,x_{3},\fX)
    &=\frac{(1)_{q}(2)_{q}}{(0)_{q}(1)_{q}}
      \left( \cF_{3,\tp{1}}^{2}-
      \cF_{3,\tp{0}}\cF_{3,\tp{2}}\right)
    \\&=(2)_{q}\left(\left(\frac{2}{(1-q)}\right)^{2}
    - \frac{q+3}{(1-q)(1-q^{2})}\right).
  \end{align}
\end{example}

The following  gives
polynomials with positive integer coefficients
by
Theorem \ref{thm:conv}.

\begin{corollary}\label{cor:conv}
  Let $l=1$. Consider
  $\cF_{i}
  =\Lam(s_{i},l,w,\scc_{i},f_{i,s_{i}},\phi,\rho_{i},x_{i},\fX)$
  for $i\in\oi(3)$ with the parcel convolution
  $\cFr=\cFo*\cFt$.  Also, for each $i\in\oi(2)$ and fitting
  $\tp{s_{i},l,m,n,k}$, suppose $\ka_{1}\in \Qgo$ and
  $\ka_{2}\in \Q$ such that
  \begin{align}
    q^{\ka_{2}(\mn+\nn)}
    \Delta(\cF_{i})(s_{i},l,w,m,n,k,\phi,\rho_{i},x_{i},\fX)
    &>_{q^{\ka_{1}}}0.
  \end{align}
  Then, for each fitting $\tp{s_{3},l,m,n,k}$, we have
  \begin{align}
    q^{\ka_{2}(\mn+\nn)}
    \Delta(\cFr)(s_{3},l,w,m,n,k,\phi,\rhor,x_{3},\fX)
    &
      >_{q^{\ka_{1}}}0.
  \end{align}
\end{corollary}
\begin{proof}
  In equation~\eqref{eq:conv-ext-cb} of
  Theorem~\ref{thm:conv}, we have $m\ccn n=
  \bta-\alp$,
  $m_{\alp,\gam}\ccn n_{\alp,\gam} =\gam-\alp$, and
  $m_{\gam,\bta}\ccn n_{\gam,\bta} =\bta-\gam$.  Thus,
  the assertion follows, because
  $\ka_{2}(m_{\alp,\gam,1}+n_{\alp,\gam,1})
  +\ka_{2}(m_{\gam,\bta,1}+ n_{\gam,\bta,1})
  =\ka_{2}(\btao-\alpo+\btat -\alpt) =\ka_{2}(\mn+\nn)$.
\end{proof}

\section{Explicit parcel convolutions, critical points, phase transitions, and merged
  determinants}
\label{sec:examples-conv}
We compute explicit parcel convolutions of finite gates
in this section,
as ones of infinite gates appear later.

\subsection{A parcel convolution of weight
  one}
\label{sec:conv-ex-weight-one}
Let $s_{1}=\tp{1,4}$, $l=1$, $w=\tp{1}$, $\scc=\llq$,
$\rho=\tp{1}$, $x=\tp{q}$, and $\fX=\{q\}$.
Then, we have
the
 $q$-Stirling parcel
$\cFo= \Lam(s_{1},l,\scc,c_{s_{1},l,q},\fX)$ and
monomial parcel
$\cFt=\Lam(s_{1},l,w,\scc,\Psi_{s_{1},\tp{\tp{0,1,0}},q},\rho,x,\fX)$.
Moreover,
Claim~\ref{c:ext-int-hadam-more-deform} of
Corollary~\ref{cor:ext-int-hadam} gives the
$\scc$-merged-log-concave parcel
$\cFr=\Lam(s_{1},l,w,\scc,f_{3,s_{1}},\rho,x,\fX)$ such that
$f_{3,s_{1},m}= q^{\mn} c_{s_{1},l,q,m}$ for $m\ldZl s_{1}$.

Let $s_{2}=2s_{1}=\tp{2,8}$.  Then,
Claims~\ref{c:con-ind-sp-cases-equators-gates-equators}
and~\ref{c:con-ind-sp-cases-rho-mediator} of
Lemma~\ref{lem:con-ind-sp-cases} and Theorem~\ref{thm:conv}
provides the $\scc$-merged-log-concave parcel
\begin{align}
  \cF_{4}
  &=\cFr*\cFr=\Lam(s_{2},l,w,\scc,
    f_{4,s_{2}},\rho,x,\fX).
\end{align}
Explicitly, we have the following $q$-polynomials:
\begin{dgroup*}
  \begin{dmath}
    f_{4,s_{2},\tp{2}}
    =
    q^{9} + 5q^{8} + 12q^{7} + 18q^{6} + 18q^{5} + 12q^{4}
    + 5q^{3} + q^{2};
    \label{eq:conv-ex-weight-one-2}
  \end{dmath}
  \begin{dmath}
    f_{4,s_{2},\tp{3}}
    = 2q^{11} + 12q^{10} + 36q^{9} + 70q^{8} +
    94q^{7} + 90q^{6} \\
    + 60q^{5} + 26q^{4} + 6q^{3};
    \label{eq:conv-ex-weight-one-3}
  \end{dmath}
  \begin{dmath*}
    f_{4,s_{2},\tp{4}}
    = q^{14} + 7q^{13} + 27q^{12} + 70q^{11} + 133q^{10} + 191q^{9} + 212q^{8} + 183q^{7} + 120q^{6} + 55q^{5} + 15q^{4};
  \end{dmath*}
  \begin{dmath*}
    f_{4,s_{2},\tp{5}}
    = 2q^{16} + 12q^{15} + 40q^{14} + 90q^{13} + 154q^{12} + 214q^{11} + 244q^{10} + 236q^{9} + 186q^{8} + 122q^{7} + 60q^{6} + 20q^{5};
  \end{dmath*}
  \begin{dmath*}
    f_{4,s_{2},\tp{6}}
    = q^{19}+ 5q^{18}+ 18q^{17}+ 41q^{16}+ 74q^{15}+ 108q^{14}+ 139q^{13}+ 157q^{12}+ 154q^{11}+ 136q^{10}+ 101q^{9}+ 66q^{8}+ 35q^{7} +  15q^{6};
  \end{dmath*}
  \begin{dmath}
    f_{4,s_{2},\tp{7}}
    =
    2q^{21}+ 6q^{20}+ 14q^{19}+ 20q^{18}+ 32q^{17}+ 42q^{16}+ 50q^{15}+ 52q^{14}+ 54q^{13}+ 46q^{12}+ 40q^{11}+ 28q^{10}+ 18q^{9}+ 10q^{8}+
    6q^{7};
    \label{eq:conv-ex-weight-one-7}
  \end{dmath}
  \begin{dmath}
    f_{4,s_{2},\tp{8}}
    =
    q^{24}+ q^{23}+ 2q^{22}+ 3q^{21}+ 5q^{20}+ 5q^{19}+ 7q^{18}+ 7q^{17}+ 8q^{16}+ 7q^{15}+ 7q^{14}+ 5q^{13}+ 5q^{12}+ 3q^{11}+ 2q^{10}+ q^{9} + q^{8}.
    \label{eq:conv-ex-weight-one-8}
  \end{dmath}
\end{dgroup*}

\subsubsection{On critical points and phase transitions}
For $\theta=\tp{3,7}$, let
$P=r_{\tp{1,\infty},\theta}(\pi(1))=\{\tp{s_{1},l,m_{i},n_{i},k_{i}}\}_{i\ldZ
  \theta}$ such that $m_{i}=n_{i}=\tp{i}$ and
$k_{i}=\tp{0,1}$.  Then, we have the ideal merged pair
$\zeta=\tp{P,\cF_{4}}$.  First, $\zeta$ has the unique front
critical point $\FC(\zeta)=0.181093\dots$ that solves
\begin{align}
  \cF_{4,\tp{2}}(q)=\cF_{4,\tp{3}}(q)
\end{align}
by equations~\eqref{eq:conv-ex-weight-one-2}
and~\eqref{eq:conv-ex-weight-one-3}.  If
$\Bo_{1}(q)=(2)_{q}$, $\Bo_{2}(q)=(3)_{q}$,
$\To_{1}(q)=f_{4,s_{2},\tp{2}}(q)$, and
$\To_{2}(q)=f_{4,s_{2},\tp{3}}(q)$, then we obtain the front
phase transition of $\zeta$ at $\FC(\zeta)$ by
Lemma~\ref{lem:BoTo-crit-points}.  Second, $\zeta$ has the
unique rear critical point $\RC(\mu)=0.978644\dots$ that
solves
\begin{align}
  \cF_{4,\tp{7}}(q)=\cF_{4,\tp{8}}(q)
\end{align}
by equations~\eqref{eq:conv-ex-weight-one-7}
and~\eqref{eq:conv-ex-weight-one-8}.  Then,
Lemma~\ref{lem:BoTo-crit-points} gives the rear phase
transition of $\zeta$ at $\RC(\zeta)$.  Third, $\zeta$ has no
asymptotic critical points for the finite gate $\theta$.

\subsubsection{Polynomials
  with positive integer coefficients of an
  ideal merged pair}
Let $\mt=\nt=\tp{2}$ and $\kt=\tp{0,1}$.  Then, the
following $q$-polynomials are merged determinants of
$\zeta$:
\begin{dgroup*}
  \begin{dmath*}
    \Delta(\cF_{4})(s_{2},l,w,\mt,\nt,\kt,\rho,x,\fX)
    =
    q^{20}+ 11q^{19}+ 60q^{18}+ 215q^{17}+ 565q^{16}+ 1152q^{15}+ 1882q^{14}+ 2510q^{13}+ 2760q^{12}+ 2510q^{11}+ 1882q^{10}+ 1152q^{9}+ 565q^{8}+ 215q^{7}+ 60q^{6}+ 11q^{5}+ q^{4          };
  \end{dmath*}
  \begin{dmath*}
    \Delta(\cF_{4})(s_{2},l,w,\mr,\nr,\kr,\rho,x,\fX)=
    3q^{25}+ 39q^{24}+ 253q^{23}+ 1091q^{22}+ 3500q^{21}+ 8862q^{20}+ 18351q^{19}+ 31793q^{18}+ 46772q^{17}+ 58972q^{16}+ 64038q^{15}+ 59954q^{14}+ 48259q^{13}+ 33163q^{12}+ 19215q^{11}+ 9201q^{10}+ 3526q^{9}+ 1024q^{8}+ 203q^{7}+ 21q^{6};
  \end{dmath*}
  \begin{dmath*}
    \Delta(\cF_{4})(s_{2},l,w,m_{4},n_{4},k_{4},\rho,x,\fX)=
    q^{32}+ 15q^{31}+ 114q^{30}+ 584q^{29}+ 2263q^{28}+ 7054q^{27}+ 18368q^{26}+ 40957q^{25}+ 79554q^{24}+ 136253q^{23}+ 207590q^{22}+ 283136q^{21}+ 347192q^{20}+ 383691q^{19}+ 382378q^{18}+ 343199q^{17}+ 276483q^{16}+ 198741q^{15}+ 126308q^{14}+ 70013q^{13}+ 33168q^{12}+ 13010q^{11}+ 4008q^{10}+ 875q^{9} + 105q^{8}.
  \end{dmath*}
\end{dgroup*}

In particular, 
$\Delta(\cF_{4})(s_{2},l,w,m_{i},n_{i},k_{i},\rho,x,\fX)$ is a
 log-concave
 $q$-polynomial for each
$i\ldZ e(\theta)$. Also,
it is
palindromic for each $i= 2, 8$, but not for 
$i\ldZ \theta$.

\subsection{A parcel convolution of weight
  two}
\label{sec:conv-ex-weight-two}

Let $s_{1}=\tp{0,2}$, $s_{2}=\tp{1,2}$, $l=1$, $w=\tp{2}$,
$\scc=\llq$, $\rho=\tp{1}$, $x=\tp{q}$, and $\fX=\{q\}$. First, we
have the $q$-number
parcel $\cFo= \Lam(s_{2},l,\scc,\chi_{s_{2},q},\fX)$
and monomial parcels
$\cFt=\Lam(s_{2},l,w,\scc,
\Psi_{s_{2},\tp{\tp{0,0,0}},q},\rho,x,\fX)$
and  $ \cFr=\Lam(s_{1},l,w,\scc,
\Psi_{s_{1},\tp{\tp{0,1,0}},q},\rho,x,\fX)$,
 which are $\scc$-merged-log-concave.
 Then, Theorem~\ref{thm:merged-int-hadam} gives the
 $\scc$-merged-log-concave 
 $\cF_{4}
  =\cFo\dd_{\tp{1,1}}\cFt
    =\Lam(s_{2},l,w,\scc, \chi_{s_{2},q},\rho,x,\fX)$.
    Moreover,
    let $s_{3}=s_{1}+s_{2}=\tp{1,4}$.  Then,
Claims~\ref{c:con-ind-sp-cases-equators-gates-equators}
and~\ref{c:con-ind-sp-cases-rho-mediator} of
Lemma~\ref{lem:con-ind-sp-cases} and Theorem~\ref{thm:conv}
yields the 
$\scc$-merged-log-concave 
\begin{align}
  \cF_{5}
  &=\cFr*\cF_{4}
    =\Lam(s_{3},l,w,\scc,f_{5,s_{3}},\rho,x,\fX)
    \end{align}
with the following $q$-polynomials:
\begin{align}
  f_{5,s_{3},\tp{1}}&= 1;
                   \label{eq:conv-ex-weight-two-1}
  \\
  f_{5,s_{3},\tp{2}}&= q^{3}+ 2q^{2}+ 2q+ 1;
                   \label{eq:conv-ex-weight-two-2}
  \\
  f_{5,s_{3},\tp{3}}&= 2q^{6}+ 5q^{5}+ 8q^{4}+ 7q^{3}+ 4q^{2}+ q;
                   \label{eq:conv-ex-weight-two-3}
  \\
  f_{5,s_{3},\tp{4}}&= q^{11}+ 3q^{10}+ 7q^{9}+ 11q^{8}+ 14q^{7}+ 14q^{6}+ 11q^{5}+ 7q^{4}+ 3q^{3} + q^{2}.
                   \label{eq:conv-ex-weight-two-4}
\end{align}

\subsubsection{On critical points and phase transitions}
For $\theta=\tp{2,3}$ and
$P=r_{\tp{1,\infty},\theta}(\pi(1))$, we have the ideal merged
pair $\zeta=\tp{P,\cF_{5}}$.
First, $\zeta$ has no front critical points, because
$0<q<1$ does not solve
\begin{align}
  \cF_{5,\tp{1}}(q)=\cF_{5,\tp{2}}(q)
\end{align}
by equations~\eqref{eq:conv-ex-weight-two-1}
and~\eqref{eq:conv-ex-weight-two-2}.  Second, $\zeta$ has the
unique rear critical point $\RC(\zeta)=0.618034\dots$ that
solves
\begin{align}
  \cF_{5,\tp{3}}(q)=\cF_{5,\tp{4}}(q)
\end{align} 
by equations~\eqref{eq:conv-ex-weight-two-3}
and~\eqref{eq:conv-ex-weight-two-4}.  Then, we have the rear
phase transition by Lemma~\ref{lem:BoTo-crit-points}.
Third, there are no asymptotic critical points of $\zeta$
for the finite gate $\theta$.

\subsubsection{Polynomials with positive
  integer coefficients of an ideal merged pair}
Let $P=\{\tp{s_{1},l,m_{i},n_{i},k_{i}}\}_{i\ldZ \theta}$
with $\mn=\nn=\tp{1}$,
$m_{4}=n_{4}=\tp{4}$, and $\kn=k_{4}=\kt$.
Then, we have the following merged determinants of 
$\zeta$:
\begin{dgroup*}
  \begin{dmath*}
    \Delta(\cF_{5})(s_{3},l,w,\mn,\nn,\kn,\rho,x,\fX)
    =q^{2}+ 2q+ 1;
  \end{dmath*}
  \begin{dmath*}
    \Delta(\cF_{5})(s_{3},l,w,\mt,\nt,\kt,\rho,x,\fX)=
    q^{10}+ 6q^{9}+ 17q^{8}+ 31q^{7}+ 41q^{6}+ 42q^{5}+ 35q^{4}+ 24q^{3}+ 13q^{2}+ 5q+ 1;
  \end{dmath*}
  \begin{dmath*}
    \Delta(\cF_{5})(s_{3},l,w,\mr,\nr,\kr,\rho,x,\fX)=
    3q^{18}+ 21q^{17}+ 81q^{16}+ 219q^{15}+ 456q^{14}+ 768q^{13}+ 1074q^{12}+ 1266q^{11}+ 1266q^{10}+ 1074q^{9}+ 768q^{8}+ 456q^{7}+ 219q^{6}+ 81q^{5}+ 21q^{4}+ 3q^{3};
  \end{dmath*}
  \begin{dmath*}
    \Delta(\cF_{5})(s_{3},l,w,m_{4},n_{4},k_{4},\rho,x,\fX)=
    q^{30}+ 8q^{29}+ 38q^{28}+ 132q^{27}+ 369q^{26}+ 870q^{25}+ 1782q^{24}+ 3232q^{23}+ 5260q^{22}+ 7754q^{21}+ 10423q^{20}+ 12836q^{19}+ 14527q^{18}+ 15136q^{17}+ 14527q^{16}+ 12836q^{15}+ 10423q^{14}+ 7754q^{13}+ 5260q^{12}+ 3232q^{11}+ 1782q^{10}+ 870q^{9}+ 369q^{8}+ 132q^{7}+ 38q^{6}+ 8q^{5} + q^{4}.
  \end{dmath*}
\end{dgroup*}

In particular, 
$\Delta(\cF_{5})(s_{3},l,w,m_{i},n_{i},k_{i},\rho,x,\fX)$
is a log-concave
$q$-polynomial for each $i\ldZ e(\theta)$ and
palindromic for each $i\ldZ e(\theta)$ except $i\neq 2$.

\section{Graded monomial products}
\label{sec:mono-conv-eta-prod-multinom}

We discuss graded monomial products
 in Definition~\ref{def:monomial-conv} with
certain deformations of eta products.

\subsection{Merged-log-concavity of
   graded
  monomial products}
\label{sec:monomial-conv}
We first prove that all monomial convolutions are generating
functions of merged-log-concave parcels by parcel
convolutions.  Then, we prove the same for all graded
monomial products by external Hadamard products.

For a multimonomial index $\tp{d,w,\alp,\bta,\gam}$, the
monomial convolution involves the change of variables
$q\mapsto q^{\alp_{\lam}}$.  Hence, we discuss the following
change of parcel parameters.

\begin{proposition}
  \label{prop:merged-param-change}
  Let $l\in \Zgeo$ and $\ka,\rhoo,\rhot\in
  \Zgeo^{l}$. Consider a parcel
  $\cF=\Lam(s,l,w,\scco,\fs,\phi,\rhoo,x_{1},\fX)$ and
  $x_{2}=x_{1}^{\ka}$.  Assume the following.
  \begin{enumerate}
  \item\label{a:merged-param-change-sq}
    There are squaring
    orders
    $O_{2}=\{\scet,\scct\}\Sup O_{1}=\{\sceo,\scco\}$.
  \item \label{c:merged-param-change-med-x}
    $\phi$ is a
    $\tp{s,l,w,\scct,\ka,x_{1},\fX}$-mediator.
  \item \label{c:merged-param-change-med-y}
    $\phi$ is a
    $\tp{s,l,w,\scct,\rhot,x_{2},\fX}$-mediator.
  \end{enumerate}
  Then, we have the following.
  \begin{enumerate}[label=(\alph*)]
  \item \label{c:merged-param-change-parcel-existence}
    There
    exists a parcel
    $\cG=\Lam(s,l,w,\scct,\gs,\phi,\rhot,x_{2},\fX)=\cF$ such
    that
    $g_{s,m}=f_{m} \cdot B(s,l,w,m,\phi,\ka,x_{1},\fX)$ for each
    $m\ldZl s$.
  \item \label{c:merged-param-change-merged}
    For squaring orders $O'=\{\sce',\scc'\}\Sup
    O_{2}$,
    let $\mu_{\scc',i}=\tp{s,l,w,\scc',\phi,\rho_{i},x_{i},\fX}$
    and $\mu_{\sce',i}=\tp{s,l,w,\sce',\phi,\rho_{i},x_{i},\fX}$
    for $i\in \oi(2)$.
    Assume
    $\ka\rc \rhot=\rhoo$.
    Then, $\cF$ is
    $\mu_{\scc',2}$-merged-log-concave
    if and only if $\cF$ is
    $\mu_{\scc',1}$-merged-log-concave.
    Also, $\cF$ is
    $\mu_{\sce',2}$-merged-log-concave
     if and
     only if $\cF$ is
     $\mu_{\sce',1}$-merged-log-concave.
  \end{enumerate}
\end{proposition}
\begin{proof}
  We prove
Claim~\ref{c:merged-param-change-parcel-existence}.  For
  each $m\ldZl s$,
  \begin{align}
    \cF_{m}
    =\frac{f_{s,m}}{\prod \phi(x_{2})^{m\rc w}\cdot [m]!_{x_{2}}^{w}}
    \frac{\prod \phi(x_{2})^{m\rc w}\cdot [m]!_{x_{2}}^{w}}{
    \prod \phi(x_{1})^{m\rc w}\cdot [m]!_{x_{1}}^{w}}
    =\frac{g_{s,m}}{\prod \phi(x_{2})^{m\rc w}\cdot [m]!_{x_{2}}^{w}}.
  \end{align}
  Also, $\gs$ is pairwise $\tp{s,l,\scct}$-positive by
  Assumptions~\ref{a:merged-param-change-sq}
  and~\ref{c:merged-param-change-med-x}, since $\fs$
  is
  pairwise $\tp{s,l,\scco}$-positive.  Hence,
  Assumption~\ref{c:merged-param-change-med-y}
  gives
  Claim~\ref{c:merged-param-change-parcel-existence}.

  Claim~\ref{c:merged-param-change-merged}
  holds by
  Claim~\ref{c:merged-param-change-parcel-existence},
  because
   $x_{2}^{\rhot}=x_{1}^{\rhoo}$ implies
  $\Ups(s,l,w,m,n,k,\phi,\rhoo,x_{1},\fX)
  =\Ups(s,l,w,m,n,k,\phi,\rhot,x_{2},\fX)$
    for each $m,n\in \Zl$ and $k\in \Ztl$.
\end{proof}

We introduce the following notations of multimonomial
indices. We use quadratic polynomials $t_{\alp,\gam}$ and
rational numbers $u_{\alp,\gam}$ in
Definition~\ref{def:intro-quad-poly} and graded monomial
products $M(d,w,\alp,\bta, \gam, q,z)$ and monomial
convolutions $\cM(d,w,\alp,\bta, \gam, q,v)$ in
Definition~\ref{def:multimonomial-ind}.

\begin{definition}\label{def:succ-mono-conv}
  Suppose a multimonomial index
  $\tp{d,w,\alp,\bta,\gam}$.
  \begin{enumerate}
  \item Let
    $\fUq{\alp,\gam,q}=\left\{q^{u_{\alp,\gam}}\right\}$.
  \item On $\Q( \fUq{\alp,\gam,q})$, we put squaring
    orders
    $\tgq{\alp,\gam,q}=>_{q^{u_{\alp,\gam}}}$ and
    $\tgeq{\alp,\gam,q}=\geq_{q^{u_{\alp,\gam}}}$ if
    $t_{\alp,\gam}(m)\in \Q_{\geq 0}$ for each
    $m \in \Zdgez$, or
    $\tgq{\alp,\gam,q}=>_{q^{\pm u_{\alp,\gam}}}$ and
    $\tgeq{\alp,\gam,q}=\geq_{q^{\pm u_{\alp,\gam}}}$
    otherwise.
  \item Suppose $\bta=\iota^{d}(1)$. Then, we call
  $\tp{d,w,\alp,\gam}$ reduced multimonomial index (or
  multimonomial index for short).  Also, let
  \begin{align}
    M(d,w,\alp, \gam, q,z)
    &=
      M(d,w,\alp,\bta, \gam, q,z),\\
    \cM(d,w,\alp, \gam, q,v)
    &=\cM(d,w,\alp,\bta, \gam, q,v).
  \end{align}
  We call $M(d,w,\alp, \gam, q,z)$ and
  $\cM(d,w,\alp, \gam, q,v)$ reduced graded monomial
  product
  and reduced monomial convolution (or graded
  monomial
  product and monomial convolution for short).
  \end{enumerate}
\end{definition}

In particular, for each monomial convolution
$\cM(d,w,\alp,\bta,\gam,q,z)$, there exists a reduced
monomial convolution
$\cM(d',w,\alp',\gam',q,z)=\cM(d,w,\alp,\bta,\gam,q,z)$.

We use the following notation of 
 parcels to write graded monomial products
and monomial convolutions. The notation generalizes one in
equation~\eqref{eq:intro-gen-functions}.
\begin{definition}\label{def:gen-fun-parcel}
  Let $\cF=\Lam(s,l,w,\scc,\fs,\phi,\rho,x,\fX)$ and
  $v_{1},\dots,v_{l}$ be  indeterminates.
  \begin{enumerate}
  \item We write the generating function
    $\cZ_{\cF}(v_{1},\dots,v_{l})$ of $\cF$ such that
    \begin{align}
      \cZ_{\cF}(v_{1},\dots,v_{l})=
      \sum_{\lam\in\Zlgez}\cF_{\lam}\cdot \prod_{i\in \oi(l)} v_{i}^{\lam_{i}}
      \in \Q(\fX)[[v_{1},\dots,v_{l}]].
    \end{align}
  \item In particular, if $s=\tp{0,\infty}$ and
    $\cF=\Lam(s,l,w,\scc,\Psi_{s,\gam,q},\rho,x,\fX)$, then let
    \begin{align}
      \cZ_{w,\gam,q}(v_{1},\dots,v_{l})
      = \cZ_{\cF}(v_{1},\dots,v_{l}).      
    \end{align}
    \end{enumerate}
  \end{definition}
  
Then, we obtain the merged-log-concavity of all
monomial convolutions.  For parcels $\cFo,\dots,\cF_{d}$,
we write $*_{i\in\oi(d)} \cF_{i}$ for the parcel convolution
$\cFo*\dots* \cF_{d}$.

\begin{theorem}\label{thm:mono-conv}
  Consider a multimonomial index
  $\tp{d,w,\alp,\bta,\gam}$ with
   $\del\in \Zgeo$ and $\ka\in \Zgeo^{d}$ such that
  $\ka \rc \alp=\iota^{d}(\del)$.  Let
  $s=\tp{0,\infty}$, $l=1$,
$\scc=\tgq{\alp,\gam,q}$,
$\rho\in \Zgeo^{l}$,  $y=\tp{q^{\del}}$, and
  $\fX=\fUq{\alp,\gam,q}$.
  For each $i\in\oi(d)$, put $x_{i}=\tp{q^{\alp_{i}}}$ and
  $\cF_{i}= \Lam(s,l,w,\scc,\Psi_{s,\tp{\gam_{i}},q^{\alp_{i}}},
  p_{i},x_{i},\fX)$.  Then, we have the following.
  \begin{enumerate}
  \item \label{c:mono-conv-conv} There is the
    $\scc$-merged-log-concave parcel
    $\cG=\Lam(s,l,w,\scc,\fs, \rho,y,\fX)
    =*_{i\in\oi(d)}\cF_{i}^{* \bta_{i}}$.
  \item For an indeterminate $v$,
    $\cZ_{\cG}(v)= \cM(d,w,\alp,\bta,\gam,q,v)$.
    \label{c:mono-conv-mono-conv}
\end{enumerate}
\end{theorem}
\begin{proof}
  Since Claim~\ref{c:mono-conv-mono-conv} follows from
  Claim~\ref{c:mono-conv-conv}, we prove
  Claim~\ref{c:mono-conv-conv}.  Without loss of generality,
  assume a reduced multimonomial index
  $\tp{d,w,\alp,\gam}$.  First, let $d=1$.  Then, we have
  $y=x_{1}^{\ka}$ and
  $\cF_{1}=
  \Lam(s,l,w,\scc,\Psi_{s,\tp{\gam_{1}},q^{\alp_{1}}}, \phi,
  p_{1},x_{1},\fX)$ for the $l$-canonical mediator
  $\phi$. Hence, Claim~\ref{c:merged-param-change-parcel-existence} of
  Proposition~\ref{prop:merged-param-change} gives the
  parcel
  \begin{align}
    \cG=\Lam(s,l,w,\scc,\fs,\rho,y,\fX)=\cF_{1}.    
  \end{align}
  By Theorem~\ref{thm:monomial-poly},
  $\cG$ is
  $\tp{s,l,w,\scc,\phi,p_{1},x_{1},\fX}$-merged-log-concave.
  Hence,
  Claim~\ref{c:mono-conv-conv} holds,
  since $\cG$ is
  $\tp{s,l,w,\scc,\phi,\rho,y,\fX}$-merged-log-concave by
  Claim~\ref{c:merged-param-change-merged} of
  Proposition~\ref{prop:merged-param-change}.

  Second, let $d\geq 2$. Then, the induction gives the
  $\scc$-merged-log-concave parcel
  \begin{align}
    \cH=\Lam(s,l,w,\scc,\hs, \rho,y,\fX)=*_{i\in\oi(d-1)} \cF_{i}.
  \end{align}
  Also, we have
  $\cF_{d}=
  \Lam(s,l,w,\scc,\Psi_{s,\tp{\gam_{i}},q^{\alp_{i}}},\phi,
  \ka_{d}\rho ,x_{i},\fX)$ for the canonical $l$-mediator
  $\phi$.  Let $\lam=\tp{\del,\alp_{d},\del}$,
  $\tau=\tp{\rho,\ka_{d}\rho,\rho}$,
  $O=\tp{O_{i}=\{\sce,\scc\}}_{i\in \oi(3)}$, and
  $o=\tp{1,\ka_{d}}$.  Hence, $\tp{\cH,\cF_{d}}$ carries the
  convolution index
  $\tp{\iota^{3}(s),l,w,O,\phi,\tau,\tp{y,x_{d},y},\fX,q,\lam,o}$ by the
  exponent equator $o_{1}\lamo= o_{2}\lamt=\del=\lamr$, the
  base-shift equator
  $o_{1}^{-1}\tau_{1}=o_{2}^{-1}\tau_{2}=\rho=\tau_{3}$, and
  Claim~\ref{c:con-ind-sp-cases-mediators} of
  Lemma~\ref{lem:con-ind-sp-cases}.  Theorem~\ref{thm:conv}
  now gives the $\scc$-merged-log-concave parcel
  \begin{align}
    \cG
    =\Lam(s,l,w,\scc,\fs,\rho,y,\fX)=
    \cH*\cF_{d}=*_{i\in\oi(d)}\cF_{i}.
  \end{align}
\end{proof}
Furthermore, we prove the  merged-log-concavity
of graded monomial products.
\begin{theorem}\label{thm:merged-graded-monom-prod}
  Consider a graded monomial product
  $M(d,w,\alp,\bta, \gam, q,z)$.  Let $s=\tp{0,\infty}$,
  $\scc=\tgq{\alp,\gam,q}$,
  $\rho\in \Zgeo^{d}$, $x=\tp{q^{\alp_{i}}}_{i\in \oi(d)}$,
  and  $\fX=\fUq{\alp,\gam,q}$.
  Then, we have the merged-log-concave parcel
  $\cF= \Lam(s,d,w,\scc,f_{s},\rho,x,\fX)$ such that
  \begin{align}
    \cZ_{\cF}(z_{1},\dots,z_{d})= M(d,w,\alp,\bta, \gam,
    q,z).
    \label{eq:merged-graded-monom-prod}
  \end{align}

\end{theorem}
\begin{proof}
  For each $i\in \oi(d)$, Theorem~\ref{thm:mono-conv} gives
  the merged-log-concave 
  $\cH_{i}= \Lam(s,1,w,\scc,f_{i,s},\rho,x_{i},\fX)$ such that
  $\cM(1,w,\tp{\alp_{i}},\tp{\bta_{i}},\tp{\gam_{i}},
  q,z_{i})=\cZ_{\cF_{i}}(z_{i})$.  Then,
  Theorem~\ref{thm:merged-ext-hadam} yields the
  merged-log-concave 
  $\cF= \Lam(s,d,w,\scc,f_{s},\rho,x,\fX) =\cH_{1}\ehd \dots
  \ehd \cH_{d}$ in
  equation~\eqref{eq:merged-graded-monom-prod}.
\end{proof}

\subsection{Graded monomial products
  and eta products}\label{sec:etaprod}

We first realize the eta function $\eta(\tau)$ in
Definition~\ref{def:eta-fun-prod} and its inverse by
weight-one linear and quadratic monomial parcels.

\begin{lemma}\label{lem:euler}
  Let $s=\tp{0,\infty}$, $l=1$, and $w=\tp{1}$.  Consider
  $\ka_{1},\ka_{2}\in \Q$ such that $\ka_{1}+\ka_{2}=1$.
  Let
  $\cFo=\Lam(s,l,w,\scc,
  \Psi_{s,\tp{\tp{0,\ka_{1},-\frac{1}{24}}},q},\rho,x,\fX)$
  and
  $\cFt=\Lam(s,l,w,\scc,
  \Psi_{s,\tp{\tp{\fraa,\frac{\ka_{1}}{2},
        \frac{1}{24}}},q},\rho,x,\fX)$.
  Then,
  $\cZ_{\cFo}(q^{\ka_{2}})
      =\eta(\tau)^{-1}$ and
      $\cZ_{\cFt}(-q^{\frac{\ka_{2}}{2}})
  =\eta(\tau)$,
  provided
       $q=e^{2\pi \i \tau}$ of
  $\Ima \tau>0$.
\end{lemma}
\begin{proof}
  Since $\eta(\tau)=q^{\frac{1}{24}}(q;q)$, Euler's
  binomial identities
   imply
  $\eta(\tau)^{-1} =q^{-\frac{1}{24}}\sum_{\lam\in\Zgez}
  \frac{q^{\lam}}{(\lam)_{q}}$ and
  $\eta(\tau) =q^{\frac{1}{24}}\sum_{\lam\in\Zgez}
  \frac{q^{\frac{\lam(\lam-1)}{2}}}{(\lam)_{q}}(-q)^{\lam}
  =q^{\frac{1}{24}}\sum_{\lam\in\Zgez} (-1)^{\lam}
  \frac{q^{\frac{\lam(\lam+1)}{2}}}{(\lam)_{q}}$.  Thus,
  \begin{align}
      \eta(\tau)^{-1}
      &=q^{-\frac{1}{24}}\sum_{\lam\in\Zgez}
    \frac{q^{\lam(\ka_{1}+\ka_{2})}}{(\lam)_{q}}
    =
    q^{-\frac{1}{24}}\sum_{\lam\in\Zgez}
    \frac{q^{\lam\ka_{1}}}{(\lam)_{q}}
    (q^{\ka_{2}})^{\lam}
    =
    \cZ_{\cFt}(q^{\ka_{2}}),\\
    \eta(\tau)
    &=q^{\frac{1}{24}}\sum_{\lam\in\Zgez}
      (-1)^{\lam}
      \frac{q^{\frac{\lam^{2}}{2}+\lam\frac{\ka_{1}+\ka_{2}}{2}}}
      {(\lam)_{q}}
      =
      q^{\frac{1}{24}}\sum_{\lam\in\Zgez}
      \frac{q^{\frac{\lam^{2}}{2}+\lam\frac{\ka_{1}}{2}}}
      {(\lam)_{q}}
      (-q^{\frac{\ka_{2}}{2}})^{\lam}
      =
      \cZ_{\cFo}(-q^{\frac{\ka_{2}}{2}}).
  \end{align}
\end{proof}
Hence, we have the following by $\etg(\bta,\ka)$ and
$J(z,q,\alp,\bta,\ka)$ in
Definition~\ref{def:intro-hat-bta-gam}.  Then, graded
monomial products
$M(d,w, \alp,\abs{\bta}, \etg(\bta,\ka), q,z)$ give the
merged-log-concavity of infinitely-many
$J(z,q,\alp,\bta,\ka)$-analogs of each eta product
$\ep_{d,\alp,\bta}(\tau)$ by choices of $\ka\in \Qd$.

\begin{proposition}\label{prop:lifted-eta-pro}
  Let $d\in \Zgeo$, $w=\tp{1}$, $\alp\in \Zgeo^{d}$,
  $\bta\in \Zd_{\neq 0}$, and $\ka\in \Qd$.  Then, we have the
  following.
  \begin{enumerate}
  \item
    \label{c:lifted-eta-pro-parcel}
    There exists a width-$d$ merged-log-concave parcel $\cF$
    such that
    $\cZ_{\cF}(z_{1},\dots,z_{d})=M(d,w, \alp,\abs{\bta},
    \etg(\bta,\ka), q,z)$.
    \item
    \label{c:lifted-eta-pro-analog}
    For each $q=e^{2\pi \i \tau}$ of $\Ima(\tau)>0$,
    $M(d,w, \alp,\abs{\bta}, \etg(\bta,\ka), q,z)$ is a
    $J(z,q,\alp,\bta,\ka)$-analog of 
    $\ep_{d,\alp,\bta}(\tau)$.
  \end{enumerate}
\end{proposition}
\begin{proof}
  Claim~\ref{c:lifted-eta-pro-parcel} holds by
  Theorem~\ref{thm:merged-graded-monom-prod}, since
  $\tp{1,w,\tp{\etg(\bta,\ka)_{i}}}$ is a monomial index for
  each $i\in \oi(d)$.

  Let us prove Claim~\ref{c:lifted-eta-pro-analog}.
  Consider $J(z,q,\alp,\bta,\ka)\mapsto \iota^{d}(1)$ as $z$ varies
  in $\bC^{d}$.  Then, when
  $J(z,q,\alp,\bta,\ka)=\iota^{d}(1)$, each $z_{i}$ has to be
  $q^{(1-\ka_{i})\alp_{i}}$ if $ \etg(\bta,\ka)_{i,1}=0$, or
  $ -q^{\frac{(1-\ka_{i})\alp_{i}}{2}}$ if
  $\etg(\bta,\ka)_{i,1}=\fraa$.  Hence, Lemma~\ref{lem:euler}
  gives Claim~\ref{c:lifted-eta-pro-analog}, since
  $\cZ_{w,\tp{\etg(\bta,\ka)_{i}},q^{\alp_{i}}}(
  q^{(1-\ka_{i})\alp_{i}}) =\eta(\tau^{\alp_{i}})^{-1}$ if
  $\etg(\bta,\ka)_{i,1}=0$ and
  $\cZ_{w,\tp{\etg(\bta,\ka)_{i}},q^{\alp_{i}}}
  (-q^{\frac{(1-\ka_{i})\alp_{i}}{2}})= \eta(\tau^{\alp_{i}})$
  if $\etg(\bta,\ka)_{i,1}=\fraa$.
\end{proof}

\begin{example}\label{ex:gradedmonomprod}
  Let $s=\tp{0,\infty}$, $d=3$, $w=\tp{1}$,
  $\alp_{1}=\tp{1,1,1}$, $\alp_{2}=\tp{2,1,1}$,
  $\bta=\tp{1,2,2}$, $\rho=\iota^{d}(1)$, and
  $\gam=\tp{\tp{0,0,0},\tp{0,\fraa,0}, \tp{0,0,0}}$.  Also,
  let $x_{i}=\iota^{d}(q)^{\alp_{i}}$,
  $\scc_{i}=\tgq{\alp_{i},\gam,q}$, and
  $\fX_{i}=\fUq{\alp_{i},\gam,q}$ for $i\in \oi(2)$.  Then,
  Claim~\ref{c:lifted-eta-pro-parcel} in
  Proposition~\ref{prop:lifted-eta-pro} gives width-$d$
  parcels
  $\cF_{i}= \Lam(s,d,w,\scc_{i},f_{s},\rho,x_{i},\fX_{i})$
  such that
  $\cZ_{\cF_{i}}(z_{1},\dots,z_{d})=M(d,w,
  \alp_{i},\bta,\gam, q,z)$ for $i\in \oi(2)$.

  Hence, suppose
  $\bP_{s,d,\xi,h} =\{\bP_{s,d,\xi,h,i}
  =\tp{s,d,m_{i},n_{i},k_{i}}\}_{i\ldZ \theta}$ in
  Example~\ref{ex:f-y-2} with
  $k_{i}=\tp{0,1,0,2,1,0}$.  For
  instance,
  $\cF_{1,m_{1}} =
  \frac{4q^{\fraa}}{(q^2-1)(q^4-1)(q-1)^{2}}$.  Moreover, we
  have the following unimodal $q$-polynomial:
    \begin{dmath*}
    \Delta(\cF_{1})(s,d,w,m_{1},n_{1},k_{1},\rho,x_{1},\fX_{1})
    =
    q^{38}+2 q^{37}+7 q^{36}+12 q^{35}+25 q^{34}+38 q^{33}+63 q^{32}+88 q^{31}+128 q^{30}+168 q^{29}+221 q^{28}+274 q^{27}+331 q^{26}+388 q^{25}+437 q^{24}+486 q^{23}+515 q^{22}+544 q^{21}+544 q^{20}+544 q^{19}+515 q^{18}+486 q^{17}+437 q^{16}+388 q^{15}+331 q^{14}+274 q^{13}+221 q^{12}+168 q^{11}+128 q^{10}+88 q^9+63 q^8+38 q^7+25 q^6+12 q^5+7 q^4+2 q^3+q^2.
  \end{dmath*}
  However,
  the following is not a unimodal $q$-polynomial:
  \begin{dmath*}
    \Delta(\cF_{2})(s,d,w,m_{1},n_{1},k_{1},\rho,x_{2},\fX_{2})
    =
    q^{30}+q^{29}+5 q^{28}+5 q^{27}+14 q^{26}+13 q^{25}+29 q^{24}+25 q^{23}+49 q^{22}+40 q^{21}+70 q^{20}+54 q^{19}+86 q^{18}+62 q^{17}+92 q^{16}+62 q^{15}+86 q^{14}+54 q^{13}+70 q^{12}+40 q^{11}+49 q^{10}+25 q^9+29 q^8+13 q^7+14 q^6+5 q^5+5 q^4+q^3+q^2.
    \end{dmath*}

  \end{example}
  Since we are interested in the unimodality and eta
  products, we conjecture the following (see also
  Conjecture~\ref{conj:unimodal-width-two} for the
  width-one
  case
  $M(1,\tp{1},\tp{1},\tp{1},\tp{\tp{0,0,0}},q,\tp{z_{1}}
  )$).
\begin{conjecture}\label{conj:gradedmonomprod}
  Under the notation of Example~\ref{ex:gradedmonomprod},
  $\Delta(\cF_{1})(s,d,w,m_{i},n_{i},k_{i},\rho,x_{1},\fX_{1})$ is
  a unimodal $q$-polynomial for each $i\in \Zgeo$.

\end{conjecture}

\subsection{Weighted
  \texorpdfstring{$q$}{q}--multinomial
  coefficients}
\label{sec:explicit-monom-conv}

We introduce the following  notion of
 weighted $q$-multinomial coefficients
\footnote{The author could not find a reference
  to the notion in literature.
  He also asked to the audience during
  his talk at International Workshop P-positivity in Matroid
  Theory and Related Topics of RIMS 2021.}.
This is
to describe  monomial convolutions, which construct graded
monomial products.

\begin{definition}
\label{def:weighted-gauss}
  Let $d\in \Zgeo$, $\alp\in \Zgeo^{d}$, and
  $\del=\lcm(\alp)$.  Suppose $i\in \Z$ and $j\in \Zd$.
  Then, the following is a weighted $q$-multinomial
  coefficient:
  \begin{align}
    {i \brack j}_{\alp,q}
    =
    \begin{dcases}
      \frac{(i)_{q^{\del}}}
      { \prod_{\lam\in\oi(d)}(j_{\lam})_{q^{\alp_{\lam}}}}
      \mif j\geq 0 \mand \sum j=i,\\
      0 \melse.
    \end{dcases}
  \end{align}
  If $d=2$, then ${i \brack j}_{\alp,q}$ is a weighted
  $q$-binomial coefficient.
\end{definition}
In particular, if $\del=1$, then
${i \brack j}_{\alp,q} ={i \brack j}_{q}$.  Moreover, we
employ the following notation.
\begin{definition}
  For $d\in \Zgeo$, $\lam\in\oi(d)$, and $j\in \Zd$,
  let
  $o(j,\lam)=
  \tp{j_{1},\dots,j_{\lam}-1,\dots, j_{d}}\in
  \Zd$.
\end{definition}
For instance, $j\in \Zd$ and $i\in \Zgeo$ give the
$q$-Pascal identity
${i \brack j}_{q} = \sum_{\lam\in\oi(d)} q^{
  \sum_{\ka\in\oi(\lam-1)} j_{\ka}} {i-1 \brack
  o(j,\lam)}_{q}$. We extend this by weighted
$q$-multinomial coefficients.

\begin{proposition}\label{prop:weighted-gauss}
  Let $d\in \Zgeo$, $\alp\in \Zgeo^{d}$,
  $\del_{l}=\lcm(\alp)$, and $\del_{g}=\gcd (\alp)$.  Assume
  $j\in \Zd$.  Then, when $i\in \Zgeo$, we have a weighted
  $q$-Pascal identity
  \begin{align}
      {i \brack j}_{\alp,q}
      = \sum_{\lam\in\oi(d)}
      q^{\del_{l} \sum_{\ka\in\oi(\lam-1)} j_{\ka}}
      {i-1 \brack o(j,\lam)}_{\alp,q}
      \left[\frac{\del_{l}}{\alp_{\lam}}\right]_{q^{\alp_{\lam}
          j_{\lam}}}.
      \label{eq:weighted-gauss-weighted-q-pascal}  
    \end{align}
    In particular,
    ${i \brack j}_{\alp,q}\geq_{q^{\del_{g}}}0$
    for each $i\in \Z$.
\end{proposition}
\begin{proof}
  Assume $j\in \Zgez^{d}$ and $\sum j=i$; otherwise, the
  weighted $q$-Pascal identity holds by $0=0$.  Then, by the
  unweighted
  $q$-Pascal identity, we have
  \begin{align}
    {i \brack j}_{\alp,q}
    &=
    {i \brack j}_{q^{\del_{l}}}
    \cdot \frac{\prod_{\lam\in\oi(d)}(j_{\lam})_{q^{\del_{l}}}}
    {\prod_{\lam\in\oi(d)}(j_{\lam})_{q^{\alp_{\lam}}}}\\
      &=
    \sum_{\lam\in\oi(d)}
    q^{\del_{l} \sum_{\ka\in\oi(\lam-1)} j_{\ka}}
    {i-1 \brack o(j,\lam)}_{q^{\del_l}}
    \cdot 
    \frac{\prod_{\lam\in\oi(d)}(j_{\lam})_{q^{\del_{l}}}}
    {\prod_{\lam\in\oi(d)}(j_{\lam})_{q^{\alp_{\lam}}}}.
  \end{align} 
  On each summand above, if $j_{\lam}< 1$, then
  ${i-1 \brack o(j,\lam)}_{q^{\del_{l}}}=0= {i-1 \brack
    o(j,\lam)}_{\alp,q}$. If not, then since each
  $\alp_{\lam}$ divides $\del_{l}$, we have
  \begin{align}
    {i-1 \brack o(j,\lam)}_{q^{\del_{l}}}
    \cdot
    \frac{\prod_{\lam\in\oi(d)}(j_{\lam})_{q^{\del_{l}}}}
    {\prod_{\lam\in\oi(d)}(j_{\lam})_{q^{\alp_{\lam}}}}
    &=
    \frac{(i-1)_{q^{\del_{l}}}}{
      \prod_{\ka\in\oi(d)}(o(j,\lam)_{\ka})_{q^{\del_{l}}}}
    \cdot 
    \frac{\prod_{\ka\in\oi(d)}(j_{\ka})_{q^{\del_{l}}}}{
      \prod_{\ka\in\oi(d)}(j_{\ka})_{q^{\alp_{\ka}}}}\\
    &=
    {i-1 \brack o(j,\lam)}_{\alp,q}
    \cdot
    \frac{ (1-q^{\del_{l}j_{\lam}})}{
      (1-q^{\alp_{\lam}j_{\lam}})}\\
    & =
    {i-1 \brack o(j,\lam)}_{\alp,q}
    \cdot
    \left[\frac{\del_{l}}{\alp_{\lam}}
    \right]_{q^{\alp_{\lam}j_{\lam}}}.
  \end{align}
  Thus, we obtain the weighted $q$-Pascal identity.

  Since $ {i \brack j}_{\alp,q}=1$ or $0$ when
  $i\in \Z_{\leq 0}$, the weighted $q$-Pascal identity gives the
  latter statement by the induction on $i\in \Zgez$.
\end{proof}

\begin{example}\label{ex:weighted-gauss}
  Let $\alp_{\lam}=\tp{1,\lam,\lam}$ for $\lam\in \oi(2)$.
  Then, we have the following unimodal $q$-polynomials:
  \begin{align}
    {3 \brack \tp{1,1,1}}_{\alpo,q}
    &  =\frac{(3)_{q}}{(1)_{q}(1)_{q}(1)_{q}}
      =q^{3}+ 2q^{2}+ 2q+ 1;\\
    {3 \brack \tp{1,1,1}}_{\alpt,q}
    &=\frac{(3)_{q^{2}}}{(1)_{q}(1)_{q^2}(1)_{q^2}}
      =
      q^{7}+ q^{6}+ 2q^{5}+ 2q^{4}
      + 2q^{3}+ 2q^{2}+ q+ 1.
    \end{align}
  However,
  ${10 \brack \tp{4,3,3}}_{\alp_{2},q}$ is
  \begin{dmath*}
    q^{76}+ q^{75}+ 3q^{74}+ 4q^{73}+ 9q^{72}+ 11q^{71}+ 21q^{70}+ 26q^{69}+ 43q^{68}+ 53q^{67}+ 80q^{66}+ 97q^{65}+ 138q^{64}+ 165q^{63}+ 221q^{62}+ 262q^{61}+ 336q^{60}+ 392q^{59}+ 485q^{58}+ 559q^{57}+ 668q^{56}+ 761q^{55}+ 884q^{54}+ 993q^{53}+ 1125q^{52}+ 1248q^{51}+ 1380q^{50}+ 1512q^{49}+ 1637q^{48}+ 1769q^{47}+ 1879q^{46}+ 2004q^{45}+ 2088q^{44}+ 2198q^{43}+ 2252q^{42}+ 2336q^{41}+ 2355q^{40}+ 2409q^{39}+ 2390q^{38}+ 2409q^{37}+ 2355q^{36}+ 2336q^{35}+ 2252q^{34}+ 2198q^{33}+ 2088q^{32}+ 2004q^{31}+ 1879q^{30}+ 1769q^{29}+ 1637q^{28}+ 1512q^{27}+ 1380q^{26}+ 1248q^{25}+ 1125q^{24}+ 993q^{23}+ 884q^{22}+ 761q^{21}+ 668q^{20}+ 559q^{19}+ 485q^{18}+ 392q^{17}+ 336q^{16}+ 262q^{15}+ 221q^{14}+ 165q^{13}+ 138q^{12}+ 97q^{11}+ 80q^{10}+ 53q^{9}+ 43q^{8}+ 26q^{7}+ 21q^{6}+ 11q^{5}+ 9q^{4}+ 4q^{3}+ 3q^{2} + q + 1,
  \end{dmath*}
  which is not a unimodal $q$-polynomial by
  $\dots +2409q^{39}+ 2390q^{38}+ 2409q^{37}+\cdots$.
  
  Hence, it would be interesting to determine $\bta\in
  \Zgeo^{d}$
  and $i\in \Zgeo$ so
  that ${i\brack j}_{\bta,q}$ is a unimodal $q$-polynomial
  for each $j\in \Zgeo^{d}$ and describe how these
  parameters change
  the unimodality of ${i\brack j}_{\bta,q}$,
  as that of
  ${i\brack j}_{q}$ is  important
  \cite{Cay,Oha,Syl}.
  For instance, if
  $i=30$ and $\bta=\tp{1,4,8}$, then one can check
  that
  ${i\brack j}_{\bta,q}$ is a unimodal $q$-polynomial for
  each $j\in \Zgeo^{3}$.
\end{example}

\subsection{Monomial convolutions
and merged determinants}

We write monomial convolutions by the weighted
$q$-multinomial coefficients and the following rational
functions.
\begin{definition}\label{def:small-monom-num}
  Consider a multimonomial index
  $\tp{d,w,\alp,\gam}$.  We
  define the $t$-monomials
  $\psi_{\alp,\gam,q}=\{\psi_{\alp,\gam,q,j}\in
  \Q(\fUq{\alp,\gam,q})\}_{j\in \Zd}$ such that
  $\psi_{\tp{\alp_{\lam}},\tp{\gam_{\lam}},q,\tp{j_{\lam}}}
  =\Psi_{\tp{0,\infty},\tp{\gam_{\lam}},q^{\alp_{\lam}},\tp{j_{\lam}}}$
  and 
  $\psi_{\alp,\gam,q,j}= \prod_{\lam\in\oi(d)}
  \psi_{\tp{\alp_{\lam}},\tp{\gam_{\lam}},q,\tp{j_{\lam}}}$.
\end{definition}
In particular,
$\psi_{\alp,\gam,q,j}=
  \prod_{\lam\in\oi(d)}
      q^{\alp_{\lam}(
      \gam_{\lam,1}j_{\lam}^{2}+
      \gam_{\lam,2}j_{\lam}
       +              \gam_{\lam,3})}$
     if  $j\in \Zgez^{d}$.

\begin{proposition}\label{prop:monom-multinom} 
  Suppose a multimonomial index $\tp{d,w,\alp,\gam}$ with
  $\del=\lcm(\alp)$. Then, we have the following.
  \begin{enumerate}
  \item \label{c:monom-multinom-explicit} 
    $\cM(d,w,\alp,\gam,q,v)
      =
      \sum_{i\in \Zgez}
      \frac{\sum_{j\in \Zd}\psi_{\alp,\gam,q,j}
      {i \brack j}^{w_{1}}_{\alp,q}}
      {(i)^{w_{1}}_{q^{\del}}}\cdot v^{i}$.
    \item \label{c:monom-multinom-parcel} Let
      $s=\tp{0,\infty}$,
    $l=1$, and $\rho\in \Zgeo^{l}$.
    Then, there
    is the $\scc$-merged-log-concave parcel
    $\cF
    =
    \Lam(s,l,w,\tgq{\alp,\gam,q},
    \fs,\rho,\tp{q^{\del}},\fUq{\alp,\gam,q})$
    such that
    $f_{s,m}= \sum_{ j\in \Zd} \psi_{\alp,\gam,q,j} {\mn \brack
      j}^{w_{1}}_{\alp,q}$
    and
 \begin{align}
   \cZ_{\cF}(v)=\cM(d,w,\alp,\gam,q,v).   
 \end{align}
  \end{enumerate}
\end{proposition}
\begin{proof} 
  We have
  $\cM(d,w,\alp,\gam,q,v)
  =\prod_{\lam\in\oi(d)}
  \left(
    \sum_{i\in\Zgez}
    \frac{\psi_{\tp{\alp_{\lam}},\tp{\gam_{\lam}},q,\tp{i}}}
  { (i)_{q^{\alp_{\lam}}}^{w_{1}}}\cdot v^{i}\right)$.
  Thus, we obtain
  Claim~\ref{c:monom-multinom-explicit}, because
  \begin{align}
    \cM(d,w,\alp,\gam,q,v)
    &=\sum_{i\in\Zgez}\left(
      \sum_{j\in \Zgez^{d},\sum j=i}
      \frac{\psi_{\alp,\gam,q,j}}{
      \prod_{\lam\in\oi(d)} (j_{\lam})^{w_{1}}_{q^{\alp_{\lam}}}
      }\right) v^{i}
    \\&=\sum_{i\in\Zgez}\left(
    \sum_{j\in \Zgez^{d},\sum j=i}
    \psi_{\alp,\gam,q,j} \cdot
    \frac{ (i)_{q^{\del}}^{w_{1}}}{
    \prod_{\lam\in\oi(d)} (j_{\lam})^{w_{1}}_{q^{\alp_{\lam}}}
    }\cdot \frac{1}{(i)_{q^{\del}}^{w_{1}}}
    \right) v^{i}
    \\&=\sum_{i\in\Zgez}
    \frac{\sum_{j\in \Zgez^{d},\sum j=i} \psi_{\alp,\gam,q,j}
    {i \brack j}^{w_{1}}_{\alp,q}}{(i)_{q^{\del}}^{w_{1}}}
    \cdot v^{i}.
  \end{align}
  
  Claim~\ref{c:monom-multinom-parcel} follows from
  Theorem~\ref{thm:mono-conv} and
  Claim~\ref{c:monom-multinom-explicit}.
\end{proof}

We define the following merged determinants
of monomial convolutions.
For
 $d\in \Zgeo$ and
$\alp\in \Zdgeo$, let
$\lcm(\alp)$ and $\gcd(\alp)$ be the least common multiple
and greatest common divisor of $\alp_{1},\dots,\alp_{d}$.

\begin{definition}\label{def:delta-eta-pro}
  Let $l=1$, $\rho\in \Zlgeo$, $m,n\in \Zl$, and
  $k\in \Z^{2l}$.
  Consider a multimonomial index
  $\tp{d,w,\alp,\bta,\gam}$ with
   $x=\tp{q^{\lcm(\alp)}}$ and
  $\fX=\fUq{\alp,\gam,q}$. Then,
  in $\Q(\fX)$,
  we define the merged determinant
  \begin{align}
    &\Delta(d,w,\alp,\bta,\gam,
        m,n,k,\rho,q)
    \\&=
    \Ups(s,l,w,m,n,k,\rho,x,\fX)
    \det
    \begin{bmatrix}
      \cM(d,w,\alp,\bta,\gam,q)_{\mn}   &
                                          \cM
      (d,w,\alp,\bta,\gam,q)_{\nn+\kt}\\
                                          \cM
        (d,w,\alp,\bta,\gam,q)_{\mn-\kt}
                                              &
                              \cM
                             (d,w,\alp,\bta,\gam,q)_{\nn}\\
      \end{bmatrix}.
    \end{align}
    We put
    $\Delta(d,w,\alp, \gam,m,n,k,\rho,q)
    =
    \Delta(d,w,\alp,\iota^{d}(1), \gam,m,n,k,\rho,q)$
    for a reduced
    multimonomial
    index $\tp{d,w,\alp,\gam}$. 
  \end{definition}
  Moreover, we define the following multimonomial
  indices.
  \begin{definition}
    \label{def:tame-multimonomial-index}
    We call a multimonomial index
    $\tp{d,w,\alp,\bta,\gam}$
    tame, if $\gam_{i,1}+\gam_{i,2}\in \Zgez$ and
    $\gam_{i,3}=0$ for each $i\in \oi(d)$.
  \end{definition}
  This realizes the following squaring orders
  $\tgq{\alp,\gam,q}$.
  \begin{lemma}\label{lem:scc-alp-gam-q}
    Suppose a tame multimonomial index
       $\tp{d,w,\alp,\bta,\gam}$ with
    $\delta=\gcd(\alp)$.
    Then,
    $\tgq{\alp,\gam,q}=>_{q^{\delta}}$.
  \end{lemma}
  \begin{proof}
    By the monomial conditions,
    $\gam=\tp{
      \tp{\frac{\lam_{i}}{2},-\frac{\lam_{i}}{2}+
         \lam_{i}',0}}_{i\in\oi(d)}$
     for some $\lam_{i},\lam_{i}'\in \Zgez$.
     Thus, the assertion holds by
    $t_{\alp,\gam}(m)=
    \sum_{i\in \oi(d)}\frac{m_{i}\lam_{i}(m_{i}\lam_{i}-1)}{2}
    +m_{i}\lam_{i}'\in \Zgez$ for $m\ldZd s$.
  \end{proof}

Then, we obtain the following  strict
inequality and $q$-polynomials with
positive integer coefficients by monomial convolutions and
weighted $q$-multinomial coefficients.

\begin{theorem}\label{thm:monom-conv-mult}
  Suppose a multimonomial index $\mu=\tp{d,w,\alp,\gam}$.
  Let $\del=\lcm(\alp)$, $s=\tp{0,\infty}$, $l=1$,
  $\scc=\tgq{\alp,\gam,q}$, and $\rho\in \Zgeo^{l}$.  Then,
  for each fitting $\tp{s,l,m,n,k}$ with $a=\nu(k)$ and
  $b=\nu(m,n,k)$, we have the following strict inequality:
  \begin{align}
    &\Delta(d,w,\alp, \gam,m,n,k,\rho,q)
    \\&=
    {b \brack a}_{q^{\del\rhoo}}^{w_{1}}
    \cdot \frac{(\mn)_{q^{\del\rhoo}}^{w_{1}}}
    {(\mn)_{q^{\del}}^{w_{1}}}
    \cdot \frac{(\nn)_{q^{\del\rhoo}}^{w_{1}}}
    {(\nn)_{q^{\del}}^{w_{1}}}
    \cdot
    \sum_{j_{1}\in \Zd}
    \psi_{\alp,\gam,q,j_{1}}
    {\mn \brack j_{1}}^{w_{1}}_{\alp,q}
    \cdot
    \sum_{j_{2}\in \Zd}
    \psi_{\alp,\gam,q,j_{2}}
    {\nn \brack j_{2}}^{w_{1}}_{\alp,q}
    \\&-
    {b \brack a^{\ve}}_{q^{\del\rhoo}}^{w_{1}}
    \cdot \frac{(\mn-\kt)_{q^{\del\rhoo}}^{w_{1}}}
    {(\mn-\kt)_{q^{\del}}^{w_{1}}}
    \cdot \frac{(\nn+\kt)_{q^{\del\rhoo}}^{w_{1}}}
    {(\nn+\kt)_{q^{\del}}^{w_{1}}}
    \\              &\cdot
                      \sum_{j_{1}\in \Zd}
                      \psi_{\alp,\gam,q,j_{1}}
                      {\mn +\kt\brack j_{1}}_{\alp,q}^{w_{1}}
                      \cdot
                      \sum_{j_{2}\in \Zd}
                      \psi_{\alp,\gam,q,j_{2}}
                      {\nn-\kt \brack j_{2}}^{w_{1}}_{\alp,q}
    \\        &       \scc 0.
  \end{align}

  In particular, if $\mu$ is tame, then
  $\Delta(d,w,\alp, \gam,m,n,k,\rho,q)$ is a
  $q^{\del_{g}}$-polynomial with positive integer
  coefficients.
\end{theorem}
\begin{proof}
  Claim~\ref{c:monom-multinom-parcel} of
  Proposition~\ref{prop:monom-multinom} provides the
  $\scc$-merged-log-concave parcel
  $\cF =
  \Lam(s,l,w,\scc,\fs,\rho,\tp{q^{\del}},\fUq{\alp,\gam,q})$
  such that $\cM(d,w,\alp,\gam,q,v) =\cZ_{\cF}(v)$.
  Hence,  the strict inequality holds by
  Claim~\ref{c:merged-binom-bshift-left-right-equations}
  of
  Theorem~\ref{thm:merged-binom-bshift}.
  Furthermore, the
  latter assertion follows from the strict inequality,
  Proposition~\ref{prop:weighted-gauss}, and
  Lemma~\ref{lem:scc-alp-gam-q}.
\end{proof}
In particular, let $k=\tp{0,1}$, $a=\nu(k)$,
$b=\nu(\tp{h},\tp{h},k)$, $\alp=\iota^{d}(1)$, and
$w=\rho=\tp{1}$ in Theorem~\ref{thm:monom-conv-mult}.  Then,
we obtain Theorem~\ref{thm:intro-poly-positivity} in
Section~\ref{sec:intro}, since $a=\tp{0,1}$ and
$b=\tp{h,h+1}$ imply
${b \brack a}_{q} = {h \brack 0}_{q} {h+1 \brack 1}_{q}
=[h+1]_{q}$ and
${b \brack a^{\ve}}_{q} ={h \brack 1}_{q} {h+1 \brack 0}_{q}
=[h]_{q}$.

\begin{example}\label{ex:lifted-eta-1}
  Consider the multimonomial index
  $\tp{d,w,\alp,\gam}$ such
  that $d=2$, $w=\tp{1}$, $\alp=\tp{2,3}$, and
  $\gam=\tp{\tp{0,0,0},\tp{0,1,0}}$.  Then,
  $\tgq{\alp,\gam,q}=\llq$ by
  Lemma~\ref{lem:scc-alp-gam-q}.
  Also, by
  $x=\tp{q^{\lcm(\alp)}}$,
  we have $(m)_{x}^{w} \cM(d,w,\alp,\gam,q)_{m}\llq 0$ for
  each $m\in \Zgez$.  For example,
  \begin{dgroup*}
    \begin{dmath*}
      (0)_{x}\cM(d,w,\alp,\gam,q)_{0}=
      1,
    \end{dmath*}
    \begin{dmath*} 
      (1)_{x}
      \cM(d,w,\alp,\gam,q)_{1}=
      q^{6}+ q^{4}+ q^{3}+ q^{2}+ 1,
    \end{dmath*}
    \begin{dmath*}
      (2)_{x} \cM(d,w,\alp,\gam,q)_{2}
      =
      q^{16}+ q^{15}+ q^{14}+ q^{13}+ 3q^{12}+ q^{11}+ 2q^{10}+ 2q^{9}+ 3q^{8}+ q^{7}+ 3q^{6}+ q^{5}+ 2q^{4}+ q^{3}+ q^{2}+ 1,
    \end{dmath*}
    \begin{dmath*}
      (3)_{x}  \cM(d,w,\alp,\gam,q)_{3}
      =
      q^{31}+ q^{30}+ q^{29}+ 2q^{28}+ 3q^{27}+ 3q^{26}+ 3q^{25}+ 5q^{24}+ 4q^{23}+ 6q^{22}+ 5q^{21}+ 7q^{20}+ 5q^{19}+ 9q^{18}+ 5q^{17}+ 8q^{16}+ 6q^{15}+ 8q^{14}+ 4q^{13}+ 8q^{12}+ 4q^{11}+ 6q^{10}+ 4q^{9}+ 5q^{8}+ 2q^{7}+ 4q^{6}+ q^{5}+ 2q^{4}+ q^{3}+ q^{2} + 1.
    \end{dmath*}
  \end{dgroup*}
  
  Furthermore, let $k=\tp{0,1}$ and $\rho=\tp{1}$.  Then,
  Theorem~\ref{thm:monom-conv-mult} gives the following
  $q$-polynomials with positive coefficients:
  \begin{dgroup*}
    \begin{dmath*}
      \Delta(d,w,\alp,\gam,\tp{0},\tp{0},k,\rho,q)
      =  1;
    \end{dmath*}
    \begin{dmath*}
      \Delta(d,w,\alp,\gam,\tp{1},\tp{1},k,\rho,q)
      =  q^{18}+ q^{16}+ q^{15}+ 2q^{14}+ q^{13}+ 3q^{12}+ q^{11}+ 3q^{10}+ 2q^{9}+ 2q^{8}+ q^{7}+ 3q^{6}+ q^{5}+ q^{4}+ q^{3}+ q^{2};
    \end{dmath*}
    \begin{dmath*}
      \Delta(d,w,\alp, \gam,\tp{2},\tp{2},k,\rho,q)=
      q^{44}+ q^{43}+ 2q^{42}+ 2q^{41}+ 5q^{40}+ 4q^{39}+ 7q^{38}+ 7q^{37}+ 13q^{36}+ 10q^{35}+ 17q^{34}+ 16q^{33}+ 22q^{32}+ 19q^{31}+ 28q^{30}+ 22q^{29}+ 31q^{28}+ 28q^{27}+ 32q^{26}+ 27q^{25}+ 37q^{24}+ 27q^{23}+ 32q^{22}+ 28q^{21}+ 31q^{20}+ 22q^{19}+ 28q^{18}+ 19q^{17}+ 22q^{16}+ 16q^{15}+ 17q^{14}+ 10q^{13}+ 13q^{12}+ 7q^{11}+ 7q^{10}+ 4q^{9}+ 5q^{8}+ 2q^{7}+ 2q^{6}+ q^{5}+ q^{4};
    \end{dmath*}                                                         
    \begin{dmath*}
      \Delta(d,w,\alp, \gam,\tp{3},\tp{3},k,\rho,q)
      =
      q^{80}+ q^{79}+ 2q^{78}+ 3q^{77}+ 6q^{76}+ 6q^{75}+ 12q^{74}+ 14q^{73}+ 21q^{72}+ 25q^{71}+ 37q^{70}+ 40q^{69}+ 56q^{68}+ 64q^{67}+ 83q^{66}+ 92q^{65}+ 118q^{64}+ 129q^{63}+ 158q^{62}+ 171q^{61}+ 206q^{60}+ 216q^{59}+ 259q^{58}+ 271q^{57}+ 310q^{56}+ 320q^{55}+ 373q^{54}+ 367q^{53}+ 418q^{52}+ 420q^{51}+ 464q^{50}+ 452q^{49}+ 508q^{48}+ 478q^{47}+ 529q^{46}+ 504q^{45}+ 539q^{44}+ 499q^{43}+ 550q^{42}+ 490q^{41}+ 525q^{40}+ 475q^{39}+ 502q^{38}+ 435q^{37}+ 467q^{36}+ 396q^{35}+ 416q^{34}+ 352q^{33}+ 364q^{32}+ 296q^{31}+ 313q^{30}+ 248q^{29}+ 252q^{28}+ 200q^{27}+ 205q^{26}+ 152q^{25}+ 156q^{24}+ 115q^{23}+ 115q^{22}+ 82q^{21}+ 81q^{20}+ 55q^{19}+ 56q^{18}+ 36q^{17}+ 34q^{16}+ 22q^{15}+ 21q^{14}+ 12q^{13}+ 11q^{12}+ 6q^{11}+ 6q^{10}+ 3q^{9}+ 2q^{8}+ q^{7} + q^{6}.
    \end{dmath*}
  \end{dgroup*}
\end{example}

In Example~\ref{ex:lifted-eta-1},
$\Delta(d,w,\alp,\gam, \tp{i},\tp{i},k,\rho,q)$ are not
log-concave
  $q$-polynomials for $i\in\oi(3)$.
Even 
$\Delta(d,w,\iota^{2}(1),\gam,\tp{1},\tp{1},k,\rho,q)= q^{3} + q^{2}
+ 2q$ is not a log-concave
 $q$-polynomial
either.  However, we state the
following conjecture.

\begin{conjecture}\label{conj:delta-eta-1}
  Let  $\lamo,\lamt \in\Zgeo$ and
  $\lamr\in \Zgez$.  Consider the multimonomial index
  $\tp{d,w,\alp,\bta_{\lamo},\gam}$ such that $d=2$,
  $w=\tp{1}$, $\alp=\iota^{d}(1)$,
  $\bta_{\lamo}=\iota^{d}(\lamo)$, and
  $\gam=\tp{\tp{\fraa,-\fraa,0},\tp{0,1,0}}$.  Then,
  \begin{align}
    \Delta(d,w,\alp,\bta_{\lamo},
    \gam,\tp{\lamr},\tp{\lamr},\tp{0,1},\tp{\lamt},q)
  \end{align}
  is a log-concave $q$-polynomial.
\end{conjecture}

\begin{remark}\label{rmk:eta-merged-lc}
  Suppose $w=\tp{1}$ and a tame $\mu$ in
  Theorem~\ref{thm:monom-conv-mult}.  Then, we first
  obtain the
  $q^{\frac{\del_{g}}{12}}$-polynomial
  $q^{\frac{\sum\alp}{12}}\cdot \Delta(d,w,\alp, \gam,m,n,k,\rho,q)$
  with positive integer coefficients.  In particular,
  a
  multimonomial index
  $\tp{d,w,\alp, \tp{
      \tp{\fraa,\gam_{i,2},\frac{1}{24} }}_{i\in
      \oi(d)}}$
  such that $\fraa+\gam_{i,2}\in \Zgez$
  yields the
  $q^{\frac{\del_{g}}{24}}$-polynomial for the eta product
  $\prod_{i\in\oi(d)} \eta(\alp_{i}\tau)$ in
  Proposition~\ref{prop:lifted-eta-pro}.
  
  Second, we obtain the Laurent
  $q^{\frac{\del_{g}}{12}}$-polynomial
  $q^{ -\frac{\sum \alp}{12}}\cdot \Delta(d,w,\alp, \gam,m,n,k,\rho,q)$
  with positive integer coefficients.  Hence, a
  multimonomial index
  $\tp{d,w,\alp, \tp{\tp{0,\gam_{i,2},-\frac{1}{24}}}_{i\in
      \oi(d)}}$ such that $\gam_{i,2}\in \Zgez$ gives the
  Laurent $q^{\frac{\del_{g}}{12}}$-polynomial for the eta
  product $\prod_{i\in\oi(d)} \eta(\alp_{i}\tau)^{-1}$ in
  Proposition~\ref{prop:lifted-eta-pro}.
\end{remark}

\subsection{Examples
and conjectures
  for merged determinants
  of monomial convolutions}

We consider multimonomial indices $\tp{d,w,\alp,\gam}$
such that $\gam=\iota^{d}(\tp{\fraa,-\fraa,0})$,
$\iota^{d}(\tp{0,0,0})$, or $\iota^{d}(\tp{0,1,0})$.

\begin{example}
  Let $d\in \Zgeo$ and $\alp\in \Zgeo^{d}$.  Then,
  $\gamo=\iota^{d}(\tp{0,0,0})$ and
  $\gamt=\iota^{d}(\tp{0,1,0})$ give the 
  same eta products $\prod_{i\in\oi(d)}\eta(\alp_{i}\tau)^{-1}$
   in
   Remark~\ref{rmk:eta-merged-lc}.  But, at the
   level of merged determinants, 
   $\frac{\Delta(d,w,\alp, \gamo,m,n,k,\rho,q)}{
    \Delta(d,w,\alp, \gamt,m,n,k,\rho,q)}$ are not
  even $q$-monomials in the following.
  
  Let $s=\tp{0,\infty}$, $l=1$, $d=2$, $w=\tp{1}$,
  $\alp=\tp{1,2}$, $k=\tp{0,1}$, and $\rho=\tp{1}$.  Also, let
  $m_{i}=\tp{i}\in \Zlgez$ so that $\tp{s,l,m_{i},m_{i},k}$ is
  fitting for each $i\in \Zgez$.  First, $\gamo$ gives 
  \begin{dgroup*}
    \begin{dmath*}
      \Delta(d,w,\alp, \gamo,\mn,\mn,k,\rho,q)=
      1 + 2 q + 3 q^{2} + 2 q^{3} + q^{4},
    \end{dmath*}
    \begin{dmath*}
      \Delta(d,w,\alp, \gamo,\mt,\mt,k,\rho,q)
      =
      1 + 2 q + 6 q^{2}+ 8 q^{3}+ 13 q^{4}
      + 14 q^{5}+ 11 q^{6}+ 12 q^{7}
      + 7 q^{8}+  4 q^{9}+ 3 q^{10},
    \end{dmath*}
    \begin{dmath*}
      \Delta(d,w,\alp, \gamo,\mr,\mr,k,\rho,q)
      =
      1 + 2 q + 6 q^{2}+ 12 q^{3}+ 21 q^{4}
      + 32 q^{5}+ 49 q^{6}+ 58 q^{7}+ 
      69 q^{8}+ 78 q^{9}+ 77 q^{10}+ 76 q^{11}
      + 68 q^{12}+ 58 q^{13}+ 44 q^{14}+ 
      34 q^{15}+ 22 q^{16}+ 12 q^{17}+ 7 q^{18}
      + 2 q^{19}+ q^{20}.
    \end{dmath*}
  \end{dgroup*}
  Second,  $\gamt$ gives 
  \begin{dgroup*}
    \begin{dmath*}
      \Delta(d,w,\alp, \gamt,\mn,\mn,k,\rho,q)=
      2 q^{3}+ 2 q^{4}+ 2 q^{5}+ 3 q^{6},
    \end{dmath*}
    \begin{dmath*}
      \Delta(d,w,\alp, \gamt,\mt,\mt,k,\rho,q)
      =
      3 q^{6}+ 4 q^{7}+ 7 q^{8}+ 12 q^{9}
      + 11 q^{10}+ 14 q^{11}+ 13 q^{12}+ 
      8 q^{13}+ 6 q^{14}+ 2 q^{15}+ q^{16},
    \end{dmath*}
    \begin{dmath*}
      \Delta(d,w,\alp, \gamt,\mr,\mr,k,\rho,q)=
      4 q^{9}+ 6 q^{10}+ 12 q^{11}+ 23 q^{12}+ 30 q^{13}+ 46 q^{14}+ 58 q^{15}+ 
      69 q^{16}+ 76 q^{17}+ 80 q^{18}+ 78 q^{19}+ 67 q^{20}+ 60 q^{21}+ 
      44 q^{22}+ 32 q^{23}+ 22 q^{24}+ 12 q^{25}+ 7 q^{26}+ 2 q^{27} + q^{28}.
    \end{dmath*}
  \end{dgroup*}
\end{example}

We put
the following differences by merged
determinants of monomial convolutions.

\begin{definition}
  Suppose a multimonomial index
  $\tp{d,w,\alp,\bta,\gam}$. Let $k=\tp{0,1}$,
  $\del=\lcm(\alp)$,
  and $\lamo,\lamt,\lamr\in \Zgez$.
  Then, in $\Q[q^{\pm u_{\alp,\gam}}]$,
  we put
  \begin{align}
    \xi(d,w,\alp,\bta,
    \gam,\lamo,\lamt,\lamr,\rho,q)             
    &=\Delta(d,w,\alp,\bta, \gam,
      \tp{\lamo+\lamr},\tp{\lamo+\lamt+\lamr},k,\rho,q)
      \\&      - \Delta(d,w,\alp,\bta,
      \gam,\tp{\lamo},\tp{\lamo+\lamt},k,\rho,q).
    \end{align}
  In particular, if $\bta=\iota^{d}(1)$, then let
  $\xi(d,w,\alp, \gam,\lamo,\lamt,\lamr,\rho,q)
      =
      \xi(d,w,\alp,\bta, \gam,\lamo,\lamt,\lamr,\rho,q)$.
\end{definition}

We then conjecture the following positivity of $\xi$.
\begin{conjecture}\label{conj:merged-diff-q-pos}
  Let $d\in \Zget$, $w,\rho\in \Zgeo^{1}$, and
  $\alp\in \Zgeo^{d}$.  Consider the multimonomial
  index $\tp{d,w,\alp,\gam}$ such that
  $\gam=\iota^{d}(\tp{0,0,0})$.  Then, for each
  $\lamo,\lamt\in \Zgez$, we have
  \begin{align}
    \xi(d,w,\alp, \gam,\lamo,\lamt,1,\rho,q)\llq 0.
  \end{align}
\end{conjecture}

\begin{example}\label{ex:merged-diff-q-pos}
  Suppose the multimonomial index
  $\tp{d,w,\alp,\gam}$
  such that $d=2$, $w=\tp{1}$,
  $\alp=\tp{1,1}$, and
  $\gam=\iota^{d}(\tp{0,0,0})$.
  Let $k=\tp{0,1}$
  and $\rho=\tp{1}$.  Then,
  \begin{dgroup*}
    \begin{dmath*}
      \Delta(d,w,\alp, \gam,\tp{0},\tp{0},k,\rho,q)=
      1,
    \end{dmath*}
    \begin{dmath*}
      \Delta(d,w,\alp, \gam,\tp{1},\tp{1},k,\rho,q)=
      3q + 1,
    \end{dmath*}
    \begin{dmath*}
      \Delta(d,w,\alp, \gam,\tp{2},\tp{2},k,\rho,q)=
      q^{4}+ 3q^{3}+ 8q^{2}+ 3q+ 1,
    \end{dmath*}
    \begin{dmath*}
      \Delta(d,w,\alp, \gam,\tp{3},\tp{3},k,\rho,q)=
      3q^{7}+ 5q^{6}+ 12q^{5}+ 14q^{4}+ 18q^{3}+
      8q^{2}+ 3q + 1,
    \end{dmath*}
    \begin{dmath*}
      \Delta(d,w,\alp, \gam,\tp{4},\tp{4},k,\rho,q)=
      q^{12}+ 3q^{11}+ 12q^{10}+ 18q^{9}+ 30q^{8}+ 39q^{7}+ 46q^{6}+ 39q^{5}+ 38q^{4}+ 18q^{3}+ 8q^{2} + 3q + 1.
    \end{dmath*}
  \end{dgroup*}
  Thus, the following are $q$-polynomials
  with positive integer coefficients:
  \begin{dgroup*}
    \begin{dmath*}
      \xi(d,w,\alp, \gam,0,0,1,\rho,q)= 3q;
    \end{dmath*}
    \begin{dmath*}
      \xi(d,w,\alp, \gam,1,0,1,\rho,q)=
      q^{4}+ 3q^{3}+ 8q^{2};
    \end{dmath*}
    \begin{dmath*}
      \xi(d,w,\alp, \gam,2,0,1,\rho,q)=
      3q^{7}+ 5q^{6}+ 12q^{5}+ 13q^{4}+ 15q^{3};
    \end{dmath*}
    \begin{dmath*}
      \xi(d,w,\alp, \gam,3,0,1,\rho,q)=
      q^{12}+ 3q^{11}+ 12q^{10}+ 18q^{9}+ 30q^{8}+
      36q^{7}+
      41q^{6}+ 27q^{5} + 24q^{4}.
    \end{dmath*}
  \end{dgroup*}
\end{example}

\begin{example}
  Conjecture~\ref{conj:merged-diff-q-pos} does not
  extend to
  the case of $d=1$.
  Suppose the multimonomial index $\tp{d,w,\alp,\gam}$
  such that $d=1$, $w=\tp{1}$, $\alp=\tp{2}$, and
  $\gam=\tp{\tp{0,0,0}}$.  Let $k=\tp{0,1}$ and $\rho=\tp{1}$.  Then, we
  have
  $\xi(d,w,\alp, \gam,2,1,1,\rho,q)\not\llq 0$, since
  \begin{align}
    \Delta(d,w,\alp, \gam,\tp{2},\tp{3},k,\rho,q)
    &=  q^6 + q^4,\\
    \Delta(d,w,\alp, \gam,\tp{3},\tp{4},k,\rho,q)
    &=  q^8 + q^6.
  \end{align}
\end{example}

\begin{example}\label{ex:merged-diff-non-q-pos}
Conjecture~\ref{conj:merged-diff-q-pos} does not extend to
the case of $\gam=\iota^{d}(\tp{0,1,0})$ either.
Consider the multimonomial index $\tp{d,w,\alp,\gam}$
with $d=2$, $w=\tp{1}$, $\alp=\tp{2,1}$, and
$\gam=\iota^{d}(\tp{0,1,0})$.  Let
$k=\tp{0,1}$ and $\rho=\tp{1}$.  Then, 
$
\xi(d,w,\alp, \gam,0,0,1,\rho,q)
\not\llq 0$, because
\begin{align}
  \Delta(d,w,\alp, \gam,\tp{0},\tp{0},k,\rho,q)
  &=1,\\
  \Delta(d,w,\alp, \gam,\tp{1},\tp{1},k,\rho,q)
  &=3q^{6}+ 2q^{5}+ 2q^{4}+ 2q^{3}.
\end{align}
\end{example}

In Claim~\ref{c:mono-conv-mono-conv} of
Theorem~\ref{thm:mono-conv},
$\Delta(d,w,\alp,\bta, \gam,m,n,k,\rho,q) \in \Q[q^{\pm
  u_{\alp,\gam}}]$. Hence, we define the following notation
by $C_{q^{u_{\alp,\gam}}}$ in
Definition~\ref{def:poly-w-const} to state another
conjecture for $d\in \Zget$ and
$\gam=\iota^{d}\left(\tp{0,1,0}\right)$.

\begin{definition}
  Consider a multimonomial index
  $\tp{d,w,\alp,\bta,\gam}$. Let $l=1$,
  $m,n\in \Zl$, $k=\tp{0,1}$,
  and  $\lamo,\lamt,\lamr\in \Zgez$.
  Then,   we put
  \begin{align}
    \Delta_{C}(d,w,\alp,\bta, \gam,m,n,k,\rho,q)
    =
    C_{q^{u_{\alp,\gam}}}(\Delta(d,w,\alp,\bta,
    \gam,m,n,k,\rho,q))  \in\Q[q^{u_{\alp,\gam}}].
  \end{align}
  In particular, let
  \begin{align}
    \xi_{C}(d,w,\alp,\bta, \gam,\lamo,\lamt,\lamr,\rho,q)             
    &=\Delta_{C}(d,w,\alp,\bta, \gam,
      \tp{\lamo+\lamr},\tp{\lamo+\lamt+\lamr},k,\rho,q)
    \\ &-\Delta_{C}(d,w,\alp,\bta, \gam,\tp{\lamo},
         \tp{\lamo+\lamt},k,\rho,q).
  \end{align}
  Also, if $\bta=\iota^{d}(1)$, then let
  $\xi_{C}(d,w,\alp, \gam,\lamo,\lamt,\lamr,\rho,q) =
  \xi_{C}(d,w,\alp,\bta, \gam,\lamo,\lamt,\lamr,\rho,q)$.
\end{definition}
Moreover, we recall the following
{\it generalized Narayana  numbers}~\cite{Guy}.
\begin{definition}\label{def:narayana}
  Let $d\in \Zgeo$ and $h,\lam\in \Zgez$.  Then, the
  generalized Narayana number $N(\lam,d,h)$ satisfies
  \begin{align}
    N(\lam,d,h)=
    \frac{\lam+1}{d}{d \choose h}{d \choose h-1-\lam}.
  \end{align}
\end{definition}

For instance, $N(0,d,h)$ are Narayana numbers, which refine
Catalan numbers.  We also use the following numbers.

\begin{definition}
  Let $d\in \Zget$, $h_{1},h_{2}\in \Z$, and
  $\alp\in \Zgeo^{d}$.  Then, we define
  \begin{align}
    p(d,h_{1},h_{2},\alp)
    =(d- \len_{h_{2}}(\alp))h_{2}+h_{1} \len_{h_{2}}(\alp)
    \in \Z.
  \end{align}
\end{definition}
Notice that  $h_{2}>h_{1}\geq 1$ implies
$p(d,h_{1},h_{2},\alp)\geq 2$.  We now conjecture the following
positivity on $\xi_{C}$, periodicity on $\Delta_{C}$ and
$p(d,h_{1},h_{2},\alp)$, and equality on $\Delta_{C}$ and
$N(\lam,d,h)$.

\begin{conjecture}
  \label{conj:merged-diff-q-pos-period-nara}
  Let $d\in \Zget$, $w,\rho\in \Zgeo^{1}$, 
  $\alp\in \Zgeo^{d}$, and $k=\tp{0,1}$.
  For $i\in \oi(2)$,
  consider
  multimonomial indices $\tp{d,w,\alp,\gam_{i}}$
such that
  $\gamo=\iota^{d}(\tp{0,1,0})$ and
  $\gamt=\iota^{d}(\tp{\fraa,-\fraa,0})$.
  \begin{enumerate}
  \item \label{c:merged-diff-q-pos-period-nara-q-pos}
    If $\lamo,\lamt\in \Zgez$, then we have
    \begin{align}
      \xi_{C}(d,w,\alp, \gamo,\lamo,\lamt,1,\rho,q)             
      \llq 0.
    \end{align}
  \item \label{c:merged-diff-q-pos-period-nara-period}
    Let $h_{1},h_{2}\in \Zgeo$ with $h_{2}>h_{1}$.
    Suppose $\alp_{i}=h_{1}$ or $\alp_{i}=h_{2}$ for
    each $i\in\oi(d)$.  Then, whenever
    $\lamo,\lamt\in\Zgez$, we have the following
    periodicity
    of $\Delta_{C}$: first,
    if $\len_{h_{1}}(\alp)>0$ and $\len_{h_{2}}(\alp)>0$,
    then
    \begin{align}
      \xi_{C}(d,w,\alp,
      \gamt,\lamo,\lamt,p(d,h_{1},h_{2},\alp),\rho,0)
      =0;
    \end{align}
    second,     if $\len_{h_{1}}(\alp)=d$ or
    $\len_{h_{2}}(\alp)=d$, then
    \begin{align}
      \xi_{C}(d,w,\alp, \gamt,\lamo,\lamt,d,\rho,0)=0.
    \end{align}
    
  \item \label{c:merged-diff-q-pos-period-nara-nara} Assume
    $\alp=\iota^{d}(1)$ and $\lam,h\in \Zgez$ such that
    $ h\leq d-\lam$. Then, we have
    \begin{align}
      \Delta_{C}(d,w,\alp, \gamt,\tp{h},\tp{h+\lam},k,\rho,0)
      =N(\lam,d+1,h+1+\lam).
    \end{align}
  \end{enumerate}
\end{conjecture}
In particular, weighted $q$-multinomial coefficients would
extend the generalized Narayana numbers
$N(\lam,d+1,h+1+\lam)$ by
$\Delta_{C}(d,w,\alp,\gamt,\tp{h},\tp{h+\lam},k,\rho,0)>0$.

Example~\ref{ex:merged-diff-q-pos} supports
Claim~\ref{c:merged-diff-q-pos-period-nara-q-pos} in
Conjecture~\ref{conj:merged-diff-q-pos-period-nara}.
To see this, suppose
 multimonomial indices
$\tp{d,w,\alp,\gam}$ of $w=\tp{1}$, $d=2$,
$\alp=\iota^{d}(1)$, and $\gam=\iota^{d}(\tp{0,1,0})$.
Also, let
$k=\tp{0,1}$, $\rho=\tp{1}$, and
$\gam'=\iota^{d}(\tp{0,0,0})$.  Then,
we have
$\Delta_{C}(d,w,\alp,\gam,\tp{i},\tp{i},k,\rho,q)
  =
  \Delta_{C}(d,w,\alp,\gam',\tp{i},\tp{i},k,\rho,q)$, since $\alp$ is flat.
Example~\ref{ex:merged-diff-non-q-pos} is also consistent
with Claim~\ref{c:merged-diff-q-pos-period-nara-q-pos} in
Conjecture~\ref{conj:merged-diff-q-pos-period-nara}, as
$  \Delta_{C}(d,w,\alp, \gam,\tp{1},\tp{1},k,\rho,q)
  =3q^{3}+ 2q^{2}+ 2q+ 2$.

  The following supports
  Claims~\ref{c:merged-diff-q-pos-period-nara-period}
  and~\ref{c:merged-diff-q-pos-period-nara-nara} in
  Conjecture~\ref{conj:merged-diff-q-pos-period-nara}.
\begin{example}\label{ex:difference-2-1}
  Consider the multimonomial index $\tp{d,w,\alp,\gam}$ such
  that $d=3$, $w=\tp{1}$, $\alp=\iota^{d}(1)$, and
  $\gam=\iota^{d}(\tp{\fraa,-\fraa,0})$.  Let $k=\tp{0,1}$ and
  $\rho=\tp{1}$.  Then,
  \begin{dgroup*}
    \begin{dmath*}
      \Delta(d,w,\alp, \gam,\tp{0},\tp{0},k,\rho,q)=
      1,
    \end{dmath*}
    \begin{dmath*}
      \Delta(d,w,\alp, \gam,\tp{1},\tp{1},k,\rho,q)=
      3q + 6,
    \end{dmath*}
    \begin{dmath*}
      \Delta(d,w,\alp, \gam,\tp{2},\tp{2},k,\rho,q)
      = 6q^{4}+  18q^{3}+ 33q^{2}+ 18q+ 6,
    \end{dmath*}
    \begin{dmath*}
      \Delta(d,w,\alp, \gam,\tp{3},\tp{3},k,\rho,q)
      = 10q^{9}+
      35q^{8}+ 88q^{7}+ 155q^{6}+ 162q^{5}+ 144q^{4}+
      83q^{3}+ 43q^{2}+ 8q+ 1,
    \end{dmath*}
    \begin{dmath*}
      \Delta(d,w,\alp, \gam,\tp{4},\tp{4},k,\rho,q)=
      15q^{16}+ 57q^{15}+ 162q^{14}+ 357q^{13}
      + 642q^{12}+
      858q^{11}+ 1041q^{10}+ 1041q^{9}+ 912q^{8}+
      678q^{7}+ 447q^{6}+ 222q^{5}+ 96q^{4}+ 27q^{3}
      + 6q^{2},
    \end{dmath*}
    \begin{dmath*}
      \Delta(d,w,\alp, \gam,\tp{5},\tp{5},k,\rho,q)=
      21q^{25}+ 84q^{24}+ 255q^{23}+ 618q^{22}+ 1296q^{21}+ 2328q^{20}+ 3528q^{19}+ 4905q^{18}+ 6105q^{17}+ 6951q^{16}+ 7161q^{15}+ 6882q^{14}+ 5958q^{13}+ 4791q^{12}+ 3450q^{11}+ 2280q^{10}+ 1308q^{9}+ 690q^{8}+ 294q^{7}+ 111q^{6}+ 27q^{5}+ 6q^{4},
    \end{dmath*}
    \begin{dmath*}
      \Delta(d,w,\alp, \gam,\tp{6},\tp{6},k,\rho,q)=
      28q^{36}+ 116q^{35}+ 367q^{34}+ 938q^{33}+ 2114q^{32}+ 4229q^{31}+ 7638q^{30}+ 12248q^{29}+ 18374q^{28}+ 25457q^{27}+ 33106q^{26}+ 40206q^{25}+ 46214q^{24}+ 49983q^{23}+ 51342q^{22}+ 49690q^{21}+ 45663q^{20}+ 39534q^{19}+ 32462q^{18}+ 24969q^{17}+ 18121q^{16}+ 12227q^{15}+ 7745q^{14}+ 4463q^{13}+ 2382q^{12}+ 1124q^{11}+ 476q^{10}+ 164q^{9}+ 52q^{8}+ 8q^{7} +
      q^{6}.
    \end{dmath*}
  \end{dgroup*}
  We have $\len_{1}(\alp)=d$.  Thus,
  Claim~\ref{c:merged-diff-q-pos-period-nara-period}
  and~\ref{c:merged-diff-q-pos-period-nara-nara} in
  Conjecture~\ref{conj:merged-diff-q-pos-period-nara} agree
  with $\xi_{C}(d,w,\alp,\gam,i,0,3,\rho,0) =0$ for
  $i\in \oi(0,3)$, and with
  \begin{align}
    \Delta_{C}(d,w,\alp, \gam,\tp{0},\tp{0},k,\rho,0)
    &=1=              N(0,4,1),\\
    \Delta_{C}(d,w,\alp, \gam,\tp{1},\tp{1},k,\rho,0)
    &=6=              N(0,4,2),\\
    \Delta_{C}(d,w,\alp, \gam,\tp{2},\tp{2},k,\rho,0)
    &= 6=             N(0,4,3),\\ 
    \Delta_{C}(d,w,\alp, \gam,\tp{3},\tp{3},k,\rho,0)
    &=  1=            N(0,4,4).
  \end{align}
\end{example}

\begin{example}\label{ex:difference-2-2}
  Suppose the multimonomial index $\tp{d,w,\alp,\gam}$ such
  that $d=4$, $w=\tp{1}$, $\alp=\iota^{d}(1)$, and
  $\gam=\iota^{d}\tp{\tp{\fraa,-\fraa,0}}$.  Let
  $\rho=\tp{1}$, $k=\tp{0,1}$, and $\lam=2$.  Then, as in
  Claim~\ref{c:merged-diff-q-pos-period-nara-period} of
  Conjecture~\ref{conj:merged-diff-q-pos-period-nara}, we
  have the following periodicity:
  \begin{align}
    \Delta_{C}(d,w,\alp,
    \gam,\tp{0},\tp{0+\lam},k,\rho,0)
    &=
    \Delta_{C}(d,w,\alp,
    \gam,\tp{4},\tp{4+\lam},k,\rho,0)
    =                   6;\\
    \Delta_{C}(d,w,\alp,
    \gam,\tp{1},\tp{1+\lam},k,\rho,0)
    &=
    \Delta_{C}(d,w,\alp,
    \gam,\tp{5},\tp{5+\lam},k,\rho,0)
    =15;\\
    \Delta_{C}(d,w,\alp,
    \gam,\tp{2},\tp{2+\lam},k,\rho,0)
    &=
    \Delta_{C}(d,w,\alp,
    \gam,\tp{6},\tp{6+\lam},k,\rho,0)
    =6;\\
    \Delta_{C}(d,w,\alp,
    \gam,\tp{3},\tp{3+\lam},k,\rho,0)
    &=\Delta_{C}(d,w,\alp,
    \gam,\tp{7},\tp{7+\lam},k,\rho,0)
    =16.
  \end{align}
  In particular, the generalized Narayana numbers
  in Claim~\ref{c:merged-diff-q-pos-period-nara-nara}
  of Conjecture~\ref{conj:merged-diff-q-pos-period-nara}
  give the first three numbers in the above,
  since
  $0\leq h \leq d-\lam=2$
  and
  \begin{align}
    N(\lam,d+1,0+1+\lam)=
    N(2,5,3)&=6,\\
    N(\lam,d+1,1+1+\lam)=
    N(2,5,4)&=15,\\
    N(\lam,d+1,2+1+\lam)=
    N(2,5,5)&=6.
  \end{align}
  
  For another example, let $\lam=3$.  Then,
  we obtain the
  following periodicity:
  \begin{align}
    \Delta_{C}(d,w,\alp,
    \gam,\tp{0},\tp{0+\lam},k,\rho,0)
    &=
    \Delta_{C}(d,w,\alp,
    \gam,\tp{4},\tp{4+\lam},k,\rho,0)
    =  4;\\
    \Delta_{C}(d,w,\alp,
    \gam,\tp{1},\tp{1+\lam},k,\rho,0)
    &=
    \Delta_{C}(d,w,\alp,
    \gam,\tp{5},\tp{5+\lam},k,\rho,0)
    =               4;\\
    \Delta_{C}(d,w,\alp,
    \gam,\tp{2},\tp{2+\lam},k,\rho,0)
    &=
    \Delta_{C}(d,w,\alp,
    \gam,\tp{6},\tp{6+\lam},k,\rho,0)
    =    24;\\
    \Delta_{C}(d,w,\alp,
    \gam,\tp{3},\tp{3+\lam},k,\rho,0)
    &= 
    \Delta_{C}(d,w,\alp,
    \gam,\tp{7},\tp{7+\lam},k,\rho,0)
    =    24.
  \end{align}
  Then, the generalized Narayana numbers in
  Claim~\ref{c:merged-diff-q-pos-period-nara-nara} of
  Conjecture~\ref{conj:merged-diff-q-pos-period-nara}
  account the first two
  numbers, because $0\leq h\leq d-\lam=1$
  and
  \begin{align}
    N(\lam,d+1,0+1+\lam)=
    N(3,5,4)&=4,\\
    N(\lam,d+1,1+1+\lam)=
    N(3,5,5)&=4.
  \end{align}
\end{example}

We compute the following for some non-flat $\alp$ in
Claim~\ref{c:merged-diff-q-pos-period-nara-period} of
Conjecture~\ref{conj:merged-diff-q-pos-period-nara}.
\begin{example}\label{ex:difference-2-3}
  Consider the multimonomial index $\tp{d,w,\alp,\gam}$ such
  that $d=2$, $w=\tp{1}$, $\alp=\tp{1,2}$, and
  $\gam=\iota^{d}(\tp{\fraa,-\fraa,0})$.  Let $\rho=\tp{1}$
  and
  $k=\tp{0,1}$.  Also, let $h_{1}=1$ and $h_{2}=2$.
  Then,
  $\len_{2}(\alp)=1$ gives
  $p(d,h_{1},h_{2},\alp)=2\cdot 2
    -1\cdot 1=3$.
    Furthermore,
    Claim~\ref{c:merged-diff-q-pos-period-nara-period} of
  Conjecture~\ref{conj:merged-diff-q-pos-period-nara} is
  consistent with
  \begin{align}
    \Delta_{C}(d,w,\alp,
    \gam,\tp{0},\tp{0},k,\rho,0)
    &=
    \Delta_{C}(d,w,\alp,
      \gam,\tp{3},\tp{3},k,\rho,0)
    =               1,\\
    \Delta_{C}(d,w,\alp,
    \gam,\tp{1},\tp{1},k,\rho,0)
    &=
      \Delta_{C}(d,w,\alp,
      \gam,\tp{4},\tp{4},k,\rho,0)
    =               3,\\
    \Delta_{C}(d,w,\alp,
    \gam,\tp{2},\tp{2},k,\rho,0)
    &=              
    \Delta_{C}(d,w,\alp,
      \gam,\tp{5},\tp{5},k,\rho,0)
    =               1.
  \end{align}
\end{example}

\section{Primal monomial parcels}
\begin{definition}\label{def:primal}
  Suppose a parcel $\cF=\Lam(s,l,w,\scc,\fs,\rho,x,\fX)$.  We
  call $\cF$ primal, if $s=\tp{0,\infty}$, $l=1$,
  $w=\tp{1}$, $\cF_{\tp{0}}=1$, and $x$ is fully
  $\scc$-admissible.
\end{definition}

The primal parcels are important for the theory of the
merged-log-concavity.  First, they construct parcels of
arbitrary gates, widths, and weights by cut and shift
operators and Hadamard products.  Second, they allow
arbitrary base shifts by canonical mediators.  Third, they
consist of {\it primal monomial parcels}
$\Lam(s,l,w,\scc,\Psi_{s,\gam,q},\rho,x,\fX)$ such that
$\gam_{1,3}=0$.

\subsection{Primal monomial parcels
  and quantum dilogarithms}
\label{sec:q-dilog}

In this section, we claim that primal monomial parcels are
essentially quantum
dilogarithms discussed in~\cite{FadKas,FadVol,Kir,KonSoi,Rom,Schu,Zag}.
For the claim, we recall the notion of quantum dilogarithms not by the
merged-log-concavity but by the pentagon identities in
the following theorem, which is due
to~\cite{FadKas,FadVol,Schu} (see
also~\cite{Kir,Zag}).  For the completeness of this
manuscript, we provide a full proof of the theorem.

\begin{theorem}[\cite{FadKas,FadVol,Schu}]
  \label{thm:lin-quad-pent-series}
  Let $Q=\Q(\fX)$ and $q\in Q$. For $u_{0}=1$ and
  $u_{1}\in \Qti$, consider a formal power series
  $k_{q}(t)= \sum_{\lam\in\Zgez}u_{\lam}t^{\lam} \in Q[[t]]$.
  Let
  $\til{Q}=Q\langle z_{1},z_{2}\rangle/\langle
  z_{1}z_{2}=qz_{2}z_{1}\rangle$.  Also, let $p=1-q$ and
  $\mu=u_{1}p$.
  \begin{enumerate} 
  \item \label{c:lin-quad-pent-series-lin}
    The following statements are equivalent.
    \begin{enumerate}
    \item  \label{s:lin-quad-pent-series-lin-pent}
      We have the $(-\mu)$-pentagon
      identity
      $k_{q}(z_{1})k_{q}(z_{2})=
        k_{q}(z_{2})k_{q}(- \mu z_{2}z_{1})
        k_{q}(z_{1})$
        in $\til{Q}$.
      \item  \label{s:lin-quad-pent-series-lin-series}
      $k_{q}(t)= \sum_{\lam\in\Zgez}
        \frac{\mu^{\lam}}{(\lam)_{q}}\cdot t^{\lam}$.
      \end{enumerate}
    
  \item \label{c:lin-quad-pent-series-quad} The following
    statements are equivalent.
    \begin{enumerate}
    \item   \label{s:lin-quad-pent-series-quad-pent}
      We have the $\mu$-pentagon
      identity
      $k_{q}(z_{2})k_{q}(z_{1})=
        k_{q}(z_{1})
        k_{q}(\mu  z_{2}z_{1})
        k_{q}(z_{2})$
      in $\til{Q}$.
    \item \label{s:lin-quad-pent-series-quad-series}
      $k_{q}(t)= \sum_{\lam\in\Zgez}
        \frac{q^{\frac{\lam(\lam-1)}{2}}\cdot \mu^{\lam}}{(\lam)_{q}} \cdot t^{\lam}$.
    \end{enumerate}
  \end{enumerate}
\end{theorem}
\begin{proof}   
  Let us prove Claim~\ref{c:lin-quad-pent-series-lin}.
  First, assume the $(-u_{1})$-pentagon identity.  By 
  $z_{1}z_{2}=qz_{2}z_{1}$, both sides of
  the $(-\mu)$-pentagon identity read
  \begin{align}
    k_{q}(z_{1})k_{q}(z_{2})
    &=  
      \sum_{v_{1},v_{2}\in\Zgez}
      u_{v_{1}}u_{v_{2}}q^{v_{1}v_{2}}
      z_{2}^{v_{2}}z_{1}^{v_{1}},
      \label{eq:lin-quad-pent-series-lin-pent-lhs-weyl}\\
    k_{q}(z_{2})k_{q}(- \mu z_{2}z_{1})k_{q}(z_{1})
    &=\sum_{v'_{1},v'_{2},v'_{3}\in\Zgez}
    u_{v'_{1}}u_{v'_{2}}u_{v'_{3}}
    z_{2}^{v'_{1}}(- \mu z_{2}z_{1})^{v'_{2}}z_{1}^{v'_{3}}
    \\&=\sum_{v'_{1},v'_{2},v'_{3}\in\Zgez}
    u_{v'_{1}}u_{v'_{2}}u_{v'_{3}}
    (- \mu)^{v'_{2}}
    q^{\frac{v'_{2}(v'_{2}-1)}{2}}
    z_{2}^{v'_{1}+v'_{2}}z_{1}^{v'_{2}+v'_{3}}.
    \label{eq:lin-quad-pent-series-lin-pent-rhs-weyl}
  \end{align}
  To obtain the term $z_{2}^{v_{2}}z_{1}^{1}$ 
  in equation~\eqref{eq:lin-quad-pent-series-lin-pent-rhs-weyl},
  we only have  the following two cases:
  first, $v'_{1}=v_{2}$, $v'_{2}=0$, and
  $v'_{3}=1$; second,
  $v'_{1}=v_{2}-1$, $v'_{2}=1$, and $v'_{3}=0$.
  Both cases
  satisfy $q^{\frac{v'_{2}(v'_{2}-1)}{2}}=1$.  Hence,
  we have
  $u_{1}u_{v_{2}}q^{v_{2}}=u_{v_{2}}u_{1}-u_{v_{2}-1}
  u_{1}\mu $,
  comparing coefficients of $z_{2}^{v_{2}}z_{1}^{1}$ in
  equations~\eqref{eq:lin-quad-pent-series-lin-pent-lhs-weyl}
  and~\eqref{eq:lin-quad-pent-series-lin-pent-rhs-weyl}.
  Then, since 
  $\mu u_{v_{2}-1}=u_{v_{2}}(1-q^{v_{2}})$
  by $u_{1}\in \Qti$,
  $ u_{v_{2}-1}
  \frac{\mu}{1-q^{v_{2}}}=u_{v_{2}}$
  implies Statement~\ref{s:lin-quad-pent-series-lin-series}.
  
  Second, suppose
  Statement~\ref{s:lin-quad-pent-series-lin-series}.  Then,
  both sides of the $(-\mu)$-pentagon identity read
  \begin{align}
    k_{q}(z_{1})k_{q}(z_{2})
    &=
      \sum_{v_{1},v_{2}\in\Zgez}
      \frac{\mu^{v_{1}+v_{2}}}
      {(v_{1})_{q}(v_{2})_{q}}q^{v_{1}v_{2}}
      z_{2}^{v_{2}}z_{1}^{v_{1}}, 
      \label{eq:lin-quad-pent-series-lin-pent-lhs-series}\\
    k_{q}(z_{2})k_{q}(- \mu z_{2}z_{1})
    k_{q}(z_{1})
    &=\sum_{v'_{1},v'_{2},v'_{3}\in\Zgez}
      (-1)^{v'_{2}}  q^{\frac{v'_{2}(v'_{2}-1)}{2}}
      \frac{\mu^{v'_{1}+2v'_{2}+v'_{3}}}{(v'_{1})_{q}
      (v'_{2})_{q} (v'_{3})_{q}}
      z_{2}^{v'_{1}+v'_{2}}z_{1}^{v'_{2}+v'_{3}}.
      \label{eq:lin-quad-pent-series-lin-pent-rhs-series}
  \end{align}
  Now,
  $v'_{1}+v'_{2}=v_{2}$, $v'_{2}+v'_{3}=v_{1}$, and
  $v'_{2}=\lam$ imply $v'_{1}=v_{2}-\lam\geq 0$ and
  $v'_{3}=v_{1}-\lam\geq 0$.  Thus, comparing coefficients of
  $z_{2}^{v_{2}}z_{1}^{v_{1}}$ in
  equations~\eqref{eq:lin-quad-pent-series-lin-pent-lhs-series}
  and~\eqref{eq:lin-quad-pent-series-lin-pent-rhs-series},
  the $(-\mu)$-pentagon identity follows from
  \begin{align}
    \frac{q^{v_{1}v_{2}}}{(v_{1})_{q}(v_{2})_{q}}
    &=
      \sum_{\lam\in\oi(0,\min(v_{1},v_{2}))}
      \frac{    (-1)^{\lam}
      q^{\frac{\lam(\lam-1)}{2}}}{ (v_{1}-\lam)_{q}(\lam)_{q}
      (v_{2}-\lam)_{q}}.
  \end{align}  
  Multiplied by $(v_{1})_{q}(v_{2})_{q}$ on both sides,
  this is equivalent to
  \begin{align}
    \label{eq:lin-quad-pent-series-lin-pent-triple-binom}
    q^{v_{1}v_{2}}
    &=
      \sum_{\lam\in\oi(0,\min(v_{1},v_{2}))}
      {v_{2} \brack \lam}{v_{1} \brack \lam}(\lam)_{q}
      (-1)^{\lam}  q^{\frac{\lam(\lam-1)}{2}}.
  \end{align}
  
  Let us prove
  equation~\eqref{eq:lin-quad-pent-series-lin-pent-triple-binom}.
  Assume $v_{1}\geq v_{2}$ without loss of the
  generality. When $v_{2}=0$,
  equation~\eqref{eq:lin-quad-pent-series-lin-pent-triple-binom}
  clearly holds.  Furthermore, the induction on $v_{2}$
  gives
  \begin{align}
    \sum_{\lam\in\oi(0,v_{2}-1)}
    {v_{2} -1 \brack \lam}{v_{1} \brack \lam}
    (\lam)_{q}(-1)^{\lam}  q^{\frac{\lam(\lam-1)}{2}}
    &=q^{v_{1}(v_{2}-1)},\\
    \sum_{\lam\in\oi(0,v_{2}-1)}
    {v_{2} -1 \brack \lam}{v_{1}-1 \brack \lam}
    (\lam)_{q}(-1)^{\lam}  q^{\frac{\lam(\lam-1)}{2}}
    &=q^{(v_{1}-1)(v_{2}-1)}.
  \end{align}
  Thus, we obtain
  equation~\eqref{eq:lin-quad-pent-series-lin-pent-triple-binom},
  as the $q$-Pascal identity gives
\begin{dmath*}
    \sum_{\lam\in\oi(0,v_{2})}
    {v_{2} \brack \lam}{v_{1} \brack \lam}
    (\lam)_{q}(-1)^{\lam}  q^{\frac{\lam(\lam-1)}{2}}
    =   \sum_{\lam\in\oi(0,v_{2})}
    {v_{2} -1 \brack \lam}{v_{1} \brack \lam}
    (\lam)_{q}(-1)^{\lam}  q^{\frac{\lam(\lam-1)}{2}}
    +\sum_{\lam\in\oi(v_{2})}q^{v_{2}-\lam}
    {v_{2} -1 \brack \lam-1}{v_{1} \brack \lam}
    (\lam)_{q}(-1)^{\lam}  q^{\frac{\lam(\lam-1)}{2}}
    =
    \sum_{\lam\in\oi(0,v_{2}-1)}
    {v_{2} -1 \brack \lam}{v_{1} \brack \lam}
    (\lam)_{q}(-1)^{\lam}  q^{\frac{\lam(\lam-1)}{2}}
    +\sum_{\lam\in\oi(v_{2})}q^{(v_{2}-1)-(\lam-1)}
    \frac{(1-q^{v_{1}})}{(1-q^{\lam})}
    {v_{2} -1 \brack \lam-1}{v_{1} -1\brack \lam-1}
    (\lam)_{q}(-1)^{\lam}  q^{\frac{\lam(\lam-1)}{2}}
    =q^{v_{1}(v_{2}-1)}
    +(1-q^{v_{1}})q^{v_{2}-1}
    \sum_{\lam\in\oi(v_{2})}
    {v_{2} -1 \brack \lam-1}{v_{1}-1 \brack \lam-1}
    (\lam-1)_{q}(-1)^{\lam}  q^{\frac{(\lam-1)(\lam-2)}{2}}
    =q^{v_{1}(v_{2}-1)}
    +(q^{v_{1}}-1)q^{v_{2}-1}
    \sum_{\lam\in\oi(0,v_{2}-1)}
    {v_{2} -1 \brack \lam}{v_{1}-1 \brack \lam}
    (\lam)_{q}(-1)^{\lam}  q^{\frac{\lam(\lam-1)}{2}}
    =q^{v_{1}(v_{2}-1)}
    +(q^{v_{1}}-1)q^{v_{2}-1}q^{(v_{1}-1)(v_{2}-1)}
    =q^{v_{1}v_{2}}.
  \end{dmath*} 
  
  Let us prove Claim~\ref{c:lin-quad-pent-series-quad}.
  First, suppose
  Statement~\ref{s:lin-quad-pent-series-quad-pent}.  By
  $u_{0}=1$,
  $k_{q}(t)\in Q[[t]]^{\times}$.  Let
  $\omg_{q}(t)=k_{q}(t)^{-1}=1-u_{1}t+\dots$.
  Then, inverting both sides of
  the $\mu$-pentagon identity, we have
  \begin{align}
    \label{eq:lin-quad-pent-series-lin-pent-omg}
    \omg_{q}(z_{1})\omg_{q}(z_{2})=
    \omg_{q}(z_{2})
    \omg_{q}(-(-\mu)p z_{2}z_{1})
    \omg_{q}(z_{1}).
  \end{align}
  Hence, Claim~\ref{c:lin-quad-pent-series-lin} implies
  $\omg_{q}(t)=
  \sum_{\lam\in\Zgez}\frac{(-\mu)^{\lam}}{(\lam)_{q}}t^{\lam}$.
  By Euler's
  binomial identities, 
  Statement~\ref{s:lin-quad-pent-series-quad-series}
  follows from $k_{q}(t)=\omg_{q}(t)^{-1}
  =\sum_{\lam\in\Zgez}\frac{(-1)^{\lam}q^{\frac{\lam(\lam-1)}{2}}
    (-\mu)^{\lam}}{(\lam)_{q}}t^{\lam} =
  \sum_{\lam\in\Zgez}\frac{q^{\frac{\lam(\lam-1)}{2}}
    \mu^{\lam}}{(\lam)_{q}}t^{\lam}$.
  
  Second, suppose
  Statement~\ref{s:lin-quad-pent-series-quad-series}.
  Then,
  we go backward the discussion above, as
  Claim~\ref{c:lin-quad-pent-series-lin} is an
  equivalence.
  This gives
  equation~\eqref{eq:lin-quad-pent-series-lin-pent-omg} for
  $\omg_{q}(t)=k_{q}(t)^{-1}$.  Thus, $k_{q}(t)$ satisfies
  the $\mu$-pentagon identity.
\end{proof}

We now recall the 
dilogarithms below by the pentagon identities in
Theorem~\ref{thm:lin-quad-pent-series}.

\begin{definition}\label{def:q-dilog}
  Suppose $\ka\in \Q$ and $q\in \Q(\fX)$.
  If
  \begin{numcases}{k_{q}(t)=}
    \sum_{\lam\in\Zgez}
    \frac{q^{\ka \lam}}{(\lam)_{q}}\cdot t^{\lam}\in \Q(\fX)[[t]], \mor    \label{eq:first-kind}\\
    \sum_{\lam\in\Zgez}
    \frac{q^{\frac{\lam(\lam-1)}{2}+\ka\lam}}{(\lam)_{q}}
    \cdot t^{\lam}\in \Q(\fX)[[t]],
    \label{eq:second-kind}
  \end{numcases} 
  then $k_{q}(t)$ is called a quantum dilogarithm.
  For our convenience, we refer to $k_{q}(t)$ in
  equation~\eqref{eq:first-kind} or~\eqref{eq:second-kind}
  as a linear or quadratic quantum dilogarithm,
   respectively.
\end{definition}

If $k_{q}(t)$ is a
linear quantum dilogarithm, then
$(1-q)\cdot
\Log(k_{q}(t))|_{q=1}=\sum_{\lam\in\Zgeo}\frac{t^{\lam}}{\lam^{2}}=\Li_{2}(t)$,
which is the dilogarithm.
Similarly, if $k_{q}(t)$ is a
quadratic quantum
dilogarithm, then
$-(1-q)\cdot \Log(k_{q}(t))|_{q=1}=\Li_{2}(-t)$.

For instance,
$1+\sum_{\lam\in\Zgeo}
\frac{q^{\frac{\lam^{2}}{2}}}{(q^{\lam}-q^{\lam-1}) \dots
  (q^{\lam}-q)(q^{\lam}-1)} t^{\lam} =\sum_{\lam\in\Zgez}
\frac{q^{\frac{\lam}{2}}}{(\lam)_{q}} (-t)^{\lam} $
of~\cite[Section 6.4]{KonSoi} is a
linear quantum dilogarithm for the indeterminate $-t$.

Let us state the following lemma.
\begin{lemma}\label{lem:euclid}
  Let $\ka\in \Q$ and $q\in \Q(\fX)$.  Then, the following
  statements hold.
  \begin{enumerate}
  \item \label{c:euclid-linear} Let $\ka=\frac{u}{d}$ for
    coprime $d,u\in \Z$.  If $q^{\ka \lam}\in\Q(\fX)$ for each
    $\lam\in \Zgez$, then $q^{\frac{1}{d}}\in \Q(\fX)$.
  \item \label{c:euclid-quad} Let $\fraa+\ka=\frac{u}{d}$
    for coprime $d,u\in \Z$.  Suppose
    $q^{\frac{\lam^{2}}{2}+\ka \lam}\in\Q(\fX)$ for each
    $\lam\in \Zgez$.  Then, $q^{\frac{1}{d}}\in \Q(\fX)$.
    Also,
    $q^{\frac{\lam^{2}}{2}+\ka \lam}\in\Q(q^{\frac{1}{d}})$
    for each $\lam\in \Zgez$.
  \end{enumerate}
\end{lemma}
\begin{proof}
  Claim~\ref{c:euclid-linear} holds by
  the Euclid algorithm, since
  $q^{m}(q^{\ka})^{n}=q^{\frac{md+n u}{d}}\in \Q(\fX)$
  for
  each $m,n\in \Z$.
  
  Let us prove Claim~\ref{c:euclid-quad}.
  Because $q^{\fraa+\ka}\in \Q(\fX)$ for $\lam=1$,
  $q^{m}\left(q^{\frac{u}{d}}\right)^{n}\in \Q(\fX)$ for
  each $m,n\in
  \Z$. Hence, the former statement holds as above. Since
  $\frac{\lam^{2}}{2}+\ka \lam =\lam(\fraa+\ka)+{\lam
    \choose 2}$, the latter statement follows from
  ${\lam \choose 2}\in \Z$.
\end{proof}
Hence, we state the following by the
linear and quadratic
primal monomial parcels and 
quantum dilogarithms.
\begin{theorem}
  \label{thm:inv-lim}
  If there is a primal monomial parcel
  $\cF=\Lam(s,l,w,\scc,\Psi_{s,\gam,q},\rho,x,\fX)$, then the
  generating function $\cZ_{w,\gam,q}(t)$ of $\cF$ is a
  quantum dilogarithm.  Conversely, each quantum
  dilogarithm
  coincides with the generating function of a primal
  monomial parcel.
\end{theorem}
\begin{proof}
  Let us prove the former part.  The monomial conditions
  of
  $\mu=\tp{l,w,\gam}$ imply $\gam_{1,1}=0$ or
  $\gam_{1,1}=\fraa$, If $\gam_{1,1}=0$, then
  $\cZ_{w,\gam,q}(t) =\sum_{\lam\in\Zgez}
  \frac{q^{\gam_{1,2}\lam}}{(\lam)_{q}} \cdot t^{\lam}\in
  \Q(\fX)[[t]]$.  Also, if $\gam_{1,1}=\fraa$, then
  $\cZ_{w,\gam,q}(t) =\sum_{\lam\in\Zgez}
  \frac{q^{\frac{\lam^{2}}{2}+\gam_{1,2}\lam}}{(\lam)_{q}} \cdot
  t^{\lam}\in \Q(\fX)[[t]]$.  Moreover, since
  $\gam_{1,2}\in \Q$ for the monomial index $\mu$,
  $\cZ_{w,\gam,q}(t)$ is a quantum dilogarithm.
  
  Let us prove the latter part. Consider
   a
  linear quantum
  dilogarithm
  $k_{q}(t)=\sum_{\lam\in\Zgez} \frac{q^{\ka\lam}}{(\lam)_{q}}
  \cdot t^{\lam}\in \Q(\fXo)[[t]]$.  Then,
  Claim~\ref{c:euclid-linear} of Lemma~\ref{lem:euclid}
  gives $d\in \Zgeo$ such that
  $\frac{q^{\ka\lam}}{(\lam)_{q}}
  \in \Q(q^{\frac{1}{d}}) $. 
  Let $\fXt\subset \Q(\fXo)$
  be a set
  of  indeterminates
  with
  $q^{\frac{1}{d}}\in \fXt$.   Then, there is a primal
  monomial parcel
  $\cF=\Lam(s,l,w,\scc,\Psi_{s,\gam,q},\rho,x,\fXt)$ such that
  $\gam=\tp{\tp{0,\ka,0}}$ and
  $\cZ_{w,\gam,q}(t)=k_{q}(t)$.  A similar discussion holds
  for quadratic quantum dilogarithms by
  Claim~\ref{c:euclid-quad} of Lemma~\ref{lem:euclid}.
\end{proof}

\subsection{Phase transitions and the golden ratio}
\label{sec:primal-monom}

We now study the almost strictly unimodal
sequences and phase transitions.
This gives the golden
ratio as a critical point.
We  use the
following notion.
\begin{definition}
  \label{def:van-prob}
  Suppose the merged pair $\zeta_{P}=\tp{P,\cF}$ of a
  parcel
  $\cF=\Lam(s,l,w,\scc,\fs,\phi,\rho,x,\fX)$ and fitting path
  $P=\{\tp{s,l,m_{i},n_{i},k_{i}}\}_{i\ldZ \theta}$.
  \begin{enumerate}
 \item   We call
   $\zeta_{P}$ vanishing,
   if
   $\lim_{i\to \infty}u(\zeta,r)_{i}=0$ for each $r\in \AfX$.
 \item\label{def:prob}
   We call
   $\zeta_{P}$ probabilistic, if
   $\sum_{i\ldZ e(\theta) }u(\zeta,r)_{i}<\infty$ for each
   $r\in \AfX$. 
 \end{enumerate}
 If $\zeta_{P}$ is vanishing (resp. probabilistic)
 for each fitting
 path $P$, then we call $\cF$ vanishing (resp.
 probabilistic).
 \end{definition}

 Unimodal sequences are important in probability.  On
 the
 terminology in Item~\ref{def:prob} in
 Definition~\ref{def:van-prob}, assume a probabilistic
 $\zeta_{P}$ and $r\in \AfX$.  Then, $u(\zeta,r)_{i}>0$ for each
 $i\ldZ e(\theta)$ by Lemma~\ref{lem:total-pos-neg}.  Then,
 Theorem~\ref{thm:path-parcel-seq} gives the almost
 strictly
 unimodal sequence $u(\zeta,r)$ with the {\it probability
   generating function}
 $(\sum_{i\ldZ e(\theta) }u(\zeta,r)_{i})^{-1} \sum_{i\ldZ e(\theta)
 }u(\zeta,r)_{i} z^{i}$ for an indeterminate $z$.

We have the following equivalence.

\begin{proposition}
  \label{prop:van-prob-tgam}
  Suppose a primal monomial parcel
  $\cF=\Lam(s,l,w,\scc,\Psi_{s,\gam,q},\rho,x,\fX)$.  Let
  $r\in \AfX$.  Then, the following statements are
  equivalent.
  \begin{enumerate}
  \item\label{s:van-prob-tgam-van}
    $\cF$ is vanishing.
  \item\label{s:van-prob-tgam-tgam}
    There exists some
    $\lam\in \Zgeo$ such that
    $t_{\gam}(\tp{\mu})>0$ 
    for each $\mu\in \Z_{\geq\lam}$.
  \item\label{s:van-prob-tgam-geom}
    There exist
    real numbers $N_{\gam}(r)\geq 1$ and
    $0<S_{\gam}(r)<1$ such that 
    $\cF_{\tp{\mu}}(r)\leq  N_{\gam}(r)
    S_{\gam}(r)^{\mu}$ for each $\mu\in \Zgez$.
  \item $\cF$ is probabilistic.
    \label{s:van-prob-tgam-prob}
    \end{enumerate}
\end{proposition}
\begin{proof} 
  Let us prove
  Statement~\ref{s:van-prob-tgam-tgam}
  from Statement~\ref{s:van-prob-tgam-van}.
  Since $w=\tp{1}$, the monomial conditions of
  $\ka=\tp{l,w,\gam}$ imply $\gam_{1,1}=\fraa$ or 0.
  Suppose $\gam_{1,1}=0$.  If $\gam_{1,2}\leq 0$, then
  $\cF_{\tp{0}}(r) =1< \cF_{\tp{1}}(r) =
  \frac{q(r)^{\gam_{1,2}}}{1-q(r)} < \cF_{\tp{2}}(r)
  =\frac{q(r)^{2\gam_{1,2}}}{(1-q(r))(1-q(r)^{2})}<\dots $,
  since
  $0<q(r)<1$ in Claim~\ref{c:adm-bounds} of
  Lemma~\ref{lem:adm-bounds-poring-inclusions}.
  However, this is
  against Statement~\ref{s:van-prob-tgam-van},
  which makes
$\Omg(\cF)$ vanishing.
Hence,
Statement~\ref{s:van-prob-tgam-tgam} follows from
  $\gam_{1,2}>0$.  If $\gam_{1,1}=\fraa$, then
  Statement~\ref{s:van-prob-tgam-tgam} holds for
  any
  $\gam_{1,2}\in \Q$.
  
  Let us prove
  Statement~\ref{s:van-prob-tgam-geom} from
  Statement~\ref{s:van-prob-tgam-tgam}.
  Since 
  $\gam_{1,1}=0$ and $ \gam_{1,2}>0$, or $\gam_{1,1}=\fraa$,
  we have $U_{\gam}\in \Qgo$ and
  $V_{\gam}(r)\in \Zgez$ such that
  each integer $\mu\geq V_{\gam}(r)+1$ satisfies
  \begin{align}
    \frac{q(r)^{U_{\gam}}}{1-q(r)^{{V_{\gam}(r)+1}}}
    &< 1,
      \label{ineq:tgam-geom-bounds-ups-upper}\\
    t_{\gam}(\tp{\mu})
    &\geq U_{\gam}  \mu.
      \label{ineq:tgam-geom-bounds-tgam-ugam}
  \end{align}
  Let
  $\ups(r)
    =  \frac{1}{1-q(r)}$,
    $S_{\gam}(r)
     =\frac{q(r)^{U_{\gam}}}
     {1-q(r)^{V_{\gam}(r)+1}}$, and
      $\chi_{\gam}(r)
      =\frac{\ups(r)}{S_{\gam}(r)}  > 1$.
      In particular,   $\mu\in \Z_{\geq 1}$ gives
      \begin{align}
        \ups(r)
        &
          >\frac{q(r)^{U_{\gam}}}{1-q(r)^{\mu}},
          \label{ineq:tgam-geom-bounds-ups-lower}\\
        S_{\gam}(r)
        &\geq
          \frac{q(r)^{U_{\gam}}}{1-q(r)^{V_{\gam}(r)+\mu}}.
          \label{ineq:tgam-geom-bounds-kagam-lower}
      \end{align}
      Also, let
      $\omg_{\gam}(r)=\max(q(r)^{t_{\gam}(\tp{0})}=1,
      \cdots, q(r)^{t_{\gam}(\tp{V_{\gam}(r)})})$ and
        $N_{\gam}(r)
        =\left(\omg_{\gam}(r)
          \chi_{\gam}(r)\right)^{V_{\gam}(r)} \geq 1$.
        
        First, assume $\mu\geq V_{\gam}(r)+1$.  Then,
        inequality~\eqref{ineq:tgam-geom-bounds-tgam-ugam}
gives
  $\cF_{\tp{\mu}}(r)
  = \frac{q(r)^{t_{\gam}(\tp{\mu})}}
  {(\mu)_{q}|_{q=q(r)}}
  \leq \prod_{i\in \oi(\mu)}\frac{q(r)^{U_{\gam}}}{1-q(r)^{i}}$.
  Hence,
  $\cF_{\tp{\mu}}(r)
    \leq \ups(r)^{V_{\gam}(r)}
    S_{\gam}(r)^{\mu-V_{\gam}(r)}
    =\chi_{\gam}(r)^{V_{\gam}(r)}S_{\gam}(r)^{\mu}
    \leq       N_{\gam}(r)
    S_{\gam}(r)^{\mu}$
    by inequalities~\eqref{ineq:tgam-geom-bounds-ups-lower}
and~\eqref{ineq:tgam-geom-bounds-kagam-lower}.
    Second, suppose $0\leq \mu\leq V_{\gam}(r)$.  Then,
    $\cF_{\tp{\mu}}(r)
    \leq
    \left(\frac{\omg_{\gam}(r)}{1-q(r)}
    \right)^{\mu}\leq 
    \left(\omg_{\gam}(r)\ups(r)\right)^{V_{\gam}(r)}$.
    Thus,
inequality~\eqref{ineq:tgam-geom-bounds-ups-upper}
    gives
    $\cF_{\tp{\mu}}(r)
    \leq   (\omg_{\gam}(r)  \ups(r))^{V_{\gam}(r)}
    \cdot
    \frac{S_{\gam}(r)^{\mu}}{S_{\gam}(r)^{V_{\gam}(r)}}
    =N_{\gam}(r)      S_{\gam}(r)^{\mu}$.

    Statement~\ref{s:van-prob-tgam-prob} follows
    from
    Statement~\ref{s:van-prob-tgam-geom},
    as
    $\sum_{m\ldZl s}\cF_{m}(s)
  \leq N_{\gam}(r)
  \sum_{i \in \Zgez} S_{\gam}(r)^{i}<\infty$.
  Also,
  Statement~\ref{s:van-prob-tgam-prob}
  implies
  Statement~\ref{s:van-prob-tgam-van}
  by Lemma~\ref{lem:total-pos-neg}. 
\end{proof}

We  discuss phase transitions
first by asymptotic critical points.

\begin{lemma}\label{lem:asymp-hill-crit-pts}
  Consider the merged pair $\zeta=\tp{P,\cF}$ of an
  infinite-length $P$  and
   a primal monomial parcel
   $\cF=\Lam(s,l,w,\scc,\Psi_{s,\gam,q},\rho,x,\fX)$.
   When $r\in \AfX$,
    we have the following.
   \begin{enumerate}
   \item
     \label{c:asymp-hill-crit-pts-hill}
     If
     $
     \gam_{1,1}=\gam_{1,2}=0$, then
     $u(\zeta,r)$ is asymptotically
     hill-shape.
   \item
     \label{c:asymp-hill-crit-pts-non-exist}
     If not, then
     $\zeta$ has no asymptotic
  critical points.
   \end{enumerate}
\end{lemma}
\begin{proof}
  Let $P=\{\tp{s,l,m_{i},n_{i},k_{i}}\}_{i\ldZ \theta}$ of an
  infinite $\theta$.  Let $\theta_{1}=1$ for simplicity.
  Then, Lemma~\ref{lem:finest} provides
  $\ka\in \Zgeo$ and $h\in \Z$ such that
  $m_{i}=n_{i}=\tp{i\ka+h}$ for 
  $i\ldZ \theta$.  
  Claim~\ref{c:asymp-hill-crit-pts-hill} now holds by
  Proposition~\ref{prop:path-parcel-asymp-hill-shape},
  since
  Claim~\ref{c:adm-bounds} of
  Lemma~\ref{lem:adm-bounds-poring-inclusions}
  gives
  $\lim_{j\to \infty}\frac{\cF_{m_{j+1}}(r)}{\cF_{m_{j}}(r)}
  =\lim_{j\to \infty}
  \frac{(j\ka+h)_{q}|_{q=q(r)}}
  {((j+1)\ka+h)_{q}|_{q=q(r)}}
    =\lim_{j\to\infty}
    \prod_{i\in\oi(\lam)}
    \frac{1}{1-q(r)^{j\ka+h+i}}
    =1$.

    Let us prove
    Claim~\ref{c:asymp-hill-crit-pts-non-exist}.  First,
    suppose $\lam\in \Zgeo$ such that each $\mu\geq \lam$
    satisfies
  \begin{align}
    t_{\gam}(\tp{\mu})>0.
    \label{ineq:non-exist-asympt-crit-pts-tgam-pos}
  \end{align}
  Then, Lemma~\ref{lem:fitting-path-mn-positivity} and
  Proposition~\ref{prop:van-prob-tgam} yield
  $\lim_{j\to \infty}u(\zeta,r)_{j}=0$.
  Thus,
  Claim~\ref{c:asymp-hill-crit-pts-non-exist} follows,
  since
  $u(\zeta,r)>0$ is strictly decreasing or hill-shape by
  Claim~\ref{c:path-parcel-seq-unim} of
  Theorem~\ref{thm:path-parcel-seq}.

  Second, assume that
  inequality~\eqref{ineq:non-exist-asympt-crit-pts-tgam-pos}
  does not hold. This gives $\gam_{1,1}=0$ and
  $\gam_{1,2}\leq 0$, since $\gam_{1,1}=\fraa$ or
  $0$ by the monomial conditions.  Thus,
  the  assumption of
  Claim~\ref{c:asymp-hill-crit-pts-non-exist} implies
  $\gam_{1,1}=0$ and $\gam_{1,2}<0$.
  Then, $q(r)^{\gam_{1,2}}> 1$ by
  Claim~\ref{c:adm-bounds} of
  Lemma~\ref{lem:adm-bounds-poring-inclusions}.  
  Claim~\ref{c:asymp-hill-crit-pts-non-exist} now holds,
  since
  $\lim_{j\to \infty}\frac{\cF_{m_{j+1}}}{\cF_{m_{j}}}
  =\lim_{j\to \infty}
  \frac{q(r)^{ ((j+1)\ka+h)\cdot \gam_{1,2}}}
  {((j+1)\ka+h)_{q}|_{q=q(r)}}
  \cdot \frac{(j\ka+h)_{q}|_{q=q(r)}}
  {q(r)^{ (j\ka+h)    \cdot \gam_{1,2}}}
  =\lim_{j\to\infty}\frac{q(r)^{\ka\gam_{1,2}}}{
    \prod_{i\in\oi(\ka)}(1-q(r)^{j\ka+h+i})}
  >1$.
\end{proof}
By the following
notation, we obtain front and rear phase transitions.
\begin{definition}
  Let $\theta=\tp{1,\infty}$ and 
  $\chi\geq 1$ be a gate.
  Suppose a merged pair $\Omg_{\lam}(\cF)$.
  Then, we write $\Omg_{\lam}^{\chi}(\cF)$
  for the restricted merged
  pair $\tp{r_{\theta,\chi}(\pi(\lam)),\cF}$.
    \end{definition}

\begin{lemma}\label{lem:front-rear-pt}
  Let $\cF=\Lam(s,l,w,\scc,\Psi_{s,\gam,q}, \rho,x,\fX)$
  be a primal monomial parcel.
  Suppose some $\lam\in \Zgeo$ such that
  each $\mu\in \Z_{\geq \lam}$ satisfies
  \begin{align}
    t_{\gam}(\tp{\mu})>0.
    \label{ineq:front-rear-pt-tgam}
  \end{align}
  Then,
  we have the following.
  \begin{enumerate}
  \item \label{c:front-rear-pt-front} The merged pair
    $\zeta_{\lam}=\Omg_{\lam}(\cF)$ has a front phase
    transition.
  \item \label{c:front-rear-pt-rear} The merged pair
    $\zeta_{\lam,\chi}=
    \Omg_{\lam}^{\chi}(\cF)$
    has a rear phase transition for each finite gate
    $\chi$.
  \end{enumerate}
\end{lemma}
\begin{proof}
  For $\theta=\tp{1,\infty}$, suppose the fitting path
  $\pi(\lam)=\{\tp{s,l,m_{i},n_{i},k_{i}}\}_{i\ldZ \theta}$ with
  $m_{0}=n_{0}=\tp{0}$.  Let us prove
  Claim~\ref{c:front-rear-pt-front}.  We use
  Lemma~\ref{lem:BoTo-crit-points} as follows.  There is
  $d\in \Zgeo$ such that $T=q^{\frac{1}{d}}\in \fX$ for the
  fully admissible $q$.  Then, let
  $\Bo_{1}(T)=\To_{1}(T)=1$, $\Bo_{2}(T)=(\lam)_{q}$, and
  $\To_{2}(T)=q^{t_{\gam}(\tp{\lam})}$.  In particular,
  $\Qu(\To)(T)= q^{t_{\gam}(\tp{\lam})}$ and
  $\Qu(\Bo)(T)=(\lam)_{q}$.  Then, $\Qu(\To)(0)=0$ by
  inequality~\eqref{ineq:front-rear-pt-tgam} and
  $\Qu(\To)(1)=1$. Also, we have $\Qu(\Bo)(0)=1$ and
  $\Qu(\Bo)(1)=0$.  Thus, Lemma~\ref{lem:BoTo-crit-points}
  gives Claim~\ref{c:front-rear-pt-front} by
  $\cF_{\tp{0}}(q)=\frac{\To_{1}(T)}{\Bo_{1}(T)}$ and
  $\cF_{\tp{\lam}}(q)=\frac{\To_{2}(T)}{\Bo_{2}(T)}$.
  
  Set $\ka\in \Zgeo$ and $\chi=\tp{1,\ka}$.  Let us prove
  Claim~\ref{c:front-rear-pt-rear} by
  Lemma~\ref{lem:BoTo-crit-points}.  Hence, consider
  $\Bo_{1}(T)=(\ka\lam)_{q}$,
  $\Bo_{2}(T)=((\ka+1)\lam)_{q}$,
  $\To_{1}(T)=q^{t_{\gam}(\tp{\ka\lam})}$, and
  $\To_{2}(T)=q^{t_{\gam}(\tp{(\ka+1)\lam})}$.  Then, we
  have $\Qu(\Bo)(0)=1$, $\Qu(\Bo)(1)=0$, and
  $\Qu(\To)(1)=1$. Also,
  $\Qu(\To)(0)= 0$, because
  inequality~\eqref{ineq:front-rear-pt-tgam}
  and $\gam_{1,1}\geq 0$ imply
  \begin{dmath*}
    t_{\gam}(\tp{(\ka+1)\lam})- t_{\gam}(\tp{\ka\lam})
    =(((\ka+1)\lam)^{2}-(\ka\lam)^{2})
    \gam_{1,1} +
    ((\ka+1)\lam-\ka\lam) \gam_{1,2}
    =(2\ka\lam+\lam^2)\gam_{1,1} +\lam\gam_{1,2}
    =2\ka\gam_{1,1} +t_{\gam}(\tp{\lam}) >0.
  \end{dmath*}
\end{proof}

We have the following positive values of $t_{\gam}$ by
semi-phase transitions.
\begin{lemma}\label{lem:front-rear-tgam}
  Consider a primal monomial parcel
  $\cF=\Lam(s,l,w,\scc, \Psi_{s,\gam,q},\rho,x,\fX)$.
  Let
    $\lam\in \Zgeo$.
  \begin{enumerate}
  \item \label{c:front-rear-tgam-front} If
    $\Omg_{\lam}(\cF)$ has a front semi-phase transition,
    then $t_{\gam}(\tp{\mu}) > 0$ for each $\mu\in\Z_{\geq \lam}$.
  \item\label{c:front-rear-tgam-rear}
    For a finite gate $\chi$,
    if
  $\Omg_{\lam}^{\chi}(\cF)$ has a
  rear semi-phase transition, then there is $\lam'\in \Zgeo$
  such that $t_{\gam}(\tp{\mu}) > 0$ for each
  $\mu\in \Z_{\geq \lam'}$.
  \end{enumerate}
\end{lemma}
\begin{proof}
  For $\theta=\tp{1,\infty}$, 
  let
  $\pi(\lam)=\{\tp{s,l,m_{i},n_{i},k_{i}}\}_{i\ldZ \theta}$ with
  $m_{0}=n_{0}=\tp{0}$.
  Also, let $\ka=\tp{l,w,\gam}$.
  First, let us prove Claim~\ref{c:front-rear-tgam-front}.
  We have $0<(\lam)_{q}|_{q=q(r)}<1$ by
  Claim~\ref{c:adm-bounds} of
  Lemma~\ref{lem:adm-bounds-poring-inclusions}.
  Also, the
  front semi-phase transition gives $r\in \AfX$ such that
  $\cF_{m_{0}}(r)
    =1
    =\frac{q(r)^{t_{\gam}(\tp{\lam})}}
    {(\lam)_{q}|_{q=q(r)}}
    =
    \cF_{m_{1}}(r)$.
    Hence,
    $0<q(r)^{t_{\gam}(\tp{\lam})}<1$,
    which implies
  $t_{\gam}(\tp{\lam})=    \gam_{1,1}\lam^2+
  \gam_{1,2}\lam>0$.
  Claim~\ref{c:front-rear-tgam-front}
   now holds, since
   $\gam_{1,1}\geq 0$
  by the monomial conditions of $\ka$.
  
  Second, let us prove
  Claim~\ref{c:front-rear-tgam-rear}.
  The rear semi-phase
  transition demands $r\in \AfX$ such that
  $\cF_{m_{\chi_{2}}}(r)
    = \frac{q(r)^{t_{\gam}(\tp{\chi_{2}\lam})}}
    {(\chi_{2}\lam)_{q}|_{q=q(r)}}
    =
    \frac{q(r)^{t_{\gam}(\tp{(\chi_{2}+1)\lam})}}
    {((\chi_{2}+1)\lam)_{q}|_{q=q(r)}}
    = \cF_{m_{\chi_{2}+1}}(r)$.
    Then, 
      $0<q(r)^{t_{\gam}(\tp{(\chi_{2}+1)\lam})-
    t_{\gam}(\tp{\chi_{2}\lam})}<1$ by
  $0<
  \frac{((\chi_{2}+1)\lam)_{q}|_{q=q(r)}}
  {(\chi_{2}\lam)_{q}|_{q=q(r)}} <1$.
    In particular,
  $t_{\gam}(\tp{(\chi_{2}+1)\lam})-
    t_{\gam}(\tp{\chi_{2}\lam})
    =(2\chi_{2}\lam+\lam^2)\gam_{1,1}
    +\lam\gam_{1,2}
    >0$.
    Thus, the monomial conditions of $\ka$ imply
  $\gam_{1,1}=\fraa$, or $\gam_{1,1}=0$ and $\gam_{1,2}>0$.
\end{proof}

Then, we obtain the following on the phase transitions of
primal monomial parcels.  This extends
Proposition~\ref{prop:monom-front-pt} by vanishing and probabilistic parcels.

\begin{theorem}\label{thm:prim-monom-pt}
  Suppose a primal monomial parcel
  $\cF=\Lam(s,l,w,\scc,\Psi_{s,\gam,q},\rho,x,\fX)$.  Then, we have the following.
  \begin{enumerate}
  \item \label{c:prim-monom-pt-no-asymp-pt} Each
    merged pair
    $\tp{P,\cF}$ has no asymptotic semi-phase transitions.
  \item \label{c:prim-monom-pt-vanishing-pcl-pt}
    The following statements are equivalent.
    \begin{enumerate}
    \item $\cF$ is vanishing.
      \label{s:prim-monom-pt-vanishing-pcl-pt-vanishing-pcl}
    \item $\cF$ is probabilistic. 
      \label{s:prim-monom-pt-vanishing-pcl-pt-prob}
    \item \label{s:prim-monom-pt-vanishing-pcl-pt-tgam}
      There exists some $\lam\in \Zgeo$ such that
      $t_{\gam}(\tp{\mu}) > 0$ for each $\mu\in\Z_{\geq \lam}$.
    \item For some $\lam\in \Zgeo$, $\Omg_{\lam}(\cF)$ has a
      front phase transition.
      \label{s:prim-monom-pt-vanishing-pcl-pt-front-pt}  
    \item For some $\lam\in \Zgeo$, $\Omg_{\lam}(\cF)$ has a
      front semi-phase transition.
      \label{s:prim-monom-pt-vanishing-pcl-pt-front-semi-pt}
    \item \label{s:prim-monom-pt-vanishing-pcl-pt-rear-pt} For
      some $\lam\in \Zgeo$ and finite gate $\chi$,
      $\Omg_{\lam}^{\chi}(\cF)$ has a rear
      phase transition.
    \item \label{s:prim-monom-pt-vanishing-pcl-pt-rear-semi-pt}
      For some $\lam\in \Zgeo$ and finite gate $\chi$,
      $\Omg_{\lam}^{\chi}(\cF)$ has a rear
      semi-phase transition.
    \end{enumerate}
  \item \label{c:prim-monom-pt-tgam-pt} If one of
    Statements~\ref{s:prim-monom-pt-vanishing-pcl-pt-tgam}
    --~\ref{s:prim-monom-pt-vanishing-pcl-pt-front-semi-pt}
    holds by some $\lam\in \Zgeo$, then each of
    Statements~\ref{s:prim-monom-pt-vanishing-pcl-pt-tgam}
    --~\ref{s:prim-monom-pt-vanishing-pcl-pt-rear-semi-pt} holds
    by the same $\lam\in \Zgeo$.
  \end{enumerate}
\end{theorem}
\begin{proof}
  Claim~\ref{c:prim-monom-pt-no-asymp-pt} holds by
  Lemma~\ref{lem:asymp-hill-crit-pts}.
    
  Let us prove Claim~\ref{c:prim-monom-pt-vanishing-pcl-pt}.
  First,
  Statements~\ref{s:prim-monom-pt-vanishing-pcl-pt-vanishing-pcl},
  \ref{s:prim-monom-pt-vanishing-pcl-pt-prob},
  and~\ref{s:prim-monom-pt-vanishing-pcl-pt-tgam} are
  equivalent by Proposition~\ref{prop:van-prob-tgam}.
  Second,
  Statements~\ref{s:prim-monom-pt-vanishing-pcl-pt-tgam},
  \ref{s:prim-monom-pt-vanishing-pcl-pt-front-pt},
  and~\ref{s:prim-monom-pt-vanishing-pcl-pt-front-semi-pt}
  are equivalent by
  Claim~\ref{c:front-rear-pt-front} of
  Lemma~\ref{lem:front-rear-pt} and
  Claim~\ref{c:front-rear-tgam-front} of
  Lemma~\ref{lem:front-rear-tgam}.  Third,
  the equivalence
  among
  Statements~\ref{s:prim-monom-pt-vanishing-pcl-pt-tgam},
  \ref{s:prim-monom-pt-vanishing-pcl-pt-rear-pt},
  and~\ref{s:prim-monom-pt-vanishing-pcl-pt-rear-semi-pt}
  follows from Claim~\ref{c:front-rear-pt-rear} of
  Lemma~\ref{lem:front-rear-pt} and
  Claim~\ref{c:front-rear-tgam-rear} of
  Lemma~\ref{lem:front-rear-tgam}.

  Let us prove Claim~\ref{c:prim-monom-pt-tgam-pt}.
  Claim~\ref{c:front-rear-pt-front} of
  Lemma~\ref{lem:front-rear-pt} and
  Claim~\ref{c:front-rear-tgam-front}
  of Lemma~\ref{lem:front-rear-tgam}
  imply that
  Statements~\ref{s:prim-monom-pt-vanishing-pcl-pt-tgam},
  \ref{s:prim-monom-pt-vanishing-pcl-pt-front-pt},
  and~\ref{s:prim-monom-pt-vanishing-pcl-pt-front-semi-pt}
  hold by the same $\lam$, which give
  Statements~\ref{s:prim-monom-pt-vanishing-pcl-pt-rear-pt}
  and~\ref{s:prim-monom-pt-vanishing-pcl-pt-rear-semi-pt}
  by Claim~\ref{c:front-rear-pt-rear} of
  Lemma~\ref{lem:front-rear-pt}.
\end{proof}

Furthermore, Theorem~\ref{thm:prim-monom-pt} gives the
following polynomials with positive integer coefficients by
the finest fitting path $\pi(1)$.

\begin{proposition}\label{prop:ideal-front-pt-gam}
  Let $\cF=\Lam(s,l,w,\scc,\Psi_{s,\gam,q},\rho,x,\fX)$
  be a primal monomial parcel.
  Then,
  Statements~\ref{s:ideal-front-pt-gam-front-pt}
  and~\ref{s:ideal-front-pt-gam-gam} below are equivalent.
  \begin{enumerate}
  \item \label{s:ideal-front-pt-gam-front-pt} The merged
    pair $\xi=\Omg(\cF)$ is ideal with a front phase
    transition.
  \item \label{s:ideal-front-pt-gam-gam} We have
    some $\lam\in \Zgeo$ such that
    \begin{align}
      \gam=\tp{\tp{0,\frac{\lam}{2},0}},
      \label{eq:ideal-front-pt-gam-gam11-zero}
    \end{align}
    or some $\lam\in \Zgez$ such that
    \begin{align}
      \gam=\tp{\tp{\fraa,\frac{\lam}{2},0}}.
      \label{eq:ideal-front-pt-gam-gam11-half}
    \end{align}
  \end{enumerate}
\end{proposition}
\begin{proof}
  For
  $\theta=\tp{1,\infty}$,
  let 
  $\pi(1)=\{\tp{s,l,m_{i},n_{i},k_{i}}\}_{i\ldZ \theta}$
  with $a_{i}=\nu(k_{i})$ and
  $b_{i}=\nu(m_{i},n_{i},k_{i})$.
  Suppose $\phi(x)= \tp{1-q}\in \Q(\fX)^{l}$
  and $i\ldZ \theta$. 
  Then, Theorem~\ref{thm:monomial-poly} gives
  \begin{align}
    q^{-2(\gam_{1,1} i^2+\gam_{1,2}i)}
    \Delta(\cF)(s,l,w,m_{i},n_{i},k_{i},\phi,\rho,x,\fX)
    &=
    \dt(V)_{a_{i}}^{b_{i}}(s,l,w,\phi,\rho,
    t_{\gam,\Del}(m_{i},n_{i},k_{i}),x)\\
    &\llq 0.\label{eq:ideal-front-pt-gam-quasi-merged}
  \end{align}
  
  Let us prove Statement~\ref{s:ideal-front-pt-gam-gam} from
  Statement~\ref{s:ideal-front-pt-gam-front-pt}.  Since
  $\Delta(\cF)(s,l,w,m_{i},n_{i},k_{i},\phi,\rho,x,\fX)\llq 0$ for the
  ideal $\xi$,
  inequality~\eqref{eq:ideal-front-pt-gam-quasi-merged}
  implies
  \begin{align}
    2(\gam_{1,1} i^2+\gam_{1,2}i)\in \Z.
    \label{inc:ideal-front-pt-gam-exp}
  \end{align}
  First, assume $\gam_{1,1} =0$. 
  Then,
  inclusion~\eqref{inc:ideal-front-pt-gam-exp} gives
  $2\gam_{1,2}\in \Z$ by $i=1$.  
  Thus,
  equation~\eqref{eq:ideal-front-pt-gam-gam11-zero} holds by
    $t_{\gam}(\tp{1}) =\gam_{1,2}>0$
  in   Claim~\ref{c:prim-monom-pt-tgam-pt} of
  Theorem~\ref{thm:prim-monom-pt}.
  Second, assume
   $\gam_{1,1}=\fraa$. 
   Then, inclusion~\eqref{inc:ideal-front-pt-gam-exp}
   implies $1+2\gam_{1,2}\in \Z$
   by $i=1$.
   Hence,
   equation~\eqref{eq:ideal-front-pt-gam-gam11-half}
   follows
   from $t_{\gam}(\tp{1}) =\fraa+\gam_{1,2}>0$
   in Claim~\ref{c:prim-monom-pt-tgam-pt} of
  Theorem~\ref{thm:prim-monom-pt}.

    Let us prove
    Statement~\ref{s:ideal-front-pt-gam-front-pt} from
    Statement~\ref{s:ideal-front-pt-gam-gam}.  First,
    suppose
    equation~\eqref{eq:ideal-front-pt-gam-gam11-zero}.
    Since $2(\gam_{1,1} i^2+\gam_{1,2}i) =\lam i\in \Zgeo$, $\xi$ is ideal by
    inequality~\eqref{eq:ideal-front-pt-gam-quasi-merged}.
    Thus, Statement~\ref{s:ideal-front-pt-gam-front-pt}
    holds by Claim~\ref{c:prim-monom-pt-tgam-pt} of
    Theorem~\ref{thm:prim-monom-pt}, since
    $t_{\gam}(\tp{\mu}) =\frac{\lam}{2} \mu>0$ for 
    $\mu\in\Zgeo$.  Second, suppose
    equation~\eqref{eq:ideal-front-pt-gam-gam11-half}.
    Since
    $2(\gam_{1,1} i^2+\gam_{1,2}i) =i^{2}+\lam i\in \Zgeo$, $\xi$ is ideal by
inequality~\eqref{eq:ideal-front-pt-gam-quasi-merged}.
Claim~\ref{c:prim-monom-pt-tgam-pt} of
Theorem~\ref{thm:prim-monom-pt}
now gives
Statement~\ref{s:ideal-front-pt-gam-front-pt},
since
    $t_{\gam}(\tp{\mu})
    =\frac{\mu^{2}}{2}+\frac{\lam\mu}{2}>0$
    for  $\mu\in\Zgeo$.
\end{proof}

We introduce the following notation to compare merged pairs
by bases and almost strictly unimodal sequences.

\begin{definition}\label{def:covering}
  For $i\in\oi(2)$, suppose parcels
  $\cF_{i}=
  \Lam(s,l,w,\scc_{i},f_{i,s},\phi_{i},\rho_{i},x,\fX_{i})$ so
  that each $x_{j}\in \Q(\fXo)\cap \Q(\fXt)$.
  \begin{enumerate}
  \item If $r_{1}\in A_{\fXo}$ and $r_{2}\in A_{\fXt}$ satisfy
    $x(r_{1})=x(r_{2})\in \R^{l}$, then we write
    $\tp{\fXo,r_{1}}\equiv^{x}\tp{\fXt,r_{2}}$.
  \item Consider $\tp{\theta,\fX_{i}}$-merged pairs
    $\zeta_{i}=\tp{P,\cF_{i}}$ for $i\in\oi(2)$.  If
    $u(\zeta_{1},r_{1})\geq u(\zeta_{2},r_{2})$ whenever
    $\tp{\fXo,r_{1}}\equiv^{x}\tp{\fXt,r_{2}}$, then we say
    that
    $\zeta_{1}$ covers $\zeta_{2}$ on $x$ and write
    $\zeta_{1}\geq^{x} \zeta_{2}$.
  \end{enumerate}
\end{definition}

Let us state the following reflexivity and transitivity of
the covering relation $\geqx$.

\begin{lemma}\label{lem:covering}
  For $i\in\oi(3)$, consider the
  $\tp{\theta,\fX_{i}}$-merged pairs
  $\zeta_{i}=\tp{P,\cF_{i}}$ of
  $\cF_{i}=
  \Lam(s,l,w,\scc_{i},f_{i,s},\phi_{i},\rho_{i},x,\fX_{i})$.
  Then, we have the following.
  \begin{enumerate}
  \item \label{c:covering-reflexivity}
    $\zetao\geqx \zetao$.
  \item\label{c:covering-transitivity}
    $\zetao\geqx \zetat$
    and $\zetat\geqx \zetar$ imply $\zetao \geqx \zetar$.
  \end{enumerate}
\end{lemma}
\begin{proof}
  Claim~\ref{c:covering-reflexivity} holds by
  $u(\zetao,r_{1})\geq u(\zetao,r_{1})$ in $\R$ for
  $r_{1}\in \AfXo$.  Claim~\ref{c:covering-transitivity}
  follows, since
  $u(\zetao,r_{1})\geq u(\zetat,r_{2}) \geq u(\zetar,r_{3})$ if
  $\tp{\fXo,r_{1}}\equiv^{x}
  \tp{\fXt,r_{2}}\equiv^{x}\tp{\fXr,r_{3}}$.
\end{proof} 

The antisymmetricity of $\geqx$ does
not hold in general.  However, we state the following.

\begin{lemma}\label{lem:covering-uniqueness}
  For $i\in\oi(2)$, consider the
  $\tp{\theta,\fX_{i}}$-merged pairs
  $\zeta_{i}=\Omg(\cF_{i})$ of primal monomial parcels
  $\cF_{i}=
  \Lam(s,l,w,\scc_{i},\Psi_{s,\gam_{i},q},\rho_{i},x,\fX_{i})$.
  If $\zetao\geqx \zetat\geqx \zetao$, then $\cFo=\cFt$.
\end{lemma}
\begin{proof}
  By the primal assumption, $q^{\frac{1}{d_{1}}}\in \fXo$ and
  $q^{\frac{1}{d_{2}}}\in \fXt$ for some
  $d_{1},d_{2}\in \Zgeo$.  Then, the covering relations imply
  $\cF_{1,\tp{i}}(q^{\frac{1}{d_{1}}})
  =\cF_{2,\tp{i}}(q^{\frac{1}{d_{2}}})$ for each
  $i\in \Zgez$ and $0<q<1$.  This implies
  $\cF_{1,\tp{i}}=\cF_{2,\tp{i}}\in
  \Q\left(q^{\frac{1}{d_{1}d_{2}}}\right)$ so that
  $\cF_{1,\tp{i}}=\cF_{2,\tp{i}}$ have infinite solutions.
\end{proof} 

In the following corollary, we identify the linear-half
monomial $\cL$ in Definition~\ref{def:linear-half} as the
extremal parcel of primal monomial parcels by polynomials
with positive integer coefficients and phase transitions.
Also,
we obtain the golden ratio as the critical point of
$\Omg(\cL)$ by the fully optimal coordinate of $\cL$.
Thus, by Theorem~\ref{thm:inv-lim}, it can be seen that we
find the golden ratio of quantum dilogarithms in the theory
of the merged-log-concavity.  Also, the Hadamard product
$\cL^{\dd 2}$ gives the golden angle as the critical point
of the ideal merged pair
$\Omg(r_{\tp{1},\tp{2}}(\cL^{\dd 2}))$ in
Section~\ref{sec:golden-angle}.

\begin{corollary}\label{cor:linear-half-maximum}
  Consider 
 the
 ideal merged pair $\zeta=\Omg(\cL)$
 of the
  primal monomial parcel
  $\cL=\Lam(s,l,w,\scc,\Psi_{s,\gam,q},\rho,x,\fX)$.
  \begin{enumerate}
  \item
    \label{c:linear-half-maximum-maximum}
    For a primal monomial parcel
    $\cF=\Lam(s,l,w,\scc',\Psi_{s,\gam',q},\rho',x,\fX')$,
    suppose
    the ideal merged pair $\zeta'=\Omg(\cF)$ with a front
    phase transition.  Then, $\zeta\geqx \zeta'$.
  \item
    \label{c:linear-half-maximum-unique}
    Among all primal monomial parcels
    $\Lam(s,l,w,\scc',\Psi_{s,\gam',q},\rho',x,\fX')$, $\cL$
    gives the unique maximum ideal merged pair
    $\zeta$ on the
    covering relation $\geqx$.
  \item
    \label{c:linear-half-maximum-crit}
    The single critical point
    of $\zeta$ is
    $\FC(\zeta)=
      \frac{-1+\sqrt{5}}{2}$, which is
      the golden ratio.
  \end{enumerate}
\end{corollary}
\begin{proof}
  Let us prove Claim~\ref{c:linear-half-maximum-maximum}.
  Suppose
  $\tp{\fX,r}\equiv^{x}\tp{\fX',r'}$,
  which gives
  $0<q(r)=q(r')<1$
  in Claim~\ref{c:adm-bounds}
  of Lemma~\ref{lem:adm-bounds-poring-inclusions}.
  Moreover,
  Proposition~\ref{prop:ideal-front-pt-gam} implies
  $\gam'=\tp{\tp{0,\frac{\lam}{2},0}}$
  for some $\lam\in \Zgeo$, 
  or
  $\gam'=\tp{\tp{\fraa,\frac{\lam}{2},0}}$
  for some $\lam\in \Zgez$.
  Thus,
  $t_{\gam}(\tp{i})=\frac{i}{2}\leq t_{\gam'}(\tp{i})$
  for each $\gam'$ and $i\in \Zgez$.
  In particular,
  \begin{align}
    u(\zeta,r)_{i}
    =
    \frac{q(r)^{t_{\gam}(\tp{i})}}{(i)_{q}|_{q=q(r)}}
    \geq
    \frac{q(r')^{t_{\gam'}(\tp{i})}}{(i)_{q}|_{q=q(r')}}
    =u(\zeta',r')_{i}.
  \end{align}

  Claim~\ref{c:linear-half-maximum-unique} follows from
  Lemma~\ref{lem:covering-uniqueness} and
  Claim~\ref{c:linear-half-maximum-maximum}.
  
  Let us prove Claim~\ref{c:linear-half-maximum-crit}.  A
  real number $0<\qq<1$ is a front critical point of $\zeta$
  if and only if
  $\cL_{\tp{0}}=1=\frac{\qq}{1-q}=\cL_{\tp{1}}$,
  whose solution
  is
  $\FC(\zeta) =\frac{-1+\sqrt{5}}{2}$.
  Since $s$ is infinite and $\cL$ is
  vanishing by
  Claim~\ref{c:prim-monom-pt-vanishing-pcl-pt} of
  Theorem~\ref{thm:prim-monom-pt},
   $\zeta$ has no
  other critical points. 
\end{proof}

We put the following  parcel $\cQ$,
which is in
Section~\ref{sec:intro-merged-rational-functions} by a
different notation.

\begin{definition}\label{def:quadratic-half}
  Let $s=\tp{0,\infty}$, $l=1$, $w=\tp{1}$, and
  $\fX=\{\qq\}$.  We
  define the quadratic-half monomial parcel
  $\cQ=\Lam(s,l,w,\scc,
  \Psi_{s,\tp{\tp{\fraa,0,0}}},x,\fX)$.
\end{definition}                     

By quadratic primal monomial parcels, $\cQ$ has the
following analog of Corollary~\ref{cor:linear-half-maximum}.

\begin{corollary}
  \label{cor:quadratic-half-maximum-maximum}
  Consider the primal monomial parcel
  $\cQ=\Lam(s,l,w,\scc,\Psi_{s,\gam,q},\rho,x,\fX)$ and the
  ideal merged pair $\zeta=\Omg(\cQ)$.
  \begin{enumerate}
  \item
    \label{c:quadratic-half-maximum-maximum}
    For 
    a quadratic primal monomial parcel
    $\cF=\Lam(s,l,w,\scc',\Psi_{s,\gam',q},\rho',x,\fX')$,
    suppose
    the  ideal merged pair
    $\zeta'=\Omg(\cF)$
    with a front phase transition.
    Then,
    $\zeta\geqx \zeta'$.    
  \item
    \label{c:quadratic-half-maximum-unique}
    Among all quadratic primal monomial parcels
    $\cF=\Lam(s,l,w,\scc',\Psi_{s,\gam',q},\rho',x,\fX')$,
    $\cQ$
     gives the  unique
    maximum ideal merged pair $\zeta$ on the covering relation
    $\geqx$.
  \item
    \label{c:quadratic-half-maximum-crit}
    The single critical point
    of $\zeta$ is
    $\FC(\zeta)=
      \frac{-1+\sqrt{5}}{2}$.
  \end{enumerate}
\end{corollary}
\begin{proof}
Claims
  hold as in
  Corollary~\ref{cor:linear-half-maximum},
  since
 $\cQ_{\tp{1}}=\frac{\qq}{1-q}$ and
  $\gam'_{1,1}=\fraa$ by the monomial conditions of
  $\tp{l,w,\gam'}$. 
\end{proof}

\subsection{Convolutions of vanishing monomial parcels and
  phase transitions}
\label{sec:conv-prim}

Theorem~\ref{thm:prim-monom-pt} gives vanishing
monomial parcels $\cF$ with phase transitions.
 We obtain more vanishing parcels
 with phase transitions by convolutions.
 We first state the following,
since convolutions of vanishing sequences
are not necessarily vanishing.

\begin{lemma}
  \label{lem:conv-prob}
  Let $\del\in \Zgeo$.  For each $i\in\oi(\del)$, assume a
  primal monomial parcel
  $\cF_{i}=\Lam(s,l,w,\scc,\Psi_{s,\gam_{i},q},\rho,x, X)$ such
  that $\Omg(\cF_{i})$ is vanishing.  Consider the parcel
  convolution $\cH=*_{i\in\oi(\del)}\cF_{i}$.  Then,
  $\Omg(\cH)$ is vanishing.
\end{lemma}
\begin{proof}
  Let $r\in \AfX$ and $\mu\in \Zgez$.  Then,
  Proposition~\ref{prop:van-prob-tgam} gives
  $N_{i}(q(r))\geq 1$ and $0<S_{i}(q(r))<1$ such that
  $\cF_{\tp{\mu}}(r)\leq N_{i}(q(r)) \cdot S_{i}(q(r))^{\mu}$ for
  each $i\in \oi(\del)$.  Consider real numbers
  $\til{N}(q(r))$ and $\til{S}(q(r))$ such that
  $\til{N}(q(r))\geq N_{i}(q(r))$ and
  $1>\til{S}(q(r))\geq S_{i}(q(r))>0$ for each
  $i\in \oi(\del)$.  The assertion holds by
  $\cH_{\tp{\mu}}(r)\leq \til{N}(q(r))^{\del} \cdot (\mu+1)^{\del-1}
  \cdot \til{S}(q(r))^{\mu}$.
\end{proof}
Then,  the following
monomial indices  give
not only monomial
convolutions of eta products in
Proposition~\ref{prop:lifted-eta-pro},
but also the convolutions with front critical
points.

\begin{proposition}
  \label{prop:conv-prob-pt}
  Suppose integers $0\leq d_{1}\leq d_{2}$ such that
  $d_{2}\in \Zgeo$.  Let $l=1$ and $w=\tp{1}$.  Suppose
  monomial indices $\tp{l,w,\gam_{i}}$ for $i\in\oi(d_{2})$
  such that $i\in\oi(d_{1})$ gives
  $\gam_{i}= \tp{0,\gam_{i,2},0}$ with
  \begin{align}
    \gam_{i,2}&> 0, \label{ineq:conv-prob-pt-linear-gam}
  \end{align}
  and  $i\in\oi(d_{1}+1,d_{2})$ gives
  $\gam_{i}= \tp{\fraa,\gam_{i,2},0}$  with
  \begin{align}
    \gam_{i,2}
    &> -\fraa.
      \label{ineq:conv-prob-pt-quad-linear-gam}
  \end{align} 
  For  $i\in\oi(d_{2})$,
  let $\cF_{i}=\Lam(s,l,w,\scc,\Psi_{s,\gam_{i},q},\rho,x,\fX)$
  be primal monomial parcels
   such that $\fX=\{q^{\ka}\}$
  of some $\ka^{-1}\in \Zgeo$.
  Moreover, consider the parcel
  convolution
  $\cH=*_{i\in\oi(d_{2})}\cF_{i}=\Lam(s,l,w,\scc,\hs,\rho,x,\fX)$
  and the merged pair $\zeta=\Omg(\cH)$.  Then, we have
  the following.
  \begin{enumerate}
  \item
    \label{c:conv-prob-pt-pt}
    $\zeta$ has the unique phase transition at the front
    critical point $0<\FC(\zeta)=q^{\ka}<1$
    that solves
    $1-q=\sum_{i\in\oi(d_{1})}
      q^{\gam_{i,2}}
      +\sum_{i\in\oi(d_{1}+1,d_{2})}q^{\fraa+\gam_{i,2}}$.
  \item
    \label{c:conv-prob-pt-no-rear-no-front-crit}
    $\zeta$ has neither rear nor front critical points, if one
    of
    inequalities~\eqref{ineq:conv-prob-pt-linear-gam}
    and~\eqref{ineq:conv-prob-pt-quad-linear-gam}
     fails.
  \end{enumerate}
\end{proposition}
\begin{proof}
  Let us prove Claim~\ref{c:conv-prob-pt-pt}.  Since $s$ is
  infinite, $\zeta$ has no rear critical points.
  Moreover, each
  $\cF_{i}$ satisfies
  Statement~\ref{s:prim-monom-pt-vanishing-pcl-pt-tgam} of
  Theorem~\ref{thm:prim-monom-pt}. Hence, $\zeta$ has no asymptotic
  critical points by
  Claim~\ref{c:prim-monom-pt-vanishing-pcl-pt}of
  Theorem~\ref{thm:prim-monom-pt} and
  Lemma~\ref{lem:conv-prob}.
  
  Let us prove that $\zeta$ has the single front critical
  point.
  First, 
  $p_{1}(q^{\ka})=1-q$
  strictly decreases
  over
  $0\leq q^{\ka}\leq 1$ with
   $p_{1}(0)=1$ and
  $p_{1}(1)=0$.
  Second,
  inequalities~\eqref{ineq:conv-prob-pt-linear-gam}
  and~\eqref{ineq:conv-prob-pt-quad-linear-gam}
  imply that
  $p_{2}(q^{\ka})= \sum_{i\in\oi(d_{1})}q^{\gam_{i,2}}
  +\sum_{i\in\oi(d_{1}+1,d_{2})}q^{\fraa+\gam_{i,2}}$
  strictly increases over
  $0\leq q^{\ka}\leq 1$ with
   $ p_{2}(0)=0$ and $ p_{2}(1)=d_{2}\geq 1$.
    Hence, there is the
  unique solution $0<\FC(\zeta)=q^{\ka}<1$ for
  $1=\cH_{\tp{0}}
  = \cH_{\tp{1}}=\frac{ p_{2}(q^{\ka}) }{p_{1}(q^{\ka})}$
  such that
  $p_{2}(\FC(\zeta))=p_{1}(\FC(\zeta))$.  Also, there are
  $r_{1},r_{2}\in \AfX$ such that
  $\cH_{\tp{0}}(r_{1})< \cH_{\tp{1}}(r_{1})$ and
  $\cH_{\tp{0}}(r_{2})> \cH_{\tp{1}}(r_{2})$.
  
  Let us confirm
  Claim\ref{c:conv-prob-pt-no-rear-no-front-crit}.  If one
  of inequalities~\eqref{ineq:conv-prob-pt-linear-gam}
  and~\eqref{ineq:conv-prob-pt-quad-linear-gam} fails, then
  we have no front critical points as $p_{2}(q^{\ka})>1$ for
  $0<q^{\ka}<1$.  We have no rear critical points either for
  the infinite $s$.
\end{proof}
In particular, Proposition~\ref{prop:conv-prob-pt} determines front
phase transitions of $\Omg(\cH)$ for
each convolution $\cH$ of primal monomial parcels.
Moreover, the following gives explicit front critical points
by metallic ratios.

\begin{corollary}\label{cor:conv-prob-linear-half-quad-half}
  Consider the linear-half
  $\cL=\Lam(s,l,w,\scc, \Psi_{s,\gam,q},\rho,x,\fX)$ and
  quadratic-half
  $\cQ=\Lam(s,l,w,\scc,\Psi_{s,\gam',q},\rho,x,\fX)$
  of $\fX=\{\qq\}$.
    Let
  $n\in \Zgeo$. Suppose merged pairs
  $\zeta_{n}=\Omg(\cL^{*n})$ and
  $\til{\zeta}_{n}=\Omg(\cQ^{*n})$.  Then, we have the
  following.
  \begin{enumerate}
  \item \label{c:conv-prob-linear-half-quad-half-fully-opt}
    $\fX$ is fully optimal for $\cL^{*n}$ and
    $\cQ^{*n}$.
  \item \label{cor:conv-prob-linear-half-quad-half-ideal-pt}
    $\zeta_{n}$ is ideal and has the unique phase transition
    at the front critical point
    $\FC(\zeta_{n})=\frac{-n+\sqrt{n^2+4}}{2}$,
    which is a metallic ratio.
  \item \label{cor:conv-prob-linear-half-quad-half-quad-half}
    The same holds for $\til{\zeta}_{n}$ with the front critical
    point $\FC(\tilde{\zeta}_{n})= \FC(\zeta_{n})$.
  \end{enumerate}
\end{corollary}
\begin{proof}
  Claim~\ref{c:conv-prob-linear-half-quad-half-fully-opt}
  holds by $\cL^{*n}_{\tp{1}}=\frac{n\qq}{1-q}$.
  Furthermore, $\zeta_{n}$ is ideal by
  Corollary~\ref{cor:conv} and
  Proposition~\ref{prop:ideal-front-pt-gam}.
Thus,
  Claim~\ref{cor:conv-prob-linear-half-quad-half-ideal-pt}
  follows, since $\zeta_{n}$ has the front phase transition
  with the critical point $\FC(\zeta_{n})$  
  by
  Proposition~\ref{prop:conv-prob-pt}.
  Claims~\ref{cor:conv-prob-linear-half-quad-half-quad-half}
  holds similarly, as
  $\cQ^{*n}_{\tp{1}}=\cL^{*n}_{\tp{1}}$.
\end{proof}

\section{Monomial convolutions and
  ideal boson-fermion gases}
\label{sec:casimir}
Generalizing Section~\ref{sec:intro-stat-pt},
Section~\ref{sec:casimir} considers
some ideal boson-fermion gases with or without Casimir energies (Ramanujan summation) by
monomial convolutions. Thus, we obtain
statistical-mechanical phase transitions by the mathematical
merged-log-concavity.

Unless stated otherwise, Section~\ref{sec:casimir} assumes
the following.
Let $s=\tp{0,\infty}$, $l=1$, $w=\tp{1}$, and
$v\in \Q$. Also, as in
Section~\ref{sec:intro-stat-pt},
let $q=e^{-\bta}$ by the thermodynamic beta
$\bta>0$
and $t=e^{-\mu'}$ by
$\mu'=-\mu\bta$ of the chemical potential
$\mu<0$.

\subsection{Monomial convolutions
  without Casimir energies}\label{sec:eta-mixed}

We  look at the following systems of ideal boson,
fermion, and boson-fermion gases without Casimir energies by
monomial convolutions.

\subsubsection{Ideal boson  gases}

From Definition~\ref{def:boson-hamil-numb},
the bosonic  operators $a_{b,\lam}, a_{b,\lam}^{\dagger}$
of $\lam\in \Zgeo$ 
give  
the bosonic system $ B(1,v)$ of the Hamiltonian $H_{b,v}$
and the number operator $N_{b}$.  Also, $B(1,v)$
has the grand canonical partition function
\begin{align}
  \cZ_{B(1,v)}(q,t)=\Tr{e^{-\bta H_{b,v}}  \cdot e^{-\mu' N_{b}}}.
\end{align}
Moreover, we have the monomial convolution
$\cZ_{w,\tp{\tp{0,1-v,0}},q}(t)$ of the primal monomial
parcel
$\Lam(s,l,w,\scc,\Psi_{s,\tp{\tp{0,1-v,0}},q},x,\fX)$.  Also,
$(1-q^{\lam-v}t)^{-1}= \sum_{n_{\lam}\in\Zgez}
e^{- \bta n_{\lam} \eps_{v,\lam}} \cdot
e^{-\mu' n_{\lam}}$ and
$(q^{1-v}t;q)^{-1}_{\infty}=\sum_{r\in \Zgez}
\frac{q^{r(1-v)}} {(r)_{q}}t^{r}$ by Euler's linear
binomial identity.  Thus,  
as in Section~\ref{sec:intro-ideal-boson-gas},
we obtain
\begin{align}
  \cZ_{B(1,v)}(q,t)=\cZ_{w,\tp{\tp{0,1-v,0}},q}(t).
  \label{eq:ideal-boson-partition-monom-conv}
\end{align}

\subsubsection{Ideal fermion gases}         
Similarly from
Definition~\ref{def:ferm-hamil-numb},
the operators
$a_{f,\lam}, a_{f,\lam}^{\dagger}$
 give the fermionic system
$F(1,v)$ of  the Hamiltonian $H_{f,v}$ and the
number operator $N_{f}$.  Thus, $F(1,v)$
carries the grand canonical partition function
\begin{align}
  \cZ_{F(1,v)}(q,t)=
  \Tr{e^{-\bta H_{b,v}}\cdot e^{-\mu' N_{b}}}.
\end{align}
Moreover, there is the monomial convolution
$\cZ_{w,\tp{\tp{\fraa,\fraa-v,0}},q}(t)$ of
the primal monomial parcel
$\Lam(s,l,w,\scc,
\Psi_{s,\tp{\tp{\fraa,\fraa-v,0}},q},x,\fX)$.
Also,
$1+q^{\lam-v}t=
\sum_{n_{\lam}}e^{-\bta n_{\lam}
   \eps_{v,\lam}}\cdot
 e^{-\mu' n_{\lam}}
$ and
 $(-q^{1-v}t;q)_{\infty}= \sum_{r\in\Zgez}\frac{q^{\frac{r^{2}}{2}+
    r\left(\fraa-v\right)}}{(r)_{q}}t^{r}$ 
by
Euler's quadratic binomial identity.
Hence,
as in Section~\ref{sec:intro-ideal-fermion-gas}, 
we have
\begin{align}
  \cZ_{F(1,v)}(q,t)
  =\cZ_{w,\tp{\tp{\fraa,\fraa-v,0}},q}(t).
  \label{eq:ideal-fermion-partition-monom-conv}
\end{align}

\subsubsection{Ideal boson-fermion
  gases}
Suppose integers $0\leq d_{1}\leq d_{2}$
such that $d_{2}\in \Zgeo$ and
$\ka\in \Q^{d_{2}}$.  Consider the boson-fermion system
$M_{d_{1},d_{2}}(\ka)$ that has sub-systems
$B(1,\ka_{\lam})$ for $\lam\in\oi(d_{1})$ and
sub-systems $F(1,\ka_{\lam})$ for
$\lam\in\oi(d_{1}+1,d_{2})$ with negligible interactions
among sub-systems.  Thus,
$M_{d_{1},d_{2}}(\ka)$ has the
grand canonical partition function
\begin{dmath*}
  \cZ_{M_{d_{1},d_{2}}(\ka)}(q,t)
  =
  \prod_{\lam\in\oi(d_{1})}
  \Tr{e^{-\bta H_{b,\ka_{\lam}}}
    \cdot e^{-\mu' N_{b}}}
  \cdot \prod_{\lam\in\oi(d_{1}+1,d_{2})}
  \Tr{e^{-\bta H_{f,
        \ka_{\lam}}}  \cdot e^{-\mu' N_{f}}}.
\end{dmath*}

Moreover,  $\alp=\iota^{d_{2}}(1)$ gives the
multimonomial index $\tp{d_{2},w,\alp,\gam_{\ka}}$ such that
\begin{align}
  \gam_{\ka,\lam}=
  \begin{dcases}
    & \tp{0,1-\ka_{\lam},0} \mforeach
      \lam\in\oi(d_{1}),\\
    &\tp{\fraa,\fraa-\ka_{\lam},0}
      \mforeach
      \lam\in\oi(d_{1}+1,d_{2}).
  \end{dcases}
\end{align}
Then,
by  equations~\eqref{eq:ideal-boson-partition-monom-conv} 
and~\eqref{eq:ideal-fermion-partition-monom-conv},
we realize
\begin{align}
  \cZ_{M_{d_{1},d_{2}}(\ka)}(q,t)=
  \cM(d_{2},w,\alp,\gam_{\ka},q,t),
\end{align}
 which is a monomial convolution.

\subsection{Monomial convolutions with Casimir
  energies}
\label{sec:casimir-q}

Sections~\ref{sec:intro-stat-pt} and~\ref{sec:eta-mixed}
ignore {\it zero-point energies}. This is often done in
statistical mechanics~\cite{KapGal}.  However, let us
attempt not to ignore.  In particular, by monomial
convolutions, we consider the Casimir energies as the
following Ramanujan summation of the divergent zero-point
energy sums of $B(1,v)$ and $F(1,v)$.

Suppose an entire function $a(y)$ of $y\in \bC$ such that the
exponential type of $a(y)$ is less than $\pi$. Then, by
the Bernoulli numbers $B_{\lam+1}$, we recall the Ramanujan
summation~\cite{Can}:
\begin{align}
  \sum^{\cR}_{\lam\in\Zgeo}a(\lam)
  =\int_{0}^{1}a(y)dy-\fraa a(0)
  -\sum_{\lam\in\Zgeo}a^{(\lam)}(0)
  \frac{B_{\lam+1}}{(\lam+1)!}.
  \label{eq:rama-sum}
\end{align}
For example, $a_{v}(y)=y-v$ gives
$\sum^{\cR}_{\lam\in\Zgeo}a_{v}(\lam)
=\frac{5}{12}-\frac{v}{2}$, which is
$-\frac{1}{12}$~\cite{Ram} when $v=1$.

\subsubsection{Ideal   boson gases}
\label{sec:prim-monom-ideal-boson-casimir}

Let us take the following operators
with zero-point energies.
\begin{definition}
  \label{def:boson-hamil-numb-casimir}
  For $v\in \Q$, let $\cH_{b,v}$ denote the Hamiltonian
  operator such that
  \begin{align}
    \cH_{b,v}
    &=\sum_{\lam\in\Zgeo}
    \frac{\eps_{v,\lam}}{2}
    \{a_{b,\lam}^{\dagger},a_{b,\lam}\}\\
    &=\sum_{\lam\in\Zgeo}
    \left(\eps_{v,\lam}
    a_{b,\lam}^{\dagger}a_{b,\lam}
    +\frac{\eps_{v,\lam}}{2}\right).
    \label{eq:boson-hamil-zero-inf}
  \end{align}
  Also, if $u\in \Zgeo$, then let
  $\cH_{b,v,u}$ and $N_{b,u}$
  denote the Hamiltonian and number operators
  such that
  \begin{align}
    \cH_{b,v,u}
    &=\sum_{\lam\in\oi(u)}
      \left(
      \eps_{v,\lam}
        a_{b,\lam}^{\dagger}a_{b,\lam}
    +\frac{  \eps_{v,\lam}}{2}
    \right),\\
    N_{b,u} &=\sum_{\lam\in\oi(u)}  a_{b,\lam}^{\dagger}a_{b,\lam}.
  \end{align}
\end{definition}

Consider the bosonic system $\cB(1,v)$ of $\cH_{b,v}$ and
$N_{b}$ in Definition~\ref{def:boson-hamil-numb}.  Unlike in
equation~\eqref{eq:hamil-bv}, $\frac{  \eps_{v,\lam}}{2}$ in
equation~\eqref{eq:boson-hamil-zero-inf} represent
zero-point energies.  Thus, $\cB(1,v)$ extends $B(1,v)$ by
the zero-point energies. Also, $\cB(1,v)$ proposes the grand
canonical partition function
\begin{align}
  \cZ_{\cB(1,v)}(q,t)=
  \Tr{e^{-\bta \cH_{b,v}}  \cdot e^{-\mu' N_{b}}}.
\end{align}
However, $\cZ_{\cB(1,v)}(q,t)$ has 
$q^{\frac{\sum_{1\leq \lam}  \eps_{v,\lam}}{2}}$ of the divergent
zero-point energy sum $\frac{\sum_{1\leq \lam}  \eps_{v,\lam}}{2}$. This
would make $\cZ_{\cB(1,v)}(q,t)=0$ by $0<q<1$.  Thus, we
consider the  regularized grand canonical partition function
$\cZ^{\cR}_{\cB(1,v)}(q,t)$ of $\cB(1,v)$ by the Ramanujan
summation and a monomial convolution as follows.

For each $u\in \Zgeo$, there exists the bosonic system
$\cB(1,v,u)$ of
$\cH_{b,v,u}$ and
$N_{b,u}$.
Also, for each
$\lam\in \Zgeo$ and $n_{\lam}\in\Zgez$, we have
$\langle n_{\lam}| e^{-\bta\eps_{v,\lam}
  \left(a^{\dagger}_{\lam}\alp_{\lam}
  +\fraa\right)}\cdot e^{-\mu' a^{\dagger}_{\lam}\alp_{\lam}}|n_{\lam}\rangle
  =
  e^{-\bta\eps_{v,\lam}n_{\lam}}\cdot e^{-\mu' n_{\lam}}
  \cdot e^{-\bta \frac{\eps_{v,\lam}}{2}}$.
  Thus,
  by $\prod_{\lam \in \oi(u)}e^{-\bta \frac{\eps_{v,\lam}}{2}}
  =q^{\frac{\sum_{\lam\in \oi(u)}\eps_{v,\lam}}{2}}$,
$\cB(1,v,u)$ has the grand canonical partition function
\begin{align}
  \cZ_{\cB(1,v,u)}(q,t)=
  \Tr{e^{-\bta \cH_{b,v,u}}  \cdot e^{-\mu' N_{b,u}}}
  =(q^{1-v}t;q)_{u}^{-1}
    \cdot
    q^{\sum_{\lam\in\oi(u)}\frac{\eps_{v,\lam}}{2}}.
    \label{eq:casimir-boson-finite}
\end{align}

In equations~\eqref{eq:casimir-boson-finite}, $u\to\infty$ yields
the divergent zero-point energy sum
\begin{align}
  \lim_{u\to\infty}\fraa
  \sum_{\lam\in\oi(u)}\eps_{v,\lam}.
  \label{eq:ram-cas-ferm}
\end{align}  
But, the Ramanujan summation~\eqref{eq:rama-sum} gives the
following regularization (c.f. \cite[Section 1.3]{Pol} for
$v=1$):
\begin{align}
  \lim_{u\to\infty}\fraa
  \sum_{\lam\in\oi(u)}\eps_{v,\lam}
  \rightarrow &
       \fraa\sum_{\lam\in\Zgeo}^{\cR}
       \eps_{v,\lam}
       =\frac{5}{24}-\frac{v}{4},
\end{align}
which is  the Casimir
energy
of $\cB(1,v)$.
Thus,
we have the regularized grand canonical partition
function $\cZ^{\cR}_{\cB(1,v)}(q,t)$ of $\cB(1,v)$
 such that
 \begin{align}
   \cZ^{\cR}_{\cB(1,v)}(q,t)
   :=\cZ_{B(1,v)}(q,t)
   q^{\fraa\sum_{\lam\in\Zgeo}^{\cR}
   \eps_{v,\lam}}
  =
   \cZ_{w,\tp{\tp{0,1-v,
  \frac{5}{24}-\frac{v}{4}}},q}(t)
  \label{eq:casimir-boson-partition}
 \end{align}
by
the monomial convolution in equation~\eqref{eq:ideal-boson-partition-monom-conv}.
In particular, this
$\cZ^{\cR}_{\cB(1,v)}(q,t)$ is non-trivial.

\subsubsection{Ideal fermion gases}
\label{sec:pifc}

Let us take the following operators
with zero-point energies.
\begin{definition}\label{def:fermion-hamil-numb-casimir}
  For $v\in \Q$,
  let $\cH_{f,v}$
  denote the Hamiltonian operator such
    that
    \begin{align}
      \cH_{f,v}&=\sum_{\lam\in\Zgeo}
      \frac{\eps_{v,\lam}}{2}   [a^{\dagger}_{f,\lam},a_{f,\lam}]\\
      &=\sum_{\lam\in\Zgeo}
      \left(
      \eps_{v,\lam}a_{f,\lam}^{\dagger}a_{f,\lam}
      -\frac{\eps_{v,\lam}}{2}
      \right).
      \label{eq:fermion-hamil-zero-inf}
    \end{align}
    If $u\in \Zgeo$, then let $\cH_{f,v,u}$ and
    $N_{f,u}$ be the
   Hamiltonian and number operators such that
    \begin{align}
      \cH_{f,v,u}
      &=\sum_{\lam\in\oi(u)}
        \left(
        \eps_{v,\lam}a_{f,\lam}^{\dagger}a_{f,\lam}
        -\frac{\eps_{v,\lam}}{2}
        \right),\\
      N_{f,u}
      &=\sum_{\lam\in\oi(u)}  a_{f,\lam}^{\dagger}a_{f,\lam}.
    \end{align}
  \end{definition}

  Consider the fermionic system $\cF(1,v)$ of
  $\cH_{f,v}$ and
  $N_{b}$, which extends $F(1,v)$ by the
  zero-point energies $-\frac{\eps_{v,\lam}}{2}$ in
equation~\eqref{eq:fermion-hamil-zero-inf}.  Also,
$\cF(1,v)$ suggests the grand canonical partition function
\begin{align}
  \cZ_{\cF(1,v)}(q,t)=
  \Tr{e^{-\bta \cH_{f,v}}  \cdot e^{-\mu' N_{f}}}.
\end{align}
However,
this has
 $q^{-\frac{\sum_{\lam\in \Zgeo}\eps_{v,\lam}}{2}}$, which
 would make $\frac{1}{0}$.  Hence, we consider the regularized grand canonical
partition function $\cZ^{\cR}_{\cF(1,v)}(q,t)$ of $\cF(1,v)$
by the  Ramanujan summation and a monomial convolution as follows.

For each $u\in \Zgeo$, there is the fermionic system
$\cF(1,v,u)$ of $\cH_{f,v,u}$ and
$N_{f,u}$.
Also, $\lam\in \Zgeo$ and  $n_{\lam}\in \{0,1\}$  imply
$\langle n_{\lam}| e^{-\bta\eps_{v,\lam}
  \left(a^{\dagger}_{\lam}\alp_{\lam}
  -\fraa\right)}\cdot
  e^{-\mu' a^{\dagger}_{\lam}\alp_{\lam}}|n_{\lam}\rangle 
  =
  e^{-\bta\eps_{v,\lam}n_{\lam}}
  \cdot e^{-\mu' n_{\lam}}\cdot e^{-\bta \cdot \frac{-\eps_{v,\lam}}{2}}$.
Then,
$\cF(1,v,u)$ has the grand canonical partition function
\begin{align}
  \cZ_{\cF(1,v,u)}(q,t)=
  \Tr{e^{-\bta \cH_{f,v,u}}  \cdot e^{-\mu' N_{f,u}}}
    =(-q^{1-v}t;q)_{u}
    \cdot
    q^{\frac{-1}{2}
    \sum_{\lam\in\oi(u)}\eps_{v,\lam}
    }.
\end{align}
Thus, the Ramanujan summation \eqref{eq:rama-sum} gives
the
regularization
\begin{align}
  \lim_{u\to\infty}-\fraa
  \sum_{\lam\in\oi(u)}\eps_{v,\lam}
  \rightarrow &
       -\fraa\sum_{\lam\in\Zgeo}^{\cR}
       \eps_{v,\lam}
       =-\frac{5}{24}+ \frac{v}{4},
\end{align}
which is  the Casimir energy
of $\cF(1,v)$.
Hence,
we have the regularized
grand canonical partition function
$\cZ^{\cR}_{\cF(1,v)}(q,t)$
of $\cF(1,v)$ such that
\begin{align}
  \cZ^{\cR}_{\cF(1,v)}(q,t)
  :=\cZ_{F(1,v)}(q,t)
  q^{-\fraa\sum_{\lam\in\Zgeo}^{\cR}
  \eps_{v,\lam}}
  =
  \cZ_{w,\tp{\tp{\fraa,\fraa-v, -\frac{5}{24}+
  \frac{v}{4}}},q}(t)
  \label{eq:casimir-fermion-partition}
\end{align}
by
the monomial convolution in equation~\eqref{eq:ideal-fermion-partition-monom-conv}.
In particular, this
$\cZ^{\cR}_{\cF(1,v)}(q,t)$
is mathematically defined.

\subsubsection{Ideal boson-fermion gases}
\label{sec:eta-mixed-casimir}
Suppose integers $0\leq d_{1}\leq d_{2}$
such that $d_{2}\in \Zgeo$ and
$\ka\in \Q^{d_{2}}$.  Consider the boson-fermion system
$\cM_{d_{1},d_{2}}(\ka)$ that has sub-systems
$\cB(1,\ka_{\lam})$ for $\lam\in\oi(d_{1})$ and sub-systems
$\cF(1,\ka_{\lam})$ for $\lam\in\oi(d_{1}+1,d_{2})$ with
negligible interactions among sub-systems.
Also,
$\cM_{d_{1},d_{2}}(\ka)$ proposes the grand canonical
partition function
(see Remark~\ref{rmk:ignore}):
\begin{dmath*}
  \cZ_{\cM_{d_{1},d_{2}}(\ka)}(q,t)
  =
  \prod_{\lam\in\oi(d_{1})}
  \Tr{e^{-\bta \cH_{b,\ka_{\lam}}}
    \cdot e^{-\mu' N_{b}}}
  \cdot \prod_{\lam\in\oi(d_{1}+1,d_{2})}
  \Tr{e^{-\bta \cH_{f,
        \ka_{\lam}}}  \cdot e^{-\mu' N_{f}}},
\end{dmath*}
which have divergent zero-point energy sums.  However,
equations~\eqref{eq:casimir-boson-partition}
and~\eqref{eq:casimir-fermion-partition} give the
regularized grand canonical partition function
$\cZ^{\cR}_{\cM_{d_{1},d_{2}}(\ka)}(q,t)$ of
$\cM_{d_{1},d_{2}}(\ka)$ such that
\begin{align}
  \cZ^{\cR}_{\cM_{d_{1},d_{2}}(\ka)}(q,t)
  :=
  \prod_{\lam\in\oi(d_{1})}
  \cZ^{\cR}_{\cB(1,\ka_{\lam})}(q,t)
  \cdot \prod_{\lam\in\oi(d_{1}+1,d_{2})}
  \cZ^{\cR}_{\cF(1,\ka_{\lam})}(q,t),
\end{align}
which is mathematically defined and non-trivial. 

Furthermore,  $\alp=\iota^{d_{2}}(1)$ yields the
multimonomial index $\tp{d_{2},w,\alp,\gam^{C}_{\ka}}$ such
that
\begin{align}
  \gam^{C}_{\ka,\lam}=
  \begin{dcases}
    &\tp{0,1-\ka_{\lam},
      \frac{5}{24}-\frac{\ka_{\lam}}{4}}
    \mif \lam\in\oi(d_{1}),\\
    & \tp{\fraa,\fraa-\ka_{\lam},
      -\frac{5}{24}+\frac{\ka_{\lam}}{4}}
    \mif
    \lam\in\oi(d_{1}+1,d_{2}).
  \end{dcases}
\end{align}
Hence,
by equations~\eqref{eq:casimir-boson-partition}
and~\eqref{eq:casimir-fermion-partition},
we obtain
\begin{align}
  \cZ^{\cR}_{\cM_{d_{1},d_{2}}(\ka)}(q,t)=
  \cM(d_{2},w,\alp,\gam^{C}_{\ka},q,t),
  \label{eq:casimir-boson-fermion-partition-monom-conv}
\end{align}  
which is a monomial convolution. 
Moreover, $\cM_{d_{1},d_{2}}(\ka)$ has the Casimir
energy
\begin{dmath*}
  \sum_{\lam\in\oi(d_{1})}\left(-\frac{5}{24}
    +\frac{\ka_{\lam}}{4}\right)
  +
  \sum_{\lam\in\oi(d_{1}+1,d_{2})}
  \left(\frac{5}{24}-\frac{\ka_{\lam}}{4}\right)
  =
  -\frac{5(2d_{1}-d_{2})}{24}
  +\frac{\sum\ka(1,d_{1})
    -\sum\ka(d_{1}+1,d_{2})}{4}.
\end{dmath*}
Thus, we obtain 
\begin{align}
  \cZ^{\cR}_{\cM_{d_{1},d_{2}}(\ka)}(q,t) =
  \cZ_{M_{d_{1},d_{2}}(\ka)}(q,t) \cdot
  q^{
    -\frac{5(2d_{1}-d_{2})}{24}
  +\frac{\sum\ka(1,d_{1})
    -\sum\ka(d_{1}+1,d_{2})}{4}
  }.
  \label{eq:explicit-casimir-bf-partition-monom-conv}
\end{align} 

Now,
equations~\eqref{eq:casimir-boson-fermion-partition-monom-conv}
and~\eqref{eq:explicit-casimir-bf-partition-monom-conv} give
the explicit description of
$\cZ^{\cR}_{\cM_{d_{1},d_{2}}(\ka)}(q,t)$ by $q$-multinomial
coefficients in Claim~\ref{c:monom-multinom-explicit} of
Proposition~\ref{prop:monom-multinom} (see
Proposition~\ref{prop:lifted-eta-pro} for eta products, when
$\ka_{\lam}=1$ for each $\lam\in\oi(d_{2})$).

\subsection{Phase transitions of
  ideal boson-fermion  gases}
\label{sec:sta-math}

As in Section~\ref{sec:intro-stat-pt}, the $t$-power series
of $\cZ_{M_{d_{1},d_{2}}(\ka)}(q,t)$ realize the Helmholtz
free energies of the ideal boson-fermion systems
$M_{d_{1},d_{2}}(\ka)$ without the Casimir energies.  Hence,
$M_{d_{1},d_{2}}(\ka)$ obtain some statistical-mechanical
phase transitions by the merged-log-concavity in
Proposition~\ref{prop:conv-prob-pt}.  So do
$\cM_{d_{1},d_{2}}(\ka)$, even with the Casimir energies in
equation~\eqref{eq:explicit-casimir-bf-partition-monom-conv}.
In particular,
inequalities~\eqref{ineq:conv-prob-pt-linear-gam}
and~\eqref{ineq:conv-prob-pt-quad-linear-gam} in
Proposition~\ref{prop:conv-prob-pt} imply $1-v>0$ and
$\fraa-v>-\fraa$ in
equations~\eqref{eq:ideal-boson-partition-monom-conv}
and~\eqref{eq:ideal-fermion-partition-monom-conv}.
Therefore, we have the energy positivity
$\eps_{v,\lam}=\lam-v>0$ of $\lam\in\Zgeo$ in the
Hamiltonians $H_{b,v}$ and $H_{f,v}$ by the phase
transitions.  Moreover, the phase transitions give not only
the variation of unimodal sequences by continuous parameters
and polynomials with positive integer coefficients, but also
non-zero particle vacua as the temperature increases.

\bibliographystyle{amsalpha}

\noindent
So Okada \\
National Institute of Technology, Oyama College  \\
Oyama Tochigi, Japan 323-0806 \\ 
e-mail: so.okada@gmail.com
\end{document}